\documentclass[10pt]{amsbook}

\newtheorem{theorem}{{\bf Theorem}}[section]
\newtheorem{lemma}[theorem]{{\bf Lemma}}           
\newtheorem{cor}[theorem]{{\bf Corollary}}
\newtheorem{prop}[theorem]{{\bf Proposition}}
\theoremstyle{definition}
\newtheorem{definition}[theorem]{{\bf Definition}}
\newtheorem{example}[theorem]{{\bf Example}}

\newcommand{\ntext}{\bf}

\theoremstyle{remark}
\newtheorem{remark}[theorem]{{\bf Remark}}

\numberwithin{equation}{section}

\begin{document}

\author{Hiroshi ISOZAKI}
\address{Institute of Mathematics, 
University of Tsukuba,
Tsukuba, 305-8571, Japan}
\email{isozakih@math.tsukuba.ac.jp}
\author{Yaroslav KURYLEV}
\address{Department of Mathematics, University College London, United Kingdom}
\email{y.kurylev@ucl.ac.uk}

\title[INTRODUCTION TO SPECTRAL THEORY AND INVERSE PROBLEM ON ASYMPTOTICALLY HYPERBOLIC MANIFOLDS]{\Large{INTRODUCTION TO SPECTRAL THEORY AND INVERSE PROBLEM ON ASYMPTOTICALLY HYPERBOLIC MANIFOLDS}}
\date{\today}

\maketitle

\begin{center}
{\Large {\bf Foreword}}
\end{center}


\subsection{Fourier analysis on manifolds}
The Fourier transform on $L^2({\bf R}^n)$ and its inversion formula are well-known :
\begin{equation}
 \widehat f(\xi) = (2\pi)^{-n/2}\int_{{\bf R}^n}
 e^{-ix\cdot\xi}f(x)dx,
 \label{eq:ForewordFouriertransf}
\end{equation}
\begin{equation}
f(x) = (2\pi)^{-n/2}\int_{{\bf R}^n}
 e^{ix\cdot\xi}\widehat f(\xi)d\xi.
 \label{eq:ForewordInversionformula}
\end{equation}
Since $- \Delta e^{ix\cdot\xi} = |\xi|^2e^{ix\cdot\xi}$, $e^{ix\cdot\xi}$ is an eigenfunction  of 
$- \Delta$. Therefore (\ref{eq:ForewordFouriertransf}) and (\ref{eq:ForewordInversionformula}) illustrate the expansion of arbitrary functions in terms of eigenfunctions (more appropriately {\it generalized eigenfunctions} since they do not belong to $L^2({\bf R}^n)$) of the Laplacian. 

There are two directions of development of the above fact. One is quantum mechanics, where the Schr{\"o}dinger operator $H = - \Delta + V(x)$ is the most basic tool to decribe the physical system of atoms or molecules. If $H$ has the continuous spectrum, it is known that there exists a system of generalized eigenfunctions of $H$ which play the same role as $e^{ix\cdot\xi}$. Moreover, by using these generalized eigenfunctions one can define an operator called the scattering matrix or the S-matrix, which is the fundamental object to study the physical properties of quantum mechanical particles through the scattering experiment.

The other direction is the Fourier transform on manifolds, especially on homogeneous spaces of Lie groups, which is a central theme in the representation theory. Hyperbolic manifolds, one of the deepest sources of classical mathematics, appear also in this context. In particular, hyperbolic quotient manifolds by the action of discrete subgroups of $SL(2,{\bf R})$ and the associated S-matrix are important objects in number theory.


\subsection{Perturbation of the continuous spectrum}
The aim of the perturbation theory of continuous spectrum is, given an operator $H_0$ whose spectral property is rather easy to understand, to study the spectral properties of $H_0 + V$, where $V$ is the perturbation deforming the operator $H_0$. When $H = H_0 + V$ has the continuous spectrum, an effective way of studying its spectral properties is to construct a generalized Fourier tranform associated with $H$. To accomplish this idea, it is necessary that the Fourier transform for $H_0$ can be constructed easily. For example, it is the case for the Laplacian $- \Delta$ on ${\bf R}^n$. If the perturbation term $V$ is an operator on the same Hilbert space as for $H_0$ and is not so strong, one can construct the Fourier transform associated with $H_0 + V$ by using the technique of functional analysis and partial differential equations.

This is not so easy for operators on hyperbolic manifolds. Even the construction of the Fourier transform associated with the Laplace-Beltrami operator on the hyperbolic space is no longer a trivial work. To construct the Fourier transform on hyperbolic spaces based on the upper half space model or the ball model, one needs deep knowledge of Bessel functions. Under the action of discrete subgroups, the properties of groups will reflect on the structure of manifolds or the construction of generalized eigenfunctions.  


\subsection{Spectral and scattering theory on hyperbolic manifolds}
 In the present note, we deal with the spectral theory and the associated forward and inverse problems for Laplace-Beltrami operators on hyperbolic manifolds. 
 Since we are mainly interested in its spectral properties, Selberg's work \cite{Se56} and its developments are beyond our scope. As an approach to the hyperbolic manifolds from the spectral theory, the first important paper is that of Faddeev \cite{Fa67}. Lang \cite{La75} is a detailed exposition of  Faddeev's article. There are also works of Roelcke \cite{Roe66}, Venkov \cite{Ve90} and a recent article of Iwaniec \cite{Iwa02}. The study of spectral theory, in particular, that of continuous spectrum is drastically changed in these 30 years. The book of Lax-Phillips \cite{LaPh76} has distinguished features, leaning over the analysis of wave equation. The derivation of the analytic continuation of Eisenstein series from that of the resolvent was done by Colin de Verdi{\`e}re \cite{Col81}. Agmon \cite{Ag86} used the modern spectral theory for this problem. Hislop \cite{His94} uses Mourre theory（which is a modern powerful technique to study the continuous spectrum of self-adjoint operators, see e.g. \cite{Is04a}) to prove the resolvent estimates for the Laplacian on hyperbolic spaces. 
 
The scattering metric proposed by Melrose \cite{Me95} aims at constructing a general calculus on non-compact manifolds on which the scattering theory is developed. Melrose' theory includes the following model. Let $\mathcal M$ be a compact $n$-dimensional Riemannian manifold with boundary. Assume that near the boundary, $\mathcal M$ is diffeomorphic to $M\times (0,1)$, $M$ being a compact $n-1$-dimensional manifold, and introduce the following metric
\begin{equation}
ds^2 = \frac{(dy)^2 + A(x,y,dx,dy)}{y^2}, \quad 0 < y < 1, \quad x \in M,
\nonumber
\end{equation}
where $A(x,y,dx,dy)$ is a symmetirc covariant tensor such that as $y \to 0$
\begin{equation}
A(x,y,dx,dy) \sim A_0(x,dx) + yA_1(x,dx,dy) + y^2A_2(x,dx,dy) + \cdots,
\label{eq:IntroAxydxdyexpansion}
\end{equation}
$A_0$ being the Riemannian metric on $M$. This generalizes the upper half-space model of the hyperbolic space. Spectral structures of the associated Laplace-Beltarmi operator were studied by Mazzeo \cite{Ma88} and Mazzeo-Melrose \cite{MaMe87}. Related inverse problem was studied by Joshi-Sa Barreto \cite{JoSaBa00}. In particular, Sa Barreto \cite{SaBa05} proved that the coincidence of the scattering operators gives rise to an isometry of associated metrics. Here the essential role is played by the boundary control method 
presented by Belishev \cite{Be87}, (see also \cite{BeKu87}, \cite{Be97}, \cite{BeKu92}), which makes it possible to reconstruct a Riemannian manifold from the boundary spectral data of the associated Laplace-Belrami operator. 

 A feature of Melrose theory is that it proves the analytic continuation of the resolvent of  Laplace-Beltrami operator for a broad class of metric so that it enables us to study the resonance, another important subject in spectral and scattering theory (\cite{GuZw97}), \cite{Zw99}). We do not deal with the resonance in this note. However, 
let us mention the recent article of Borthwick \cite{Bo07} which studies the inverse problem related to the resonance based on Melrose theory and includes a thorough list of references.

In the case of the Schr{\"o}dinger operator $- \Delta + V(x)$ on ${\bf R}^n$, 
the behavior of solutions to the Schr{\"o}dinger equation has a clear difference depending on the decay order of the potential at infinity. If we assume that $V(x) = O(|x|^{-\rho}), \ |x| \to \infty$, the border line is the case $\rho = 1$. This is also true on hyperbolic spaces. The difference occurs in the case $\rho = 1$ of the decay order $d_h^{-\rho}$, where $d_h$  denotes the hyperbolic distance. In (\ref{eq:IntroAxydxdyexpansion}), $y$ corresponds to $e^{-d_h}$. Hence from the view point of perturbation theory, the theory of scattering metric deals with the case in which the perturbation term is expanded as the power of $e^{-d_h}$.


\subsection{Contents of this note}
The purpose of this note is the exposition of the basic knowledge of the generalized Fourier transform on asymptotically hyperbolic manifolds and their applications to inverse scattering problem. We deal with the general short-range perturbation of the metric, namely, we consider the metric which differ from the standard hyperbolic metric with the term decaying like $d_h^{-1-\epsilon}$, where $d_h$ is the hyperbolic distance. 

More precisely we shall study an $n$-dimensional connected Riemannian manifold ${\mathcal M}$, which is written as a union of open sets:
\begin{equation}
{\mathcal M} = {\mathcal K}\cup{\mathcal M}_1\cup\cdots\cup{\mathcal M}_N.
\nonumber
\end{equation}
The basic assumptions are as follows:

\medskip
\noindent
(A-1) $\ \ \overline{\mathcal K}$ {\it is compact}.

\medskip
\noindent
(A-2) $\ \ {\mathcal M}_i\cap{\mathcal M}_j = \emptyset, \quad i \neq j$.

\medskip
\noindent
(A-3) $\ $ {\it Each} ${\mathcal M}_i$, $i = 1, \cdots, N$, {\it is diffeomorphic either to} ${\mathcal M}_0 = M\times(0,1)$ {\it or to} ${\mathcal M}_{\infty} = M\times(1,\infty)$, $M$ being a compact Riemannian manifold of dimension $n-1$. {\it Here the manifold $M$ is allowed to be different for each} $i$.

\medskip
\noindent
(A-4) {\it On each $\mathcal M_{i}$, the Riemannian metric $ds^2$ has the following form
\begin{equation}
ds^2 = y^{-2}\left((dy)^2 + h(x,dx) + A(x,y,dx,dy)\right),
\label{C0metricds2Axydxdy}
\end{equation}
\begin{equation}
A(x,y,dx,dy) = \sum_{i,j=1}^{n-1}a_{ij}(x,y)dx^idx^j + 2\sum_{i=1}^{n-1}a_{in}(x,y)dx^idy + a_{nn}(x,y)(dy)^2,
\nonumber
\end{equation}
where $h(x,dx) = \sum_{i,j=1}^{n-1}h_{ij}(x)dx^idx^j$ is a positive definite metric on $M$, 
and $a_{ij}(x,y), 1 \leq i,j \leq n$, satisfies the following condition
\begin{equation}
|\widetilde D_x^{\alpha}D_y^{\beta}\, a(x,y)| \leq C_{\alpha\beta}(1 + |\log y|)^{-{\rm min}(|\alpha|+\beta,1)-1-\epsilon_0}, \quad \forall \alpha, \beta
\label{IntroDecay}
\end{equation}
for some $\epsilon_0 > 0$.
Here $\widetilde D_x = \tilde y(y)\partial_x$, $\tilde y(y) \in C^{\infty}((0,\infty))$ such that $\tilde y(y) = y$ for $y > 2$ and $\tilde y(y) = 1$ for $0 < y < 1$.} 

\medskip
Of course this metric $ds^2$ depends on the end $\mathcal M_i$, hence should be written as $ds^2 = y^{-2}\big((dy)^2 + h_i(x,dx) + A_i(x,y,dx,dy)\big)$.

\medskip
Picking up the wave equation, we shall study the following scattering problem. Consider the initial value problem for the wave equation
\begin{equation}
\left\{
\begin{split}
& \partial_t^2u = \Delta_g u \quad {\rm on} \quad \mathcal M, \\
& u\big|_{t=0} = f, \quad \partial_tu\big|_{t=0} = - i\sqrt{-\Delta_g}f,
\end{split}
\right.
\nonumber
\end{equation}
where $f$ is orthogonal to the point spectral subspace for $- \Delta_g$. Then for any compact set $K$ on $\mathcal M$, the solution $u(t)$ behaves as
\begin{equation}
\int_K|u(t)|^2dV_g \to 0 , \quad {\rm as} \quad t \to \pm \infty.
\nonumber
\end{equation}
Namely, the wave disappears from any compact set in $\mathcal M$. On each end $\mathcal M_j$, it will behave like
\begin{equation}
\|u(t) - u_j^{(\pm)}(t)\| \to 0, \quad {\rm as} \quad t \to \pm \infty,
\nonumber
\end{equation}
where $u_j^{(\pm)}(t)$ is the solution to the free wave equation
\begin{equation}
\left\{
\begin{split}
& \partial_t^2u_{j}^{(\pm)} = \Delta_{g_j^0} u_j^{(\pm)}, \quad {\rm on} \quad \mathcal M_j, \\
& u_j^{(\pm)}\big|_{t=0} = f_j^{(\pm)}, \quad \partial_tu_j^{(\pm)}\big|_{t=0} = - i\sqrt{-\Delta_{g_j^0}}f_j^{(\pm)},
\end{split}
\right.
\nonumber
\end{equation}
$\Delta_{g_j^0}$ being the Laplace-Beltrami operator associated with the metric $y^{-2}\big((dy)^2 + h_j(x,dx)\big)$. The scattering operator $\mathcal S$ assigns to the asymptotic data in the remote past that in the remote future:
\begin{equation}
\mathcal S : \big(f_1^{(-)},\cdots,f_N^{(-)}\big) \to \big(f_1^{(+)},\cdots,f_N^{(+)}\big).
\nonumber
\end{equation}
The inverse scattering is an attempt to recover the metric of $\mathcal M$ from the scattering operator $\mathcal S$.
To study this problem, we first investigate the spectral properties of the associated Laplace-Beltrami operator $ - \Delta_g$. Namely

\begin{itemize}
\item
Location of the essential spectrum.

\item Absence of eigenvalues embedded in the continuous spectrum when one of the ends is regular, i.e. one $\mathcal M_i$ is diffeomorphic to $M\times(0,1)$.

\item Discreteness of embedded eigenvalues in the continuous spectrum when all the ends are cusps, i.e. all $\mathcal M_i$ are diffeomorphic to $M_i\times(1,\infty)$.

\item Limiting absorption principle for the resolvent and the absolute continuity of the continuous spectrum.

\end{itemize}
Our next issue is the forward problem. Namely

\begin{itemize}
\item 
Construction of the generalized Fourier transform associated with $- \Delta_g$.

\item 
Asymptotic completeness of time-dependent wave operators.

\item 
Characterization of the space of scattering solutions to the  Helmhotz equation in terms of the generalized Fourier transform.

\item
Asymptotic expansion of scattering solutions to the Helmholtz equation and the S-matrix.

\end{itemize}
As a byproduct, we also study 

\begin{itemize}
\item Representation of the fundamental solution to the wave equation in the upper-half space model.
\item Radon transform and the propagation of singularities for the wave equation.
\end{itemize}
Finally, we shall discuss the inverse problem. Namely
\begin{itemize}
\item Identification of the Riemannian metric from the scattering matrix.
\end{itemize}
We show that two asymptotically hyperbolic manifolds satisfying the above assumptions are isometric, if the metrics coincide on one regular end, and also the S-matrices coincide on that end.

\medskip
The ingredient of each chapter is as follows.

\medskip
\noindent
{\it Chapter 1} $\ $ Fourier transforms on hyperbolic spaces

We discuss the construction of the Fourier transform associated with the Laplace-Beltrami operator of ${\bf H}^n$ as well as its spectral properties. Moreover, we characterize the solution space of the Helmholtz equation in terms of the Fourier transform. We also study the fundamental solution to the wave equation and the Radon transform. We mainly use the estimates of Bessel functions. This chapter is the basis of whole arguments in this note. Main results are Theorems 3.13, 4.2, 4.3, 6.5 and 6.6.

\medskip
\noindent
{\it Chapter 2} $\ $ Perturbation of the metric

This is an exposition of spectral and scattering theory for Laplace-Beltrami operators associated with asympotically hyperbolic metrics on ${\bf R}^n_+$ and their scattering matrices. As in Chapter 1, we will discuss the generalized Fourier transform, the asymptotoic expansion of the resolvent, the Helmholtz equation and the Radon transform. This is also an introduction to the classical spectral theory. Main results are Theorems 2.3, 7.1, 7.8, 7.9, 7.10 and 8.9.

\medskip
\noindent
{\it Chapter 3} $\ $ Manifolds with hyperbolic ends

The general hyperbolic manifolds are constructed by the action of discrete groups on ${\bf H}^n$. We shall consider simple cases and study the spectral properties of the resulting quotient manifolds. We also discuss the action of $SL(2,{\bf Z})$. Main results are Theorems 3.8, 3.12, 3.13 and 3.14.

\medskip
\noindent
{\it Chapter 4} $\ $ Radon transform and propagation of singularities in ${\bf H}^n$

The Radon transform describes singularities of solutions to the wave equation. We shall discuss this classical matter in this chapter for the hyperbolic space. The goal is Theorem 5.2 which is a generalization of Theorem 6.6 in Chapter 1. 

\medskip
\noindent
{\it Chapter 5} $\ $ Introduction to inverse scattering

Local perturbations of the metric of hyperbolic manifolds are identified from the scattering matrix. We shall prove this fact by using spectral theory. Our goal is Theorem 4.8, which asserts that if the metrics coincide on one regular end of the asymptoticaly hyperbolic manifolds, and also the S-matrices coincide on that end, then two manifolds are isometric.

The method we have given here works not only for asymptotically hyperbolic ends but also for the manifolds on which the spectral representation is established. In particlular, Theorem 4.8 holds for manifolds with asymptotically Euclidean ends, or the mixture of Euclidean and hyperbolic ends. 

\medskip
\noindent
{\it Chapter 6} $\ $ Boundary control method

To identify the metric, we reduce the problem to that of the inverse spectral problem on non-compact manifolds with compact boundaries. The crucial role is played by the boundary control method developed by Belishev and Kurylev. This section is devoted to a comprehensive and self-contained exposition of this approach.
We shall give a complete proof of the BC-method except for Tataru's theorem on the uniqueness of solutions to non-characteristic Cauchy  problem for the wave equation.

\medskip
\noindent
{\it Appendix A} $\ $ Radon transform and propagation of singularities in ${\bf R}^n$

The relation between the propagation of singularities and the Radon transform is not obvious even for the case of perturbed Euclidean metric. We shall give detailed proof for this subject for the case of general short-range perturbation of the Euclidean metric. Main results are Theorems 6.7 and 6.10.

\medskip
The main part of our results
 will be proved under a weaker decay assumption on the metric. More precisely, if  we assume instead of (A-4) that in the region $0 < y < y_0$
 
\begin{equation}
ds^2 = y^{-2}\left((dy)^2 + h(x,dx) + B(x,y,dx)\right),
\label{C0ds2Bxydx}
\end{equation}
$$
B(x,y,dx) = \sum_{i,j=1}^{n-1}b_{ij}(x,y)dx^idx^j,
$$
where each $b_{ij}(x,y)$ satisfies
\begin{equation}
|\widetilde D_x^{\alpha}D_y^{\beta}\, b(x,y)| \leq C_{\alpha\beta}(1 + \rho(x,y))^{-1-\epsilon}, \quad \epsilon > 0,
\label{C0decayassumption2}
\end{equation}
 $\rho(x,y)$ being the distance of $(x,y) \in \mathcal M$ from some fixed point, we can derive the same results as those presented below. 
 In fact, we shall prove that the metric of the form (\ref{C0metricds2Axydxdy}) satisfying (\ref{IntroDecay}) is transformed to the metric of the form (\ref{C0ds2Bxydx}) satisfying (\ref{IntroDecay}) (see Theorem 1.6 in Chapter 4), and once we adopt (\ref{C0ds2Bxydx}), we only use the decay assumption  (\ref{C0decayassumption2}).
 
 Even if we start from the metric of the form (\ref{C0metricds2Axydxdy}) satisfying (\ref{C0decayassumption2}), the results below, except for Theorem 2.10, Corollary 2.11 in Chapter 2 and Theorems in Chapter 4, also hold.
 The difference is that the non-existence of eigenvalues embedded in the continuous spectrum  may not be true. However, even in this case, one can show that the embedded eigenvalues are discrete with possible accumulation points $0$ and $\infty$ just like Chapter 3, Theorem 3.5.

\medskip
We have tried to make Chapters 1, 2 and 6 as elementary as possible so that one needs little knowledge to understand the spectral theory and inverse problems. The readers interested in only the inverse problems can skip Chapter 4 and Appendix. If one wants to know the essential step of the limiting absorption principle (resolvent estimates), one should skip Chapter 1 and read subsections 2.3, 2.4 and 2.5 of Chapter 2 first. Although it is written for the upper-half space model, the same idea works for the analysis of ends. We employed the method of integration by parts to prove the limiting absorption principle, which is essentially due to Eidus \cite{Ei69}. This approach is simple and needs no preparatory tool, moreover it is flexible and applicable to various situation. For the other approaches, see e.g. \cite{EGM98}, \cite{FrHi89}, \cite{Kub73},  \cite{Mue87}, \cite{Mue92}.

To construct the generalized Fourier transform, we compute the asymptotic expansion at infinity of the resolvent. This is a classical idea, and has been frequently used (see e.g. \cite{Sa79}, or \cite{Gu92}).
We also utilize the Besov type space introduced by Agmon-H{\"o}rmander \cite{AgHo76} to construct eigenoperators, which, as has been done by Yafaev \cite{Yaf91}, makes it possible to characterize the solution space of the Helmholtz equation by the generalized Fourier transform and to derive the S-matrix from the asymptotic expansion of solutions to the Helmholtz equation.

One can deal with other types of metric by the methods employed here. 
For example, the asymptotically Euclidean ends can be treated in the same way by utilizing results in Chap. 2, \S 5, \S 6 and  Appendix A. The inverse scattering from
asymptotically (Euclidean) cylindrical ends has been studied in \cite{IKL10}. In practical situation, this problem includes that of wave guides. In \cite{IKL11}, inverse scattering from cusp of asymptotically hyperbolic manifolds (or orbifolds) in 2-dimensions is studied. The idea consists in generalizing the notion of S-matrix, which makes it possible to determine all geometrically finite hyperbolic surfaces.
One can also consider a mixture of these different types of ends.

 There are many unknown problems on spectral properties and inverse scattering for a big variety of other types of ends. We hope that the methods in this paper will be helpful for the future study of these fields.

\subsection{Remarks on notation}
\begin{itemize}
\item For two Banach spaces $X, \ Y$, ${\bf B}(X;Y)$ denotes the totality of 
bounded linear operators from $X$ to $Y$. 

\item For a self-adjoint operator $A$ 
\begin{eqnarray*}
\sigma(A) &=& {\rm the \ spectrum \ of} \ A, \\
\sigma_p(A) &=& {\rm the \ set \ of \ all \ eigenvalues \ of } \ A, \\
\sigma_{ac}(A) &=& {\rm the \ absolutely \ continuous \ spectrum \ of } \ A, \\
\sigma_{d}(A) &=& {\rm the \ dscrete \ spectrum \ of \ }\ A, \\
\sigma_{e}(A) &=& {\rm the \ essential \ spectrum \ of \ }\ A.
\end{eqnarray*}

\item For an open set $\Omega$ in a manifold, $C_0^{\infty}(\Omega)$ is the set of all infinitely differentiable functions with compact support in $\Omega$.

\item For a measure $d\mu$ on $\Omega$, $L^2(\Omega ;d\mu)$ denotes all functions $f$ such that
$$
\|f\| = \left(\int_{\Omega}|f|^2d\mu\right)^{1/2} < \infty.
$$

\item For an open set $\Omega$, $H^m(\Omega)$ is the Sobolev space of order $m$ on $\Omega$, namely the set of all functions $f$ on $\Omega$ whose all weak derivatives of order up to $m$ belong to $L^2(\Omega;d\mu)$.

\item
 $H^m_{loc}(\Omega)$ denotes the set of all $u$ such that
$u \in H^m(\omega)$ for all relatively compact open set $\omega$ in $\Omega$.

\item In the inequalities, $C$'s denote various constants. Although these constants may vary from line to line, they are denoted by the same letter $C$.

\item  Theorems, Lemmas, etc. are quoted as follows. In each chapter, Theorem $m.n$ means Theorem $m.n$ of \S $\, m$ of that chapter. Theorem $p.m.n$ means Theorem $m.n$ of Chapter $p$. 

\end{itemize}
Throughout this note, we have assumed the standard knowledge of functional analysis. We have also given a brief explanation for the basic knowledge of the spectrum of self-adjoint operators and partial differential equations when it appears. The reader should consult Kato \cite{Ka76}, 
Reed-Simon \cite{ReSi}, Isozaki \cite{Is04a} for details.

\subsection{Very short perspective} 
Let us explain the basic strategy of constructing the Fourier transform in this paper taking ${\bf R}^1$ as an example. We regard
$H = - d^2/dx^2$ as the Laplacian on the 1-dimensional manifold ${\bf R}^1$.
The resolvent $R(z) = (H - z)^{-1}$ of $H$ has the following expession:
$$
R(z)f(x) = \frac{i}{2\sqrt{z}}\int_{-\infty}^{\infty}
e^{i\sqrt{z}|x-y|}f(y)dy, \quad {\rm Im}\,\sqrt{z} > 0.
$$
Therefore assuming that $f \in L^1({\bf R}^1)$ and $z \to \lambda > 0$, and letting $x \to \pm \infty$, we have
$$
R(\lambda + i0)f(x) \sim i\sqrt{\frac{\pi}{2\lambda}}
e^{\pm i\sqrt{\lambda}x}\widehat f(\pm \sqrt{\lambda}).
$$
Let $E_H(\lambda)$ be the spectral measure for $H$. Then by Stone's formula, we have for $0 < a < b < \infty$
$$
(E_H((a,b))f,f) = \lim_{\epsilon\to0}\frac{1}{2\pi i}
\int_a^b([R(\lambda + i\epsilon) - R(\lambda - i\epsilon)]f,f)d\lambda.
$$
Letting $u = R(\lambda + i0)f$, we have by integration by parts
\begin{eqnarray*}
([R(\lambda + i0) - R(\lambda - i0)]f,f) &=& (u,f) - (f,u)\\
&=& \lim_{R\to\infty}\int_{-R}^R\left(u''\overline{u} - 
u\overline{u''}\right)dx \\
&=& \lim_{R\to\infty}[u'\overline{u} - u\overline{u'}]_{-R}^R\\
&=& \frac{\pi i}{\sqrt{\lambda}}\left(|\widehat f(\sqrt{\lambda})|^2 + 
|\widehat f(-\sqrt{\lambda})|^2\right), 
\end{eqnarray*}
which implies
$$
\|f\|^2 = \lim_{a\to0,b\to\infty}(E_H((a,b))f,f) = 
\int_{-\infty}^{\infty}|{\widehat f}(k)|^2dk.
$$
These calculations suggest that
\begin{itemize}
\item The Fourier transform is obtained from the asymptotic expansion  at infinity of the Green operator of the  Laplacian.
\item Parseval's formula is a consequence of Stone's formula and integration by parts.
\end{itemize}
We should stress that
\begin{itemize}
\item The limit $R(\lambda \pm i0)$ of the resolvent $R(\lambda \pm i\epsilon)$ as $\epsilon \downarrow 0$ plays an important role.
\end{itemize}
The procedure of taking the limit as $\epsilon \downarrow 0$ of $R(\lambda \pm i\epsilon)$ is called the {\it limiting absorption principle}.
\medskip

We shall explain these matters on asymptotically hyperbolic spaces.


\tableofcontents

\chapter{Fourier transforms on the hyperbolic space}


\section{Basic geometry in the hyperbolic space}


\subsection{Upper-half space model} 
We begin with reviewing elementary geometric properties of the hyperbolic space ${\bf H}^n$. Throughout this article ${\bf H}^n$ is the Euclidean upper-half space
\begin{equation}
{\bf R}^n_+ = \{(x,y) \; ;\; x \in {\bf R}^{n-1},\ y > 0\}
\label{eq:upperhalfspace}
\end{equation}
equipped with the metric
\begin{equation}
ds^2 = \frac{|dx|^2 + (dy)^2}{y^n}.
\label{eq:Riemannianmetric}
\end{equation}
In the following, for $v = (v_1,\cdots,v_d) \in {\bf R}^d$, $|v|$ means its Euclidean length :
$|v| = \Big(\sum_{i=1}^dv_i^2\Big)^{1/2}.$


\begin{theorem}
(1) The following 4 maps are the generators of the group of isometries on  ${\bf H}^n$ :\\
\noindent
(a) dilation : $(x,y) \to (\lambda x, \lambda y), \ \lambda > 0$,\\
\noindent
(b) translation : $(x,y) \to (x + v,y), \; v \in {\bf R}^{n-1}$,\\
\noindent
(c) rotation : $(x,y) \to (Rx,y), \; R \in O(n-1)$,\\
\noindent
(d) inversion with respect to the unit sphere centered at $(0,0)$ :
$$ 
(x,y) \to (\overline{x},\overline{y}) = \frac{(x,y)}{|x|^2 + |y|^2}.
$$
\noindent
(2) Any isometry on ${\bf H}^n$ is  a product of the above 4 isometries.
\end{theorem}
Proof. The assertion (1) follows from a direct computation. We use 
$$
d\overline{x} = \frac{dx}{r^2} - \frac{2x}{r^3}dr, \quad
d\overline{y} = \frac{dy}{r^2} - \frac{2y}{r^3}dr,
$$
where $r^2 = x^2 + y^2, \ \overline{x} = x/r^2, \ \overline{y} = y/r^2$, to prove (d). The proof of the assertion (2) is in  \cite{BePe92} pp. 21, 24. \qed

\bigskip
Recall that the inversion with respect to the sphere $\{|x - x_0| = r\}$ is the map: $x \to r^2(x - x_0)/|x - x_0|^2 + x_0$. We give  examples of the isometry in ${\bf H}^2$ and ${\bf H}^3$, which can be proved by a straightforward computation.

\subsection{${\bf H}^2$ and linear fractional transformation} When $n = 2$, it is convenient to identify a point $(x,y) \in {\bf H}^2$ with the complex number $z = x + iy$. For a matrix 
$$
\gamma = \left(
\begin{array}{cc}
a & b \\
c & d
\end{array}
\right) \in SL(2,{\bf R}),
$$
the linear fractional transformation
$$
z \to \gamma\cdot z := \frac{az + b}{cz + d}
$$
defines an isometry on ${\bf H}^2$. 

\subsection{${\bf H}^3$ and quarternions} 
Represent a point $(x_1,x_2,x_3) 
\in {\bf H}^3$ by a quarternion
$$
{\bf z} = x_1{\bf 1} + x_2{\bf i} + x_3{\bf j} = 
\left(
\begin{array}{cc}
x_1 + ix_3 & x_2 \\
- x_2 & x_1 - ix_3
\end{array}
\right)
$$
with $\bf k$-component equal to 0 ; then ${\bf H}^3 \subset {\bf Q}$.
For a matrix 
$$
\gamma = \left(
\begin{array}{cc}
a & b \\
c & d
\end{array}
\right) \in SL(2,{\bf C}),
$$
the M{\"o}bius transformation
$$
{\bf z} \to \gamma\cdot {\bf z} := (a{\bf z} + b)(c{\bf z} + d)^{-1}.
$$
acts from ${\bf H}^3$ to ${\bf Q}$. Using $ad - bc = 1$, straightforward although lengthy computations show that $\gamma\cdot{\bf z}$ actually belongs to ${\bf H}^3$.
Thus $\gamma$ defines an isometry on ${\bf H}^3$. 

\subsection{Geodesics} 
The equation of a geodesic in a Riemannian manifold with 
metric $ds^2 = g_{ij}dx^idx^j$ is 
$$
\frac{d^2x^k}{dt^2} + \Gamma^k_{ij}\frac{dx^i}{dt}\frac{dx^j}{dt} = 0,
$$
$$
\Gamma^k_{ij} = \frac{1}{2}g^{kp}\left(\frac{\partial g_{jp}}{\partial x^i} 
+ \frac{\partial g_{ip}}{\partial x^j} - \frac{\partial g_{ij}}{\partial x^p}\right),
$$
where $(g^{ij})$ is the inverse matrix of $(g_{ij})$.
It is well-known that this may be rewritten as Hamilton's canonical equation
with Hamiltonian $h(x,\xi) = \frac{1}{2}g^{ij}\xi_i\xi_j$:
$$
\dfrac{dx^i}{dt} = \dfrac{\partial h}{\partial \xi_i}, \quad
\dfrac{d\xi_i}{dt} = - \dfrac{\partial h}{\partial x^i}.
$$
(One can  check it directly by 
using the formula $\dfrac{\partial g^{ij}}{\partial x^m} = - g^{ik}\left(\dfrac{\partial g_{kr}}{\partial x^m}\right)g^{rj}$). In the case of ${\bf H}^n$, with $(\xi, \eta)$ dual to $(x, y)$, 
Hamilton's  equation turns out to be
\begin{equation}
\left\{
\begin{split}
&\frac{dx}{dt} = y^2\xi, \quad \frac{dy}{dt} = y^2\eta, \\
&\frac{d\xi}{dt} = 0, \quad \frac{d\eta}{dt} = - y(|\xi|^2 + \eta^2).
\end{split}
\right.
\nonumber
\end{equation}
Hence $\xi$ does not depend on $t$. If $\xi = 0$, the curve becomes a straight
line $\{x = x(0)\}$. When $\xi \neq 0$, $(x(t),y(t))$ moves in the 2-dimensional plane spanned by 2 vectors $(\xi,0)$ and $(0,1)$, which is denoted by $\Pi$. We use the same $(x,y)$ to denote the rectangular coordinates on $\Pi$. Since the energy $h$ is conserved, $y(t)^2(|\xi|^2 + \eta(t)^2)$ is a constant, which is denoted by $2E$. Then $\eta^2 = 2E/y^2 - |\xi|^2$, which implies
$$
y' = \frac{dy}{dx} = \frac{\eta}{|\xi|} = \pm \sqrt{\frac{A}{y^2} - 1}, \quad 
A = \frac{2E}{|\xi|^2}.
$$
Solving this equation, we get $(x + B)^2 + y^2 = A$. 
We have thus proven


\begin{lemma}
There are only two kinds of geodesics in ${\bf H}^n$ : \\
\noindent
(a) the hemi-circles with center on the hyperplane $\{y = 0\}$, \\
\noindent
(b) the straight lines perpendicular to the hyperplane $\{y = 0\}$. 
\end{lemma}

We see by Lemma 1.2 that for two points $P, Q \in {\bf H}^n$, there exists only one geodesic passing through
$P$ and $Q$.

\begin{lemma} For two points $(a,b), (a',b') \in {\bf H}^n$, there exists an isometry which maps $(a,b)$ to $(0,1)$ and $(a',b')$ to $(0,c)$, where
$$
\left(\tanh\frac{|\log c|}{2}\right)^2 = \frac{|a - a'|^2 + (b - b')^2}{|a - a'|^2 + (b + b')^2}.
$$
\end{lemma}
Proof. By the following isometries, $(a,b)$ is mapped to $(0,1)$ :
$$
(a,b) \to (\frac{a}{b},1) \quad {\rm (dilation)} \quad
\to (0,1) \quad {\rm (translation)}.
$$
Then $(a',b')$ is mapped to $(\frac{a' - a}{b},\frac{b'}{b})$. Therefore, we have only to show that for any $(x,y)$ there exists an isometry which maps $(x,y)$ to $(0,c)$ with suitable $c$ leaving $(0,1)$ invariant. 
The problem is then reduced to 2-dimensions. Consider the linear fractional transformation  by
$$
\gamma = \left(
\begin{array}{cc}
\cos\theta & - \sin\theta \\
\sin\theta & \cos\theta
\end{array}
\right),
$$
which leaves $i$ invariant. Then for given $z = x + iy$,
$$
\gamma\cdot z = \frac{\frac{|z|^2 - 1}{2}\sin2\theta + x\cos2\theta + iy}{
|z\sin\theta + \cos\theta|^2}.
$$
By choosing $\theta$ so that the real part vanishes, we get the isometry which maps $x + iy$ to $ic$. Let us compute $c$. Assuming that $x > 0$, by our choice of $\theta$,
$$
\cos2\theta = \frac{1 - |z|^2}{[(1 - |z|^2)^2 + 4x^2]^{1/2}}, \quad
\sin2\theta = \frac{2x}{[(1 - |z|^2)^2 + 4x^2]^{1/2}}. 
$$
Therefore
\begin{eqnarray*}
|z\sin\theta + \cos\theta|^2 &=& \frac{1 + |z|^2}{2} + \frac{1 - |z|^2}{2}\cos2\theta + x\sin2\theta \\
&=& 
\frac{1}{2}\left[1 + |z|^2 + \left((1 - |z|^2)^2 + 4x^2\right)^{1/2}\right],
\end{eqnarray*}
hence
$$
c = \frac{2y}{1 + |z|^2 + \left((1 - |z|^2)^2 + 4x^2\right)^{1/2}} 
= \frac{1 + |z|^2 - ((1 - |z|^2)^2 + 4x^2)^{1/2}}{2y}.
$$
This implies
$$
\left(\tanh\frac{|\log c|}{2}\right)^2 = \frac{1 + |z|^2 - 2y}{1 + |z|^2 + 2y}.
$$
Putting $x = |a - a'|/b, y = b'/b$, we complete the proof of the lemma.
\qed

\bigskip
The hyperbolic distance from $(0,1)$ to $(0,c)$ is given by
$$
\left|\int_1^c\frac{dy}{y}\right| = |\log c|.
$$
This and Lemma 1.3 imply the following formula. 


\begin{lemma}
The hyperbolic distance  $d = d \big((x,y),(x',y')\big)$ between $(x,y)$ and $(x',y')$ is given by
$$
\left(\tanh \frac{d}{2}\right)^2 = 
\frac{|x - x'|^2 + |y - y'|^2}{|x - x'|^2 + |y + y'|^2}.
$$
\end{lemma}

From this lemma, we get
\begin{equation}
\frac{1}{2}\big(\cosh d - 1\big) = \frac{|x - x'|^2 + |y - y'|^2}{4yy'}.
\label{eq:Chap1Sec1Hyperbolicdistance}
\end{equation}


\begin{lemma}
The geodesic sphere in ${\bf H}^n$ is a Euclidean sphere.
\end{lemma}

For example the geodesic sphere in ${\bf H}^n$ with center $(0,1)$ and radius $r > 0$ is written as
$$
|x|^2 + (y - (1 + 2\delta))^2 = 4\delta(1 + \delta), \quad
\delta = (\cosh r - 1)/2.
$$
This is a Euclidean sphere with center $(0,\cosh r)$ and radius $\sinh r$.

The following formula is a corollary of the previous considerations :
\begin{equation}
ds^2 = (dr)^2 + \sinh^2r(d\theta)^2,
\label{NormalCoordinatesmetric}
\end{equation}
where $(r,\theta) \in [0,\infty)\times S^{n-1}$ are geodesic polar coordinates centered at $(0,1)$, and $(d\theta)^2$ is the standard metric on $S^{n-1}$.


\subsection{Estimate of the metric} 
Let $d_h(x,y)$ be the hyperbolic distance between $(x,y)$ and $(1,0)$. For $w \in {\bf R}^d$, we put $\langle w\rangle = (1 + |x|^2)^{1/2}$, and define
\begin{equation}
\rho_0(x,y) = \log \langle x\rangle + \langle \log y\rangle.
\label{eq:C1S1Defofrho}
\end{equation}


\begin{lemma}
There exists a constant $C_0 > 0$ such that on ${\bf H}^n$
$$
C_0^{-1}\big(1 + \rho_0(x,y)\big) \leq 1 + d_h(x,y) \leq C_0\big(1 + \rho_0(x,y)\big).
$$
\end{lemma}
Proof. By (\ref{eq:Chap1Sec1Hyperbolicdistance}), $\cosh d_h = (|x|^2 + y^2 + 1)/(2y)$. If $y$ is small, $e^{d_h} \sim (|x|^2 + 1)/y$, and we obtain the lemma easily. If $y$ is large, $e^{d_h} \sim y + |x|^2/y$. The estimate from above is easy, since $e^{d_h} \leq C(y + |x|^2)$. The estimate from below is obtained by cosidering two cases $y > \sqrt{|x|}$ and $y < \sqrt{|x|}$. \qed


\section{Besov type spaces}
The Fourier transform $\widehat f(\xi)$ of a function $f(x)$ on ${\bf R}^n$ becomes smooth if $f(x)$ decays rapidly at infinity, and we can restrict $\widehat f(\xi)$ on a hypersurface in ${\bf R}^n$. The best possible space to describe the relation between the decay at infinity of ${\bf R}^n$ and the restriction of its Fourier transform on a hypersurface was found by Agmon-H{\"o}rmander \cite{AgHo76}. Let us point out that Murata (\cite{Mu74}, \cite{Mu80}) had discovered this space in his study of the asymptotic behavior at infinity of solutions of linear partial differential equations. This space furnishes a natural framework to characterize  solutions to the Helmholtz equation. We introduce this space for ${\bf H}^n$.

  
\subsection{The  Besov type space}
Let $\bf h$ be a Hilbert space endowed with inner product $(\;,\;)_{\bf h}$ and norm $\|\cdot \|_{\bf h}$. 
We decompose $(0,\infty)$ into $(0,\infty) = \cup_{k \in {\bf Z}}I_k$, where
$$
I_k = 
\left\{
\begin{array}{cc}
 \big(\exp(e^{k-1}),\exp(e^k)\big], & k \geq 1 \\
 \big(e^{-1},e\big],& k = 0 \\
 \big(\exp(- e^{|k|}),\exp(- e^{|k| -1})\big], & k \leq - 1.
\end{array}
\right.
$$
We fix a natural number $n \geq 2$ and put 
$$
d\mu(y) = \frac{dy}{y^n}.
$$


\begin{definition}
Let $\mathcal B$ be the space of $\bf h$-valued function on $(0,\infty)$ satisfying
$$
\|f\|_{\mathcal B} = \sum_{k\in{\bf Z}} 
e^{|k|/2}\left(\int_{I_k}\|f(y)\|_{\bf h}^2d\mu(y)\right)^{1/2} < \infty.
$$
\end{definition}


\begin{lemma} (1) The following inequality holds :
$$
\int_0^{\infty}y^{(n - 1)/2}\|f(y)\|_{\bf h}d\mu(y) \leq C\|f\|_{\mathcal B}, \quad
\forall f \in {\mathcal B}
$$
(2) For any $T \in {\mathcal B}^{\ast}$, there exits a unique $v_T \in L^2_{loc}((0,\infty);{\mathcal H})$ such that
$$
T(f) = \int_0^{\infty}\big(f(y),v_T(y)\big)_{\bf h}d\mu(y), \quad \forall 
f \in {\mathcal B},
$$
$$
\|T\| = \sup_{k\in{\bf Z}}e^{-|k|/2}\left(\int_{I_k}\|v_T(y)\|_{\bf h}^2
d\mu(y)
\right)^{1/2}.
$$
\end{lemma}
Proof.  By the Schwarz inequality, we have
$$
\int_0^{\infty}y^{(n-1)/2}\|f(y)\|_{\bf h}\frac{dy}{y^n} 
\leq 
\sum_k\left(\int_{I_k}\frac{dy}{y}\right)^{1/2}
\left(\int_{I_k}\frac{\|f(y)\|_{\bf h}^2}{y^n}dy\right)^{1/2}.
$$
Since $\int_{I_k}dy/y \leq Ce^{|k|}$, we get the assertion (1).

Let $T_k$ be the restriction of $T$ on $L^2(I_k;{\mathcal H})$. Then we have for $f$ which vanishes outside $I_k$
$$
|T_k(f)| = |T(f)| \leq \|T\|\|f\|_{\mathcal B} = \|T\|
e^{|k|/2}\left(\int_{I_k}\|f(y)\|_{\bf h}^2d\mu(y)\right)^{1/2}.
$$
Therefore by the theorem of Riesz, there exists $v_T^{(k)}(y) \in L^2(I_k;{\mathcal H})$ such that
$$
T(f) = \int_{I_k}\left(f(y),v_T^{(k)}(y)\right)_{\bf h}d\mu(y), \quad
\forall f \in L^2(I_k;{\bf h}),
$$
$$
\left(\int_{I_k}\|v_T^{(k)}(y)\|_{\bf h}^2d\mu(y)\right)^{1/2} \leq
\|T\|e^{|k|/2}.
$$
Putting $v_T(y) =v_T^{(k)}(y), \ y \in I_k$, we then have
$$
\sup_k e^{-|k|/2}\left(\int_{I_k}\|v_T(y)\|_{\bf h}^2d\mu(y)\right)^{1/2} \leq
\|T\|.
$$
Let $\chi_k$ be the characteristic function of $I_k$. Then for any $f \in {\mathcal B}$
\begin{eqnarray*}
T(f) &=& \sum_k T(\chi_k f) \\
&=& \sum_k\int_{I_k}\left(f(y),v_T^{(k)}(y)\right)_{\bf h}d\mu(y) \\
& =& \int_0^{\infty}\left(f(y),v_T(y)\right)_{\bf h}d\mu(y).
\end{eqnarray*}
We now put
$$
a_k = e^{|k|/2}\left(\int_{I_k}\|f(y)\|_{\bf h}^2d\mu(y)\right)^{1/2},
\quad
b_k = e^{- |k|/2}\left(\int_{I_k}\|v_T(y)\|_{\bf h}^2d\mu(y)\right)^{1/2}.
$$
Then since
\begin{eqnarray*}
|T(f)| &\leq& \sum_k\int_{I_k}\|f(y)\|_{\bf h}\|v_T(y)\|_{\bf h}
d\mu(y) \\
&\leq& \sum_ka_kb_k \ \leq \ \sum_ka_k\left(\sup_kb_k\right),
\end{eqnarray*}
we have $\|T\| \leq \sup_kb_k$. \qed

\bigskip
By this lemma, ${\mathcal B}^{\ast}$ is identified with the Banach space with norm
$$
\|v\|_{\mathcal B^{\ast}} = 
\sup_{k \in {\bf Z}} e^{-|k|/2}\left(\int_{I_k}\|v(y)\|_{\bf h}^2d\mu(y)\right)^{1/2}. 
$$
However, the following norm is easier to handle.


\begin{lemma} There exists a constant $C > 0$ such that
$$
C\|v\|_{\mathcal B^{\ast}} \leq \left(\sup_{R > e}\;\frac{1}{\log R}\int_{\frac{1}{R} < y < R}
\|v(y)\|_{\bf h}^2d\mu(y)\right)^{1/2} \leq C^{-1}\|v\|_{\mathcal B^{\ast}}.
$$
\end{lemma}
Proof. We put
$$
A = \sup_{k \in {\bf Z}}e^{-|k|}\int_{I_k}\|v(y)\|_{\bf h}^2d\mu,
\quad
B = \sup_{R > e}\;\frac{1}{\log R}\int_{\frac{1}{R} < y < R}
\|v(y)\|_{\bf h}^2d\mu.
$$
For any $\epsilon > 0$ there exists $k \in {\bf Z}$ such that
$$
e^{-|k|}\int_{I_k}\|v(y)\|_{\bf h}^2d\mu > A - \epsilon.
$$
By putting $\log R = e^{|k|}$, we have
$$
\frac{1}{\log R}\int_{\frac{1}{R} < y < R}\|v(y)\|_{\bf h}^2d\mu \geq
e^{-|k|}\int_{I_k}\|v(y)\|_{\bf h}^2d\mu.
$$
This implies $B \geq A$.
 
On the other hand for any $\epsilon > 0$ there exists $R > e$ such that
$$
\frac{1}{\log R}\int_{\frac{1}{R} < y < R}\|v(y)\|^2_{{\bf h}}d\mu > B - \epsilon.
$$
Choosing $k \in {\bf Z}$ such that $\exp(e^k) \leq R \leq \exp(e^{k+1})$  we then have
\begin{eqnarray*}
\frac{1}{\log R}\int_{\frac{1}{R} < y < R}\|v(y)\|_{\bf h}^2d\mu 
&\leq& 
\frac{1}{\log R}\sum_{|p| \leq k + 1}\int_{I_p}\|v(y)\|_{\bf h}^2d\mu \\
&\leq& \frac{A}{\log R}\sum_{|p| \leq k + 1}e^{|p|} \leq CA. \qquad 
\qed
\end{eqnarray*}


\begin{definition}
We identify ${\mathcal B}^{\ast}$ with the space equipped with norm 
$$
\|u\|_{\mathcal B}^{\ast} = 
\left(\sup_{R > e}\;\frac{1}{\log R}\int_{\frac{1}{R} < y < R}
\|u(y)\|_{\bf h}^2d\mu\right)^{1/2} < \infty.
$$
\end{definition}

The following inequality holds :

\begin{equation}
 \left|(f,v)\right| = \left|\int_0^{\infty}(f(y),v(y))_{{\bf h}}d\mu\right|
 \leq C\|f\|_{{\mathcal B}}\|v\|_{{\mathcal B}^{\ast}}.
 \nonumber
\end{equation}


\begin{lemma} \label{L:2.5}
(1) The following assertions (\ref{eq:Chap1Sec2logR1}) and (\ref{eq:Chap1Sec2logR2}) are equivalent.
\begin{equation}
\lim_{R\to\infty}\frac{1}{\log R}\int_{\frac{1}{R} < y < R}
\|u(y)\|_{\bf h}^2d\mu = 0.
\label{eq:Chap1Sec2logR1}
\end{equation}
\begin{equation}
\lim_{R\to\infty}\frac{1}{\log R}\int_0^{\infty}
\rho\big(\frac{\log y}{\log R}\big)\|u(y)\|_{\bf h}^2d\mu = 0, \quad
\forall \rho \in C_0^{\infty}({\bf R}).
\label{eq:Chap1Sec2logR2}
\end{equation}
(2)  A function $u$ belongs to ${\mathcal B}^{\ast}$ if and only if
\begin{equation}
\sup_{R>e}\frac{1}{\log R}\int_0^{\infty}
\rho\big(\frac{\log y}{\log R}\big)\|u(y)\|^2_{{\bf h}}d\mu < \infty, \quad \forall \rho \in C_0^{\infty}({\bf R})
\nonumber
\end{equation}
\end{lemma}
Proof. To prove (1), we have only to note that (\ref{eq:Chap1Sec2logR1}) is equivalent to
\begin{equation}
\lim_{R\to\infty}\frac{1}{\log R}\int_{R^{a} < y < R^b}
\|u(y)\|_{\bf h}^2d\mu = 0, \quad - \infty < \forall a < \forall b < \infty.
\label{eq:Chap1Sec2logR3}
\end{equation}
Indeed, letting $R = R^c$, $c = {\rm max}\,\{|a|,|b|\}$, in (\ref{eq:Chap1Sec2logR1}), we get (\ref{eq:Chap1Sec2logR3}). Letting $a = 1, \; b = - 1$ in (\ref{eq:Chap1Sec2logR3}), we get (\ref{eq:Chap1Sec2logR1}). Since $a$ and $b$ are arbitrary, (\ref{eq:Chap1Sec2logR3}) is equivalent to (\ref{eq:Chap1Sec2logR2}).

The assertion (2) is proved  similarly.
\qed

\bigskip
In the upper half-space model ${\bf R}^n_+$, we represent a point of ${\bf R}^n_+$ as $(x,y), x \in {\bf R}^{n-1}, y > 0$, and put ${\bf h} = L^2({\bf R}^{n-1})$. 


\subsection{Weighted $L^2$ space}  
The following spaces are also useful. 


\begin{definition}
For $s \in {\bf R}$, we define the space $L^{2,s}$ by
$$
 u \in L^{2,s} \Longleftrightarrow \|u\|_s^2 = \int_0^{\infty}
 (1 + |\log y|)^{2s}\|u(y)\|_{\bf h}^2\,d\mu(y) < \infty.
$$
\end{definition}


\begin{lemma}
 For $s > 1/2$, we have the following inclusion relations :
$$
 L^{2,s} \subset {\mathcal B} \subset L^{2,1/2} \subset L^2 \subset
 L^{2,-1/2} \subset {\mathcal B}^{\ast} \subset L^{2,-s}.
$$
\end{lemma}
Proof. We put
$$
a_{k,s} = \left(\int_{I_k}(1 + |\log y|)^{2s}\|u(y)\|_{\bf h}^2d\mu(y)\right)^{1/2}.
$$
Since
$$
 C^{-1}e^{|k|} \leq 1 + |\log y| \leq Ce^{|k|}, \quad 
 y \in I_k,
$$
we have
$$
 C^{-1}e^{|k|s}a_{k,0} \leq a_{k,s} \leq Ce^{|k|s}a_{k,0}. 
$$
This implies
$$
 \|u\|_{1/2} = \sqrt{\sum_k(a_{k,1/2})^2} 
 \leq \sum_k a_{k,1/2} 
 \leq C\sum_ke^{|k|/2}a_{k,0} = C\|u\|_{\mathcal B}.
$$
Letting $\epsilon = s - 1/2 > 0$, we have
\begin{eqnarray*}
 \|u\|_{\mathcal B} = \sum_ke^{-|k|\epsilon}e^{|k|s}a_{k,0} 
 \leq C\sum_ke^{-|k|\epsilon}a_{k,s} 
 \leq C(\sum_ka_{k,s}^2)^{1/2} = C\|u\|_s.
\end{eqnarray*}
These two relations yield
$L^{2,s} \subset {\mathcal B} \subset L^{2,1/2}$.
Passing to the dual spaces, we have
$L^{2,-1/2} \subset {\mathcal B}^{\ast} \subset L^{2,-s}$. \qed


\section{1-dimensional problem}


\subsection{Some facts from functional analysis} 
Let us recall basic terminologies.  A densely defined linear operator $A$ on a Hilbert space $\mathcal H$ is said to be {\it symmetric} if
$
(Au,v) = (u,Av), \ \forall u, v \in D(A).
$
If $A$ is symmetric, $D(A) \subset D(A^{\ast})$ and $A^{\ast}u = Au$ for $u \in D(A)$. A symmetric operator $A$ is said to be {\it self-adjoint} if
$D(A^{\ast}) = D(A)$. 
The {\it closure} $\overline{A}$ of a symmetric operator $A$ is defined as follows: $u \in D(\overline{A})$, $\overline{A}u = f$ if and only if there exists $\{u_n\} \in D(A)$ such that $u_n \to u$, $Au_n \to f$. A symmetric operator $A$ is said to be {\it essentially self-adjoint} if $\overline A$ is self-adjoint. $A$ is essentially self-adjoint if and only if
${\rm Ker}\,(A^{\ast} \pm i) = \{0\}$. This is equivalent to ${\rm Ker}\,(A^{\ast} - z) = \{0\}$ if ${\rm Im}\,z \neq 0$. For the proof of these facts, see e.g. \cite{ReSi}, Vol. 1  and Vol. 3.

 Suppose we are given a differential operator
$A = a(y)\partial_y^2 + b(y)\partial_y + c(y)$ 
on the interval $(0,\infty)$. 
We shall assume that the coefficients of $A$ is sufficiently smooth, $a(y) \neq 0$ on $(0, \infty)$,  and that there exists a function $\rho(y) > 0$ 
such that $A\big|
_{C_0^{\infty}((0,\infty))}$ is essentially self-adjoint in $\mathcal H = L^2((0,\infty);\rho(y)dy)$. For ${\rm Im}\,z \neq 0$, let $\varphi_0(y)$ and $\varphi_{\infty}(y)$ be non-trivial solutions of $(A - z)u = 0$ on $(0,\infty)$ such that
$$
\varphi_0(y) \in L^2((0,1);\rho(y)dy), \quad
\varphi_{\infty}(y) \in L^2((1,\infty);\rho(y)dy).
$$


\begin{lemma}
 $\varphi_0(y)$ and $\varphi_{\infty}(y)$ are linearly independent.
\end{lemma}

Proof. If they were linearly dependent, then $\varphi_0(y) \in {\mathcal H}$. Therefore, since $A$ is self-adjoint, $\varphi_0(y) = 0$, which is a contradiction. \qed

\bigskip
Let $W(y)$ be the Wronskian: 
$$
 W(y) = \varphi_0(y)\varphi_{\infty}'(y) - 
 \varphi_0'(y)\varphi_{\infty}(y) \neq 0
$$
 and define the {\it Green function} $G(y,y')$ by
\begin{equation}
 G(y,y') = \frac{1}{a(y')\rho(y')W(y')}\left\{
 \begin{split}
& \varphi_0(y)\varphi_{\infty}(y'), \quad 0 < y < y',  \\
& \varphi_{\infty}(y)\varphi_{0}(y'), \quad  0 < y' < y.
 \end{split}
 \right.
 \nonumber
\end{equation}
The integral operator
$$
Gf(y) = \int_0^{\infty}G(y,y')f(y')\rho(y')dy'
$$
is called the {\it Green operator} of $A - z$. Let $\|\cdot\|$  be the norm in $\mathcal H$.


\begin{lemma}
(1) If ${\rm Im}\, z \neq 0$, 
$$
 \|Gf\| \leq \frac{1}{|{\rm Im}\,z|}\|f\|.
$$
(2) For $f \in \mathcal H$, 
$(A - z)Gf = f.$
\end{lemma}
Proof. (1) is a standard fact (see e.g. \cite{ReSi} Vol 1). For $f \in C_0^{\infty}((0,\infty))$, we put  $u = Gf$. One can then find a small $\epsilon > 0$ such that $u = C\varphi_0(y)$ for $y < \epsilon$ and  $u = C'\varphi_{\infty}(y)$ for $y > 1/\epsilon$. Hence $u \in \mathcal H$. 
Using $(A - z)\varphi_0 = (A - z)\varphi_{\infty} = 0$, we have, by a direct, 
computation
$$
 (A - z)u = \left(\varphi_{\infty}'\varphi_0 -  \varphi_{0}'\varphi_{\infty}\right)
 \frac{a\rho}{a\rho W}f = f.
$$
This implies that $G = (A - z)^{-1}$ on $C_0^{\infty}((0,\infty))$, and proves (2) for such $f'$s. 
As $\|(A_z)^{-1}\| \leq |{\rm Im}\,z|^{-1}$, 
by  approximating $f \in L^2((0,\infty))$  by $f_n \in C_0^{\infty}((0,\infty))$, we obtain (1) and (2) for the whole $\mathcal H$. 
\qed

\bigskip
We  explain the elliptic regularity theorem in the 1-dimensional case. Let $I\subset {\bf R}$ be an open interval and $A = - d^2/dx^2 + a_1(x)d/dx + a_0(x)$ be a differential operator with smooth coefficients. The formal adjoint $A^{\dagger}$ is defined by
$$
A^{\dagger}\varphi(x) = - \frac{d^2}{dx^2}\varphi(x) - \frac{d}{dx}\left(\overline{a_1(x)}\varphi(x)\right) + \overline{a_0(x)}\varphi(x).
$$
A function $u(x)$ is said to be a weak solution of the equation $Au = f$ on $I$ if
$$
\int_Iu(x)\overline{A^{\dagger}\varphi(x)}dx = \int_If(x)\overline{\varphi(x)}dx, \quad 
\forall \varphi \in C_0^{\infty}(I).
$$


\begin{lemma}
If $u$ is a weak solution to the equation $Au = f$ on $I$ with $f \in C^{\infty}(I)$, 
then actually $u \in C^{\infty}(I)$ and $Au = f$ holds in the classical sense.
\end{lemma}
Proof. By Corollary 3.1.6 of \cite{Hor}, we have $u \in C^2(I)$
if, e.g. $f \in C^1(I)$. Since $u'(x)$ is a weak solution to the equation
$$
\left(- \frac{d^2}{dx^2} + (a_1 + a_1')\frac{d}{dx} + a_0\right)u' = f' - a_0'u,$$
we have $u' \in C^2(I)$, hence $u \in C^3(I)$. Repeating this procedure, we prove the lemma. \qed


\subsection{Bessel functions} 
We summarize basic knowledge of Bessel functions utilized in this note. For the details, see \cite{MUH59}, \cite{Le72} and  \cite{Wa62}.

The modified Bessel function (of 1st kind) $I_{\nu}(z)$ with parameter $\nu \in {\bf C}$ 
is defined by
\begin{equation}
I_{\nu}(z) = \left(\frac{z}{2}\right)^{\nu}\sum_{n=0}^{\infty}
\frac{(z^2/4)^n}{n!\,\Gamma(\nu + n + 1)}, \quad
z \in {\bf C}\setminus(-\infty,0].
\label{eq:Chap1Sec3DefinitionofInu}
\end{equation} 
It is related with the Bessel function $J_{\nu}(z)$ by
\begin{equation}
I_{\nu}(y) = e^{-\nu\pi i/2}J_{\nu}(iy), \quad y > 0.
\nonumber
\end{equation}
The following function $K_{\nu}(z)$ is also called the modified Bessel function, or the K-Bessel function, or sometimes the Macdonald function:
\begin{equation}
K_{\nu}(z) = \frac{\pi}{2}\frac{I_{-\nu}(z) - I_{\nu}(z)}{\sin(\nu\pi)}, \quad
\nu \notin {\bf Z},
\label{eq:KnuandInu}
\end{equation}
\begin{equation}
K_n(z) = K_{-n}(z) = \lim_{\nu \to n}K_{\nu}(z), \quad 
n \in {\bf Z}.
\nonumber
\end{equation}
These $I_{\nu}(z), K_{\nu}(z)$ solve the following equation 
\begin{equation}
z^2u'' + zu' - (z^2 + \nu^2)u = 0,
\nonumber
\end{equation}
and have the following asymptotic expansions as $|z| \to \infty$:
\begin{equation}
I_{\nu}(z) \sim \frac{e^z}{\sqrt{2\pi z}} +
 \frac{e^{-z + (\nu + 1/2)\pi i}}{\sqrt{2\pi z}}, \quad |z| \to \infty, \quad
- \frac{\pi}{2} < {\rm arg}\,z < \frac{\pi}{2},
\label{eq:Chap1Sec3Inunearinfty}
\end{equation}
\begin{equation}
K_{\nu}(z) \sim \sqrt{\frac{\pi}{2z}}e^{-z}, \quad 
|z| \to \infty, \quad - \pi < {\rm arg}\,z < \pi.
\label{eq:Chap1Sec3Knunearinfinity}
\end{equation}
The asymptotics as $z \to 0$ are as follows:
\begin{equation}
I_{\nu}(z) \sim \frac{1}{\Gamma(\nu + 1)}\left(\frac{z}{2}\right)^{\nu},
\label{eq:Chap1Sec3Inunear0}
\end{equation}
\begin{equation}
K_{\nu}(z) \sim \frac{\pi}{2\sin(\nu\pi)}
\left(\frac{1}{\Gamma(1 - \nu)}\left(\frac{z}{2}\right)^{-\nu} - 
\frac{1}{\Gamma(1 + \nu)}\left(\frac{z}{2}\right)^{\nu}\right), \quad 
\nu \not\in {\bf Z}
\label{eq:Chap1Sec3Knunear0}
\end{equation}
\begin{equation}
K_n(z) \sim \left\{
\begin{split}
& - \log z, \quad n = 0, \\
& 2^{n-1}(n - 1)!z^{-n}, \quad n = 0, 1, 2, \cdots
\end{split}
\right.
\nonumber
\end{equation}

Let $n \geq 2$ be an integer, and a parameter $\zeta \in {\bf C}$ satisfy ${\rm Re}\,\zeta \geq 0$. We consider the differential operator
\begin{equation}
L_0(\zeta) = y^2(- \partial_y^2 + \zeta^2) + (n - 2)y\partial_y - \frac{(n-1)^2}{4}
\label{eq:DiffOpL0zeta}
\end{equation}
on the interval $(0,\infty)$. Let $(\;,\,)$ be the inner product of $L^2((0,\infty);dy/y^n)$. We have
\begin{equation}
(L_0(\zeta)u,v) = (u,L_0(\overline{\zeta})v), \quad \forall u, v
\in C_0^{\infty}((0,\infty)).
\label{eq:Chap1Sec3L0zetauv}
\end{equation}
When $\zeta \neq 0$, the equation $(L_0(\zeta) + \nu^2)u = 0$ has
two linearly independent solutions
$$
y^{(n-1)/2}I_{\nu}(\zeta y), \ y^{(n-1)/2}K_{\nu}(\zeta y),
$$ 
and when $\zeta = 0$ and $\nu \neq 0$, these two linearly independent solutions are $y^{\frac{n-1}{2} \pm \nu}$.


\begin{theorem}
If $\zeta \geq 0$, $L_0(\zeta)\big|_{C_0^{\infty}((0,\infty))}$ is essentially self-adjoint.
\end{theorem}
Proof. We have only to show that
$$
(u,(L_0(\zeta) \pm i)\varphi) = 0, \quad \forall \varphi \in C_0^{\infty}((0,\infty)) \Longrightarrow u = 0.
$$
Suppose $(u,(L_0(\zeta) + i)\varphi) = 0, \ \forall \varphi \in C_0^{\infty}((0,\infty))$. Then by Lemma 3.3, $u \in C^{\infty}((0,\infty))$ and 
$(L_0(\zeta) - i)u = 0$ holds in the classical sense. Picking $\nu = \exp(-\pi i/4)$, we have
$$
u = ay^{(n-1)/2}I_{\nu}(\zeta y) + by^{(n - 1)/2}K_{\nu}(\zeta y).
$$
Since $u \in L^2((1,\infty);dy/y^n)$, we have $a = 0$ by (\ref{eq:Chap1Sec3Inunearinfty}). Since ${\rm Re}\,\nu > 0$ and $u \in L^2((0,1);dy/y^n)$, we also have $b = 0$ by (\ref{eq:Chap1Sec3Knunear0}). When $\zeta = 0$, $u$ is written as 
$$
u = ay^{(n -1)/2 + \alpha - i\beta} + by^{(n -1)/2 - \alpha + i\beta}
, \quad 
\alpha, \beta > 0
$$
As above $a = 0$, since $u \in L^2((1,\infty));dy/y^n)$, and $b = 0$ since
$u \in L^2((0,1));dy/y^n)$. \qed


\subsection{Green function}
We construct the Green function of $L_0(\zeta) + \nu^2$ when ${\rm Re}\,\zeta > 0$. In the following we always assume that
\begin{equation}
\nu \not\in {\bf Z}, \quad {\rm Re}\,\nu \geq 0.
\nonumber
\end{equation}


\begin{definition} \label{DEF:3.5}
We put
\begin{equation}
G_0(y,y';\zeta,\nu) = \left\{
\begin{split}
& (yy')^{(n-1)/2}K_{\nu}(\zeta y)I_{\nu}(\zeta y'), \quad 
y > y' > 0, \\
& (yy')^{(n-1)/2}I_{\nu}(\zeta y)K_{\nu}(\zeta y'), \quad 
y' > y > 0
\end{split} 
\right.
\nonumber
\end{equation}
and define the integral operator $G_0(\zeta,\nu)$ by
\begin{equation}
G_0(\zeta,\nu)f(y) = \int_0^{\infty}G_0(y,y';\zeta,\nu)f(y')\frac{dy'}{(y')^n}.
\nonumber
\end{equation}
\end{definition}


\begin{lemma} $
(L_0(\zeta) + \nu^2)G_0(\zeta,\nu)f = f, \quad \forall f \in C_0^{\infty}((0,\infty))$.
\end{lemma}
Proof. Using the equality
\begin{equation}
I_{\nu}(z)K'_{\nu}(z) - I'_{\nu}(z)K_{\nu}(z) = - \frac{1}{z},
\nonumber
\end{equation}
we have
\begin{equation}
\begin{split}
& \left(y^{(n-1)/2}I_{\nu}(\zeta y)\right)
\left(y^{(n-1)/2}K_{\nu}(\zeta y)\right)'  \\
& - \left(y^{(n-1)/2}I_{\nu}(\zeta y)\right)'
\left(y^{(n-1)/2}K_{\nu}(\zeta y)\right) = - y^{n-2}.
\end{split}
\nonumber
\end{equation}
We then compute as in the proof of Lemma 3.2 (2). \qed


\begin{lemma}
The Green function $G_0(y,y';\zeta,\nu)$ is analytic with respect to $\zeta$ when ${\rm Re}\,\zeta > 0$, and the following inequalities hold.
\begin{equation}
|G_0(y,y';\zeta,\nu)| \leq C(yy')^{(n- 1)/2},
\label{eq:Chap1Sec3G01}
\end{equation}
\begin{equation}
|G_0(y,y';\zeta,\nu)| \leq \frac{C}{|\zeta|}(yy')^{(n- 2)/2},
\label{eq:Chap1Sec3G02}
\end{equation}
\begin{equation}
\Big|\frac{\partial}{\partial\zeta}G_0(y,y';\zeta,\nu)\Big| \leq \frac{C}{|\zeta|}(yy')^{(n- 2)/2}(y + y').
\label{eq:Chap1Sec3delzetaG0}
\end{equation}
Here the constant $C$  depends on $\nu$, but is independent of $\zeta$ when ${\rm Re}\,\zeta > 0$.
\end{lemma}
Proof. By virtue of (\ref{eq:Chap1Sec3Inunearinfty}) $\sim$ (\ref{eq:Chap1Sec3Knunear0}), we have
\begin{equation}
|I_{\nu}(z)| \leq C\left(\frac{|z|}{1 + |z|}\right)^{{\rm Re}\,\nu}
(1 + |z|)^{-1/2}e^{{\rm Re}\,z},
\label{eq:Chap1Sec3AbsInu}
\end{equation}
\begin{equation}
|K_{\nu}(z)| \leq C\left(\frac{|z|}{1 + |z|}\right)^{-{\rm Re}\,\nu}
(1 + |z|)^{-1/2}e^{-{\rm Re}\,z}
\label{eq:Chap1Sec3AbsKnu}
\end{equation}
Since $t/(1 + t)$ is monotone increasing for $t \geq 0$ , we have for $y > y' > 0$ 
\begin{equation}
|K_{\nu}(\zeta y)I_{\nu}(\zeta y')| \leq C\frac{e^{- {\rm Re}\,\zeta(y - y')}}{(1 + |\zeta y|)^{1/2}(1 + |\zeta y'|)^{1/2}}.
\nonumber
\end{equation}
Hence,
\begin{equation}
|G_0(y,y';\zeta,\nu)| \leq C(yy')^{(n-1)/2}
\frac{e^{-{\rm Re}\zeta|y - y'|}}{(1 + |\zeta y|)^{1/2}(1 + |\zeta y'|)^{1/2}},
\label{C1S3Goexpdecay}
\end{equation}
which implies (\ref{eq:Chap1Sec3G01}), (\ref{eq:Chap1Sec3G02}). By the following formulas
\begin{equation}
\label{3.14'}
2I_{\nu}'(z) = I_{\nu-1}(z) + I_{\nu+1}(z),
\end{equation}
\begin{equation}
- 2K_{\nu}'(z) = K_{\nu-1}(z) + K_{\nu+1}(z)
\nonumber
\end{equation}
(see e.g. \cite{MUH59} p. 173) and  (\ref{eq:Chap1Sec3Inunearinfty}) $\sim$ (\ref{eq:Chap1Sec3Knunear0}), we have
\begin{equation}
|zI_{\nu}'(z)| \leq C\left(\frac{|z|}{1 + |z|}\right)^{{\rm Re}\,\nu}
(1 + |z|)^{1/2}e^{{\rm Re}\,z},
\nonumber
\end{equation}
\begin{equation}
|zK_{\nu}'(z)| \leq C\left(\frac{|z|}{1 + |z|}\right)^{-{\rm Re}\,\nu}
(1 + |z|)^{1/2}e^{-{\rm Re}\,z}.
\nonumber
\end{equation}
Therefore we have
\begin{equation}
\Big|\frac{\partial}{\partial \zeta}I_{\nu}(\zeta y)\Big| \leq 
\frac{C}{|\zeta|}\left(\frac{|\zeta y|}{1 + |\zeta y|}\right)^{{\rm Re}\,\nu}
(1 + |\zeta y|)^{1/2}e^{{\rm Re}\,\zeta y}, 
\nonumber
\end{equation}
\begin{equation}
\Big|\frac{\partial}{\partial \zeta}K_{\nu}(\zeta y)\Big| \leq 
\frac{C}{|\zeta|}\left(\frac{|\zeta y|}{1 + |\zeta y|}\right)^{-{\rm Re}\,\nu}
(1 + |\zeta y|)^{1/2}e^{-{\rm Re}\,\zeta y}. 
\nonumber
\end{equation}
Using the straightforward inequality
\begin{equation}
\left(\frac{1 + |\zeta y'|}{1 + |\zeta y|}\right)^{1/2} \leq
\frac{y + y'}{(yy')^{1/2}},
\nonumber
\end{equation}
we obtain (\ref{eq:Chap1Sec3delzetaG0}). \qed

\bigskip
One can check that the constants $C$ in (\ref{eq:Chap1Sec3G01}) $\sim$ (\ref{eq:Chap1Sec3delzetaG0}) may be chosen independently of $\nu$ when $\nu$ varies over a compact set in $\{{\rm Re}\,\nu \geq 0\}\setminus{\bf Z}$.

\bigskip
We define ${\mathcal B}, {\mathcal B}^{\ast}$ by putting ${\bf h} = {\bf C}$ in \S 2.


\begin{lemma} We have
\begin{equation}
\|G_0(\zeta,\nu)f\|_{{\mathcal B}^{\ast}} \leq C\|f\|_{{\mathcal B}},
\nonumber
\end{equation}
where the constant $C$ is independent of $\nu$ when $\nu$ varies over a compact set in $\{{\rm Re}\,\nu \geq 0\}\setminus{\bf Z}$, and also of $\zeta$ when ${\rm Re}\,\zeta > 0$.
\end{lemma}
Proof. We put $u = G_0(\zeta,\nu)f$. By (\ref{eq:Chap1Sec3G01}), we have
$$
\frac{|u(y)|^2}{y^n} \leq \frac{C}{y}\left(\int_0^{\infty}
\frac{|f(y')|}{(y')^{1/2}}\frac{dy'}{(y')^{n/2}}\right)^2.
$$
Hence we have
\begin{eqnarray*}
\|u\|_{{\mathcal B}^{\ast}} &\leq& C\int_0^{\infty}
\frac{1}{(y')^{1/2}}\frac{|f(y')|}{(y')^{n/2}}dy' \\
&=& \sum_k\int_{I_k}\frac{1}{(y')^{1/2}}\frac{|f(y')|}{(y')^{n/2}}dy'\\
&\leq& \sum_k\left(\int_{I_k}\frac{dy}{y}\right)^{1/2}
\left(\int_{I_k}|f(y)|^2d\mu(y)\right)^{1/2} 
\leq C\|f\|_{{\mathcal B}}. 
\qed
\end{eqnarray*}


\begin{lemma} (1)
If  $u \in {\mathcal B}^{\ast}$ satisfies $(L_0(\zeta) - z)u = 0$ for $\zeta > 0,\ {\rm Im}\,z \neq 0$, then $u = 0$. \\
 \noindent
(2) If  $u \in L^2((0,\infty))$ satisfies $(L_0(\zeta) - t)u = 0$ for $\zeta > 0,\ t \in {\bf R}$, then $u = 0$.
\end{lemma}
Proof. We prove the assertion (1). Letting $\nu = \pm i{\sqrt z},\ {\rm Re}\,\nu > 0$, $u$ is written as 
$u = ay^{(n-1)/2}I_{\nu}(\zeta y) + by^{(n- 1)/2}K_{\nu}(\zeta y)$.
Since $u \in {\mathcal B}^{\ast}$, letting $y \to \infty$, we see that  $a = 0$. Letting
$y \to 0$, we also see $b = 0$. The assertion (2) is proved in a similar way. 
\qed


\begin{cor} \l
If $\zeta > 0, \ z = - \nu^2, \ {\rm Im}\,z \neq 0$,  then
\begin{equation}
G_0(\zeta,\nu) = (L_0(\zeta) - z)^{-1}
\end{equation}
holds, where the right-hand side is the resolvent of $L_0(\zeta)$ in $L^2((0, \infty); \frac{d y}{y^n}).$
\end{cor}


\subsection{Limiting absorption principle}
Let $X$ be a Banach space and $X^{\ast}$ its dual. A sequence $\{u_n\}_{n=1}^{\infty} \subset X^{\ast}$ is said to converge to $u \in X^{\ast}$ in $\ast$-weak sense if
$$
\langle u_n,v\rangle \to \langle u,v\rangle, \quad \forall v \in X.
$$


\begin{theorem}
For $\zeta > 0$, $\lambda > 0$, $f \in {\mathcal B}$,
\begin{equation}
(L_0(\zeta) - \lambda \mp i\epsilon)^{-1}f \to 
G_0(\zeta,\mp i\sqrt{\lambda})f, \quad \epsilon \to 0
\nonumber
\end{equation}
in $\ast$-weak sense.
\end{theorem}
Proof. We put $u(\nu) = G_0(\zeta,\nu)f$, where $\nu = -i\sqrt{\lambda + i\epsilon}$ for $\lambda + i\epsilon$, and $\nu = i\sqrt{\lambda - i\epsilon}$ for $\lambda - i\epsilon$.
By Corollary 3.10, $u(\nu) = (L_0(\zeta) - \lambda \mp i\epsilon)^{-1}f$. 
Since, by Lemma 3.8, $u(\nu)$ are bounded in ${\mathcal B}^{\ast}$, by
Lebesgue's convergence theorem 
$(u(\nu),g) \to (G_0(\zeta,\mp i\sqrt{\lambda})f,g)$, $\forall g \in C_0^{\infty}((0,\infty))$. As $C_0^{\infty}((0,\infty))$ is dense in $\mathcal B$, applying again Lemma 3.8 proves the theorem. \qed

\bigskip
In the following, we write
\begin{equation}
(L_0(\zeta) - \lambda \mp i0)^{-1} = G_0(\zeta,\mp i\sqrt{\lambda}).
\nonumber
\end{equation}
By Lemma 3.8, we have the following uniform, with respect to $\zeta > 0$, estimate
\begin{equation}
\sup_{\zeta \geq 0}\|(L_0(\zeta) - \lambda \mp i0)^{-1}\|_{{\bf B}(
{\mathcal B};{\mathcal B}^{\ast})} = C(\lambda) < \infty,
\label{eq:Chap1Sec3EstimteL0}
\end{equation}
where, for  $0 < a < b < \infty$,
\begin{equation}
\sup_{a < \lambda < b}C(\lambda) < \infty.
\label{eq:UniformEstimate}
\end{equation}
Later we will also  prove ({\ref{eq:UniformEstimate})  by using techniques from partial differential equations.


\subsection{Eigenfunction expansions} 


\begin{lemma}
For $\zeta > 0$, $\sigma(L_0(\zeta)) = [0,\infty)$ and $\sigma_p(L_0(\zeta)) = \emptyset.$
\end{lemma}
Proof. We have for $u \in C_0^{\infty}((0,\infty))$
\begin{equation}
(L_0(\zeta)u,u) + \frac{(n-1)^2}{4}\|u\|^2 = \zeta^2\int_0^{\infty}|u(y)|^2\frac{dy}{y^{n-2}} + \int_0^{\infty}|u'(y)|^2\frac{dy}{y^{n-2}}.
\nonumber
\end{equation}
By integration by parts and  Cauchy-Schwarz' inequality, we have
\begin{eqnarray*}
(n-1)\int_0^{\infty}|u|^2\frac{dy}{y^n} &=& 
2{\rm Re}\int_0^{\infty}\left(\partial_y u\right)\overline{u}\frac{dy}{y^{n-1}} \\
&\leq& 2\left(\int_0^{\infty}\frac{|\partial_y u|^2}{y^{n-2}}dy\right)^{1/2}
\left(\int_0^{\infty}\frac{|u|^2}{y^n}dy\right)^{1/2}.
\end{eqnarray*}
This implies that
\begin{equation}
\int_0^{\infty}\frac{|\partial_y u|^2}{y^{n-2}}dy \geq
\frac{(n- 1)^2}{4}(u,u).
\nonumber
\end{equation}
Therefore,
\begin{equation}
(L_0(\zeta)u,u) \geq \zeta^2\int_0^{\infty}|u(y)|^2\frac{dy}{y^{n-2}}.
\nonumber
\end{equation}
Therefore $\sigma(L_0(\zeta)) \subset [0,\infty)$. 

Let us recall that for $\lambda > 0$,
$(L_0(\zeta) - \lambda)\left[y^{(n-1)/2}K_{i\sqrt{\lambda}}(\zeta y)\right] = 0$. 
Take $\chi(t) \in 
C^{\infty}((0,\infty))$ such that $\chi(t) = 0 \ (t < 1), \ \chi(t) = 1 \ (t > 2)$, and put
$$
u_N(y) = \chi(Ny)y^{(n-1)/2}K_{i\sqrt{\lambda}}(\zeta y)
$$
By (\ref{eq:Chap1Sec3Knunear0})
\begin{equation}
\begin{split}
\|u_N\|^2 &= \int_0^{\infty}\chi\big(\frac{Nt}{\zeta}\big)|K_{i\sqrt{\lambda}}
(t)|^2\frac{dt}{t} \\
&\geq \int_1^{\infty}|K_{i\sqrt{\lambda}}(t)|^2\frac{dt}{t} + 
C\int_{\zeta/N}^1\frac{dt}{t} \\
&\geq C(\log N + 1).
\end{split}
\label{C1S3unnorm^2}
\end{equation}
We put $\varphi_N(y) = u_N(y)/\|u_N\|$. Then $\|\varphi_N\| = 1$, and
\begin{eqnarray*}
& &(L_0(\zeta) - \lambda)\varphi_N = \frac{1}{\|u_N\|}
\Big\{- (Ny)^2\chi''(Ny)y^{(n-1)/2}K_{i\sqrt{\lambda}}(\zeta y) \\
& & - 2Ny\chi'(Ny)y\partial_y\big(y^{(n-1)/2}K_{i\sqrt{\lambda}}(\zeta y)\big) 
 + (n-2)Ny\chi'(Ny)y^{(n-1)/2}K_{i\sqrt{\lambda}}(\zeta y)\Big\}.
\end{eqnarray*}
Taking into account (\ref{3.14'}) and  (\ref{C1S3unnorm^2}) and facts   that
$$
\int_0^{\infty}(Ny)^2\chi'(Ny)^2\frac{dy}{y} = \int_0^{\infty}t^2\chi'(t)^2
\frac{dt}{t} < \infty,
$$
and also $\int_0^{\infty}(Ny)^4\chi'(Ny)^2dy/y < \infty,\, \int_0^{\infty}(Ny)^4\chi''(Ny)^2dy/y < \infty,$
we have $\|(L_0(\zeta) - \lambda)\varphi_N\| \to 0$. By Weyl's method of singular sequence (see \cite{Is04a} p. 25), we have $\lambda \in \sigma(L_0(\zeta))$.
Lemma 3.9 proves that
$L_0(\zeta)$ has no eigenvalues. \qed

\bigskip
Let us recall Stone's formula (\cite{Is04a} p. 17). Let $H$ be a self-adjoint operator, $R(z) = (H - z)^{-1}$ the resolvent of $H$, $E_H(\lambda)$ the spectral decomposition for $H$. If $a, b \not\in \sigma_p(H)$, letting $I = (a,b)$, we have 
\begin{equation}
\begin{split}
 \left(E_H(I)f,g\right) 
 & = 
([E_H(b) - E_H(a)]f,g)  \\
&= \lim_{\epsilon\to 0}\frac{1}{2\pi i}\int_a^b
\left([R(\lambda + i\epsilon) - R(\lambda - i\epsilon)]f,g\right)d\lambda.
\end{split}
\label{eq:Chap1Sec3Stone}
\end{equation}

Using $K_{\nu}(z) = K_{-\nu}(z)$ and (\ref{eq:KnuandInu}, we have
\begin{equation}
K_{-\nu}(z)I_{-\nu}(z') - K_{\nu}(z)I_{\nu}(z') = 
\frac{2\sin(\nu\pi)}{\pi}K_{\nu}(z)K_{\nu}(z'), \quad 
\nu \not\in {\bf Z}.
\nonumber
\end{equation}
Therefore, the integral kernel of $(L_0(\zeta) - \lambda - i0)^{-1} - 
(L_0(\zeta) - \lambda + i0)^{-1}$ is given by
\begin{equation}
\frac{2i\sinh(\sqrt{\lambda}\pi)}{\pi}(yy')^{(n-1)/2}K_{i\sqrt{\lambda}}(\zeta y)K_{i\sqrt{\lambda}}(\zeta y').
\label{eq:Chap1Sec3KiKi}
\end{equation}
We now put, for $f \in C_0^{\infty}((0,\infty))$ and $k > 0$
\begin{equation}
\left({\mathcal F}_{\zeta}f\right)(k) = 
\frac{\big(2k\sinh(k\pi)\big)^{1/2}}{\pi}
\int_0^{\infty}y^{(n-1)/2}K_{ik}(\zeta y)f(y)\frac{dy}{y^n}.
\label{eq:Chap1Sec3Fzeta}
\end{equation}


\begin{theorem}
(1) ${\mathcal F}_{\zeta}$ is uniquely extended to a unitary operator from \\
\noindent $L^2((0,\infty);dy/y^n)$ to $L^2((0,\infty);dk)$. \\
\noindent
(2) If $f \in D(L_0(\zeta))$, then 
$\left({\mathcal F}_{\zeta}L_0(\zeta)f\right)(k) = 
k^2\left({\mathcal F}_{\zeta}f\right)(k).$ \\
\noindent
(3) For $f \in L^2((0,\infty);dy/y^n)$, the inversion formula
\begin{eqnarray}
f &=& {\mathcal F}_{\zeta}^{\ast}{\mathcal F}_{\zeta}f
\label{eq:FzetaastFzeta} \\
&=& y^{(n-1)/2}\int_0^{\infty}
\frac{\big(2k\sinh(k\pi)\big)^{1/2}}{\pi}
K_{ik}(\zeta y)({\mathcal F}_{\zeta}f)(k)dk
\label{eq:InversionFzeta}
\end{eqnarray} 
holds.
\end{theorem}
Proof. It follows from (\ref{eq:Chap1Sec3Stone}) and (\ref{eq:Chap1Sec3Fzeta}) that for $0 < a < b < \infty$ 
\begin{equation}
\left([E_{L_0(\zeta)}(b) - E_{L_0(\zeta)}(a)]f,g\right) = \int_{\sqrt{a}}^{\sqrt{b}}\left({\mathcal F}_{\zeta}f(k),{\mathcal F}_{\zeta}g(k)\right)dk,
\label{eq:Chap1Sec3Spectremeasure}
\end{equation}
where we have used
\begin{equation}
\overline{K_{ik}(y)} = K_{ik}(y) = K_{-ik}(y).
\label{C1S3KikK-ik}
\end{equation}
Letting $a \to 0, b \to \infty$, we see that ${\mathcal F}_{\zeta}$ is an isometry from  
$L^2(0,\infty);dy/y^n)$ to $L^2((0,\infty);dk)$. We show the surjectivity later. For $f \in C_0^{\infty}((0,\infty))$, by part integration, we have
$$
 \int_0^{\infty}y^{(n-1)/2}K_{ik}(\zeta y)\left(L_0(\zeta)f(y)\right)\frac{dy}{y^n}  = k^2 \int_0^{\infty}y^{(n-1)/2}K_{ik}(\zeta y)f(y)\frac{dy}{y^n}.
 $$
This proves (2), if we take into account the density of $C_0^{\infty}((0,\infty))$ in $D(L_0(\zeta))$ (see Theorem 3.4). 

The isometric property of $\mathcal F_{\zeta}$ entails (\ref{eq:FzetaastFzeta}). However, the integral formula (\ref{eq:InversionFzeta}) requires a subtle analysis. Since ${\mathcal F}_{\zeta}$ is bounded from $L^2((0,\infty);
dy/y^n)$ to  $L^2(0,\infty);dk)$,
for any $f \in L^2((0,\infty);dy/y^n)$ the strong limit
\begin{equation}
\lim_{a\to 0, b \to \infty}
\frac{(2k\sinh(k\pi))^{1/2}}{\pi}\int_{\sqrt{a}}^{\sqrt{b}}
y^{(n-1)/2}K_{ik}(\zeta y)f(y)\frac{dy}{y^n} =:
\left({\mathcal F}_{\zeta}f\right)(k)
\nonumber
\end{equation}
exists in $L^2((0,\infty);dk)$. To study the inverse transformation, we define an operator ${\mathcal F}_{\zeta}(k)$ by
\begin{equation}
{\mathcal F}_{\zeta}(k)f = \left({\mathcal F}_{\zeta}f\right)(k) \quad
{\rm for} \quad k > 0 \quad {\rm  and} \quad f \in C_0^{\infty}((0,\infty)).
\nonumber
\end{equation}


\begin{remark}
In the following we often use such a notation. Namely,
let a given be an operator ${\mathcal F}$  from a Hilbert space ${\mathcal H}$ to another Hilbert space $L^2((0,\infty) ; {\bf h} ; dk)$, where $\bf h$ is an auxiliary Hilbert space. For $k > 0$ we define an operator ${\mathcal F}(k)$ from a suitable subspace $S$ of ${\mathcal H}$ to ${\bf h}$ by
$$
{\mathcal F}(k)f = ({\mathcal F}f)(k), \quad 
f \in S.
$$
Conversely if we are given a family of operators $\{{\mathcal F}(k)\}_{k > 0}, $ with range in $\bf h$, we define an operator ${\mathcal F}$ with range in $L^2((0,\infty) ; {\bf h} ; dk)$ by the above formula.
\end{remark}


\begin{lemma} For any $k > 0$, there exists a constant $0 < C(k) < \infty$ such that
\begin{equation}
\sup_{\zeta > 0}\|{\mathcal F}_{\zeta}(k)\|
_{{\bf B}({\mathcal B};{\bf C})} \leq C(k),
\nonumber
\end{equation}
where $C(k)$ is uniformly bounded on any compact in $(0,\infty)$.
\end{lemma}
Proof. Using Lemma 3.8 and Theorem 3.11, and differentiating (\ref{eq:Chap1Sec3Stone}) and (\ref{eq:Chap1Sec3Spectremeasure})
by $b$, we have, in view of (\ref{eq:Chap1Sec3KiKi}),
\begin{equation}
|{\mathcal F}_{\zeta}(k)f|^2 = \frac{k}{i\pi}\left(\big[(L_0(\zeta) - k^2 - i0))^{-1} - (L_0(\zeta) - k^2 + i0))^{-1}\big]f,f\right). 
\nonumber
\end{equation}
Using (\ref{eq:Chap1Sec3EstimteL0}), we prove the lemma. \qed

\medskip
By (\ref{eq:Chap1Sec3Fzeta}), $\mathcal F_{\zeta}(k)^{\ast}$ is simply a multiplication operator :
$$
{\bf C} \ni \alpha \to \frac{(2k\sinh(k\pi))^{1/2}}{\pi}y^{(n-1)/2}K_{ik}(\zeta y)\alpha.
$$
Lemma 3.15 implies
\begin{equation}
\sup_{\zeta > 0}\|{\mathcal F}_{\zeta}(k)^{\ast}\|
_{{\bf B}({\bf C} ; {\mathcal B}^{\ast})} \leq C(k),
\nonumber
\end{equation}
By (\ref{eq:UniformEstimate}), this $C(k)$ is bounded when $k$ varies over a compact set in $(0,\infty)$.
Hence, for any $g \in L^2((0,\infty);dk)$,
\begin{equation}
\int_{1/N}^N{\mathcal F}_{\zeta}(k)^{\ast}g(k)dk \in {\mathcal B}^{\ast}, \quad \forall N > 0.
\nonumber
\end{equation}
Letting $\chi_N(\lambda)$ be the characteristic function of $(1/N,N)$, we have for $h \in C_0^{\infty}((0,\infty))$
\begin{equation}
\Big(\int_{1/N}^N{\mathcal F}_{\zeta}(k)^{\ast}g(k)dk,h\Big) 
= \int_{1/N}^Ng(k)\overline{\big({\mathcal F}_{\zeta}(k)h\big)}dk = (\chi_Ng,{\mathcal F}_{\zeta}h).
\label{C1S3innerproductfzeta}
\end{equation}
Here the left-hand side is the coupling between ${\mathcal B}^{\ast}$ and ${\mathcal B}$, the right-hand side is the inner product of $L^2((0,\infty);dk)$. However, since $\mathcal F_{\zeta}$ is an isometry between $L^2((0,\infty);dy/y^n)$ and $L^2((0,\infty);dk)$, the right-hand side makes sense for all
$h \in L^2((0, \infty); d y/y^n)$ Thus,
the left-hand side can be extended by continuity to $h \in L^2((0,\infty));dy/y^n)$.
This implies, by Riesz' theorem, that 
$$
\int_{1/N}^N{\mathcal F}_{\zeta}(k)^{\ast}g(k)dk 
= \mathcal F_{\zeta}^{\ast}(\chi_Ng) \in 
L^2((0,\infty);{dy}/{y^n}).
$$
Since $\mathcal F_{\zeta}^{\ast}$ is partial isometry, 
in the sense of strong convergence in
$L^2((0,\infty);dy/y^n)$,
$$
\lim_{N\to\infty}\int_{1/N}^N{\mathcal F}_{\zeta}(k)^{\ast}g(k)dk = {\mathcal F}_{\zeta}^{\ast}g
$$
holds. Taking $g={\mathcal F}_{\zeta} f$ and using again that ${\mathcal F}_{\zeta}$ is a partial
isometry, we see that,
 in the sense of strong convergence in $L^2((0,\infty);dy/y^n)$,
\begin{equation}
f = \lim_{N\to\infty}\int_{1/N}^N{\mathcal F}_{\zeta}(k)^{\ast}
\left({\mathcal F}_{\zeta}f\right)(k)dk.
\nonumber
\end{equation}
This is the meaning of the inversion formula (\ref{eq:InversionFzeta}).

\bigskip
Let us prove the surjectivity of ${\mathcal F}_{\zeta}$. Denote by $C_0((0, \infty))$
the class of continuous functions with compact support in $(0, \infty)$.


\begin{lemma} \label{L:3.16}
For $f \in C_0((0,\infty))$
\begin{equation}
{\mathcal F}_{\zeta}(k)f = C_{\pm}(k)
\lim_{y\to 0}y^{-(n-1)/2\pm ik}(L_0(\zeta)- k^2 \mp i0)^{-1}f,
\nonumber
\end{equation}
\begin{equation}
C_{\pm}(k) = \frac{1}{\pi}\left(\frac{\zeta}{2}\right)^{\pm ik}
\Gamma(1 \mp ik)(2k\sinh(k\pi))^{1/2}.
\nonumber
\end{equation}
\end{lemma}
Proof. By the definition of Green's function, it follows from the asymptotics
(\ref{eq:Chap1Sec3Inunear0}) that, for small $y >0$,
\begin{equation}
G_0(\zeta,\mp i\sqrt{k})f(y) \sim 
\frac{(\zeta/2)^{\mp ik}}{\Gamma(1 \mp ik)}
y^{(n-1)/2 \mp ik}\int_0^{\infty}(y')^{(n-1)/2}
K_{ik}(\zeta y')f(y')\frac{dy'}{(y')^n},
\nonumber
\end{equation}
from which the lemma follows. \qed


\begin{lemma}
Suppose $u \in {\mathcal B}^{\ast}$ satisfies $(L_0(\zeta) - k^2)u = 0$ for $\zeta > 0, k > 0$ and $\lim_{y\to 0}y^{-(n-1)/2+ik}u$ exists. Then $u = 0$.
\end{lemma}
Proof. Since $u$ is written as $u = ay^{(n-1)/2}I_{ik}(\zeta y) + 
by^{(n-1)/2}I_{-ik}(\zeta y)$,
$$
y^{-(n-1)/2+ik}u \sim ac_+(k)y^{2ik} + bc_-(k) \quad 
{\rm as} \quad y \to 0
$$
with constants $c_{\pm}(k) \neq 0$. If the limit of the right-hand side exists, $a = 0$. Hence $u = 
by^{(n-1)/2}I_{-ik}(\zeta y)$. Looking at the behavior as $y \to \infty$, we have $b = 0$. \qed


\begin{lemma}
(1) Suppose $\zeta > 0, k > 0$, and $f \in C_0((0,\infty)), u \in {\mathcal B}^{\ast}$ satisfy $(L_0(\zeta) - k^2)u = f$. Furthermore assume that as $y \to 0$, $u \sim Cy^{(n-1)/2 - ik}$. Then $u = (L_0(\zeta) - k^2 - i0)^{-1}f$.
\\
\noindent
(2) Suppose $\zeta > 0, k > 0$, and $f \in C_0((0,\infty)), u \in {\mathcal B}^{\ast}$ satisfy $(L_0(\zeta) - k^2)u = f$. Furthermore assume that as $y \to 0$, $u \sim Cy^{(n-1)/2 + ik}$. Then $u = (L_0(\zeta) - k^2 + i0)^{-1}f$.\end{lemma}
Proof. By Theorem 3.11, $(L_0(\zeta) - k^2 - i0)^{-1}f \in \mathcal B^{\ast}$ and behaves like $Cy^{(n-1)/2-ik}$ near $0$. To prove (1), we put $u - (L_0(\zeta) - k^2 - i0)^{-1}f = v$, and apply the previous lemma. Taking the complex conjugate of (1), we obtain (2). \qed


\begin{lemma}
${\rm Ran}\,{\mathcal F}_{\zeta} = L^2((0,\infty);dk)$. 
\end{lemma}
Proof. For $\psi(k) \in L^1_{loc}((0,\infty))$, let $\mathcal L(\psi)$ be the set of Lebesgue points of $\psi$, i.e. the set of $\ell > 0$ such that 
$$
\psi(\ell) = \lim_{\epsilon\to 0}\frac{1}{2\epsilon}
\int_{\ell-\epsilon}^{\ell+\epsilon}\psi(k)dk.
$$
It is well-known that $(0,\infty)\setminus\mathcal L(\psi)$ is measure 0 for any $\psi \in L^1_{loc}((0,\infty))$.
Let  $\varphi(k) \in L^2((0,\infty); dk)$ be othogonal to the range of $\mathcal F_{\zeta}$, and take 
$$
\ell \in {\mathcal L}(\varphi(k))\cap{\mathcal L}(|\varphi(k)|^2).
$$
 We take $\chi(y) \in C^{\infty}((0,\infty))$, $\chi(y) = 1 \ (y < 1), \chi(y) = 0 \ (y > 2)$, and  put
$$
u_{\ell}(y) = \chi(y)y^{(n-1)/2}I_{i\ell}(\zeta y),
$$
$$
g_{\ell}(y) = (L_0(\zeta) - \ell^2)u_{\ell} = [L_0(\zeta),\chi]I_{i\ell}(\zeta y).
$$
 Since $g_{\ell}(y) \in C_0^{\infty}((0,\infty))$, $u_{\ell} = (L_0(\zeta) - {\ell}^2 + i0)^{-1}g_{\ell}$ by Lemma 3.18. The formula (\ref{eq:Chap1Sec3Fzeta}) and Lemma 3.16 imply that ${\mathcal F}_{\zeta}(k)g_{\ell} =: C(k)$ is a continuous function of $k > 0$ such that $C(\ell) \neq 0$. 
For the characteristic function $\chi_I$ of an interval $I \subset (0,\infty)$, we have 
$$
(\mathcal F_{\zeta}\chi_I(L_0(\zeta))g_{\ell})(k) = \chi_I(k^2)(\mathcal F_{\zeta}g_{\ell})(k) = \chi_I(k^2)C(k),
$$
which implies
\begin{equation}
\int_I\varphi(k)\overline{C(k)}dk = 0
\nonumber
\end{equation}
for any  interval $I \subset (0,\infty)$. We then have
\begin{equation}
\begin{split}
\varphi(\ell)\overline{C(\ell)} &= \varphi(\ell)\overline{C(\ell)} - \frac{1}{2\epsilon}\int_{\ell-\epsilon}^{\ell+\epsilon}\varphi(k)\overline{C(k)}dk  \\
&= \overline{C(\ell)}\left(\varphi(\ell) - \frac{1}{2\epsilon}\int_{\ell-\epsilon}^{\ell+\epsilon}\varphi(k)dk\right)-
\frac{1}{2\epsilon}\int_{\ell-\epsilon}^{\ell+\epsilon}\varphi(k)\left(\overline{C(k)}- \overline{C(l)}\right)dk.
\end{split}
\nonumber
\end{equation}
When $\epsilon \to 0$, the 1st term of the right-hand side tends to 0
since $\ell \in L(\varphi(k))$. The 2nd term also tends to 0 by the Schwarz inequality,
\begin{equation}
\begin{split}
& \left|\frac{1}{2\epsilon}\int_{\ell-\epsilon}^{\ell+\epsilon}\varphi(k)\left(\overline{C(k)}- \overline{C(l)}\right)dk\right| \\
& \leq \left(\frac{1}{2\epsilon}\int_{\ell-\epsilon}^{\ell+\epsilon}|\varphi(k)|^2dk\right)^{1/2}\times
\left(\frac{1}{2\epsilon}\int_{\ell-\epsilon}^{\ell+\epsilon}|C(k)-C(\ell)|^2dk\right)^{1/2},
\end{split}
\nonumber
\end{equation}
  the assumption that $\ell \in \mathcal L(|\varphi(k)|^2)$, and continuity of $C(k)$. Therefore $\varphi(\ell) = 0$, which proves the lemma due to the density of $L(\varphi(k)) \cap L(|\varphi(k)|^2)$. 
\qed


\subsection{Kontrovich-Lebedev's inversion formula}
By ${\mathcal F}_{\zeta}^{\ast}{\mathcal F}_{\zeta} = 1$,
\begin{equation}
f(y) = \int_0^{\infty}\!\!\int_0^{\infty}\frac{2\sigma\sinh(\sigma\pi)}{\pi^2}
(yy')^{-1/2}K_{i\sigma}(y)K_{i\sigma}(y')f(y')dy'd\sigma,
\nonumber
\end{equation}
and from ${\mathcal F}_{\zeta}{\mathcal F}_{\zeta}^{\ast} = 1$,
\begin{equation}
g(\sigma) = \int_0^{\infty}\!\!\int_0^{\infty}
\frac{2(\tau\sigma)^{1/2}\left(\sinh(\sigma\pi)\sinh(\tau\pi)\right)^{1/2}}{\pi^2}
\frac{K_{i\sigma}(y)K_{i\tau}(y)}{y}g(\tau)d\tau dy,
\nonumber
\end{equation}
which are called  Kontrovich-Lebedev's inversion formulae. The convergence of the integral in $L^2$ is proven above. Conditions for the pointwise convergence are given in \cite{Le72} p. 132.



\section{The upper-half space model}
\subsection{Laplace-Beltrami operator}
We return to the upper-half space model (\ref{eq:upperhalfspace})
with the Riemannian metric (\ref{eq:Riemannianmetric}).
The volume element is $dxdy/y^n$. Therefore,
\begin{equation}
L^2({\bf H}^n) = L^2({\bf R}^n_+ ;\frac{dxdy}{y^n}).
\nonumber
\end{equation}
The Laplace-Beltrami operator is given by
\begin{equation}
 - \Delta_g =   y^2(- \partial_y^2 - \Delta_x) + (n - 2)y\partial_y,
\quad
\Delta_x = \sum_{i=1}^{n-1}(\partial/\partial x_i)^2.
\nonumber
\end{equation}
We put 
\begin{equation}
H_0 = - \Delta_g - \frac{(n-1)^2}{4}.
\nonumber
\end{equation}
The partial Fourier transform $\hat f(\xi,y)$ of $f(x,y)$ is defined by
\begin{equation}
F_0f(\xi,y) = \hat f(\xi,y) = (2\pi)^{-(n-1)/2}\int_{{\bf R}^{n-1}}e^{-ix\cdot\xi}f(x,y)dx.
\nonumber
\end{equation}
Letting $L_0(\zeta)$ be as in (\ref{eq:DiffOpL0zeta}), we have
\begin{equation}
\widehat{(H_0 f)}(\xi,y) = \left(L_0(|\xi|)\hat f(\xi,\cdot)\right)(y).
\nonumber
\end{equation}


\begin{lemma}
$H_0\big|_{C_0^{\infty}({\bf R}^n_+)}$ is essentially self-adjoint.
\end{lemma}
Proof. We have only to prove that, for $u \in L^2({\bf H}^n)$,
\begin{equation}
\left((H_0 - i)\varphi,u\right) = 0 \quad \forall \varphi \in C_0^{\infty}({\bf R}^n_+) \Longrightarrow u = 0,
\nonumber
\end{equation}
and the same assertion with $i$ replaced by $- i$. Passing to the partial Fourier transform
and choosing $\varphi(x, y)=\varphi_x(x) \varphi_y(y),\,\varphi_x \in C^\infty_0({\bf R}^{n-1}),\,
\varphi_y \in C^\infty_0((0, \infty))$, for almost all $\xi \in {\bf R}^{n-1}$, 
we have 
$$
\left((L_0(|\xi|) - i)\varphi_y(y),  \hat u(\xi,y) \right)_{L^2((0, \infty); dy/y^n)}=  0.
$$
By the result for the 1-dimensional case (Theorem 3.4), we have $\hat u(\xi,y) = 0$. \qed


\subsection{Limiting absorption principle and  Fourier transform} We put
\begin{equation}
R_0(z) = (H_0 - z)^{-1}, \quad z \in {\bf C}\setminus{\bf R},
\nonumber
\end{equation}
and define the spaces ${\mathcal B}, {\mathcal B}^{\ast}$ by taking ${\bf h} = L^2({\bf R}^{n-1};dx)$ in Subsection 2.1.


\begin{theorem} \label{Th:4.2}
(1) $\ \sigma(H_0) = [0,\infty).$ \\
\noindent
(2) $\ \sigma_p(H_0) = \emptyset.$ \\
\noindent
(3) For $\lambda > 0$ and $f \in {\mathcal B}$, the following limits exist in ${\mathcal B}^{\ast}$ in the weak $\ast$-sense
$$
\lim_{\epsilon \to 0} R_0(\lambda \pm i \epsilon)f =: 
R_0(\lambda \pm i0)f,
$$
and the following inequality holds
\begin{equation}
\|R_0(\lambda \pm i0)f\|_{{\mathcal B}^{\ast}} \leq C\|f\|_{\mathcal B},
\label{eq:ResolventEstimate}
\end{equation}
where the constant $C$ does not depend on $\lambda$ if it varies over a compact set in $(0,\infty)$. \\
\noindent
(4) We put for $k > 0$, $k^2 = \lambda$, $f \in C_0^{\infty}({\bf R}^n_+)$,
\begin{equation}
\begin{split}
\left({\mathcal F}^{(\pm)}_0(k)f\right)(x) & = \frac{\big(2k\sinh(k\pi)\big)^{1/2}}{\pi}
(2\pi)^{-(n-1)/2} \\
& \times \iint\limits_{{\bf R}^{n-1}\times(0,\infty)}
e^{ix\cdot\xi} 
\Big(\frac{|\xi|}{2}\Big)^{\mp ik}
y^{(n - 1)/2}K_{ik}(|\xi|y)\widehat f(\xi,y)
\frac{d\xi dy}{y^n}.
\end{split}
\label{eq:Chap1Sec4FormulaF0k}
\end{equation}
Then we have
\begin{equation}
\frac{k}{\pi i}\left([R_0(k^2 + i0) - R_0(k^2 - i0)]f,f\right) =  \|{\mathcal F}^{(\pm)}_0(k)f\|^2_{L^2({\bf R}^{n-1})},
\label{eq:Chap1Sec4Parseval}
\end{equation}
and
\begin{equation}
\|{\mathcal F}^{(\pm)}_0(k)f\|_{L^2({\bf R}^{n-1})} \leq C\|f\|_{\mathcal B},
\label{eq:Chap1Sec4F0estimate}
\end{equation}
where the constant $C$ is independent of $k$ if it varies over a compact set in $(0,\infty)$. \\
\noindent
(5) We put $({\mathcal F}^{(\pm)}_0f)(k) = {\mathcal F}^{(\pm)}_0(k)f$. Then ${\mathcal F}^{(\pm)}_0$ is uniquely extended to a unitary operator from $L^2({\bf H}^n)$ to $L^2((0,\infty);L^2({\bf R}^{n-1});dk)$. For $f \in D(H_0)$, we have
\begin{equation}
({\mathcal F}^{(\pm)}_0H_0f)(k) = k^2({\mathcal F}^{(\pm)}_0f)(k).
\label{eq:Chap1Sec4Eigenoperator}
\end{equation}
\end{theorem}
Proof. (1) Since Lemma 3.12 implies $\sigma(L_0(|\xi|) )= [0,\infty)$, for
$z \not\in [0,\infty)$ the operator
\begin{equation}
(2\pi)^{-(n-1)/2}\int_{{\bf R}^{n-1}}e^{ix\cdot\xi}\Big(
\big(L_0(|\xi|) - z\big)^{-1}\widehat f(\xi,\cdot)\Big)(y)d\xi
\label{C1S4Resolventndim}
\end{equation}
is bounded on $L^2((0,\infty);L^2({\bf R}^{n-1});y^{-n}dy)$ and is equal to $R_0(z)$. Therefore $\sigma(H_0) \subset [0,\infty)$. The converse inclusion relation is proved by the method of singular sequence as in Lemma 3.12. Namely we take $\chi \in C_0^{\infty}({\bf R})$ such that $\chi(t) = 1 \ (|t| < 1), \ \chi(t) = 0 \ (|t| > 2)$, and normalize
$$
\chi\big(\frac{|x|}{N}\big)\chi\big(\frac{\log y}{\log N}\big)
e^{ix\cdot\xi}y^{(n-1)/2}K_{i\sqrt{\lambda}}(|\xi|y).
$$
We omit the computation.

(2) If there exists an $L^2$-solution of $(H_0 - \lambda)u = 0$, we have $(L_0(|\xi|) - \lambda)\widehat u(\xi,\cdot) = 0$, where, for almost all $\xi$,  $\widehat u(\xi,\cdot) \in L^2((0, \infty); dy/y^n)$.
 Lemma 3.9 yields $\widehat u(\xi,y) = 0$.

(3) We shall prove this statement in Chap. 2, \S 2 (see Lemma 2.2.9). In this section we confine ourselves to $f \in L^{2,s}$, $\forall s > 1/2$. We start with estimates
\begin{equation}
\|R_0(\lambda \pm i0)f\|_{{\mathcal B}^{\ast}} \leq C_s\|f\|_s, 
\label{eq:Chap1Sec4R0lambda}
\end{equation}
where the constant $C_s$ is independent of $\lambda$ when $\lambda$ varies over a compact set in $(0,\infty)$ and $\|\cdot\|_s$ is the norm in  Definition 2.6 with ${\bf h} = L^2({\bf R}^{n-1};dx)$.
Observe that
$$
\sup_{R>e}\frac{1}{\log R}\int_{1/R}^R\left[\int_{{\bf R}^{n-1}}|F(\xi,y)|^2d\xi\right]\frac{dy}{y^n} \leq 
\int_{{\bf R}^{n-1}}\left[\sup_{R>e}\frac{1}{\log R}\int_{1/R}^R|F(\xi,y)|^2\frac{dy}{y^n}\right]d\xi.
$$
Taking $F(\xi,y) = (L_0(\xi) - \lambda \mp i0)^{-1}\widehat f (\xi,y)$
and using  (\ref{eq:Chap1Sec3EstimteL0}), (\ref{eq:UniformEstimate} ), and  Lemmata 2.3 and  2.7 
\begin{eqnarray*}
& & \|R_0(\lambda \pm i0)f\|_{{\mathcal B}^{\ast}}^2 
\leq
\int_{{\bf R}^{n-1}}\|(L_0(|\xi|) - \lambda \mp i0)^{-1}\widehat f(\xi,\cdot)\|^2_{{\mathcal B}^{\ast}}d\xi \\
& &\leq C 
\int_{{\bf R}^{n-1}}
\|\widehat f(\xi,\cdot)\|^2_{{\mathcal B}}
d\xi 
\leq C_s 
\int_{{\bf R}^{n-1}}
\|\widehat f(\xi,\cdot)\|^2_{s}
d\xi = C_s\|f\|^2_s, 
\end{eqnarray*}
which proves (\ref{eq:Chap1Sec4R0lambda}).

Returning to formula (\ref{C1S4Resolventndim}), where $\hat f \in C_0^{\infty}({\bf H}^n)$ and using Theorem 3.11, we see that there exist limits $R_0(\lambda \pm i0)f = \lim_{\epsilon \to 0}R_0(\lambda \pm i\epsilon)f$. Using (\ref{eq:Chap1Sec4R0lambda}), we extend them to $f \in L^{2,s}$.

(4) The equality (\ref{eq:Chap1Sec4Parseval}) follows from (\ref{eq:Chap1Sec3Spectremeasure}), which together with (\ref{eq:ResolventEstimate}) proves (\ref{eq:Chap1Sec4F0estimate}).

(5) Taking into account of the 1-dimensional result, we have only to prove the unitarity. 
Restricting ourselves to ${\mathcal F}_0^{(-)}$, we obtain
by the Parseval formula (\ref{eq:Chap1Sec4Parseval})  that
 ${\mathcal F}_0^{(-)}$ is isometric. We take $\varphi(k,x) \in L^2((0,\infty)\times{\bf R}^{n-1})$,  $\chi(y) \in C^{\infty}(0,\infty)$ such that $\chi(y) = 1 \ (y < 1), \ \chi(y) = 0 \ (y > 2)$, and put
$$
 u_l(x,y) = \chi(y)y^{(n-1)/2}F_0^{\ast}\left[\Big(\frac{|\xi|}{2}\Big)^{-il}I_{il}(|\xi|y)\widehat\varphi(l,\xi)\right],
$$
where for any $\psi$
\begin{equation}
F_0^{\ast}\psi = (2\pi)^{-(n-1)/2}\int_{{\bf R}^{n-1}}e^{ix\cdot\xi}
\psi(\xi)d\xi.
\label{eq:Chap1Sect4F0ast}
\end{equation} 
Let $(H_0 - l^2)u_l = f_l$. When $y \to 0$, $\xi \neq 0$, 
$$
 \widehat u_l(\xi,y) \sim \frac{1}{\Gamma(1 + il)} y^{(n-1)/2+il}
 \widehat\varphi(l,\xi).
$$
Since for any fixed $\xi \in {\bf R}^{n-1}$ we have $\hat f(\xi, \cdot) \in C^\infty_0((0, \infty))$,
$\hat u_l(\xi, \cdot) \in {\mathcal B}^{\ast}$,
 by Lemma 3.18, $\hat u_l(\xi, \cdot)=(L_0(|\xi|) -l^2-i0\big)^{-1}$ and
$u_l = R_0(l^2 - i0)f$. 
Therefore, by Lemma 3.16 ${\mathcal F}_0^{(-)}(l)f = 
C(l)\varphi(l,\cdot)$, with some constant $C(l) \neq 0$. Therefore by the same argument as in the proof of Lemma 3.19, ${\mathcal F}_0^{(-)}$ is onto.
 \qed


\subsection{Helmholtz equation} 
Theorem 4.2 implies
\begin{equation}
{\mathcal F}^{(\pm)}_0(k)^{\ast} \in {\bf B}(L^2({\bf R}^{n-1});{\mathcal B}^{\ast}),
\label{eq:Chap1Sec4F0kastB}
\end{equation}
\begin{equation}
\begin{split}
\left({\mathcal F}^{(\pm)}_0(k)^{\ast}\varphi\right)(x,y)& = 
\frac{\big(2k\sinh(k\pi)\big)^{1/2}}{\pi}
 \\
& \ \ \ \times F_0^{\ast}\left[
\left(\frac{|\xi|}{2}\right)^{\pm ik}
y^{(n-1)/2}
K_{ik}(|\xi|y)\widehat \varphi(\xi)\right],
\end{split}
\label{eq:Chap1Sect4F0kastpsi}
\end{equation}
and by (\ref{eq:Chap1Sec4Eigenoperator}) in the weak sense
\begin{equation}
(H_0 - k^2){\mathcal F}^{(\pm)}_0(k)^{\ast}\varphi = 0, \quad 
\forall \varphi \in L^2({\bf R}^{n-1}).
\nonumber
\end{equation}
The aim of this subsection is to prove the following theorem 
(Modified {\it Poisson-Herglotz} formula).


\begin{theorem}
For $k > 0$
\begin{equation}
\{u \in {\mathcal B}^{\ast} ; (H_0 - k^2)u = 0\}
 = {\mathcal F}^{(\pm)}_0(k)^{\ast}\left(L^2({\bf R}^{n-1})
\right).
\nonumber
\end{equation}
\end{theorem}

Namely, any solution in $\mathcal B^{\ast}$ to the Helmholtz equation 
can be written as a Poisson integral of some $L^2$-function on the boundary at infinity. 
As will be shown later, the space $\mathcal B^{\ast}$ is, in some sense, 
the smallest  space for the solutions to the Helmholtz equation. 
Namely, recall the inclusion relations in Lemma 2.7. One can show that if $u \in L^{2,-1/2}$ satisfies the 
Helmholtz equation $(H_0 - k^2)u = 0$ for $k > 0$, then $u = 0$. Therefore, all the 
non-zero solutions to the Helmholtz equation decays at most like 
or slower than the functions in $\mathcal B^{\ast}$. The largest solution space was characterized 
by Helgason \cite{Hel70}, who proved that {\it all} solutions of the Helmholtz equation 
$(H_0 - \lambda)u = 0$ is written by a Poisson integral of a Sato's hyperfunction on the boundary. This result was extended to general symmetric spaces by \cite{Mine75}, \cite{KKMOOT78}. This was also extended to the Euclidean space using more general analytic functionals by \cite{HKMO72}.

In the Euclidean case, Theorem 4.3 was proved by Agmon-H{\"o}rmander \cite{AgHo76}. It was also extended to 2-body Schr{\"o}dinger operators by Yafaev \cite{Yaf91}, and for the 3-body problem by the author \cite{Is01}.

The proof of Theorem 4.3 requires a series of Lemmas.


\begin{lemma} \label{L:4.4} (A-priori estimate)

\noindent
(1) If $u \in {\mathcal B}^{\ast}$ satisfies
$(H_0 - z)u = f \in {\mathcal B}^{\ast}, z \in {\bf C}$, 
\begin{equation}
\|y\partial_y u\|_{{\mathcal B}^{\ast}} + 
\|y\partial_xu\|_{{\mathcal B}^{\ast}} 
\leq C(\|u\|_{{\mathcal B}^{\ast}} + \|f\|_{{\mathcal B}^{\ast}}).
\nonumber
\end{equation}
(2) If $u \in {\mathcal B}^{\ast}$ satisfies $(H_0 - z)u = f \in {\mathcal B}^{\ast}$ and
\begin{equation}
\lim_{R\to\infty}\frac{1}{\log R}\int_{1/R}^R\left[
\|u(y)\|^2_{L^2({\bf R}^{n-1})} + 
\|f(y)\|^2_{L^2({\bf R}^{n-1})} \right]\frac{dy}{y^n} = 0,
\nonumber
\end{equation}
we have
\begin{equation}
\lim_{R\to\infty}\frac{1}{\log R}\int_{1/R}^R
\left[\|y\partial_yu(y)\|^2_{L^2({\bf R}^{n-1})} + 
\|y\partial_xu\|^2_{L^2({\bf R}^{n-1})}\right]\frac{dy}{y^n} = 0.
\nonumber
\end{equation}
\end{lemma}
Proof. We put $D_y = y\partial_y, \ D_x = y\partial_x$. Then 
\begin{equation}
H_0 = - D_y^2 + (n - 1)D_y - D_x^2 - \frac{(n-1)^2}{4},
\nonumber
\end{equation}
and for $u, v \in C_0^{\infty}({\bf R}^n_+)$
\begin{equation}
(H_0u,v) = (D_yu,D_yv) + (D_xu,D_xv) - \frac{(n-1)^2}{4}(u,v).
\label{eq:Chap1Sec4H0uv}
\end{equation}
We pick $\rho \in C_0^{\infty}({\bf R})$ such that $\rho(t) = 1$ for $|t| < 1$, and put
$$
\rho_{r,R}(x,y) = \rho\left(\frac{|x|}{r}\right)\rho\left(
\frac{\log y}{\log R}\right), \quad
\rho_{R}(y) = \rho\left(\frac{\log y}{\log R}\right),
$$
for large parameters $r, R >> 1$. If $u \in {\mathcal B}^{\ast}$ satisfies $(H_0 - z)u = f \in {\mathcal B}$, 
we have, cf (\ref{eq:Chap1Sec4H0uv}),
\begin{equation}
(f,\rho_{r,R} u) = (D_xu,D_x(\rho_{r,R} u)) + (D_yu,D_y(\rho_{r,R} u)) - E(z)(u,\rho_{r,R} u),
\label{eq:Chap1Sec4frhorR}
\end{equation}
with $E(z) = (n-1)^2/4 + z$. Let us note that putting
$\widetilde \rho_x = D_x\rho_{r,R}, \ \widetilde \rho_y = D_y\rho_{r,R}$, we have
\begin{equation}
{\rm Re}\,(D_xu,\widetilde\rho_xu) = - \frac{1}{2}(u,(D_x\widetilde\rho_x)u),
\nonumber
\end{equation}
\begin{equation}
{\rm Re}\,(D_yu,\widetilde\rho_yu) = - \frac{1}{2}(u,y^n\big(\frac{\widetilde \rho_y}{y^{n-1}}\big)'u), \quad ' = \partial_y.
\nonumber
\end{equation}
We take the real part of (\ref{eq:Chap1Sec4frhorR}) and let $r \to \infty$. Since, pointwise
$$
D_x\widetilde\rho_x \to 0, \quad 
\widetilde\rho_y \to \frac{1}{\log R}\rho'\big(\frac{\log y}{\log R}\big),
$$
we obtain
$$
{\rm Re}\,(f,\rho_Ru) = (\rho_RD_xu,D_xu) + (\rho_RD_yu,D_yu) - 
\frac{1}{2}(u,\psi_Ru) - {\rm Re}\,E(z)(u,\rho_Ru),
$$
$$
\psi_R = y^n\partial_y\left(\frac{1}{y^{n-1}\log R}
\rho'\big(\frac{\log y}{\log R}\big)\right).
$$
Using Cauchy-Schwarz inequality and dividing by $\log R$, we obtain
\begin{equation}
\begin{split}
&\frac{1}{\log R}\int_{0}^{\infty}\left[(\rho_RD_xu,D_xu) + (\rho_RD_yu,D_yu)\right]\frac{dy}{y^n} \\
&\ \ \ \ \leq \frac{1}{\log R}\int_0^{\infty}\left[(\phi_Ru,u) + (\phi_Rf,f)\right]\frac{dy}{y^n},
\end{split}
\label{eq:IntegralrhoRDxuDxu}
\end{equation}
where $\phi_R$ has the form $\displaystyle{\phi_R(y) = C(R)
\phi(\frac{\log y}{\log R})}$ for some $\phi \in C_0^{\infty}({\bf R})$ and $C(R)$ 
is bounded on $(e,\infty)$. 
Taking the supremum with respect to $R$, we obtain,
by Lemma \ref{L:2.5},   the assertion (1).

 Letting $R \to \infty$ in (\ref{eq:IntegralrhoRDxuDxu}) and using Lemma 2.5 (1), we obtain (2).
\qed

\bigskip
\begin{lemma}
For $\varphi \in L^2({\bf R}^{n-1})$,
\begin{equation}
\lim_{R\to\infty}\frac{1}{\log R}\int_{1/R}^R
\|({\mathcal F}^{(\pm)}_0(k)^{\ast}\varphi)(\cdot,y)\|^2_{L^2({\bf R}^{n-1})}
\frac{dy}{y^n} = C\|\varphi\|^2_{L^2({\bf R}^{n-1})},
\nonumber
\end{equation}
where $C = C(k) > 0$.
\end{lemma}
Proof. 
By (\ref{eq:Chap1Sect4F0kastpsi}) and (\ref{eq:Chap1Sec3AbsKnu}) and  Lebesgue's convergence theorem, we have 
\begin{equation}
\begin{split}
& \frac{1}{y^{n-1}}\|({\mathcal F}^{(\pm)}_0(k)^{\ast}\varphi)(\cdot,y)\|^2_{L^2({\bf R}^{n-1})} = 
\tilde C(k)\int_{{\bf R}^{n-1}}|K_{ik}(|\xi|y)\varphi(\xi)|^2d\xi 
\\
&\leq \tilde C(k) \int_{{\bf R}^{n-1}}e^{-2|\xi|y}|\varphi(\xi)|^2d\xi. 
\end{split}
\nonumber
\end{equation}
Thus, 
$$
 \frac{1}{y^{n-1}}\|({\mathcal F}^{(\pm)}_0(k)^{\ast}\varphi)(\cdot,y)\|^2 \rightarrow 0, \quad 
 \hbox{as}\,\, y \to \infty.
 $$
 This implies that, as $R \to \infty,$
\begin{equation}
\frac{1}{\log R}\int_1^R\|({\mathcal F}^{(\pm)}_0(k)^{\ast}\varphi)(\cdot,y)\|^2_{L^2({\bf R}^{n-1})}
\frac{dy}{y^n} \to 0.
\label{eq:Chap1Sect4F0yplusinfty}
\end{equation}
To compute the limit as $y \to 0$, we first use  (\ref{eq:Chap1Sec3Knunear0}) to see that
\begin{eqnarray*}
  \frac{1}{y^{n-1}}\|({\mathcal F}^{(\pm)}_0(k)^{\ast}\varphi)(\cdot,y)\|^2_{L^2({\bf R}^{n-1})} 
&=& C(k)\int_{{\bf R}^{n-1}}|K_{ik}(|\xi|y)\varphi(\xi)|^2d\xi \\
\sim C(k)\|\varphi\|^2_{L^2({\bf R}^{n-1})} &+& {\rm Re}\left[C(\varphi)y^{-2ik}\right], 
\end{eqnarray*}
where $C(k) >0$ and
$$
C(\varphi) = C_0\int_{{\bf R}^n}|\xi|^{-2ik}|\varphi(\xi)|^2d\xi.
$$
Hence, 
$$
\frac{1}{\log R}\int_{1/R}^1\|({\mathcal F}^{(\pm)}_0(k)^{\ast}\varphi)(\cdot,y)\|^2_{L^2({\bf R}^{n-1})}
\frac{dy}{y^n} \to C(k)\|\varphi\|^2_{L^2({\bf R}^{n-1})}. \qed
$$

\bigskip
The above lemma and (\ref{eq:Chap1Sec4F0kastB}) imply the following corollary.


\begin{cor}
There exists a constant $C = C(k) > 0$ such that
\begin{equation}
C\|\varphi\|_{L^2({\bf R}^{n-1})} \leq \|{\mathcal F}^{(\pm)}_0(k)^{\ast}\varphi\|_{{\mathcal B}^{\ast}} \leq C^{-1}\|\varphi\|_{L^2({\bf R}^{n-1})}.
\nonumber
\end{equation}
\end{cor}

Next we show that the Fourier transform ${\mathcal F}^{(\pm)}_0(k)$ is derived from the asymptotic expansion of the resolvent as $y \to 0$, cf. Lemma \ref{L:3.16}.


\begin{lemma} \label{L:4.7}
For $f \in {\mathcal B}$ we put
\begin{equation}
u_{\pm} = R_0(k^2 \pm i0)f,
\nonumber
\end{equation}
\begin{equation}
v_{\pm}(x,y)  = \omega_{\pm}(k)y^{(n-1)/2 \mp ik}
\Big({\mathcal F}^{(\pm)}_0(k)f\Big)(x),
\nonumber
\end{equation}
\begin{equation}
\omega_{\pm}(k) = \frac{\pi}{\big(2k\sinh(k\pi)\big)^{1/2}
\Gamma(1 \mp ik)}
\label{eq:Chap1Sect4omegaplusminusk}
\end{equation}
Then we have as $R \to \infty$
\begin{equation}
\frac{1}{\log R}\int_{1/R}^1\|u_{\pm}(\cdot,y) - v_{\pm}(\cdot,y)\|^2_{L^2({\bf R}^{n-1})}
\frac{dy}{y^n} \to 0.
\nonumber
\end{equation}
\end{lemma}
Proof. First we show the lemma for $\hat f \in C_0^{\infty}({\bf R}^n_{+})$. Since
supp$\,\hat f$ is compact, we have as $y \to 0$ 
\begin{equation}
\begin{split}
\hat u_{\pm}(\xi,y) &= y^{(n-1)/2}I_{\mp ik}(|\xi|y)
\int_0^{\infty}(y')^{(n-1)/2}K_{ik}(|\xi|y')
\hat f(\xi,y')\frac{dy'}{(y')^{n}} \\
&\sim \omega_{\pm}(k)
y^{(n-1)/2 \mp ik}F_0{\mathcal F}^{(\pm)}_0(k)f.
\end{split}
\label{C1S4limitR0zinfty}
\end{equation}
It then follows from (\ref{eq:Chap1Sec4F0estimate}) and (\ref{eq:Chap1Sec3EstimteL0}) that

$$
\frac{1}{\log R}\int_{1/R}^1\|u_{\pm}(\cdot,y) - 
v_{\pm}(\cdot,y)\|^2_{L^2({\bf R}^{n-1})}
\frac{dy}{y^n} \to 0,
$$
as $R \to \infty$.
The general case is proved if we note that by (\ref{eq:ResolventEstimate}) and
(\ref{eq:Chap1Sec4F0estimate})
$$
\frac{1}{\log R}\int_{1/R}^1\|u_{\pm}(\cdot,y) - 
v_{\pm}(\cdot,y)\|^2_{L^2({\bf R}^{n-1})}
\frac{dy}{y^n} \leq C\|f\|^2_{{\mathcal B}},
$$
and approximate $f$ by $f_n$ with $\hat f_n \in C_0^{\infty}({\bf R}^n_+)$.  \qed

\bigskip
By the well-known formula
$$
 \Gamma(1 + s)\Gamma(1 - s) = s\Gamma(s)\Gamma(1 - s)
 = \frac{\pi s}{\sin(\pi s)},
$$
we have
\begin{equation}
 |\Gamma(1 + i\sigma)|^2 = \frac{\pi\sigma}{\sinh \pi\sigma}, \quad \sigma > 0,
 \label{eq:Chap1Sec4Gammasquuare}
\end{equation}
which implies
\begin{equation}
|\omega_{\pm}(k)|^2 = \frac{\pi}{2k^2}.
\label{eq:modulusomegasquaredpi2k}
\end{equation}
The function (\ref{eq:Chap1Sect4omegaplusminusk}) and the formulas (\ref{eq:Chap1Sec4Gammasquuare}), (\ref{eq:modulusomegasquaredpi2k}) will be used frequently throughout these notes.


\begin{cor} For $u_{\pm} = R_0(k^2 \pm i0)f$, with $f \in {\mathcal B}$, we have
\begin{equation}
\lim_{R\to\infty}\frac{1}{\log R}\int_{1/R}^1
\|u_{\pm}(\cdot,y)\|^2_{L^2({\bf R}^{n-1})}
\frac{dy}{y^n} = 
\frac{\pi}{2k^2}\|{\mathcal F}^{(\pm)}_0(k)f\|^2_{L^2({\bf R}^{n-1})},
\label{C1S4limnormupm}
\end{equation}
\begin{equation}
\lim_{R\to\infty}\frac{1}{\log R}\int_{1/R}^1
\|\big(y\partial_y - \frac{n-1}{2} \pm ik\big)u_{\pm}(\cdot,y)\|^2_{L^2({\bf R}^{n-1})}
\frac{dy}{y^n} =  0.
\label{C1S4radcond}
\end{equation}
\end{cor}
Proof. Let $u_{\pm}, v_{\pm}$ be as in the previous lemma, and denote them by $u, v$. Let $\|\cdot\| = \|\cdot\|_{L^2({\bf R}^{n-1})}$. Since $\|u\|^2 - \|v\|^2 = (u - v,u) + (v,u - v)$, we have
$\big|\|u\|^2 - \|v\|^2\big| \leq (\|u\| + \|v\|)\|u - v\|$. Thus, by (\ref{eq:ResolventEstimate}), 
(\ref{eq:Chap1Sec4F0estimate}) and Lemma \ref{L:4.7} that, as $R \to \infty$,
\begin{equation}
\begin{split}
 & \frac{1}{\log R} \left|\int_{1/R}^1(\|u\|^2 - \|v\|^2)\frac{dy}{y^n}\right|  \\ 
\leq& 
\frac{1}{\log R}  \left(\int_{1/R}^1(\|u\|^2 + \|v\|^2)\frac{dy}{y^n}\right)^{1/2}\times
\frac{1}{\log R}  \left(\int_{1/R}^1(\|u - v\|^2)\frac{dy}{y^n}\right)^{1/2} \to 0.
\end{split}
\nonumber
\end{equation}
We then obtain (\ref{C1S4limnormupm}) by using
$$
|\omega_{\pm}(k)|^2\|{\mathcal F}^{(\pm)}_0(k)f\|^2_{L^2({\bf R}^{n-1})} 
= \frac{1}{\log R}\int_{1/R}^1\|v_{\pm}(\cdot,y)\|^2_{L^2({\bf R}^{n-1})}
\frac{dy}{y^n}.
$$
Noting Lemma \ref{L:4.4} (1) and differentiating (\ref{C1S4limitR0zinfty}), we obtain (\ref{C1S4radcond}). 
\qed


\begin{lemma} \label{L:4.9}
For $f \in {\mathcal B}$, let $u = R_0(k^2 \pm i0)f$,  
$D_x = y\partial_x, \ D_y = y\partial_y$. Then we have
\begin{equation}
\lim_{R\to\infty}\frac{1}{\log R}\int_1^R\|u(\cdot,y)\|^2_{L^2({\bf R}^{n-1})}
\frac{dy}{y^n} = 0.
\label{eq:Chap1Sec4LogRunorm}
\end{equation}
\begin{equation}
\lim_{R\to\infty}\frac{1}{\log R}\int_1^R\big[
\|D_xu(\cdot,y)\|^2_{L^2({\bf R}^{n-1})} +
\|D_yu(\cdot,y)\|^2_{L^2({\bf R}^{n-1})}\big]\frac{dy}{y^n} = 0.
\label{eq:Chap1Sec4LogRunorm2}
\end{equation}
\end{lemma}
Proof. We first prove (\ref{eq:Chap1Sec4LogRunorm}) for $\widehat f \in C_0^{\infty}({\bf R}^n_+), \ 
u = R_0(k^2 - i0)f$. If
 $f(x,y) = 0$ for $y < C^{-1}$ and $y > C$, $\hat u(\xi,y)$ is written as for $y > C$ 
$$
\hat u(\xi,y) = y^{(n-1)/2}K_{i\sqrt{\lambda}}(|\xi|y)\int_{C^{-1}}^C
(y')^{(n-1)/2}\hat h(\xi,y')\frac{dy'}{(y')^n},
$$
where, due to (\ref{eq:Chap1Sec3AbsInu}), (\ref{C1S4Resolventndim}) and Definition \ref{DEF:3.5},
 $h \in L^2({\bf R}^n_+)$. Denoting 
$$
g(\xi) = \int_{C^{-1}}^C
(y')^{(n-1)/2}\hat h(\xi,y')\frac{dy'}{(y')^n},
$$
we have by  (\ref{eq:Chap1Sec3AbsKnu})
$$
|\hat u(\xi,y)| \leq Cy^{(n-1)/2}e^{-|\xi|y}g(\xi), \quad 
g \in L^2({\bf R}^{n-1}).
$$
Hence,
$$
\frac{1}{\log R}\int_1^R\|u(\cdot,y)\|_{L^2({\bf R}^{n-1})}^2\frac{dy}{y^n} 
\leq \frac{C}{\log R}\int_1^R\|e^{-|\xi|y}g(\xi)\|^2_{L^2({\bf R}^{n-1})}
\frac{dy}{y}.
$$
Therefore, (\ref{eq:Chap1Sec4LogRunorm}) for $\hat f \in C^\infty_0({\bf R}_+^n)$
follows from Lebesgue's convergence theorem. Taking note of 
$$
 \frac{1}{\log R}\int_1^R\|u(\cdot,y)\|_{L^2({\bf R}^{n-1})}^2\frac{dy}{y^n} 
\leq C\|f\|_{{\mathcal B}}^2,
$$
we have only to approximate $\hat f$ by  functions from $C_0^{\infty}({\bf R}^n_+)$ to prove 
(\ref{eq:Chap1Sec4LogRunorm})  for the general case.

We put
$$
\langle u,v\rangle = \int_{1}^{\infty}\big(u(\cdot),v(\cdot)\big)
d\mu, \quad d\mu = dy/y^n,
$$
where $(\cdot,\cdot)$ is the inner product of $L^2({\bf R}^{n-1})$. Take
$\rho \in C^{\infty}({\bf R})$ such that $\rho(t) = 0 \ (|t| > 3), \ \rho(t) = 1 \ 
(|t| < 2)$, and put $\rho_R(y) = \rho((\log y)/(\log R))$.
We multiply the equation $(H_0 - k^2)u = f$ by $\rho_R(y)\overline{u}$ and integrate by parts to see
\begin{eqnarray*}
& &\langle D_yu,\rho_RD_yu\rangle + 
\langle D_yu,y^n\big(\frac{\rho_R}{y^{n-1}}\big)'u\rangle + 
(D_yu,\rho_Ru)\big|_{y= 1}  \\
& &- \frac{n-1}{2}(u,\rho_Ru)\big|_{y=1} 
- \frac{n-1}{2}\langle u,y^n\big(\frac{\rho_R}{y^{n-1}}\big)'u
\rangle \\
& & + \langle D_xu,\rho_RD_xu\rangle - E(k^2)\langle u,\rho_Ru\rangle = 
\langle f,\rho_Ru\rangle.
\end{eqnarray*}
(We should insert $\rho(|x|/r)$, and let $r \to \infty$ using Theorem \ref{Th:4.2}(3) and 
Lemma \ref{L:4.4}(1)). 
We now put $\psi(t) = 1 \ (t < 3), \ \psi(t) = 0 \ (t > 4)$, $\psi_R(y) = 
\psi((\log y)/(\log R))$, and $\|\cdot\| = \|\cdot\|_{L^2({\bf R}^{n-1})}$ 
to obtain
\begin{eqnarray*}
& &\langle D_yu,\rho_RD_yu\rangle + \langle D_xu,\rho_RD_xu\rangle \\
&\leq& C(k)\Big(\int_{1}^{\infty}\psi_R(y)\|D_yu(y)\|\cdot\|u(y)\|d\mu 
+ \int_{1}^{\infty}\psi_R(y)\|u(y)\|^2d\mu  \\
& & + \int_{1}^{\infty}\psi_R(y)\|f(y)\|\cdot\|u(y)\|d\mu 
+ (\|D_yu(1)\| + \|u(1)\|)\|u(1)\|\Big).
\end{eqnarray*}
We divide both sides by $\log R$. Then the first term of the right-hande side is dominated from above by
\begin{equation}
\left(\frac{1}{\log R}\int_1^{\infty}\psi_R(y)\|D_yu\|^2d\mu\right)^{1/2}
\left(\frac{1}{\log R}\int_1^{\infty}\psi_R(y)\|u\|^2d\mu\right)^{1/2}.
\nonumber
\end{equation}
By Lemma 4.4 (1), we have
\begin{equation}
\sup_{R>2}\frac{1}{\log R}\int_1^{\infty}\psi_R(y)\|D_yu\|^2d\mu < \infty.
\nonumber
\end{equation}
Using (\ref{eq:Chap1Sec4LogRunorm}), we see that
$$
\lim_{R\to\infty}\frac{1}{\log R}\int_{1}^{\infty}\psi_R(y)\|u(y)\|^2d\mu = 0.
$$
Using the same considerations to estimate 
$\frac{1}{\log R} \int_{1}^{\infty}\psi_R(y)\|f(y)\|\cdot\|u(y)\|d\mu$,
we arrive at (\ref{eq:Chap1Sec4LogRunorm2}). \qed

\bigskip


\begin{lemma}
If $u \in {\mathcal B}^{\ast}, \ (H_0 - k^2)u = 0$, $f \in {\mathcal B}$, and  either ${\mathcal F}^{(+)}_0(k)f = 0$ or ${\mathcal F}^{(-)}_0(k)f = 0$ holds, then $(u,f) = 0$.
\end{lemma}
Proof. Assume that ${\mathcal F}^{(-)}_0(k)f = 0$. Take $\rho(t) \in C_0^{\infty}({\bf R})$ such that $\rho(t) = 1 \ (|t| < 1)$, and put 
$$
\rho_R(y) = \rho\big(\frac{\log y}{\log R}\big), \quad
\rho_{R,r}(y) = \chi\big(\frac{\log y}{\log R}\big)
\rho\big(\frac{\log y}{\log r}\big), \quad
\chi(t) = \int_{-\infty}^t\rho(s)ds.
$$
Letting $v = R_0(k^2- i0)f$, we then have
\begin{eqnarray*}
0 &=& (\rho_{R,r}(y)v,(H_0 - k^2)u) \\
&=& (\rho_{R,r}f,u) - 
((D_y^2\rho_{R,r})v,u) - 2((D_y\rho_{R,r})D_yv,u) 
+(n-1)((D_y\rho_{R,r})v,u).
\end{eqnarray*}
Let $r \to \infty$. Then, for any $R>0$ and sufficiently large $r$,
$$
\rho\left(\frac{\log y}{\log R}\right) \rho\left(\frac{\log y}{\log r}\right) =\rho\left(\frac{\log y}{\log R}\right).
$$
Using this formula, together with the fact that
$$
D_y \rho\left(\frac{\log y}{\log r}\right)=\frac{1}{\log r} \rho'\left(\frac{\log y}{\log r}\right),
$$
so that we obtain an extra factor $\frac{1}{\log r}$, we can use Lemma \ref{L:4.9} to show
that it is possible to replace $\rho_{R, r}$ in the above equation by $\chi_R(y)=\chi(\log y/\log R)$.
Thus,
\begin{equation}
(\chi_Rf,u) = 
((D_y^2\chi_R)v,u) + 2((D_y\chi_R)D_yv,u) - (n-1)((D_y\chi_R)v,u).
\label{eq:ChapSec4psiRf}
\end{equation}
Observe that, due to the assumption ${\mathcal F}^{(-)}_0(k)f = 0$,
it follows from Corollary 4.8 and Lemma 4.9  that
$$
 \frac{1}{\log R}\int_{1/R < y < R}\|v(\cdot,y)\|^2_{L^2({\bf R}^{n-1})}
 \frac{dy}{y^n} \to 0, \quad \hbox{as}\,\, R \to \infty.
$$
Since
$$
D_y\chi_R(y) = \frac{1}{\log R}\rho\big(\frac{\log y}{\log R}\big), \quad
D^2_y\chi_R(y) = \frac{1}{\log^2 R}\rho'\big(\frac{\log y}{\log R}\big) ,
$$
it then follows that the 1st and 3rd terms in the right-hand side of (\ref{eq:ChapSec4psiRf})
tend to $0$ as $R\to \infty$.  Integrating by parts in the 2nd term of the right-hand side of 
(\ref{eq:ChapSec4psiRf}) and using the fact that, by Lemma \ref{L:4.4}(1), $D_y u \in  \mathcal B^{\ast}$,
the same considerations show that this term also tends to $0$. Thus,  $(f,u) = 0$. \qed
 
\bigskip


\begin{lemma} Let $X, Y$ be Banach spaces, and $T \in {\bf B}(X,Y)$. Then the following 4 assertions are equivaent. \\
\noindent
(1) $\ {\rm Ran}\,T$  is closed. \\
(2) $\ {\rm Ran}\,T^{\ast}$ is closed. \\
(3) $\ {\rm Ran}\,T = N(T^{\ast})^{\perp} = 
\{y \in Y ; \langle y,y^{\ast}\rangle = 0 \ \forall y^{\ast} \in N(T^{\ast})\}$. \\
(4) $\ {\rm Ran}\,T^{\ast} = N(T)^{\perp} = 
\{x \in X^{\ast} ; \langle x,x^{\ast}\rangle = 0 \ \forall x \in N(T)\}$. 
\end{lemma}
　
For the proof, see e.g. \cite{Yo66} p. 205. 

\bigskip
\noindent
{\bf Proof of Theorem 4.3}. We put $X = {\mathcal B}, \ Y = L^2({\bf R}^{n-1}), \ T = {\mathcal F}^{(\pm)}_0(k)$ in the above lemma. By Corollary 4.6, ${\rm Ran}\,T^{\ast}$  is closed. Hence ${\rm Ran}\,T$ is closed. Corollary 4.6 also implies $N(T^{\ast}) = \{0\}$. Therefore ${\rm Ran}\,T = Y$, and ${\rm Ran}\,T^{\ast} = N(T)^{\perp}$. 
Lemma 4.10 shows that if $u \in {\mathcal B}^{\ast}$ and $(H_0 - k^2)u = 0$, then $u \in N(T)^{\perp}$. Therefore $u \in {\rm Ran}\,T^{\ast}$. 
\qed

\bigskip

\begin{cor}
$ \ \ {\mathcal F}^{(\pm)}_0(k){\mathcal B} = L^2({\bf R}^{n-1}). $
\end{cor}



\section{Modified Radon transform}

\subsection{Modified Radon transform on ${\bf H}^n$} 
The Radon transform is usually defined as an integral over some submanifolds (see e.g. \cite{Hel99}). In this section, we define the Radon transform in terms of the Fourier transform. For this purpose it is convenient to change its definition slightly. 

\begin{definition}
For $k \in {\bf R}\setminus\{0\}$ we define operators $\mathcal F^0(k)$ and $\mathcal F_0(k)$ by
\begin{equation}
\begin{split}
\mathcal F^0(k)f(x)& = \sqrt{\frac{2}{\pi}}\,k\, \sqrt{\frac{\sinh(k\pi)}{k\pi}}\\
& \ \ \ \times F_0^{\ast}\left(\Big(\frac{|\xi|}{2}\Big)^{-ik}
\int_0^{\infty}y^{\frac{n-1}{2}}K_{ik}(|\xi|y)\widehat f(\xi,y)\frac{dy}{y^n}\right), \\
\mathcal F_0(k) & = \frac{\Omega(k)}{\sqrt2}\mathcal F^0(k), \\
\Omega(k) & = \frac{-i}{\Gamma(1-ik)}
\sqrt{\frac{k\pi}{\sinh(k\pi)}}.
\end{split}
\nonumber
\end{equation}
Here $g(k) := (k\pi/\sinh(k\pi))^{1/2}$ is defined on ${\bf C}\setminus\{i\tau \,; \,\tau \in (- \infty,1]\cup[1,\infty)\}$ as a single-valued analytic function. In particular, $g(k) = g(-k)$ for $k > 0$.
\end{definition}

Note that by (\ref{eq:Chap1Sec4FormulaF0k}), $\mathcal F^0(k) = \mathcal F_0^{(+)}(k)$ for $k > 0$, and by (\ref{eq:Chap1Sec4Gammasquuare}), $|\Omega(k)|= 1$.
The following lemma follows easily from this definition and Theorem 4.2.


\begin{lemma}
(1) $\mathcal F_0$ is uniquely extended to an isometry from $L^2({\bf H}^n)$ to $\widehat{\mathcal H}$ $ := L^2({\bf R};L^2({\bf R}^{n-1});dk)$, and it diagonalizes $H_0$ :
$$
\left(\mathcal F_0H_0f\right)(k,x) = k^2\left(\mathcal F_0f\right)(k,x).
$$
(2) Let $r_+$ be the projection onto the subspace $\widehat{\mathcal H}_+ := L^2((0,\infty);L^2({\bf R}^{n-1});dk)$. Then the range of $r_+\mathcal F_0$ is $\widehat{\mathcal H}_+$. \\
\noindent
(3) $g \in \widehat{\mathcal H}$ belongs to the range of $\mathcal F_0$ if and only if
$$
\widehat g(-k,\xi) = \frac{\Gamma(1 -ik)}{\Gamma(1 + ik)}\left(\frac{|\xi|}{2}\right)^{2ik}\widehat g(k,\xi), \quad \forall k > 0.
$$
\end{lemma}

We then define the modified Radon transform associated with $H_0$ by


\begin{definition}
For $s \in {\bf R}$, we define
$$
\big(\mathcal R_0f\big)(s,x) = \frac{1}{\sqrt{2\pi}}\int_{-\infty}^{\infty}
e^{iks}\big(\mathcal F_0f\big)(k,x)dk.
$$
\end{definition}

Recall that $\mathcal F_0$ is written explicitly as 
\begin{equation} \label{5.a}
\mathcal F_0(k)f(x) = \frac{-ik}{\sqrt{\pi}\, \Gamma(1 - ik)}
F_0^{\ast}\Big(\big(\frac{|\xi|}{2}\big)^{-ik}\int_0^{\infty}
y^{\frac{n-1}{2}}K_{ik}(|\xi|y)\widehat f(\xi,y)\frac{dy}{y^n}\Big).
\end{equation}
Lemma 5.2 implies the following theorem.


\begin{theorem}\  $\mathcal R_0$ is an isometry from $L^2({\bf H}^n)$ to $\widehat{\mathcal H}$. Moreover we have
$$
\mathcal R_0 H_0 = - \partial_s^2\mathcal R_0.
$$
\end{theorem}


\subsection{Asymptotic profiles of solutions to the wave equation} 
The Radon transform thus defined describes the behiavior of solutions to the wave equation at infinity. Recall that the solution to  the wave equation
\begin{equation}
\left\{
\begin{split}
& \partial_t^2u + H_0u = 0, \\
& u\big|_{t=0}= f, \ \partial_tu\big|_{t=0} = g
\end{split}
\right.
\nonumber
\end{equation}
is written as 
$$
u(t) = \cos(t\sqrt{H_0})f + \sin(t\sqrt{H_0})\sqrt{H_0}^{-1}g.
$$


\begin{theorem}
For any $f \in L^2({\bf H}^n)$, we have as $t \to \pm\infty$
\begin{equation}
\left\|\cos(t\sqrt{H_0})f - \frac{y^{(n-1)/2}}{\sqrt2}(\mathcal R_0f)(-\log y \mp t,x)\right\|_{L^2({\bf H}^n)} \to 0,
\nonumber
\end{equation}
\begin{equation}
\left\|\sin(t\sqrt{H_0})f \mp \frac{iy^{(n-1)/2}}{\sqrt2}(\mathcal R_0\,{\rm sgn}(-i\partial_s)f)(-\log y \mp t,x)\right\|_{L^2({\bf H}^n)} \to 0,
\nonumber
\end{equation}
where 
\begin{equation}
{\rm sgn}\,(- i\partial_s)\phi(s) = \frac{1}{2\pi}\iint_{{\bf R}^1\times{\bf R}^1}
e^{ik(s- s')}{\rm sgn}\,(k)\phi(s')ds'dk,
\nonumber
\end{equation}
and where  ${\rm sgn}\,(k) = 1 \ (k > 0), \ {\rm sgn}\,(k) = - 1 \ (k < 0)$.
\end{theorem}

Proof. We  prove this theorem only for the case $t \to \infty$. Since the map : $f(k,x) \to y^{(n-1)/2}f(\log y,x
)$ is unitary from $\widehat{\mathcal H}$ onto $L^2({\bf H}^n)$, it follows from Theorem 5.4 that we have only to prove the case when 
$\phi(k,\xi) := (F_0\mathcal F_0^{(+)}f)(k,\xi) \in C_0^{\infty}({\bf R}_+ \times{\bf R}^{n-1})$. 
Let supp$\,\phi(k,\xi) \subset \{\delta_0 < k < \delta_0^{-1}\}\times\{R^{-1} < |\xi| < R\}$ for some $\delta_0, R > 0$. We put
\begin{equation} 
\begin{split}
\label{1.5.5}
&u(t,\xi,y) 
= F_0e^{-it\sqrt{H_0}}f \\
&= F_0\big(\mathcal F_0^{(+)}\big)^{\ast}e^{-itk}\mathcal F_0^{(+)}f \\
&= \int_0^{\infty}\frac{(2k\sinh(k\pi))^{1/2}}{\pi}\left(\frac{|\xi|}{2}\right)^{ik}y^{\frac{n-1}{2}}K_{ik}(|\xi|y)e^{-itk}\phi(k,\xi)dk.
\end{split}
\end{equation}
By the well-known integral representation
\begin{equation}
K_{\nu}(z) = \frac{1}{2}\int_{-\infty}^{\infty}e^{-z\cosh(s)}e^{\nu s}ds,
\nonumber
\end{equation}
(see e.g. \cite{Wa62}, Chap. 6, formula (7) or \cite{Le72}, formula (5.10.23)),
one can show that if $z > \delta_0$ for some $\delta_0 > 0$,
$$
|\partial_k^mK_{ik}(z)| \leq C_me^{-z/2}, \quad \forall m \geq 0,
$$
where the constant $C_m$ is independent of $k$.
Therefore, for any $\delta > 0$, by using $(-it)^{-1} \partial_k e^{-itk}=e{-itk}$
and integrating by parts, we see that, for any $N>0$,
\begin{equation}
\int_{\delta}^{\infty}\|u(t,\cdot,y)\|^2_{L^2({\bf R}^{n-1})}\frac{dy}{y^n} \leq  \frac{C_N}{(1+|t|)^N}.
\label{C1S5deltainfty}
\end{equation}
In the region $0 < y < \delta$, $K_{ik}(|\xi|y)$ is expanded as 
\begin{eqnarray*}
K_{ik}(|\xi|y) &=& \frac{\pi}{2i\sinh(k\pi)}\left(\frac{1}{\Gamma(1-ik)}
\Big(\frac{|\xi|y}{2}\Big)^{-ik} - \frac{1}{\Gamma(1+ik)}
\Big(\frac{|\xi|y}{2}\Big)^{ik}\right) \\
& & \hskip 5mm +\, r_1(k,|\xi|y),
\end{eqnarray*}
where $|r_1(k,|\xi|y)| \leq C|\xi|y$ uniformly for $\delta_0 < k < \delta_0^{-1},\, R^{-1}<|\xi| <R$. We put
$$
u_1(t,\xi,y) = \int_0^{\infty}\frac{(2k\sinh(k\pi))^{1/2}}{\pi}
\left(\frac{|\xi|}{2}\right)^{ik}y^{\frac{n-1}{2}}r_1(k,|\xi|y)e^{-itk}\phi(k,\xi)dk.
$$
Then 
$$
|u_1(t,\xi,y)| \leq C(\xi)y^{\frac{n+1}{2}}\int_{\delta_0}^{1/\delta_0}|\phi(k,\xi)|dk,
$$
hence
\begin{equation}
\int_0^{\delta}\|u_1(t,\cdot,y)\|_{L^2({\bf R}^{n-1})}^2\frac{dy}{y^n} \leq C_{\phi}\delta^2,
\label{C1S50deltaint}
\end{equation}
where $C_{\phi}$ is independent of $t \in {\bf R}$. We  put
\begin{equation}
\begin{split}
u_0(t,\xi,y) & = \frac{1}{i}\int_0^{\infty}\sqrt{\frac{k}{2\sinh(k\pi)}}
\left(\frac{1}{\Gamma(1 - ik)}\Big(\frac{|\xi|y}{2}\Big)^{-ik} - 
\frac{1}{\Gamma(1 + ik)}\Big(\frac{|\xi|y}{2}\Big)^{ik}\right) \\
& \hskip 20 mm \times\Big(\frac{|\xi|}{2}\Big)^{ik}y^{\frac{n-1}{2}}e^{-itk}\phi(k,\xi)dk.
\end{split}
\nonumber
\end{equation}
Then,
\begin{equation}
u_0(t,\xi,y) = u_0^{(+)}(t,\xi,y) + u_0^{(-)}(t,\xi,y).
\label{eq:Chap1Sect5u0txiy}
\end{equation}
Here
\begin{equation}
\begin{split}
& u_0^{(+)}(t,\xi,y) = \frac{1}{i}\int_0^{\infty}\sqrt{\frac{k}{2\sinh(k\pi)}}
\frac{1}{\Gamma(1 - ik)}y^{\frac{n-1}{2}}e^{-ik(t + \log y)}\phi(k,\xi)dk \\
&= \frac{y^{(n-1)/2}}{\sqrt{\pi}}\int_0^{\infty}e^{ik(-\log y - t)}
\left(F_0\mathcal F_0(k)f\right)(\xi)dk,
\end{split}
\nonumber
\end{equation}
\begin{equation} 
\begin{split}
 u_0^{(-)}(t,\xi,y) & = \frac{-1}{i}\int_0^{\infty}\sqrt{\frac{k}{2\sinh(k\pi)}}
\frac{1}{\Gamma(1 + ik)}
\Big(\frac{|\xi|}{2}\Big)^{2ik}y^{\frac{n-1}{2}}e^{-ik(t - \log y)}\phi(k,\xi)dk\\
& = \frac{y^{(n-1)/2}}{\sqrt{\pi}}\int_{-\infty}^{0}e^{ik(-\log y + t)}
\left(F_0\mathcal F_0(k)f\right)(\xi)dk.
\end{split}
\nonumber
\end{equation}
In the last equation we have used that, in view of (\ref{5.a}), (\ref{C1S3KikK-ik}),
 $(|\xi|/2)^{2ik}F_0\mathcal F^0(k)f = - F_0\mathcal F^0(-k)f.$
Rewriting $u_0^{(\pm)}(t,\xi,y)$ as
$$
u_0^{(\pm)}(t,\xi,y) = g_{\pm}(- \log y \mp t,\xi)y^{(n-1)/2}
$$
with $g_{\pm} \in L^2({\bf R}\times{\bf R}^{n-1})$, we have
$$
\int_0^{\delta}\|u_0^{(+)}(t,\cdot,y)\|^2_{L^2({\bf R}^{n-1})} 
\frac{dy}{y^n} = \int_{-\log\delta - t}^{\infty}\|g_+(\rho,\cdot)\|_{L^2({\bf R}^{n-1})}^2d\rho,
$$
which tends to 0 as $t \to - \infty$. Similarly
$$
\int_0^{\delta}\|u_0^{(-)}(t,\cdot,y)\|^2_{L^2({\bf R}^{n-1})} 
\frac{dy}{y^n} = \int^{\infty}_{- \log\delta + t}\|g_-(\rho,\cdot)\|_{L^2({\bf R}^{n-1})}^2d\rho,
$$
which tends to 0 as $t \to \infty$. 
In view of (\ref{C1S5deltainfty}), (\ref{C1S50deltaint}), we have thus proven that
\begin{equation}
\|u(t,\cdot) - u_0^{(\pm)}(t,\cdot)\|_{L^2({\bf H}^n)} \to 0 
\quad {\rm as} \quad t \to \pm \infty.
\nonumber
\end{equation}
In other words
\begin{equation}
\|F_0e^{-it\sqrt{H_0}}f - u_0^{(+)}(t)\|_{L^2({\bf H}^n)} \to 0 \quad (t \to \infty),
\nonumber
\end{equation} 
\begin{equation}
\|F_0e^{-it\sqrt{H_0}}f - u_0^{(-)}(t)\|_{L^2({\bf H}^n)} \to 0 \quad (t \to - \infty),
\nonumber
\end{equation}
\begin{equation}
\|F_0e^{it\sqrt{H_0}}f - u_0^{(-)}(-t)\|_{L^2({\bf H}^n)} \to 0 \quad (t \to \infty).
\nonumber
\end{equation}
The theorem follows from these formulas together with Definition 5.3 and (\ref{eq:Chap1Sect5u0txiy}).
\qed

\bigskip
By the change of variable $s = - \log y - t$, we get the following corollary.


\begin{cor}
For any $f \in L^2({\bf H}^n)$, we have as $t \to \infty$
$$
\sqrt{2}e^{(n-1)(s+t)/2}\left(\cos(t\sqrt{H_0})f\right)(x,e^{-s-t}) \to  
\big(\mathcal R_0f\big)(s,x) \quad {\rm in } \quad 
L^2({\bf R}^n).
$$
\end{cor}


\section{Radon transform and the wave equation}


\subsection{Radon transform and horosphere} 
As is seen in Theorem 5.5, the modified Radon transform is closely connected with the wave equation. We shall also study its geometrical feature in this section.
The fundamental solution for the wave equation on ${\bf H}^n$ is written explicitly in terms of spherical mean. For $n =3$, it has the following form (see e.g. \cite{Hel84} or \cite{ChVe96}):
\begin{equation}
\cos(t\sqrt{H_0})f(z) = \frac{\partial}{\partial t}\left(\frac{1}{4\pi\sinh(t)}\int_{S(z,t)}f(z')dS\right),
\label{eq:Chap1Sect6FundamentalSol}
\end{equation}
where $S(z;t) = \{z' ; d_h(z',z) = t\}$, and $d_h(z',z)$ is the hyperbolic 
distance. It follows from (\ref{eq:Chap1Sec1Hyperbolicdistance}) that
\begin{equation}
S(z,t) = \left\{(x',y') ; |x' -x|^2 + |y' - \cosh(t)y|^2 = 
\sinh^2(t)y^2\right\}.
\nonumber
\end{equation}
Therefore, $dS = \sinh^2(t)y^2d\omega$,  $d\omega$ being the Euclidean surface element on $S^2$, and
\begin{equation}
\cos(t\sqrt{H_0})f(z) = \frac{\partial}{\partial t}\left(\frac{\sinh(t)y^2}{4\pi}
\int_{S^2}f((x,\cosh(t)y) + \sinh(t)y\omega)d\omega\right).
\nonumber
\end{equation}
Let $t \to \infty$ and $y \to 0$ keeping $t + \log y = -s$. Then
$$
(x,\cosh(t)y) + \sinh(t)y\omega \to \big(x,\frac{e^{-s}}{2}\big) + \frac{e^{-s}}{2}\omega,
$$
Therefore, the sphere $S(z,t)$ converges to the sphere
\begin{equation}
\Sigma(s,x) = \big\{(x',y') ; \big|x' - x\big|^2 + \big|y' - \frac{e^{-s}}{2}\big|^2 = \frac{e^{-2s}}{4}\big\}.
\nonumber
\end{equation}
This is the horosphere tangent to $\{y' = 0\}$. We then have
$$
\cos(t\sqrt{H_0})f(z) \sim \frac{-y}{8\pi}\frac{\partial}{\partial s}\left(e^{-s}\int_{\Sigma(s,x)}fd\omega\right),
$$
which, compared with Theorem 5.5 with $n =3$, implies that
\begin{equation}
\mathcal R_0f(s,x) = \frac{-\sqrt2}{8\pi}\frac{\partial}{\partial s}\left(e^{-s}\int_{\Sigma(s,x)}fd\omega\right).
\nonumber
\end{equation}
From this formula, one can easily see that, if $f$ is supported in the region $y > \delta > 0$, 
then $\mathcal R_0f(s,x) = 0$ for $e^{-s} < \delta$. The converse is also true. Namely, if $\mathcal R_0f(s,x) = 0$ for $e^{-s} < \delta$, $f(x,y)$ vanishes for $y < \delta$. This is the {\it support theorem} for the Radon transform. See \cite{LaPh79} and \cite{SaBa05}. 


\subsection{1-dimensional wave equation} 
In the Euclidean space, there are 3 ways of constructing fundamental solutions to the wave equation : (1) the method of spherical means, (2) the method of plane waves and (3) the method of Fourier transforms. In the hyperbolic space, the first method is usually adopted. For example, in the work of Helgason \cite{Hel84}, a generalization of Asgeirsson's mean value theorem on two-point homogeneous space is used to derive the formula (\ref{eq:Chap1Sect6FundamentalSol}). In the following we shall apply the Fourier analysis to the fundamental solution. Let us start with the 1-dimensional case. The basic formula is


\begin{lemma}
\begin{equation}
K_{\nu}(x)K_{\nu}(y) = \frac{\pi}{2\sin(\nu\pi)}\int_{\log(y/x)}^{\infty}J_0(\sqrt{2xy\cosh t  - x^2 - y^2})\sinh(\nu t)dt
\nonumber
\end{equation}
($x > 0, y > 0, |{\rm Re}\,\nu| <  1/4$).
\end{lemma}
Proof. See \cite{DiFe33}, p. 302 and \cite{Le72} p. 140.\qed

\bigskip
 For $x > 0$ and $k \in {\bf R}$, we have by  (\ref{eq:Chap1Sec3DefinitionofInu}) and (\ref{eq:KnuandInu})
\begin{equation}
\overline{I_{ik}(x)} = I_{-ik}(x), \quad
\overline{K_{ik}(x)} = K_{ik}(x) = K_{-ik}(x),
\label{eq:Chap1Sec6Inubar}
\end{equation}
Let $\theta(t)$  be the Heaviside function: $\theta(t) = 1 \ (t > 0)$, $\theta(t) = 0 \ (t \leq 0)$.
By Lemma 6.1 and (\ref{eq:Chap1Sec6Inubar}), we have  for $x, y > 0$
\begin{equation}
\begin{split}
& \int_{-\infty}^{\infty}\sinh(\pi k)K_{ik}(x)K_{ik}(y)
\sin(tk)dk \\
& =
\frac{\pi^2}{2}\left(\theta\big(t - \log(\frac{y}{x})\big) - 
\theta\big(- t - \log(\frac{y}{x}\big)\big)\right)
J_0\big(\sqrt{2xy\cosh t - x^2 - y^2}\big).
\end{split}
\label{eq:Chap1Sec6intsinhKikKik}
\end{equation}
We put
$$
\rho(k) = \frac{2k\sinh(\pi k)}{\pi^2},
$$
and  define for $\zeta > 0$
\begin{equation}
U_{adv}(t,y,y';\zeta) = \frac{(yy')^{\frac{n-1}{2}}}{2\pi}\int_{{\bf R}^2}
\frac{K_{ik}(\zeta y)K_{ik}(\zeta y')}{k^2 - (\omega + i0)^2}\rho(k)e^{-it\omega}dk d\omega,
\nonumber
\end{equation}
\begin{equation}
U_{ret}(t,y,y';\zeta) = \frac{(yy')^{\frac{n-1}{2}}}{2\pi}\int_{{\bf R}^2}
\frac{K_{ik}(\zeta y)K_{ik}(\zeta y')}{k^2 - (\omega - i0)^2}\rho(k)e^{-it\omega}dk d\omega.
\nonumber
\end{equation}
The subscripts $adv$ and $ret$ mean advanced and retarded, respectively.


\begin{lemma} (1) For $t > 0$ and $y, y' > 0$, we have
\begin{equation}
U_{adv}(t,y,y';\zeta)  =
(yy')^{\frac{n-1}{2}}\theta\big(t - \big|\log\big(\frac{y}{y'}\big)\big|\big)
J_0\big(\zeta\sqrt{2yy'\cosh t - y^2- (y')^2}\big),
\nonumber
\end{equation}
and for $t < 0$,
\begin{equation}
U_{adv}(t,y,y';\zeta) = 0.
\nonumber
\end{equation}
(2) For $t \in {\bf R}$,
\begin{equation}
U_{ret}(t,y,y';\zeta) =  U_{adv}(-t,y,y';\zeta).
\nonumber
\end{equation}
\end{lemma}
Proof. Let us recall that if $a > 0$
\begin{equation}
\int_{-\infty}^{\infty}\frac{e^{iax}}{x - b \mp i0}dx = 
\left\{
\begin{array}{cc}
2\pi i e^{iab} & (-) \\
0 & (+),
\end{array}
\right.
\label{eq:Chap1Sec6intx-bmpi0}
\end{equation}
and if $a < 0$
\begin{equation}
\int_{-\infty}^{\infty}\frac{e^{iax}}{x - b \mp i0}dx = 
\left\{
\begin{array}{cc}
0 & (-) \\
- 2\pi ie^{iab} & (+).
\end{array}
\right.
\label{eq:Chap1Sec6intx-bmpi02}
\end{equation}
Using
\begin{equation}
\frac{1}{k^2 - (\omega + i0)^2} = \frac{1}{2k}\left(
\frac{1}{\omega + k + i0} - \frac{1}{\omega - k + i0}\right),
\nonumber
\end{equation}
we then have
\begin{equation}
\int_{-\infty}^{\infty}\frac{e^{-it\omega}}{k^2 - (\omega + i0)^2}d\omega = \left\{
\begin{split}
 &2\pi\,\frac{\sin(tk)}{k} \quad (t > 0) \\
 & 0 \hskip 14mm \quad (t < 0).
\end{split}
\right.
\nonumber
\end{equation}
Therefore by (\ref{eq:Chap1Sec6intsinhKikKik}) we have if $y, y'>0$
\begin{equation}
\int\!\!\int \frac{K_{ik}(\zeta y)K_{ik}(\zeta y')}
{k^2 - (\omega + i0)^2}e^{-it\omega}\rho(k)dk d\omega
\nonumber
\end{equation}
\begin{equation}
= \left\{
\begin{split}
&2\pi\left(\theta(t - \log(\frac{y'}{y})) - \theta(-t - \log(\frac{y'}{y}))\right)J_0\big(\zeta\sqrt{2yy'\cosh t - y^2- (y')^2}\big)  \hskip 4mm (t > 0) \\
&\hskip 40mm 0 \hskip 68mm(t< 0),
\end{split}
\right.
\nonumber
\end{equation}
which proves (1). Using (\ref{eq:Chap1Sec6Inubar}), we prove (2).
\qed


\begin{lemma} (1)
For $f \in C_0^{\infty}((0,\infty))$, we put
\begin{equation}
u_+(t,y,\zeta) = \int_0^{\infty}U_{adv}(t,y,y';\zeta)f(y')\frac{dy'}{(y')^n}.
\nonumber
\end{equation}
Then the following formulas hold:
\begin{equation}
(L_0(\zeta) + \partial_t^2)u_+(t,y,\zeta) = f(y)\delta(t),
\label{eq:Chap1Sec6L0zeta-deltsquared}
\end{equation}
\begin{equation}
u_+(t,y,\zeta) = 0 \quad {\rm for} \quad t < 0,
\label{eq:Chap1Sec6utnegativet}
\end{equation}
\begin{equation}
(\partial_tu_+)(+0,y,\zeta) = f(y).
\label{eq:Chap1Sec6deltu0}
\end{equation}
\end{lemma}
Proof. 
 Observe that, due to Lemma 6.2, for $f \in C^\infty_0((0, \infty))$,
$u_+(t, y, \zeta)$ is a well-defined smooth function of $(y, t),\, y, t>0$. 
The formula (\ref{eq:Chap1Sec6utnegativet}) is obvious. 
Consider now, for $t>0$,
\begin{equation} \label{6.8a}
\begin{split}
&  (L_0(\zeta) + \partial_t^2)u_+(t,y,\zeta) \\
= &\ 
 \frac{1}{2\pi}\int_0^{\infty}\!\int_{{\bf R}^2}(yy')^{\frac{n-1}{2}}
K_{ik}(\zeta y)K_{ik}(\zeta y')\rho(k)e^{-it\omega}\frac{f(y')}{(y')^n}
dkd\omega dy=0,
\end{split}
\end{equation}
where we have used Theorem 3.13 (2) and (3). Using (\ref{eq:Chap1Sec6intx-bmpi0}) and (\ref{eq:Chap1Sec6intx-bmpi02}), we have
\begin{eqnarray*}
& &\int_{-\infty}^{\infty}\frac{2\omega}{k^2 - (\omega + i0)^2}e^{-it\omega}d\omega \\
&=& \int_{-\infty}^{\infty}\frac{e^{-it\omega}}{k - \omega - i0}d\omega - 
\int_{-\infty}^{\infty}\frac{e^{-it\omega}}{k + \omega + i0}d\omega \\
&=& \left\{
\begin{array}{cc}
4\pi i\cos(tk) & ( t > 0), \\
0 & (t < 0).
\end{array}
\right.
\end{eqnarray*}
Therefore, we have
\begin{equation}
\partial_tu_+(t,y,\zeta) = \int_0^{\infty}\!\int_{{\bf R}^2}(yy')^{\frac{n-1}{2}}
K_{ik}(\zeta y)K_{ik}(\zeta y')\cos(tk)\rho(k)f(y')\frac{dkdy'}{(y')^n},
\nonumber
\end{equation}
which proves (\ref{eq:Chap1Sec6deltu0}). 

Formula (\ref{eq:Chap1Sec6L0zeta-deltsquared}) follows from (\ref{eq:Chap1Sec6utnegativet})
and (\ref{6.8a}).
\qed

\bigskip
We now define
\begin{equation}
U(t,y,y';\zeta) = U_{adv}(t,y,y';\zeta) - U_{ret}(t,y,y';\zeta).
\nonumber
\end{equation}
The following lemma is an easy consequence of Lemma 6.2 $(2)$ and Lemma 6.3.


\begin{lemma}
For $f \in C_0^{\infty}((0,\infty))$, we put
$$
u(t,y,\zeta) = \int_0^{\infty}U(t,y,y';\zeta)f(y')\frac{dy'}{(y')^n}.
$$
Then we have
\begin{equation}
(\partial_t^2 - L_0(\zeta))u(t,y,\zeta) = 0,
\nonumber
\end{equation}
\begin{equation}
u(0,y,\zeta) = 0,
\nonumber
\end{equation}
\begin{equation}
\partial_tu(0,y,\zeta) = f(y).
\nonumber
\end{equation}
\end{lemma}
 Note that $U_{adv}(t,y,y';\xi)$ is the Scwartz kernel of the operator
$\frac1t \sin(t L_0(\xi))$ and, therefore, defines a bounded operator in $L^2((0, \infty); dy/y^n)$.
This can be also directly observed from Theorem 3.13 (1) and (3), if we take 
$f \in L^2((0, \infty); dy/y^n)$.


\subsection{Wave equation in ${\bf H}^n$} 
We define an operator $P(t,y,y')$ by
\begin{equation}
P(t,y,y')f(x)  = (2\pi)^{-\frac{n-1}{2}}\int_{{\bf R}^{n-1}}e^{ix\cdot\xi}
p(\xi;t,y,y')
\widehat f(\xi)d\xi,
\label{eq:Chap1Sect6Defof P(t)}
\end{equation}
\begin{equation}
p(\xi;t,y,y') = J_0(|\xi|\sqrt{2yy'\cosh(t) - y^2 - (y')^2}),
\nonumber
\end{equation}
which is a Fourier multiplier acting on functions of $x \in {\bf R}^{n-1}$, depending on parameters $t, y, y'$.
Since $J_0(z)$ is an even function of $z$, $p(\xi;t,y,y')$ is smooth with respect to $\xi$ and all the other parameters $y, y'$ and $t$. By Lemma 6.4, the solution of the Cauchy problem
\begin{equation}
\left\{
\begin{split}
&\partial_t^2u + H_0u = 0, \\
& u(0) = 0, \quad \partial_tu(0) = f
\end{split}
\right.
\nonumber
\end{equation}
is written as
\begin{equation}
\begin{split}
&u(t,x,y) = \int_0^{\infty}(yy')^{\frac{n-1}{2}}
\left(\theta(t - |\log\frac{y}{y'}|) - \theta(-t - |\log\frac{y}{y'}|)\right)
\\
& \hskip 30mm \times \left(P(t,y,y')f(\cdot,y')\right)(x)\frac{dy'}{(y')^n}.
\end{split}
\nonumber
\end{equation}
Differentiating this formula with respect to $t$, we get the fundamental solution.

\begin{theorem}
Let $P$ be defined by (\ref{eq:Chap1Sect6Defof P(t)}). Then we have the following formula:
\begin{equation}
\begin{split}
\cos(t\sqrt{H_0})f(x,y) &= \int_0^{\infty}\big(yy'\big)^{\frac{n-1}{2}}\left(\delta\big(t - |\log\frac{y}{y'}|\big) + \delta\big(t + |\log\frac{y}{y'}|\big)\right) \\
& \hskip 10mm \times P(t,y,y')f(\cdot,y')(x)\frac{dy'}{(y')^n} \\
& + \int_0^{\infty}\big(yy'\big)^{\frac{n-1}{2}}\left(\theta\big(t - |\log\frac{y}{y'}|\big) - \theta\big(- t - |\log\frac{y}{y'}|\big)\right) \\
& \hskip 10mm \times \partial_tP(t,y,y')f(\cdot,y')(x)\frac{dy'}{(y')^n}.
\end{split}
\nonumber
\end{equation}
\end{theorem}

In view of Corollary 5.6, we can derive an explicit form of the modified Radon transform $\mathcal R_0f$.
Take $f \in C_0^{\infty}({\bf H}^n)$ and $s \in {\bf R}$.  
We let $t \to \infty$ and $y \to 0$ keeping $- t - \log y = s$. 
Then we have $y = e^{-s-t}$, $t - |\log(y/y')| = - s - \log y'$, and $t + |\log(y/y')|  \to \infty$. Moreover,
under these conditions, 
\begin{equation}
p(\xi;t,y,y') \to J_0(|\xi|\sqrt{e^{-s}y' - (y')^2}),
\nonumber
\end{equation}
\begin{equation}
\partial_tp(\xi;t,y,y') \to - \frac{e^{-s}|\xi|^2y'}{2}
\frac{J_1(|\xi|\sqrt{e^{-s}y' - (y')^2})}{|\xi|\sqrt{e^{-s}y' - (y')^2}},
\nonumber
\end{equation}
where we have used $J_0'(z) = - J_1(z)$.
Note that the right-hand side is again a smooth function of $s, \xi$ and $y'$, and when $y' = e^{-s}$, this $p(\xi,t,y,y') = 1$. Therefore the modified Radon transform has the following expression.


\begin{theorem}
For $f \in C_0^{\infty}({\bf H}^n)$ and $s \in {\bf R}$, we have
\begin{equation}
\mathcal R_0f(s,x) = \sqrt2e^{(n-1)s/2} f(x,e^{-s}) 
- \sqrt2e^{-s}\int_0^{e^{-s}}y^{-\frac{n-1}{2}}
A(s,y)f(\cdot,y){dy},
\nonumber
\end{equation}
where $A(s,y)f(\cdot,y)$ is  defined by
\begin{equation}
A(s,y)f(\cdot,y) = (2\pi)^{-(n-1)/2}\int_{{\bf R}^{n-1}}
e^{ix\cdot\xi}A(\xi;s,y)\widehat f(\xi,y)d\xi,
\nonumber
\end{equation}
\begin{equation}
A(\xi;s,y) = \frac{|\xi|^2}{2}\frac{J_1(|\xi|\sqrt{e^{-s}y - y^2})}{|\xi|\sqrt{e^{-s}y - y^2}}.
\nonumber
\end{equation}
\end{theorem}

Passing to the Fourier transform in Theorem 6.6 and using Definition 5.3, we have
\begin{equation}
\begin{split}
& \frac{1}{\pi}\int_{-\infty}^{\infty}\int_0^{\infty}
e^{iks}\left(\frac{|\xi|}{2}\right)^{-ik}
\frac{-ik}{\Gamma(1 - ik)}y^{-\frac{n+1}{2}}K_{ik}(|\xi|y)\widehat f(\xi,y)dydk \\& = 2e^{\frac{(n-1)s}{2}}\widehat f(\xi,e^{-s}) - e^{-s}|\xi|^2
\int_0^{e^{-s}}y^{-\frac{n-1}{2}}\frac{J_1(|\xi|\sqrt{e^{-s}y - y^2})}{|\xi|\sqrt{e^{-s}y - y^2}}\widehat f(\xi,y)dy.
\end{split}
\nonumber
\end{equation}
Taking $\widehat f(\xi,y)$ to be of the form $\varphi(\xi)\psi(y)$, and then letting $|\xi|=1$, we have
\begin{equation}
\begin{split}
& \frac{1}{\pi}\int_{-\infty}^{\infty}\int_0^{\infty}e^{iks}\,2^{ik}
\frac{-ik}{\Gamma(1-ik)}y^{-\frac{n+1}{2}}K_{ik}(y)\psi(y)dydk \\
&= 2e^{\frac{(n-1)s}{2}}\psi(e^{-s}) - 
e^{-s}\int_0^{e^{-s}}y^{-\frac{n-1}{2}}\frac{J_1(\sqrt{e^{-s}y - y^2})}{\sqrt{e^{-s}y - y^2}}\psi(y)dy.
\nonumber
\end{split}
\end{equation}
Since this holds for any $C^{\infty}_0((0,\infty))$-function $\psi(y)$, we have proven the following lemma.


\begin{lemma}
For $y > 0$
\begin{equation}
\begin{split}
& \frac{1}{\pi}\int_{-\infty}^{\infty}e^{iks}\,\frac{-ik}{\Gamma(1 - ik)}2^{ik}K_{ik}(y)dk \\
&= 2e^{-s}\delta(e^{-s}-y) 
- e^{-s} y \, \theta(e^{-s} - y)\frac{J_1(\sqrt{e^{-s}y - y^2})}{\sqrt{e^{-s}y - y^2}},
\nonumber
\end{split}
\end{equation}
where $\theta$ is the Heaviside function.
\end{lemma}

Letting $s + \log 2 = t$, one can rewrite the above formula as follows
\begin{equation}
\begin{split}
& \frac{1}{2\pi}\int_{-\infty}^{\infty}e^{ikt}\,\frac{-ik}{\Gamma(1 - ik)}K_{ik}(y)dk \\
&= 2e^{-t}\delta(2e^{-t}-y) 
- e^{-t} y\, \theta(2e^{-t} - y)\,\frac{J_1(\sqrt{2e^{-t}y - y^2})}{\sqrt{2e^{-t}y - y^2}}.
\nonumber
\end{split}
\end{equation}


\chapter{Perturbation of the metric}

We shall study in this chapter spectral properties of $- \Delta_g$, where $\Delta_g$ is the Laplace-Beltrami operator associated with a Riemannian metric, which is a perturbation of the hyperbolic metric on ${\bf H}^n$.
We shall prove the limiting absorption principle, construct the generalized Fourier transform and introduce the scattering matrix. To study ${\bf H}^n$ in an invariant manner, it is better to employ the ball model and geodesic polar coordinates centered at the origin. However, we use the upper-half space model, since it is of independent interest, necessary in order to make the arguments in Chapter 1 complete by the method adopted here, and also of a preparatory character  to deal with hyperbolic ends in Chapter 3. 


\section{Preliminaries from  elliptic partial differential equations}

\subsection{Regularity theorem} 
In this section, for the notational convenience, we denote points $x \in {\bf R}^n$ by $x = (x_1,\cdots,x_n)$. We consider the differential operator
\begin{equation}
A = \sum_{|\alpha| \leq 2}a_{\alpha}(x)(-i \partial_x)^{\alpha}
\nonumber
\end{equation}
defined on ${\bf R}^n$. The coefficients $a_{\alpha}(x)$ are assumed to satisfy
\begin{equation}
a_{\alpha}(x) \in C^{\infty}({\bf R}^n), \quad \partial_x^{\beta}
a_{\alpha}(x) \in L^{\infty}({\bf R}^n), \quad \forall \beta,
\nonumber
\end{equation}
\begin{equation}
\sum_{|\alpha| = 2}a_{\alpha}(x)\xi^{\alpha} \geq C|\xi|^2, \quad 
\forall x \in {\bf R}^n, \quad \forall \xi \in {\bf R}^n,
\nonumber
\end{equation}
$C$ being a positive constant. A function
$u \in L^2_{loc}({\bf R}^n)$ is said to be a weak solution of $Au = f$ if it satisfies
\begin{equation}
\int_{{\bf R}^n}u(x)\overline{A^{\dagger}\varphi(x)}dx = \int_{{\bf R}^n}f(x)\overline{\varphi(x)}dx, \quad \forall \varphi \in C_0^{\infty}
({\bf R}^n),
\nonumber
\end{equation}
where $A^{\dagger}$ is the formal adjoint of $A$. 


\begin{theorem}
If $u \in L^2({\bf R}^n)$ is a weak solution of $Au = f$ and $f \in H^m({\bf R}^n)$ for some $m \geq 0$, then $u \in H^{m+2}({\bf R}^n)$, and
\begin{equation}
\|u\|_{H^{m+2}({\bf R}^n)} \leq C(\|u\|_{L^2({\bf R}^n)} + 
\|f\|_{H^m({\bf R}^n)}).
\nonumber
\end{equation}
\end{theorem}

 For the proof see e.g. \cite{Mi73}. By using Theorem 1.1, one can prove the following inequality. Let $\Omega$ be a bounded open set in ${\bf R}^n$ with smooth boundary, and $\Omega_{\epsilon}$ an $\epsilon$-neighborhood of $\Omega$. Then
\begin{equation}
\|u\|_{H^{m+2}(\Omega)} \leq C_{\epsilon}(\|u\|_{L^2(\Omega_{\epsilon})} + 
\|f\|_{H^m(\Omega_{\epsilon})}).
\label{eq:Chap2Sec1Ellipticestimate}
\end{equation}


\subsection{A-priori estimates in ${\bf H}^n$} 
We next consider ${\bf R}^n_+$. We put
\begin{equation}
D_i = x_n\partial_i, \quad 1 \leq i \leq n, \quad
D = (D_1,\cdots,D_n),
\nonumber
\end{equation}
and let $(\;,\,), \ \|\cdot\|$ be the following inner product and the norm:
\begin{equation}
(u,v) = \int_{{\bf R}^n_+}u(x)\overline{v(x)}\frac{dx}{(x_n)^n}, \quad 
\|u\|^2 = (u,u).
\nonumber
\end{equation}
For operators $A$ and $B$,  $[A,B]$ denotes the commutator $AB - BA$. Straightforward computations show the following lemma.


\begin{lemma}
(1) For $\ j \neq n, \ 1 \leq i \leq n$,
\begin{equation}
[D_i,D_j] = \delta_{in}D_j.
\nonumber
\end{equation}
(2) For $\ u, v \in C_0^{\infty}({\bf R}^n_+)$, 
\begin{equation}
(D_iu,v) = - (u,D_iv) + \delta_{in}(n - 1)(u,v).
\nonumber
\end{equation}
\end{lemma}

We use the following weight
\begin{equation}
\rho(x) = \log(1 + |x|^2) + \sqrt{1+ (\log x_n)^2}.
\label{C2S1rhoxy}
\end{equation}
Comparing $\rho$ with $\rho_0$ in Lemma 1.1.6, there exists a constant $C > 0$ such that
\begin{equation}
C^{-1}(1 + d_h(x)) \leq \rho(x) \leq C(1 + d_h(x)),
\label{C2S1dhandrho1}
\end{equation}
where $d_h(x)$ is the geodesic distance between $x$ and $(0,1)$
in the metric $ds^2=dx^2/x_n^n$, cf. (\ref{eq:Riemannianmetric}) of Ch.1. 
We put
\begin{equation}
\widetilde D_i = \widetilde y(x_n)\partial_{x_i}, \quad (i = 1,\cdots, n-1), \quad
\widetilde D_n = D_n,
\label{C2S1tildeDx}
\end{equation}
where $\widetilde y(x_n) \in C^{\infty}({\bf R})$, $\widetilde y(x_n) = 1$ for $x_n < 1$, $\widetilde y(x_n) = x_n$ for $x_n > 2$. Then we have for $s \in {\bf R}$ and $|\alpha| \geq 1$
\begin{equation}
|\widetilde D^{\alpha}\rho(x)^s| +  |D^{\alpha}\rho(x)^s|
\leq C_s\rho(x)^{s-1}.
\label{CS1rhoestimate2}
\end{equation}

\medskip
We consider the differential operator
$A = A_0 + A_1$ with
\begin{equation}
A_0 = - D_n^2 + (n - 1)D_n - \sum_{i=1}^{n-1}D_i^2,
\nonumber
\end{equation}
\begin{equation}
A_1 = \sum_{i,j=1}^na_{ij}(x)D_iD_j + \sum_{i=1}^nb_i(x)D_i + c(x).
\nonumber
\end{equation}
We rewrite $A$ as
\begin{equation}
A = P_2(x,D) + P_1(x,D), \quad 
D = (D_1,\cdots,D_n),
\nonumber
\end{equation}
where
\begin{equation}
P_2(x,\xi) = |\xi|^2 + \sum_{i,j=1}^na_{ij}(x)\xi_i\overline{\xi_j},
\nonumber
\end{equation}
\begin{equation}
P_1(x,\xi) = (n - 1)\xi_n + \sum_{i=1}^nb_i(x)\xi_i + c(x).
\nonumber
\end{equation}
We assume that the coefficients $a_{ij}(x), \ b_i(x), \ c(x)$ are in $C^{\infty}({\bf R}^n_+ ; {\bf R})$ and satisfy
\begin{equation}
|\widetilde D^{\alpha} a(x)| \leq C_{\alpha}\,\rho(x)^{-\epsilon}, \quad 
\forall \alpha,
\label{C2S1axdecaty}
\end{equation}
for some $\epsilon > 0$, where $a(x)$ represents any of 
$a_{ij}(x), \ b_i(x), \ c(x)$.
Moreover, $a_{ij}$ is real and symmetric : $a_{ij} = a_{ji}$, and $P_2(x,\xi)$ is {\it uniformly elliptic}, namely, there exists a constant $C_0 > 0$ such that
\begin{equation}
P_2(x,\xi) \geq C_0|\xi|^2, \quad \forall \xi \in {\bf C}^n, \quad 
\forall x \in {\bf R}^n_+.
\label{eq:Chap2Sect1Assumption3}
\end{equation}

Let $\mathcal B$ and $\mathcal B^{\ast}$ be defined as in Chap. 1, \S 2, 
with ${\bf h} = L^2({\bf R}^{n-1})$.
For $s \in {\bf R}$, we introduce the function space $\mathcal X^s$ as follows
\begin{equation}
\mathcal X^s \ni u \Longleftrightarrow \rho(x)^su(x) \in L^2({\bf H}^n) = L^2\Big({\bf R}^n_+;\frac{dx}{x_n^n}\Big),
\label{C2S1Xs}
\end{equation}
equipped with the norm
\begin{equation}
\|u\|_{\mathcal X^s} = \|\rho^su\|_{L^2({\bf H}^n)}.
\label{C2S1normofXs}
\end{equation}


\begin{theorem} (1) If $u \in {\mathcal B}^{\ast}$ satisfies 
$(A - z)u = f \in {\mathcal B}^{\ast}$ with $z \in {\bf C}$, then
\begin{equation}
\|D_iu\|_{{\mathcal B}^{\ast}} \leq C(1 + |z|)^{1/2}(\|u\|_{{\mathcal B}^{\ast}}  + \|f\|_{{\mathcal B}^{\ast}} ), \quad 1 \leq i \leq n.
\nonumber
\end{equation}
(2) Furthermore, if
\begin{equation}
\lim_{R\to\infty}\frac{1}{\log R}\int_{1/R}^R\left[
\|u(\cdot,x_n)\|^2_{L^2({\bf R}^{n-1})} + 
\|f(\cdot,x_n)\|^2_{L^2({\bf R}^{n-1})}\right]\frac{dx_n}{(x_n)^n} = 0
\nonumber
\end{equation}
holds, then, for $1 \leq i \leq n$, we have
\begin{equation}
\lim_{R\to\infty}\frac{1}{\log R}\int_{1/R}^R
\|D_iu(\cdot,x_n)\|^2_{L^2({\bf R}^{n-1})}\frac{dx_n}{(x_n)^n} = 0.
\nonumber
\end{equation}
(3) Assertion (2) also holds with $\lim$ replaced by $\liminf$.\\
\noindent
(4) If $u, f \in L^2({\bf H}^n)$, then
\begin{equation}
\|D_iu\| \leq C(1 + |z|)^{1/2}(\|u\| + \|f\|), \quad 1 \leq i \leq n,
\label{eq:Chap2Sec1Diu}
\end{equation}
\begin{equation}
\|D_iD_ju\| \leq C(1 + |z|)(\|u\| + \|f\|), \quad 1 \leq i, j  \leq n.
\label{eq:Chap2Sec1DiDju}
\end{equation}
(5) If $u, f \in \mathcal B^{\ast}$,
\begin{equation}
\|D_iD_ju\|_{\mathcal X^{-s}} \leq C_s(1 + |z|)(\|u\|_{\mathcal B^{\ast}} + \|f\|_{\mathcal B^{\ast}}), \quad 1 \leq i, j  \leq n,
\label{eq:Chap2Sec1DiDjuBast}
\end{equation}
for any $s > 1/2$. \\
\noindent
(6) If $u, f \in {\mathcal X}^s$ for some $s \in {\bf R}$, then 
\begin{equation}
\|D_iu\|_{\mathcal X^s} \leq C(1 + |z|)^{1/2}(\|u\|_{\mathcal X^s} + \|f\|_{\mathcal X^s}), \quad 1 \leq i \leq n,
\label{eq:Chap2Sec1Dius}
\end{equation}
\begin{equation}
\|D_iD_ju\|_{\mathcal X^s} \leq C(1 + |z|)(\|u\|_{\mathcal X^s} + \|f\|_{\mathcal X^s}), \quad 1 \leq i, j  \leq n.
\label{eq:Chap2Sec1DiDjus}
\end{equation}

In the above estimates in (1), (4), (5) and (6), the constants $C$ and $C_s$ are independnet of $z \in {\bf C}$.
\end{theorem}

We note that assertion (4) is a particular case of assertion (6) with $s=0$,
while assertion (5) follows from (6), if we take into the account that 
$\mathcal B^{\ast} \subset \mathcal X^{-s},\, s>1/2$.  

\medskip

Proof. We take $\chi(t) \in C_0^{\infty}({\bf R})$ such that $\chi(t) = 1 \ (|t| < 1), \ 
\chi(t) = 0 \ (|t| > 2)$, and put
\begin{equation}
\chi_{R,r}(x) = \chi\big(\frac{\log x_n}{\log R}\big)
\chi\big(\frac{|x'|}{r}\big), 
\quad
\chi_{R}(x_n) = \chi\big(\frac{\log x_n}{\log R}\big),
\nonumber
\end{equation}
where $x' = (x_1,\cdots,x_{n-1})$. Since with $g_{ij} = \delta_{ij} + a_{ij}$,
\begin{eqnarray*}
(g_{ij}D_iD_ju,\chi_{R,r}^2u) &=& 
- (g_{ij}D_iu,\chi_{R,r}^2D_iu)  \\
& &   - (D_ju,\left(D_i(g_{ij}\chi_{R,r}^2)\right)u) + 
\delta_{in}(n-1)(D_ju,g_{ij}\chi_{R,r}^2u).
\end{eqnarray*}
Thus, we have
\begin{eqnarray*}
- \;\sum_{i,j=1}^n(g_{ij}D_iD_ju,\chi_{R,r}^2u) &=& 
\sum_{i,j=1}^n (g_{ij}\chi_{R,r}D_ju,\chi_{R,r}D_iu)  \\
& &   + \sum_{i,j=1}^n(D_ju,\left(D_i(g_{ij}\chi_{R,r}^2)\right)u) \\
& &- \sum_{j=1}^n
\delta_{in}(n-1)(D_ju,g_{nj}\chi_{R,r}^2u).
\end{eqnarray*}
We split the 2nd term of the right-hand side into
$$
- \sum(\chi_{R,r}D_ju,(D_ig_{ij})\chi_{R,r}u) - 
2\sum(\chi_{R,r}D_ju,g_{ij}(D_i\chi_{R,r})u)
$$
and use the uniform ellipticity (\ref{eq:Chap2Sect1Assumption3}) to see that
\begin{equation}
\begin{split}
C_0\|\chi_{R,r}Du\|^2 & \leq  {\rm Re}\,(Au,\chi_{R,r}^2u) + \epsilon\|\chi_{R,r}Du\|^2 \\
&\ \ \ \ \ + C_{\epsilon}(\|\psi_Ru\|^2 + 
\|(D\chi_{R,r})u\|^2).
\end{split}
\nonumber
\end{equation}
Here  $\psi_R$ is defined by
$$
\psi_{R}(x_n) = \psi\big(\frac{\log x_n}{\log R}\big), 
$$
where $\psi \in C_0^{\infty}({\bf R})$, $\psi = 1$ on the support of $\chi$.
For small $\epsilon > 0$, the term $\epsilon\|\chi_{R,r}Du\|^2$ is absorbed by the left-hand side. Therefore, by using the equation $(A - z)u = f$, we have
$$
\|\chi_{R,r}Du\|^2  \leq  C (1+|z|) (\|\psi_Ru\|^2 + 
\|(D\chi_{R,r})u\|^2 + \|\psi_Rf\|^2).
$$
We fix $R$ and let $r \to \infty$ to see that $\chi_{R,r}$ can be replaced by $\chi_R$. Moreover 
$$
|(D\chi_R)(x_n)| \leq \frac{C}{\log R}\psi_R(x_n) \leq C\psi_R(x_n)
$$ 
for $R > e$.
Therefore, we have
\begin{equation}
\|\chi_{R}Du\|^2  \leq C(1 + |z|)(\|\psi_Ru\|^2 + \|\psi_Rf\|^2).
\label{eq:Chap2Sec1chiRD}
\end{equation}
Dividing this inequality by $\log R$ and taking the supremum with respect to $R$, we obtain the assertion (1). Letting $R \to \infty$, we obtain (2) and (3). 

Letting $R \to \infty$ in (\ref{eq:Chap2Sec1chiRD}), we prove (\ref{eq:Chap2Sec1Diu}).
To prove (\ref{eq:Chap2Sec1DiDju}), we first observe that the previous considerations do not require (\ref{C2S1axdecaty}) in full
generality, just that $a \in L^\infty(\bf R_+^n)$. This makes it possible to consider
only  the case when $u$ is compactly supported. 
 In fact, in the general case putting
$\chi_{R,r}u = v$ we have
\begin{equation}
(A-z)v = \chi_{R,r}f + [A,\chi_{R,r}]u.
\nonumber
\end{equation}
Since $[A,\chi_{R,r}] = \sum_ic_i(x)D_i + d(x)$ and $c_i(x), d(x)$  and $c_i, d \in L^\infty$ independently on $R, r >e$, we can apply (\ref{eq:Chap2Sec1Diu}) and (\ref{eq:Chap2Sec1DiDju}) 
to see that the right-hand side is in $L^2(\bf R_+^n)$ uniformly with respect to $R, r$. 

Now assuming that $u$ is compactly supported, we split $u$ as $u = u_1 + u_2 + u_3$, where $u_i = \chi_i(\frac{\log x_n}{\log R})u$ so that 
${\rm supp}\,u_1 \subset \{x_n < 2/R\}$, ${\rm supp}\,u_2 \subset \{1/R < x_n < 2R\}$, ${\rm supp}\,u_3 \subset \{x_n > R\}$.
 Using
$$
\|D_iD_ju\|^2 = (D_j^2u,D_i^2u) + (D_ju,[D_j,D_i]D_ju),
$$
we have
\begin{equation}
\sum_{i,j}\|D_iD_ju\|^2 \leq C(\|\sum_iD_i^2u\|^2 + \sum_i\|D_iu\|^2).
\nonumber
\end{equation}
We have
\begin{equation}
A_0 u_i = - A_1u_i +z u_i+ f_i, \quad i = 1, 3,
\label{C2S1A0ui}
\end{equation}
where 
$$
\|f_i\| \leq C(\|f\| + \|D_nu\| + \|u\|) \leq C (1+|z|)^{1/2} (\|f\| + \|u\|),
$$
with the last inequality following from (\ref{eq:Chap2Sec1Diu}).
Since
$
\|A_0u_i\|^2 = \sum_{j,k}(D_j^2u_i,D_k^2u_i),
$
taking the $L^2$-norm of the both sides of (\ref{C2S1A0ui}),
and using condition (\ref{C2S1axdecaty}), we have, for $i=1, 3$,
$$
\sum_{j,k}\|D_jD_ku_i\| \leq \epsilon\sum_{j,k}\|D_jD_ku_i\| + 
C_{\epsilon} (1+|z|) \big(\sum_{j}\|D_ju_i\| + \|u\| + \|f\|\big), 
$$
where $\epsilon = \epsilon(R) \to 0$ as $R \to \infty$. 
Therefore (\ref{eq:Chap2Sec1DiDju}) holds for $i = 1, 3$ with sufficiently large $R$. For $i = 2$, we have only to note that $u_2$ satisfies the following 2nd order elliptic equation with bounded coefficients:
\begin{equation}
\sum_{i,j}\widetilde{a}_{ij}(x)\partial_i\partial_j u_2 + \sum_i
\widetilde{a_i}(x)\partial_iu_2 + \widetilde{c}(x)u_2 = f_2
\nonumber
\end{equation}
and use Theorem 1.1.

To prove (5), we put $v = \rho(x)^{-s}u$ and $g = (A-z)v$. Then Lemma 1.2.7,  estimate (\ref{CS1rhoestimate2}) and  assertion (1) imply that $v, g \in L^2({\bf H}^n)$. By  assertion (4), 
we then have $D_iv, D_iD_j v \in L^2({\bf H}^n)$, which, in turn, implies that
$D_iD_ju \in \mathcal X^{-s}$ and the inequality (\ref{eq:Chap2Sec1DiDjuBast}).

The proof of (\ref{eq:Chap2Sec1Dius}) is similar to the proof of (\ref{eq:Chap2Sec1Diu}) if w use $\rho(x)^s\chi_{R,r}(x)$ instead of $\chi_{R,r}(x)$.

 To prove (\ref{eq:Chap2Sec1DiDjus}), we again consider $v= \rho(x)^{-s} u$,
which, due to (\ref{eq:Chap2Sec1Dius}) satisfies $(A-z) v=g \in L^2({\bf H}^n)$.
Using (\ref{eq:Chap2Sec1Diu}) together with  (\ref{eq:Chap2Sec1Dius}) and (\ref{CS1rhoestimate2}),
we arrive at (\ref{eq:Chap2Sec1DiDjus}).
\qed

\subsection{Essential self-adjointness}

 On the upper space $\bf R_+^n$, we introduce the Riemannian metric
\begin{equation} \label{E2.1}
ds^2= \frac{1}{x_n^n} \sum_{i, j=1}^n g_{ij}(x) dx_i dx_j,
\end{equation}
where $g_{ij}=\delta_{ij}+ a_{ij}$. Assume that $A$ is symmetric on $C_0^\infty(\bf R_+^n)$.
}

\begin{theorem}
$A\big|_{C_0^{\infty}({\bf R}^n_+)}$ is essentially self-adjoint.
\end{theorem}
Proof. We show that for $u \in L^2({\bf H}^n)$ 
\begin{equation}
(u,(A - i)\varphi) = 0, \quad \forall \varphi \in C_0^{\infty}({\bf H}^n)
\Longrightarrow u = 0
\nonumber
\end{equation}
and the same assertion holds with $i$ replaced by $- i$. Applying (\ref{eq:Chap2Sec1Ellipticestimate}), we see that $u \in H^2_{loc}({\bf R}^n_+)$, and $(A + i)u = 0$ holds, moreover, by Theorem 1.3 (4),
\begin{equation}
D_iu, \ D_iD_ju \in L^2({\bf H}^n).
\nonumber
\end{equation}
Letting
\begin{equation}
\Omega_{r,R} = \{|x'| < r, \ 1/R < x_n < R\}, \quad
\Omega_{R} = \{1/R < x_n < R\},
\nonumber
\end{equation}
we then have
\begin{equation}
\int_{\Omega_{r,R}}Au\overline{u}d\mu = - i\int_{\Omega_{r,R}}|u|^2d\mu, \quad
d\mu = dx/(x_n)^n.
\nonumber
\end{equation}
Integrating by parts and taking the imaginary part,
\begin{equation}
\int_{\Omega_{r,R}}|u|^2d\mu \leq C\sum_i\int_{\partial\Omega_{r,R}}
|u||D_iu|dS,
\nonumber
\end{equation}
where $dS$ is the surface measure associated with hyperbolic metric.
Noting that
\begin{equation}
\int_{1/R < x_n < R}|uD_iu|d\mu < \infty,
\nonumber
\end{equation}
there is a sequence 
 $r_n \to \infty$ such that,
 $$\sum_i \int_{\Sigma_{R, n}}|u| |D_i u| dS \to 0 \quad \hbox{as}\,\, n \to \infty,
 $$
 where $\Sigma_{R, n}=\{(x', x_n): \,|x'| =r_n,\,  R^{-1}< x_n < R\}$.
 Using these $r_n'$s, we see that
\begin{equation}
\int_{\Omega_{R}}|u|^2d\mu \leq C\sum_{i=1}^n\left(\int_{x_n = 1/R} + \int_{x_n=R}\right)
|u||D_iu|\frac{dx'}{(x_n)^{n-1}}.
\label{eq:Chap2Sect1OmegaRu}
\end{equation}
We next put
\begin{equation}
f(x_n) = \sum_{i=1}^n\int_{{\bf R}^{n-1}}|uD_iu(x',x_n)|\frac{dx'}{(x_n)^{n-1}}.
\nonumber
\end{equation}
Then, since $u, D_i u \in L^2(\bf H^n)$, we have
\begin{equation}
\int_0^{\infty}f(x_n)\frac{dx_n}{x_n} < \infty.
\nonumber
\end{equation}
Hence, $\liminf_{x_n\to\infty}f(x_n) = 0$ and $\liminf_{x_n\to0}f(x_n) = 0$. Using this fact, letting
$R_n$ tend to infinity along a suitable sequence in (\ref
{eq:Chap2Sect1OmegaRu}), we have $u = 0$. \qed

\subsection{Rellich's theorem}
It is well-known that,  for a bounded open set $\Omega \subset {\bf R}^n$, the inclusion $H^1(\Omega) \subset L^2(\Omega)$ is compact. This is often stated in the following form and is called Rellich's theorem.


\begin{theorem}
Let $\Omega$ be a bounded open set in  ${\bf R}^n$, and $m \geq 1$. Then for any bounded sequence $\{f_k\}$ in $H^m(\Omega)$, there exists a subsequence $\{f_{k'}\}$ convergent in $H^{m-1}(\Omega)$.
\end{theorem}

For the proof, see e.g. \cite{Mi73}.

 
\subsection{Unique continuation theorem} 
Let us assume that on a connected open set $\Omega \subset {\bf R}^n$, we are given a differential operator
\begin{equation}
A = \sum_{|\alpha| \leq 2}a_{\alpha}(x)\partial_x^{\alpha},
\nonumber
\end{equation}
where for $|\alpha| = 1, 2$, $a_{\alpha}(x) \in C^{\infty}$, and for $|\alpha| = 0$, $a_{\alpha}(x) \in L^{\infty}$, moreover for $|\alpha| = 2$, $a_{\alpha}(x)$ is real-valued and satisfies
\begin{equation}
 \sum_{|\alpha|=2}a_{\alpha}(x)\xi^{\alpha} \geq C|\xi|^2, \quad 
 \forall x \in {\Omega}, \quad \forall \xi \in {\bf R}^n,
 \nonumber
\end{equation}
for a constant $C > 0$. Then, if  $u$ satisfies 
$Au = 0$ on $\Omega$, and vanishes on an open subset of $\Omega$, then $u$ vanishes identically on $\Omega$. For the proof, see e.g. \cite{Mi73} for a $C^{\infty}$-coefficient case, and  \cite{Ar57} for the general case.  


\section{Basic spectral properties for Laplace-Belrami operators on ${\bf H}^n$}

\subsection{Assumption on the metric} 
In the sequel, we denote points in ${\bf H}^n = {\bf R}^n_+$ as $(x,y)$, where $x \in {\bf R}^{n-1}$, $y > 0$, and put
\begin{equation}
D_x = y\partial_x, \quad \widetilde D_x = \tilde y(y) \partial_x, \quad  \widetilde D_y=D_y = y\partial_y,
\label{eq:DxDywidetildeDx}
\end{equation}
where $\tilde y(y) \in C^{\infty}((0,\infty))$ is a positive function such that $\tilde y(y) = 1$ for $y < 1$, $\tilde y(y) = y$ for $y > 2$. Recall that we put\begin{equation}
\rho(x,y) = \log\left(1+ |x|^2 + y^2\right) + \sqrt{1 +  |\log y|^2},
\nonumber
\end{equation}
and have the following inequality
$$
C^{-1}(1 + \rho(x,y)) \leq 1 + d_h(x,y) \leq C(1 + \rho(x,y)),
$$
$$
|\widetilde D^{\alpha}\rho(x,y)^s| + |D^{\alpha}\rho(x,y)^s| \leq C_s\rho(x,y)^{s-1}, \quad |\alpha| \geq 1, \quad s \in {\bf R},
$$
where $d_h(x,y)$ is the distance between $(x,y)$ and $(0,1)$ with respect to the standard hyperbolic metric (Lemma 1.1.6).

To describe the space of metric, we introduce the following class of functions. 

\begin{definition}
{\it For $s \in {\bf R}$, let $\mathcal W^s$ be the set of real-valued $C^{\infty}$-functions $f(x,y)$ defined on ${\bf R}^{n-1}\times(0,\infty)$ such that for any (multi) index $\alpha$, $\beta$, there exists a constant $C_{\alpha\beta} > 0$ such that}
\begin{equation}
|(\widetilde D_x)^{\alpha}(D_y)^{\beta}\,f(x,y)| \leq C_{\alpha\beta}\, \rho(x,y)^{s-{\rm min}(|\alpha|+\beta,1)}.
\label{eq:eq:C2S2Wsy>1}
\end{equation}
\end{definition}

On the upper half-space ${\bf R}^n_+$, we consider the Riemannian metric
\begin{equation}
ds^2 = y^{-2}\Big((dx)^2 + (dy)^2 + A(x,y,dx,dy)\Big),
\label{eq:Chap2Sect2ds2generalform}
\end{equation}
where $A(x,y,dx,dy)$ is a symmetric covariant tensor of the form
\begin{equation}
A(x,y,dx,dy) = \sum_{i,j=1}^{n-1}a_{ij}(x,y)dx^idx^j + 2\sum_{i=1}^{n-1}a_{in}(x,y)dx^idy + 
a_{nn}(x,y)(dy)^2.
\nonumber
\end{equation}
Here each $a_{ij}(x,y) \ (1 \leq i, j \leq n)$ is assumed to satisfy the following condition:

\medskip
\noindent
{\bf (C)} \  {\it There exists a constant $\epsilon > 0$ such that}
$a_{ij} \in {\mathcal W}^{-1-\epsilon}$ \ {\it for}  $\ y >1$.

\medskip
Let us look at the Laplace-Beltrami operator associated with the above metric $ds^2$. 
 Let  $\mathcal P$ the set of differential operators $P$ defined by
\begin{equation}
\mathcal P \ni P \Longleftrightarrow
P = \sum_{\alpha,\beta}(c_{\alpha\beta}+a_{\alpha\beta}) D_x^{\alpha}D_y^{\beta}, 
\nonumber
\end{equation}
where  $c_{\alpha\beta}$ are constants, $a_{\alpha\beta} \in  \mathcal {\mathcal W}^{-1-\epsilon}$ and the above sum is finite.
Then by a direct computation using Lemma 1.2 one can show that
$\mathcal P$ is an algebra.

We rewrite (\ref{eq:Chap2Sect2ds2generalform}) into $ds^2 = g_{ij}(X)dX^idX^j$, $X = (X_1,\cdots,X_n) = (x,y)$, where $g_{ij}(X) = y^{-2}\big(\delta_{ij} + a_{ij}(x,y)\big)$ and we assume that 
$a_{ij} \xi_i \xi_j > -|\xi|^2$. Letting
$\big(g^{ij}\big) = \big(g_{ij}\big)^{-1}$, we have
\begin{equation}
g^{ij}(x,y) = y^2\big(\delta^{ij} + \widehat g^{ij}(x,y)\big),
\nonumber
\end{equation}
where $\widehat g^{ij}(x,y) \in {\mathcal W}^{-1-\epsilon}$. The associated Laplace-Beltrami operator $\Delta_g$ is then written as
\begin{equation}
- \Delta_g = D_y^2 - (n-1)D_y + D_x^2 + \sum_{i,j=1}^{n}a^{ij}(x,y)D_iD_j + 
\sum_{i=1}^nb^i(x,y)D_i,
\nonumber
\end{equation}
where $(D_1,\cdots,D_{n}) = (y\partial_x,y\partial_y)$ and $a^{ij}(x,y),  b^i(x,y) \in {\mathcal W}^{-1-\epsilon}$. Hence $\Delta_g \in \mathcal P$.

The operator $- \Delta_g$ is symmetric in $L^2({\bf R}^n_+;\sqrt{g}\;dxdy)$, where $g = \det(g_{ij})$. In order to compare it with the Laplace-Beltrami operator for the standard hyperbolic metric, it is convenient to use the unitary gauge transformation from $L^2({\bf R}_+^n; {\sqrt g}dx dy)$ onto
$L^2({\bf R}_+^n; dx dy/y^n)$:
$$ u \to (y^{2n}g)^{1/4}u,$$ so that 
\begin{equation}
- \Delta_g - \frac{(n-1)^2}{4} \to - (y^{2n}g)^{1/4}\Delta_g(y^{2n}g)^{-1/4} - \frac{(n-1)^2}{4}
\nonumber
\end{equation}
in $L^2({\bf R}^n_+;dxdy/y^n)$.


\subsection{Transformed Laplace-Beltrami operators} 
We are thus led to the differential operators
\begin{equation}
H = - (y^{2n}g)^{1/4}\Delta_g(y^{2n}g)^{-1/4} - \frac{(n-1)^2}{4} = H_0 + V,
\nonumber
\end{equation}
\begin{equation}
 H_0 = - D_y^2 + (n-1)D_y - D_x^2 - \frac{(n-1)^2}{4},
\quad
 V = \sum_{|\alpha|\leq2}a_{\alpha}(x,y)D^{\alpha}
 \nonumber
\end{equation}
in $L^2({\bf R}^n_+;dxdy/y^n)$, with the inner product  denoted by $(\cdot,\cdot)$. $H\big|_{C_0^{\infty}({\bf H}^n)}$ is symmetric, 
\begin{equation}
(Hf,g) = (f,Hg), \quad \forall f, g \in C_0^{\infty}({\bf H}^n),
\label{eq:Chap2Sect2Symmetric}
\end{equation}
 and uniformly elliptic in the sense of \S 1. By our assumption $a_{\alpha\beta}$ satisfies the condition (C).

\medskip
One should keep in mind that our operator $- H$ is unitarily equivalent to the Riemannian Laplacian $\Delta_g$ associated with the metric $ds^2$ of (\ref{eq:Chap2Sect2ds2generalform}) which is shifted by $(n-1)^2/4$. The arguments to be developed in Chapters 2 and 3 are also applicable to the more general operators with perturbation of 1st order differential operators, except for Theorem 2.10. Even in this case, however, Theorem 2.10 still holds except for a discrete set of $\lambda$'s, which can be proved by the same way as in  Theorems 3.3.5 and 3.3.6.

\medskip
 By Theorem 1.4, $H\big|_{C_0^{\infty}({\bf H}^n)}$ is essentially self-adjoint. Let
\begin{equation}
R_0(z) = (H_0 - z)^{-1}, \ R(z) = (H - z)^{-1}.
\nonumber
\end{equation}

\begin{lemma}
For $z \not\in {\bf C}\setminus{\bf R}$, $R_0(z)VR(z)$ is compact. Hence 
\begin{equation}
\sigma_d(H) \subset (- \infty,0), \quad 
\sigma_e(H) = [0,\infty).
\nonumber
\end{equation}
\end{lemma}
Proof. By Theorem 1.3 (4),  $VR(z) \in {\bf B}(L^2;L^2)$, and $R_0(z)V = (VR_0(\overline{z}))^{\ast} \in 
{\bf B}(L^2;L^2)$. We take $\chi(t) \in C_0^{\infty}({\bf R})$ satisfying $\chi(t) = 1 \ (|t| < 1), \ \chi(t) = 0 \ (|t| > 2)$, and put
\begin{equation}
 \chi_R(x,y) = \chi\left(\frac{|x|}{R}\right)\chi
 \left(\frac{\log y}{\log R}\right).
 \nonumber
\end{equation}
Then $\chi_R R(z)$, and henceforth $R_0(z)V\chi_RR(z)$ are compact and,
due to the decay assumption of the coefficients,  $\|R_0(z)V(1 - \chi_R)R(z)\| \to 0
 \ (R \to \infty)$. Hence $R_0(z)VR(z)$ is also compact. Since $\sigma(H_0) = \sigma_e(H_0) = [0,\infty)$, the lemma follows from Weyl's theorem (\cite{Is04a}, p. 26). \qed
 
\bigskip
The main purpose of this section is to prove the following theorem.


\begin{theorem} 
(1) $\ \sigma_p(H)\cap(0,\infty) = \emptyset$. \\
\noindent
(2) For any $\lambda > 0$, $\lim_{\epsilon \to 0}
R(\lambda \pm i\epsilon) =: R(\lambda \pm i0)$ exists in the weak-$\ast$ sense, namely
\begin{equation}
\exists\lim_{\epsilon \to 0}(R(\lambda \pm i\epsilon)f,g) =: (R(\lambda \pm i0)f,g), \quad 
\forall f, g \in \mathcal B.
\nonumber
\end{equation} 
\noindent
(3) For any compact interval $I\subset (0,\infty)$ there exists a constant $C > 0$ such that
\begin{equation}
\|R(\lambda \pm i0)f\|_{{\mathcal B}^{\ast}} \leq C\|f\|_{\mathcal B}, \quad \forall \lambda \in I.
\label{C2S2RlambdafEstimate}
\end{equation}
(4) For any $f, g \in {\mathcal B}$, $(0,\infty) \ni \lambda \to (R(\lambda \pm i0)f,g)$ is continuous. \\
\noindent
(5) Let $E_H(\cdot)$ be the resolution of the identity for $H$. Then $E_H((0,\infty))L^2({\bf H}^n)$ is equal to the absolutely continuous subspace for $H$.
\end{theorem}

Note that the proof of the estimate (\ref{C2S2RlambdafEstimate}) implies the following inequality
\begin{equation}
\|R(z)f\|_{{\mathcal B}^{\ast}} \leq C\|f\|_{\mathcal B}, \quad \forall {\rm Re}\,z \in I.
\label{R(z)UniformEstimate}
\end{equation}


\subsection{Resolvent estimates} 
We shall prove Theorem 2.3 by first establishing some a-priori estimates for solutions to the equation $(H -z)u = f$, and then passing to limiting procedures. Although our method seems to be tricky, the basic idea consists in the following observation.
Let us note that by virtue of Lemma 1.4.7, $u_{\pm}^0 = R_0(\lambda \pm i0)f$ behaves like
$$
 \hat u_{\pm}^0(\xi,y) \sim C_{\pm}(\xi)y^{(n - 1)/2 \mp i\sqrt{\lambda}} \quad 
 (y \to 0).
$$
Therefore, we infer
$$
 \left(y{\partial}_y - (\frac{n-1}{2} \mp i\sqrt{\lambda})\right)u_{\pm}^0 = 
 o(y^{(n-1)/2}) \quad (y \to 0).
$$
This suggests the importance of the term $\left(y{\partial}_y - (\frac{n-1}{2} \mp i\sqrt{\lambda})\right)u_{\pm}^0$ to derive the estimates for $u_{\pm}^0$. We put
\begin{equation}
 \sigma_{\pm} = \frac{n-1}{2} \mp i\sqrt{z}.
 \nonumber
\end{equation}
Here for $z = re^{i\theta}, r > 0, - \pi < \theta < \pi$, we take the branch of $\sqrt{z}$ as $\sqrt{r}e^{i\theta/2}$. 

We begin by estimating $u^0 = 
R_0(\lambda + i0)f$. 
Let $(\;,\;)_{\bf h}, \ \|\cdot\|_{\bf h}$ denote the inner product and  norm of $L^2({\bf R}^{n-1})$, respectively.


\begin{lemma} Suppose $u$ satisfies $(H_0 - z)u = f$, and let $w_{\pm} = (D_y - \sigma_{\pm})u$.  Let $\varphi(y) \in C^1((0,\infty);{\bf R})$ and $0 <a < b < \infty$. Then we have
\begin{eqnarray*}
& &  \int_a^b(D_y\varphi + 2\varphi)\|D_xu\|_{\bf h}^2\frac{dy}{y^n} + 
\left[\frac{\varphi(\|w_{\pm}\|_{\bf h}^2 - \|D_xu\|_{\bf h}^2)}{y^{n-1}}\right]_{y=a}^{y=b}\\
& =& \mp 2\,{\rm Im}\,\sqrt{z}\int_a^b\varphi\left(\|w_{\pm}\|_{\bf h}^2 + 
\|D_xu\|_{\bf h}^2\right)\frac{dy}{y^n} \\
& & + \int_a^b(D_y\varphi)\|w_{\pm}\|_{\bf h}^2\frac{dy}{y^n} 
- 2\,{\rm Re}\int_a^b\varphi(f,w_{\pm})_{\bf h}\frac{dy}{y^n}.
\end{eqnarray*}
\end{lemma}

Proof. 
We rewrite the equation $(H_0 - z)u = f$ as 
\begin{equation}
D_y(D_y - \sigma_{\pm})u = \sigma_{\mp}(D_y - \sigma_{\pm})u - 
D_x^2u - f.
 \label{eq:Dyminussigmaplusu}
\end{equation}  
Taking the inner product of (\ref{eq:Dyminussigmaplusu}) and $\varphi w_{\pm}$, we have
\begin{equation}
\begin{split}
& \int_a^b\varphi(D_yw_{\pm},w_{\pm})_{\bf h}\frac{dy}{y^n} \\
= 
&\  \sigma_{\mp}\int_a^b\varphi\|w_{\pm}\|_{\bf h}^2\frac{dy}{y^n} 
-\int_a^b\varphi(D_x^2u,w_{\pm})_{\bf h}\frac{dy}{y^n} - 
\int_a^b\varphi(f,w_{\pm})_{\bf h}\frac{dy}{y^n}.
\end{split}
\label{C2S2ProofLemma2.4}
\end{equation}
Take the real part.
By integration by parts, the left-hand side is equal to
\begin{equation}
\begin{split}
& {\rm Re}\int_a^b\varphi(D_yw_{\pm},w_{\pm})_{\bf h}\frac{dy}{y^n} \\ 
= & \ 
\left[\frac{\varphi\|w_{\pm}\|_{\bf h}^2}{2y^{n-1}}\right]_{y=a}^{y=b} 
- \frac{1}{2}\int_a^b(D_y\varphi)\|w_{\pm}\|_{\bf h}^2\frac{dy}{y^n} + 
\frac{n-1}{2}\int_a^b\varphi\|w_{\pm}\|_{\bf h}^2\frac{dy}{y^n}.
\end{split}
\label{C2S2identity1}
\end{equation}
Let us note that using
$$                                    
(-D_x^2u,D_yu)_{\bf h} = (v,D_yv)_{\bf h} - \|v\|_{\bf h}^2, 
\quad v = \sqrt{D_x^2}u = y\sqrt{-\Delta_x}u,
$$
we have
\begin{eqnarray*}
& &- {\rm Re}\int_a^b\varphi\big(D_x^2u,w_{\pm}\big)_{\bf h}\frac{dy}{y^n} 
\\
& =& \left[\frac{\varphi\|D_xu\|_{\bf h}^2}{2y^{n-1}}\right]_{y=a}^{y=b} 
- \frac{1}{2}\int_a^b(D_y\varphi)\|D_xu\|_{\bf h}^2\frac{dy}{y^n} 
+ \left(\dfrac{n-3}{2} - {\rm Re}\,\sigma_{\pm}\right)\int_a^b
\varphi\|D_xu\|_{\bf h}^2\frac{dy}{y^n}.
\end{eqnarray*}
Apply this to the 2nd term of the right-hand side of (\ref{C2S2ProofLemma2.4}).  We then have
\begin{equation}
\begin{array}{rcl}
& & {\displaystyle {\rm Re}\int_a^b\varphi(D_yw_{\pm},w_{\pm})_h\frac{dy}{y^n}} \\
& =& {\displaystyle ({\rm Re}\, \sigma_{\mp})\int_a^b\varphi\|w_{\pm}\|_h^2\frac{dy}{y^n} 
-{\rm Re}\int_a^b\varphi(y^2\Delta_hu,w_{\pm})_h\frac{dy}{y^n} 
- {\rm Re}\int_a^b\varphi(f,w_{\pm})_h\frac{dy}{y^n}} \\
&= &{\displaystyle 
\left(\frac{n-1}{2} \mp {\rm Im}\,\sqrt{z}\right)\int_a^b\varphi\|w_{\pm}\|_h^2\frac{dy}{y^n}
+ \left[\frac{\varphi\|D_xu\|_h^2}{2y^{n-1}}\right]_{y=a}^{y=b} } \\
& &{\displaystyle  - \frac{1}{2}\int_a^b(D_y\varphi)\|D_xu\|_h^2\frac{dy}{y^n}
-  (1 \pm {\rm Im}\,\sqrt{z})\int_a^b\varphi
\|D_xu\|_h^2\frac{dy}{y^n} - {\rm Re}\int_a^b\varphi(f,w_{\pm})_h\frac{dy}{y^n}.}
\end{array}
\label{C2S2identity2}
\end{equation}
Equating (\ref{C2S2identity1}) and (\ref{C2S2identity2}), we obtain the lemma. 
 \qed

\medskip
We shall derive estimates of the resolvent $R_0(z) = (H_0 - z)^{-1}$, when $z \in {\bf C}\setminus{\bf R}$ approaches the real axis.


\begin{lemma} 
Let $u = R_0(z)f$. Let $w_{\pm} = (D_y - \sigma_{\pm})u$, and put for  $C^1 \ni \varphi \geq 0$ and constants $0 < a < b$, 
\begin{equation}
L_{\pm} =  \int_a^b\big(D_y\varphi + 2\varphi\big)\|D_x u\|_{\bf h}^2\frac{dy}{y^n}    
+ \left[\frac{\varphi(\|w_{\pm}\|_{\bf h}^2 - \|D_xu\|_{\bf h}^2)}{y^{n-1}}\right]_{y=a}^{y=b},
\label{C2S"Apm}
\end{equation}
\begin{equation}
R_{\pm} = \int_{a}^b(D_y\varphi)\|w_{\pm}\|_{\bf h}^2\frac{dy}{y^n}
- 2{\rm Re}\int_a^b\varphi(f,w_{\pm})_{\bf h}\frac{dy}{y^n}.
\label{C2S"Bpm}
\end{equation}
Then we have the following inequality. 
\begin{equation}
L_+ \leq R_+, \quad L_- \geq R_-, \quad {\rm if} \quad {\rm Im}\,\sqrt z \geq 0,
\label{C2S2ABinequality1}
\end{equation}
\begin{equation}
L_+ \geq R_+, \quad L_- \leq R_-, \quad {\rm if} \quad {\rm Im}\,\sqrt z \leq 0,
\label{C2S2ABinequality2}
\end{equation}
 \end{lemma}

Proof.
Using Lemma 2.4,  $\varphi \geq 0$, and the sign of ${\rm Im}\,\sqrt{z}$, we obtain the lemma.
\qed

\medskip
In the following, $z$ varies over the region
\begin{equation}
J_{\pm} = \{z \in {\bf C} \, ; \, a \leq {\rm Re}\,z \leq b, \ 
0 < \pm {\rm Im}\,z < 1\},
\label{C2S2Jplusminus}
\end{equation}
where $0 < a < b$ are arbitrarily chosen constants.


\begin{lemma}
Let $u = R_0(z)f$ with $f \in \mathcal B$. Then, for any $\epsilon > 0$, there exists a constant $C_{\epsilon} > 0$ such that
$$
\int_0^{\infty}\|D_xu\|_{\bf h}^2\frac{dy}{y^n} \leq \epsilon\|u\|_{\mathcal B^{\ast}}^2 + 
C_{\epsilon}\|f\|_{\mathcal B}^2, \quad \forall z \in J_{\pm}.
$$
\end{lemma}
Proof. Assume that $z \in J_+$. Letting $\varphi = 1$ and using (\ref{C2S2ABinequality1}), we have
$$
\int_a^b\|D_xu\|_{\bf h}^2\frac{dy}{y^n} \leq \left[\frac{\|D_xu\|_{\bf h}^2 - \|w_{+}\|_{\bf h}^2}{2y^{n-1}}\right]_{y=a}^{y=b} +\left|\int_a^b(f,w_{+})_{\bf h}\frac{dy}{y^n}\right|.
$$
By Theorem 1.3 (4),  $w_{+}, D_xu \in L^2$ for $z \not\in {\bf R}$. Hence
\begin{equation}
\liminf_{y\to0}\frac{\|w_{+}\|_{\bf h}^2 + \|D_xu\|_{\bf h}^2}{y^{n-1}} = 0, \quad
\liminf_{y\to\infty}\frac{\|w_{+}\|_{\bf h}^2 + \|D_xu\|_{\bf h}^2}{y^{n-1}} = 0.
\label{C2S2liminf}
\end{equation}
Therefore letting $a \to 0$ and $b \to \infty$ along suitable sequences, we have$$
\int_0^{\infty}\|D_xu\|_{\bf h}^2\frac{dy}{y^n} \leq \left|\int_0^{\infty}(f,w_{+})_{\bf h}\frac{dy}{y^n}\right| \leq \epsilon\|w_{+}\|^2_{\mathcal B^{\ast}} + C_{\epsilon}\|f\|_{\mathcal B}^2.
$$
Theorem 1.3 (1) yields $\|w_{+}\|_{\mathcal B^{\ast}} \leq C(\|u\|_{\mathcal B^{\ast}} + \|f\|_{\mathcal B^{\ast}})$, which proves the lemma when $z \in J_{+}$. The case for $z \in J_-$ is proved similarly by using $w_-$. \qed


\begin{lemma}  Let $u$, $f$ be as in the previous lemma,  and $w_{\pm} = (D_y - \sigma_{\pm})u$. Then for any $\epsilon > 0$, there exists a constant $C_{\epsilon} > 0$ such that, for any $y > 0$, 
\begin{equation}
\frac{\|w_{+}\|_{\bf h}^2 - \|D_xu\|_{\bf h}^2}{y^{n-1}} \leq \epsilon\|u\|_{\mathcal B^{\ast}}^2 +
C_{\epsilon}\|f\|_{\mathcal B}^2, \quad \forall z \in J_+,
\nonumber
\end{equation}
\begin{equation}
\frac{\|w_{-}\|_{\bf h}^2 - \|D_xu\|_{\bf h}^2}{y^{n-1}} \leq \epsilon\|u\|_{\mathcal B^{\ast}}^2 +
C_{\epsilon}\|f\|_{\mathcal B}^2, \quad \forall z \in J_-.
\nonumber
\end{equation}
\end{lemma}
Proof. As in the previous lemma, assume that $z \in J_+$. Letting $\varphi = 1$ and using (\ref{C2S2ABinequality1}), we have
$$
\frac{\|w_{+}\|_{\bf h}^2 - \|D_xu\|_{\bf h}^2}{y^{n-1}}\Big|_{y = b} \leq 
\frac{\|w_{+}\|_{\bf h}^2  - \|D_xu\|_{\bf h}^2}{y^{n-1}}\Big|_{y = a} + C\|f\|_{\mathcal B}\|w_+\|_{{\mathcal B}^{\ast}}.
$$
Using (\ref{C2S2liminf}) and
[letting $a \to 0$ along a suitable sequence,   we obtain the lemma by  Theorem 1.3 (1). \qed


\begin{lemma} Let $u$, $f$, $w_{\pm}$ be as in the previous lemma.
Then, for any $\epsilon > 0$, there exists a constant $C_{\epsilon} > 0$ such that
\begin{equation}
\|w_{+}\|_{{\mathcal B}^{\ast}} \leq \epsilon\|u\|_{\mathcal B^{\ast}} + C_{\epsilon}\|f\|_{\mathcal B}, \quad \forall z \in J_{+},
\nonumber
\end{equation}
\begin{equation}
\|w_{-}\|_{{\mathcal B}^{\ast}} \leq \epsilon\|u\|_{\mathcal B^{\ast}} + C_{\epsilon}\|f\|_{\mathcal B}, \quad \forall z \in J_{-}.
\nonumber
\end{equation}
\end{lemma}

Proof. We divide the inequality in Lemma 2.7 by $y$ and integrate on $(1/R,R)$. We then use Lemma 2.6 to estimate the integral of $\|D_xu\|_{\bf h}^2$, and obtain 
the lemma.
\qed


\begin{lemma} There exists a constant $C > 0$ such that
$$
\|R_0(z)f\|_{{\mathcal B}^{\ast}} \leq C\|f\|_{\mathcal B}, \quad 
\forall z \in J_{\pm}.
$$
\end{lemma}
Proof. We consider the case that $z \in J_+$, and put $\sqrt{z} = k + i\epsilon$ for $z \in J_{+}$. Then $\epsilon > 0$ and $k > C$ for some constant $C > 0$. Letting $w_{+} = (D_y - \sigma_{+})u$, we then have 
\begin{equation}
 {\rm Im}\,D_y(w_{+},u)_{\bf h} = {\rm Im}\,(n - 1 + 2ik)(w_{+},u)_{\bf h} - {\rm Im}\,(f,u)_{\bf h}.
\label{C2S2ImDywu}
\end{equation}
This is a consequence of the formula
$$
D_y(w_{+},u)_{\bf h} = (D_yw_{+},u)_{\bf h} + \|w_{+}\|_{\bf h}^2 + \left(\frac{n-1}{2} + \epsilon + ik\right)(w_{+},u)_{\bf h}
$$
and (\ref{eq:Dyminussigmaplusu}). We integrate (\ref{C2S2ImDywu}). Since

$$
 \int_a^bD_y(w_{+},u)_{\bf h}\frac{dy}{y^n}  = \left[\frac{(w_{+},u)_{\bf h}}{y^{n-1}}\right]_a^b + (n-1)\int_a^b(w_{+},u)_{\bf h}\frac{dy}{y^n},
$$
we then have
\begin{equation}
 {\rm Im}\,\left[\frac{(w_{+},u)_{\bf h}}{y^{n-1}}\right]_a^b = 2k\,{\rm Re}
 \int_a^b(w_{+},u)_{\bf h}\frac{dy}{y^n} - {\rm Im}\int_a^b(f,u)_{\bf h}\frac{dy}{y^n}.
 \label{C2S2Imwyyn-1inte}
\end{equation}
Using $w_{+} = D_yu - \sigma_{+}u$ and integrating by parts, we have
$$
 {\rm Re}\int_a^b(w_{+},u)_{\bf h}\frac{dy}{y^n} = \frac{1}{2}
 \left[\frac{\|u\|_{\bf h}^2}{y^{n-1}}\right]_a^b - \epsilon
 \int_a^b\|u\|_{\bf h}^2\frac{dy}{y^n}.
$$
Therefore (\ref{C2S2Imwyyn-1inte}) is computed as
\begin{equation}
  {\rm Im}\,\left[\frac{(w_{+},u)_{\bf h}}{y^{n-1}}\right]_a^b =  k\,
 \left[\frac{\|u\|_{\bf h}^2}{y^{n-1}}\right]_a^b - 2\epsilon k
 \int_a^b\|u\|_{\bf h}^2\frac{dy}{y^n}
 - {\rm Im}\int_a^b(f,u)_{\bf h}\frac{dy}{y^n},
 \nonumber
\end{equation}
which implies
 \begin{equation}
  {\rm Im}\,\left[\frac{(w_{+},u)_{\bf h}}{y^{n-1}}\right]_a^b \leq k\,
 \left[\frac{\|u\|_{\bf h}^2}{y^{n-1}}\right]_a^b 
 + C\|f\|_{\mathcal B}\|u\|_{{\mathcal B}^{\ast}}.
 \nonumber
\end{equation}
Note that for $z \not\in {\bf R}$, $w_{+}$ and $u$ are in $L^2((0,\infty);L^2({\bf R}^{n-1});dy/y^n)$. Hence, there exists a sequence $b_1 < b_2 <\cdots \to \infty$ such that
$$
\frac{|(w_{+},u)_{\bf h}(b_m)| + \|u(b_m)\|_{\bf h}^2}{b_m^{n-1}} \to 0.
$$
For $w_+$, we take $a = y < b = b_m$  to have
$$
\frac{\|u(y)\|_{\bf h}^2}{y^{n-1}} \leq C_k\left(\frac{\|w_+(y)\|_{\bf h}^2}{y^{n-1}} 
+ \frac{|(w_{+},u)_{\bf h}(b_m)| + \|u(b_m)\|_{\bf h}^2}{b_m^{n-1}} 
+  \|f\|_{\mathcal B}\|u\|_{{\mathcal B}^{\ast}}\right).
$$
Letting $m \to \infty$, we see that
$$
 \frac{\|u(y)\|_{\bf h}^2}{y^{n-1}} \leq C\left(\frac{\|w_+(y)\|_{\bf h}^2}{y^{n-1}}
 + \|f\|_{\mathcal B}\|u\|_{{\mathcal B}^{\ast}}\right).
$$
Dividing by $y$ and integrating from $1/R$ to $R$, we have
$$
 \frac{1}{\log R}\int_{1/R}^R\|u(y)\|_{\bf h}^2 \frac{dy}{y^n}
 \leq \frac{C}{\log R}\int_{1/R}^R\|w_+(y)\|_{\bf h}^2
 \frac{dy}{y^n} + 
  C\|f\|_{\mathcal B}\|u\|_{{\mathcal B}^{\ast}},
$$
which implies 
$$
\|u\|_{{\mathcal B}^{\ast}}^2\leq C\|w_+\|_{{\mathcal B}^{\ast}}^2 + C\|f\|_{\mathcal B}\|u\|_{{\mathcal B}^{\ast}}.
$$
This, together with Lemma 2.8, yields 
$$
\|u\|_{{\mathcal B}^{\ast}} \leq C\|f\|_{\mathcal B}, \quad \forall z \in J_+.
$$
Similarly, we can prove the lemma for $z \in J_-$.
 \qed

\bigskip
Lemma 2.9 completes the proof of Theorem 1.4.2.


\subsection{Radiation conditions and uniqueness theorem} 
The following theorem specifies the fastest decay order of non-trivial solutions to the Helmholtz equation $(H - \lambda)u = 0$.


\begin{theorem} Let $\lambda > 0$.
If $u \in {\mathcal B}^{\ast}$ satisfies $(H - \lambda)u = 0$ for $0 < y < y_0$ with some $y_0 > 0$, and
\begin{equation}
  \liminf_{R\to\infty}\frac{1}{\log R}\int_{1/R}^1\|u(y)\|^2_{L^2({\bf R}^{n-1})}
  \frac{dy}{y^n} = 0,
 \nonumber
\end{equation}
then $u = 0$ for $0 < y < y_0$.
\end{theorem}

We should stress that we have only to assume the equation $(H - \lambda)u = 0$ to be satisfied near $y = 0$.
The proof is given in the next section.

 
\begin{cor}
 $\ \ \ \sigma_p(H)\cap\big(0,\infty\big) = \emptyset$.
\end{cor}

We say that $u \in {\mathcal B}^{\ast}$ satisfies the {\it outgoing radiation condition} (for $\sigma_+$), or {\it incoming radiation condition} (for $\sigma_-$), if
 the following two conditions (\ref{eq:Chap2Sect2RadCond1/R1}) and (\ref{eq:Chap2Sect2RadCondplusinfty}) are fulfilled:
\begin{equation}
\lim_{R\to\infty}\frac{1}{\log R}\int_{1/R}^1
\|(D_y - \sigma_{\pm}(\lambda))u(y)\|^2_{L^2({\bf R}^{n-1})}\frac{dy}{y^n} = 0,
\label{eq:Chap2Sect2RadCond1/R1}
\end{equation}
\begin{equation}
\sigma_{\pm}(\lambda) = \frac{n-1}{2} \mp i\sqrt{\lambda}.
\nonumber
\end{equation}
\begin{equation}
\lim_{R\to\infty}\frac{1}{\log R}\int_{1}^R
\|u(y)\|^2_{L^2({\bf R}^{n-1})}\frac{dy}{y^n} = 0.
\label{eq:Chap2Sect2RadCondplusinfty}
\end{equation}

 
\begin{lemma}
 Assume that $\lambda > 0$ and $u \in \mathcal B^{\ast}$ satisfies the equation $(H - \lambda)u = 0$, and the outgoing or incoming radiation condition. Then
 $u = 0$.
\end{lemma}
Proof. We assume that $u$ satisfies the outgoing radiation condition. We take $0 \leq \rho(t) \in C_0^{\infty}({\bf R})$ satisfying ${\rm supp}\,\rho \subset (-1,1)$, $\int_{-1}^1\rho(t)dt = 1$, and put 
$$
\varphi_R(y) = \chi\big(\frac{\log y}{\log R}\big), \quad 
\chi(t) = \int_{-\infty}^t\rho(s)ds.
$$
Let $(\;,\;)_{\bf h}$ and $\|\cdot\|_{\bf h}$ denote the inner product and the norm of $L^2({\bf R}^{n-1})$, respectively.  We multiply the equation $(H - \lambda)u = 0$ by $\varphi_R(y)\overline{u}$ and integrate over ${\bf R}^{n-1}\times
(0,R)$ to obtain
\begin{equation} \label{2.20c}
\begin{split}
0 &= {\rm Im}\int_{0}^R\left((- D_y^2 + (n - 1)D_y + V)u,\varphi_Ru \right)_{\bf h}\frac{dy}{y^n}  \\
&= - {\rm Im}\frac{(D_yu,u)_{\bf h}}{y^{n-1}}\Big|_{y=R} + {\rm Im}\frac{1}{\log R}\int_{0}^R\rho\big(\frac{\log y}{\log R}\big)\left(D_yu,u\right)_{\bf h}\frac{dy}{y^{n}} \\
&\ \ \  + 
{\rm Im}\int_{0}^R\left(Vu,\varphi_Ru\right)_{\bf h}\frac{dy}{y^n}.
\end{split}
\end{equation}
 Observe that (\ref{eq:Chap2Sect2RadCondplusinfty}) implies, due to Theorem 1.3 (2), that
\begin{equation} \label{2.20a}
\lim_{R \to \infty} \frac{1}{\log R} \int_1^R  ||D_y u||^2_{L^2({\bf R}^{n-1})} \frac{dy}{y^n}=0.
\end{equation}
Indeed, let $\psi(y) \in C^\infty({\bf R}_+),\, \psi=1$ for $y>1$ and $\psi=0$ for $y <1/2$. Then, with
$v= \psi u$,
$$
(H-\lambda) v=f:= \left[H, \psi \right] u \in \mathcal B,
$$
due to Theorem 1.3 (1) and the fact, that $\hbox{supp}(f) \subset \{1/2 <y <1\}.$
Thus, $v$ satisfies conditions of Theorem 1.3 (2), which implies (\ref{2.20a}).

Conditions (\ref{eq:Chap2Sect2RadCondplusinfty}), (\ref{2.20a}) yield that 
\begin{equation} \label{2.20b}
\lim_{R \to \infty} \frac{1}{\log R} \int_{1/R}^R ||(D_y-\sigma_{\pm}) u(y)||_{L^2({\bf R}^{n-1})} 
\frac{dy}{y^n} =0.
\end{equation}
Also (\ref{eq:Chap2Sect2RadCondplusinfty}), (\ref{2.20a}) imply that 
$$
\lim \inf_{y \to \infty} \frac{|(D_y u, u)_{y=a}|}{y^{n-1}} =0.
$$
 We also see that
$$
{\rm Im}\int_{0}^R\left(Vu,\varphi_Ru\right)_{\bf h}\frac{dy}{y^n} \to {\rm Im}
\int_0^{\infty}(Vu,u)_{\bf h}\frac{dy}{y^n} = 0.
$$
Indeed, $\int_0^{\infty}\big|(Vu,u)_{\bf h}d\big|y/y^n < \infty$, since $Vu \in \
\mathcal X^{s}$, $1/2 < s < (1 + \epsilon)/2$ due to (2.4) and
Theorem 1.3 (5). As $V$ is symmetric, this gives the result.

 Hence, by (\ref{2.20c}),  there is a subsequence $R_1 < R_2 < \cdots \to \infty$ such that
\begin{equation}
 {\rm Im}\frac{1}{\log R_j}\int_0^{\infty}
\rho\big(\frac{\log y}{\log R_j}\big) (D_yu,u)_{\bf h}\frac{dy}{y^n} \to 0.
 \nonumber
\end{equation}


Combining this equation with (\ref{2.20b}), we have 
$$
\lim_{j\to\infty}\frac{{\sqrt \lambda}}{\log R_j}\int_0^{\infty}\big(\rho(\frac{\log y}{\log R_j})u,u\big)_{\bf h}\frac{dy}{y^n} = 0, \quad \forall \rho \in C_0^{\infty}({\bf R}). 
$$
This implies that
$$
\lim_{j\to\infty}\frac{1}{\log R'_j}\int_{1/R'_j}^{R_j'}\|u(y)\|_{\bf h}^2\frac{dy}{y^n} = 0
$$
along a suitable sequence $R_1' < R_2' < \cdots \to \infty$.
The lemma then follows from Theorem 2.10. \qed 


\subsection{Proof of Theorem 2.3}
The assertion (1) has been proved in Corollary 2.11. Let $\epsilon$ be as in the condition (C) in Subsection 2.1, and take $s$ such that
\begin{equation}
\frac{1}{2} < s < \frac{1 + \epsilon}{2}.
\nonumber
\end{equation}
Take a compact interval $I \subset (0,\infty)$ arbitrarily, and put
$$
J = \{\lambda \pm i\epsilon \, ; \, \lambda \in I, \ 0 < \epsilon < 1\}.
$$


\begin{lemma}
(1) There exists a constant $C > 0$ such that
\begin{equation}
\sup_{z \in J}\|R(z)f\|_{{\mathcal X}^{-s}}  \leq C\|f\|_{\mathcal B},
\label{Lemma2.2.13.Inequality1}
\end{equation}
\begin{equation}
\sup_{z \in J}\|R(z)f\|_{\mathcal B^{\ast}} \leq C\|f\|_{\mathcal B}.
\label{Lemma2.2.13.Inequality2}
\end{equation}
(2) For any $\lambda > 0$ and $f \in {\mathcal B}$, the strong limit $\lim_{\epsilon \to 0}R(\lambda \pm i\epsilon)f$ exists in ${\mathcal X}^{-s}$. Also, the weak limit $\lim_{\epsilon \to 0}R(\lambda \pm i\epsilon)f$ exists in $\mathcal B^{\ast}$.
\\
\noindent
(3) $\ R(\lambda \pm i0)f$ is an ${\mathcal X}^{-s}$-valued strongly continuous function of $\lambda > 0$, and also a $\mathcal B^{\ast}$-valued weakly continuous function of $\lambda > 0$. In particular, 
$$
\lim_{\epsilon\to 0}(R(\lambda\pm i\epsilon)f,g) = (R(\lambda\pm i0)f,g),
\quad \forall g\in \mathcal B.
$$
\end{lemma}
Proof. 
 If (1) does not hold, there exist $z_n \in J$ and 
$f_n \in {\mathcal B}$ satisfying 
\begin{equation}
 \|f_n\|_{\mathcal B} \to 0, \quad 
\|u_n\|_{{\mathcal X}^{-s}} = 1, \quad u_n = R(z_n)f_n. 
\nonumber
\end{equation}
These imply that 
\begin{equation} \label{E2.2}
(H_0-z_n) u_n=f_n-V u_n
\end{equation}
and we can assume, without loss of generality, that $z_n \to \lambda \in I$.  By Theorem 1.3 (6), 
$$
||D^{\alpha} u_n||_{\chi^{-s}} \leq C,\quad |\alpha| \leq 2.
$$
Therefore, by the condition (C), $V u_n \in \mathcal B$ and 
$$
||V u_n||_{\mathcal B} \leq C.
$$
Returning to (\ref{E2.2}), this implies, due to Lemma 2.9, that 
\begin{equation} \label{E2.3}
||u_n||_{{\mathcal B}^*} \leq C.
\end{equation}
Therefore, there exists a subsequence, which we continue to denote by $u_n$, such that
$$
u_n \rightarrow u,
$$
in the sense of the weak convergence.

On the other hand, applying Theorem 1.3 (4), we see that, with $|\alpha| \leq 2$ and $1/2 <t, t'<s$,
\begin{equation} \label{E2.4}
||D^{\alpha} u_n||_{\chi^{-t}} \leq C \left(||u_n||_{\chi^{-t}}+ ||f_n||_{\chi^{-t}} \right) \leq C;
\end{equation}
\begin{equation} \label{E2.5}
||D^{\alpha} (u_n-u_m)||_{\chi^{-t'}} \leq C \left(||u_n-u_m||_{\chi^{-t'}}+ ||f_n-f_m||_{\chi^{-t'}} +|z_n-z_m|\right).
\end{equation}
These imply, using Rellich's theorem, that there exists a subsequence such that 
$D^{\alpha} u_n \to D^{\alpha} u$ in $\chi^{-s}, \, |\alpha| \leq 2$ and, in particular, $||u||_{\chi^{-s}}=1$.
Also $u_n \to u$ in ${\mathcal B}^*$,  as follows from Lemma 2.9 together with (\ref{E2.2}), (\ref{E2.5}).

Then
$$
u=-R_0(\lambda \pm i0) Vu, \quad Vu \in {\mathcal B},
$$
and, by Corollary 1.4.8 (2) and Lemma 1.4.9, $u$ satisfies the radiation condition. Thus, by Lemma 2.12,
$u=0$, contradicting $||u||_{\chi^{-s}}=1$.
 This completes the proof of (\ref{Lemma2.2.13.Inequality1}). 

To prove (\ref{Lemma2.2.13.Inequality2}),
we observe that $\mathcal B^{\ast}$ is reflexive and, therefore, 
sequentially weakly compact by Theorem V.2.1 of \cite{Yo66}.
We then use (\ref{E2.2}) with $z_n, u_n, f_n$ replaced by $z, R(z)f, f$ and follow the same arguments.
 
 The assertion (2), (3) can be proved by the similar manner.
 \qed

\bigskip
The assertions (2), (3), (4) of Theorem 2.3 are now easily derived from Lemma 2.13 and the resolvent equation $R(z) = R_0(z) - R_0(z)VR(z)$. 
To this end, we use Theorem 1.3 (6) with $s < (1+\epsilon)/2$, (C) in the decay assumption of the metric in subsection 2.1 and Theorem 1.4.2 (3).

For the proof of (5), see  \cite{IkSa72} or \cite{Is04a}, p. 49.  \qed

\bigskip
The following lemma is a consequence of the above proof.


\begin{lemma}
 For any $f \in {\mathcal B}$ and $\lambda > 0$, $u = R(\lambda \pm i0)f$ satisfies the equation $(H - \lambda)u = f$, 
and the radiation condition. 
Conversely, any solution $u \in {\mathcal B}^{\ast}$ of the above equation 
satisfying the radiation condition is unique and is given by $u = R(\lambda \pm i0)f$. 
 \end{lemma}


\section{Growth order of solutions to reduced wave equations}


\subsection{Abstract differential equations}
 Let $X$ be a Hilbert space and consider the following differential equation for an $X$-valued function $u(t)$: 
\begin{equation}
- u^{\prime\prime}(t) + B(t)u(t) + V(t)u(t) - Eu(t) = P(t)u(t), \quad t > 0,
\label{eq:Chap2Sect3Diffeq}
\end{equation}
$E > 0$ being a constant. The following assumptions are imposed. 

\bigskip
\noindent
{\it (A-1)   $B(t)$ is a non-negative self-adjoint operator valued function with domain $D(B(t)) = D \subset X$ independent of $t > 0$. For each $x \in D$, the map $(0,\infty) \ni t \to B(t)x \in X$ is $C^1$, and there exist constants $t_0 > 0$ and $\delta > 0$ such that
\begin{equation}
t\frac{dB(t)}{dt} + (1 + \delta)B(t) \leq 0, \quad \forall t  > t_0.
\label{eq:Chap2Sect3tdBtdt}
\end{equation}
(A-2) For any fixed $t$, $V(t)$ is bounded self-adjoint on $X$ and satisfies
\begin{equation}
V(t) \in C^1((0,\infty);{\bf B}(X)),
\label{eq:Chap2Sect3VtC1}
\end{equation}
\begin{equation}
\frac{1}{t}\|V(t)\| + 
\big\|\frac{dV(t)}{dt}\big\| \leq C(1 + t)^{-1-\epsilon},
\quad \forall t \geq 1,
\label{eq:Chap2Sect3normVt}
\end{equation}
for some constants $C, \epsilon > 0$.\\
\noindent
(A-3) For any fixed $t$, $P(t)$ is a closed (not necessarily self-adjoint) operator on $X$ with domain $D(P(t)) \supset D$ satisfying
\begin{equation}
P(t)^{\ast}P(t) \leq C(1 + t)^{-2-2\epsilon}\big(B(t) + 1\big).
\label{eq:Chap2Sect3B1tastB1t}
\end{equation}
Moreover,
$$
{\rm Re}\,P(t) := \frac{1}{2}\left(P(t) + P(t)^{\ast}\right)
$$
is a bounded operator on $X$ and satsifies}
\begin{equation}
\|{\rm Re}\,P(t)\| \leq C(1 + t)^{-1-\epsilon}, \quad \forall t > 0.
\label{eq:Chap2Sect3RealB1t}
\end{equation}


\begin{theorem}
Under the above assumptions (A-1), (A-2), (A-3), if
$$
\liminf_{t\to\infty}(\|u'(t)\|_X + \|u(t)\|_X) = 0
$$
holds, there exists $t_1 > 0$ such that $u(t) = 0$,  $\forall t > t_1$.
\end{theorem}

The proof below is a modification of the method in \cite{Sa79} p. 29. In the following, $\|\cdot\|_X$ is simply written as $\|\cdot\|$. We put
\begin{equation}
(Ku)(t) = \|u'(t)\|^2 + E\|u(t)\|^2 - (B(t)u(t),u(t)) - (V(t)u(t),u(t)).
\nonumber
\end{equation}


\begin{lemma}
There exist constants $C_1, T_1 > 0$ such that
$$
\frac{d}{dt}(Ku)(t) \geq - {C_1}{(1 + t)^{-1-\epsilon}}(Ku)(t), \quad 
\forall t > T_1.
$$
\end{lemma}
Proof. By choosing $\epsilon$ small enough, we can assume that, in addition to  (A-2) and (A-3), 
\begin{equation}
 \|V'(t)\| \leq C(1 + t)^{-1-2\epsilon}.
 \label{eq:Chap2Sect3Vprimet}
\end{equation}
By the equation (\ref{eq:Chap2Sect3Diffeq})
\begin{eqnarray*}
\frac{d}{dt}(Ku)(t) & = & 2{\rm Re}\,\Big[(u'',u') + E(u,u') - (Bu,u') 
- (Vu,u')\Big]   - ((B' + V')u,u) \\
& =&  - 2{\rm Re}\,(Pu,u') - ((B' + V')u,u). 
\end{eqnarray*}
By (\ref{eq:Chap2Sect3B1tastB1t})
\begin{equation}
\|Pu\| \leq C\big(1 + t)^{-1-\epsilon}(\sqrt{(Bu,u)} + \|u\|\big).
\label{eq:Chap2Sect3EstimatePt}
\end{equation}
By (\ref{eq:Chap2Sect3Vprimet}), there exists $t_0 = t_0(\epsilon) > 0$ such that for $t > t_0$ 
$$
|(V'(t)u,u)| \leq \frac{\epsilon}{2}(1 + t)^{-1-\epsilon}\|u\|^2.
\nonumber
$$
By (\ref{eq:Chap2Sect3tdBtdt})
$$
- (B'u,u) \geq \frac{1+\delta}{t}(Bu,u).
$$
Putting the above estimates together we have that there is $C_{\epsilon} > 0$ such that 
for $t > t_0$
\begin{eqnarray*}
\frac{d}{dt}(Ku)(t) 
&\geq& - Ct^{-1-\epsilon}(\|u'\|^2 + \|u\|\|u'\| + \frac{\epsilon}{2}\|u\|^2)
 + \frac{1}{t}(Bu,u) \\
& \geq& - C_{\epsilon}t^{-1-\epsilon}\|u'\|^2  - C\epsilon t^{-1-\epsilon}\|u\|^2
 + \frac{1}{t}(Bu,u). 
\end{eqnarray*}
We rewrite the right-hand side as
\begin{eqnarray*}
& &-C_{\epsilon}t^{-1-\epsilon}(\|u'\|^2 + E\|u\|^2) + (C_{\epsilon}E- C\epsilon)t^{-1-\epsilon}\|u\|^2 + \frac{1}{t}(Bu,u) \\
&=& - C_{\epsilon}t^{-1-\epsilon}(Ku)(t) \\
& &+ (C_{\epsilon}E- C\epsilon)t^{-1-\epsilon}\|u\|^2 - 
C_{\epsilon}t^{-1-\epsilon}(Vu,u)  + 
\big(\frac{1}{t} - \frac{C_{\epsilon}}{t^{1+\epsilon}}\big)(Bu,u).
\end{eqnarray*}
Choose $C_{\epsilon}$ large enough so that $C_{\epsilon}E - C\epsilon \geq \frac{1}{2}C_{\epsilon}E$. Using (\ref{eq:Chap2Sect3normVt}), choose $t_0 = t_0(\epsilon,C_{\epsilon})$ such that, for $t > t_0$, $\frac{E}{2}\|u\|^2 - (Vu,u) \geq 0$, and $1 - Ct^{-\epsilon} > 0$. Thus, the 3rd line is non-negative for  $t > t_0$. Hence the lemma is proved. \qed

\bigskip
Let $m > 0$ be an integer and put
\begin{equation}
(Nu)(t) = t\left[K(e^{d(t)}u) + \frac{m^2 - \log t}{t^{2\alpha}}\|e^{d(t)}u\|^2\right],
\nonumber
\end{equation}
\begin{equation}
\frac{1}{3} < \alpha < \frac{1}{2}, \quad d(t) = \frac{m}{1 - \alpha}
t^{1-\alpha}.
\nonumber
\end{equation}

\begin{lemma} If ${\rm supp}\,u(t)$ is unbounded, there exist constants
 $m_1 \geq 1$, $T_2 \geq T_1$ such that
\begin{equation}
(Nu)(t) \geq 0, \quad \forall t \geq T_2, \quad \forall m \geq m_1.
\nonumber
\end{equation}
\end{lemma}
Proof. Letting $w(t) = e^{d(t)}u(t)$, we have
\begin{equation}
\begin{split}
\frac{d}{dt}(Nu) & = Kw + t\frac{d}{dt}(Kw) + (1 - 2\alpha)\frac{m^2 - \log t}{t^{2\alpha}}\|w\|^2 \\
& \ \ \ - t^{-2\alpha}\|w\|^2 + 2(m^2 - \log t)t^{1-2\alpha}{\rm Re}\,(w',w)\\
&= \|w'\|^2 + \Big(E + (1 - 2\alpha)\frac{m^2 - \log t}{t^{2\alpha}} - t^{-2\alpha}\Big)\|w\|^2 \\
& \ \ \  - (Bw,w) - (Vw,w) + t\frac{d}{dt}(Kw) \\
& \ \ \ + 2t^{1-2\alpha}(m^2 - \log t){\rm Re}\,(w',w).
\end{split}
\label{eq:Chap2Sect3ddtNt}
\end{equation}
By  direct computation, 
\begin{equation}
\begin{split}
w' & = e^du' + mt^{-\alpha}w, \\
w'' & = e^du'' + mt^{-\alpha}e^du' + mt^{-\alpha}w' - \alpha mt^{-\alpha - 1}w 
\\ 
&= Bw + Vw - Ew + 2mt^{-\alpha}w' \\
&  - \left[P
 + (\alpha m t^{-\alpha -1} +m^2t^{-2\alpha})\right]w.
\end{split}
\nonumber
\end{equation}
Hence,
\begin{equation}
\begin{split}
\frac{d}{dt}(Kw) & = 2{\rm Re}\,(w'' + Ew - Vw - Bw,w') - 
(\big(B' + V'\big)w,w) \\
&=  4mt^{-\alpha}\|w'\|^2 
 - 2(\alpha mt^{-\alpha - 1} + m^2t^{-2\alpha}){\rm Re}\,(w,w') \\
&\ \ \  - (\big(B' + V'\big)w,w) - 2 {\rm Re}\,(Pw,w'). 
\end{split}
\label{eq:Chpa2Sect3ddtKw}
\end{equation}
By (\ref{eq:Chap2Sect3ddtNt}) and (\ref{eq:Chpa2Sect3ddtKw}) we have
\begin{equation}
\begin{split}
&\ \ \ \ \  \frac{d}{dt}(Nu) \\
&= (4mt^{1-\alpha} + 1)\|w'\|^2 + 
\{E + (1 - 2\alpha)t^{-2\alpha}(m^2 - \log t) - t^{-2\alpha}\}\|w\|^2 \\
&\ \ \ - 2(\alpha mt^{-\alpha} + t^{1-2\alpha}\log t)\;{\rm Re}\;(w,w') 
 - ((V + tV')w,w) \\
 & \ \ \ - ((tB' + B)w,w)  - 2t{\rm Re}\,(Pw,w') \\
& =: I_1 + I_2 + I_3.
\end{split}
\nonumber
\end{equation} 
For large $t > 0$, $I_1$ is estimated from below as
\begin{equation}
 I_1 \geq (4mt^{1-\alpha} + 1)\|w'\|^2 + \big(\frac{E}{2} + 
 (1 - 2\alpha)t^{-2\alpha}m^2\big)\|w\|^2.
 \nonumber
\end{equation}
By (\ref{eq:Chap2Sect3normVt}), 
$I_2$ is estimated from below as
\begin{eqnarray*}
 I_2 &\geq& - 2(\alpha m t^{-\alpha} + t^{1-2\alpha}\log t)\|w\|\|w'\| - Ct^{-\epsilon}\|w\|^2 \\
 &\geq& - \epsilon m^2 t^{-2\alpha}\|w\|^2 - C_\epsilon\|w'\|^2 \\
 & & - 2t^{1-2\alpha}\log t\|w\|\|w'\| - Ct^{-\epsilon}\|w\|^2.
 \nonumber
\end{eqnarray*}
By (\ref{eq:Chap2Sect3tdBtdt}),
$I_3$ is stimated from below as
\begin{eqnarray*}
I_3 \geq \delta(Bw,w) - 
2t\|Pw\|\cdot\|w'\|.
\end{eqnarray*}
Using (\ref{eq:Chap2Sect3EstimatePt}), we estimate the 2nd term as 
\begin{eqnarray*}
  2t\|Pw\|\cdot\|w'\| \leq  \frac{1}{2}\|w'\|^2 + Ct^{-\epsilon}((Bw,w) + \|w\|^2).
\end{eqnarray*}
Therefore for large $t$, we have
\begin{equation}
 I_3 \geq - \frac{1}{2}\|w'\|^2 - Ct^{-\epsilon}\|w\|^2.
 \nonumber
\end{equation}
Putting the above estimates together, we then have
\begin{equation}
 \frac{d}{dt}(Nu) \geq \frac{7}{2}mt^{1-\alpha}\|w'\|^2 + \frac{E}{3}\|w\|^2 - 
 2t^{1-2\alpha}\log t\,\|w\|\|w'\|.
 \nonumber
 \end{equation}
Finally, we use the inequality
\begin{equation}
 t^{1-2\alpha}\log t\|w\|\|w'\| \leq \epsilon t^{1-\alpha}\|w'\|^2 + 
 C_{\epsilon}t^{1-3\alpha}(\log t)^2\|w\|^2
 \nonumber 
\end{equation}
and $1 - 3\alpha < 0$. Then there is $t_0 > 0$ independent of $m$ such that
\begin{equation}
  \frac{d}{dt}(Nu)(t) \geq 3mt^{1-\alpha}\|w'\|^2 + \frac{E}{4}\|w\|^2 \geq 0 
 \label{eq:Chap2SEct3ddtNutgeq3m}
\end{equation}
for $t > t_0$. 

On the other hand, $Nu(t)$ can be	 rewritten as
\begin{equation}
 \begin{split}
  (Nu)(t) & = te^{2d}\big[\|mt^{-\alpha}u + u'\|^2 + E\|u\|^2 \\
  & \ \ \ - (Bu,u) - (Vu,u) + t^{-2\alpha}(m^2 - \log t)\big]\|u\|^2 \\
  &= te^{2d}\big[2t^{-2\alpha}\|u\|^2m^2 + 2t^{-\alpha}{\rm Re}\,(u,u')m \\
  & + (Ku - t^{-2\alpha}\|u\|^2\log t)\big].
 \end{split}
 \label{eq:Chap2Sect3Nutte2d}
\end{equation}
By the assumption of the lemma, ${\rm supp}\,u(t)$ is unbounded. Therefore, there is $T_2 > t_0$ such that $\|u(T_2)\| > 0$.
By choosing $m_1$ large enough, we then have 
\begin{equation}
(Nu)(T_2) > 0, \quad \forall m > m_1.
\label{eq:Chap2Sect3NuT2}
\end{equation}
The inequalities (\ref{eq:Chap2SEct3ddtNutgeq3m}) and (\ref{eq:Chap2Sect3NuT2}) prove the lemma. \qed

\bigskip
\noindent
{\it Proof of Theorem 3.1}. We show that if ${\rm supp}\,u(t)$ is unbounded,
\begin{equation}
\liminf_{t\to\infty}(\|u'(t)\|^2 + \|u(t)\|^2) > 0
\label{eq:Chap2Sect3liminfuprimet}
\end{equation}
holds. We first consider the case in which there exists a sequence $t_n \to \infty$ such that $(Ku)(t_n) > 0 \ (n = 1, 2, \cdots)$. Let $T_1$ be as in  Lemma 3.2. Then for some $T > T_1$, $(Ku)(T) > 0$. We show that  $(Ku)(t) \geq 0, \ \forall t > T$.
In fact Lemma 3.2 implies
\begin{equation}
\frac{d}{dt}\left\{\exp\left(C_1\int_T^t(1 + s)^{-1-\epsilon}ds\right)
(Ku)(t)\right\} \geq 0, \quad \forall t > T.
\nonumber
\end{equation}
Hence,
\begin{equation}
(Ku)(t) \geq \exp\left(- C_1\int_T^t(1 + s)^{-1-\epsilon}ds\right)
(Ku)(T), \quad \forall t > T.
\nonumber
\end{equation}
This then implies that, for $t >t(E)$,
\begin{equation}
\begin{split}
\|u'(t)\|^2 + E\|u(t)\|^2 &= Ku(t) + (B(t)u(t),u(t)) + (V(t)u(t),u(t)) \\
&\geq \exp\left(- C_1\int_T^t(1 + s)^{-1-\epsilon}ds\right)
(Ku)(T) \\
& \ \ \  - CEt^{-\epsilon}\|u(t)\|^2.
\end{split}
\nonumber
\end{equation}
Therefore, we arrive at
\begin{equation}
\liminf_{t\to\infty}(\|u'(t)\|^2 + \|u(t)\|^2) \geq 
\frac{1}{2}\exp\left(- C_1\int_T^{\infty}(1 + s)^{-1-\epsilon}ds\right)
(Ku)(T) > 0.
\nonumber
\end{equation}
 
We next consider the case in which $(Ku)(t) \leq 0$ for all $ t$ large enough.
Lemma 3.3 and (\ref{eq:Chap2Sect3Nutte2d}) show that, for large $t$,
\begin{equation}
2t^{-2\alpha}\|u(t)\|^2m^2 + 2t^{-\alpha}{\rm Re}\,(u(t),u'(t))m
- t^{-2\alpha}\|u(t)\|^2\log t \geq 0,
\nonumber
\end{equation}
which together with
\begin{equation}
 \frac{d}{dt}\|u(t)\|^2 = 2{\rm Re}\,(u(t),u'(t)),
 \nonumber
\end{equation}
yields, for large $t > 0$, that 
\begin{equation}
\frac{d}{dt}\|u(t)\|^2 \geq t^{-\alpha}\left(\frac{1}{m}\log t - 
2m\right)\|u(t)\|^2 \geq 0.
\label{eq:Chap2Sect3ddtnormut}
\end{equation}
Since the support of $u(t)$ is unbounded, by choosing $T$ large enough so that $\|u(T)\| > 0$.
In view of (\ref{eq:Chap2Sect3ddtnormut}), we then have
\begin{equation}
 \|u(t)\| \geq \|u(T)\| > 0, \quad \forall t > T,
 \nonumber
\end{equation}
which proves (\ref{eq:Chap2Sect3liminfuprimet}). \qed


\subsection{Canonical form} 
In order to apply Theorem 3.1 to the operator $H$ in the previous section, we transform the metric $ds^2$ into the following canonical form.

\begin{theorem}
Let $ds^2$ be the Riemannian metric satisfying the condition (C). Choose a sufficiently small $y_0 > 0$. Then there exists a diffeomorphism $(x,y) \to (\overline{x},\overline{y})$ in the region $0 < y < y_0$ such that
$$
|\partial_{\overline x}^{\alpha}D_{\overline y}^{\beta}(\overline{x} - x)| \leq C_{\alpha\beta}(1 + d_h(x,y))^{-{\rm min}(|\alpha|+\beta,1)-1-\epsilon/2}, \quad \forall \alpha, \beta,
$$
$$
|\partial_{\overline x}^{\alpha}D_{\overline y}^{\beta}\Big(\frac{\overline{y} - y}{\overline y}\Big)| \leq C_{\alpha\beta}(1 + d_h(x,y))^{-{\rm min}(|\alpha|+\beta,1)-1-\epsilon/2}, \quad \forall \alpha, \beta,
$$
and in the $(\overline{x},\overline{y})$ coordinate system, the Riemannian metric takes the form
\begin{equation}
ds^2 = (\overline y)^{-2}\left((d\overline{x})^2 + (d\overline{y})^2 + \sum_{i,j=1}^{n-1}b_{ij}(\overline{x},\overline{y})d\overline{x}^id\overline{x}^j\right).
\nonumber
\end{equation}
Here $b_{ij}(\overline{x}^i,\overline{x}^j)$ satisfies the condition (C) with $\epsilon$ replaced by $\epsilon/2$.
\end{theorem}

The point is that there is no cross term $d\overline x^id\overline{y}$. The proof is a slight modification of the one given in Chap. 4, \S 2. This theorem also holds for the asymptotically hyperbolic ends with regular infinity to be discussed in Chap. 3, \S 2.

Let us prove Theorem 2.10. In the coordinate system of Theorem 3.4, (denoting $(\overline x, \overline y)$ by $(x,y)$), the equation $(- \Delta_g- \frac{(n-1)^2}{4} - \lambda)u = 0$ becomes
$$
\Big(-\frac{1}{\sqrt{g}}\partial_y\big(\sqrt{g}g^{nn}\partial_y\big) - 
\sum_{i,j=1}^{n-1}\frac{1}{\sqrt g}\partial_{x_i}\big(\sqrt{g}g^{ij}\partial_{x_j}\big) - \frac{(n-1)^2}{4}
- \lambda\Big)u = 0.
$$
This is rewritten as
$$
\Big(- D_y^2 + hD_y - \sum_{i,j=1}^{n-1}D_{x_i}h^{ij}D_{x_j} - \frac{(n-1)^2}{4} + Q - \lambda
\Big)u = 0,
$$
where $Q= \sum_{i=1}^{n-1} b_i(x, y) D_i+c(x,y)$. Here
$h - (n-1), h^{ij} - \delta^{ij}$ and $Q$ satisfy the condition (C), since for $y$ close to 0, $d_h(x,y)$ and $\rho(x,y)$ are equivalent. Putting $t = - \log y$ and $u = v\exp(-\frac{1}{2}\int_{t_0}^th(x,e^s)ds)$, we have
\begin{equation}
(- \partial_t^2 + B(t) - \lambda)v = P(t)v,
\nonumber
\end{equation}
where
$$
B(t) = - e^{-2t}\sum_{i,j=1}^{n-1}\partial_{x_i}({\delta_{ij}} + a_{ij}(t,x))\partial_{x_j},
$$
$$
P(t) = - e^{-t}\sum_{i=1}^{n-1}b_i(t,x)\partial_{x_i} + c(t,x),
$$
and, for large $t>0$,  $a_{ij}, b_i, c$ satisfy
$$
|\partial_x^{\alpha}\partial_t^{\beta} m(t,x)| \leq C_{\alpha\beta}(1 + t)^{-\beta -1-\epsilon}, \quad \forall \alpha, \beta.
$$
We have, therefore, for large $t > 0$
$$
tB'(t) + 2B(t) = - \sum{i, j=1}^{n-1}\partial_{x_i}e^{-2t}\{(-2t + 2)(\delta_{ij} + a_{ij}) + \partial_ta_{ij})\}\partial_{x_j} \leq 0,
$$
Hence, with $X = L^2({\bf R}^{n-1})$, the assumption (\ref{eq:Chap2Sect3tdBtdt}) is satisfied. Rewriting $P(t)^{\ast}P(t)$ as 
$$
P(t)^{\ast}P(t) = \sum_{|\alpha|\leq2}a_{\alpha}(t,x)(D_x)^{\alpha}, \quad
D_x = e^{-t}\partial_x,
$$
we have, for any $\varphi \in C_0^{\infty}({\bf R}^{n-1})$,
\begin{eqnarray*}
(P(t)^{\ast}P(t)\varphi,\varphi) &\leq& C(1 + t)^{-2-2\epsilon}
\Big(e^{-2t}\|\partial_x\varphi\|^2 + \|\varphi\|^2\Big) \\
&\leq & C(1 + t)^{-2-2\epsilon}\Big((B(t)\varphi,\varphi) + (\varphi,\varphi)\Big),
\end{eqnarray*}
which proves (\ref{eq:Chap2Sect3B1tastB1t}).  Note that as $t \to \infty$, $y \to 0$ and 
\begin{equation}
\exp(-\frac{1}{2}\int_{t_0}^th(x,e^s)ds) = y^{(n-1)/2}\Big(1 + O(|\log y|^{-1})\Big).
\label{C2S3expequalyn-1}
\end{equation}

Our next goal is to show that the condition in Theorem 3.1 is satisfied. 
To this end, we return to the proof of Theorem 1.3 (2). Take $\chi(t) \in C_0^{\infty}({\bf R})$ such that $\chi(t) = 1$ for $-1 < t < -1/2$, and $\chi(t) = 0$ for $t < -2$ or $t > - 1/4$. Take $\psi \in C_0^{\infty}({\bf R})$ such that $\psi = 1$ on ${\rm supp}\,\chi$, and $\psi(t)= 0$ for $t > 0$ or $t < -3$. Then 
the estimate (\ref{eq:Chap2Sec1chiRD}) is valid for this choice of $\chi$ and $\psi$. Following the arguments after this inequality, we obtain
$$
\liminf_{R\to\infty}\frac{1}{\log R}\int_{R^{-1}}^{R^{-1/2}}
\|D_iu(y)\|^2\frac{dy}{y^n} = 0
$$
if the condition of Theorem 2.10 is satisfied. This implies that
$$
\liminf_{y\to 0}\frac{\|D_yu(y)\|^2 + \|u(y)|^2}{y^{n-1}} = 0.
$$
Since $t = - \log y$, it follows from this formula together with (\ref{C2S3expequalyn-1}) that 
$$
\liminf_{t\to\infty}\left(\|v'(t)\| + \|v(t)\|\right) = 0.
$$
Therefore, by Theorem 3.1, $v(t) = 0$ for large $t$, i.e. $u(y) = 0$ for small $y$. By the unique continuation theorem, this in turn imples that $u(y) = 0$ for $y < y_0$. \qed

\medskip


\subsection{Asymptotically Euclidean metric} Let us remark that Theorem 3.1 also applies to asymptotically Euclidean metrics on ${\bf R}^n$. In fact, given a metric $g_{ij}(x)$ satisfying
$$
|\partial_x^{\alpha}(g_{ij}(x) - \delta_{ij})| \leq C_{\alpha}(1 + |x|)^{-|\alpha|-1-\epsilon_0}, \quad \forall \alpha,
$$
one can construct a diffeomorphism near infinity such that this metric is transformed into
$$
(dr)^2 + r^{2}h(r,\omega,d\omega), \quad r > r_0, \quad \omega \in S^{n-1},
$$
where $h(r,\omega,d\omega)$ is a positive definite metric on $S^{n-1}$, and behaves like $h_0(\omega,d\omega)$ at infinity, where $h_0(\omega,d\omega)$ is the standard metric on $S^{n-1}$ (see Appendix A, \S 2).


\section{Abstract theory for spectral representations}  


\subsection{Basic ideas} 
Let $H = \int_{-\infty}^{\infty}\lambda dE(\lambda)$ be a self-adjoint operator on a Hilbert space $\mathcal H$, and $I$  an open interval contained in $\sigma_{ac}(H)$. Let $\bf h$ be an auxiliary Hilbert space and  
$\widehat{\mathcal H} = L^2(I;{\bf h};\rho(\lambda)d\lambda)$ the Hilbert space of all $\bf h$-valued $L^2$-functions on $I$ with respect to the measure $\rho(\lambda)d\lambda$. By a spectral representation of $H$ on $I$, we mean a unitary operator $U : E(I)\mathcal H \to \widehat{\mathcal H}$ such 
that
\begin{equation}
(UHf)(\lambda) = \lambda (Uf)(\lambda), \quad \forall f \in D(H), \quad \forall \lambda \in I.
\nonumber
\end{equation}
We mainly consider the following situation. There exist Banach spaces ${\mathcal H_+, \mathcal H}_-$ such that ${\mathcal H}_+ \subset {\mathcal H} \subset {\mathcal H}_-$ and for $\lambda \in I$, $\lim_{\epsilon \downarrow 0} (H - \lambda \mp i\epsilon)^{-1}$ exists as a bounded operator in ${\bf B}({\mathcal H}_+;{\mathcal H}_-)$. For the limits $(H - (\lambda \pm i0))^{-1}$ one can associate the
operators $U_{\pm}(\lambda) \in {\bf B}({\mathcal H}_+;{\bf h})$ and the spectral representations 
$U_{\pm}$ satisfying
\begin{equation}
(U_{\pm}f)(\lambda) = U_{\pm}(\lambda)f, \quad \forall \lambda \in I, \quad \forall f \in {\mathcal H}_+.
\nonumber
\end{equation}
Then there is a unitary operator $\widehat S(\lambda)$ on $\bf h$ such that
\begin{equation}
 U_+(\lambda) = \widehat S(\lambda)U_-(\lambda) , \quad \forall \lambda \in I.
 \nonumber
\end{equation}
This $\widehat S(\lambda)$ is called the scattering matrix or S-matrix. The two limits $\lim_{\epsilon \downarrow 0} (H - \lambda \mp i\epsilon)^{-1}$ appear naturally in computing the limit $\lim_{t\to\pm\infty} e^{-itH}$. Hence, the S-marix is closely related with the asymptotic behavior of solutions to the time-dependent Schr{\"o}dinger equation $i\partial_tu = Hu$. However, the scattering matrix depends on the spectral representations $U_{\pm}$ so that there exist apparently different S-matrices for the same operator $H$. In this and the next sections, we shall introduce three kinds of S-matrices and study their relationships in the case of ${\bf R}^n$ and ${\bf H}^n$. We begin with an abstract framework.


\subsection{Stationary wave operators}
Assume that we are given a Hilbert space ${\mathcal H}$ and Banach spaces ${\mathcal H}_{\pm}$ with norms $\|\cdot\|$, 
and $\|\cdot\|_{\pm}$ satisfying
\begin{equation}
 {\mathcal H}_+ \subset {\mathcal H} \subset {\mathcal H}_-,
\quad
\|f\|_{-} \leq \|f\| \leq \|f\|_{+} , \quad 
 \forall f \in {\mathcal H}_+.
 \nonumber
\end{equation}
We also assume that the above inclusions are dense, and that the inner product $(\,,\,)$  of ${\mathcal H}$ is naturally identified with the coupling of ${\mathcal H}_+$ and ${\mathcal H}_-$. This means that there exists an isometry
 $T : {\mathcal H}_- \to ({\mathcal H}_+)^{\ast}$ such that 
\begin{equation}
 \langle f,Tu\rangle = (f,u), \quad \forall f \in {\mathcal H}_+, 
 \quad \forall u \in {\mathcal H},
 \nonumber
\end{equation}
where $\langle f,v\rangle$ denotes the value $v(f)$ of $v \in \left({\mathcal H}_+\right)^{\ast}$ for $f \in {\mathcal H}_+$. In this case we simply write ${\mathcal H}_- = \left({\mathcal H}_+\right)^{\ast}$.

Let $H_j, j = 1, 2$, be self-adjoint operators on ${\mathcal H}$ such that $D(H_1) = D(H_2)$. For $j = 1, 2$, we put $R_j(z) = (H_j - z)^{-1}$.
Since $D(H_1) = D(H_2)$, we have
\begin{equation}
(H_2 - H_1)R_j(z) \in {\bf B}(\mathcal H;\mathcal H), \quad z \not\in {\bf R}.
\label{C2S4V21Resolbnded}
\end{equation}

Now for $j = 1, 2$, we assume the following:

\medskip
\noindent
(A-1) {\it For any $\varphi(\lambda) \in C_0^{\infty}({\bf R})$, $\varphi(H_j)\mathcal H_+ \subset \mathcal H_+$.} 

\medskip
\noindent
(A-2) {\it There exists an open set $I \subset {\bf R}$ such that $\sigma_p(H_j)\cap I = \emptyset$, and the following strong limit exists}
\begin{equation}
\lim_{\epsilon \to 0} R_j(\lambda \pm i\epsilon) =: R_j(\lambda \pm i0) \in 
{\bf B}({\mathcal H}_+ ; {\mathcal H}_-), \quad \forall \lambda \in I.
\nonumber
\end{equation}
{\it Moreover for any $f \in {\mathcal H}_+$,
$I \ni \lambda \to R_j(\lambda \pm i0)f \in {\mathcal H}_-$ is strongly continuous.}

\medskip
\noindent
(A-3) {\it We put $G_{jk}(z) = (H_j - z)R_k(z)$ for $z \not\in {\bf R}$, and assume that for $\lambda \in I$, $\epsilon > 0$ there exists a strong limit 
$$
\lim_{\epsilon \to 0}G_{jk}(\lambda \pm i\epsilon) =: G_{jk}(\lambda \pm i0) \in {\bf B}({\mathcal H}_+ ; {\mathcal H}_+).
$$
 Furthermore for any $f \in {\mathcal H}_+$, $I \ni \lambda \to G_{jk}(\lambda \pm i0)f \in {\mathcal H}_+$ is strongly continuous.}

\bigskip
We first introduce an operator which shows the similarity of $H_1$ and $H_2$.
Let $E_j(\lambda)$ be the spectral measure for $H_j$, and for $\lambda \in I$, put
\begin{equation}
E_j'(\lambda) = \frac{1}{2\pi i}\left(R_j(\lambda + i0) - R_j(\lambda - i0)\right).
\nonumber
\end{equation}
 By the assumption (A-2), $E_j'(\lambda) \in {\bf B}({\mathcal H}_+ ; {\mathcal H}_-)$. Now for any compact interval $e \subset I$ and $f \in {\mathcal H}_+$, we define
\begin{equation}
\Omega_{jk}^{(\pm)}(e)f = \int_e E_j'(\lambda)G_{jk}(\lambda \pm i0)fd\lambda.
\nonumber
\end{equation}
This is called the stationary wave operator. By the above assumptions, $\Omega_{jk}^{(\pm)}(e) \in {\bf B}({\mathcal H}_+;{\mathcal H}_-)$. However, we have the following stronger results. Let us recall one terminology. For two Hilbert spaces $\mathcal H_1$ and $\mathcal H_2$, closed subspaces $S_1 \subset \mathcal H_1$ and $S_2 \subset {\mathcal H}_2$ and  $U \in {\bf B}({\mathcal H}_1;{\mathcal H}_2)$, we say that $U$ is a partial isometry from ${\mathcal H}_1$ to $\mathcal H_2$ with initial set $S_1$ and final set $S_2$ if $U : S_1 \to S_2$ is unitary and $U : S_1^{\perp} \to 0$. $U$ is a partial isometry if and only if $U^{\ast}U$ and $UU^{\ast}$ are orthogonal projections onto its initial set $S_1$ and final set $S_2$, respectively.


\begin{theorem} Let $e$ be any compact interval in $I$. \\
\noindent
(1) $\Omega_{jk}^{(\pm)}(e)$ is uniquely extended to a bounded operator on ${\mathcal H}$, and is a partial isometry with initial set $E_k(e){\mathcal H}$ and final set $E_j(e){\mathcal H}$. \\
\noindent
(2) $(\Omega_{jk}^{(\pm)}(e))^{\ast} = \Omega_{kj}^{(\pm)}(e)$, where $^{\ast}$ means the adjoint in ${\mathcal H}$. \\
\noindent
(3) $\Omega_{jk}^{(\pm)}(e)$ intertwines $H_j$ and $H_k$. That is, for any bounded Borel function $\varphi(\lambda)$,
$$
\varphi(H_j)\Omega_{jk}^{(\pm)}(e) = \Omega_{jk}^{(\pm)}(e)\varphi(H_k).
$$
\end{theorem}

Theorem 4.1 is proved through a series of Lemmas.


\begin{lemma}
Let $f(\lambda), g(\lambda)$ be $\mathcal H_+$-valued bounded measurable functions on $I$, and $e, e'$ compact intervals in $I$. We put
$$
\varphi = \int_eE_j'(\lambda)f(\lambda)d\lambda, \quad 
\psi = \int_{e'}E_j'(\lambda)g(\lambda)d\lambda.
$$
Then $\varphi, \psi \in {\mathcal H}$ and
$$
(\varphi,\psi) = \int_{e\cap e'}(E_j'(\lambda)f(\lambda),g(\lambda))d\lambda.
$$
\end{lemma}
Proof. If $f(\lambda), g(\lambda)$ are constant functions $f$ and $g$, by Stone's formula, $\varphi = E_j(e)f, \psi = E_j(e')g$. Hence,
$$
(\varphi,\psi) = (E_j(e\cap e')f,g) = \int_{e\cap e'}(E_j'(\lambda)f,g)d\lambda.$$
If $f(\lambda)$, $g(\lambda)$ are step functions, i.e. $f(\lambda) = \sum_{n}\chi_n(\lambda)f_n, g(\lambda) = \sum_n\chi_n(\lambda)g_n$, $\chi_n(\lambda)$ being a characteristic function of the interval $e_n$, $\varphi$ and $\psi$ are written as
$$
\varphi = \sum_nE_j(e\cap e_n)f_n, \quad 
\psi = \sum_n E_j(e'\cap e_n)g_n.
$$
Therefore,
\begin{eqnarray*}
(\varphi,\psi) &=& \sum_{m,n}(E_j(e\cap e'\cap e_m\cap e_n)f_m,g_n) \\
&=& \sum_{m,n}\int_{e\cap e'\cap e_m\cap e_n}(E_j'(\lambda)f_m,g_n)d\lambda \\
&=& \int_{e\cap e'}(E_j'(\lambda)f(\lambda),g(\lambda))d\lambda.
\end{eqnarray*}
Hence, the lemma holds for step functions. 

Let $f(\lambda), g(\lambda)$ be  bounded measurable functions, i.e. 
we can approximate them by step functions $f_m(\lambda), g_n(\lambda)$ such that\begin{equation}
\lim_{m\to\infty}\|f(\lambda) - f_m(\lambda)\|_+ = 0 \quad {a.e.}
\label{C2S4Strongmeasurable}
\end{equation}
and similarly for $g$.
We put 
$$
\varphi_m = \int_eE_j'(\lambda)f_m(\lambda)d\lambda, \quad
\psi_n = \int_{e'}E_j'(\lambda)g_n(\lambda)d\lambda.
$$
Then we have
$$
\|\varphi_m - \varphi_{m'}\|^2 = \int_e(E_j'(\lambda)(f_m(\lambda) - f_{m'}(\lambda),
f_m(\lambda) - f_{m'}(\lambda))d\lambda \to 0,
$$
when $m , m' \to \infty$. Indeed, assumption (A-2) and  boundedness of 
$f$ imply that the integrand is uniformly bounded with respect to $m, m'$. 
Also (\ref{C2S4Strongmeasurable}) implies that this integrand tends to $0$ a.e.
By Lebesgue's theorem, the result follows. 

Thus,
the sequence $\{\varphi_m\}$ converges to $\varphi$ in ${\mathcal H}$ and similaly, $\{\psi_n\}$ converges to $\psi$. Moreover, letting $m, n \to \infty$ in the formula
$$
(\varphi_m,\psi_n) = \int_{e\cap e'}(E_j'(\lambda)f_m(\lambda),g_n(\lambda))d\lambda,
$$
we complete the proof of the lemma. \qed


\begin{lemma}
If $f, g \in {\mathcal H}_+$ and $e, e'$ are compact intervals in $I$, we have
$$
\Omega_{jk}^{(\pm)}(e)f, \ \Omega_{jk}^{(\pm)}(e')g \in \mathcal H,
$$
$$
(\Omega_{jk}^{(\pm)}(e)f,\Omega_{jk}^{(\pm)}(e')g) = 
(E_k(e\cap e')f,g).
$$
\end{lemma}
Proof. By Lemma 4.2
$$
(\Omega_{jk}^{(\pm)}(e)f,\Omega_{jk}^{(\pm)}(e')g) = \int_{e\cap e'}(E_j'(\lambda)G_{jk}(\lambda \pm i0)f,G_{jk}(\lambda \pm i0)g)d\lambda.
$$
Using the resolvent equation, we have
\begin{equation}
\begin{split}
G_{jk}^{\ast}(\lambda \pm i\epsilon)\frac{1}{2\pi i}[R_j(\lambda + i\epsilon) - R_j(\lambda - i\epsilon)]G_{jk}(\lambda \pm i\epsilon) \\
 = \frac{1}{2\pi i}[R_k(\lambda + i\epsilon) - R_k(\lambda - i\epsilon)].
 \end{split}
 \label{eq:Chap2Sect4GjkastRj}
\end{equation}
Hence, 
\begin{equation}
\begin{split}
\left(\frac{1}{2\pi i}[R_j(\lambda + i\epsilon) - R_j(\lambda - i\epsilon)]G_{jk}(\lambda \pm i\epsilon)f,G_{jk}(\lambda \pm i\epsilon)g\right) \\
= 
\left(\frac{1}{2\pi i}[R_k(\lambda + i\epsilon) - R_k(\lambda - i\epsilon)]f,g\right).
\end{split}
\nonumber
\end{equation}
Letting $\epsilon \to 0$, we finally obtain
\begin{equation}
(E_j'(\lambda)G_{jk}(\lambda \pm i0)f,G_{jk}(\lambda \pm i0)g) = 
(E_k'(\lambda)f,g),
\label{eq:Chap2Sect4innerproductEGG}
\end{equation}
which proves the lemma. \qed

\bigskip
By Lemma 4.3, $\Omega_{jk}^{(\pm)}(e)$ is a partial isometry on  ${\mathcal H}$ with initial set $E_k(e){\mathcal H}$.


\begin{lemma} For any compact interval $e \subset I$, we have $(\Omega_{jk}^{(\pm)}(e))^{\ast} = \Omega_{kj}^{(\pm)}(e)$.
\end{lemma}
Proof. Since $G_{kj}^{\ast}(z)G_{jk}^{\ast}(z) = 1$, by multiplying (\ref{eq:Chap2Sect4GjkastRj}) by $G_{kj}^{\ast}(\lambda \pm i\epsilon)$, we have
\begin{equation}
\begin{split}
\frac{1}{2\pi i}[R_j(\lambda + i\epsilon) - R_j(\lambda - i\epsilon)]G_{jk}(\lambda \pm i\epsilon) \\
= G_{kj}^{\ast}(\lambda \pm i\epsilon)\frac{1}{2\pi i}[R_k(\lambda + i\epsilon) - R_k(\lambda - i\epsilon)]
\end{split}
\nonumber
\end{equation}
Letting $\epsilon \to 0$, we have for $f, g \in {\mathcal H}_+$
\begin{equation}
(f,E_j'(\lambda)G_{jk}(\lambda \pm i0)g) = (E_k'(\lambda)G_{kj}(\lambda \pm i0)f,g),
\label{eq:Chap2Sect4fEjprimelambdaGjk}
\end{equation}
which proves the lemma. \qed

\bigskip
This lemma implies that the final set of $\Omega_{jk}^{(\pm)}(e)$ is the initial set of $\Omega_{kj}^{(\pm)}(e)$, i.e. $\Omega_{jk}^{(\pm)}(e)$ is a partial isometry with initial set $E_k(e)\mathcal H$ and final set $E_j(e)\mathcal H$.


\begin{lemma}
For any compact intervals $e, e' \subset I$, we have
 $E_j(e')\Omega_{jk}^{(\pm)}(e) = \Omega_{jk}^{(\pm)}(e)E_k(e')$.
\end{lemma}
Proof. Lemma 4.2 yields for $f, g \in {\mathcal H}_+$
\begin{equation}
\begin{split}
(E_j(e')\Omega_{jk}^{(\pm)}(e)f,g) = (\Omega_{jk}^{(\pm)}(e)f,E_j(e')g) \\
= \int_{e\cap e'}(E_j'(\lambda)G_{jk}(\lambda \pm i0)f,g)d\lambda.
\end{split}
\nonumber
\end{equation}
By (\ref{eq:Chap2Sect4fEjprimelambdaGjk}) the right-hand side is equal to
\begin{eqnarray*}
\int_{e\cap e'}(f,E_k'(\lambda)G_{kj}(\lambda \pm i0)g)d\lambda 
&=& (f,E_k(e')\Omega_{kj}^{(\pm)}(e)g)\\
&=& (\Omega_{jk}^{(\pm)}(e)E_k(e')f,g),
\end{eqnarray*}
which proves the lemma. \qed

\bigskip
The assertion (3) of Theorem 4.1 is a direct consequence of the above lemma. Approximating $I$ by compact intervals, we define $\Omega_{jk}^{(\pm)}(I)$. 


\subsection{Time-dependent wave operators} 
 We consider the relation between stationary and time-dependent wave operators. We impose a new assumption.

\medskip
\noindent
(A-4) {\it For any open set $e \subset I$, there is a set $\mathcal D_e \subset {\mathcal H}_+\cap E_1(e)\mathcal H$, which is assumed to be dense in $E_1(e)\mathcal H$, such that for any $f \in \mathcal D_e$} 
$$
\int_{-\infty}^{\infty}\|(H_2 - H_1)e^{-itH_1}f\|dt < \infty.
$$


\begin{theorem} Under the assumptions (A-1) $\sim$ (A-4), for any open set $e \subset I$, the strong limit
$$
{\mathop{\rm s-lim}_{t\to\pm\infty}}\,e^{itH_2}e^{-itH_1}E_1(e) =: W_{21}^{(\pm)}(e)
$$
in $\mathcal H$ exists and $\Omega_{21}^{(\pm)}(e) = W_{21}^{(\pm)}(e)$.
\end{theorem}
Proof. The assumption (A-4) implies that, for $f \in {\mathcal D}_e$,
$$
\int_{-\infty}^{\infty}\|\frac{d}{dt}\left(e^{itH_2}e^{-itH_1}f\right)\|dt < \infty
$$
holds. Hence there exist the limits  ${\mathop{\rm s-lim}_{t\to\pm\infty}}\,e^{itH_2}e^{-itH_1}E_1(e)f$ and, therefore, by the density of $\mathcal D_e$ the existence of 
$W_{21}^{(\pm)}(e)$. 

To prove $\Omega_{21}^{(\pm)}(e) = W_{21}^{(\pm)}(e)$ for any open set $e \subset I$, it suffices, due to Lemma 4.3, to consider relatively compact sets $e$. 

Let $V_{21} = H_2 - H_1$. For $f \in \mathcal D_e$ we have
\begin{equation}
W_{21}^{(+)}(e)f 
= f + i\int_0^{\infty}e^{itH_2}V_{21}e^{-itH_1}fdt.
\nonumber
\end{equation}
Hence, for $f\in \mathcal D_e, g \in \mathcal H$,
\begin{equation}
(W_{21}^{(+)}(e)f,g) = (f,g) + \lim_{\epsilon\to0}i
\int_0^{\infty}(e^{itH_2}V_{21}e^{-itH_1}f,e^{-2\epsilon t}g)dt.
\label{C2S4W21intepsilon}
\end{equation}
Using the following relations
\begin{equation}
R(\lambda + i\epsilon) = i\int_0^{\infty}e^{it(\lambda + i\epsilon - H)}dt,
\quad
R(\lambda - i\epsilon) = - i\int_{-\infty}^{0}e^{it(\lambda - i\epsilon - H)}dt
\nonumber
\end{equation}
and  Plancherel's formula for the Fourier transform, we have for $f \in \mathcal D_e$ and $g \in \mathcal H$
\begin{equation}
\begin{split}
&  i\int_0^{\infty}(e^{itH_2}V_{21}e^{-itH_1}f,e^{-2\epsilon t}g)dt \\
& = - \frac{1}{2\pi i}\int_{-\infty}^{\infty}(V_{21}R_1(\lambda + i\epsilon)f,
R_2(\lambda + i\epsilon)g)d\lambda.
\end{split} 
\label{C2S4int=int}
\end{equation}
Here we should note that $\|V_{21}R(\cdot + i\epsilon)f\|_{\mathcal H},  
\|R_2(\cdot + i\epsilon)g\|_{\mathcal H} \in L^2({\bf R})$, hence the integral of the right-hand side is absolutely convergent. To see this, we have only to note that
$$
\|R_j(\lambda + i\epsilon)h\|^2 = \int_{-\infty}^{\infty}
\frac{1}{(\mu-\lambda)^2+\epsilon^2}d_{\mu}(E_j(\mu)h,h),
$$
$$
V_{21}R_1(\lambda + i\epsilon)f = V_{21}(H_1+i)^{-1}R_1(\lambda + i\epsilon)
(H_1 + i)f,
$$
and $(H_1 + i)f \in E_1(e)\mathcal H$, also $V_{21}(H_1 + i)^{-1} \in {\bf B}(\mathcal H;\mathcal H)$ by (\ref{C2S4V21Resolbnded}).

We now let
$$
\delta_2(\lambda,\epsilon) = \frac{1}{2\pi i}\big(R_2(\lambda + i\epsilon) - R_2(\lambda - i\epsilon)\big),
$$ 
and prove that, if $f \in \mathcal D_e$ and $g$ is such that $d_{\mu}(E_2(\mu)g,g)$ is compactly supported,
\begin{equation}
  i\int_0^{\infty}(e^{itH_2}V_{21}e^{-itH_1}f,e^{-2\epsilon t}g)dt 
= \lim_{N\to\infty}\int_{-N}^{N}(\delta_2(\lambda,\epsilon)
V_{21}R_1(\lambda + i\epsilon)f,g)d\lambda.
\label{C2S4int=int2}
\end{equation}

Indeed, by using the identity $R_2(z) - R_1(z) = - R_2(z)V_{21}R_1(z)$, we have
\begin{equation}
\begin{split}
& - \frac{1}{2\pi i}\int_{-N}^{N}(V_{21}R_1(\lambda + i\epsilon)f,R_2(\lambda + i\epsilon)g)d\lambda  \\
&  = \int_{-N}^{N}(\delta_2(\lambda,\epsilon)V_{21}R_1(\lambda + i\epsilon)f,g)d\lambda \\
& + \frac{1}{2\pi i}\int_{-N}^N((R_2(\lambda + i\epsilon) - R_1(\lambda + i\epsilon))f,g)d\lambda.
\nonumber
\end{split}
\end{equation}
However, $\displaystyle{\frac{1}{2\pi i}\int_{-N}^N(R_j(\lambda +i\epsilon)f,g)d\lambda \to \frac{1}{2}(f,g)}$ when $N \to \infty$. In fact, 
\begin{equation}
(R_j(z)f,g) = \int_{-\infty}^{\infty}\frac{1}{\mu -z}d_{\mu}(E_j(\mu)f,E_j(\mu)g),
\nonumber
\end{equation} 
where the domain of integration is bounded by our assumptions on $f$ and $g$.
Therefore
$$
\frac{1}{2\pi}\int_{-N}^N(R_j(\lambda +i\epsilon)f,g)d\lambda = 
\frac{1}{2\pi}\int_{-\infty}^{\infty}
\ln\left(\frac{-N-\mu+i\epsilon}{N-\mu+i\epsilon}\right)d_{\mu}(E_j(\mu)f,E_j(\mu)g).
$$
Since
$\displaystyle{
\ln\left(\frac{-N-\mu+i\epsilon}{N-\mu+i\epsilon}\right) \to \pi i}$ for any $\mu$, the result follows by Lebesgue's dominated convergence theorem.

Let us take bounded open intervals $J, J_1$ such that
\begin{equation}
e \subset \overline{e} \subset J \subset \overline{J} \subset J_1 \subset 
\overline{J_1} \subset I,
\label{C2S4eandJ}
\end{equation}
and $g = \varphi(H_2)h$ for some $\varphi(\lambda) \in C_0^{\infty}(J)$ and $h \in \mathcal H_+$. Such $g$'s are dense in $E_2(I)\mathcal H$. Then we have
\begin{equation}
(\delta_2(\cdot,\epsilon)V_{21}R_1(\cdot + i\epsilon)f,g) \in L^1({\bf R}), \quad \epsilon > 0,
\label{C2S4innerprodcutL1}
\end{equation}
\begin{equation}
\lim_{\epsilon\to0}\int_{-\infty}^{\infty}(\delta_2(\lambda,\epsilon)V_{21}R_1(\lambda + i\epsilon)f,g)d\lambda = (\Omega_{21}^{(+)}(e)f,g) - (f,g).
\label{C2S4W=Omega}
\end{equation}
 
 In fact, since $V_{21}R_1(\lambda + i\epsilon) = G_{21}(\lambda + i\epsilon) - 1$, we have
$$
(\delta_2(\lambda,\epsilon)V_{21}R_1(\lambda + i\epsilon)f,g) = 
(\delta_2(\lambda,\epsilon)G_{21}(\lambda + i\epsilon)f,g) - 
(f,\delta_2(\lambda,\epsilon)g).
$$
Then the 2nd term of the right-hand side is written as
\begin{equation}
(f,\delta_2(\lambda,\epsilon)g) = \frac{\epsilon}{\pi}\int_{-\infty}^{\infty}
\frac{1}{(\mu-\lambda)^2+\epsilon^2}d_{\mu}(f,E_2(\mu)g).
\label{C2S4fdeltalambdag}
\end{equation}
If $\lambda \not\in \overline{J_1}$, the right-hand side is dominated from above by $C\epsilon(1 + |\lambda|^2)^{-1}$. On the other hand, assumptions (A-1), (A-2) imply that the left-hand side is bounded for $\lambda \in \overline{J_1}$ uniformly with respect to $\epsilon$. Therefore $(f,\delta_2(\cdot,\epsilon)g) \in L^1({\bf R})$, and by Stone's theorem
\begin{equation}
\lim_{\epsilon\to0}\int_{-\infty}^{\infty}(f,\delta_2(\lambda,\epsilon)g)d\lambda = (f,g).
\label{C2S4(f,g)}
\end{equation}
By the resolvent equation, $R_1(z) = R_1(i)(1 + (z-i)R_1(z))$. Then we have
$$
G_{21}(z) = (H_2 - H_1)R_1(i)\left(1 + (z - i)R_1(z)\right) + 1.
$$
Since $f \in E_1(e)\mathcal H$, we have $\|(\lambda + i\epsilon - i)R_1(\lambda + i\epsilon)f\| \leq C_{f}$ uniformly for $\lambda \not\in \overline{J_1}$ and $\epsilon > 0$. Hence so is $\|G_{21}(\lambda + i\epsilon)f\|$.
Then formula (\ref{C2S4fdeltalambdag})
 implies that if $\lambda \not\in \overline{J_1}$,
$$
|(G_{21}(\lambda + i\epsilon)f,\delta_2(\lambda,\epsilon)g)| \leq 
C\epsilon (1 + |\lambda|^2)^{-1},
$$
which implies
\begin{equation}
\int_{{\bf R}\setminus\overline{J_1}}(\delta_2(\lambda,\epsilon)G_{21}(\lambda + i\epsilon)f,g)d\lambda \to 0, \quad 
\epsilon \to 0.
\label{C2S4J1cto0}
\end{equation}
Since $f \in E_1(e) H$, 
$$
\int_{J_1} E_2'(\lambda) G_{21}(\lambda+i\epsilon) f d\lambda 
\rightarrow \Omega^{(+)}_{21}(e) f.
$$
Together with (\ref{C2S4J1cto0}), this implies that
\begin{equation}
\label{C2S4E2lambdaG21lambdaint}
\int_{{\bf R}} \left(\delta_2(\lambda, \epsilon)G_{21}(\lambda+i\epsilon)f, g \right)d\lambda
\rightarrow \left(\Omega^{(+)}_{21}(e) f, g  \right)
\end{equation}
Equations (\ref{C2S4(f,g)}) and (\ref{C2S4E2lambdaG21lambdaint}) prove (\ref{C2S4W=Omega}).
By (\ref{C2S4W21intepsilon}), (\ref{C2S4int=int2}) and (\ref{C2S4W=Omega}) we get 
$W_{21}^{(+)}(e)= \Omega_{21}^{(+)}$ when $e$ is a relatively compact interval in $I$.



For an open subset $e \subset I$, we have only to appriximate $e$ by a finite number of relatively compact intervals. The proof for $W_{21}^{(-)}(e) = \Omega_{21}^{(-)}(e)$ is similar. \qed


\subsection{Spectral representation}
Let us recall that for a self-adjoint operator ${H = \int_{-\infty}^{\infty}\lambda dE(\lambda)}$, we take an open interval $I$ in $\sigma_{ac}(H)$. We take an auxiliary Hilbert space $\bf h$ and a measure $\rho(\lambda)d\lambda$ on $I$, 
$\rho(\lambda) \in L^1(I;d\lambda)$, and put
$$
\widehat{\mathcal H}(I) = L^2(I;{\bf h};\rho(\lambda)d\lambda).
$$
A unitary operator $U$ from $E(I){\mathcal H}$ onto $\widehat{\mathcal H}(I)$ satisfying
$$
(UHf)(\lambda) = \lambda(Uf)(\lambda), \quad \lambda \in I, \quad f \in D(H)
$$
is called a spectral representaion of $H$ on $I$. By the functional calculus, 
\begin{equation}
(U\varphi(H)f)(\lambda) = \varphi(\lambda)(Uf)(\lambda)
\label{C2S4FuncCalcu}
\end{equation}
holds for any bounded Borel function $\varphi$ and $f \in E(I)\mathcal H$. In fact, (\ref{C2S4FuncCalcu}) is first proven for the resolvent $\varphi(H) = (H - z)^{-1}$, next for the spectral measure $E(\mu)$ by using Stone's formula, and then for  any bounded Borel function.

Let ${\mathcal H}_+, {\mathcal H}_-$ be Banach spaces satisfying the assumptions in subsection 4.2. We assume that, for $\lambda \in I$, there exists a bounded operator $U(\lambda) \in {\bf B}({\mathcal H}_+;{\bf h})$, which is stronlgy continuous in $\lambda$,  such that
\begin{equation}
(Uf)(\lambda) = U(\lambda)f, \quad \lambda \in I, \quad f \in {\mathcal H}_+.
\nonumber
\end{equation}
Then $U(\lambda)^{\ast} \in {\bf B}({\bf h};{\mathcal H}_-)$. 
Let us show that for $\Phi \in \widehat H(I)$
\begin{equation}
U^{\ast}\Phi = \int_I U(\lambda)^{\ast}\Phi(\lambda)\rho(\lambda)d\lambda \in E(I){\mathcal H}.
\label{eq:Chap2Sect4intIUlambdaast}
\end{equation}
Indeed, let us first assume that ${\rm supp}\,\Phi \subset J$, where $\overline{J}$ is a compact set in $I$.
Then, for $f \in {\mathcal H}_+$, we have
\begin{eqnarray*}
\Big(\int_I U(\lambda)^{\ast}\Phi(\lambda)\rho(\lambda)d\lambda,f\Big) =
\int_I(\Phi(\lambda),U(\lambda)f)_{\bf h}\rho(\lambda)d\lambda
= (\Phi,Uf)_{\widehat{\mathcal H}} 
= (U^{\ast}\Phi,f).
\end{eqnarray*}
As $U^{\ast}$ is partial isometry, the right-hand side can be extended to 
$f \in {\mathcal H}$,
which together with Riesz' theorem implies (\ref{eq:Chap2Sect4intIUlambdaast}) 
for $\Phi$ with ${\rm supp}\,\Phi \subset J$. 
Since $J$ is arbitrary, and $I \subset \sigma_{ac}(H)$, (\ref{eq:Chap2Sect4intIUlambdaast}) is exteded onto $\widehat{\mathcal H}(I)$.

As a consequence, we have the inversion formula for $f \in E(I){\mathcal H}$ 
\begin{equation}
f = \int_IU(\lambda)^{\ast}(Uf)(\lambda)\rho(\lambda)d\lambda.
\label{eq:Chap2Sect4InversionFormula}
\end{equation}
In fact, for $g \in {\mathcal H}_+$, 
$$
(f,g)_{\mathcal H} = (Uf,Ug)_{\widehat{\mathcal H}} = \int_I((Uf)(\lambda),U(\lambda)g)_{\bf h}\rho(\lambda)d\lambda.
$$
Hence we have
$$
(f,g)_{\mathcal H} = \int_I(U(\lambda)^{\ast}(Uf)(\lambda),g)\rho(\lambda)d\lambda,
$$
which proves ({\ref{eq:Chap2Sect4InversionFormula}) by virtue of (\ref{eq:Chap2Sect4intIUlambdaast}). 

We need a new assumption:

\medskip
\noindent
(A-5) {\it There exists a  subspace $\mathcal D \subset D(H)\cap \mathcal H_+$ such that ${\mathcal D}$ as well as $H{\mathcal D}$ are dense in $\mathcal H_+$ and $D(H)$. }

\medskip
Then, for $\psi \in {\bf h}, f \in {\mathcal D}$,
$$
(U(\lambda)^{\ast}\psi,(H - \lambda)f) = 0
$$
holds, since  $U(\lambda)Hf = \lambda U(\lambda)f$. Therefore, $U(\lambda)^{\ast} \in {\bf B}(\bf h;{\mathcal H}_-)$ satisfies the equation
\begin{equation}
(H - \lambda)U(\lambda)^{\ast} = 0,
\nonumber
\end{equation}
and is called the {\it eigenoperator} of $H$. Here the self-adjoint operator $H$ in $\mathcal H$ is extended to $\mathcal H_-$ via the relation
\begin{equation}
(Hu,f) = (u,Hf), \quad u \in \mathcal H_-, \quad f \in \mathcal D.
\label{C2S4HinmathcalH-}
\end{equation}

We now discuss the perturbation theory for spectral representations. For $H_1$ we assume that

\bigskip
\noindent
(A-6)$\ $ {\it For any $\lambda \in I$ there exists $U_1(\lambda) \in {\bf B}({\mathcal H}_+;{\bf h})$ such that for $f, g \in {\mathcal H}_+$ 
\begin{equation}
(E_1'(\lambda)f,g) = \rho(\lambda)(U_1(\lambda)f,U_1(\lambda)g)_{\bf h}.
\nonumber
\end{equation}
Moreover, $U_{1}$ defined by $(U_1f)(\lambda) = U_1(\lambda)f$ is uniquely extended to a unitary operator from $E_1(I){\mathcal H}$ to $\widehat H(I)$. }

\bigskip
By this assumption, we have for $f \in D(H_1)$
\begin{equation}
(U_1H_1f)(\lambda) = \lambda(U_1f)(\lambda) a.e..
\label{eq:Chap2Sect4F1H1fequallambdaF}
\end{equation}
In fact, let $f \in {\mathcal D}$. Since $R_1(z)H_1 = 1 + zR_1(z)$, we have $E_1'(\lambda)H_1f = \lambda E_1'(\lambda)f$. The assumption (A-6) then implies
\begin{eqnarray*}
(E_1(I)H_1f,g) &=& \int_I\lambda\left((U_1f)(\lambda),(U_1g)(\lambda)\right)_{\bf h}\rho(\lambda)d\lambda \\
&=& \int_I\left((U_1H_1f)(\lambda),(U_1g)(\lambda))_{\bf h}\rho(\lambda\right)d\lambda,
\end{eqnarray*}
which proves (\ref{eq:Chap2Sect4F1H1fequallambdaF}) for $f \in {\mathcal D}$.  Since
${\mathcal D}$ is dense in $D(H_1)$ we obtain (\ref{eq:Chap2Sect4F1H1fequallambdaF}).
Therefore, $U_1(\lambda)^{\ast} \in {\bf B}({\bf h};{\mathcal H}_-)$ is an eigenoperator of $H_1$:
\begin{equation}
(H_1 - \lambda)U_1(\lambda)^{\ast} = 0.
\nonumber
\end{equation}

We construct the spectral representation of $H_2$ by using that of $H_1$.
Starting from $U_1$, we will construct two operators $U_2^{(+)}, U_2^{(-)}$ corresponding to wave operators $W_{21}^{(+)}, W_{21}^{(-)}$. For $\lambda \in I$, we define
\begin{equation}
U_2^{(\pm)}(\lambda) = U_1(\lambda)G_{12}(\lambda \pm i0).
\nonumber
\end{equation}
For  $f \in {\mathcal H}_+$, we put $(U_2^{(\pm)}f)(\lambda) = U_2^{(\pm)}(\lambda)f$. Then we have the following theorem.


\begin{theorem} Under the assumptions (A-1) $\sim$ (A-6), we have
\begin{equation}
(E_2'(\lambda)f,g) = \rho(\lambda)(U_2^{(\pm)}(\lambda)f,U_2^{(\pm)}(\lambda)g)_{\bf h}, \quad f, g \in {\mathcal H}_+.
\nonumber
\end{equation}
Moreover $U_2^{(\pm)} = U_1\left(W_{21}^{(\pm)}(I)\right)^{\ast}$, and $U_2^{(\pm)}$ is a spectral representation for $H_2$. 
\end{theorem}
\noindent
Proof. The first half of the theorem follows from (\ref{eq:Chap2Sect4innerproductEGG}) and (A-6). By virtue of (\ref{eq:Chap2Sect4fEjprimelambdaGjk}) and (A-6) we have
\begin{eqnarray*}
(E_2'(\lambda)G_{21}(\lambda \pm i0)f,g) &=& 
(f,E_1'(\lambda)G_{12}(\lambda \pm i0)g) \\
&=& \rho(\lambda)(U_1(\lambda)f,U_2^{(\pm)}(\lambda)g)_{\bf h}.
\end{eqnarray*}
Integration with respect to $\lambda$ then yields, in view of Theorem 4.6, that
$$
(W_{21}^{(\pm)}(I)f,g) = (U_1f,U_2^{(\pm)}g)_{\widehat{\mathcal H}},
$$
hence $W_{21}^{(\pm)}(I) = (U_2^{(\pm)})^{\ast}U_1$. We have, therefore, $U_2^{(\pm)} = U_1\left(W_{21}^{(\pm)}(I)\right)^{\ast}$. Since ${\rm Ran}\, W_{21}^{(\pm)}(I) = E_2(I){\mathcal H}$ and $W_{21}^{(\pm)}\varphi(H_1) = \varphi(H_2)W_{21}^{(\pm)}$ for any bounded Borel function $\varphi(\lambda)$, $U_2^{(\pm)}$ is a partial isometry with initial set $E_2(I)\mathcal H$ and final set $\widehat{\mathcal H}(I)$. Moreover $U_2^{(\pm)}\varphi(H_2) = \varphi(\lambda)U_2^{(\pm)}$ for any bounded Borel function. Therefore, $U_2^{(\pm)}$ is a spectral representation for $H_2$. \qed

\bigskip
By the relation $U_2^{(\pm)}(\lambda)^{\ast} = (1 - R_2(\lambda \mp i0)V)U_1(\lambda)^{\ast}$, $V = H_2 - H_1$, we have
\begin{equation}
(H_2 - \lambda)U_2^{(\pm)}(\lambda)^{\ast} = 0.
\nonumber
\end{equation}
Hence $U_2^{(\pm)}(\lambda)^{\ast}$ is an eigenoperator of $H_2$.
Let us summarize the results obtained so far.
Let $E_2(\lambda)$ be the resolution of identity for $H_2$.


\begin{theorem}
(1) Let $V_{21} = H_2 - H_1$ and put 
$$
U_2^{(\pm)}(\lambda) = U_1(\lambda)(1 - V_{21}R(\lambda \pm i0))=
U_1(\lambda) G_{12}(\lambda \pm i0).
$$
Then $U_2^{(\pm)}(\lambda) \in {\bf B}({\mathcal H}_+ ; {\bf h})$ for $\lambda  \in I$. \\
\noindent
(2) $U_2^{(\pm)}(\lambda)^{\ast} \in {\bf B}({\bf h};\mathcal H_-)$ is an eigenoperator of $H$ with eigenvalue $\lambda \in I$ in the following sense 
$$
((H_2 - \lambda)f,U_2^{(\pm)}(\lambda)^{\ast}\varphi) = 0
$$
for any $f \in {\mathcal H}_+$ such that $H_2f \in {\mathcal H}_+$ and $\varphi \in {\bf h}$. Moreover, 
$$
(U_2^{(\pm)}H_2f)(\lambda) = \lambda(U_2^{(\pm)}f)(\lambda), \quad 
f \in D(H_2), \quad \lambda \in I.
$$ 
\noindent
\noindent
(3) The operator $U_2^{(\pm)}$ defined by $(U_2^{(\pm)}f)(\lambda) = U_2^{(\pm)}(\lambda)f$ for $f \in {\mathcal H}_+$ is uniquely extended to a partial isometry with 
the initial set $E_2(I)\mathcal H$ and the final set $\widehat{\mathcal H}(I)$. \\
\noindent
(4) For any $\Phi \in \widehat{\mathcal H}(I)$ and any compact interval $e \subset I$,
$$
\int_e U_2^{(\pm)}(\lambda)^{\ast}\Phi(\lambda)\rho(\lambda)d\lambda \in {\mathcal H}.
$$
(5) For any $f \in E_2(I)\mathcal H$, the following inversion formula holds:
$$
f = \mathop{\rm s-lim}_{n\to\infty}\int_{I_n}U_2^{(\pm)}(\lambda)^{\ast}(U_2^{(\pm)}f)(\lambda)\rho(\lambda)d\lambda,
$$
where $I_n = [a_n,b_n], a < a_n < b_n < b$, $a_n \to a, b_n \to b$ and $I = (a,b)$.
\end{theorem}
Proof. We have only to show the assertions (4) and (5). Let $I_e(\Phi)$ be the integral in (4). 
We first assume that supp$\,\Phi(\lambda)$ is a compact set $e$ in $I$. We take $f \in E_2(I)\mathcal H$ such that $U_2^{(\pm)}f = \Phi$. Then for any $g \in {\mathcal H}_+$ , we have
\begin{eqnarray*}
(f,g) &=& (U_2^{(\pm)}f,U_2^{(\pm)}g) \\
  &=& \int_e((U_2^{(\pm)}f)(\lambda),(U_2^{(\pm)}g)(\lambda))\rho(\lambda)d\lambda 
  =\int_e \left( \Phi(\lambda), (U_2^{(\pm)} g)(\lambda) \right) \rho(\lambda)d\lambda
  \\
 &=& \int_e(U_2^{(\pm)}(\lambda)^{\ast}\Phi(\lambda),g)\rho(\lambda)d\lambda 
 = (I_e(\Phi),g).
\end{eqnarray*}
We have, therefore, $I_e(\Phi) = f \in \mathcal H$. 
This implies also that, for any $f \in E_2(I)\mathcal H$ and a compact interval $e \subset I$,
$$
E_2(e)f = \int_eU_2^{(\pm)}(\lambda)^{\ast}(U_2^{(\pm)}f)(\lambda)\rho(\lambda)d\lambda,
$$
since $(U_2^{(\pm)}E_2(e)f)(\lambda) = \chi_e(\lambda)(U_2^{(\pm)}f)(\lambda)$, where $\chi_e(\lambda)$ is the characteristic function of $e$. Therefore 
$$
\left\|\int_eU_2^{(\pm)}(\lambda)^{\ast}(U_2^{(\pm)}f)(\lambda)\rho(\lambda)d\lambda\right\| \to 0
$$
if the measure of $e$ tends to 0. This proves (5). \qed


\subsection{S-matrix}
The scattering operator for $H_1, H_2$  (on $I$) is defined by
\begin{equation}
S = (W_{21}^{(+)}(I))^{\ast}W_{21}^{(-)}(I).
\nonumber
\end{equation}
This is unitary on $E_1(I){\mathcal H}$. Let us rewrite it by using the spectral representation. We define
\begin{equation}
\widehat S = U_1SU_{1}^{\ast}.
\nonumber
\end{equation}
Letting $V_{21} = H_2 - H_1$, we also put
\begin{equation}
\widehat S(\lambda) = 1 - 2\pi i\rho^2(\lambda)A(\lambda), 
\nonumber
\end{equation}
\begin{equation}
A(\lambda) = U_1(\lambda)V_{21}U_1(\lambda)^{\ast}
- U_1(\lambda)V_{21}R_2(\lambda + i0)V_{21}U_1(\lambda)^{\ast}.
\nonumber
\end{equation}
Then  $\widehat S(\lambda) \in {\bf B}(\bf h;\bf h)$ and is called the S-matrix or the scattering matrix.

\begin{theorem}
$\widehat S(\lambda)$ is unitary on $\bf h$, and for any $\hat f \in \widehat{\mathcal H}$ 
\begin{equation}
(\widehat S \hat f)(\lambda) = \widehat S(\lambda)\hat f(\lambda)
\nonumber
\end{equation}
holds. Here the right-hand side means that we fix $\lambda$ arbitrarily, regard $\hat f(\lambda)$ as an element of $\bf h$ and apply $\widehat S(\lambda)$. 
\end{theorem}
Proof. Noting that
$$
W_{21}^{(\pm)}(I) = E_1(I) + i\int_0^{\pm\infty}e^{isH_2}V_{21}e^{-isH_1}E_1(I)ds,
$$
we have
$$
W_{21}^{(+)}(I) - W_{21}^{(-)}(I) = i\int_{-\infty}^{\infty}
e^{itH_2}V_{21}e^{-itH_1}E_1(I)dt.
$$
By the definition of $S$, we have
$$
(S - 1)E_1(I) = (W_{21}^{(+)})^{\ast}(W_{21}^{(-)}(I) - 
W_{21}^{(+)}(I)).
$$
Letting $f = E_1(I)f, g = E_1(I)g$, we then have
\begin{equation}
\begin{split}
&  (Sf,g) - (f,g) \\
&=  - i\int_{-\infty}^{\infty}(e^{itH_2}V_{21}e^{-itH_1}f,W_{21}^{(+)}(I)g)dt \\
&= - i\int_{-\infty}^{\infty}(V_{21}e^{-itH_1}f,e^{-itH_1}g)dt \\
 & - \int_0^{\infty}ds\int_{-\infty}^{\infty}(V_{21}e^{-itH_1}f,
e^{isH_2}V_{21}e^{-i(s+t)H_1}g)dt,
\end{split}
\label{eq:Chap2Sect4FormulaofSminus1}
\end{equation}
where we have used $e^{-itH_2}W_{21}^{(+)}(I) = W_{21}^{(+)}(I)e^{-itH_1}$. Letting $\hat f(\lambda) = U_1(\lambda)f, \hat g(\lambda) = U_1(\lambda)g$, we obtain
\begin{eqnarray*}
& &\int_{-\infty}^{\infty}(V_{21}e^{-isH_2}V_{21}e^{-itH_1}f,e^{-i(s+t)H_1}g)dt \\
&=& \int_{-\infty}^{\infty}dt\int_I(U_1(\lambda)V_{21}e^{-isH_2}V_{21}e^{-itH_1}f,e^{-i(s+t)\lambda}\hat g(\lambda))_{\bf h}\rho(\lambda)d\lambda.
\end{eqnarray*}
Inserting $e^{-\epsilon|t|}$ and letting $\epsilon \to 0$, this converges to
\begin{eqnarray*}
& &2\pi\int_I(U_1(\lambda)V_{21}e^{-is(H_2 - \lambda)}V_{21}E_1'(\lambda)f,\hat g(\lambda))_{\bf h}\rho(\lambda)d\lambda\\
&=&2\pi\int_I(U_1(\lambda)V_{21}e^{-is(H_2 - \lambda)}V_{21}U_1(\lambda)^{\ast}\hat f(\lambda),\hat g(\lambda))_{\bf h}\rho(\lambda)^2d\lambda,
\end{eqnarray*}
where we have used $E_1'(\lambda) = \rho(\lambda)U_1(\lambda)^{\ast}U_1(\lambda)$. Therefore, 
the last term of the most right-hand side of  (\ref{eq:Chap2Sect4FormulaofSminus1}) is equal to
$$
- 2\pi\int_0^{\infty}ds\int_I(U_1(\lambda)V_{21}e^{-is(H_2-\lambda)}
V_{21}U_1(\lambda)^{\ast}\hat f(\lambda),\hat g(\lambda))_{\bf h}
\rho(\lambda)^2d\lambda.
$$
Inserting $e^{-\epsilon s}$ and letting $\epsilon \to 0$, this converges to
$$
2\pi i\int_I(U_1(\lambda)V_{21}R_2(\lambda + i0)
V_{21}U_1(\lambda)^{\ast}\hat f(\lambda),\hat g(\lambda))_{\bf h}
\rho(\lambda)^2d\lambda.
$$
Similarly the first term of the most right-hand side of (\ref{eq:Chap2Sect4FormulaofSminus1}) is rewritten as
$$
- 2\pi i\int_I(U_1(\lambda)V_{21}U_1(\lambda)^{\ast}\hat f(\lambda),\hat g(\lambda))_{\bf h}
\rho(\lambda)^2d\lambda.
$$
This proves the representation of $\widehat S$. Since $\widehat S$ is unitary on $\widehat{\mathcal H}$, so is $\widehat S(\lambda)$ on $\bf h$. \qed


\begin{lemma}
For any $\lambda \in I$, we have
\begin{equation}
U_2^{(+)}(\lambda) = \widehat S(\lambda)U_2^{(-)}(\lambda).
\nonumber
\end{equation}
\end{lemma}
Proof. By Theorem 4.7, we have
$$
\left(W_{12}^{(+)}(I)\right)^{\ast} = \left(U_1\right)^{\ast}U_2^{(+)}, \quad
W_{12}^{(-)} = \left(U_2^{(-)}\right)^{\ast}U_1.
$$
Therefore by the definition of $\widehat S$, we have
$$
\widehat S U_2^{(-)} = U_2^{(+)},
$$
which proves the lemma. \qed


\section{Examples of spectral representations}


\subsection{Spectral representation on ${\bf R}^n$} 
Let us apply the results in the previous section to Schr{\"o}dinger operators $H_0 = - \Delta$ and 
$$
H = - \sum_{i,j=1}^na_{ij}(x)\partial_i\partial_j + \sum_{i=1}^na_i(x)\partial_i + a_0(x)
$$ 
on ${\bf R}^n$, where $\partial_i = \partial/\partial x_i$. Let $\mathcal H = L^2({\bf R}^n;dx)$ and assume that $H$ is formally self-adjoint and uniformly elliptic on ${\bf R}^n$, i.e. there exists a constant $C_0 > 0$ such that
$$
C^{-1}|\xi|^2 \leq \sum_{i,j=1}^na_{ij}(x)\xi_i\xi_j \leq C|\xi|^2, \quad 
\forall x, \xi \in {\bf R}^n.
$$
The coefficients $a_{ij}(x) - \delta_{ij}$ and $a_i(x)$ of $H$ are assumed to be smooth and satisfy
\begin{equation}
\left|\partial^{\alpha}a(x)\right| \leq C_{\alpha}(1 + |x|)^{-1-\epsilon-|\alpha|}, \quad \forall \alpha, \quad \forall x \in {\bf R}^n
\nonumber
\end{equation}
for a constant $\epsilon > 0$. For $s \in {\bf R}$ we define the space $L^{2,s}$ by
\begin{equation}
L^{2,s} \ni f \Longleftrightarrow \|f\|_s^2 = \int_{{\bf R}^n}(1 + |x|^2)^s|f(x)|^2dx < \infty.
\nonumber
\end{equation}
Let $s > 1/2$ be arbitraily fixed. Then, by choosing $\mathcal H_{\pm} = L^{2,\pm s}$, the assumptions (A-1) $\sim$ (A-3) are satisfied for $H_1 = H_0$, $H_2 = H$ and $I = (0,\infty)$. We should remark that by this choice of $\mathcal H_{\pm}$, the boundary value of the resolvent $R_{j}(\lambda \pm i0)f$ is strongly continuous in $L^{2,-s}$ as a function of $\lambda > 0$. These facts are well-known and are proved in e.g. \cite{Is04a}, where they are proved for the potential perturbation of $- \Delta$, however, the proof also works for the case of the 2nd order variable coefficients. Let us also note that Theorem 3.1 can also be applied in this case.

As a spectral representation for $H_0$, we employ the usual Fourier transformation:
\begin{equation}
\left(U_0(\lambda)f\right)(\omega) = (2\pi)^{-1/2}\int_{{\bf R}^n}e^{-i\sqrt\lambda\omega\cdot x}f(x)dx,
\nonumber
\end{equation}
and ${\bf h} = L^2(S^{n-1})$ and $\rho(\lambda) = \frac{1}{2}\lambda^{(n-2)/2}$. Then the assumption (A-4) is also satisfied. Let $R(z) = (H - z)^{-1}$ and $V = H - H_0$. Then 
\begin{equation}
U_{\pm}(\lambda) = U_0(\lambda)(1 - VR(\lambda \pm i0))
\nonumber
\end{equation}
gives the spectral representation for $H$.


\subsection{Spectral representations on ${\bf H}^n$}
 Let $\mathcal H = L^2({\bf R}^n_+ ;dxdy/y^n)$ and consider the operators $H_0$ and $H$ introduced in \S 2.
Let $L^{2,s}$ be defined by Definition 2.6 of Chap. 1. Let ${\mathcal H}_{\pm} = L^{2,\pm s}$ with $1/2 < s < (1+\epsilon)/2$ and $H_1 = H_0$, $H_2 = H$ and $I = (0,\infty)$. First we check (A-1). Let $\langle \log y\rangle^s = (1 + |\log y|^2)^{s/2}$. We show that there exists a constant $C_s$ independent of $z \not\in {\bf R}$ such that
\begin{equation}
\|\langle\log y\rangle^s(H_j-z)^{-1}\langle \log y\rangle^{-s}\| 
\leq C_s|{\rm Im}\, z|^{-2}(1+|z|).
\label{C2S5commuteHlogy}
\end{equation}
Once we have proven (\ref{C2S5commuteHlogy}), we can use an abstract theorem from functional analysis (see Lemma 3.1 in Chap. 3, where $\sigma$ can be an arbitrary negative number) to show 
$$
\langle\log y\rangle^s\varphi(H_j)\langle \log y\rangle^{-s} \in {\bf B}(\mathcal H;\mathcal H), \quad \forall \varphi \in C_0^{\infty}({\bf R}),
$$ 
which yields (A-1).

Let us prove (\ref{C2S5commuteHlogy}). We have
\begin{equation}
\begin{split}
& \langle\log y\rangle^s(H_j-z)^{-1}\langle\log y\rangle^{-s} \\
& =
(H_j-z)^{-1} + (H_j-z)^{-1}[H_j,\langle\log y\rangle^s](H_j-z)^{-1}\langle \log y\rangle^{-s}.
\end{split}
\nonumber
\end{equation}
Since $[H_j,\langle\log y\rangle^s]$ is a 1st order differential operator with respect to $D_x$, $D_y$ with bounded coefficients, one can show
$$
\|[H_j,\langle\log y\rangle^s](H_j-z)^{-1}\| \leq C_s|{\rm Im}\,z|^{-1}(1 + |z|)
$$
by using Theorem 1.3 (4) and the standard estimate of the resolvent. The inequality (\ref{C2S5commuteHlogy}) imediately follows from this.

 Theorem 2.3 
together with Lemma 1.2.7 justify  (A-2). As above, by this choice of $L^{2,\pm s}$ the strong continuity of $R_j(\lambda \pm i0)f$ with respect to $\lambda$ is guaranteed.

To prove (A-4) for a proper ${\mathcal D}_e$, $e=(a^2, b^2),\, 0<a<b< \infty$, we first 
observe that it is sufficient to show that, for $1 < s<1+\epsilon$ and $f \in {\mathcal D}_e$,
$$
\int_{-\infty}^{\infty}\Big(||e^{-itH_0} f||_{\chi^{-s}}+
\sum_j||D_j e{^{-itH_0}} f ||_{\chi^{-s}}+\sum_{j,l}||D_j D_l e^{-itH_0} f ||_{\chi^{-s}} \Big) dt < \infty.
$$
Assuming that $H_0\mathcal D_e \subset \mathcal D_e$, and utilising Theorem 1.3 (6), we can confine to the proof that
$$
\int_{-\infty}^{\infty}\left(||e^{-itH_0} f||_{\chi^{-s}}+
||e^{-itH_0} H_0 f ||_{\chi^{-s}}\right) dt < \infty,
\quad f \in {\mathcal D}_e.
$$
Let
$$
{\mathcal D}_e= \left\{f:\, \phi(k, \xi)=(F_0 \mathcal F_0^{(+)}f)(k, \xi) \in C^\infty_0((a, b) \times {\bf R}^{n-1})  \right\}.
$$
Since then $(F_0 \mathcal F_0^{(+)} H_0 f)(k, \xi)= k^2 \phi(k, \xi)\in C^\infty_0((a, b) \times {\bf R}^{n-1})$,
we have $H_0 {\mathcal D}_e \subset {\mathcal D}_e$, it suffices to show that
\begin{equation} \label{E2.9}
\int_{-\infty}^{\infty}||e^{-itH_0} f||_{\chi^{-s}} dt < \infty
\end{equation}
This is proved in the same way as in Theorem 1.5.5. In fact, letting $u(t,\xi,y) = F_0e^{-itH_0}f$, we have $$
u(t,\xi,y) =\int_0^\infty\frac{\left( 2k \hbox{sinh}(k \pi) \right)^{1/2}}{\pi} \left( \frac{|\xi|}{2} \right)^{ik}
y^{(n-1)/2} K_{ik}(|\xi|y) e^{-itk^2} \phi(k, \xi) dk,
$$
(cf. Chap. 1, (\ref{1.5.5})). Then, similar to Chap. 1, (\ref{C1S5deltainfty}), we show that, for any $\sigma >0$,
\begin{equation}
\int_{\delta}^{\infty}\|u(t,\cdot,y)\|^2_{L^2({\bf R}^{n-1})}\frac{dy}{y^n} \leq C_N(1 + |t|)^{-N}, \quad \forall N > 0
\label{E2.8}
\end{equation}
To consider the behavior of $u(t, \cdot, y)$ for $0 < y<\sigma$, we, similar to the proof of Theorem 1.5.5,
use the decomposition
$$
u(t,\xi,y) = u_0^{(+)}(t,\xi,y) + u_0^{(-)}(t,\xi,y)+u_1(t, \xi, y),
$$
which have the same representations as in Theorem 1.5.5 with, however,
$e^{-ikt}$ replaced by $e^{-ik^2t}$. Since, for $k \in (a, b)$ and bounded $|\xi|, y,$
we have
$$
|r_1(k, |\xi|, y)| \leq C |\xi| y,\quad | \partial_k^2 r_1(k, |\xi|, y)| \leq C \log (|\xi| y)|\xi| y,
$$
(see (\ref{eq:Chap1Sec3DefinitionofInu}), (\ref{eq:KnuandInu})), we see that, for $ y < \sigma,$
$$
|u_1(t, \xi, y)| \leq C_{\phi} y^{(n+1)/2} (1+|\log (y)|) (1+|t|)^{-2}.
$$
This implies that
\begin{equation} \label{E2.7}
\int_{-\infty}^\infty \left( \int_0^\sigma ||u_1(t, \cdot, y)||_{L^2({\bf R}^(n-1)}  \frac{dy}{y^n} \right)^{1/2} dt 
<\infty.
\end{equation}
Using (\ref{eq:Chap1Sect5u0txiy}), we see that, for $t >\frac{2 |\log(y)|}{a}$ and $t <\frac{ |\log(y)|}{2b}$,
\begin{equation} \label{E2.6}
|u_0^{\pm}(t, \xi, y)|  \leq C_{\phi} y^{(n-1)/2} (1+|t|)^{-2},
\end{equation}
which implies that
$$
\int_{-\infty}^\infty \left( \int_0^\infty ||u_0^{\pm}(t, \cdot, y)
 \Theta_{a, b}(y, t)||_{L^2({\bf R}^{n-1})} (1+|\log(y)|)^{-2s} \frac{dy}{y^n} \right)^{1/2} dt 
< \infty
$$
for $s > 1$. 
Here $\Theta_{a,b}(y, t)=1$, if  $t >\frac{2 |\log(y)|}{a}$ and $t <\frac{ |\log(y)|}{2b}$, and $0$ otherwise.
As for the remaining part, we have, by the stationary phase method, that, for $\frac{1}{2b} <
\frac{|t|}{|\log(y)|}<\frac{2}{a}$,
$$
|u_0^{\pm}(t, \xi, y)|  \leq C_{\phi} y^{(n-1)/2} \left(|t|+|\log(y)| \right)^{-1/2}
$$
Taking into account that the domain of integration with respect to $\xi$ is bounded, we obtain that
$$
\int_0^\sigma|\log(y)|^{-2s} ||u_0^{\pm}(t, \cdot, y)||_{L^2({\bf R}^{n-1})} 
\left( 1- \Theta_{a,b}(y, t) \right) \frac{dy}{y^n} \leq C_{\phi} (1+|t|)^{-2s}.
$$
This estimate, together with (\ref{E2.7}), shows that
$$
\int_{-\infty}^{\infty}||u(t, \cdot, y) H(\sigma-y)||_{\chi^{-s}} dt < \infty,
$$
which, due to (\ref{E2.8}), implies (\ref{E2.9}).

As for the spectral representation, we put
\begin{equation}
\begin{split}
\left(U_0(\lambda)f\right)(x) & = \frac{\big(2\sqrt{\lambda}\sinh(\sqrt{\lambda}\pi)\big)^{1/2}}{\pi}
(2\pi)^{-(n-1)/2} \\
& \times \iint\limits_{{\bf R}^{n-1}\times(0,\infty)}
e^{ix\cdot\xi} 
\Big(\frac{|\xi|}{2}\Big)^{-i\sqrt{\lambda}}
y^{(n - 1)/2}K_{i\sqrt{\lambda}}(|\xi|y)\hat f(\xi,y)
\frac{d\xi dy}{y^n}.
\end{split}
\nonumber
\end{equation}
and ${\bf h} = L^2({\bf R}^{n-1})$, $\rho(\lambda) = \frac{1}{2}\lambda^{-1/2}$. Then the assumptions
 (A-5), (A-6) are fulfilled. Taking
\begin{equation}
U_{\pm}(\lambda) = U_0(\lambda)(1 - VR(\lambda \pm i0))
\nonumber
\end{equation}
gives, due to Theorem 4.7, the spectral representation for $H$, where $R(z) = (H - z)^{-1}$ and $V = H - H_0$.


\subsection{Absolutely continuous subspace}
Let us recall the well-known classification of the spectra of self-adjoint operators. Let $H = \int_{-\infty}^{\infty}\lambda dE_H(\lambda)$ be a self-adjoint operator in a Hilbert space $\mathcal H$. Then for any $u \in \mathcal H$, $(E_H(I)u,u)$, where $I$ is any Borel set in ${\bf R}$, defines a Borel measure on ${\bf R}$. Then the {\it absolutely continuous subspace} for $H$ is defined by
\begin{equation}
\mathcal H_{ac}(H) = \{u \in {\mathcal H}\, ; \, (E_H(\cdot)u,u) \ {\rm is \ absolutely \ continuous \ with \ respect \ to } \ d\lambda\}.
\label{C2S5Absocontisp}
\end{equation}
This is a closed subspace in $\mathcal H$.
The importance of this subspace is that it is usually stable under the perturbation appearing in scattering phenomena (see e.g. \cite{Ka76}). 

Let $R_H(z) = (H-z)^{-1}$, and $I$ be an open interval in $\sigma(H)$. If the limiting absorption principle holds on $I$, i.e. the condition (A-2) in \S 4 is guaranteed on $I$, we have
\begin{equation}
E_H(I)\mathcal H \subset {\mathcal H}_{ac}(H).
\label{C2S5EHIsubsetHacH}
\end{equation}
In fact, for $u$ in a dense subset of $\mathcal H$, we have by Stone's formula
$$
(E_H(B)u,u) = \frac{1}{2\pi i}\int_B\left((R_H(\lambda + i0)- R_H(\lambda - i0))u,u\right)d\lambda,
$$
for any Borel set $B$ in $I$, which yields (\ref{C2S5EHIsubsetHacH}). Therefore, for our case of $H = - \Delta_g$ for the asymptotically Euclidean metric, or 
$H = - \Delta_g - (n-1)^2/4$ for the asymptotically hyperbolic metric,
$$
E_H((0,\infty))\mathcal H = {\mathcal H}_{ac}(H).
$$
In this case, we often say that the continuous spectrum of $H$ is absolutely continuous, or $H$ has no singular continuos spectrum.

The spectral representation $U^{(\pm)}$ is then a unitary operator from $\mathcal H_{ac}(H)$ to the representation space $L^2((0,\infty));{\bf h};\rho(\lambda)d\lambda)$, where ${\bf h} = L^2(S^{n-1})$ for the Euclidean metric, and ${\bf h} = L^2({\bf R}^{n-1})$ for the hyperbolic metric.


\section{Geometric S-matrix}
In \S 4 and \S 5, we have constructed two Fourier transforms $U_{\pm}$ for $H = H_0 + V$, however only one Fourier transform $U_0$ is adopted for $H_0$. As a matter of fact, it is natural to associate two kinds of Fourier transforms also with $H_0$. To see this let us recall that the Green operator for $- \Delta - \lambda$ on  ${\bf R}^3$ is written as
$$
(- \Delta - \lambda \mp i0)^{-1}f = 
\frac{1}{4\pi} \int_{{\bf R}^3}
\frac{e^{\pm i\sqrt{\lambda}|x-y|}}{|x-y|}f(y)dy.
$$
Noting the asymptotic expansion 
$|x - y| \sim r - \omega\cdot y \ (\omega = x/r)$ 
as $r = |x| \to \infty$, we have for $f \in C_0^{\infty}({\bf R}^3)$
$$
 (- \Delta - \lambda \mp i0)^{-1}f \sim
\frac{e^{\pm i\sqrt{\lambda}r}}{4\pi r} \int_{{\bf R}^3}
e^{\mp i\sqrt{\lambda}\omega\cdot y}f(y)dy,
\quad (r \to \infty).
$$
This suggests that we have two Fourier transforms 
$$
\left(U_0^{(\pm)}(\lambda)f\right)(\omega) = (2\pi)^{-n/2}
\int_{{\bf R}^n}
e^{\mp i\sqrt{\lambda}\omega\cdot y}f(y)dy
$$
for $H_0 = - \Delta$ in ${\bf R}^n$. They are related as
$$
U_0^{(+)}(\lambda) = JU_0^{(-)}(\lambda),
$$
where $J$ is the unitary operator on $L^2(S^{n-1})$ defined by
\begin{equation}
 J : \varphi(\omega) \to \varphi(- \omega).
 \label{eq:Chap2Sect6Antipodal}
\end{equation}

In the case of the hyperbolic space ${\bf H}^n$, two Fourier transforms for $H_0 = - \Delta_g$ are defined by
\begin{equation}
\begin{split}
\left(U_0^{(\pm)}(\lambda)f\right)(x) & = \frac{\big(2\sqrt{\lambda}\sinh(\sqrt{\lambda}\pi)\big)^{1/2}}{\pi}
(2\pi)^{-(n-1)/2} \\
& \times \iint\limits_{{\bf R}^{n-1}\times(0,\infty)}
e^{ix\cdot\xi} 
\Big(\frac{|\xi|}{2}\Big)^{\mp i\sqrt{\lambda}}
y^{(n - 1)/2}K_{i\sqrt{\lambda}}(|\xi|y)\hat f(\xi,y)
\frac{d\xi dy}{y^n}.
\end{split}
\nonumber
\end{equation}
They are related as
$$
U_0^{(+)}(\lambda) = J(\lambda)U_0^{(-)}(\lambda),
$$
\begin{equation}
J(\lambda) = F_0^{\ast}\Big(\frac{|\xi|}{2}\Big)^{-2i\sqrt{\lambda}}
F_0.
\nonumber
\end{equation}

Let us return to the abstract theory in \S 4. Assume that we have two spectral representatios  ${\mathcal F}_0^{(\pm)}$  for $H_0$. Define
\begin{eqnarray}
 {\mathcal F}^{(\pm)}(\lambda) &=& {\mathcal F}^0(\lambda)(1 - VR(\lambda \pm i0)), \nonumber \\
  {\mathcal F}^{0}(\lambda) &=& {\mathcal F}_0^{(+)}(\lambda), 
  \nonumber \\
   \nonumber
  {\mathcal G}^{(\pm)}(\lambda) &=& {\mathcal F}_0^{(\pm)}(\lambda)
  (1 - VR(\lambda \pm i0)). \nonumber 
  \end{eqnarray}
Note that
$$
{\mathcal G}^{(+)}(\lambda) = {\mathcal F}^{(+)}(\lambda).
$$
Then by Theorem 4.7, ${\mathcal F}^{(\pm)}, \ {\mathcal G}^{(\pm)}$ give  spectral representations for $H$. The S-matrix in \S 4 is defined through ${\mathcal F}^{(\pm)}(\lambda)$. Namely
\begin{eqnarray*}
 \widehat S &=& {\mathcal F}^{(+)} \Big({\mathcal F}^{(-)}\Big)^{\ast}, \\
 \widehat S(\lambda) &=& 1 - 2\pi i\,\rho(\lambda){\mathcal F}^0(\lambda)(V - VR(\lambda + i0)V){\mathcal F}^0(\lambda)^{\ast}, \\
 &=& 1 - 2\pi i\,\rho(\lambda){\mathcal F}^{(+)}(\lambda)V{\mathcal F}^0(\lambda)^{\ast}.
\end{eqnarray*}
 
Here we introduce a new assumption.

\medskip
\noindent
(A-7) {\it There exists a unitary operator $J(\lambda)$ on  ${\bf h}$ 
satisftying}
\begin{equation}
 {\mathcal F}_0^{(+)}(\lambda) = J(\lambda){\mathcal F}_0^{(-)}(\lambda).
 \nonumber
\end{equation}

 We define a unitary operator $J$ on 
$L^2(I;{\bf h};d\lambda)$ by
\begin{equation}
 \big(Jf\big)(\lambda) = J(\lambda)f(\lambda).
 \nonumber
\end{equation}
Then we have
\begin{equation}
{\mathcal F}^{(-)}(\lambda) = J(\lambda){\mathcal G}^{(-)}(\lambda), \quad 
 {\mathcal F}^{(-)} = J{\mathcal G}^{(-)}.
 \nonumber
\end{equation}
We define a new scattering operator by
\begin{equation}
 \widehat S_{geo} = {\mathcal G}^{(+)}\Big({\mathcal G}^{(-)}\Big)^{\ast},
 \nonumber
\end{equation}
and a new scattering matrix by
\begin{equation}
\begin{split}
 \widehat S_{geo}(\lambda) & = \widehat S(\lambda)J(\lambda)  \\
 & = J(\lambda) - 2\pi i\,\rho(\lambda){\mathcal F}^{(+)}(\lambda)V{\mathcal F}_0^{(-)}(\lambda)^{\ast}.
 \end{split}
 \label{eq:Chap2Sect6GeometricSmatrix}
\end{equation}
We call $\widehat S_{geo}(\lambda)$ the {\it geometric scattering matrix}. Since ${\mathcal F}^{(+)} = {\mathcal G}^{(+)}$, we have
\begin{equation}
\widehat S_{geo} = \widehat S J,
\nonumber
\end{equation}
and the following theorem holds.


\begin{theorem}
 $\widehat S_{geo}(\lambda)$ is unitary on ${\bf h}$, and
 $$
 \Big(\widehat S_{geo}\hat f\Big)(\lambda) = 
 \widehat S_{geo}(\lambda)\hat f(\lambda), \quad 
 \forall \hat f \in \widehat{\mathcal H}, \quad \forall \lambda \in I.
  $$
\end{theorem}

The reason why $\widehat S_{geo}(\lambda)$ is called the {\it geometric} S-matrix is as follows. Suppose we are given a Schr{\"o}dinger operator $H$ on a Riemannian manifold ${\mathcal M}$. In some cases, we can associate a {\it boundary at infinity} $\partial_{\infty}{\mathcal M}$ for ${\mathcal M}$, and construct the spectral representation ${\mathcal F}^{(\pm)}(\lambda)$ as above with ${\bf h} = L^2(\partial_{\infty}{\mathcal M})$, and prove the asymptotic expansion
$$
R(\lambda \pm i0)f \simeq C_{\pm}(\lambda)a(\rho)e^{\pm iS(\rho,\lambda)}
{\mathcal F}^{(\pm)}(\lambda)f, 
\quad (\rho \to \infty)
$$
at infinity in an appropriate topology. Here, $R(z) = (H - z)^{-1}$ and  $\rho$ is a geodesic distance from a fixed point $x_0$ of $\mathcal M$. Moreover the solutions of the equation $(H - \lambda)u  = 0$ belonging to a certain class admit the following asymptotic expansion at infinity
$$
u \simeq C_-(\lambda)a(\rho)e^{-iS(\rho,\lambda)}\varphi_-  + 
C_+(\lambda)a(\rho)e^{+iS(\rho,\lambda)}\varphi_+,
$$
$$
\varphi_+ = \widehat S_{geo}(\lambda)\varphi_-,
$$
(see e.g. \cite{Me95}).
The geometric S-matrix is non-trivial even for the case $V = 0$, since $\widehat S_{geo}(\lambda) = J(\lambda)$. We shall discuss these facts in the next section for the case of ${\bf R}^n$ and ${\bf H}^n$.


\section{Helmholtz equation and geometric S-matrix}


\subsection{The case of ${\bf H}^n$} 
We incoporate the results in Chap. 1 \S 4 and Chap. 2 \S 5. 
For $k > 0$ we define $\mathcal F_0^{(\pm)}(k)$ by Chap. 1 (\ref{eq:Chap1Sec4FormulaF0k}) and put
\begin{equation}
\mathcal F^0(k) = \mathcal F_0^{(+)}(k),
\nonumber
\end{equation}
\begin{equation}
\mathcal F^{(\pm)}(k) = \mathcal F^0(k)\big(1 - VR((k \pm i0)^2)\big),
\label{eq:Chap2Sect7Fplusminusk}
\end{equation}
and ${\mathcal H}_{\pm} = L^{2,\pm s}$ for $s > 1/2$. Note that we write $(k \pm i0)^2$ instead of $k^2 \pm i0$. Later this choice will turn out to be convenient. Then $\mathcal F^0(k) \in {\bf B}(L^{2,s};L^2({\bf R}^{n-1}))$, and Theorem 4.7, together with the results of section 5.2, implies
\begin{equation}
\frac{k}{\pi i}\left(\left[R(k^2 + i0) - R(k^2 - i0)\right]f,f\right) = \|\mathcal F^{(\pm)}(k)f\|_{L^2({\bf R}^{n-1})}^2,
\nonumber
\end{equation}
where $R(z) = (H - z)^{-1}$.
Therefore by Theorem 2.3, for any $0 < a < b < \infty$ there exists a constant $C > 0$ such that
\begin{equation}
\|{\mathcal F}^{(\pm)}(k)f\|_{L^2({\bf R}^{n-1})} \leq C\|f\|_{{\mathcal B}}.
\quad a < \forall k < b,
\label{7.1'}
\end{equation}
By the argument in \S 4, we have the following theorem. Let $E(\lambda)$ be the resolution of identity for $H$.


\begin{theorem}
(1) $\mathcal F^{(\pm)}$ defined by $\big({\mathcal F}^{(\pm)}f\big)(k) = \mathcal F^{(\pm)}(k)f$ is uniquely extended to a unitary operator from $E((0,\infty))L^2({\bf H}^n)$ to $L^2((0,\infty);L^2(R^{n-1});dk)$. Moreover, 
\begin{equation}
\left(\mathcal F^{(\pm)}Hf\right)(k) = k^2\left(\mathcal F^{(\pm)}f\right)(k), \quad \forall k > 0, \quad 
\forall f \in D(H).
\nonumber
\end{equation}
(2) For $f \in E((0,\infty))L^2({\bf H}^n)$, the inversion formula holds:
$$
f = \mathop{\rm s-lim}_{N\to\infty}\int_{1/N}^N
\mathcal F^{(\pm)}(k)^{\ast}(\mathcal F^{(\pm)}f)(k)dk.
$$
(3) $\mathcal F^{(\pm)}(k)^{\ast} \in {\bf B}(L^2({\bf R}^{n-1});\mathcal B^{\ast})$ is an eigenoperator of $H$ in the sense that
$$
(H - k^2)\mathcal F^{(\pm)}(k)^{\ast}\phi = 0,\quad \forall \phi \in L^2({\bf R}^{n-1}).
$$
(4) The wave operators
\begin{equation}
W_{\pm} = \mathop{\rm s-lim}_{t\to\pm\infty}e^{itH}e^{-itH_0}
\nonumber
\end{equation}
exist and $W_{\pm} = (\mathcal F^{(\pm)})^{\ast}\mathcal F^0$. \\
\noindent
(5) The S-matrix is written as
\begin{equation}
\widehat S(k) = 1 - \frac{\pi i}{k}\mathcal F^{(+)}(k)V\mathcal F^0(k)^{\ast},
\label{eq:Chap2Sect7Smatrx}
\end{equation}
and  satisfies
\begin{equation}
\mathcal F^{(+)}(k) = \widehat S(k)\mathcal F^{(-)}(k).
\label{eq:Chap2Sect7FplusSmatrixFminus}
\end{equation}
\end{theorem}

\bigskip
We next consider the geometric scattering matrix for $H$. For $k > 0$ we define
\begin{equation}
{\mathcal G}^{(\pm)}(k) = \mathcal F_0^{(\pm)}(k)\big(1 - VR((k \pm i0)^2)\big).
\label{eq:Chap2Sect7Gplusminus}
\end{equation}
As above, $\mathcal G^{(\pm)}(k) \in {\bf B}({\mathcal B};L^2({\bf R}^{n-1}))$ and $\mathcal G^{(\pm)}$ give other spectral representations for $H$. Note that, letting $F_0$ be the Fourier transform on ${\bf R}^{n-1}$, we have
\begin{equation}
 {\mathcal F}_0^{(+)}(k) = J(k){\mathcal F}_0^{(-)}(k),
 \nonumber
\end{equation}
\begin{equation}
J(k) = F_0^{\ast}\Big(\frac{|\xi|}{2}\Big)^{-2ik}F_0.
\label{eq:Chap2Sect7Jk}
\end{equation}

We extend Theorem 1.4.7 for $H$. For $u, v \in {\mathcal B}^{\ast}$, we define
\begin{equation}
u \simeq v
 \Longleftrightarrow 
 \lim_{R\to\infty}\frac{1}{\log R}
\int_{1/R}^R\|u(y) - v(y)\|^2_{L^2({\bf R}^{n-1})}
\frac{dy}{y^{n}} = 0.
\nonumber
\end{equation}


\begin{lemma} Let $\chi(y) = 1 \ (y < 1/2), \ \chi(y) = 0 \ (y > 1)$, and $\omega_{\pm}(k)$ be as in Chap. 1 (\ref{eq:Chap1Sect4omegaplusminusk}). Then for any $\varphi \in L^2({\bf R}^{n-1})$ and $k > 0$
\begin{equation}
\begin{split}
 {\mathcal F}_0^{(+)}(k)^{\ast}\varphi 
 & \ \simeq \frac{k}{\pi i}\omega_+(k)\chi(y)y^{(n-1)/2-ik}
\varphi \\
& \ \ \ \ \ - \frac{k}{\pi i}\omega_-(k)\chi(y)y^{(n-1)/2+ik}J(k)^{\ast}\varphi,
\end{split}
\nonumber
\end{equation}
\begin{equation}
\begin{split}
 {\mathcal F}_0^{(-)}(k)^{\ast}\varphi 
 & \ \simeq \frac{k}{\pi i}\omega_+(k)\chi(y)y^{(n-1)/2-ik}
J(k)\varphi \\
& \ \ \ \ \ - \frac{k}{\pi i}\omega_-(k)\chi(y)y^{(n-1)/2+ik}\varphi.
\end{split}
\nonumber
\end{equation}
\end{lemma}
Proof. In view of Chap. 1
(\ref{eq:Chap1Sect4F0yplusinfty}), we have only to compute the behavior of the left-hand side as $y \to 0$ for $\widehat{\varphi} \in C_0^{\infty}({\bf R}^{n-1})$. We use Chap.1 (\ref{eq:Chap1Sec3Knunear0}) in the expression Chap.1
(\ref{eq:Chap1Sect4F0kastpsi}) and compute directly to prove the lemma. \qed


\begin{lemma} Let $\chi(y)$ and $\omega_{\pm}(k)$ be as in the previous lemma. Then, for
$f \in {\mathcal B}$ and $k > 0$,
\begin{equation}
R(k^2\pm i0)f \simeq 
\omega_{\pm}(k)\chi(y)y^{(n-1)/2\mp ik}{\mathcal G}^{(\pm)}(k)f.
\nonumber
\end{equation}
\end{lemma}
Proof. The lemma follows from the resolvent equation
$$
R(k^2 \pm i0) = R_0(k^2 \pm i0) - R_0(k^2 \pm i0)V
R(k^2\pm i0),
$$
Lemmas 4.7, 4.9 of Chap.1 and (\ref{eq:Chap2Sect7Gplusminus}). \qed

\bigskip
By (\ref{eq:Chap2Sect6GeometricSmatrix}), the geometric scattering matrix is defined to be
\begin{equation}
\widehat S_{geo}(k) = J(k) - \frac{\pi i}{k} {\mathcal F}^{(+)}(k)V{\mathcal F}_0^{(-)}(k)^{\ast}.
\nonumber
\end{equation}


\begin{lemma}
For $\varphi \in L^2({\bf R}^{n-1})$ 
\begin{equation}
\begin{split}
{\mathcal G}^{(-)}(k)^{\ast}\varphi & \simeq
\frac{k}{\pi i}\omega_+(k)\chi(y) y^{(n-1)/2-ik}{\widehat S}_{geo}(k)\varphi \\
& \ \ \ \ \ - \frac{k}{\pi i}\omega_-(k)\chi(y) y^{(n-1)/2+ik}\varphi.
\end{split}
\nonumber
\end{equation}
\end{lemma}
Proof.  By (\ref{eq:Chap2Sect7Gplusminus})
\begin{equation}
{\mathcal G}^{(-)}(k)^{\ast} \varphi= {\mathcal F}_0^{(-)}(k)^{\ast} \varphi
- R(k^2 + i0)V{\mathcal F}_0^{(-)}(k)^{\ast} \varphi.
\nonumber
\end{equation}
Since ${\mathcal F}^{(+)}(k)={\mathcal G}^{(+)}(k)$, we obtain, by Lemmas 7.2 and 7.3, that
\begin{equation}
\begin{split}
\nonumber
&{\mathcal G}^{(-)}(k)^{\ast} \varphi
\simeq
\frac{k}{\pi i}\omega_+(k)\chi(y) y^{(n-1)/2-ik} J(k)\varphi
\\ \nonumber
& - \frac{k}{\pi i}\omega_-(k)\chi(y) y^{(n-1)/2+ik} \varphi-
\omega_+(k)\chi(y) y^{(n-1)/2-ik} \left[J(k)- {\widehat S}_{geo}(k) \right] \varphi
\\ \nonumber
& \simeq
\frac{k}{\pi i}\omega_+(k)\chi(y) y^{(n-1)/2-ik}{\widehat S}_{geo}(k)\varphi
- \frac{k}{\pi i}\omega_-(k)\chi(y) y^{(n-1)/2+ik}\varphi.
\qed
\end{split}
\end{equation}
}


\begin{lemma} There exists a constant $C= C(k) > 0$ such that for any $\varphi \in L^2({\bf R}^{n-1})$
\begin{equation}
\lim_{R\to\infty}\frac{1}{\log R}\int_{1/R}^R\|{\mathcal G}^{(-)}(k)^{\ast}\varphi\|^2_{L^2({\bf R}^{n-1})}\frac{dy}{y^n} = 
C\|\varphi\|_{L^2({\bf R}^{n-1})}^2.
\nonumber
\end{equation}
\end{lemma}
Proof. We put
$a_+ = {\widehat S}_{geo}(k)\varphi, \ a_- = \varphi$. Then by Lemma 7.4 $\|{\mathcal G}^{(-)}(k)^{\ast}\varphi\|^2_{L^2({\bf R}^{n-1})}$ behaves like 
\begin{eqnarray*}
& |C_+(k)|^2y^{n-1}\|a_+\|_{\bf h}^2 + |C_-(k)|^2y^{n-1}\|a_-\|_{\bf h}^2 \\
& + C_+(k)\overline{C_-(k)}y^{n-1 -2ik}(a_+,a_-)_{\bf h} + 
C_-(k)\overline{C_+(k)}y^{n-1+2ik}(a_-,a_+)_{\bf h},
\end{eqnarray*}
 where $C_{\pm}(k)$ are constants. Simple computation shows that the 3rd and 4th terms tend to 0. As $\widehat S_{geo}(\lambda)$ is unitary, the lemma follows. \qed

Together with (\ref{7.1'}), this implies

\begin{cor} There is a constant $C > 0$ such that
\begin{equation}
C^{-1}\|\varphi\|_{L^2({\bf R}^{n-1})} \leq \|{\mathcal G}^{(\pm)}(k)^{\ast}\varphi\|_{{\mathcal B}^{\ast}} \leq C\|\varphi\|_{L^2({\bf R}^{n-1})}.
\nonumber
\end{equation}
\end{cor}


\begin{lemma}
If $u \in {\mathcal B}^{\ast}, \ (H - k^2)u = 0, \ f \in {\mathcal B}$, and either ${\mathcal G}^{(+)}(k)f = 0$ or ${\mathcal G}^{(-)}(k)f = 0$ holds, then $(u,f) = 0$.
\end{lemma}
Proof. The same as Lemma 1.4.10. \qed

\bigskip
These preparations are sufficient to extend Theorem 1.4.3 to $H$.


\begin{theorem}
For $k > 0$
\begin{equation}
\{u \in {\mathcal B}^{\ast} \,; \, (H - k^2)u = 0\} = 
{\mathcal G}^{(\pm)}(k)^{\ast}\big(L^2({\bf R}^{n-1})\big).
\nonumber
\end{equation}
\end{theorem}


\begin{theorem} If
$u \in {\mathcal B}^{\ast}$ satisfies $(H - k^2)u = 0$ for $k > 0$, there exist
$\varphi_{\pm} \in L^2({\bf R}^{n-1})$ such that
\begin{equation}
u \simeq \frac{k}{\pi i}\omega_+(k)\chi(y)y^{(n-1)/2-ik}\varphi_+ - 
\frac{k}{\pi i}\omega_-(k)\chi(y)y^{(n-1)/2 + ik}\varphi_-.
\nonumber
\end{equation}
Moreover,
\begin{equation}
\varphi_+ = \widehat S_{geo}(k)\varphi_-.
\nonumber
\end{equation}
\end{theorem}
Proof. By Theorem 7.8, $u$ can be written as $u = {\mathcal G}^{(-)}(k)^{\ast}\psi$. Using Lemma 7.4, we prove the theorem. \qed


\begin{theorem}
For any $\varphi_- \in L^2({\bf R}^{n-1})$, there exist unique $u \in {\mathcal B}^{\ast}$ and $\varphi_+ \in L^2({\bf R}^{n-1})$ such that the equation $(H - k^2) u = 0$ and the expansion in Theorem 7.9 hold.
\end{theorem}
Proof. The existence of such $\varphi_+$ and $u$ follows from Theorem 7.9. We prove the uniqueness. If $\varphi_- = 0$, we have $u \simeq C(k)\chi(y)y^{(n-1)/2 - ik}\varphi_+$, hence $u$ satisfies the radiation conditions (\ref{eq:Chap2Sect2RadCond1/R1}), (\ref{eq:Chap2Sect2RadCondplusinfty}). 
Then $u = 0$ by Lemma 2.12, which also proves $\varphi_+ = 0$.\qed


\subsection{The case of ${\bf R}^n$} 
It is worthwhile to give a brief look at the case of ${\bf R}^n$. We define the weighted $L^2$ space $L^{2,s}$ and the Besov type space $\mathcal B$ by
\begin{equation}
L^{2,s} \ni u \Longleftrightarrow 
\|u\|_s^2 = \int_{{\bf R}^n}(1 + |x|)^{2s}|u(x)|^2dx < \infty,
\nonumber
\end{equation}
\begin{equation}
\|u\|_{\mathcal B} = \mathop\sum_{j=0}^{\infty}2^{j/2}\|u\|_{L^2(\Omega_j)} < \infty,
\nonumber
\end{equation}
\begin{equation}
\Omega_j = \{x \in {\bf R}^n ; r_{j-1} < |x| < r_j\},
\nonumber
\end{equation}
where $r_j = 2^j \ (j \geq 0), \ r_{-1} = 0$. The dual space of $\mathcal B$ has the following equivalent norm
\begin{equation}
\|u\|^2_{{\mathcal B}^{\ast}} = \sup_{R>1}\frac{1}{R}
\int_{|x|<R}|u(x)|^2dx.
\nonumber
\end{equation}
Let $H$ be as in subsection 5.1, ${\bf h} = L^2(S^{n-1})$, and put for $k > 0$
\begin{equation}
\left(\mathcal F_0^{(\pm)}(k)f\right)(\omega) = (2\pi)^{-n/2}\int_{{\bf R}^n}e^{\mp ik\omega\cdot x}
f(x)dx,
\nonumber
\end{equation}
\begin{equation}
\mathcal F^0(k) = \mathcal F_0^{(+)}(k),
\nonumber
\end{equation}
\begin{equation}
\mathcal F^{(\pm)}(k) = \mathcal F^0(k)\big(1 - VR((k \pm i0)^2)\big),
\nonumber
\end{equation}
\begin{equation}
\mathcal G^{(\pm)}(k) = \mathcal F_0^{(\pm)}(k)\big(1 - VR((k \pm i0)^2)\big).
\nonumber
\end{equation}
 Then the results in \S 5 and \S 6 can be applied to $H$. 
 Let $E(\lambda)$ be the resolution of identity for $H$.


\begin{theorem}
(1) $\mathcal F^{(\pm)}$ defined by $\big({\mathcal F}^{(\pm)}f\big)(k) = \mathcal F^{(\pm)}(k)f$ is uniquely extended to a unitary operator from $E((0,\infty))L^2({\bf R}^n)$ to $L^2((0,\infty);L^2(S^{n-1});k^{n-1}dk)$. Moreover 
\begin{equation}
\left(\mathcal F^{(\pm)}Hf\right)(k) = k^2\left(\mathcal F^{(\pm)}f\right)(k), \quad \forall k > 0, \quad 
\forall f \in D(H).
\nonumber
\end{equation}
(2) For $f \in E((0,\infty))L^2({\bf R}^n)$, the inversion formula holds:
$$
f = \mathop{\rm s-lim}_{N\to\infty}\int_{1/N}^N
\mathcal F^{(\pm)}(k)^{\ast}(\mathcal F^{(\pm)}f)(k)k^{n-1}dk.
$$
(3) $\mathcal F^{(\pm)}(k)^{\ast} \in {\bf B}(L^2({S}^{n-1});\mathcal B^{\ast})$ is an eigenoperator of $H$ in the sense that
$$
(H - k^2)\mathcal F^{(\pm)}(k)^{\ast}\phi = 0,\quad \forall \phi \in L^2(S^{n-1}).
$$
\end{theorem}


\section{Modified Radon transform}


\subsection{Extension of the Fourier transform}
In order to construct the modified Radon transform associated with $H$ in \S 2, we extend the definition of the generalized Fourier transform for all $k \in {\bf R}$. 
Let us repeat the definitions of the Fourier transforms introduced so far:
\begin{equation}
\begin{split}
\left({\mathcal F}^{(\pm)}_0(k)f\right)(x) & = \sqrt{\frac{2}{\pi}}\,k\,
\sqrt{\frac{\sinh(k\pi)}{k\pi}} \\
& \ \ \ \ \times F_0^{\ast}\left(\Big(\frac{|\xi|}{2}\Big)^{\mp ik}\int\limits_0^{\infty}
y^{(n - 1)/2}K_{ik}(|\xi|y)\widehat f(\xi,y)
\frac{dy}{y^n}\right), \\
\mathcal F^0(k) &= \mathcal F_0^{(+)}(k), \\
\mathcal F_0(k) & = \frac{1}{\sqrt2}\Omega(k)\mathcal F^0(k),\\
\Omega(k) & = \frac{-i}{\Gamma(1-ik)}\sqrt{\frac{k\pi}{\sinh(k\pi)}},\\
J(k) &= F_0^{\ast}\left(\frac{|\xi|}{2}\right)^{-2ik}F_0,
\end{split}
\label{E2.15}
\end{equation}
$F_0$ being the Fourier transformation on ${\bf R}^{n-1}$. We have also defined
\begin{equation}
\mathcal F^{(\pm)}(k) = \mathcal F^{(+)}_0(k)(1 - VR((k \pm i0)^2)).
\nonumber
\end{equation}
Note that the operators $\mathcal F_0^{(\pm)}(k)$, $\mathcal F^{(\pm)}(k)$ can be extended using the above formulae for $0 \neq k \in {\bf R}$ and, by 
(\ref{C1S3KikK-ik}) of Chap. 1,
\begin{equation}
\mathcal F_0^{(+)}(k) = - \mathcal F_0^{(-)}(-k) = J(k)\mathcal F_0^{(-)}(k) = - J(k)\mathcal F_0^{(+)}(-k),
\nonumber
\end{equation}
\begin{equation}
\mathcal F^{(+)}(k) = -J(k)\mathcal F^{(-)}(-k).
\nonumber
\end{equation}
We now define a new Fourier transformation $\mathcal F_{\pm}(k)$ by
\begin{equation}
\mathcal F_{\pm}(k) = \frac{1}{\sqrt2}\Omega(\pm k)\mathcal F^{(\pm)}(k), \quad 0 \neq k \in {\bf R},
\label{E2.1^}
\end{equation}
and put $(\mathcal F_{\pm}f)(k) = \mathcal F_{\pm}(k)f$. 
Let $\widehat S(k)$ be the S-matrix defined by (\ref{eq:Chap2Sect7Smatrx}). Then by (\ref{eq:Chap2Sect7FplusSmatrixFminus}), we have
\begin{equation}
\mathcal F_+(k) = \frac{\Gamma(1 + ik)}{\Gamma(1 - ik)}\widehat S(k)\mathcal F_-(k), \quad k > 0.
\nonumber
\end{equation}
By definition we also have
\begin{equation}
\mathcal F_+(-k) = - J(-k)\mathcal F_-(k).
\nonumber
\end{equation}
The following Theorem can be proved easily from the above formulas.


\begin{theorem}
(1) ${\mathcal F}_{\pm} : L^2({\bf H}^n) \to L^2({\bf R};L^2({\bf R}^{n-1});dk)$ is a partial isometry with initial set $E((0,\infty))L^2({\bf H}^n)$, $E(\lambda)$ being the resolution of identity for $H$, and
$$
(\mathcal F_{\pm}Hf)(k) = k^2(\mathcal F_{\pm}f)(k), \quad k \in {\bf R}, \quad f \in D(H).
$$
(2) For $k > 0$, we have
$$
\mathcal F_+(k) = - \frac{\Gamma(1 + ik)}{\Gamma(1 - ik)}\widehat S(k)J(k)\mathcal F_+(-k).
$$
Consequently, the range of $\mathcal F_{\pm}$ has the following characterization:
$$
g \in {\rm Ran}\,\mathcal F_+ \Longleftrightarrow  g(k) = - \frac{\Gamma(1 + ik)}{\Gamma(1 - ik)}\widehat S(k)J(k)g(-k), \quad k > 0,
$$
$$
g \in {\rm Ran}\,\mathcal F_- \Longleftrightarrow  J(k)g(-k) = - \frac{\Gamma(1 + ik)}{\Gamma(1 - ik)}\widehat S(k)g(k), \quad k > 0.
$$
\end{theorem}

Note that the above relation is rewritten as
$$
g \in {\rm Ran}\,\mathcal F_+ \Longleftrightarrow  g(k) = - \frac{\Gamma(1 + ik)}{\Gamma(1 - ik)}\widehat S_{geo}(k)g(-k), \quad k > 0.
$$

We put
\begin{equation}
\mathcal H_{>0} = L^2((0,\infty);L^2({\bf R}^{n-1});dk), \quad
\mathcal H_{<0} = L^2((-\infty,0);L^2({\bf R}^{n-1});dk),
\label{C2S8Hplusminus}
\end{equation}
and let $r_{+}$ and $r_-$ be the projections onto $\mathcal H_{>0}$ and $\mathcal H_{<0}$, respectively.


\begin{lemma}
\begin{equation}
W_+ = 2(\mathcal F_+)^{\ast}r_+\mathcal F_0, \quad 
W_- = 2(\mathcal F_+)^{\ast}r_-\mathcal F_0,
\label{eq:Chap2Sect8Fplusrplus}
\end{equation}
\begin{equation}
W_+ = 2(\mathcal F_-)^{\ast}Gr_-\mathcal F_0, \quad
W_- = 2(\mathcal F_-)^{\ast}Gr_+\mathcal F_0,
\label{eq:Chap2Sect8Fminustminus}
\end{equation}
where $G$ is the operator of multiplication by $\dfrac{\Gamma(1 -ik)}{\Gamma(1 + ik)}$.
\end{lemma}
Proof. Recall that $|\Omega(k)| = 1$ and $J(k)$ is unitary on $L^2({\bf R}^{n-1})$.
By Theorem 7.1(4),
using $\mathcal F^0(-k) = - J(-k)\mathcal F^0(k)$ and $\mathcal F^{(-)}(-k) = - J(-k)\mathcal F^{(+)}(k)$, we have, for $f, g \in {\mathcal B}$,
\begin{eqnarray*}
(W_-f,g) &=& (\mathcal F^0f,\mathcal F^{(-)}g) \\
&=& \int_0^{\infty}(\mathcal F^0(k)f,\mathcal F^{(-)}(k)g)dk \\
&=& \int_{-\infty}^0(J(-k)\mathcal F^0(k)f,J(-k)\mathcal F^{(+)}(k)g)dk \\
&=& \int_{-\infty}^0(\Omega(k)\mathcal F^0(k)f,\Omega(k)\mathcal F^{(+)}(k)g)dk \\
&=& 2\int_{-\infty}^0(\mathcal F_0(k)f,\mathcal F_+(k)g)dk \\
&=& (2(\mathcal F_+)^{\ast}r_-\mathcal F_0f,g),
\end{eqnarray*}
which proves (\ref{eq:Chap2Sect8Fplusrplus}) for $W_-$. By the similar and simpler manner, one can prove (\ref{eq:Chap2Sect8Fplusrplus}) for $W_+$. 
Using $\mathcal F^0(-k) = - J(-k)\mathcal F^0(k)$ and $\mathcal F^{(+)}(-k) = - J(-k)\mathcal F^{(-)}(k)$, we have for $f, g \in {\mathcal B}$
\begin{eqnarray*}
(W_+f,g) &=& (\mathcal F^0f,\mathcal F^{(+)}g) \\
&=& \int_0^{\infty}(\mathcal F^0(k)f,\mathcal F^{(+)}(k)g)dk \\
&=& \int_{-\infty}^0(J(-k)\mathcal F^0(k)f,J(-k)\mathcal F^{(-)}(k)g)dk \\
&=& \int_{-\infty}^0\frac{\Omega(-k)}{\Omega(k)}(\Omega(k)\mathcal F^0(k)f,\Omega(-k)\mathcal F^{(-)}(k)g)dk \\
&=& 2\int_{-\infty}^0\frac{\Omega(-k)}{\Omega(k)}(\mathcal F_0(k)f,\mathcal F_-(k)g)dk \\
&=& (2(\mathcal F_-)^{\ast}Gr_-\mathcal F_0f,g),
\end{eqnarray*}
which proves (\ref{eq:Chap2Sect8Fminustminus}) for $W_+$. Similarly, we can prove (\ref{eq:Chap2Sect8Fminustminus}) for $W_-$. \qed

We define operators $\hat I$ and $U$ on $L^2({\bf R};L^2({\bf R}^{n-1};dk)$ by
\begin{equation}
\begin{split}
(\hat If)(k) & = f(-k), \\
(Uf)(k) & = \frac{\Gamma(1 - ik)}{\Gamma(1 + ik)}
\big(F_0^{\ast}\left(\frac{|\xi|}{2}\right)^{2ik}F_0f\big)(k).
\end{split}
\nonumber
\end{equation}
Direct computation shows the following relations:
\begin{equation}
\begin{split}
\hat Ir_+ &= r_- \hat I, \\
\hat IU\hat I &= U^{-1}, \\
Ur_{\pm} & = r_{\pm}U.
\end{split}
\label{eq:Chap2Sect8IUrplusminus}
\end{equation}


\begin{lemma} 
\begin{equation} 
\mathcal F_0(\mathcal F_0)^{\ast} = \frac{1}{2}(I + \hat IU).
\label{eq:Chap2Sect8F0F0ast}
\end{equation}
\end{lemma}
Proof. Let $\Pi = (I + \hat IU)/2$. Then by (\ref{eq:Chap2Sect8IUrplusminus}), one can show $\Pi^{\ast} = \Pi^2 = \Pi$. Moreover, $g = \Pi f$ satisfies $\hat Ig = Ug$. Therefore by Lemma 1.5.2 (3), $\Pi$ is the projection onto the range of $\mathcal F_0$. \qed


\begin{lemma}
\begin{equation}
\mathcal F_+ = r_+\mathcal F_0(W_+)^{\ast} + r_-\mathcal F_0(W_-)^{\ast},
\label{eq:Chap2Sect8FplusWplusminus}
\end{equation}
\begin{equation}
\mathcal F_- = Gr_+\mathcal F_0(W_-)^{\ast} + Gr_-\mathcal F_0(W_+)^{\ast}.
\label{eq:Chap2Sect8FminusWplusminus}
\end{equation}
\end{lemma}
Proof. By (\ref{eq:Chap2Sect8Fplusrplus}) and (\ref{eq:Chap2Sect8F0F0ast}),
\begin{eqnarray*}
\mathcal F_0(W_+)^{\ast} &=& 2\mathcal F_0(\mathcal F_0)^{\ast}r_+\mathcal F_+ \\
&=& r_+\mathcal F_+ + \hat IUr_+\mathcal F_+.
\end{eqnarray*}
Since $\hat IUr_+ = r_-\hat IU$ by (\ref{eq:Chap2Sect8IUrplusminus}), multiplying both sides by $r_+$, we obtain
$$
r_+\mathcal F_0(W_+)^{\ast} = r_+\mathcal F_+.
$$
Similarly, we have
$$
r_-\mathcal F_0(W_-)^{\ast} = r_-\mathcal F_+.
$$
Adding these two equalities, we obtain (\ref{eq:Chap2Sect8FplusWplusminus}). The formula (\ref{eq:Chap2Sect8FminusWplusminus}) is proved in a similar manner. \qed



\subsection{Modified Radon transform}
We now define the modified Radon transform for $H$.

\begin{definition}
For $s \in {\bf R}$, we define
\begin{equation}
\left(\mathcal R_{\pm}f\right)(s) = \frac{1}{\sqrt{2\pi}}\int_{-\infty}^{\infty}e^{iks}\left(\mathcal F_{\pm}f\right)(k)dk.
\nonumber
\end{equation}
\end{definition}


\begin{theorem}
$\mathcal R_{\pm}$ is a partial isometry from $L^2({\bf H}^n)$ to $L^2({\bf R};L^2({\bf R}^{n-1});dk)$ with initial set $E((0,\infty))L^2({\bf H}^n)$. The Fourier transform of the final set of $\mathcal R_{\pm}$ is characterized by Theorem 8.1 (2). Moreover
$$
\mathcal R_{\pm}H = - \partial_s^2\mathcal R_{\pm}.
$$
\end{theorem}

The scattering operator can also be defined by the Radon transform.


\begin{definition}
We define the scattering operator $\mathcal S_R$ by
$$
\mathcal S_R = \mathcal R_+(\mathcal R_-)^{\ast}.
$$
\end{definition}


\begin{lemma}
The scattering operator $\mathcal S_R$ is a partial isometry with initial set {\rm Ran}$\,{\mathcal R}_-$ and final set {\rm Ran}$\,\mathcal R_+$. The relation between $S = (W_+)^{\ast}W_-$ and $\mathcal S_R$ is given by the following formula. Let $\mathcal F_1$ be the 1-dimensional Fourier transformation. Then 
\begin{equation}
\mathcal F_1\mathcal S_R(\mathcal F_1)^{\ast} = r_+\mathcal F_0S(\mathcal F_0)^{\ast}r_+ G^{\ast} + r_-\mathcal F_0S^{\ast}(\mathcal F_0)^{\ast}r_- G^{\ast} 
+ \frac{1}{2}{\hat IU}G^{\ast}.
\nonumber
\end{equation}
\end{lemma}
Proof. The first half of the lemma follows from the definition. Since $\mathcal F_1\mathcal S_R(\mathcal F_1)^{\ast} = \mathcal F_+(\mathcal F_-)^{\ast}$, the second half follows from Lemma 8.4 and direct computation. \qed


\subsection{Asymptotic profiles of solutions to the wave equation}
 We compute the asymptotic profile of the solution 
$$
u(t) = \cos(t\sqrt{H})f + \sin(t\sqrt{H})\sqrt{H}^{-1}g
$$
to the wave equation
\begin{equation}
\left\{
\begin{split}
&\partial_t^2u + Hu = 0, \\
& u\big|_{t=0} = f, \quad \partial_tu\big|_{t=0} = g.
\end{split}
\right.
\nonumber
\end{equation}


\begin{theorem}
For any $f \in E((0,\infty))L^2({\bf H}^n)$, we have as $t \to \infty$
\begin{equation}
\left\|\cos(t\sqrt{H})f - \frac{y^{(n-1)/2}}{\sqrt2}({\mathcal R}_+f)(-\log y - t,x)\right\|_{L^2({\bf H}^n)} \to 0,
\nonumber
\end{equation}
\begin{equation}
\left\|\sin(t\sqrt{H})f - \frac{iy^{(n-1)/2}}{\sqrt2}({\mathcal R}_+
\hbox{sgn} (-i\partial_s)f)(-\log y - t,x)\right\|_{L^2({\bf H}^n)} \to 0,
\nonumber
\end{equation}
where $\hbox{sgn}$ is defined in Theorem 1.5.5.
\end{theorem}
Proof.  
Using the relations
\begin{equation}
\mathcal F^{(+)}(k)^{\ast}
= \mathcal F^0(k)^{\ast} - R((k-i0)^2)V\mathcal F^{(+)}(k)^{\ast},
\nonumber
\end{equation} 
 we have by the spectral representation theorem
\begin{equation}
\begin{split}
e^{-it\sqrt{H}}f &= \int_0^{\infty}e^{-itk}\mathcal F^{(+)}(k)^{\ast}
\left(\mathcal F^{(+)}f\right)(k)dk \\
&= \int_0^{\infty}e^{-itk}\mathcal F^0(k)^{\ast}
\left(\mathcal F^{(+)}f\right)(k)dk \\
&\ \ -  \int_{0}^{\infty}e^{-itk}
R(k^2 - i0)V\mathcal F^0(k)^{\ast}\left(\mathcal F^{(+)}f\right)(k)dk. 
\end{split}
\label{eq:Chap2Sect8e-itsqrtHFplus}
\end{equation}
By the same computation as in the proof of Theorem 1.5.5, the first  term of the right-hand side of (\ref{eq:Chap2Sect8e-itsqrtHFplus}) tends to 
\begin{equation}
\frac{y^{(n-1)/2}}{\sqrt{\pi}}\int_0^{\infty}
e^{ik(-\log y - t)}\left(\mathcal F_+f\right)(k)dk
\nonumber
\end{equation}
as $t \to \infty$.

We need the following lemma to deal with 
the 2nd term of the right-hand side of (\ref{eq:Chap2Sect8e-itsqrtHFplus}).


\begin{lemma}
Let $A$ be a self-adjoint operator on a Hilbert space $\mathcal H$. For $\psi(k) \in C_0((0,\infty);\mathcal H)$ we put
\begin{equation}
\Psi_{\pm}(t) = \int_0^{\infty}e^{\pm ikt}\psi(k)dk.
\nonumber
\end{equation}
Then for any $\epsilon > 0$ 
\begin{equation}
\left\|\int_0^{\infty}(A - k \mp i\epsilon)^{-1}e^{\pm ikt}\psi(k)dk\right\|
\leq \int_t^{\infty}\|\Psi_{\pm}(s)\|ds
\nonumber
\end{equation}
holds. Similarly letting 
\begin{equation}
\Phi_{\pm}(t) = \int_{-\infty}^0e^{\mp ikt}\psi(k)dk
\nonumber
\end{equation}
for $\psi(k) \in C_0((-\infty,0);\mathcal H)$, we have for any $\epsilon > 0$
\begin{equation}
\left\|\int_{-\infty}^{0}(A + k \pm i\epsilon)^{-1}e^{\mp ikt}\psi(k)dk\right\|
\leq \int_{-\infty}^{t}\|G_{\mp}(s)\|ds.
\nonumber
\end{equation}
\end{lemma}

Proof. By virtue of the identity
$$
(A - k \mp i\epsilon)^{-1} = \pm i\int_0^{\infty}
e^{\mp is(A - k \mp i\epsilon)}ds,
$$
we have
$$
\int_0^{\infty}(A - k \mp i\epsilon)^{-1}e^{\pm ikt}\psi(k)dk = 
\pm i \int_0^{\infty}e^{\mp is(A \mp i\epsilon)}
\Psi_{\pm}(s + t)ds,
$$
which proves the first half of the lemma. We also have
$$
(A + k \mp i\epsilon)^{-1} = \pm i\int_{-\infty}^0
e^{\pm is(A + k \mp i\epsilon)}ds
$$
which proves the second half. \qed

\bigskip
\noindent
{\it Proof of Theorem 8.9 (continued).}
Letting $\sqrt{H} = A$, we have
$$
(H - k^2 \mp i0)^{-1} = (A - k \mp i0)^{-1}(A + k)^{-1}.
$$
Therefore, to show that the 2nd term of the right-hand side of (\ref{eq:Chap2Sect8e-itsqrtHFplus}) tends to 0, letting
$$
\psi(k) = (A + k)^{-1}V{\mathcal F}_0(k)^{\ast}
\left(\mathcal F^{(+)}f\right)(k),
$$
$$
\Psi(t) = \int_0^{\infty}e^{-ikt}\psi(k)dk,
$$
we have only to prove
\begin{equation}
\int_0^{\infty}\|\Psi(t)\|dt < \infty.
\nonumber
\end{equation}
Take $g \in L^2({\bf H}^n)$, and consider
$$
(\Psi(t),g) = \int_0^{\infty}e^{-ikt}(V\mathcal F_0(k)^{\ast}\big(\mathcal F^{(+)}f\big)(k),(A+k)^{-1}g)dk.
$$
Arguing in the same way as the proof of (A-4) in Subsection 5.2. we have
$$
|(\Psi(t),g)| \leq C(1 + t)^{-1-\epsilon}\|g\|,
$$
implying that 
 $\|\Psi(t)\| \leq C (1+t)^{-1-\epsilon}$.
We have thus derived that
\begin{equation}
\left\|e^{-it\sqrt{H}}f - \frac{y^{(n-1)/2}}{\sqrt\pi}\int_0^{\infty}
e^{ik(-\log y - t)}\left(\mathcal F_+f\right)(k)dk\right\| \to 0
\label{eq:Chap2Sect8eminusitsqrtH0inftyint}
\end{equation}
as $t \to \infty$.

By using the relation
\begin{equation}
\mathcal F^{(-)}(k)^{\ast} \mathcal F^{(-)}(k) = 
\mathcal F^{(+)}(-k)^{\ast} \mathcal F^{(+)}(-k),
\nonumber 
\end{equation}
we have as above
\begin{equation}
\begin{split}
e^{-it\sqrt{H}}f &= \int_0^{\infty}e^{-itk}\mathcal F^{(-)}(k)^{\ast}
\left(\mathcal F^{(-)}f\right)(k)dk \\
 &= \int_{-\infty}^{0}e^{itk}\mathcal F^{(+)}(k)^{\ast}
\left(\mathcal F^{(+)}f\right)(k)dk \\
&= \int_{-\infty}^0e^{itk}\mathcal F^0(k)^{\ast}
\left(\mathcal F^{(+)}f\right)(k)dk \\
&\ \ -  \int_{-\infty}^{0}e^{itk}
R(k^2 + i0)V\mathcal F^0(k)^{\ast}\left(\mathcal F^{(+)}f\right)(k)dk. 
\end{split}
\nonumber
\end{equation}
Arguing as above, we can derive
\begin{equation}
\left\|e^{-it\sqrt{H}}f - \frac{y^{(n-1)/2}}{\sqrt\pi}\int_{-\infty}^{0}
e^{ik(-\log y + t)}\left(\mathcal F_+f\right)(k)dk\right\| \to 0
\label{eq:Chap2Sect8eminusitsqrtHminusinfty0int}
\end{equation}
as $t \to -\infty$. Theorem 8.9 then follows from (\ref{eq:Chap2Sect8eminusitsqrtH0inftyint}) and (\ref{eq:Chap2Sect8eminusitsqrtHminusinfty0int}). \qed


\subsection{Invariance principle}
Suppose for two self-adjoint operators $A$ and $B$, the wave operator
$$
W_{\pm} = {\mathop{\rm s-lim}_{t\to\pm\infty}}\,e^{itA}e^{-itB}P_{ac}(B),
$$
exists, where $P_{ac}(B)$ denotes the projection onto the absolutely continuous subspace for $B$. Then, for a suitable Borel function $\phi(s)$ on ${\bf R}$, the wave operator
$$
W^{(\phi)}_{\pm} = {\mathop{\rm s-lim}_{t\to\pm\infty}}\,e^{it\phi(A)}e^{-it\phi(B)}P_{ac}(B),
$$ 
exists and $W_{\pm} = W_{\pm}^{(\phi)}$. This fact is called {\it invariance principle}, and is proved in a general setting (see e.g. pp. 545, 579  of \cite{Ka76}). We are interested in the case where $\phi(s) = \sqrt{s}$. Then $W_{\pm}$ is the wave operator for the Schr{\"o}dinger equation, and $W_{\pm}^{(\phi)}$ is the wave operator for the wave equation. 

Under the assumptions in the present chapter, we can prove this invariance principle directly for the above operators $H$ and $H_0$ on ${\bf H}^n$. In fact, letting
$$
H_{+} = E_H((0,\infty))H,
$$
where $E_H(\lambda)$ is the spectral resolution for $H$, the existence of the strong limit
\begin{equation}
{\mathop{\rm s-lim}_{t\to\pm\infty}}\, e^{it\sqrt{H_{+}}}e^{-it\sqrt{H_0}}P_{ac}(H_0)
\label{C2S8sqrtwaveop}
\end{equation}
can be proven by the same argument as that for the wave operator
$$
W_{\pm} = {\mathop{\rm s-lim}_{t\to\pm\infty}}\, e^{itH}e^{-itH_0}.
$$
Observing the proof of Theorem 8.9 (see the arguments after (\ref{eq:Chap2Sect8e-itsqrtHFplus})), we see that for $f \in {\mathcal H}_{ac}(H) = E_H((0,\infty))L^2({\bf H}^n)$ (see Chap. 2, Subsection 5.3)
$$
\Big\|e^{-it\sqrt{H}}f - \int_0^{\infty}e^{-itk}\mathcal F^0(k)^{\ast}(\mathcal F^{(\pm)}f)(k)dk\Big\| \to 0,
$$
as $t \to \infty$, which implies that
$$
{\mathop{\rm s-lim}_{t\to\pm\infty}}\, e^{it\sqrt{H_+}}e^{-it\sqrt{H_0}}P_{ac}(H_0) = \big(\mathcal F^{(+)}\big)^{\ast}\mathcal F^0 = W_+.
$$
Note that, since $E_H((0,\infty)) = P_{ac}(H)$, we have
\begin{equation}
{\mathop{\rm s-lim}_{t\to\pm\infty}}\, e^{it\sqrt{H_+}}e^{-it\sqrt{H_0}}P_{ac}(H_0) = {\mathop{\rm s-lim}_{t\to\pm\infty}}\, P_{ac}(H)e^{it\sqrt{H}}e^{-it\sqrt{H_0}}P_{ac}(H_0).
\label{C2S8WaveopH+HacH}
\end{equation}
We have thus proven the following theorem.


\begin{theorem}
Let $H$ and $H_0$ be as in Subsection 2.2. Then the wave operator for the wave equation
\begin{equation}
{\mathop{\rm s-lim}_{t\to\pm\infty}}\, e^{it\sqrt{H_+}}e^{-it\sqrt{H_0}}P_{ac}(H_0)
\nonumber
\end{equation}
exists and is equal to the wave operator for the Schr{\"o}dinger equation
\begin{equation}
{\mathop{\rm s-lim}_{t\to\pm\infty}}\,e^{itH}e^{-itH_0}P_{ac}(H_0).
\nonumber
\end{equation}
\end{theorem}

In particular, this theorem implies that the scattering matrix for the Schr{\"o}dinger equation and that for the wave equation coincide.


\chapter{Manifolds with hyperbolic ends}


\section{Classification of 2-dimensional hyperbolic manifolds}
The hyperbolic manifold is, by definition, a complete Riemannian manifold with all sectional curvatures equal to $- 1$. General hyperbolic manifolds are constructed by the action of discrete groups on the upper-half space. The resulting quotient manifold is either compact, or non-compact but of finte volume, or non-compact with infinite volume. In the latter two cases, the manifold can be split into bounded part and unbounded part, this latter being called the end. To study the general structure of ends is beyond our scope. We briefly look at the 2-dimensional case.

\subsection{M{\"o}bius transformation}
Recall that ${\bf C}_+ = \{z = x + iy \, ; y > 0\}$ is a 2-dimensional hyperbolic space equipped with the metric
\begin{equation}
ds^2 = \frac{(dx)^2 + (dy)^2}{y^2}.
\label{eq:Chap3Sec1metric}
\end{equation}
Let $\partial{\bf C}_+ = \partial{\bf H}^2 = \{(x,0) \, ; x \in {\bf R}\}\cup{\infty} = 
{\bf R}\cup{\infty}$.
For a matrix 
\begin{equation}
\gamma = \left(
\begin{array}{cc}
a & b \\
c & d
\end{array}
\right) \in SL(2,{\bf R})
\nonumber
\end{equation}
the M{\"o}bius transformation is defined by
\begin{equation}
{\bf C}_+ \ni z \to \gamma \cdot z:= \frac{az + b}{cz + d},
\label{eq:Chap3Action}
\end{equation}
which is an isometry on ${\bf H}^2$. 
Since $\gamma$ and $-\gamma$ define the same action, one usually identifies them and  considers the factor group:
$$
PSL(2,{\bf R}) := SL(2,{\bf R})/\{\pm I\}.
$$

The non-trivial M{\"o}bius transformations $\gamma$ are classified into 3 categories :
\begin{equation}
\begin{split}
elliptic &\Longleftrightarrow {\rm there \ is \ only \ one \ fixed \ point \ in} \ {\bf C}_+ \\
&\Longleftrightarrow |{\rm tr}\, \gamma| < 2, \\
parabolic &\Longleftrightarrow {\rm there \ is \ only \ one \ degenerate \ fixed \ point \ on} \ \partial{\bf C}_+ \\
&\Longleftrightarrow |{\rm tr}\, \gamma| = 2, \\
hyperbolic &\Longleftrightarrow {\rm there \ are \ two \  fixed \ points \ on} \ \partial{\bf C}_+ \\
&\Longleftrightarrow |{\rm tr}\, \gamma| > 2.
\end{split}
\nonumber
\end{equation}

\subsection{Fuchsian group}
Let $\Gamma$ be a discrete subgroup of $SL(2,{\bf R})$, which is usually called a {\it Fuchsian} group. As a short introduction to the theory of Fuchsian groups, we refer \cite{Kat92}. Let $\mathcal M = \Gamma\backslash{\bf H}^2$ be the fundamental domain by the action (\ref{eq:Chap3Action}). $\Gamma$ is said to be {\it geometrically finite} if $\mathcal M$ is chosen to be a finite-sided convex polygon. The sides are then geodesics of ${\bf H}^2$. The geometric finiteness is equivalent to that $\Gamma$ is finitely generated. 


\subsection{Examples}
As a simple example, consider the cyclic group $\Gamma$ which generates the action $z \to z + 1$. This is parabolic with fixed point $\infty$. The associated fundamental domain is $\mathcal M = (-1/2,1/2]\times(0,\infty)$,
with which one can endow the metric (\ref{eq:Chap3Sec1metric}). 
It has two infinities : $(-1/2,1/2]\times \{0\}$ and $\infty$. The part 
$(-1/2,1/2]\times(0,1)$ has an infinite volume. Let us call it {\it regular infinity} in this note. 
The part $(-1/2,1/2]\times(1,\infty)$ has a finite volume, and is called {\it cusp}. The sides $x = \pm 1/2$ are geodesics.

Another simple example is the cyclic group generated by the hyperbolic action $z \to \lambda z$, $\lambda > 1$. The sides of the fundamental domain 
$\mathcal M = \{1 \leq |z| \leq \lambda\}$ are semi-circles orthogonal to $\{y = 0\}$, which are geodesics. The quotient manifold is diffeomorphic to $S^1\times(-\infty,\infty)$. It is parametrized by $(t,r)$, where $t \in {\bf R}/\log\lambda{\bf Z}$ and $r$ is the signed distance from the segment $\{(0,t)\, ; 1 \leq t \leq \lambda\}$. The metric is then written as
\begin{equation}
ds^2 = (dr)^2 + \cosh^2r\,(dt)^2.
\label{eq:Chap3Sect1Funnelmetric}
\end{equation}
The part $x > 0$ (or $x < 0)$ of $\mathcal M$ is called {\it funnel}. Letting $y = 2e^{-r}$, one can rewrite (\ref{eq:Chap3Sect1Funnelmetric}) as
\begin{equation}
ds^2 = \Big(\frac{dy}{y}\Big)^2 + \Big(\frac{1}{y} + \frac{y}{4}\Big)^2(dt)^2.
\nonumber
\end{equation}
This means that the funnel can be regarded as a perturbation of the regular infinity.

\subsection{Classification}
The set of limit points of a Fuchsian group $\Gamma$, denoted by $\Lambda(\Gamma)$, is defined as follows : $w \in \Lambda(\Gamma)$ if there exist $z_0 \in {\bf C}_+$ and distinct $\gamma_n \in \Gamma$, $n = 1,2,\cdots$,  such that $\gamma_n\cdot z_0 \to w$.
 Since $\Gamma$ acts discontinuously on ${\bf C}_+$, $\Lambda(\Gamma) \subset \partial{\bf H}^2$. There are only 3 possibilities.
\begin{itemize}
\item 
({\it Elementary}) : 
$\Lambda(\Gamma)$ is a finite set.

\item
({\it The 1st kind}) :
$\ \Lambda(\Gamma) = \partial{\bf H}$.

\item
({\it The 2nd kind}) :
$\ \Lambda(\Gamma)$ is a perfect (i.e. every point is an accumulation point), nowhere dense set of $\partial{\bf H}$.
\end{itemize}

If $\Lambda(\Gamma)$ is a finite set, $\Gamma$ is said to be {\it elementary}.  
Any elementary group is either cyclic or is conjugate in $PSL(2,{\bf R})$ to a group generated by $\gamma\cdot z = \lambda z$, $(\lambda > 1)$, and $\gamma'\cdot z = - 1/z$.

 For non-elementary case, we have the following theorem.


\begin{theorem}
Let ${\mathcal M} = \Gamma\backslash{\bf H}^2$ be a non-elementary geometrically finite hyperbolic manifold. Then there exists a compact subset $\mathcal K$ such that $\mathcal M \setminus {\mathcal K}$ is a finite disjoint union of cusps and funnels.
\end{theorem}

For the proof of this theorem, see \cite{Bo07}, p. 27, Theorem 2.13. 

One more explanation is necessary about Theorem 1.1. Let $\Gamma$ be a Fuchsian group. For a point $z_0 \in \overline{{\bf R}^2_+}$, we put
$$
\Gamma_{z_0} = \{\gamma \in \Gamma \, ; \gamma\cdot z_0 = z_0\}.
$$
If $\Gamma_{z_0} \neq \{1\}$, $z_0$ is called a fixed point of $\Gamma$.  A fixed point in ${\bf R}^2_+$ is called an elliptic fixed point. Let $\mathcal M_{sing}$ be the set of elliptic fixed points of $\Gamma$. By a suitable choice of local coordinates,  $\mathcal M = \Gamma\backslash{\bf H^2}$ becomes a Riemann surface, moreover by introducing the metric $y^{-2}\left((dx)^2 + (dy)^2)\right)$, $\mathcal M\setminus{\mathcal M}_{sing}$ is a hyperbolic manifold. However, this metric is singular around the points from $\mathcal M_{sing}$. 
In this case, there exists a neighborhood $U$ of $z_0 \in {\mathcal M}_{sing}$ such that $U = \Gamma_{z_0}\backslash B$, where $B$ is a ball in ${\bf H}^2$.
Then $\mathcal M$ turns out to be an {\it orbifold}. Theorem 1.1 also holds for the orbifold case. However,
in this note, we do not enter into the orbifold structure in detail. The case $\Gamma = SL(2,{\bf Z})$ will be explained in \S 5.


\section{Model space}
By the above classification, it is natural to consider the manifold whose ends are asymptotically equal to either $\mathcal M_{reg} = M\times(0,1)$, or $\mathcal M_{cusp} = M\times(1,\infty)$, where $M$ is a compact manifold, and the metrics of $\mathcal M_{reg}$ and $\mathcal M_{cusp}$ have the form
\begin{equation}
ds^2 = \frac{(dy)^2 + h(x,dx)}{y^2},
\label{eq:Chap3Sec1Modelmetric}
\end{equation}
where $h(x,dx) = \sum_{i,j=1}^{n-1}h_{ij}(x)dx^idx^j$ is the metric on $M$, $x$ being local coordinates on $M$.
Let $\Delta_M$ be the Laplace-Beltrami operator on $M$, $0 = \lambda_0 < \lambda_1 \leq \cdots$ the eigenvalues, and $\varphi_m(x)$, $m = 0,1,2,\cdots$, the associated complete orthonormal system of eigenvectors of $- \Delta_M$. We define for $\phi \in L^2(M)$
\begin{equation}
P_m\phi = (\phi,\varphi_m)_{L^2(M)}\,\varphi_m,
\label{C3S2Pmphi}
\end{equation}
\begin{equation}
\Pi_m\phi = (\phi,\varphi_m)_{L^2(M)}.
\label{C3S2Paimphi}
\end{equation}

We now let $\mathcal M = M\times(0,\infty)$ equipped with the metric (\ref{eq:Chap3Sec1Modelmetric}). 
The Laplace-Beltrami operator on $\mathcal M$ is $ y^2(\partial_y^2 + \Delta_M) - (n - 2)y\partial_y$. We put 
\begin{equation}
 H_{free} = - y^2(\partial_y^2 + \Delta_M) + 
 (n - 2)y\partial_y - \frac{(n-1)^2}{4}=-\Delta_{\mathcal M}- \frac{(n-1)^2}{4}.
 \label{C3S2Hfree}
\end{equation}

Here we need to explain the change of usage of suffix. In Chapters 1 and 2, we used the subscript 0 to denote {\it unperturbed} operators. However, in the sequel, we use the suffix {\it free} for that purpose. The suffix 0 will be used to distinguish the case in which the eigenvalue $\lambda_0 = 0$ is involved.

Spectral properties of $H_{free}$ can be studied in essentially the same way as in Chap. 2. We have only to replace the space $L^2({\bf R}^{n-1})$ by $L^2(M)$ and the Fourier transform by the eigenfunction expansion associated with $- \Delta_M$. The expansion coefficient of $f(x,y)$ is denoted by
\begin{equation}
 {\widehat f}_m(y) = (f(\cdot,y),\varphi_m)_{L^2(M)} = \big(\Pi_m\, f\big)(y).
 \label{eq:Chap3Sect2FourierCoeffi}
\end{equation}
For $f \in C_0^{\infty}(\mathcal M)$, we have
\begin{equation}
\big(\Pi_mH_{free}f\big)(y) = 
  L_{free}(\sqrt{\lambda_m})
\widehat f_m(y),
\nonumber
\end{equation}
where $L_{free}(\zeta)$ is defined by Chap. 1. (\ref{eq:DiffOpL0zeta}). As in Corollary 1.3.10, for $\lambda_m \neq 0$, the Green operator of $L_{free}(\sqrt{\lambda_m}) - \lambda \mp i\epsilon$ is
\begin{equation}
\big(L_{free}(\sqrt{\lambda_m}) - \lambda \mp i\epsilon)\big)^{-1} 
= G_{free}(\sqrt{\lambda_m},\mp i
\sqrt{\lambda \pm i\epsilon}),
\nonumber
\end{equation}
where $G_{free}(\zeta,\nu)$ is defined by Definition 1.3.5.
The Fourier transformation associated with $L_{free}(\sqrt{\lambda_m})$ is given in 
Chap.1, (\ref{eq:Chap1Sec3Fzeta}):
\begin{equation}
\left({ F}_{free,m}\psi\right)(k) = 
\frac{\big(2k\sinh(k\pi)\big)^{1/2}}{\pi}
\int_0^{\infty}y^{(n-1)/2}K_{ik}(\sqrt{\lambda_m}\,y)
\psi(y)\frac{dy}{y^n}.
\label{C3S2Ffreem}
\end{equation}
Letting $\zeta = \sqrt{\lambda_m}$ in Theorem 1.3.13, we obtain the following theorem.


\begin{theorem}  Let $\lambda_m \neq 0$. \\
\noindent
(1) ${F}_{free,m}$ is a unitary operator from $L^2((0,\infty);dy/y^n)$ onto $L^2((0,\infty);dk)$. \\
\noindent
(2) For $\psi \in D(L_{free}(\sqrt{\lambda_m}))$
$$
({F}_{free,m}L_{free}(\sqrt{\lambda_m})\psi)(k) = 
k^2({ F}_{free,m}\psi)(k).
$$
(3) For $\psi \in L^2((0,\infty);dy/y^n)$ the inversion formula holds :
\begin{equation}
\begin{split}
  \psi  & = \big({ F}_{free,m}\big)^{\ast}{ F}_{free,m}\psi \\
   & = y^{(n-1)/2}\int_0^{\infty}
   \frac{(2k\sinh(k\pi))^{1/2}}{\pi}
  K_{ik}(\sqrt{\lambda_m}\,y)
  ({F}_{free,m}\psi)(k)dk.
 \end{split}
 \nonumber
 \end{equation}
\end{theorem}

We consider the case $\lambda_m = 0$, i.e. $m =0$: 
\begin{equation}
L_{free}(0) = - y^2\partial_y^2 + (n-2)y\partial_y - \frac{(n-1)^2}{4}.
\nonumber
\end{equation} 
Since this is Euler's operator, we have
\begin{equation}
 (L_{free}(0) - \lambda \mp i\epsilon))^{-1} = G_{free,0}(\mp i
 \sqrt{\lambda \pm i\epsilon}),
 \nonumber
\end{equation}
\begin{equation}
 G_{free,0}(\nu)\psi(y) = \int_0^{\infty}G_{free,0}(y,y';\nu)\psi(y')\frac{dy'}{(y')^n},
 \label{eq:Chap3Sect1G0nu}
\end{equation}
\begin{equation}
 G_{free,0}(y,y',\nu) = \frac{1}{2\nu}\left\{
  \begin{split}
  y^{\frac{n-1}{2}+\nu}(y')^{\frac{n-1}{2}-{\nu}}, \quad 
   0 < y < y', \\
      y^{\frac{n-1}{2}-\nu}(y')^{\frac{n-1}{2}+{\nu}}, \quad 
   0 < y' < y. 
  \end{split}
 \right.
 \label{eq:Chap3Sect1GreenofG0}
\end{equation}
In the same way as in Lemma 1.3.8, we can prove
\begin{equation}
 \|G_{free,0}(\nu)\psi\|_{{\mathcal B}^{\ast}} \leq \frac{C}{|\nu|}\|\psi\|_{{\mathcal B}},
 \nonumber
\end{equation}
where the constant $C$ is independent of $\nu$. The Fourier transform $F_{free,0}$ associated with $L_{free}(0)$ has 2 components:
\begin{equation}
 F_{free,0} = \big(F_{free,0}^{(+)},F_{free,0}^{(-)}\big),
 \label{C3S2Ffree0}
\end{equation}
\begin{equation}
 \big(F_{free,0}^{(\pm)}\psi\big)(k) = \frac{1}{\sqrt{2\pi}}
 \int_0^{\infty}y^{\frac{n-1}{2} \pm ik} \psi(y)\frac{dy}{y^n}.
 \label{C3S2Ffreeopm}
\end{equation}
Let us check this fact. By (\ref{eq:Chap3Sect1G0nu}), we have for $\psi \in C_0^{\infty}((0,\infty))$ 
\begin{equation}
 G_{free,0}(\mp ik)\psi(y) \sim 
 \pm \frac{i}{k}\sqrt{\frac{\pi}{2}}\left\{
 \begin{split}
\,y^{\frac{n-1}{2}\mp ik}
 F_{free,0}^{(\pm)}(k)\psi, \quad 
 y \to 0, \\
\,y^{\frac{n-1}{2}\pm ik}
   F_{free,0}^{(\mp)}(k)\psi, \quad 
 y \to \infty.
 \end{split}
 \right.
 \nonumber
\end{equation}
On the other hand, we have
\begin{eqnarray*}
& & \frac{1}{2\pi i}\big(G_{free,0}(-ik) - G_{free,0}(ik)\big)\psi \\
&=& \frac{1}{4\pi k}
\int_0^{\infty}(yy')^{\frac{n-1}{2}}\left\{
\big(\frac{y'}{y}\big)^{ik}
+ \big(\frac{y}{y'}\big)^{ik}
\right\}\psi(y')\frac{dy'}{(y')^n} \\
&=& \frac{1}{2k\sqrt{2\pi}}\left(
y^{\frac{n-1}{2} - ik}F_{free,0}^{(+)}(k)\psi + 
 y^{\frac{n-1}{2} + ik}F_{free,0}^{(-)}(k)\psi
 \right).
\end{eqnarray*}
Hence we have
\begin{equation}
\frac{1}{2\pi i}\left(\big[G_{free,0}(-ik) - 
G_{free,0}(ik)\big]\psi,\psi
\right) = \frac{1}{2k}\left|(F_{free,0}\psi)(k)\right|^2.
\nonumber
\end{equation}
Integrating this equality and arguing as in Chap. 1, \S 3, we obtain the following Theorem 2.2. Alternatively, one can use the fact that
$$
 (F_{free,0}\psi)(k) = \left(
 \widetilde \psi(-k), \widetilde \psi(k) \right),
$$
where 
$\widetilde \psi$ is the Fourier transform of $U\psi(t) = e^{-(n-1)t/2}\psi(e^t)$. In fact,  $U$ is unitary from $L^2((0,\infty);dy/y^n)$ to  
$L^2({\bf R};dt)$, and we have
\begin{equation}
U\left(- y^2\partial_y^2 + (n-2)y\partial_y - \frac{(n-1)^2}{4}
\right)U^{\ast} = 
- \partial_t^2.
\label{eq:Chap3Sect1Eulerunitrary}
\end{equation}


\begin{theorem}
(1) $F_{free,0} : L^2((0,\infty);dy/y^n) \to (L^2((0,\infty);dk))^2$ is unitary. \\
\noindent
(2) For $f \in D(L_{free,0}(0))$,
\begin{equation}
 (F_{free,0}L_{free,0}(0)f)(k) =  k^2 (F_{free,0}f)(k).
 \nonumber
\end{equation}
(3) For $f \in L^2((0,\infty);dy/y^n)$, the inversion formula holds:
\begin{equation}
\begin{split}
 f &= \left(F_{free,0}\right)^{\ast}F_{free,0}f \\
 & = \frac{1}{\sqrt{2\pi}}\int_0^{\infty}y^{(n-1)/2}
\left(y^{-ik}
 F_{free,0}^{(+)}(k)f + y^{ik}
 F_{free,0}^{(-)}(k)f\right) dk.
\end{split}
\nonumber
\end{equation}
\end{theorem}

We now return to the operator $H_{free}$ whose resolvent is written as
\begin{equation}
 (H_{free} - \lambda \mp i0)^{-1}f 
 = 
\sum_{m=0}^{\infty}\varphi_m(x)\left(
G_{free}(\sqrt{\lambda_m},\mp i\sqrt{\lambda})\widehat f_m\right)(y).
\label{eq:Chap3Sect1ResolventofH0}
\end{equation}
Here $G_{free}(\sqrt{\lambda_0},\mp i\sqrt{\lambda}) = G_{free,0}(\mp i\sqrt{\lambda})$.
Repeating the proof of Lemma 1.4.1, we can show the following lemma.


\begin{lemma}
$H_{free}\Big|_{C_0^{\infty}(\Omega)}$ is essentially self-adjoint.
\end{lemma}

Recall that the generalized Fourier transform is derived from the asymptotic behavior of the resolvent at infinity. For $M\times(0,\infty)$, there are two infinities ; $y = 0$ and $y = \infty$, the former corresponding to the regular infinity, the latter to the cusp. We put the suffix {\it reg} or {\it c} for the Fourier transforms associated with regular infinity or cusp.


\begin{definition} 
Let ${\mathcal D}(M\times(0,\infty))$ be the set of functions $f(x,y) \in C^{\infty}(M\times(0,\infty))$ such that $\widehat f_m \in C_0^{\infty}((0,\infty))$, moreover $\widehat f_m = 0$ except for a finite number of $m$. We put
\begin{equation}
{\bf h} = L^2(M) \oplus {\bf C},\quad 
\widehat{\mathcal H} = L^2((0,\infty) ; {\bf h} ; dk),
\nonumber
\end{equation} 
\begin{equation}
{\mathcal F}_{free}^{(\pm)}
 = \left({\mathcal F}_{free, reg}^{(\pm)},{\mathcal F}_{free, c}^{(\pm)}\right), \nonumber
\end{equation}
and define on ${\mathcal D}(M\times(0,\infty))$
\begin{equation}
 {\mathcal F}_{free, reg}^{(\pm)}
  =  \sum_{m=0}^{\infty}
  C_{m}^{(\pm)}(k)P_m\otimes
F_{free,m}^{(\pm)},
\label{eq:Chap3Sect2F0regpm}
\end{equation}
\begin{equation}
F_{free,m}^{(\pm)} = 
\left\{
\begin{array}{cc}
F_{free,m} & (\lambda_m \neq 0) \\
F_{free,0}^{(\pm)} & (\lambda_m = 0),
\end{array}
\right.
\label{eq:Chap3Sect2F0gammaastpm}
\end{equation}
\begin{equation}
C_{m}^{(\pm)}(k) = 
\left\{
\begin{split}
& \left(\frac{\sqrt{\lambda_m}}{2}\right)^{\mp ik} \quad (\lambda_m \neq 0) \\
& \dfrac{\pm i}{k\omega_{\pm}(k)}\sqrt{\dfrac{\pi}{2}} \quad (\lambda_m = 0),
\end{split}
\right.
\label{eq:Chap3Sect2constCgammaastk}
\end{equation}
\begin{equation}
{\mathcal F}_{free, c}^{(\pm)} = 
P_{0}\otimes F_{free,0}^{(\mp)}.
\label{eq:Chap3Sect2Focplusminus}
\end{equation}
\end{definition}

We define ${\mathcal B}, {\mathcal B}^{\ast}$, and $L^{2,s}$ by putting ${\bf h} = L^2(M)\oplus{\bf C}$ in Chap. 1, \S 2. Note that, geometrically, $\mathcal B$ corresponds to the diadic decomposition with respect to the geodesic distance, and ${\mathcal B}^{\ast}$ to the integral mean over the geodesic ball. Let
\begin{equation}
  R_{free}(z) = (H_{free} - z)^{-1}.
  \nonumber
\end{equation}
Then Theorem 2.1.3 remains valid for $H_{free}$ if  $\mathcal X^s$ is replaced by $L^{2,s}$.


\begin{theorem} 
(1) $\sigma(H_{free}) = [0,\infty)$. \\
\noindent
(2) $\sigma_p(H_{free}) = \emptyset$. \\
\noindent
(3) For $\lambda > 0$ and $f,g \in {\mathcal B}$, the following weak limit exists
$$
\lim_{\epsilon \to 0}(R_{free}(\lambda \pm i\epsilon)f,g) =: 
(R_{free}(\lambda \pm i0)f,g).
$$
Moreover
\begin{equation}
 \|R_{free}(\lambda \pm i0)f\|_{{\mathcal B}^{\ast}} \leq C\|f\|_{\mathcal B},
 \nonumber
\end{equation}
where the constant $C$ does not depend on $\lambda$ if $\lambda$ varies over a compact set in $(0,\infty)$. \\
\noindent
(4) Letting ${\mathcal F}_{free}^{(\pm)}(k)f = 
({\mathcal F}_{free}^{(\pm)}f)(k)$ for $f \in {\mathcal D}(M\times(0,\infty))$, we have
\begin{equation}
 \|{\mathcal F}_{free}^{(\pm)}(k)f\|_{{\bf h}} \leq 
 C\|f\|_{\mathcal B},
 \nonumber
\end{equation}
where the constant $C$ does not depend on $k$ if $k$ varies over a compact set in $(0,\infty)$.\\
\noindent
(5) ${\mathcal F}_{free}^{(\pm)}$ is uniquely extended to a unitary operator from $L^2(M\times(0,\infty);\sqrt{g_M}dxdy/y^n)$ to
$\widehat{\mathcal H}$. Moreover if $f \in D(H_{free})$
\begin{equation}
({\mathcal F}_{free}^{(\pm)}H_{free}f)(k) = k^2
({\mathcal F}_{free}^{(\pm)}f)(k).
\nonumber
\end{equation}
\end{theorem}
Proof. The assertions (1), (2) follow from Lemma 1.3.2. Note that $L_{free}(0)$ should be treated separately, however, it is easy by (\ref{eq:Chap3Sect1Eulerunitrary}). The proof of (3) is almost the same as Theorem 2.2.3 (2), (3), the term $L_{free}(0)$ requires a small change, though. In the next section, we shall give the proof for the more general case (see Theorem 3.8). Applying Stone's formulas for each $L_{free}(\sqrt{\lambda_m})$, we have
\begin{equation}
 \frac{1}{2\pi i}\left([R_{free}(\lambda + i0) - 
 R_{free}(\lambda - i0)]f,f\right) = 
 \|{\mathcal F}_{free}^{(\pm)}(k)f\|^2,
  \nonumber
\end{equation}
which implies (4). Since each
${\mathcal F}_{free,m}$ is unitary, (5) follows. \qed

\medskip
The relation of $\mathcal F_{free}^{(\pm)}$ and the asymptotic behavior of the resolvent is as follows.

\begin{theorem}
For $k > 0$ and $f \in {\mathcal B}$, we have
\begin{equation}
\lim_{R\to\infty}\frac{1}{\log R}\int_{1/R<y<1}
\|R_{free}(k^2 \pm i0)f - v^{(\pm)}_{reg}\|_{L^2(M)}^2\frac{dy}{y^n} = 0,
\label{eq:Chap3Sect1R0nearregular}
\end{equation}
\begin{equation}
v^{(\pm)}_{reg} = \omega_{\pm}(k)\, y^{(n-1)/2  \mp ik}
{\mathcal F}^{(\pm)}_{free, reg}(k)f,
\nonumber
\end{equation}
\begin{equation}
\lim_{R\to\infty}\frac{1}{\log R}\int_{1<y<R}
\|R_{free}(k^2 \pm i0)f - v_{c}^{(\pm)}\|_{L^2(M)}^2\frac{dy}{y^n} = 0,
\label{eq:Chap3Sect1R0nearcusp}
\end{equation}
\begin{equation}
v_{c}^{(\pm)} = \omega_{\pm}^{(c)}(k)\,
y^{(n-1)/2 \pm ik}
{\mathcal F}_{free, c}^{(\pm)}(k)f.
\nonumber
\end{equation}
Here $\omega_{\pm}(k)$ is defined by Chap. 1 (\ref{eq:Chap1Sect4omegaplusminusk}), and
\begin{equation}
\omega_{\pm}^{(c)}(k) = \pm \frac{i}{k}\sqrt{\frac{\pi}{2}}.
\nonumber
\end{equation}
\end{theorem}

Proof. By Theorem 2.5(3) and (4),
we have only to prove the theorem for $f \in {\mathcal D}(M\times(0,\infty))$. Assume that
$f = 0$ for $y < \epsilon$ and $y > 1/\epsilon$. Then if $y < \epsilon$, we have by (\ref{eq:Chap3Sect1ResolventofH0}), (\ref{eq:Chap3Sect1GreenofG0}) and Chap.1 Definition 3.5
\begin{equation}
\begin{split}
  & R_{free}(k^2 \pm i0)f\\
 =\,&  \pm \frac{i}{k}\sqrt{\frac{\pi}{2}}\frac{1}{\sqrt{|M|}}
 y^{(n-1)/2 \mp ik}F_{free,0}^{(\mp)}(k)\widehat f_0 \\
 &   + \frac{\pi}{\big(2k\sinh(k\pi)\big)^{1/2}}
 \sum_{m \geq 1}\varphi_m(x)
y^{(n-1)/2}I_{\mp ik}(\sqrt{\lambda_m}y)
F_{free,m}(k)\widehat f_m.
\end{split}
\nonumber
\end{equation}
Using Definition 2.4 and Chap. 1 (\ref{eq:Chap1Sec3Inunear0}), we obtain (\ref{eq:Chap3Sect1R0nearregular}). 

For $y > 1/\epsilon$, we have by using Chap. 1 (\ref{eq:Chap1Sec3G02})
\begin{eqnarray*}
& & \|R_{free}(k^2 \pm i0)f - \frac{1}{\sqrt{|M|}}
G_{free,0}(\mp ik)\widehat f_0\|_{L^2(M)}^2 \\
&\leq& Cy^{n-2}\sum_{m\geq 1}
\left(\int_0^{\infty}|\widehat f_m(y)|\frac{dy}{y^{(n+2)/2}}\right)^2,
\end{eqnarray*}
which proves (\ref{eq:Chap3Sect1R0nearcusp}). \qed


\section{Manifolds with hyperbolic ends}


\subsection{The formula of Helffer-Sj{\"o}strand}
We prepare a useful tool from functional analysis introduced by Helffer-Sj{\"o}strand \cite{HeSj89}. Let $\sigma \in {\bf R}$, and suppose $f(t) \in C^{\infty}({\bf R})$ satisfies
\begin{equation}
|f^{(k)}(t)| \leq C_k(1 + |t|)^{\sigma-k}, \quad \forall k, \quad \forall t \in {\bf R}.
\label{C3S3HelfferSjostland}
\end{equation}
Then there exists $F(z) \in C^{\infty}({\bf C})$ such that
\begin{equation}
\left\{
\begin{split}
& F(t) = f(t), \quad t \in {\bf R}, \\
& |F(z)| \leq C(1 + |z|)^{\sigma}, \\
& |\overline{\partial_z}F(z)| \leq C_n|{\rm Im}\,z|^n(1 + |z|)^{\sigma-n-1}, \quad \forall n, \\
& {\rm supp}\,F(z) \subset \{|{\rm Im}\,z| \leq 2 + 2|{\rm Re}\,z|\}.
\end{split}
\right.
\label{C3S3Fzestimate}
\end{equation}
Here $\overline{\partial_z} = \frac{1}{2}(\partial_x + i\partial_y)$.
This function $F$ is called an {\it almost analytic extension} of $f$. If $f \in C_0^{\infty}({\bf R})$, we can construct $F(z) \in C_0^{\infty}({\bf C})$.

Let us explain the idea of the proof. For $z \in {\bf C}$, let $\langle z\rangle = (1 + |z|^2)^{1/2}$. Take $\chi(y) \in C_0^{\infty}({\bf R})$ such that $\chi(y) = 1 \ (|y| < 1)$, $\chi(y) = 0 \ (|y| > 2)$, and put
$$
F(z) = \sum_{n=0}^{N-1}\frac{i^n}{n!}f^{(n)}(x)y^n\chi\Big(\frac{y}{\langle x\rangle}\Big).
$$
Then we have
\begin{equation}
\begin{split}
2\overline{\partial_z}F(z) = & \frac{i^{N-1}}{(N-1)!}f^{(N)}(x)y^{N-1}\chi\Big(\frac{y}{\langle x\rangle}\Big) \\
& + \sum_{n=0}^{N-1} \frac{i^n}{n!}
f^{(n)}(x)y^n\chi'\Big(\frac{y}{\langle x\rangle}\Big)
\Big(\frac{i}{\langle x\rangle} - \frac{xy}{\langle x\rangle^3}\Big).
\end{split}
\nonumber
\end{equation}
On the support of the first term of the right-hand side, $|y| \leq 2\langle x\rangle$. Hence for $1 \leq n \leq N-1$, it is dominated by $C\langle x\rangle^{\sigma-N}|y|^{N-1} \leq C|y|^n\langle z\rangle^{\sigma-n-1}$. On the support of the 2nd term, $\langle x\rangle \leq |y| \leq 2\langle x\rangle$. Hence, it is dominated by 
$$
C\sum_{n=0}^{N-1} \frac{1}{n!}
\langle x\rangle^{\sigma-n-1}|y|^n\Big|\chi'\Big(\frac{y}{\langle x\rangle}\Big)\Big| \leq 
C \langle x\rangle^{\sigma-1}\exp{\frac{|y|}{\langle x\rangle}} \leq
C_n |y|^n \langle z\rangle^{\sigma-n-1}.
$$
Hence, $|\overline{\partial_z}F(z)| \leq C_n|{\rm Im}\,z|^n(1 + |z|)^{\sigma-n-1}$ holds for $1 \leq n \leq N-1$. By the similar computation, one can show $|F(z)| \leq C(1 + |z|)^{\sigma}$. For the general construction of $F(z)$, see e.g. \cite{Is04a} p. 363.


\begin{lemma}
Let $f(t)$ and $F(z)$ be as above. Suppose $\sigma < 0$. Then for any self-adjoint operator $A$, the following formula holds
$$
f(A) = \frac{1}{2\pi i}\int_{{\bf C}}\overline{\partial_z}F(z)(z - A)^{-1}dzd\overline{z}.
$$
\end{lemma}

Proof. For $\lambda \in {\bf R}$, we have by the generalized Cauchy formula
$$
F(\lambda) = \frac{1}{2\pi i}\int_{|z|=R}\frac{F(z)}{z-\lambda}dz + 
 \frac{1}{2\pi i}\int_{|z|<R}\frac{\overline{\partial_z}F(z)}{z-\lambda}dzd\overline{z}.
$$
Letting $R \to \infty$, we have
$$
F(\lambda) =  
 \frac{1}{2\pi i}\int_{{\bf C}}\frac{\overline{\partial_z}F(z)}{z-\lambda}dzd\overline{z},
$$
where the integral is absolutely convergent.
Let $E(\lambda)$ be the spectral decomposition of $A$. Then we have
\begin{equation}
\begin{split}
f(A) &= \int_{-\infty}^{\infty}f(\lambda)dE(\lambda) \\
& = \frac{1}{2\pi i}\int_{-\infty}^{\infty}\int_{{\bf C}}\frac{\overline{\partial_z}F(z)}{z-\lambda}dzd\overline{z}dE(\lambda) \\
& = \frac{1}{2\pi i}\int_{{\bf C}}\overline{\partial_z}F(z)(z - A)^{-1}dzd\overline{z}. \qed
\end{split}
\nonumber
\end{equation}

Let us mention here useful formulas to compute the commutator of functions of self-adjont operators.
For two operators $P, A$, we put
$$
{\rm ad}_0(P,A) = P,
$$
$$
{\rm ad}_n(P,A) = [{\rm ad}_{n-1}(P,A),A], \quad \forall n \geq 1.
$$
If $A$ is self-adjoint and $f(s)$ satisfies $|f^{(k)}(s)| \leq C_k(1 + |s|)^{\sigma-k}, \ \forall k \geq 0$, we have
\begin{equation}
[P,f(A)] = \sum_{k=1}^{n-1}\frac{(-1)^{k-1}}{k!}{\rm ad}_k(P,A)f^{(k)}(A) + R_{n,l},
\label{C3S3[P,f(A)]left}
\end{equation}
\begin{equation}
R_{n,l} = \frac{1}{2\pi i}\int_{\bf C}\overline{\partial_z}F(z)(A-z)^{-1}
{\rm ad}_n(P,A)(A-z)^{-n}dzd\overline{z}.
\label{C3S3Rnl}
\end{equation}
\begin{equation}
[P,f(A)] = \sum_{k=1}^{n-1}\frac{1}{k!}f^{(k)}(A){\rm ad}_k(P,A) + R_{n,r},
\label{C3S3[P,f(A)]right}
\end{equation}
\begin{equation}
R_{n,r} = \frac{(-1)^{(n+1}}{2\pi i}\int_{\bf C}\overline{\partial_z}F(z)(A-z)^{-n}
{\rm ad}_n(P,A)(A-z)^{-1}dzd\overline{z}.
\label{C3S3Rnr}
\end{equation}
Here, $F(z)$ is an almost analytic extension of $f$, and we assume that
$$
\|(A-z)^{-n}{\rm ad}_n(P,A)(A-z)^{-1}\| \leq C|{\rm Im}\,z|^{-n-1}\langle z\rangle^{\mu(n)},
$$
$$
\sigma - n + \mu(n) < 0,
$$
in order to guarantee the convergence of the integrals (\ref{C3S3Rnl}), (\ref{C3S3Rnr}). Formal derivation of (\ref{C3S3[P,f(A)]left}), (\ref{C3S3[P,f(A)]right}) is rather easy. However, rigorous derivation requires examination of the domain of ${\rm ad}_n(P,A)$. When $P$ and $A$ are differential operators, this domain question boils down to  the regularity estimate for $(A-z)^{-1}$.


\subsection{Assumptions on ends}
Now we consider an $n$-dimensional connected Riemannian manifold ${\mathcal M}$, which is written as a union of open sets:
\begin{equation}
{\mathcal M} = {\mathcal K}\cup{\mathcal M}_1\cup\cdots\cup{\mathcal M}_N.
\nonumber
\end{equation}
We assume that

\medskip
\noindent
{\bf(A-1)} $\ \ \overline{\mathcal K}$ {\it is compact}.

\medskip
\noindent
{\bf (A-2)} $\ \ {\mathcal M}_p\cap{\mathcal M}_q = \emptyset, \quad p \neq q$.

\medskip
\noindent
{\bf (A-3)} $\ $ {\it Each} ${\mathcal M}_p$, $p = 1, \cdots, N$, {\it is diffeomorphic either to} ${\mathcal M}_{reg} = M_p\times(0,1)$ {\it or to} ${\mathcal M}_{c} = M_p\times(1,\infty)$, $M_p$ being a compact Riemannian manifold of dimension $n-1$, {\it which is allowed to be different for each} $p$.

\medskip
\noindent
{\bf (A-4)} {\it On each $\mathcal M_{p}$, the Riemannian metric $ds^2$ has the following form
\begin{equation}
ds^2 = y^{-2}\left((dy)^2 + h_p(x,dx) + A_p(x,y,dx,dy)\right),
\label{C3S3ds2onMp}
\end{equation}
\begin{equation}
A_p(x,y,dx,dy) = \sum_{i,j=1}^{n-1}a_{p,ij}(x,y)dx^idx^j + 2\sum_{i=1}^{n-1}a_{p,in}(x,y)dx^idy + a_{p,nn}(x,y)(dy)^2,
\nonumber
\end{equation}
where $h_p(x,dx) = \sum_{i,j=1}^{n-1}h_{p,ij}(x)dx^idx^j$ is a positive definite metric on $M_p$, 
and $a_{p,ij}(x,y), 1 \leq i,j \leq n$, satisfies the following condition
\begin{equation}
|\widetilde D_x^{\alpha}D_y^{\beta}\, a(x,y)| \leq C_{\alpha\beta}(1 + |\log y|)^{-{\rm min}(|\alpha|+\beta,1)-1-\epsilon}, \quad \forall \alpha, \beta
\label{eq:Chap3Sect2DecayMetric1}
\end{equation}
for some $\epsilon > 0$.
Here $D_y = y \partial_y$, $\widetilde D_x = \tilde y(y)\partial_x$, $\tilde y(y) \in C^{\infty}((0,\infty))$ such that $\tilde y(y) = y$ for $y > 2$ and $\tilde y(y) = 1$ for $0 < y < 1$.}

\medskip
Following Example 1.3, we call $\mathcal M_p = M_p\times(0,1)$ a {\it regular end} and ${\mathcal M}_p = M_p\times(1,\infty)$ a {\it cusp}.

Let us note that the above model in particular contains ${\bf H}^n$. In fact, we take $\mathcal K = B_2(0,1)$, and $\mathcal M_1 = {\bf H}^n\setminus B_{\log2}(0,1)$, where $B_r(0,1)$ is the geodesic ball of radius $r$ centered at $(0,1)$.
Using geodesic polar coordinates, $\mathcal M_1$ is isometric to $S^{n-1}\times(\log 2,\infty)$ equipped with the metric $(dr)^2 + \sinh^2r(d\theta)^2$. Taking $e^r = 2/y$, we see that $\mathcal M_1 = \mathcal M_{reg} = S^{n-1}\times(0,1)$ equiped with the metric $y^{-2}\Big((dy)^2 + (d\theta)^2 + (y^4/16 - y^2/2)(d\theta)^2\Big)$.

The 2nd important remark is that, if $\mathcal M_p$ is equal to $\mathcal M_{reg}$, one can assume that the above metric (\ref{C3S3ds2onMp}) takes the form
\begin{equation}
ds^2 = y^{-2}\Big((dy)^2 + h_p(x,dx) +  \sum_{i,j=1}^{n-1}a_{p,ij}(x,y)dx^idx^j\Big)
\label{eq:Chap3Sect2WarpedProduct}
\end{equation}
and each $a_{p,ij}(x,y)$ satisfies the condition (\ref{eq:Chap3Sect2DecayMetric1}). This can be proved in the same way as in Theorem 4.1.6 to be given in Chap. 4. Therefore in the following we consider the metric of the form (\ref{eq:Chap3Sect2WarpedProduct}) for such ends.

Let $\Delta_g$ be the Laplace-Beltrami operator on ${\mathcal M}$. 
As has been discussed in Chap. 2, \S 2, we pass to the gauge transformation
\begin{equation}
- \Delta_g - \frac{(n-1)^2}{4} \to H =:- \rho^{1/4}\Delta_g\rho^{-1/4} - \frac{(n-1)^2}{4},
\label{C3S3gaugetransf}
\end{equation}
where $\rho \in C^{\infty}(\mathcal M)$ is a positive function such that
on each end $\mathcal M_p$
\begin{equation}
\rho = \det g^{(p)}/\det g^{(p)}_{free},
\label{C3S3rhodefine}
\end{equation}
$g^{(p)}_{free}$ and $g^{(p)}$ being the unperturbed and perturbed metrics
\begin{equation}
g^{(p)}_{free} = y^{-2}\left((dy)^2 + h_p(x,dx)\right),
\label{C3S3gpfree}
\end{equation}
\begin{equation}
g^{(p)} = y^{-2}\left((dy)^2 + h_p(x,dx) + A_p(x,y;dx, dy)\right)
\label{C3S3gp}
\end{equation}
satisfying the above assumptions. Then $H$ is written as
\begin{equation}
H = - \Delta_g + L_2 - \frac{(n-1)^2}{4},
\label{C3S3L2define}
\end{equation}
$L_2$ being a 2nd order differential operator on $\mathcal M$, and satisfies the following conditions. 

\medskip
\noindent
{\bf (A-5)} {\it $H$ is formally self-adjoint. Namely,
$$
(H\varphi,\psi) = (\varphi,H\psi), \quad \forall \varphi,\psi \in C_0^{\infty}({\mathcal M}),
$$
where $(\,\, ,\,)$ is the inner product of $L^2({\mathcal M})$, i,e, 
$$
(\varphi,\psi) = \int_{\mathcal M}\varphi\overline{\psi}d\mathcal M,
$$
$d\mathcal M$ being the measure which coincides with the unperturbed metric on each 
$\mathcal M_p$.}

\noindent
{\bf (A-6)} {\it $L_2$ is short-range on each ${\mathcal M}_p \ (1 \leq p \leq N)$. Namely, if $L_2$ is represented as
$$
L_1 = \sum_{|\alpha|\leq 2}a_{\alpha}(x,y)D^{\alpha}, \quad
D = (D_x,D_y) = (y\partial_x,y\partial_y),
$$
there exists a constant $\epsilon > 0$ such that}
\begin{equation}
|\widetilde D_x^{\beta}D_y^{k}a_{\alpha}(x,y)| \leq C_{\beta,k}(1 + |\log y|)^{-|\beta|-k-1-\epsilon},  \quad \forall \beta, \quad \forall k.
\nonumber
\end{equation}

\medskip

We use the following partition of unity. Fix $x_0 \in {\mathcal K}$ arbitrarily, and pick $\chi_0 \in C_0^{\infty}({\mathcal M}),$ such that
$$
\chi_0(x) = \left\{
\begin{array}{ll}
1, & {\rm dist}\,(x,x_0) < R, \\
0, & {\rm dist}\,(x,x_0) > R + 1,
\end{array}
\right.
$$
where ${\rm dist}(x,x_0)$ is the distance between $x$ and $x_0$.
Taking $R$ large enough, we define $\chi_j \in C^{\infty}(\mathcal M)\, j=1, \dots, N,$ such that$$
\chi_j(x) = \left\{
\begin{array}{ll}
1 - \chi_0(x), & x \in {\mathcal M}_j, \\
0, & x \notin {\mathcal M}_j.
\end{array}
\right.
$$
Then we have
\begin{equation}
\left\{
\begin{array}{l}
\sum_{j=0}^N\chi_j = 1, \\
{\rm supp}\,\chi_j \subset {\mathcal M}_j, \quad 1 \leq j \leq N, \\
\chi_0 = 1 \quad {\rm on} \quad {\mathcal K}.
\end{array}
\right.
\label{eq:Chap3Sec2Partitionofunity}
\end{equation}
For $1 \leq j \leq N$, we construct $\widetilde\chi_j \in C^{\infty}({\mathcal M})$ such that
\begin{equation}
\quad 
{\rm supp}\,\widetilde\chi_j \subset {\mathcal M}_j, 
\quad 
\widetilde\chi_j = 1 \quad {\rm on} \quad
{\rm supp}\,\chi_j.
\nonumber
\end{equation}


\begin{theorem}
(1) $H\big|_{C_0^{\infty}({\mathcal M})}$ is essentially self-adjoint. \\
\noindent
(2) $\sigma_{e}(H) = [0,\infty)$. 
 \end{theorem}
 
Proof. To prove assertion (1), we first observe that Theorem 2.1.3(4) and (6) remain valid
for $H$, if we substitute the spaces $\mathcal X^s$ with 
$$
L^{2, s}=\{U \in L^2_{loc}:\,\, \int_{\mathcal M}
\left(1+\log^2{(d(x, x_0))}  \right)^{s}  |u(x)|^2 < \infty \}.
$$
Using this analog of Theorem 2.1.3 (4), assertion (1) is proven in the same way as in Theorem 2.1.4.
 
 To show (2), we derive a formula for the resolvent by using the partition of unity (\ref{eq:Chap3Sec2Partitionofunity}). Recall that ${\mathcal M}_j$ is diffeomorphic to $M_j\times(0,1)$ or $M_j\times(1,\infty)$. Let $H_{free(j)}$ be defined by (\ref{C3S2Hfree}) with $M$ replaced by $M_j$, and put
\begin{equation}
R(z) = (H -z)^{-1}, \quad 
R_{free(j)}(z) = (H_{free(j)} - z)^{-1}.
\label{C3S3RzRfreejz}
\end{equation} 
Note that we are using the suffix {\it free(j)} to specify unperturbed operators with respect to the model space $M_j\times (0,\infty)$.
Since
$$
(H - z)\chi_jR_{free(j)}(z)\widetilde\chi_j = \chi_j + \chi_j(H - H_{free(j)})R_{free(j)}(z)\widetilde\chi_j + 
[H,\chi_j]R_{free(j)}(z)\widetilde\chi_j,
$$
we have
\begin{equation}
\chi_jR_{free(j)}(z)\widetilde\chi_j = R(z)\chi_j + R(z)A_j(z)\widetilde\chi_j,
\nonumber
\end{equation}
\begin{equation}
 A_j(z) = [H,\chi_j]R_{free(j)}(z) + \chi_j(H - H_{free(j)}) {\widetilde \chi}_j R_{free(j)}(z).
\nonumber
\end{equation}
Letting
\begin{equation}
A(z) = \sum_{j=1}^NA_j(z)\widetilde\chi_j,
\end{equation}
we then have
\begin{equation}
R(z) = \sum_{j=1}^N\chi_jR_{free(j)}(z)\widetilde\chi_j + R(z)(\chi_0 - A(z)).
\nonumber
\end{equation}
By the assumption (A-4), $R(z)(\chi_0 - A(z))$ is compact. Indeed, for $z \not\in {\bf R}$,  $A_j(z)$ is bounded from $W^{2, 2}({\mathcal M})$ to $L^{2,s}$ with $0 < s < 1 + \epsilon$. Since $R(z)$ is locally smoothening, this implies the desired compactness if one considers the adjoint $(A(z)^{\ast} - \chi_0)R(z)^{\ast}$.

To prove (2), we first show $(- \infty,0) \subset \sigma_d(H)$. It is sufficient to prove that $f(H)$ is compact for any $f \in C_0^{\infty}((-\infty,0))$. Let $F$ be an almost analytic extension of $f$. Then, by Lemma 3.1, we have  
$$
f(H) = \sum_{j=1}^N\chi_jf(H_{free(j)}) {\widetilde \chi}_j- K,
$$
$$
K = \frac{1}{2\pi i}\int_{{\bf C}}\overline{\partial_{z}}F(z)R(z)\left(\chi_0 - A(z)\right)dzd\overline{z}.
$$
Note that $K$ is compact, since $|\partial_{\overline z} F(z)| \leq C_l(1+|z|)^{-l}$, for all $l >0$, and
so is $R(z)(\chi_0 - A(z))$. Since $\sigma(H_{free(j)}) = [0,\infty)$, we have $f(H_{free(j)}) = 0$. Therefore $f(H)$ is compact, which proves
 $\sigma_e(H) \subset [0,\infty)$. The converse inclusion relation is proven by Weyl's method of singular sequence as in Lemma 1.3.12. \qed


\subsection{Limiting absorption principle}　

\begin{lemma}
Let $f(x) \in L^1(0,\infty;dx)$ and put
$$
u(x) = \int_x^{\infty}f(t)dt.
$$
Then for $s >1/2$
$$
\int_0^{\infty}x^{2(s-1)}|u(x)|^2dx \leq 
\frac{4}{(2s-1)^2}\int_0^{\infty}x^{2s}|f(x)|^2dx.
$$
\end{lemma}
Proof. We use the following inequality of Hardy : For
$p > 1, \ g(x) \in L^1(0,\infty)$, we put
$$
F(x) = \int_x^{\infty}g(t)dt.
$$
Then we have
\begin{equation}
\int_0^{\infty}|F(x)|^pdx \leq p^p\int_0^{\infty}|x g(x)|^pdx
\nonumber
\end{equation}
(\cite{HLP52}, p. 244). Letting $\epsilon = 2s - 1 > 0, \ y = x^{\epsilon}$ for $u(x)$ in the Lemma, we  have
$$
(2s - 1)\int_0^{\infty}x^{2(s-1)}|u(x)|^2dx = \int_0^{\infty}
|u(y^{1/\epsilon})|^2dy,
$$
$$
u(y^{1/\epsilon}) = \frac{1}{\epsilon}\int_y^{\infty}
f(z^{1/\epsilon})z^{1/\epsilon - 1}dz.
$$
By Hardy's inequality, with $g(z)= \frac{1}{\epsilon} f(z^{1/\epsilon}) z{(1-\epsilon)/\epsilon}$
and $p=2$,
\begin{eqnarray*}
\int_0^{\infty}|u(y^{1/\epsilon})|^2dy &\leq& 
\frac{4}{\epsilon^2}\int_0^{\infty}|f(y^{1/\epsilon})|^2y^{2/\epsilon}dy \\
&=& \frac{4}{\epsilon}\int_0^{\infty}|f(x)|^2x^{2s}dx,
\end{eqnarray*}
which implies the Lemma. \qed

\bigskip
On each end $\mathcal M_j$ of ${\mathcal M}$, the spaces $L^{2,s}, \ {\mathcal B}, \ 
{\mathcal B}^{\ast}$ are defined in the same way as before with ${\bf h} = L^2(M_j)$. Using the above partition of unity $\chi_j$, we put
\begin{equation}
\|u\|_s = \|\chi_0u\|_{L^2} + \sum_{j=1}^N\|\chi_ju\|_s,
\nonumber
\end{equation}
\begin{equation}
\|u\|_{\mathcal B} = \|\chi_0u\|_{L^2} + \sum_{j=1}^N\|\chi_ju\|_{\mathcal B},
\nonumber
\end{equation}
\begin{equation}
\|u\|_{{\mathcal B}^{\ast}} = \|\chi_0u\|_{L^2} + \sum_{j=1}^N\|\chi_ju\|_{{\mathcal B}^{\ast}},
\nonumber
\end{equation}
where $\|\chi_ju\|_s$ is defined by
\begin{equation}
\|\chi_ju\|_s = \left(\int_0^{\infty}(1 + |\log y|)^{2s}\|\chi_ju(y)\|^2_{L^2(M_j)}
\frac{dy}{y^n}\right)^{1/2},
\nonumber
\end{equation}
and $\|\chi_ju\|_{\mathcal B}, \ \|\chi_ju\|_{{\mathcal B}^{\ast}}$ are defined similarly.

Let us  note that many a-priori estimates and preliminary results which are proven in Chapter 2 for ${\bf H}^n$ may be straightforwardly generalized for
 $\mathcal M$. For example, 
Theorem 2.1.3 remains valid if we use $L^{2,s}$ instead of $\mathcal X^s$. Similarly, Theorem 2.2.10 can be extended to the case in which $(H - \lambda)u= 0$ in one of the regular ends $M_p\times(0,y_0)$ $(0<y_0<1)$.
Analogous extensions are true for Lemmas 2.2.4 $\sim$ 2.2.8 and so on.


\begin{lemma}
Suppose all ${\mathcal M}_j \ (1 \leq j \leq N)$ have a cusp. If $u \in {\mathcal B}^{\ast}$ satisfies $(H - \lambda)u = 0$ for some $\lambda > 0$ and, on each $\mathcal M_j$,
\begin{equation}
\lim_{R\to\infty}\frac{1}{\log R}\int_2^R\|u(y)\|_{L^2(M_j)}^2\frac{dy}{y^n} = 0, 
\nonumber
\end{equation}
then $u \in L^{2,s}, \ \forall s > 0$. Moreover, for any $s > 0$ and any compact interval $I \subset (0,\infty)$, there exists a constant $C_s > 0$ such that
\begin{equation}
\|u\|_s \leq C_s\|u\|_{{\mathcal B}^{\ast}}, \quad \forall \lambda \in I.
\label{eq:Chap3SEct2udecaysfast}
\end{equation}
\end{lemma}

Proof. 
For simplicity's sake, we assume that $N = 1$. Letting $U = \chi_1u$ and $M = M_1$, we have for $\epsilon > 0$ given in the assumption (A-4)
\begin{equation}
\left\{
\begin{split}
& \Big(- y^2(\partial_y^2 + \Delta_M) + (n - 2)y\partial_y- \frac{(n-1)^2}{4} - \lambda\Big)U = F, \\
&U \in {\mathcal B}^{\ast}, \quad F \in L^{2,(1 + \epsilon)/2}. 
\end{split}
\right.
\label{eq:Chap3Sect2equationofU}
\end{equation}
In fact, $F$ consists of $U$ and its 1st and 2nd order derivatives,
which, by Theorem 2.1.3, are in $L^{2,-(1+\epsilon)/2}$,
 multiplied by functions decaying like $(1 + |\log y|)^{-1-\epsilon}, \epsilon > 0$. Therefore, $F$ is in $L^{2,(1+\epsilon)/2}$.

We apply the boot-strap arguments.
In view of Lemma 2.2.6, letting ${\bf h} = L^2(M)$ and $\Delta_M$ the Laplace-Beltrami operator on $M$, we have
\begin{equation}
\int_0^{\infty}y^2\|\sqrt{-\Delta_M}U\|_{\bf h}^2\frac{dy}{y^n} \leq 
C\left(\|U\|_{\mathcal B^{\ast}}^2 + \|F\|_{\mathcal B}^2\right).
\label{C3S3Uestimate}
\end{equation}
Let $P_0$ be the projection associated with the 0 eigenvalue of $\Delta_M$, and put
$$
U_0 = P_0U, \quad U' = U - P_0U.
$$
Then we have by (\ref{C3S3Uestimate})
\begin{equation}
\|U'\|_s \leq C_s(\|U\|_{\mathcal B^{\ast}} + \|F\|_{\mathcal B}), \quad 
\forall s > 0.
\nonumber
\end{equation}
Since $U'$ satisfies the equation
$$
(H_0 - \lambda)U' = F' \in L^{2,(1+\epsilon)/2},
$$
 we have, by Theorem 2.1.3 (6), that
\begin{equation}
U', \ D_iU',\  D_iD_jU' \in L^{2,(1+\epsilon)/2}.
\label{C3S3U'DiU'estimate}
\end{equation}

Letting
$$
t = \log y, \quad 
 u_{0}(t) = e^{-(n-1)t/2}U_{0}(e^t), \quad 
  f_{0}(t) = e^{-(n-1)t/2}F_{0}(e^t), 
$$
we see that $u_{0}(t)$ satisfies
\begin{equation}
\left\{
\begin{split}
& (- \partial_t^2  - \lambda)u_{0} = 
f_{0}, \\
& \lim_{R\to\infty}\frac{1}{R}\int_2^R|u_{0}(t)|^2dt = 0,\\
& (1 + t)^{(1 + \epsilon)/2}f_{0}(t) \in L^2((2,\infty);dt).
\end{split}
\right.
\label{C3S3u0equation}
\end{equation}
Recall that the Green function of the 1-dimensional Helmholtz equation
$$
\Big(- \frac{d^2}{dt^2} - z\Big)u = f, \quad 
{\rm Im}\,z \geq 0
$$
is given by
$\displaystyle{
\frac{i}{2\sqrt{z}}e^{i\sqrt{z}|t - s|}}$.
Hence $u_0$ is represented as
\begin{eqnarray*}
u_0(t) &=& \frac{i}{2{\sqrt{\lambda}}}\int_0^{\infty}
e^{i\sqrt{\lambda}|t-s|}f_0(s)ds +C_+e^{i\sqrt{\lambda} t} +C_-e^{-i\sqrt{\lambda} t}
\\
&=& \frac{i}{2{\sqrt{\lambda}}}\int_0^{t}
e^{i\sqrt{\lambda}(t-s)}f_0(s)ds +
 \frac{i}{2{\sqrt{\lambda}}}\int_t^{\infty}
e^{i\sqrt{\lambda}(s-t)}f_0(s)ds \\
& & + \ C_+e^{i\sqrt{\lambda} t} \ + \ C_-e^{-i\sqrt{\lambda} t}.
\end{eqnarray*}
Since $f_0(t) \in L^1((0,\infty));dt)$, we have
$$
u_0(t) \sim \Big(C_+ + \frac{i}{2\sqrt{\lambda}}\int_0^{\infty}e^{-i\sqrt{\lambda}s}f_0(s)ds\Big)e^{i\sqrt{\lambda}t} + C_-e^{-i\sqrt{\lambda}t}, 
\quad t \to \infty,
$$
$$
u_0(t) \sim C_+e^{i\sqrt{\lambda}t} + \Big(C_- + \frac{i}{2\sqrt{\lambda}}\int_0^{\infty}e^{i\sqrt{\lambda}s}f_0(s)ds\Big)e^{-i\sqrt{\lambda}t}, 
\quad t \to -\infty.
$$
They imply, by (\ref{C3S3u0equation}),
$$
C_+ = 0 = - \frac{i}{2\sqrt{\lambda}}\int_0^{\infty}e^{-i\sqrt{\lambda}s}f_0(s)ds,
$$
$$
C_- = 0 = - \frac{i}{2\sqrt{\lambda}}\int_0^{\infty}e^{i\sqrt{\lambda}s}f_0(s)ds.
$$
We then have
$$
u_0(t) = \frac{i}{2\sqrt{\lambda}}\left(
e^{-i\sqrt{\lambda}t}\int_t^{\infty}e^{i\sqrt{\lambda}s}f_0(s)ds - 
e^{i\sqrt{\lambda}t}\int_t^{\infty}e^{-i\sqrt{\lambda}s}f_0(s)ds 
\right).
$$
Using Lemma 3.3, we then have
\begin{equation}
(1 + t)^{(-1 +\epsilon)/2}u_0, \ (1 + t)^{(-1+\epsilon)/2}\frac{d}{dt}u_0 \in L^2((0,\infty);dt).
\label{C3S3u0estimate}
\end{equation}
Then by (\ref{C3S3u0equation}), we also have
\begin{equation}
 (1 + t)^{(-1+\epsilon)/2}\frac{d^2}{dt^2}u_0 \in L^2((0,\infty);dt).
\label{C3S3u0''estimate}
\end{equation}
By (\ref{C3S3U'DiU'estimate}), (\ref{C3S3u0estimate}) and (\ref{C3S3u0''estimate}), we  have
 $U, D_iU, D_iD_jU \in L^{2,(-1+\epsilon)/2}$. Hence we have $F \in L^{2,(1+2\epsilon)/2}$.
 
We return to the equation (\ref{eq:Chap3Sect2equationofU}), and apply the same arguments as above. Then we have  $U, D_iU, D_iD_jU \in L^{2,(-1+2\epsilon)/2}$, hence $F \in L^{2,(1 + 3\epsilon)/2}$. 
We repeat these procedures to obtain $U \in L^{2,(-1 +N\epsilon)/2}, \ \forall N > 0$ and the inequality (\ref{eq:Chap3SEct2udecaysfast}). \qed


\begin{theorem}
(1) If one of  ${\mathcal M}_j$ has a regular infinity, $\sigma_p(H)\cap(0,\infty) = \emptyset$. \\
\noindent
(2) If all of ${\mathcal M}_j$ have a cusp, then $\sigma_p(H)\cap(0,\infty)$ is discrete with finite multiplicities, whose possible accumulation points are $0$ and $\infty$.
\end{theorem}
Proof. We shall prove (1). Let $u$ be the eigenvector of $H$ with eigenvalue $\lambda \in (0,\infty)$. Applying Theorem 2.2.10 on ${\mathcal M}_j$ having a regular infinity, we see that $u$ vanishes in a neighborhood of infinity of ${\mathcal M}_j$. By the unique continuation theorem, $u$ vanishes identically on $\mathcal M$.

To prove (2) assume that there exist an infinite number of eigenvlaues (counting multiplicities) in a compact interval $I \subset (0,\infty)$. Let $u_n, 
n = 1, 2, \cdots$, be the associated orthonormal system of eigenvectors. Choose $x_0 \in {\mathcal K}$ arbitrarily, and let $\chi_R$ be such that $\chi(x) = 1$ for ${\rm dist}\,(x,x_0) < R$, $\chi(x) = 0$ for ${\rm dist}\,(x,x_0) > R$. By (\ref{eq:Chap3SEct2udecaysfast}), for any  $\epsilon > 0$, there exists $R > 0$ independent of $n$ such that $\|(1 - \chi_R)u_n\|_{L^2} < \epsilon$ and 
$\|\chi_R u_n  \|_{L^2} \geq 1- 2 {\sqrt \epsilon}$. 
Using Rellich's theorem, one can choose a subsequence of $\{\chi_Ru_n\}_{n\geq1}$ which converges in $L^2$,
$$
\chi_R u_n \to u, \quad \| u  \|_{L^2} \geq 1- 2 {\sqrt \epsilon}.
$$
Thus, for sufficiently large $n, m$,
\begin{eqnarray*}
& &\left(u_n, u_m \right) =\left( \chi_R u_n, \, \chi_R u_m   \right)+
\left((1 - \chi_R)u_n, \, \chi_R u_m\right) 
\\ \nonumber
& &+\left(\chi_R u_n,\,  (1 - \chi_R)u_m \right)+ 
\left( (1 - \chi_R)u_n,\,(1 - \chi_R)u_m  \right)
\\
\nonumber
& &\geq (1- 2 {\sqrt \epsilon})^2-3\epsilon >0, \quad \hbox{if}\,\, \epsilon < \frac{1}{16}.
\end{eqnarray*}
 This is a contradiction to $(u_n, u_m)=0$.
  \qed


\begin{theorem}
Suppose $\lambda > 0$, and $u \in {\mathcal B}^{\ast}$ satisfies 
$(H - \lambda)u = 0$. 
Furthermore, assume that, when ${\mathcal M}_j$ has a regular infinity,
\begin{equation}
\lim_{R\to\infty}\frac{1}{\log R}\int_{1/R}^{1/2}\|u(\cdot,y)\|^2_{L^2(M_j)}\frac{dy}{y^n} = 0,
\nonumber
\end{equation}
and when ${\mathcal M}_j$ has a cusp,
\begin{equation}
\lim_{R\to\infty}\frac{1}{\log R}\int_{2}^{R}\|u(\cdot,y)\|^2_{L^2(M_j)}\frac{dy}{y^n} = 0.
\nonumber
\end{equation}
Then:

\noindent
(1) If one of ${\mathcal M}_j$ has a regular infinity, then $u = 0$. \\
\noindent
(2) If all of ${\mathcal M}_j$ have a cusp, then $u \in L^{2,s}, \ \forall s > 0$.
\end{theorem}
Proof. Applying Theorem 2.2.10 to ${\mathcal M}_j$ with regular infinity, we see that
$u$ vanishes on an open set of ${\mathcal M}_j$, hence $u = 0$ by the unique continuation theorem. The assertion (2) follows from Lemma 3.4. \qed

\bigskip
As in Chap. 2, \S 2, we put
\begin{equation}
\sigma_{\pm}(\lambda) = \frac{n-1}{2} \mp i\sqrt{\lambda}.
\nonumber
\end{equation}
We say that a solution $u \in {\mathcal B}^{\ast}$ of the equation
\begin{equation}
(H - \lambda)u = f \in {\mathcal B}
\nonumber
\end{equation}
satisfies the outgoing radiation condition, when ${\mathcal M}_j$ has a regular infinity , if
\begin{equation}
\lim_{R\to\infty}\frac{1}{\log R}\int_{1/R}^{1/2}\|(D_y - \sigma_+(\lambda))u(\cdot,y)\|^2_{L^2(M_j)}\frac{dy}{y^n} = 0,
\end{equation}
and when ${\mathcal M}_j$ has a cusp
\begin{equation}
\lim_{R\to\infty}\frac{1}{\log R}\int_{2}^{R}\|(D_y - \sigma_-(\lambda))u(\cdot,y)\|^2_{L^2(M_j)}\frac{dy}{y^n} = 0.
\label{C3S3IncRadCond}
\end{equation}
The incoming radiation condition is defined similarly by exchanging $\sigma_+(\lambda)$ and $\sigma_-(\lambda)$.

Let us remark that, compared to the case of ${\bf H}^n$ (see Chap. 2, (\ref{eq:Chap2Sect2RadCondplusinfty})), the condition (\ref{C3S3IncRadCond}) seems to be confusing. Due to the presence of 0-eigenvalue of $\Delta_M$, there exist generalized eigenfunctions for $H_{free}$ which behave like $y^{(n-1)/2 \pm i\sqrt{\lambda}}$ as $y \to \infty$. To distinguish these two functions, we need  (\ref{C3S3IncRadCond}).


\begin{theorem}
Let $\lambda > 0$ and suppose $u \in {\mathcal B}^{\ast}$ satisfies $(H - \lambda)u = 0$ and the outgoing or incoming radiation condition. Then:

\noindent
(1) If one of ${\mathcal M}_j$ has a regular infinity, then $u = 0$. \\
\noindent
(2) If all ${\mathcal M}_j$ have a cusp, then $u \in L^{2,s}, \ \forall s > 0$.
\end{theorem}
Proof. We assume that the ends $\mathcal M_1$, $\cdots$, $\mathcal M_{\mu}$ have regular infinities, and $\mathcal M_{\mu +1}$, $\cdots$, $\mathcal M_N$ have cusps. Recall that for $1 \leq j \leq {\mu}$, $\mathcal M_j$ is diffeomorphic to $M_j\times(0,1)$, and for ${\mu} +1 \leq j \leq N$, $\mathcal M_j$ is diffeomorphic to $M_j\times(1,\infty)$. Let $\{\chi_j\}_{j=0}^N$ be a smooth partition of unity such that $\sum_{j=0}^N\chi_j = 1$ on $\mathcal M$, and ${\rm supp}\,\chi_j \subset \mathcal M_j$ for $1 \leq j \leq N$.
We shall assume that for $1 \leq j \leq \mu$, 
\begin{equation}
\chi_j(y) = \left\{
\begin{split}
1, & \quad (y < 1/2), \\
0, & \quad (y > 3/4),
\end{split}
\right.
\nonumber
\end{equation}
 and for $\mu +1 \leq j \leq N$, 
 \begin{equation}
\chi_j(y) = \left\{
\begin{split}
0, & \quad (y < 3/2), \\
1, & \quad (y > 2).
\end{split}
\right.
\nonumber
\end{equation}
We take $\rho(t) \in C_0^{\infty}({\bf R})$ such that $\rho(t) = \rho(-t)$ and
\begin{equation}
\rho(t) = \left\{
\begin{split}
c, & \quad |t| < 1, \\
0, & \quad |t| > 2,
\end{split}
\right.
\nonumber
\end{equation}
where $c$ is a positive constant such that
$$
\int_{-\infty}^0\rho(t)dt = \int_0^{\infty}\rho(t)dt = 1.
$$
We put
$$
\varphi(t) = \int_{-\infty}^{t}\rho(s)ds, \quad \psi(t) = \int_t^{\infty}\rho(s)ds,
$$
and
$$
\varphi_R(y) = \varphi\left(\frac{\log y}{\log R}\right), \quad
\psi_R(y) = \psi\left(\frac{\log y}{\log R}\right).
$$
Then we have
\begin{equation}
\begin{split}
& \chi_j(y)\varphi_R(y) \in C_0^{\infty}(\mathcal M_j) \quad {\rm for} \quad 
1 \leq j \leq \mu, \\
&\chi_j(y)\psi_R(y) \in C_0^{\infty}(\mathcal M_j) \quad {\rm for} \quad 
\mu +1 \leq j \leq N.
\end{split}
\nonumber
\end{equation}
Moreover,
\begin{equation}
\lim_{R\to\infty}\varphi_R(y) = \varphi(0) = 1, \quad 
\lim_{R\to\infty}\psi_R(y) = \psi(0) = 1.
\label{C3S3psi(0)=0}
\end{equation}
Since $(H - \lambda)u = 0$, we have
\begin{equation}
0 = ((H-\lambda)u,\chi_j\varphi_Ru) = 
(u,[H,\chi_j\varphi_R]u).
\nonumber
\end{equation}
Therefore, we have
\begin{equation}
(u,[H,\chi_j]\varphi_Ru) + (u,\chi_j[H,\varphi_R]u) = 0,
\nonumber
\end{equation}
\begin{equation}
(u,[H,\chi_j]\psi_Ru) + (u,\chi_j[H,\psi_R]u) = 0,
\nonumber
\end{equation}
\begin{equation}
(u,[H,\chi_0]u) = 0.
\nonumber
\end{equation}
We add them, and let $R \to \infty$. Then by (\ref{C3S3psi(0)=0})
$$
\sum_{j=1}^{\mu}(u,[H,\chi_j]\varphi_Ru) + 
\sum_{j=\mu +1}^N(u,[H,\chi_j]\psi_Ru) + 
(u,[H,\chi_0]u) \to \sum_{j=0}^N(u,[H,\chi_j]u) = 0.
$$
Therefore, as $R \to \infty$,
\begin{equation}
\sum_{j=1}^{\mu}(u,\chi_j[H,\varphi_R]u) + \sum_{j=\mu +1}^N(u,\chi_j[H,\psi_R]u) 
\to 0.
\label{C3S3(u,[H,psiR]u)to0}
\end{equation}
We put
$$
V_j = H - \Big(- D_y^2 + (n-1)D_y - y^2\Delta_{M_j} - \frac{(n-1)^2}{4} \Big).
\nonumber
$$
Then we have, for $1 \leq j \leq \mu$,
\begin{equation}
\begin{split}
[H,\varphi_R] &= [-D_y^2 + (n-1)D_y,\varphi_R] + [V_j,\varphi_R]\\
&= -\frac{2}{\log R}\rho\left(\frac{\log y}{\log R}\right)\big(D_y - \frac{n-1}{2}\big) + \frac{1}{\log R}L_{j,R}.
\end{split}
\label{C3S3[H,varphiR]}
\end{equation}
Here $L_{j,R}$ is a 1st order differential operator 
\begin{equation}
L_{j,R} = a_{j,R}(x,y)D_y + b_{j,R}(x,y)D_x + c_{j,R},
\label{C3S3LjR}
\end{equation}
whose coefficients satisfy, due to (\ref{eq:Chap3Sect2DecayMetric1}),
\begin{equation}
|a_{j,R}(x,y)| + |b_{j,R}(x,y)| + |c_{j,R}(x,y)| \leq C(1 + |\log y|)^{-1-\epsilon},
\label{C3S3ajRdecay}
\end{equation}
where the constant $C$ is independent of $R > 1$. Similarly, we have, for $\mu +1 \leq j \leq N$,
\begin{equation}
\begin{split}
[H,\psi_R] &= [-D_y^2 + (n-1)D_y,\psi_R] + [V_j,\varphi_R]\\
&= \frac{2}{\log R}\rho\left(\frac{\log y}{\log R}\right)\big(D_y - \frac{n-1}{2}\big) + \frac{1}{\log R}L_{j,R},
\end{split}
\label{C3S3[H,varphiR]}
\end{equation}
where $L_{j,R}$ is a 1st order differential operator having the same property as above.
In view of (\ref{C3S3(u,[H,psiR]u)to0}), we then have
\begin{equation}
\begin{split}
&- \sum_{j=1}^\mu\frac{2}{\log R}(\chi_j\rho\left(\frac{\log y}{\log R}\right)
\big(D_y - \frac{n-1}{2}\big)u,u) \\
& + \sum_{j=\mu+1}^N\frac{2}{\log R}(\chi_j\rho\left(\frac{\log y}{\log R}\right)
 \big(D_y - \frac{n-1}{2}\big)u,u) \\
& + \sum_{j=1}^N\frac{1}{\log R}(\chi_jL_{j,R}u,u) \to 0.
\end{split}
\label{C3S3Dy-(n-1)/2uto0}
\end{equation}
We consider the case when $u$ satisfies the outgoing radiation condition.
Then we have, by (\ref{C3S3Dy-(n-1)/2uto0}),
\begin{equation}
\sum_{j=1}^N\frac{2i\sqrt{\lambda}}{\log R}\, (\chi_j\rho\big(\frac{\log y}{\log R}\big)u,u) \to 0,
\label{C3S3ApplyRadCond}
\end{equation}
since  one can replace $(D_y - (n-1)/2)$ by $-i\sqrt{\lambda}$ for $1 \leq j \leq \mu$,  by $i\sqrt{\lambda}$ for $\mu +1 \leq j \leq N$, and $(\chi_jL_{j,R}u,u)/\log R \to 0$.
This shows that, for $1 \leq j \leq N$,
$$
\frac{1}{\log R}\int_0^{\infty}\chi_j(y)\rho\big(\frac{\log y}{\log R}\big)
\|u(\cdot,y)\|_{L^2(M_j)}^2\frac{dy}{y^n} \to 0.
$$
Thus, $u$ satisfies conditions of Theorem 3.6, providing the desired result.

The case in which $u$ satisfies the incoming radiation condition is proved similarly.
 \qed

\bigskip
These preparations are sufficient to prove the limiting absortion principle for $H$ as in Chap. 2, \S 2. 

 
\begin{theorem}
For $\lambda \in \sigma_e(H)\setminus\sigma_p(H)$, there exists a limit
$$
\lim_{\epsilon \to 0}R(\lambda \pm i\epsilon) \equiv 
R(\lambda \pm i0) \in {\bf B}({\mathcal B};{\mathcal B}^{\ast})
$$
in the weak $\ast$ sense. Moreover, for any compact interval $I \subset \sigma_e(H)\setminus\sigma_p(H)$, there exists a constant $C > 0$ such that
\begin{equation}
\|R(\lambda \pm i0)f\|_{{\mathcal B}^{\ast}} \leq C\|f\|_{\mathcal B},
\quad \lambda \in I.
\nonumber
\end{equation} 
For 
$f \in {\mathcal B}$, we put $u = R(\lambda \pm i0)f$. Then $u$ is a unique solution to the equation $(H - \lambda)u = f$ satisfying the outgoing (for the case $+$) or incoming (for the case $-$) radiation condition. For $f, g \in {\mathcal B}$,  $(R(\lambda \pm i0)f,g)$ is continuous with respect to $\lambda \in \sigma_e(H)\setminus\sigma_p(H)$.
\end{theorem}

In order to prove Theorem 3.8, recall that Lemmas 2.2.4 $\sim$ 2.2.9  also hold for $M_j \times (0,\infty)$ with $\bf h$ replaced by $L^2(M_j)$.
Let  $\chi_j$ be the partition of unity (\ref{eq:Chap3Sec2Partitionofunity}), and put $u = R(\lambda + i\epsilon)f, \ 
u_j = \chi_ju$. 
Then, with $\epsilon$ defined by (\ref{eq:Chap3Sect2DecayMetric1}),
\begin{equation} \label{E2.10}
\| u_j \|_{{\mathcal B}^{\ast}}\leq C_{s} \left(\| f \|_{{\mathcal B}}+\|u_j\|_{-s} \right), \quad
1/2 < s< (1+\epsilon)/2,
\end{equation}
where $C_s$ is independent of $\lambda \in I$. Indeed, we first observe that
$$
(H-\lambda-i\epsilon) u_j =\chi_jf+[H, \chi_j] u.
$$
By Theorem 2.1.3 (6),
$$
\|D_i u_j\|_{-s}, \|D_i D_l u_j\|_{-s}  \leq C_{s} \left(\| f \|_{{\mathcal B}}+\|u_j\|_{-s}\right),
$$
and as
$\left[ H, \chi_j \right],\,\left[ H_{free(j)},\, \chi_j \right]$ are compactly supported, we also have
$$
 \|\left[ H,\, \chi_j \right]u\|_{\mathcal B},\, \|\left[ H_{free(j)},\, \chi_j \right]u\|_{\mathcal B} \leq C_{s} \left(\| f \|_{{\mathcal B}}+\|u_j\|_{-s} \right).
$$
At last, rewriting the equation for $u_j$ as
$$
(H_{free(j)} -\lambda-i\epsilon) u_j= \chi_jf+[H_{free(j)}, \chi_j] u+\chi_j V u,
$$
and using (\ref{eq:Chap3Sect2DecayMetric1}), we obtain (\ref{E2.10}) by Lemma 2.2.9.
Once we have derived  estimate (\ref{E2.10}), the remaining arguments are essentially 
the same as those in Chap. 2. Namely,
arguing in the same way as in Lemma 2.2.13, we can prove the following lemma. 


\begin{lemma}
Take $s > 1/2$ sufficiently close to $1/2$. Let $I$ be any compact interval in $(0,\infty)\setminus\sigma_p(H)$, and put $J = \{\lambda \pm i\epsilon ; \lambda \in I, 0 < \epsilon < 1\}$. \\
\noindent
(1) There exists a constant $C_s > 0$ such that
\begin{equation}
\sup_{z \in J}\|R(z)f\|_{-s} \leq C_s\|f\|_{\mathcal B}.
\nonumber
\end{equation}
(2) For any $f \in {\mathcal B}$ and $\lambda \in (0,\infty)\setminus\sigma_p(H)$, the strong limit $\lim_{\epsilon\to 0}R(\lambda \pm i0)f$ exists in $L^{2,-s}$. \\
\noindent
(3) $R(\lambda \pm i0)f$ is an $L^{2,-s}$-valued continuous function of $\lambda \in (0,\infty)\setminus\sigma_p(H)$.
\end{lemma}

Since $L^{2,s}$ $(s > 1/2)$ is dense in $\mathcal B$, Theorem 3.8 follows from Lemma 3.9
and (\ref{E2.10}). \qed


\subsection{Fourier transform associated with $H$} 
One can apply the abstract theory in Chap. 2, \S 4 to $H$ after suitable modifications. However, we shall give here a direct approach to the spectral representation for $H$.

Let $H_{free(j)}$ be as above and
$\chi_j$ as in (\ref{eq:Chap3Sec2Partitionofunity}). We put
\begin{equation}
\widetilde V_j = H - H_{free(j)} \quad {\rm on} \quad \mathcal M_j.
\label{C3S3widetildeVj}
\end{equation}
This is symmetric, since so are $H$ and $H_{free(j)}$ on $C_0^{\infty}(\mathcal M_j)$.
Using
\begin{equation}
(H_{free(j)} - \lambda)\chi_jR(\lambda \pm i0) = 
\chi_j + \left([H_{free(j)},\chi_j] - \chi_j\widetilde V_j\right)
R(\lambda \pm i0),
\label{C3S3HfreechiR=}
\end{equation}
we have
\begin{equation}
\begin{split}
 \chi_jR(\lambda \pm i0) 
& =  
R_{free(j)}(\lambda \pm i0)\chi_j \\
& \ \ \ \   + 
R_{free(j)}(\lambda \pm i0)\left([H_{free(j)},\chi_j] - \chi_j\widetilde V_j \right)
R(\lambda \pm i0).
\end{split}
\label{C3S3Rfree(j)chij=}
\end{equation}
This formula suggests how the generalized Fourier transform is constructed by the perturbation method.

\subsubsection{Definition of ${\mathcal F}_{free(j)}^{(\pm)}(k)$} 
 Let $0 = \lambda_{j,0} < \lambda_{j,1} \leq \lambda_{j,2} \leq \cdots$ be the eigenvalues of the Laplace-Beltrami operator on $M_j$ and $|M_j|^{-1/2}= \varphi_{j,0}, \varphi_{j,1}, \varphi_{j,2}, \cdots$ the associated orthonormal eigenvectors, where $|M_j|$ is the volume of $M_j$.  We define, for $\phi \in L^2(M_j)$,
\begin{equation}
P_{j,m}\phi = (\phi,\varphi_{j,m})_{L^2(M_j)}\varphi_{j,m},
\label{C3S3Projjm}
\end{equation}
\begin{equation}
\Pi_{j,m}\phi = \left(\phi,\varphi_{j,m}\right)_{L^2(M_j)}.
\label{C3S3Paijm}
\end{equation}

Assume that for $1 \leq j \leq \mu$, $\mathcal M_j$ has a regular infinity, and for $\mu + 1 \leq j \leq N$, $\mathcal M_j$ has a cusp.

\medskip
\noindent 
(i) For $1 \leq j \leq \mu$ (the case of regular infinity), we define
\begin{equation}
{\mathcal F}_{free(j)}^{(\pm)}(k) = 
\sum_{m=0}^{\infty}C_{j,m}^{(\pm)}(k)\,
P_{j,m}\otimes
F_{free(j),m}^{(\pm)}(k),
\label{C3S3Ffreejpmk}
\end{equation}
where $F_{free(j),m}^{(\pm)}$ is defined by (\ref{C3S2Ffreem}), (\ref{C3S2Ffreeopm}), (\ref{eq:Chap3Sect2F0gammaastpm}) with $M$ replaced by $M_j$, and $C_{j,m}^{(\pm)}(k)$ is the constant in (\ref{eq:Chap3Sect2constCgammaastk}) with $\lambda_m$ replaced by $\lambda_{j,m}$, i.e.
\begin{equation}
C_{j,m}^{(\pm)}(k) = 
\left\{
\begin{split}
& \left(\frac{\sqrt{\lambda_{j,m}}}{2}\right)^{\mp ik}, \quad
(\lambda_{j,m} \neq 0), \\
& \frac{\pm i}{k\omega_{\pm}(k)}\sqrt{\frac{\pi}{2}}, \quad (\lambda_{j,m} = 0).
\end{split}
\right.
\label{C3S3Cjmk}
\end{equation}
Thus, in this case,  $F_{free(j)}^{(\pm)}(k)=F_{free(j), reg}^{(\pm)}(k)$, see (\ref{eq:Chap3Sect2F0regpm}).

\medskip
\noindent
(ii) For $\mu + 1 \leq j \leq N$ (the case of cusp), we define
\begin{equation}
{\mathcal F}_{free(j)}^{(\pm)}(k) =
P_{j,0}\otimes F_{free(j),0}^{(\mp)}(k).
\label{C3S3Ffreejpmkcusp} 
\end{equation}
Thus, in this case,  $F_{free(j)}^{(\pm)}(k)=F_{free(j), c}^{(\pm)}(k)$, see (\ref{eq:Chap3Sect2Focplusminus}).

\subsubsection{Definition of ${\mathcal F}^{(\pm)}(k)$}
 For $1 \leq j \leq N$, we define
\begin{equation}
\mathcal F^{(\pm)}_j(k) = \mathcal F^{(\pm)}_{free(j)}(k)Q_j(k^2 \pm i0),
\label{C3S3Fpmjkdefine}
\end{equation}
\begin{equation}
Q_j(z) = \chi_j + \left([H_{free(j)},\chi_j] - \chi_j\widetilde V_j\right)R(z) 
= (H_{free(j)}-z)\chi_jR(z).
\label{C3S3Qjz}
\end{equation}
Finally, we define the Fourier transform associated with $H$ by
\begin{equation}
{\mathcal F}^{(\pm)}(k) = \big({\mathcal F}_{1}^{(\pm)}(k), \cdots, {\mathcal F}_{N}^{(\pm)}(k)\big).
\label{C3S3Fpmkdefine}
\end{equation}

\subsubsection{Asymptotic expansion of the resolvent}
For $f, g \in {\mathcal B}^{\ast}$ on ${\mathcal M}$, by
$f \simeq g$
we mean that on each end the following expansion
\begin{equation}
\lim_{R\to\infty}\frac{1}{\log R}\int_{1/R}^{R}\rho_j(y)\|f(y) - g(y)\|^2_{L^2(M_j)}
\frac{dy}{y^n} = 0
\nonumber
\end{equation}
holds, where $\rho_j(y) = 1$ $(y<1)$, $\rho_j(y) = 0$ $(y > 1)$ when $\mathcal M_j$ has a regular infinity, and $\rho_j(y) = 0$ $(y < 1)$, $\rho_j(y) = 1$ $(y > 1)$ when $\mathcal M_j$ has a cusp. Applying Theorem 2.6 on each end, we get the following theorem. 


\begin{theorem} Let $f \in {\mathcal B}$, 
$k^2 \in \sigma_e(H)\setminus\sigma_p(H)$, and $\chi_j$ the partition of unity from (\ref{eq:Chap3Sec2Partitionofunity}). Then we have
\begin{equation}
\begin{split}
R(k^2 \pm i0) f &\simeq \omega_{\pm}(k)\sum_{j=1}^\mu
\chi_jy^{(n-1)/2\mp ik}
{\mathcal F}^{(\pm)}_j(k)f \\
& + \omega_{\pm}^{(c)}(k)\sum_{j=\mu+1}^N
\chi_jy^{(n-1)/2\pm ik}
{\mathcal F}^{(\pm)}_j(k)f.
\end{split}
\nonumber
\end{equation}
\end{theorem}

We put
\begin{equation}
{\bf h}_{\infty} = \left({\mathop\oplus_{j=1}^\mu} L^2(M_j)\right)\oplus
\left({\mathop\oplus_{j= \mu+1}^N}P_{j,0}L^2(M_j)\right), 
\label{eq:Chap3Sect2hinfty}
\end{equation}
As a matter of fact, 
$$
P_{j,0}L^2(M_j) = {\bf C}\,\varphi_{j,0} = \left\{c\,\varphi_{j,0}\,;\, c \in {\bf C}\right\}, \quad \varphi_{j,0} = |M_j|^{-1.2},
$$ 
equipped with the inner product
\begin{equation}
(c_1\varphi_{j,0},c_2\varphi_{j,0})_{{\bf C}_j} = c_1\overline{c_2}.
\label{C3S3Cvarphij0innerprod}
\end{equation}
For $\phi, \psi \in {\bf h}_{\infty}$ we define the inner product by
\begin{equation}
(\phi,\psi)_{{\bf h}_{\infty}} = \sum_{j=1}^\mu(\phi_j,\psi_j)_{L^2(M_j)} + \sum_{j=\mu+1}^N(\phi_j,\psi_j)_{{\bf C}_j}.
\label{C3S3hinftyinnerproduct}
\end{equation}
We then have the following lemma.


\begin{lemma}
For $f, g \in {\mathcal B}$ and $k^2 \in \sigma_e(H)\setminus\sigma_p(H)$,
\begin{equation}
\frac{k}{\pi i}\left(\left[R(k^2 + i0) - R(k^2 - i0)\right]f,g\right) =  \left({\mathcal F}^{(\pm)}(k)f,{\mathcal F}^{(\pm)}(k)g\right)_{{\bf h}_{\infty}}.
\nonumber
\end{equation}
\end{lemma}
Proof. Take $\chi \in C_0^{\infty}({\bf R})$ such that $\chi(t) = 1 \ (|t| < 1), \ \chi(t) = 0 \ (|t| > 2)$. Let $\chi_R \in C_0^{\infty}({\mathcal M})$ be such that $\chi_R = 1$ on a neighborhood of ${\mathcal K}$, $\chi_R = \chi(\log y/\log R)$ on each  ${\mathcal M}_j$, where $R > 0$ is a large parameter. Let $\chi_j$ be the partition of unity from (\ref{eq:Chap3Sec2Partitionofunity}).
Putting
$u = R(k^2 + i0)f, \ v = R(k^2 + i0)g$, we have 
\begin{equation}
(\chi_Ru,Hv) - (Hu,\chi_Rv) 
= ([H,\chi_R]u,v) \\ 
 = \sum_{j=1}^N(\chi_j[H,\chi_R]u,v),
\nonumber
\end{equation}
since $\chi_R = 1$ on a neighborhood of ${\rm supp}\,\chi_0$.
Next we take $\widetilde \chi_j \in C^{\infty}(\mathcal M_j)$ such that ${\rm supp}\,\widetilde \chi_j \subset \mathcal M_j$ and $\widetilde \chi_j = 1$ on ${\rm supp}\,\chi_j$. Then, by Theorems 3.8, 2.1.3 (5)
and (\ref{eq:Chap3Sect2DecayMetric1}), 
we have, as $R \to \infty$,
\begin{equation}
\begin{split}
(\chi_Ru,Hv) - (Hu,\chi_Rv) &= \sum_{j=1}^N(\chi_j[H,\chi_R]\widetilde \chi_ju,v)\\
&= \sum_{j=1}^N(\chi_j[H_{free(j)},\chi_R]\widetilde\chi_j u,v) + o(1).
\end{split}
\nonumber
\end{equation}
On each end, we have
\begin{equation}
\begin{split}
\big[- y^2\partial_y^2 + (n - 2)y\partial_y,\chi_R\big] &= 
- \frac{2}{\log R}\chi'\Big(\frac{\log y}{\log R}\Big)\Big(D_y - \frac{n-1}{2}\Big) \\
&  - \left(\frac{1}{\log R}\right)^2\chi''\Big(\frac{\log y}{\log R}\Big).
\end{split}
\nonumber
\end{equation}
Therefore,
$$
 (\chi_j[H_{free(j)},\chi_R]\widetilde\chi_j u,v) =  - \frac{2}{\log R}\left(\chi_j\chi'\Big(\frac{\log y}{\log R}\Big)\Big(D_y - \frac{n-1}{2}\Big)u,v\right) + o(1).
$$ 
Since,  by Theorem 3.8, $u$ satisfies the outgoing radiation condition, for $1 \leq j \leq \mu$, one can replace $(D_y - (n-1)/2)u$ by $- iku$. Hence,
\begin{eqnarray*}
(\chi_j[H_{free(j)},\chi_R]\widetilde\chi_ju,v) &=& \frac{2ik}{\log R}\left(\chi_j\chi'\Big(\frac{\log y}{\log R}\Big)u,v\right) + o(1) \\
 &=& \frac{2ik}{\log R}\cdot\frac{\pi}{2k^2}\left(\chi_j\chi'\Big(\frac{\log y}{\log R}\Big)y^{n-1}{\mathcal F}_j^{(+)}(k)f,
 {\mathcal F}_j^{(+)}(k)g\right) + o(1) \\
 &=& \frac{\pi i}{k}\Big({\mathcal F}_j^{(+)}(k)f,{\mathcal F}_j^{(+)}(k)g\Big)_{L^2(M_j)} + o(1),
\end{eqnarray*}
where we have used Theorem 3.10 in the 2nd line, and 
$$
\frac{1}{\log R}\int_{-\infty}^{0}\chi'\Big(\frac{\log y}{\log R}\big)\frac{dy}{y} = 1.
$$
For $\mu + 1 \leq j \leq N$, one replaces $(D_y - (n-1)/2)u$ by $iku$, and uses
$$
\frac{1}{\log R}\int_{0}^{\infty}\chi'\Big(\frac{\log y}{\log R}\big)\frac{dy}{y} = - 1
$$
to obtain
\begin{eqnarray*}
 (\chi_j[H_{free(j)},\chi_R]\widetilde\chi_j u,v) &=&- \frac{2ik}{\log R}\left(\chi_j\chi'\Big(\frac{\log y}{\log R}\Big)u,v\right) + o(1) \\
 &=& \frac{\pi i}{k} \big({\mathcal F}_j^{(+)}(k)f,
 {\mathcal F}_j^{(+)}(k)g\big)_{{\bf C}_j} + o(1).
\end{eqnarray*}
Using
\begin{eqnarray*}
(\chi_Ru,Hv) - (Hu,\chi_Rv) &\to&
(u,g) - (f,v) \\
&=& (R(k^2 + i0)f,g) - (f,R(k^2 + i0)g),
\end{eqnarray*}
we complete the proof of the lemma. \qed

\bigskip
We put
\begin{equation}
\widehat{\mathcal H} = L^2((0,\infty);{\bf h}_{\infty};dk).
\nonumber
\end{equation}


\begin{theorem}
We define $\big({\mathcal F}^{(\pm)}f\big)(k) = {\mathcal F}^{(\pm)}(k)f$ for $f \in {\mathcal B}$. Then  ${\mathcal F}^{(\pm)}$ is uniquely extended to a bounded operator from $L^2({\mathcal M})$ to $\widehat{\mathcal H}$ with the following properties.

\noindent
(1) $\ {\rm Ran}\,\,{\mathcal F}^{(\pm)} = \widehat{\mathcal H}$. \\
\noindent
(2) $\ \|f\| = \|{\mathcal F}^{(\pm)}f\|$ for $f \in {\mathcal H}_{ac}(H)$. \\
\noindent
(3) $\ {\mathcal F}^{(\pm)}f = 0$ for $f \in {\mathcal H}_p(H)$.  \\
\noindent
(4) $\ 
\left({\mathcal F}^{(\pm)}Hf\right)(k) = 
k^2\left({\mathcal F}^{(\pm)}f\right)(k)$ for 
$f \in D(H)$. \\
\noindent
(5)$\ {\mathcal F}^{(\pm)}(k)^{\ast} \in {\bf B}({\bf h}_{\infty};{\mathcal B}^{\ast})$ and
$(H - k^2){\mathcal F}^{(\pm)}(k)^{\ast} = 0$  for $k^2 \in (0,\infty)\setminus\sigma_p(H)$. \\
\noindent
(6) For $f \in {\mathcal H}_{ac}(H)$, the inversion formula holds:
\begin{eqnarray*}
f = \left({\mathcal F}^{(\pm)}\right)^{\ast}{\mathcal F}^{(\pm)}f 
= \sum_{j=1}^N\int_0^{\infty}{\mathcal F}^{(\pm)}_j(k)^{\ast}
\left({\mathcal F}^{(\pm)}_jf\right)(k)dk.
\end{eqnarray*}
\end{theorem}

\noindent
{\bf Remark} The meaning of the integral in (6) is as follows. Let $(0,\infty)\setminus\sigma_p(H) = \cup_{i=1}^{\infty}I_i$, $I_i = (a_i,b_i)$ being non-overlapping connected open interval. For $g(k) \in \widehat{\mathcal H}$, we have by (5)
$$
\int_{\sqrt{a_i}+\epsilon}^{\sqrt{b_i}-\epsilon}{\mathcal F}^{(\pm)}_j(k)^{\ast}
g(k)dk \in {\mathcal B}^{\ast}.
$$
As a matter of fact, it belongs to $L^2({\mathcal M})$, and 
$$
\lim_{\epsilon \to 0}\int_{\sqrt{a_i}+\epsilon}^{\sqrt{b_i}-\epsilon}{\mathcal F}^{(\pm)}_j(k)^{\ast}
g(k)dk \in L^2({\mathcal M})
$$
in the sense of strong convergence in $L^2({\mathcal M})$. Denoting this limit by
$$
\int_{\sqrt{I_i}}{\mathcal F}^{(\pm)}_j(k)^{\ast}
g(k)dk,
$$
we define
$$
\int_{0}^{\infty}{\mathcal F}^{(\pm)}_j(k)^{\ast}
g(k)dk = 
\sum_{i=1}^{\infty}\int_{\sqrt{I_i}}{\mathcal F}^{(\pm)}_j(k)^{\ast}
g(k)dk.
$$

\bigskip
Proof. Let $E(\lambda)$ be the spectral decomposition for $H$. Since the interval $(a_i,b_i)$ does not contain eigenvalues of $H$, we have by Lemma 3.11 and Stone's formula
$$
\frac{1}{2\pi i}\int_{a_i+\epsilon}^{b_i-\epsilon}
\left([R(\lambda + i0) - R(\lambda - i0)]f,f\right)d\lambda = 
\int_{\sqrt{a_i+\epsilon}}^{\sqrt{b_i-\epsilon}}
\|{\mathcal F}^{(\pm)}(k)f\|^2dk,
$$
for $f \in {\mathcal B}$. When $\epsilon \to 0$, the left-hand side converges to $(E((a_i,b_i))f,f)$. Therefore, so does the right-hand side and
$$
(E(I_i)f,f) = \int_{\sqrt{I_i}}\|{\mathcal F}^{(\pm)}(k)f\|^2dk.
$$
Since the end points of $(a_i.b_i)$ are eigenvalues, we have adding these formulas
$$
\left(E((0,\infty)\setminus\cup_{\lambda_n\in\sigma_p(H)}\{\lambda_n\})f,f\right) = 
\int_{0}^{\infty}\|{\mathcal F}^{(\pm)}(k)f\|^2dk.
$$
Let $P_{ac}(H)$ be the projection onto the absolutely continuous subspace for $H$. Then
$$
E((0,\infty)\setminus\cup_{\lambda_n\in\sigma_p(H)}\{\lambda_n\}) = 
P_{ac}(H).
$$
Therefore, we have
$$
\left(P_{ac}(H)f,f\right) = 
\int_{0}^{\infty}\|{\mathcal F}^{(\pm)}(k)f\|^2dk,
$$
which proves (2), (3).

Let $f \in C_0^{\infty}(\mathcal M)$. By (\ref{C3S3Fpmjkdefine}), (\ref{C3S3Qjz}) and Theorem 2.1 (2), we have
\begin{equation}
\begin{split}
\mathcal F_j^{(\pm)}(k)(H - k^2)f &= \mathcal F_{free(j)}^{(\pm)}(k)Q_j(k^2 \pm i0)(H - k^2)f \\
&= \mathcal F^{(\pm)}_{free(j)}(k)(H_{free(j)} - k^2)\chi_jf = 0.
\end{split}
\nonumber
\end{equation}
To prove (4) for $f \in D(H)$, we have only to approximate it by a sequence in $C_0^{\infty}(\mathcal M)$.

Theorem 3.8 and Lemma 3.11 imply that $\mathcal F^{(\pm)}(k) \in {\bf B}(\mathcal B;{\bf h}_{\infty})$. Therefore, $\mathcal F^{(\pm)}(k)^{\ast} \in {\bf B}({\bf h}_{\infty};\mathcal B^{\ast})$. This and (4) yield (5).

To prove (1). we have only to show that ${\rm Ran}\,{\mathcal F}^{(\pm)}$ is dense in $\widehat{\mathcal H}$, since ${\rm Ran}\,{\mathcal F}^{(\pm)}$ is closed by (2), (3). The idea is the same as the case of Lemma 1.3.19. For the sake of notational simplicity, we assume that there are only 2 ends, $\mathcal M_1$ with regular infinity and $\mathcal M_2$ with cusp. Suppose 
$$
(\varphi_1(k),\varphi_2(k)\varphi_{2,0}) \in {\bf h}_{\infty} = L^2((0,\infty);L^2(M_1);dk)\times
L^2((0,\infty);{\bf C};dk),
$$
where 
$\varphi_{2,0} = |M_2|^{-1/2}$ is the eigenfunction of $\Delta_{M_2}$ associated with zero eigenvalue,
 is orthogonal to 
${\rm Ran}\,\mathcal F^{(+)}$.
Let $\{e_1,e_2,\cdots\}$ be a complete orthnormal system of $L^2(M_1)$, and put
$$
\varphi_{1,n}(k) = (\varphi_1(k),e_n)_{L^2(M_1)}.
$$
Let $\mathcal L(\psi)$ be the set of Lebesgue points of $\psi \in L^1_{loc}((0,\infty))$ introduced in the proof of Lemma 1.3.19. We take
$$
\ell \in \Big(\cap_{n=1}^{\infty}\mathcal L(\varphi_{1,n})\Big)\cap \left(
\mathcal L(\|\varphi_{1}(k)\|_{L^2(M_1)}^2)\right)\cap\Big(\mathcal L(\varphi_{2})\Big)\cap\Big( 
\mathcal L(|\varphi_{2}|^2)\Big).
$$
Let $\{\chi_j\}_{j=0}^2$ be the partition of unity from (\ref{eq:Chap3Sec2Partitionofunity}). We fix $m$ arbitrarily, and put
\begin{equation}
\begin{split}
u_{\ell} & = 
\omega_{+}(\ell)\chi_1(y)y^{(n-1)/2-i\ell}\alpha e_m + \omega_{+}^{(c)}(\ell)\chi_j(y)y^{(n-1)/2+i\ell}\beta\varphi_{2,0},
\end{split}
\nonumber
\end{equation}
$\alpha$, $\beta$ being arbitraily chosen constants. We further put
$$
(H - \ell^2)u_{\ell} = g_{\ell}.
$$
Then, as can be checked easily, $g_{\ell} \in L^{2,(1 + \epsilon)/2}$, and by Theorems 3.8 
and 3.10, $u_{\ell}$ is written as $u_{\ell} = R(\ell^2+ i0)g_{\ell}$. Moreover, 
etting ${\mathcal F}^{(+)}(k)g_{\ell} = (C_1(k),C_2(k)\varphi_{2,0})$, we see that 
$(C_1(k),C_2(k)\varphi_{2,0})$ is an $L^2(M_1)\times{\bf C}$-valued continuous function 
of $k > 0$, satisfying
\begin{equation}
(C_1(\ell),e_n) = \delta_{mn}\alpha, \quad \quad C_2(\ell) = \beta.
\label{C3S3C1ellen}
\end{equation}
By our assumption, $(\varphi_1(k),\varphi_2(k)\varphi_{2,0})$ is orthogonal to ${\mathcal F}^{(+)}(k)E_H(I)g_{\ell} $, $I$ being any interval of $(0,\infty)$. Hence,
$$
\int_I\left((\varphi_1(k),C_1(k))_{L^2(M_1)} + \varphi_2(k)\overline{C_2(k)}\right)dk = 0
$$
for any interval $I \subset (0,\infty)$.
By the same arguments as in the proof of Lemma 1.3.19, we then have 
$$
\frac{1}{2\epsilon}\int_{\ell-\epsilon}^{\ell+\epsilon}\varphi_2(k)\overline{C_2(k)}dk \to \varphi_2(\ell)\overline{\beta}.
$$
The 1st term is computed as
\begin{equation}
\begin{split}
\frac{1}{2\epsilon}\int_{\ell-\epsilon}^{\ell + \epsilon}
(\varphi_1(k),C_1(k))_{L^2(M_1)}dk &=
\frac{1}{2\epsilon}\int_{\ell-\epsilon}^{\ell + \epsilon}
(\varphi_1(k),C_1(k)- C_1(\ell))_{L^2(M_1)}dk \\
& + \frac{1}{2\epsilon}\int_{\ell-\epsilon}^{\ell + \epsilon}
(\varphi_1(k),C_1(\ell))_{L^2(M_1)}dk.
\end{split}
\nonumber
\end{equation}
By (\ref{C3S3C1ellen}), $(\varphi_1(k),C_1(\ell))_{L^2(M_1)} = \varphi_{1,m}(k)\overline{\alpha}$, hence
$$
\frac{1}{2\epsilon}\int_{\ell-\epsilon}^{\ell + \epsilon}
(\varphi_1(k),C_1(\ell))_{L^2(M_1)}dk \to \varphi_{1,m}(\ell)\overline{\alpha}.
$$
We also have
\begin{equation}
\begin{split}
& \left|\frac{1}{2\epsilon}\int_{\ell-\epsilon}^{\ell + \epsilon}
(\varphi_1(k),C_1(k)- C_1(\ell))_{L^2(M_1)}dk\right| \\
& \leq \left(\frac{1}{2\epsilon}\int_{\ell-\epsilon}^{\ell + \epsilon}
\|\varphi_1(k)\|^2_{L^2(M_1)}dk\right)^{1/2}
 \times \left(\frac{1}{2\epsilon}\int_{\ell-\epsilon}^{\ell + \epsilon}
\|C_1(k)- C_1(\ell)\|^2_{L^2(M_1)}dk\right)^{1/2}.
\end{split}
\nonumber
\end{equation}
The right-hand side tends to 0, since $\ell$ is an Lebesgue point of $\|\varphi_1(k)\|_{L^2(M_1)}^2$, and $C_1(k)$ is an $L^2(M_1)$-valued continuous function of $k > 0$. We have, therefore, obtained that
\begin{equation}
\varphi_{1,m}(\ell)\overline{\alpha} + \varphi_2(\ell)\overline{\beta} = 0.
\nonumber
\end{equation}
Since $\alpha$, $\beta$ and $m$ are arbitrarily, we have $\varphi_1(\ell) = 0$, $\varphi_2(\ell) = 0$, which completes the proof of (1). The proof of (6) is the same as Theorem 1.3.13. \qed


\subsection{$S$ matrix} 
As in Chap. 2, we can prove the existence and completeness of time-dependent wave operators and introduce the Radon transform associated with $H$. We give a breif sketch of the proof later. Here, instead of this time-depedent approach, we construct the S-matrix by using the generalized Fourier transform.

The following theorem is proved in the same way as Theorem 1.4.3 with ${\mathcal F}_0^{(\pm)}(k)$ replaced by ${\mathcal F}^{(\pm)}(k)$, and is a generalization of the modified {\it Poisson-Herglotz} formula.


\begin{theorem}
If $k^2 \not\in \sigma_p(H)$, we have 
\begin{equation}
{\mathcal F}^{(\pm)}(k){\mathcal B} = {\bf h}_{\infty},
\nonumber
\end{equation}
\begin{equation}
\{u \in {\mathcal B}^{\ast}\,;\,(H - k^2)u = 0\} = 
{\mathcal F}^{(\pm)}(k)^{\ast}{\bf h}_{\infty}.
\nonumber
\end{equation}
\end{theorem}

We derive an asymptotic expansion of solutions to the Helmholtz equation. Let $V_j$ be the differential operator defined by
\begin{equation}
 V_j =  [H_{free(j)},\chi_j] - \chi_j\widetilde V_j \quad 
 (1 \leq j \leq N),
 \nonumber
\end{equation}
where $\widetilde V_j$ is defined by (\ref{C3S3widetildeVj}). 
We put
\begin{equation}
J_j(k) = \sum_{m=1}^{\infty}
\left(\frac{\sqrt{\lambda_{j,m}}}{2}\right)^{-2ik}
P_{j,m} = 
\left(\frac{\sqrt{-\Delta_{M_j}}}{2}\right)^{-2ik}(I - P_{j,0}),
\label{C3S3Jpk}
\end{equation}
where $\Delta_{M_j}$ is the Laplace-Beltami operator on $M_j$ and 
$P_{j,0}$ is the projection onto the zero eigenspace for $\Delta_{M_j}$.
For $1 \leq j, l \leq N$, we define 
$\widehat S_{jl}(k) \in {\bf B}(L^2(M_l);L^2(M_j))$ by
\begin{equation}
 \widehat S_{jl}(k) = 
 \left\{
 \begin{split}
 & \delta_{jl}J_j(k) + \frac{\pi i}{k}
 {\mathcal F}_j^{(+)}(k)\big(V_l\big)^{\ast}\left({\mathcal F}_{free(l)}^{(-)}(k)\right)^{\ast},  \quad 1 \leq j \leq \mu,
\\
&  \frac{\pi i}{k}{\mathcal F}_j^{(+)}(k)\big(V_l\big)^{\ast}\left({\mathcal F}_{free(l)}^{(-)}(k)\right)^{\ast}, \quad \mu + 1 \leq j \leq N.
\end{split}
\right.
\label{C3S3widehatSpqk}
\end{equation}


\begin{theorem}
For $\psi = (\psi_1,\cdots,\psi_N) \in {\bf h}_{\infty}$, the following asymptotic expansion holds:
\begin{equation}
 \begin{split}
    & \big({\mathcal F}^{(-)}(k)\big)^{\ast}\psi = \sum_{j=1}^N\big({\mathcal F}_j^{(-)}(k)\big)^{\ast}\psi_j \\
  &\simeq \ \ \frac{ik}{\pi} \omega_-(k)\sum_{j=1}^{\mu}\,\chi_j \,
  y^{(n-1)/2+ik}\,\psi_j + \frac{ik}{\pi} \omega_-^{(c)}(k)\sum_{j=\mu+1}^N
   \, \chi_j\, y^{(n-1)/2-ik}\,\psi_{j} \\
    &\  \ - \frac{ik}{\pi} \omega_+(k)\sum_{j=1}^{\mu}\sum_{l=1}^N
  \,\chi_j \,y^{(n-1)/2-ik}\,
  \widehat S_{jl}(k)\psi_l  \\
& \ \ - \frac{ik}{\pi} \omega_+^{(c)}(k)\sum_{j=\mu +1}^N\sum_{l=1}^N
  \,\chi_j \,y^{(n-1)/2+ik}\,
  \widehat S_{jl}(k)\psi_l.
   \end{split}
   \nonumber
\end{equation}
\end{theorem}
Proof. 
First note that by (\ref{C3S3Fpmjkdefine})
\begin{equation}
 \left({\mathcal F}_j^{(-)}(k)\right)^{\ast} = 
 \chi_j \left({\mathcal F}_{free(j)}^{(-)}(k)\right)^{\ast}
 + R(k^2 + i0)\big(V_j\big)^{\ast}\left({\mathcal F}_{free(j)}^{(-)}(k)\right)^{\ast}.
 \label{C3S3Fjminuskadjoint}
\end{equation}
By (\ref{C3S3Ffreejpmk}), for $1 \leq j \leq \mu$,
\begin{eqnarray*}
 \left({\mathcal F}_{free(j)}^{(-)}(k)\right)^{\ast}\phi 
 & =& \sum_{m=0}^{\infty}  \overline{C_{j,m}^{(-)}(k)}\,
\left(F_{free(j),m}^{(-)}(k)\right)^{\ast}\,P_{j,m}\phi \\
&=& \overline{C_{j,0}^{(-)}(k)}\,\frac{1}{\sqrt{2\pi}}\,y^{(n-1)/2 + ik}P_{j,0}\phi \\
&+& \sum_{m=1}^{\infty}\overline{C_{j,m}^{(-)}(k)}\,\frac{(2k\sinh(k\pi))^{1/2}}{\pi}\,
y^{(n-1)/2}K_{ik}(\sqrt{\lambda_{j,m}}\,y)\,P_{j,m}\phi,
\end{eqnarray*}
and by (\ref{C3S3Ffreejpmkcusp}), for $\mu + 1 \leq j \leq N$,
\begin{equation}
 \left({\mathcal F}_{free(j)}^{(-)}(k)\right)^{\ast}\phi 
 = \frac{1}
 {\sqrt{2\pi}}y^{(n-1)/2-ik}\phi.
\label{C3S3FfreemujN}
\end{equation}

Since $\mathcal F^{(-)}(k)^{\ast} \in {\bf B}({\bf h}_{\infty};\mathcal B^{\ast})$, we have only to prove the theorem for $\psi = (\psi_1,\cdots,\psi_N) \in {\bf h}_{\infty}$ such that 
for $1 \leq j \leq \mu$, $P_{j,m}\psi_j = 0$ except for a finite number of $m$.
By using Chap. 1, (\ref{eq:Chap1Sec3Knunear0}), (\ref{eq:Chap1Sect4omegaplusminusk}) and (\ref{eq:modulusomegasquaredpi2k}), for $1 \leq j \leq \mu$, one can show
\begin{equation}
\begin{split}
\left({\mathcal F}_{free(j)}^{(-)}(k)\right)^{\ast}\psi_j 
&\simeq
\frac{ik}{\pi} \omega_-(k)y^{(n-1)/2+ik}\psi_j \\
& \ \  - \frac{ik}{\pi} \omega_+(k)y^{(n-1)/2-ik}
\sum_{m\geq1}
\left(\frac{\sqrt{\lambda_m}}{2}\right)^{-2ik}P_{j,m}\psi_j.
\end{split}
\label{C3S3Ffree1jmu}
\end{equation}
We apply (\ref{C3S3FfreemujN}) and (\ref{C3S3Ffree1jmu}) to the 1st term of the right-hand side of (\ref{C3S3Fjminuskadjoint}). To the 2nd term, we apply
Theorem 3.10. We then have, for $1 \leq j \leq \mu$,
\begin{eqnarray*}
\left(\mathcal F_j^{(-)}(k)\right)^{\ast}\psi_j  & \simeq &
\frac{ik}{\pi}\omega_-(k)\, \chi_j\, y^{(n-1)/2+ik}\, \psi_j \\
& & - \frac{ik}{\pi}\omega_+(k)\sum_{l=1}^\mu\,\chi_l\,y^{(n-1)/2-ik}\,\widehat S_{lj}(k)\psi_j \\
& & - \frac{ik}{\pi}\omega_+^{(c)}(k)\sum_{l=\mu+1}^N\,\chi_l\,y^{(n-1)/2+ik}\,
\widehat S_{lj}(k)\psi_j.
\end{eqnarray*}
Similary, one can show, for $\mu + 1 \leq j \leq N$,
\begin{eqnarray*}
\left(\mathcal F_j^{(-)}(k)\right)^{\ast}\psi_j  & \simeq &
 \frac{ik}{\pi}\omega_-^{(c)}(k)\, \chi_j\, y^{(n-1)/2-ik}\, \psi_j \\
& & - \frac{ik}{\pi}\omega_+(k)\sum_{l=1}^\mu\,\chi_l\,y^{(n-1)/2-ik}\,\widehat S_{lj}(k)\psi_j \\
& & - \frac{ik}{\pi}\omega_+^{(c)}(k)\sum_{l=\mu+1}^N\,\chi_l\,y^{(n-1)/2+ik}\,
\widehat S_{lj}(k)\psi_j.
\end{eqnarray*}
Summing up these two formulas, we obtain the theorem. \qed

\bigskip
We define an operator-valued $N \times N$ matrix $\widehat S(k)$ by
\begin{equation}
\widehat S(k) = \Big(\widehat S_{jl}(k)\Big),
\label{C3S3Smatrix}
\end{equation}
and call it  {\it $S$-matrix}. This should be more properly called the geometric S-matrix in the context of Chap. 2, \S 6. This is a bounded operator on ${\bf h}_{\infty}$. 
Similarly to Theorem 2.7.9, we have the following asymptotic expansion.


\begin{theorem}
(1)  For any $u \in {\mathcal B}^{\ast}$ satisfying $(H - k^2)u = 0$, there exists a unique $\psi^{(\pm)} = (\psi^{(\pm)}_1,\cdots,\psi^{(\pm)}_N) \in {\bf h}_{\infty}$ such that
\begin{equation}
 \begin{split}
    u &\simeq \omega_-(k)\sum_{j=1}^{\mu}\,\chi_j\, 
  y^{(n-1)/2+ik}\,\psi_j^{(-)}  + \omega_-^{(c)}(k)\sum_{j={\mu}+1}^N
   \,\chi_j \, y^{(n-1)/2-ik}\, \psi_{j}^{(-)} \\
   & -  \omega_+(k)\sum_{j=1}^{\mu}
  \chi_j\, y^{(n-1)/2-ik}\,
  \psi_j^{(+)} -  \omega_+^{(c)}(k)\sum_{j={\mu}+1}^N
  \,\chi_j\, y^{(n-1)/2+ik}\,
  \psi_{j}^{(+)}.
   \end{split}
   \nonumber
\end{equation}
(2) For any $\psi^{(-)} \in {\bf h}_{\infty}$, there exists a unique $\psi^{(+)} \in {\bf h}_{\infty}$ and $u \in \mathcal B^{\ast}$ satisfying $(H - k^2)u = 0$, for which the expansion (1) holds. Moreover
\begin{equation}
\psi^{(+)} = \widehat S(k)\psi^{(-)}.
\nonumber
\end{equation}
\end{theorem}

Proof. By Theorem 3.13, $u \in {\mathcal F}^{(-)}(k)^{\ast}{\bf h}_{\infty}$. Using Theorem 3.14, we prove the result. \qed


\begin{theorem}
$\widehat S(k)$ is unitary on ${\bf h}_{\infty}$.
\end{theorem}
Proof. 
Let $u \in \mathcal B^{\ast}$ such that $(H - k^2)u = 0$. By Theorem 3.13, $u = \mathcal F^{(+)}(k)^{\ast}\psi^{(+)}$, $\psi^{(+)} \in {\bf h}_{\infty}$.
By similar arguments as in Theorem 3.14, with $\mathcal F^{(+)}(k)^{\ast}$ instead of $\mathcal F^{(-)}(k)^{\ast}$, one can show that there exists $\psi^{(-)} \in {\bf h}_{\infty}$  such that  the expansion in Theorem 3.15 (1) holds.  In particular, $\psi^{(+)} =\widehat S(k) \psi^{(-)} $.
This means that $\widehat S(k)$ is onto. 

Thus, we have only to prove that $\widehat S(k)$ is isometric.
Take $\psi^{(-)} = (\psi^{(-)}_1,\cdots,\psi^{(-)}_N) \in {\bf h}_{\infty}$ such that for $1 \leq j \leq \mu$, $P_{j,m}\psi^{(-)}_j = 0$ except for a finite number of $m$. We put for $1 \leq j \leq \mu$
\begin{equation}
 a_{j,m} = \left\{
 \begin{array}{lc}
 P_{j,0}\psi_j^{(-)}, & (m = 0) \\
 \displaystyle{\Big(\frac{\sqrt{\lambda_{j,m}}}{2}\Big)^{-ik}
 \Gamma(1 + ik)P_{j,m}\psi_j^{(-)}},
 & (m \neq 0)
 \end{array}
 \right.
\nonumber
\end{equation}
\begin{equation}
u_j^{(-)} =  \omega_-(k)\, \chi_j\,\Big(
y^{(n-1)/2+ ik}a_{j,0}
+ \sum_{m\geq1}
y^{(n-1)/2}I_{ik}(\sqrt{\lambda_{j,m}}\,y)a_{j,m}
\Big).
\nonumber
\end{equation}
Then, as $y \to 0$, 
\begin{equation}
 u_j^{(-)} \simeq  \omega_-(k)\, \chi_j(y)\, y^{(n-1)/2 + ik}
 \psi_j^{(-)}.
 \nonumber
\end{equation}
For $\mu+1 \leq j \leq N$, we put
\begin{equation} \label{E2.13}
 u_j^{(-)} = \omega_-^{(c)}(k)\, \chi_j\, y^{(n-1)/2 - ik}
 \psi_j^{(-)},
\end{equation}
and define
\begin{eqnarray*}
 u^{(-)} &=& \sum_{j=1}^Nu_j^{(-)},  \quad
 f = (H - k^2)u^{(-)} \in {\mathcal B},\\
 u^{(+)} &=& R(k^2 + i0)f, \quad
 u = u^{(+)} - u^{(-)}, \\
 \psi^{(+)} &=& {\mathcal F}^{(+)}(k)f.
\end{eqnarray*}
Then, by Theorem 3.10, $u$ and $\psi^{(\pm)}$ give the expansion in Theorem 3.15 (1).
Lemma 3.11 implies
\begin{equation}
\begin{split}
\frac{1}{2k}\|\psi^{(+)}\|^2 &= \frac{1}{2\pi i}\left(R(k^2 + i0)f - 
R(k^2 - i0)f,f\right) \\
&= \frac{1}{2\pi i}\left[(f,u^{(-)}) - (u^{(-)},f)\right].
\end{split}
\nonumber
\end{equation}
Here we have used the fact that
\begin{equation} \label{E2.14}
 R(k^2 - i0)f =  u^{(-)},
 \nonumber
\end{equation}
since $u^{(-)}$ is incoming. Now we do the same computation as in Lemma 3.11. Let $\chi_R$ be as in the lemma. Then,
\begin{equation}
\begin{split}
 (f,\chi_Ru^{(-)}) - (\chi_Ru^{(-)},f) & =  ([H,\chi_R]u^{(-)},u^{(-)}) \\
& =  \sum_{j=1}^N(\chi_j[H_{free(j)},\chi_R]\widetilde\chi_j\,u^{(-)},u^{(-)}) + o(1).
\end{split}
\nonumber
\end{equation}
Recall that 
$$
[H_{free(j)},\chi_R] = - \frac{2}{\log R}\chi'\big(\frac{\log y}{\log R}\big)
\big(D_y - \frac{n-1}{2}\big) + O(|\log R|^{-2}).
$$
Then, for $1 \leq j \leq \mu$, using the fact that $u^{(-)}$ has the form (\ref{E2.13}), we have  
\begin{eqnarray*}
  (\chi_j[H_{free(j)},\chi_R]\widetilde\chi_j\,u^{(-)},u^{(-)}) 
 &=& \frac{2ik}{\log R}
 \left(\chi'\Big(\frac{\log y}{\log R}\Big)u_j^{(-)},u_j^{(-)}\right) + o(1) \\
 &=& \frac{2ik}{\log R} |\omega_-(k)|^2
 \int_0^{1}\chi'\Big(\frac{\log y}{\log R}\Big)\frac{dy}{y}
 \|\psi^{(-)}_j\|^2 + o(1)  \\
 &=& \frac{\pi i}{k}\|\psi_j^{(-)}\|^2 + o(1),
\end{eqnarray*}
where, at the last step, we use equation (\ref{eq:Chap2Sect4InversionFormula})  of Ch. 1.

Similarly, for $\mu +1 \leq j \leq N$,
\begin{equation}
 (\chi_j[H_{free(j)},\chi_R]\widetilde\chi_j\,u^{(-)},u^{(-)})= \frac{\pi i}{k}
 \|\psi_j^{(-)}\|^2 + o(1).
 \nonumber
\end{equation}
Taking $R \to \infty$, we obtain $\|\psi^{(+)}\| = \|\psi^{(-)}\|$. \qed

\begin{cor}
 $\ {\mathcal F}^{(+)}(k) = \widehat S(k){\mathcal F}^{(-)}(k)$.
\end{cor}
Proof. The above $f$ satisfies $\psi^{(\pm)} = {\mathcal F}^{(\pm)}
(k)f$. Since $\psi^{(+)} = \widehat S(k)\psi^{(-)}$ and, by (\ref{E2.14}),
$\psi^{(-)}={\mathcal F}^{(-)}(k) f$, the corollary is proved.\qed


\subsection{Wave operators}
We briefly  look at the temporal asymptotics of $e^{-it\sqrt{H}}f$ for $f \in \mathcal H_{ac}(H)$. Let $\{\chi_j\}_{j=0}^N$ be the partition of unity given in Subsection 3.2. We can then show that
\begin{equation}
\|\chi_0e^{-it\sqrt{H}}f\| \to 0, \quad {\rm as} \quad t \to \pm \infty.
\label{C3S3chi0u(t)to0} 
\end{equation}
In fact, by approximating $f$, we have only to consider the case that $f \in D(H)\cap{\mathcal H}_{ac}(H)$. In this case, we have $\chi_0e^{-it\sqrt{H}}f = \chi_0(H+i)^{-1}e^{-it\sqrt{H}}(H+i)f$. 
Since $(H+i)f \in \mathcal H_{ac}(H)$, we have $\chi_0e^{-it{H}}(H+i)f \to 0$ weakly as 
$t \to \pm \infty$. 
{\ntext As also $\chi_0 (H+i)^{-1}$ is compact,} this proves (\ref{C3S3chi0u(t)to0}). It then  implies
\begin{equation}
\|e^{-it\sqrt{H}}f- \sum_{j=1}^N \chi_je^{-it\sqrt{H}}f\| \to 0, \quad {\rm as} \quad t \to \pm \infty.
\nonumber
\end{equation}
Consider the behavior of $\chi_je^{-it\sqrt{H}}f$ on the end $\mathcal M_j$. Suppose $\mathcal M_j$ is a regular end. Then the argument in Chapter 2 Subsection 8.3 works well without any essential change, and one can show that, as $t \to \infty$,
\begin{equation}
\left\|\chi_je^{-it\sqrt{H}}f - \chi_j\frac{y^{(n-1)/2}}{\sqrt{\pi}}
\int_0^{\infty}e^{ik(-\log y - t)}\Big(\mathcal F_j^{(+)}f\Big)(k)dk\right\| \to 0.\nonumber
\end{equation} 
Similarly, for $g \in L^2({\mathcal M}_j)$,
\begin{equation}
\left\|\chi_je^{-it\sqrt{H_{free(j)}}}g - \chi_j\frac{y^{(n-1)/2}}{\sqrt{\pi}}
\int_0^{\infty}e^{ik(-\log y - t)}\Big(\mathcal F_{free(j)}^{(+)}g\Big)(k)dk\right\| \to 0.
\nonumber
\end{equation}
Taking $g = \big(\mathcal F^{(+)}_{free(j)}\big)^{\ast}\mathcal F_j^{(+)}f$, 
these two limits imply
\begin{equation}
\chi_j e^{- it\sqrt{H}}f \sim 
\chi_j e^{-it\sqrt{H_{free(j)}}}
\Big(\mathcal F_{free(j)}^{(+)}\Big)^{\ast}\mathcal F_{j}^{(+)}f.
\nonumber
\end{equation} 
We can prove  similar formulae when $\mathcal M_j$ is a cusp.
This means that, in the long-run, the waves disappear from compact parts of the manifold, and, on each end, they behave like free waves. 

Similarly, we can prove
\begin{equation}
{\mathop{\rm s-lim}_{t\to\infty}}\, e^{it\sqrt{H}}\chi_je^{-it\sqrt{H_{free(j)}}}
= \big(\mathcal F_j^{(+)}\big)^{\ast}\mathcal F_{free(j)}^{(+)},
\nonumber
\end{equation}
and, therefore, there exist the wave operators,
\begin{equation}
W_{\pm} = \mathop{\rm s-lim}_{t\to\pm\infty}\sum_{j=1}^Ne^{it\sqrt{H}}\chi_je^{-it\sqrt{H_{free(j)}}} = \sum_{j=1}^N\big(\mathcal F_j^{(+)}\big)^{\ast}\mathcal F_{free(j)}^{(+)}.
\label{C3S3WaveOpandFourier}
\end{equation}
Since $\mathcal F_{free(j)}^{(+)}$ are unitary, it follows from Theorem 3.12, that and $W_{\pm}$ are complete: 
\begin{equation}
{\rm Ran}\,W_{\pm} = \mathcal H_{ac}(H).
\nonumber
\end{equation}

As in Chap. 2, \S 8, we construct $\mathcal F_{\pm}$ from $\mathcal F^{(\pm)}$, and define the Radon transform by the formula
\begin{equation}
\left(\mathcal R_{\pm}f\right)(s) = 
\frac{1}{\sqrt{2\pi}}\int_{-\infty}^{\infty}
e^{iks}\left(\mathcal F_{\pm}f\right)(k)dk.
\nonumber
\end{equation}
Then Theorem 2.8.9 also holds on $\mathcal M$.


\section{Cusps and generalized eigenfunctions}
In the following two sections, we consider the case in which ${\mathcal M}$ has only cusps as infinity. We use the same notation as in the previous section, and  for the sake of simplicity assume that $\mathcal M$ has only one cusp and the manifold at infinity $M$  satisfies $|M| = 1$. 
In this section $z$ denotes a point in ${\mathcal M}$. Moreover, we assume:

\medskip
\noindent
{\bf (C-1)} {\it The end $\mathcal M_1$ is identified with $M \times (1,\infty)$ and the metric of $\mathcal M$ is}
\begin{equation} \label{E4.1}
ds^2 = \sum_{i, j=1}^n g_{ij}(z) dz^i dz^j=\frac{(dy)^2 + h(x,dx)}{y^2} \quad {on} \quad \mathcal M_1,
\end{equation}
where we typically use local coordinates $z= (x, y)$, $x=(x_1, \dots, x_{n-1})$ being local coordinates
on $M$.


\subsection{A remark on the S-matrix}
In Theorem 3.15, we have proven that for $k > 0$ such that $k^2 \not\in \sigma_p(H)$ and $u \in {\mathcal B}^{\ast}$ satisfying $(H - k^2)u = 0$, there exist unique constant functions $\psi^{(\pm)} \in P_0L^2(M)$ such that
\begin{equation}
u \simeq \omega_-^{(c)}(k)y^{(n-1)/2-ik}\psi^{(-)} - \omega_+^{(c)}(k)y^{(n-1)/2 + ik}\psi^{(+)}, \quad \omega_{\pm}^{(c)}(k) = \pm \frac{i}{k}\sqrt{\frac{\pi}{2}}.
\label{C3S4expansion}
\end{equation}
$w_{\pm}^{(c)}(k)$ has natural extension to $k < 0$. Then taking $u(k) = \overline{u(-k)}$, we obtain, for $k < 0$, a solution to $(H-k^2)u = 0$ which also satisfies (\ref{C3S4expansion}).  With this in mind, we change the notion of the S-matrix as follows. Let
$$
\mathcal N(k) = \left\{u \in \mathcal B^{\ast}\, ; \, \Big(- \Delta_g - \frac{(n-1)^2}{4} - k^2\Big)u = 0\right\}.
$$
Then, for any $0 \neq k \in {\bf R}$, such that $k^2 \not\in \sigma_p(H)$, ${\rm dim}\,\mathcal N(k) =1$, and one can choose a basis $u(z,k) \in \mathcal N(k)$ satisfying
\begin{equation}
u \simeq y^{(n-1)/2-ik} + \widehat S(k)y^{(n-1)/2 + ik},
\label{C3S4expansion2}
\end{equation}
$\widehat S(k)$ being a complex number of modulus 1. 
Traditionally, we put
\begin{equation}
\mathcal S(s) = \widehat S(k), \quad s = (n-1)/2 - ik,
\label{C3S4mathcalSswidehatSk}
\end{equation}
and call it the S-matrix.


\subsection{Eisenstein series}  We put
\begin{equation}
\sqrt{\sigma_p(H)}  = \{\zeta \in {\bf C}\, ; \, \zeta^2 \in \sigma_p(H)\}.
\nonumber
\end{equation}
Let $\chi \in C^{\infty}((0,\infty))$ be such that $\chi(y) = 0$ for $y < 2$, $\chi(y) = 1$ for $y > 3$.
We define for $k > 0$ and $\epsilon > 0$
\begin{equation}
\varphi(z,k - i\epsilon) = \chi(y)\,y^{\frac{n-1}{2} + i(k - i\epsilon)} - R((k - i\epsilon)^2)\,[H,\chi]\,y^{\frac{n-1}{2} + i(k - i\epsilon)}.
\label{eq:varphixk}
\end{equation}
Due to (C-1), $\hbox{supp}\left([H, \chi] \right) \subset M \times (2, 3)$ and
this function $\varphi$ satisfies
$$
(H - (k -i\epsilon)^2)\varphi(z,k-i\epsilon) = 0.
$$

By the reasoning to be explained in the next section, this function is called an {\it Eisenstein series}. 
As a function of $k - i\epsilon$, this is meromorphic in the lower-half plane and has poles at $\sqrt{\sigma_p(H)}\cap{\bf C_-}$. 
Note that in the standard notation, we put $s = (n-1)/2 + i(k - i\epsilon)$ and regard $\varphi$ as a meromorphic function on $\{s \in {\bf C}\, ; \, {\rm Re}\,s > (n-1)/2\}$. 
By the limiting absorption principle, letting $\epsilon \to 0$, $\varphi(z,k - i\epsilon)$ is continuously extended to ${\bf R} \setminus \sqrt{\sigma_p(H)}$.

Using the definitions (\ref{C3S3Fpmjkdefine}), 
(\ref{eq:Chap3Sect2Focplusminus}), (\ref{C3S2Ffree0}), and (\ref{C3S3Qjz}) with $\widetilde V = 0$, we have, for $k \in (0,\infty)\setminus\sqrt{\sigma_p(H)}$,
\begin{equation}
{\mathcal F}^{(+)}(k)f = 
 \frac{1}{\sqrt{2\pi}} \int_{{\mathcal M}}
\overline{\varphi(z,k)}f(z)\,d{\mathcal M}.
\nonumber
\end{equation}
Hence, by Theorem 3.12 we have the following theorem.


\begin{theorem} $\mathcal F^{(+)}$ maps $\mathcal H_{ac}(H)$ onto $L^2((0,\infty)\,;\, P_0(L^2(M))\,;\,dk)$.
For any $f \in L^2({\mathcal M})$, the inversion formula holds:
\begin{equation}
f(z) = \frac{1}{\sqrt{2\pi}}\int_{0}^{\infty}
\varphi(z,k)\widetilde f(k)dk
+ \sum_{i}(f,\psi_i)\psi_i,
\nonumber
\end{equation}
\begin{equation}
 \widetilde f(k) =  \frac{1}{\sqrt{2\pi}} \int_{{\mathcal M}}
\overline{\varphi(z,k)}f(z)d\mathcal M,
\nonumber
\end{equation}
where $\psi_i$ is a normalized eigenvector of $H$.
\end{theorem}


\subsection{Theory of quadratic forms}
Let us recall the theory of quadratic forms associated with self-adjoint extensions of symmetric operators. For the details, see e.g.  \cite{Ka76} p. 322 or \cite{Is04a}, p. 38. Let $D$ be a dense subspace of a Hilbert space ${\mathcal H}$. A hermitian quadratic form $a(\cdot,\cdot)$ with domain $D$ is 
a mapping : $D\times D$ $\to$ $\bf C$ satisfying
$$
a(\lambda u + \mu v,w) = \lambda a(u,w) + 
\mu a(v,w), \quad \lambda, \mu \in {\bf C}, \quad 
u, v, w \in D
$$
$$
\overline{a(u,v)} = a(v,u), 
\quad u, v \in D.
$$
A hermitian quadratic form $a(\cdot,\cdot)$ is said to be positive definite if there exists a constant $C > 0$, such that
$$
a(u,u) \geq C\|u\|^2,
\quad u \in D.
$$
In this case $a(\cdot,\cdot)$ defines an inner product on
$D$. If $D$ is complete with respect to the norm $\|u\|_a = 
\sqrt{a(u,u)}$, $a(\cdot,\cdot)$  is said to be a closed form.
We say that $a(\cdot,\cdot)$ is closable if, for any sequence $u_n \in D$ such that 
$\|u_n\| \to 0, \|u_n - u_m\|_a \to 0$,  we have
$\|u_n\|_a \to 0$. 
For a closable form $a(\cdot,\cdot)$, we define a subspace $\widetilde D$ by
$$
u \in \widetilde D \Longleftrightarrow
\exists u_n \in D \ {\rm s.t.}
\|u_n - u\| \to 0, \|u_n - u_m\|_a \to 0.
$$
For $u, v \in \widetilde D$, there exist
$u_n, v_n \in D$ such that $u_n \to u$, $v_n \to v$, 
$\|u_n - u_m\|_a \to 0$, $\|v_n - v_m\|_a \to 0$. Then, the quadratic form, defined by 
$$
\widetilde a(u,v) = \lim_{m,n\to\infty},
a(u_m,v_n)
$$
can be shown to be positive defnite and closed and is called
the closed extension of $a(\cdot,\cdot)$. Then the following theorem holds.


\begin{theorem} Let
$a(\cdot,\cdot)$ be a positive definite closed form with domain $D$. 
Then there exists a unique self-adjoint operator $A$ such that
$D(A) \subset D$ and
$$
a(u,v) = (Au,v), \quad u \in D(A),\quad  v \in D.
$$
Moreover
$D = D(A^{1/2})$.
\end{theorem}

A quadratic form $a(\cdot,\cdot)$ with domain $D$ is said to be bounded from below if there exists a constant $C_0 \geq 0$ such that
$$
a(u,u) \geq - C_0\|u\|^2, \quad 
\forall u \in D.
$$
In this case the quadratic form $b(\cdot,\cdot)$ defined by
$$
b(u,v) = a(u,v) + (C_0 + 1)(u,v)
$$
is positive definite. $a(\cdot,\cdot)$ is said to be closable if so is
$b(\cdot,\cdot)$.
Let $\widetilde b(\cdot,\cdot)$ be the closed extension of $b(\cdot,\cdot)$ . By Theorem 4.2, there exists a unique self-adjoint operator $B$ such that
$D(B) \subset \widetilde D$ and
$$
\widetilde b(u,v) = (Bu,v), \quad 
u \in D(B),\quad  v \in \widetilde D.
$$
Letting
$$
\widetilde a(u,v) = \widetilde b(u,v) - (C_0 + 1)(u,v),
$$
$$
A = B - (C_0 + 1), 
$$
we have $D(A) = D(B) \subset
\widetilde D$, and
$$
A \geq - C_0, \quad \widetilde a(u,v) = (Au,v), \quad u \in D(A), \quad 
v \in \widetilde D.
$$
We call $A$ the self-adjoint operator associated with $a(\cdot,\cdot)$.


\subsection{0-mode boundary value problem} 
We show that the Eisenstein series $\varphi(z,k)$ is meromorphically extended to ${\bf C}$ with respect to $k$. Following the arguments of \cite{Col81}, we consider the boundary value problem as below. 

Recall that $\mathcal M$ is assumed to be 
\begin{equation} \label{E4.11}
\mathcal M = \mathcal K \cup \mathcal M_1, \quad \mathcal M_1 = M\times(1,\infty), \quad |M| =1,
\end{equation}
where $\overline{\mathcal K}$ is compact. We can assume that 
$$
\mathcal K \cap \big(M\times(2,\infty)\big) = \emptyset.
$$
Take $a > 3$, and put
\begin{equation}
\mathcal M_{int}^a = \mathcal K\cup \big(M\times(1,a)\big), \quad 
\mathcal M_{ext}^a =M \times (a,\infty), \quad
\Gamma^a = M\times\{a\}.
\nonumber
\end{equation}
Using the projections $P_0$ and $P'$ on $L^2(M)$,
\begin{equation}
(P_0\psi)(x) = \int_{M}\psi(x')dM_{x'}, \quad 
P' = 1 - P_0,
\nonumber
\end{equation}
we define the following Hilbert space:
\begin{equation}
\mathcal H = L^2(\mathcal M^a_{int}) \oplus (P'\otimes I^a_y)L^2(\mathcal M^a_{ext})
\subset L^2(\mathcal M),
\nonumber
\end{equation}
 with $\pi: \,L^2(\mathcal M) \to \mathcal H$ being the associated orthogonal projection. Here, for any $b>0$,
$I_y^b$ is the cut-off  projector, in the $y$-coordinate, onto $y>b$.
To define the Sobolev spaces $H^m(\mathcal M)$, we use representation (\ref{E4.11}) of 
$\mathcal M$. Namely, if $U_l,\, l=1, \dots, L,$ is a coordinate covering of $ M$, we use, as
a coordinate covering of $\mathcal M$,
$$
\mathcal M = {\mathop\cup_{l=1}^{L+P}} {\mathcal U}_l,
$$
where ${\mathcal U}_l=U_l \times (1, \infty),\, l=1, \dots, L;\, \{{\mathcal U}_l\}_{l=L+1}^{L+P}$
being a coordinate covering of $\mathcal M^2_{int}$. Using the corresponding decomposition of unity,
$$
1= \sum_{l=1}^{L+P} \Psi_l(z),\quad \hbox{supp}(\Psi_l) \subset {\mathcal U}_l,
$$
where we assume, for $y>2$,  $\Psi_l(x, y)=\psi_l(x),\, {\rm supp}(\psi_l) \subset U_l,\, l=1, \dots, L$,
we define
$$
\|f\|^2_{H^m(\mathcal M)}= \sum_{l=1}^{L+P} \| \Psi_l f \|^2_{H^m(\mathcal U_l)}.
$$
Here $H^m(\mathcal U_l),\, l = L+1,\dots, L+P,$ are usual Sobolev spaces, while 
$$
\| \Psi_l f \|^2_{H^m(\mathcal U_l)}= 
\sum_{|\alpha| \leq m} \int_1^\infty \|D^{\alpha} \left (\Psi_l f  \right)  \|^2_{L^2(M)} \, \frac{dy}{y^n},
\quad l=1, \dots, L,
$$
where $D_i= y \partial_i,\, i=1, \dots, n$.

Note that, if $m=1$,  $\|f\|_{H^1}$ is equivalent to the classical invariant definition of $H^1$ on 
a Riemannian manifold,
\begin{equation}
\|f\|^2_{H^1(\mathcal M)} \sim \|f\|^2_{L^2(\mathcal M)}+ \int_{\mathcal M} |d f|^2_g d {\mathcal M}
 =  \|f\|^2_{L^2(\mathcal M)}+\int_{\mathcal M} \,g^{ij}\,\partial_i f\, \overline{\partial_j f}\,\sqrt{g}\,dz.
\label{eq:gradugradvinnerproduct}
\end{equation}

Next we define 
$$
{\mathcal H}^m:= \pi H^m(\mathcal M), \quad m \geq 1.
$$
Note that, with $I_M$ being identity on $M$ and $b>1$,
$
\left(I_M  \otimes I^b_y \right)f \in H^m(M \times (b, \infty)) 
$ iff
$$
\sum_j \int_b^\infty y^{2m} \left[ (1+\lambda_j^2)^m | {\hat f}_j(y)|^2 
+ |\partial_y^m {\hat f}_j(y)|^2 \right]\, \frac{dy}{y^n} < \infty.
$$
Here 
$f(x, y)= \sum_{j=0}  {\hat f}_j(y) \phi_j(x)$, for $y>b$.
Thus, 
$$
(P' \otimes I^b_y) {\mathcal H}^m \to  H^m(M \times (b, \infty)),\quad b>1.
$$

 Also, if $u \in \mathcal H^m$, then $\partial_y^j(P'\otimes I^b_y)u$, $0 \leq j \leq m-1$, is continuous across $\Gamma^a,\, a>b$.

We define a quadratic form $l(\cdot,\cdot)$ with domain $\mathcal H^1$ by
\begin{equation}
\begin{split}
l(u,u) &= (du,du)_{L^2(\mathcal M^a_{int})} 
+ \|u\|^2_{L^2(\mathcal M^a_{int})} \\
&+ (du,du)_{L^2(\mathcal M^a_{ext})} 
+ \|u\|^2_{L^2(\mathcal M^a_{ext})},
\end{split}
\nonumber
\end{equation}
see (\ref{eq:gradugradvinnerproduct}). 
Then $l(\cdot,\cdot)$ is a positive definite closed form on $\mathcal H^1$, and $\sqrt{l(\cdot,\cdot)}$ is equivalent to the $\mathcal H^1$-norm. Hence, by Theorem 4.2, there exists a unique self-adjoint operator $L$ such that $L \geq 1$, 
$ D(L^{1/2}) = \mathcal H^1$ and
\begin{equation}
 l(u,v) = (Lu,v), \quad \forall u \in D(L), \quad \forall v \in \mathcal H^1.
 \nonumber
\end{equation}

We introduce the set  $D_L$ by
\begin{equation}
D_L = \{u \in \mathcal H^2\, ; \, \left(\partial_y(P_0\otimes I_y)u\right)(a - 0) = 0\}.
\label{C3S4DL}
\end{equation}
Here, for $w \in H^1(M\times(a,a+1))$ or $w \in H^1(M\times(a-1,a))$,
 $w(a\pm 0)$ is defined by
\begin{equation}
w(a \pm 0) = \lim_{\epsilon \to 0}w(\cdot, a \pm \epsilon).
\nonumber
\end{equation}


\begin{lemma} (1) $L$ has compact resolvent. \\
\noindent
(2) $\ D(L) =  D_L$. \\
\noindent
(3)  If $\zeta \not\in \sigma(L)$, for any $f \in \mathcal H$ and $\lambda \in {\bf C}$, there exists a unique solution $u \in D_L$ of the following boundary value problem
\begin{equation}
\left\{
\begin{split}
&\big(- \Delta_g  - \frac{(n-1)^2}{4}  + 1 - \zeta\big)u = f \quad {in} \quad \mathcal M^a_{int}, \\
&\big(- \Delta_g  - \frac{(n-1)^2}{4}  + 1 - \zeta\big)(P'\otimes1)u = f \quad {in} \quad \mathcal M^a_{ext},
\\
&\big(\partial_y(P_0\otimes I_y)u\big)(a - 0) = \lambda. 
\end{split}
\right.
\label{C3S40modeBVP}
\end{equation}
The solution $u = u(z,\zeta,\lambda)$ is analytic with respect to $\lambda$ and meromorphic on ${\bf C}$ with respect to $\zeta$ with possible poles at $\sigma(L)$.
\end{lemma}
Proof.  By (\ref{E4.1}), if $y>1$, the inverse to $g_{ij}$ is,
For $y > 1$, the metric takes the form
\begin{equation}
\left(g_{ij}\right) = 
\left(
\begin{array}{cc}
h_{ij}(x)/y^2 & 0 \\
0 &1/ y^2
\end{array}
\right).
\nonumber
\end{equation}
Therefore, its inverse is
\begin{equation}
\left(g^{ij}\right) = 
\left(
\begin{array}{cc}
y^2h^{ij}(x) & 0 \\
0 & y^2
\end{array}
\right).
\nonumber
\end{equation}
To show the compactness of the resolvent, we have only to show that if  $\{u_j\}$ is a bounded sequence in $\mathcal H^1$, it contains a subsequence convergent in $\mathcal H$.
Since $P_m$ is the projection onto the the eigenspace corresponding to 
$m$-th eigenvalue $\lambda_m$ of $- \Delta_M$, we have, for $u \in H^1(\mathcal M)$ 
and $R > a$,
\begin{eqnarray*}
\int_{M\times(R,\infty)}g^{ij}\,\partial_iu\,\overline{\partial_ju}\,d\mathcal M&=& \int_{M\times(R,\infty)}y^2
\left(|\partial_yu|^2 + h^{ij}\partial_{x_i}u\overline{\partial_{x_j}u}\right)\frac{dMdy}{y^n} \\
&\geq& R^2\sum_{m=0}^{\infty}\lambda_m\int_{R}^{\infty}\|P_mu(y)\|^2_{L^2(M)}\frac{dy}{y^n} \\
&\geq& \lambda_1R^2\int_{M\times(R,\infty)}|(P'\otimes I_y)u|^2d\mathcal M.
\end{eqnarray*}
By the above inequality, for any $\epsilon > 0$ there exists $R > 1$ such that
$$
\sup_j\int_{M\times(R,\infty)}|(P'\otimes I_y)u_j|^2d\mathcal M < \epsilon.
$$
On $\mathcal M\setminus M\times(R,\infty)$ we apply Rellich's theorem to extract a convergent subsequence. This proves (1). 

Any $u \in D(L)$ is written as $u = L^{-1}f$ for some $f \in \mathcal H$. 
 It satisfies
\begin{equation}
l(u,v) = (Lu,v) = (f,v), \quad \forall v \in \mathcal H^1.
\label{weaksenseHNu}
\end{equation}
Taking $v$ from $C_0^{\infty}(\mathcal M^a_{int})$ and $(P'\otimes I_y)C_0^{\infty}(\mathcal M^a_{ext})$, we see that 
$$
(- \Delta_g  - \frac{(n-1)^2}{4} + 1 - \zeta)u = f \quad {\rm weakly} \quad {\rm in}\quad  \mathcal M^a_{int}, \quad {\rm  and} \quad \mathcal M^a_{ext}.
$$
 Therefore, $u \in H^2_{loc}(\mathcal M^a_{int})$,  $(P'\otimes I_y)u \in H^2_{loc}(\mathcal M^a_{ext})$.  
 Take $v = \varphi_m(x)\chi(y)$ $(m \geq1)$, where $\chi \in C_0^{\infty}((2,\infty))$  and $\varphi_m$ is the eigenfunction associated with $\lambda_m$. Then from (\ref{weaksenseHNu}), we see that $(u(\cdot,y),\varphi_m)$ satisfies a 2nd order differential equation on $(2,\infty)$.
Therefore, we have that $ (P'\otimes
  I_y)u \in H^2_{loc}(M\times(2,\infty))$. We then have $u \in H^2(\mathcal M^a_{int})$ and, by Theorem 2.1.3, $u = (P'\otimes I_y)u \in H^2(\mathcal M^a_{ext})$.
By taking $v \in (P_0\otimes I_y)C^{\infty}(M\times(2,a])$ such that $v = 0$ for $y < 3$ in (\ref{weaksenseHNu}), and integrating by parts, we have
$$
\left(\big(y^{(n-2}\partial_y(P_0\otimes I_y)u\big)(a-0),v\right)_{L^2(\Gamma^a)} = 0.
$$
 Therefore, 
$\left(\partial_y(P_0\otimes I_y)u\right)(a-0) = 0$. These facts prove $D(L) \subset D_L$.

Take $u \in D_L$ and put $h = (- \Delta_g -(n-1)^2/4 + 1)u$ for $y \neq a$. Then by integration by parts, we have
$$
l(u,v) = (h,v)_{\mathcal H}, \quad  \forall v \in \mathcal H^1.
$$
Since $l(u,v) = (L^{1/2}u,L^{1/2}v)_{\mathcal H}$, we then have
$$
\big|(L^{1/2}u,L^{1/2}v)_{\mathcal H}\big| \leq C\|v\|_{\mathcal H}, \quad 
\forall v \in \mathcal H^1
$$
with a constant $C$ independent of $v \in \mathcal H^1 = D(L^{1/2})$.
This shows that $L^{1/2}u \in D(L^{1/2})$, which proves $D_L \subset D(L)$. 
In  particular, we have proven for $y \neq a$
$$
Lu = \big(- \Delta_g  - \frac{(n-1)^2}{4} + 1\big)u, \quad u \in D(L).
$$

The uniqueness in (3) follows from $\zeta \not\in \sigma(L)$.  Indeed, if $u_1, u_2$ be two 
different solutions, then $u_i-u_2 \in D_L$ would be an eigenfunction of $L$.
To show the existence, we take  $\eta(y) \in C^{\infty}(\mathcal M^a_{int})$ such that 
$\eta(y) = 0$ for $y < 2$, $\eta(a-0) = 0$, $(\partial_y\eta)(y-a) = 1$, and $\eta(y) = 0$ in 
$\mathcal M^a_{ext}$. Let
\begin{equation}
\tilde f = \left\{
\begin{split}
& \big(- \Delta_g  - \frac{(n-1)^2}{4} + 1 -\zeta\big)\eta \quad {\rm in} \quad \mathcal M^a_{int},\\
& 0 \quad {\rm in} \mathcal\quad M^a_{ext},
\end{split}
\right.
\nonumber
\end{equation}
and put
\begin{equation}
u = u(z,\zeta,\lambda) = \lambda\chi(y) + (L - \zeta)^{-1}f - \lambda(L - \zeta)^{-1}\tilde f.
\label{C3S4existence}
\end{equation}
This is analytic with respect to $\lambda$ and meromorphic with respect to $\zeta$. \qed

\medskip
For $0 < \alpha < \beta < \infty$, we put
$$
U_{\alpha \beta}^{(\pm)} = \{\zeta \in {\bf C}\, ; \, 
\alpha < {\rm Re}\,\zeta < \beta, \ \ 
0 \leq \pm {\rm Im}\,\zeta\}.
$$

\begin{lemma} On $M\times(0,\infty)$, we consider $H_0 = - y^2(\partial_y^2 + \Delta_M) + (n-2)y\partial_y - (n-1)^2/4$, and
 $R_0(\zeta) = (H_0 - \zeta)^{-1}$. Suppose
 $f \in C_0^{\infty}({\mathcal M})$ satisfies ${\rm supp}\,f \subset \mathcal M_1 = M\times(1,\infty)$. Let $\rho(y) \in C^{\infty}((0,\infty))$ be such that $\rho(y) = 0$ for $y < 2$, $\rho(y) = 1$ for $y > 3$. Then, for any $0 < \alpha < \beta < \infty$, there exist $\epsilon > 0, \ C > 0$ such that
\begin{equation}
\Big|\rho(y)\left((P'\otimes I_y)R_{0}(\zeta)f\right)(x,y)\Big| \leq Ce^{-\epsilon y}, \quad 
 \zeta \in U_{
 \alpha \beta}^{(\pm)}.
\nonumber
\end{equation}
\end{lemma}
Proof. 
By (\ref{eq:Chap3Sect1ResolventofH0}),
\begin{equation}
u(x,y) := (P'\otimes I^a_y)R_0(\zeta + i0)f = \sum_{m\geq1}
\varphi_m(x)\left(G_0(\sqrt{\lambda_m},\nu)\widehat f_m\right)(y),
\nonumber
\end{equation}
with $\nu = - i\sqrt{\zeta}$, where $G_0(\zeta,\nu)$ is defined by  Definition 1.3.5. Then we have 
by Chap. 1, (\ref{C1S3Goexpdecay}) 
$$
\|u(\cdot,y)\|^2_{L^2(M)} = 
\sum_{m\geq1}|G_0(\sqrt{\lambda_m},\nu)\widehat f_m(y)|^2 \leq Ce^{-\epsilon y}.
$$
Note that supp$\,\widehat f_m(y)$ is away from 0, and the singularities of $I_{\nu}(y), K_{\nu}(y)$  at $y = 0$ do no harm. Since, for any $q > 0$, $\|\Delta_x^qu(\cdot,y)\|^2$ is estimated in a similar manner, by Sobolev's inequality we have $|u(x,y)|^2  \leq Ce^{-\epsilon y}$. \qed


\subsection{Meromorphic continuation of the Eisenstein series}
Here we pass to the traditional parametrization. 
For a subset $\mathcal E \subset {\bf R}$, we write
\begin{equation}
\frac{n-1}{2}  \pm \sqrt{- \mathcal E} = \left\{s \in {\bf C}\, ;\,
s(n-1-s) - \frac{(n-1)^2}{4} \in \mathcal E\right\}.
\nonumber
\end{equation}
Let $A = L-1-\frac{(n-1)^2}{4}$, and put
\begin{eqnarray*}
\begin{split}
\Sigma(A) =\frac{n-1}{2} \pm \sqrt{-\sigma(A)}, \quad 
\Sigma(H)  = \frac{n-1}{2} \pm \sqrt{- \sigma(H)}, \\
\Sigma_d(H)  = \frac{n-1}{2} \pm \sqrt{- \sigma_d(H)}, \quad
\Sigma_p(H)  = \frac{n-1}{2} \pm \sqrt{- \sigma_p(H)}, \\
\mathcal L = \Big\{s \in {\bf C}\, ; \, {\rm Re}\, s = \frac{n-1}{2}\Big\}, 
\quad
\mathcal L_{\pm} = \Big\{s \in \mathcal L\, ; \, \pm \, {\rm Im}\, s > 0\Big\}.  
\end{split}
\nonumber
\end{eqnarray*}
Note that $\Sigma(H) = \mathcal L\cup\Sigma_d(H)$, and that
$\Sigma(A)$ is a discrete set, since $\sigma(A)$ is discrete by Lemma 4.3.

In view of 
(\ref{eq:varphixk}), we  define for $\big\{{\rm Re}\,s > (n-1)/2\big\}\setminus\Sigma_p(H)$
\begin{equation}
\begin{split}
E(z,s) & = \chi(y)\,y^{s} - \big(- \Delta_g - s(n-1-s)\big)^{-1}\left[-\Delta_g,\chi(y)\right]y^s  \\
 & = \varphi(z, k - i\epsilon),
\end{split}
\nonumber
\end{equation}
where $s = (n-1)/2 + i(k - i\epsilon)$ $(\epsilon > 0)$. 
By Theorem 3.8, $E(z,s)$ is extended to $\mathcal L\setminus\big(\Sigma_p(H) \cup \{(n-1)/2\}\big)$. 
We take    $s = (n-1)/2 + ik \in \mathcal L\setminus\big(\Sigma_p(H) \cup \{(n-1)/2\}\big)$. 
Since $\left(-\Delta_g -s(n-1-s)  \right)^{-1} f$ satisfies outgoing radiation condition,
$$
E(z, s) -y^s \sim C y^{n-1-s}.
$$
Comparing with
(\ref{C3S4expansion2}),
\begin{equation}
E(z,s) \simeq y^s + \mathcal S(s)y^{n-1-s}, \quad {\rm as}\quad y \to \infty.
\nonumber
\end{equation}

By Lemma 4.3, for $s \not\in \Sigma(A)$, there exists a unique solution $v = v(z,s) \in D_L$ of the following boundary value problem
\begin{equation}
\left\{
\begin{split}
&\big(- \Delta_g - s(n-1-s)\big)v(z,s) = 0 \quad {\rm in} \quad \mathcal M^a_{int}, 
\\
&\big(- \Delta_g - s(n-1-s)\big)(P'\otimes I^1_y)v(z,s) = 0 \quad {\rm in} \quad \mathcal M^a_{ext},\\
&\big(y\partial_y(P_0\otimes  I^1_y)v\big)(a-0,s) = 1.
\end{split}
\right.
\label{C3S4BVPforv(s)}
\end{equation}
We define
\begin{equation}
\lambda_a(s) = \Big((P_0\otimes I^1_y)v\Big)(a -0,s).
\label{C3S4lambdaas}
\end{equation}
By Lemma 4.3 (3), $\lambda_a(s)$ is meromorphic on ${\bf C}$ with respect to $s$ with poles in $\Sigma(A)$.


\begin{lemma}
(1) For $s \in \mathcal L\setminus\big(\Sigma(A)\cup\Sigma_p(H)\cup\{(n-1)/2\}\big)$, we have
\begin{equation} \label{E4.10}
\lambda_a(s) = \frac{a^s + a^{n-1-s}\mathcal S(s)}{sa^s + (n-1-s) a^{n-1-s}\mathcal S(s)}, \quad
\mathcal S(s) = a^{2s-n+1}\frac{1 - s\lambda_a(s)}{(n-1-s)\lambda_a(s) - 1}.
\end{equation}
(2) Letting $v(z,s)$ be the solution to (\ref{C3S4BVPforv(s)}), we have 
\begin{equation}
E(z,s) -\left( s a(s)-(n-1-s) \mathcal S(s) a^{(n-1-s)} \right) v= 
\left\{
\begin{split}
& y^s + \mathcal S(s)y^{n-1-s}, \quad {\rm on} \quad \mathcal M^a_{ext},  \\
& 0, \quad {\rm on} \quad \mathcal M^a_{int}.
\end{split}
\right.
\nonumber
\end{equation}
(3) $\mathcal S(s)$ and $E(z,s)$ are extended to  meromorphic functions on ${\bf C}$.
\end{lemma}

Proof. 
Lemma 4.4 implies
$$
\big|(P'\otimes I^a_y)E(z,s)\big| \leq Ce^{-\epsilon y}, \quad \epsilon > 0.
$$
Hence, we have
\begin{equation}
(P_0\otimes I^a_y)E(z,s) \simeq  y^s + \mathcal S(s)y^{n-1-s}.
\nonumber
\end{equation}
On the other hand,for $y>3$, 
$$
\left( - y^2\partial_y^2 + (n-2)y\partial_y - s(n-1-s)\right)
(P_0\otimes I^3_y)E(z,s) = 0.
$$
Therefore, we have
\begin{equation}
(P_0\otimes I^3_y)E(z,s) =   y^s + \mathcal S(s)y^{n-1-s},
\label{C3S4P0Ey>3}
\end{equation}
since any solution of the equation $\left( - y^2\partial_y^2 + (n-2)y\partial_y - s(n-1-s)\right)u(y) = 0$
 is written uniquely by a linear combination of $y^s$ and $y^{n-1-s}$. Let
\begin{equation}
u = \left\{
\begin{split}
& E(z,s) \quad {\rm in} \quad \mathcal M^a_{int}, \\
& (P'\otimes I^a_y)E(z,s) \quad {\rm in} \quad \mathcal M^a_{ext}.
\end{split}
\right.
\nonumber
\end{equation}
Then $u \in D_L$, and
\begin{equation}
\left\{
\begin{split}
&\big(- \Delta_g - s(n-1-s)\big)u = 0 \quad {\rm in} \quad \mathcal M^a_{int}, 
\\
&\big(- \Delta_g - s(n-1-s)\big)(P'\otimes1)u = 0 \quad {\rm in} \quad \mathcal M^a_{ext},\\
&\big(y\partial_y(P_0\otimes I_y)u\big)(a-0,s) = sa^s + (n-1-s)\mathcal S(s)a^{n-1-s}.
\end{split}
\right.
\nonumber
\end{equation}
Comparing with (\ref{C3S4BVPforv(s)}), we obtain, by the uniqueness, 
\begin{equation}
u = \big(sa^s + (n-1-s)\mathcal S(s)a^{n-1-s}\big)v.
\label{C3S4u=Constv(z,s)}
\end{equation}
 Using (\ref{C3S4P0Ey>3}), we obtain (1). The assertions (2) and (3) are direct consequences of Lemma 4.3 (3), (\ref{C3S4u=Constv(z,s)}) and the meromorphy of $\lambda_a(s)$. \qed


\begin{lemma}
 $\lambda_a(s) \in {\bf R}$ for $s \in \mathcal L\setminus\Sigma(A)$
 {\ntext  and $\lambda_a(s)=\lambda_a(\bar s)$.} 
 \end{lemma}
Proof. Note that if $v \in D_L$, then $\overline{v} \in D_L$, and also that
$s(n-1-s) \in {\bf R}$ if $s \in \mathcal L$. Then, if $v(z,s)$ satisfies (\ref{C3S4BVPforv(s)}), so does $\overline{v(z,s)}$. By the uniqueness, $v(z,s)$ is then real-valued. This proves that $\lambda_a(s) \in {\bf R}$. 
 As, for $s \in \mathcal L$,  $s(n-1-s)=\bar s(n-1-\bar s)$ it follows from (\ref{C3S4BVPforv(s)}) that 
$\lambda_a(s)=\lambda_a(\bar s)$. 
\qed

 
\begin{theorem}
 $\mathcal S(s)$ is holomorphic on ${\rm Re}\,s = (n-1)/2$. 
\end{theorem}
Proof. Take $s_1 = (n-1)/2 + ik_1, \ 0 \neq k_1 \in {\bf R}$, and suppose $\lambda_a(s)$ is holomorphic at $s_1$. It follows from Lemma 4.6 that $\lambda_a(s_1)$ is real. Then $(n-1-s_1)\lambda_a(s_1) - 1 \neq 0$, hence by Lemma 4.5 (1), $\mathcal S(s)$ is holomorphic at $s_1$. 

Suppose $\lambda_a(s)$ has a pole at $s_1 = (n-1)/2 + ik_1, \ 0 \neq k_1 \in {\bf R}$. Then  $\kappa_a(s) = 1/\lambda_a(s)$ is holomorphic at $s_1$, and $\kappa_a(s_1) = 0$. By the formula
\begin{equation}
\mathcal S(s) = a^{2s-n+1}\frac{\kappa_a(s) - s}{n-1-s-\kappa_a(s)},
\label{eq:Sskappa}
\end{equation}
$\mathcal S(s)$ is holomorphic at $s_1$.

Suppose $\lambda_a(s)$ is holomorphic at $s_0 = (n-1)/2$. By Lemma 4.5 (1), if $\lambda_a(s_0) \neq 2/(n-1)$, $\mathcal S(s)$ is holomorphic at $s_0$, and $\mathcal S(s_0) = -1$. 
If $\lambda_a(s_0) = 2/(n-1)$, by the Taylor expansion
$\lambda_a(s_0 + w) = 2/(n-1) + cw + O(w^2)$. We then have
$$
\mathcal S(s_0 + w) = - a^{2w}\frac{\displaystyle{\Big(c + \Big(\frac{2}{n-1}\Big)^2\Big)w + O(w^2)}}{\displaystyle{\Big(c - \Big(\frac{2}{n-1}\Big)^2\Big)w + O(w^2)}}.
$$
Since $\lambda_a(s)=\lambda_a(\bar s)$, we have $c=0$. 
Therefore, $\mathcal S(s)$ is holomorphic at $s_0$ and $\mathcal S(s_0) = 1$.

Suppose $\lambda_a(s)$ has a pole at $s_0 = (n-1)/2$. By (\ref{eq:Sskappa}), $S(s)$ is holomorphic at $s_0$ and  $\mathcal S(s_0) = -1$. \qed

\medskip
Note, since by Theorem 3.16, $\widehat S(k)$ is unitary for $k > 0$, $k^2 \not\in \sigma_p(H)$, we have $|\mathcal S(s)| = 1$ a.e. on $\mathcal L$. In particular, due to the proof of Theorem 4.7, 
$\mathcal S((n-1)/2) = 
\pm 1$.

\medskip
To prove the holomorphy of $E(z,s)$, we prepare an identity. Let $v(z,s)$ be a solution to (\ref{C3S4BVPforv(s)}), and put
\begin{equation}
\widetilde w(z,s) = \left(sa^s + (n-1-s)a^{n-1-s}\mathcal S(s)\right)v(z,s),
\nonumber
\end{equation}
and, for $k \in {\bf R}$,
\begin{equation}
w(z,k) = \widetilde w\Big(z,\frac{n-1}{2} + ik\Big).
\nonumber
\end{equation}
It satisfies the equation
\begin{equation}
(L-1 - s(n-1-s))w = 0, \quad s = \frac{n-1}{2} + ik,
\nonumber
\end{equation}
and the boundary condition
\begin{equation}
\big((P_0\otimes I^1_y)\big)w(a-0,k) = a^{(n-1)/2+ik} + a^{(n-1)/2-ik}\mathcal S\Big(\frac{n-1}{2} + ik\Big),
\nonumber
\end{equation}
where we have used the definition of $\lambda_a(s)$ and Lemma 4.5. It also satisfies
\begin{equation}
\begin{split}
\big(y\partial_y(P_0\otimes I^1_y)w\big)(a-0,k)  = & \Big(\frac{n-1}{2} 
 + ik\Big)a^{(n-1)/2+ik} \\
 &  + 
\Big(\frac{n-1}{2} - ik\Big)a^{(n-1)/2-ik}\mathcal S\Big(\frac{n-1}{2} + ik\Big).
\end{split}
\nonumber
\end{equation}


\begin{lemma} For $k, h \in {\bf R}$, the following formula holds:
\begin{equation} 
\begin{split}
&(w(\cdot,k),w(\cdot,h))_{\mathcal H} \\
& =\frac{i}{k-h}\left(a^{i(h-k)}\mathcal S\Big(\frac{n-1}{2}+ik\Big)\overline{
\mathcal S\Big(\frac{n-1}{2}+ih\Big)} - a^{i(k-h)}\right) \\
& - \frac{i}{k + h}\left(a^{i(k+h)}\overline{\mathcal S\Big(\frac{n-1}{2} + ih\Big)}
- a^{-i(k+h)}\mathcal S\Big(\frac{n-1}{2} + ik\Big)\right).
\end{split}
\label{eq:wkwhinnerproduct}
\end{equation}
\end{lemma}
Proof. Letting $ w_0(y,k) = (P_0\otimes I^1_y)w\big|_{\mathcal M_1}$, we have, by integration 
by parts and Lemma 4.4,
\begin{equation}
\begin{split}
& (L w(k),w(h))_{\mathcal H} - (w(k),L w(h))_{\mathcal H} \\
&= \displaystyle{\frac{1}{y^{n-2}}\Big(w_0(y,k)(\overline{\partial_yw_0})(y,h)} - (\partial_y w_0)(y,k)\overline{w_0(y,h)}\Big)\Big|_{y=a-0}.
\end{split}
\nonumber
\end{equation}
Using the equation and the boundary conditions, we have
\begin{equation}
\begin{split}
&(k^2 - h^2)(w(k),w(h)) \\
& =i(h+k)\left(a^{i(h-k)}\mathcal S\Big(\frac{n-1}{2}+ik\Big)\overline{\mathcal S\Big(\frac{n-1}{2}+ih\Big)} - a^{i(k-h)}\right) \\
& + i(h - k)\left(a^{i(k+h)}\overline{\mathcal S\Big(\frac{n-1}{2} + ih\Big)}
- a^{-i(k+h)}\mathcal S\Big(\frac{n-1}{2} + ik\Big)\right),
\end{split}
\nonumber
\end{equation}
which proves the lemma. \qed


\begin{theorem}
Eisenstein series $E(z,s)$ is holomorphic on ${\rm Re}\,s = (n-1)/2$.
\end{theorem}
Proof. In view of Lemma 4.5 (2), we have only to show that when $k \to k_0 \in \Sigma(A)$, 
$\|w(k)\|$ is bounded. We prove this by first letting $h \to k \neq 0$ and $k \to k_0$ in (\ref{eq:wkwhinnerproduct}).
Since $\mathcal S(s)$ is holomorphic and, by the unitarity, $|\mathcal S(s)| = 1$ 
on ${\rm Re}\,s = (n -1)/2$, the 1st term of the right-hand side of (\ref{eq:wkwhinnerproduct}) is 
bounded in this process. The second term is bounded when $k_0 \neq 0$.

 By the note after Theorem 4.7, $\mathcal S(s_0) = \pm 1$ for $s_0 = (n-1)/2$. Therefore, the 2nd term of the right-hand side of (\ref{eq:wkwhinnerproduct}) is bounded when $k, h \to k_0$. \qed


\section{$SL(2,{\bf Z})\backslash{\bf H}^2$ as a Riemann surface}
 In this section we summarize the basic properties of the quotient manifold by the action of modular group 
\begin{equation}
SL(2,{\bf Z}) = \left\{
\left(
\begin{array}{cc}
a & b \\
c & d
\end{array}
\right)\, ; \, a, b, c, d \in {\bf Z},\, ad - bc = 1
\right\},
\nonumber
\end{equation}
where the action $SL(2,{\bf Z})\times{\bf C}_+ \ni (\gamma,z) \to \gamma\cdot z \in {\bf C}_+$ is defined by (\ref{eq:Chap3Action}). 
In the following, $I_2$ denotes the $2\times2$ unit matrix.


\subsection{Fundamental domain} 
Let ${\mathcal M} = SL(2,{\bf Z})\backslash {\bf H}^2$. The fundamental domain $\mathcal M^f$ of ${\mathcal M}$ is the following set: 
\begin{equation}
{\mathcal M^f} = \{
z \in {\bf C}_+ \, ; \, |z| \geq 1, - 1/2 \leq {\rm Re}\,z \leq 1/2 \},
\nonumber
\end{equation}
\begin{equation}
\partial M^f = \partial M_1^f\cup\partial M_2^f,
\nonumber
\end{equation}
\begin{equation}
\partial M_1^f = \Big\{-\frac{1}{2} + iy\, ;\, \frac{\sqrt3}{2} \leq y < \infty\Big\}\cup\Big\{\frac{1}{2} + iy\, ; \, \frac{\sqrt3}{2} \leq y < \infty\Big\},
\nonumber
\end{equation}
\begin{equation}
\partial M_2^f = \Big\{e^{i\varphi}\, ;\, \frac{\pi}{3} \leq \varphi \leq \frac{2\pi}{3}\Big\},
\nonumber
\end{equation}
 (\cite{Ume00} p. 241).  
We put
\begin{equation}
\gamma^{(T)} = 
\left(
\begin{array}{cc}
1 & 1 \\
0 & 1
\end{array}
\right), \quad 
\gamma^{(I)} = \left(
\begin{array}{cc}
0 & - 1 \\
1 & 0
\end{array}
\right).
\nonumber
\end{equation}
Their actions are
\begin{equation}
\gamma^{(T)}\cdot z = z + 1, \quad \gamma^{(I)}\cdot z = - \frac{1}{z}.
\nonumber
\end{equation}
To get $\mathcal M$ from $\mathcal M^f$, we glue $\partial M_1^f$ by the action of $\gamma^{(T)}$, i.e. $-\frac{1}{2}+iy \to \frac{1}{2} + iy$, and $\partial M_2^f$ by the  action of $\gamma^{(I)}$, i.e. $e^{i\varphi} \to e^{i(\pi - \varphi)}$. We denote this identification by $\Pi$, i.e.
$$
\mathcal M = \mathcal M^f/\Pi.
$$
The resulting surface $\mathcal M$ has two singular points, $p_1 = \Pi(i)$ and $p_2 = \Pi(e^{i\pi/3}) = \Pi(e^{2\pi i/3})$. The nature of these singularities is clarified by the following lemmas (see \cite{Ume00}, p. 247, p. 251).
We denote by $\left\langle \gamma \right\rangle$ the cyclic group generated by $\gamma$.


\begin{lemma}  
$SL(2,{\bf Z})$ is generated by $\gamma^{(T)}$ and $\gamma^{(I)}$. 
\end{lemma}


\begin{lemma}
For $w \in {\bf C}_+$, we put
\begin{equation}
G_w = \{\gamma \in SL(2,{\bf Z}) \, ; \, \gamma\cdot w = w\}.
\nonumber
\end{equation}
That $w \in {\mathcal M}^f$ and $G_w \neq \{\pm I_2\}$ occurs only for the following three cases. \\

\smallskip
\noindent
(1) $w = i$. In this case 
$\displaystyle{G_w = \left\langle \left(
\begin{array}{cc}
0 & -1 \\
1 & 0
\end{array}
 \right)\right\rangle}$.  \\
 
\smallskip
\noindent
(2) $w = e^{\pi i/3}$. In this case
$\displaystyle{G_w = \left\langle \left(
\begin{array}{cc}
0 & -1 \\
1 & -1
\end{array}
 \right)\right\rangle}$. \\
 
 \smallskip
\noindent
(3) $w = e^{2\pi i/3}$. In this case 
$\displaystyle{G_w = \left\langle \left(
\begin{array}{cc}
-1 & -1 \\
1 & 0
\end{array}
 \right)\right\rangle}$. 
\end{lemma}

Note that in the case $w = i$, the order of the group $G_w$ is 2, while in the case $w = e^{\pi i/3}$ and $e^{2\pi i/3}$ (which are identified by $\gamma^{(T)}$ and $\gamma^{(I)}$), the order of the group $G_{w}$ is 3. As a result, the point $p_i$ has a vicinity $U_i \subset \mathcal M$, $i = 1,2$, which can be represented as $U_1 = \Gamma_1\backslash B(1/2)$, $U_2 = \Gamma_2\backslash B(1/2)$, where $\Gamma_1$, $\Gamma_2$ are the groups of rotations corresponding to $G_i$ and $G_{e^{\pi i/3}}$, and $B(r)$ is the ball of radius $r > 0$ in ${\bf C}$ centered at $0$. These introduce orbifold structure on $\mathcal M$, however, in this note, we do not issue these constructions further.


\subsection{Analytic structure} 
To introduce local coordinates on $\mathcal M$, we consider 3 different cases.

\medskip
\noindent
1. Let $V_0 =  \mathcal M^f\setminus\partial \mathcal M_2^f$, and $U_0 = \Pi(V_0)$. Define for $p \in U_0$
$$
\zeta_0(p) = \varphi_0(z) = e^{2\pi i\, z}, \quad p = \Pi(z).
$$
Then, since two points $-1/2 + iy$, $1/2 + iy$ are identified by the action of $\gamma^{(T)}$, $\zeta_0(p)$ defines analytic coordinates on $U_0$.

\medskip
\noindent
2. Let $V_1 = \mathcal M^f\setminus\partial M_1^f$, and $U_1 = \Pi(V_1)$ be a neighborhood of $p_1 = \Pi(i)$.
Define for $p \in U_1$
\begin{equation}
\zeta_1(p) = \varphi_1(z) = \left(\frac{z - i}{z + i}\right)^2, \quad \Pi(z) = p.
\nonumber
\end{equation}
Then, since two points $e^{i\varphi}$, $e^{i(\pi - \varphi)}$, where $\pi/3 \leq \varphi < \pi/2$, are identified by the action of $\gamma^{(I)}$,
$\zeta_1(p)$ defines analytic coordinates on $U_1$.

\medskip
\noindent
3.  Let 			
$V_2 = \mathcal M^f\setminus i{\bf R}$, and $U_2 = \Pi(V_2)$ be a neighborhood of $p_2 = \Pi(e^{\pi i/3}) = \Pi(e^{2\pi i/3})$.
Define for $p \in U_2$
\begin{equation}
\zeta_2(p) = \varphi_2(z) = 
\left\{
\begin{split}
\left(\frac{z - e^{\pi i/3}}{z - e^{-\pi i/3}}\right)^3, \quad 
p = \Pi(z), \quad {\rm Re}\, z > 0, \\
\left(\frac{z - e^{2\pi i/3}}{z - e^{-2\pi i/3}}\right)^3, \quad 
p = \Pi(z), \quad {\rm Re}\, z < 0.
\end{split}
\right.
\nonumber
\end{equation}
Since two points $-1/2 + iy$, $1/2 + iy$ are identified by the action of $\gamma^{(T)}$, and two points $e^{i\varphi}$, $e^{i(\pi - \varphi)}$, where $\pi/3 \leq \varphi < \pi/2$, are identified by the action of $\gamma^{(I)}$,
 this $\zeta_2(p)$ defines analytic local coordinates on  $U_2$.
 
 To check that $\varphi_1,\,\varphi_2$ satisfy the desired analytical property, 
 it is convenient to observe that $\varphi_{1}, \varphi_2$ map the corresponding sectors of the
 circle $|z|=1$ onto an interval of a ray emanating from $0$.

\medskip
Since $\zeta_{\alpha}\circ \zeta_{\beta}^{-1}$ on $\zeta_{\beta}(U_{\alpha}\cap U_{\beta})$, $\alpha, \beta= 0,1,2$, are analytic,  the local coordinate system $\{(U_{\alpha},\zeta_{\alpha})\}_{{\alpha}=0}^2$ makes ${\mathcal M}$ a Riemann surface.


\subsection{Singularities as a Riemannian manifold} 
By the metric
\begin{equation}
ds^2 = \frac{(dx)^2 + (dy)^2}{y^2} = - \frac{4dzd\overline{z}}{(z - 
\overline{z})^2} \quad {\rm on} \quad \mathcal M^f,
\nonumber
\end{equation}
$\mathcal M$ becomes a hyperbolic space. However, we must pay attention to the points $p_1, p_2$. 
By the above local coordinate $\zeta_{\alpha}(p) = \varphi_{\alpha}(z)$, $p = \Pi(z)$,
 $\alpha = 0,1,2$, this metric becomes
\begin{equation}
ds^2 = \frac{d\zeta_{\alpha} d\overline{\zeta_{\alpha}}}{({\rm Im}\,z)^2|\varphi_{\alpha}'(z)|^2}.
\nonumber
\end{equation}
Therefore, on the zeros of $\varphi_{1}'(z)$, i.e. at $i$, $\varphi_2(z)$, i.e. $e^{\pi i/3}$, $e^{2\pi i/3}$, this Riemannian metric has singularities. In these cases,
$$
\zeta_{\alpha} = \varphi_{\alpha}(z) = T(z)^{\alpha +1}, \quad T(z) = \frac{z-w}{z-\overline w},
$$
where $w =i$  for $\alpha=1$, and  $w = e^{\pi i/3}$ and $w = e^{2\pi i/3}$ for $\alpha = 2$. In these cases, 
$$
z = \frac{w - \overline{w}\zeta^{1/n}}{1 - \zeta^{1/n}} = w + (w - \overline w)\zeta^{1/n} + \cdots.
$$
Therefore, $dz/d\zeta = n^{-1}(w-\overline w)\zeta^{1/n - 1} + \cdots$, hence
\begin{equation}
\left|\varphi_{\alpha}'(z)\right|^2 = \left|\frac{dz}{d\zeta}\right|^{-2}= O(|\zeta|^{\lambda}), \quad \lambda = 2 - \frac{2}{n}.
\label{C3S5|dzdzeta|2}
\end{equation}
Note that $1 \leq  \lambda < 2$.
The volume element and the Laplace-Beltrami operator are rewritten as
\begin{equation}
\frac{dx\wedge dy}{y^2} = \frac{i}{2y^2}dz\wedge d\overline{z} = \frac{i\big|dz/d\zeta\big|^2}{2({\rm Im}\,z)^2}d\zeta\wedge d\overline{\zeta},
\label{S2dxdy/y^2dz}
\end{equation}
\begin{equation}
y^2\big(\partial_x^2 + \partial_y^2\big) = 4({\rm Im}\,z)^2\frac{\partial^2}{\partial z\partial\overline{z}} = 
\frac{4({\rm Im}\,z)^2}{\big|dz/d\zeta|^2}\frac{\partial^2}{\partial\zeta\partial\overline{\zeta}}.
\label{S2delzdelzbar}
\end{equation}
Both of them have singularities at the corresponding $w$. However,  the singularity of the volume element and that of the Laplace-Beltrami operator cancel, since we have, for $C^{\infty}$-functions $f, g$ supported near $w$, 
\begin{equation}
\int_{\mathcal M}y^2\big(\partial_x^2 + \partial_y^2\big)f\cdot g\,\frac{dxdy}{y^2} = 2i
\int \frac{\partial^2}{\partial\zeta\partial\overline{\zeta}}f\cdot g\, d\zeta d\overline{\zeta}.
\label{S2dzetadzetabar}
\end{equation}

 We take small open neighborhoods $\widetilde U_i$ of $p_i$, $i = 1,2 $ such that $\widetilde U_1\cap \widetilde U_2 = \emptyset$. We construct a partition of unity $\{\chi_{\alpha}\}_{\alpha=0}^2$ such that ${\rm supp}\,\chi_{\alpha} \subset \widetilde U_{\alpha}$, $\alpha = 1,2$, ${\rm supp}\,\chi_0 \subset U_0$, and $\sum_{\alpha=0}^3\chi_{\alpha} = 1$ on $\mathcal M$. In addition to the hyperbolic volume element, let 
\begin{equation}
dV_E^{(\alpha)} = 
 \frac{i}{2}d\zeta_{\alpha}\wedge d\overline{\zeta_{\alpha}},
\label{dVE}
\end{equation}
and define a quadratic form $a(u,v)$ by
\begin{equation}
a(u,v) = \sum_{\alpha=0}^3\int\chi_{\alpha}u\,\overline{v}\ dV_H^{(\alpha)} + 
\sum_{\alpha=0}^3\int\chi_{\alpha}\nabla u\cdot\nabla\overline{v}\ dV_E^{(\alpha)},
\label{S2Quv}
\end{equation}
where
\begin{equation}
\nabla = (\partial_t,\partial_s),  \quad (\zeta = t + is).
\nonumber
\end{equation}
We can show that the  quadratic form $a(u,v)$ with domain $C_0^{\infty}(\mathcal M)$ 
is closable in $L^2(\mathcal M, dv_H).$ Let $\widetilde a(u,v)$ be its closed extension, and
 $\widetilde H^1$  the set of $u$ such that $a(u,u) < \infty$ equipped with the inner product (\ref{S2Quv}). This is the 1st order Sobolev space on $\mathcal M$. By Theorem 4.2, we have a self-adjoint operator $A$ 
such that $a(u,v) = (Au,v)_{\mathcal M,g}$ for $u \in D(A), \ v \in \widetilde H^1$. Then $1 - A$ is a self-adjoint realization of the Laplace-Beltrami operator $\Delta_g$.

When we deal with the perturbation problem of $\Delta_g$, we should restrict ourselves to the case that the coefficients of differential of more than one order of the pertubation term vanish around  $i, e^{\pi i/3}, e^{2\pi i/3}$. The precise assumption is as follows. 

\medskip
{\it Let 
$H_0 = - \Delta_g = - y^2(\partial_y^2 + \partial_x^2)$, and  $V$ a 2nd order differential operator on ${\mathcal M}$ such that

\medskip
\noindent
{\bf (M-1)} $H= H_0 + V$ is formally self-adjoint.

\medskip
\noindent
{\bf (M-2)} Around $i, e^{\pi i/3}, e^{2\pi i/3}$,  $V$ is an operator of multiplication by a bounded real function.

\medskip
\noindent
{\bf (M-3)} Except for the neighborhoods in (M-2),  $V$ is a differential operator of the form} :
\begin{equation}
V = \sum_{i + j \leq 2}a_{ij}(x,y)(y\partial_x)^i(y\partial_y)^j
\nonumber
\end{equation}
\begin{equation}
|D^{\alpha}a_{ij}(x,y)| \leq C_{\alpha}(1 + |\log y|)^{-{\rm min}(|\alpha|,1)-1-\epsilon}, \quad 
\forall \alpha,
\nonumber
\end{equation}
\begin{equation}
D = (D_x,D_y) = (y\partial_x,y\partial_y).
\nonumber
\end{equation}

We define a self-adjoint extension of $H$ through the quadratic form discussed in \S 4.
This means that we perturb the hyperbolic metric on $\mathcal M$ except for neighborhoods of singular points so that it is asymptotically equal to the original metric at infinity.

Since the measure $dxdy/y^2$ has singularties at $i, e^{\pi/3}, e^{2\pi i/3}$, the following lemma is not obvious.


\begin{lemma}
For any $R > 1$, let $\chi_R$ be the characteristic function of ${\mathcal M}\cap\{y < R\}$. Then $\chi_R(H + i)^{-1}$ is compact in $L^2({\mathcal M};dxdy/y^2)$.
\end{lemma}
Proof. Assume that $f_n, n = 1, 2, \cdots,$ are on the unit sphere of $L^2({\mathcal M};dxdy/y^2)$, and let $u_n = (H + i)^{-1}f_n$. By Rellich's theorem, from $\{\chi_Ru_n\}_{n\geq1}$ one can extract a subsequence which converges in $L^2$ outside small neighborhoods of singular points. 

 Around  $p_1 = i$ and $p_2 = \omega$, we take local coordinate $\zeta = t + is$ as above, and for a suffiently small $r > 0$, let $B_r$ be a disc $\{t^2 + s^2 < r^2\}$. Then, 
 if $u \in L^2(\mathcal M, \frac{dx dy}{y^2})$ has a support in $B_r$, we have by (\ref{S2dxdy/y^2dz})
\begin{equation}
\int_{B_r}|u|^2dtds \leq C\int_{B_r}|u|^2dV_H^{(\alpha)},
\label{BrdtdsBrdVH}
\end{equation}
with a constant $C > 0$. By the Sobolev imbedding $H^s({\bf R}^n) \subset L^p({\bf R}^n)$, where $0 \leq s < n/2$, $p = 2n/(n-2s)$, we have
\begin{equation}
H^1({\bf R}^2) \subset L^p({\bf R}^2), \quad \forall p > 2,
\label{S2Sobolev}
\end{equation}
with continuous inclusion. 

We take $\alpha, \beta$ such that $\alpha^{-1} + \beta^{-1} = 1$, $1 < \alpha < 2/\lambda$,
where $\lambda$ is defined by (\ref{C3S5|dzdzeta|2}). 
Then, by H{\"o}lder's inequality,
$$
\int_{B_{\delta}}|u|^2dV_H^{(\alpha)} \leq C\int_{B_{\delta}}r^{-\lambda}|u|^2dtds
\leq C\left(\int_{B_{\delta}}r^{-\lambda\alpha}dtds\right)^{1/\alpha}
\left(\int_{B_{\delta}}|u|^{2\beta}dtds\right)^{1/\beta}.
$$
Since $\lambda\alpha < 2$, the 1st term of the most right-hand side tends to 0 when $\delta \to 0$. To the 2nd term of the most right-hand side we apply (\ref{S2Sobolev}). Then, for any $\epsilon > 0$, there exists $\delta > 0$ such that
$$
\int_{B_{\delta}}|u|^2dV_H^{(\alpha)} \leq \epsilon\left(\int_{B_{\delta}}|u|^2dV_H^{(\alpha)} + \int_{B_{\delta}}|\nabla u|^2dV_E^{(\alpha)}\right), \quad u \in {\widetilde H}^1.
$$
Given the bouded sequence $\{u_n\}$ in $\widetilde H^1$, the integral of $|u_n|^2$ over $B_{\delta}$ can be made small uniformly in $n$. Outside $B_{\delta}$, we use the usual Rellich theorem. This proves the lemma. \qed


\subsection{Spectrum} 
By the above Lemma 5.3, the results in  \S 3 and \S 4 also hold for $H$. Let $R(z) = (H - z)^{-1}$.

\begin{theorem}
(1) $\ \sigma_e(H) = [0,\infty)$. \\
\noindent
(2) $\ \sigma_p(H)\cap(0,\infty)$ is of finite multiplicity, discrete as a subset in ${\bf R}$, with possible accumulation points 0 and $\infty$. \\
\noindent
(3) If $\lambda \in (0,\infty)\setminus\sigma_p(H)$,
$
R(\lambda \pm i0) \in {\bf B}({\mathcal B};{\mathcal B}^{\ast}).
$
\end{theorem}


\subsection{Eisenstein series} 
We return to the case of  $H_0 = -y^2(\partial_y^2 + \partial_x^2)$.
Let
\begin{equation}
{G} = SL(2,{\bf Z}),
\quad
{ G}_0 = \left\{
\left(
\begin{array}{cc}
1 & n \\
0 & 1 
\end{array}
\right) \, ; \, n \in {\bf Z}
\right\},
\nonumber
\end{equation}
i.e. $G_0$ is the group of translations by $n$ along the $y-$axis.


\begin{lemma}
(1) For $g = 
\left(
\begin{array}{cc}
a & b \\
c & d 
\end{array}
\right), g'=
\left(
\begin{array}{cc}
a' & b' \\
c' & d' 
\end{array}
\right) \in {G}$,
$$
g'g^{-1} \in {G}_0 \Longleftrightarrow \exists n \in {\bf Z} \quad {\rm s.t.} \quad a' - a = nc, \ b' - b = nd, \ 
c' = c,\ d' = d
$$
(2)
$\ \displaystyle{\left(
\begin{array}{cc}
1 & 0 \\
0 & 1 
\end{array}
\right), 
\left(
\begin{array}{cc}
\ast & \ast \\
c & d 
\end{array}
\right)}$, $(c,d) = 1$, are the complete representative of ${G}_0\backslash{G}$. Here $(c,d) = 1$ means that $c$ and $d$ are mutually prime.
\end{lemma}
The proof is omitted.

\bigskip
Let us note that for
$z = x + iy$
\begin{equation}
{\rm Im}\,g\cdot z = \frac{y}{(cx + d)^2 + c^2y^2}
\nonumber
\end{equation}
holds. The Eisenstein series is defined by
\begin{equation}
\widetilde E(z,s) = \sum_{[g] \in {\bf G}_0\backslash{\bf G}}\left({\rm Im}\,g\cdot 
z\right)^s = y^s + 
\sum_{(c,d) = 1}\left(\frac{y}{(cx + d)^2 + c^2y^2}\right)^s.
\label{E5.1}
\end{equation}
We show that it is absolutely convergent for $ {\rm Re}\,s > 1$.


\begin{lemma} For
$\ |x| \leq 1/2, \ y \geq \sqrt3/2$, $cd \neq 0$,
$$
 \frac{y}{(cx + d)^2 + c^2y^2}  \leq \frac{2}{\sqrt3 |cd|}.
$$
\end{lemma}
Proof. Letting $r^2 = x^2 + y^2$, we have
$$
(cx + d)^2 + c^2y^2 = r^2\left(c + \frac{dx}{r^2}\right)^2 + \frac{y^2}{r^2}d^2
\geq \frac{y^2}{r^2}d^2 \geq \frac{3}{4}d^2.
$$
This together with the obvious inequality
$$
 (cx + d)^2 + c^2y^2 \geq c^2y^2
$$
proves
$$
 (cx + d)^2 + c^2y^2 \geq \frac{1}{2}\left(c^2y^2 + 
 \frac{3}{4}d^2\right) \geq \frac{\sqrt3}{2}y|cd|.
 \qed
$$

Lemma 5.6 implies the following lemma. 

\begin{lemma}
 For ${\rm Re}\,s > 1$, the series (\ref{E5.1}) is absolutely convergent and
$$
  |\widetilde E(z,s) - y^s| \leq C_s, \quad \forall z \in {\mathcal M}.
$$
\end{lemma}

Since $y^s$ satisfies on ${\bf H}^2$,
$$ 
-\Delta (y^s)-s(1-s) y^s=0,
$$
due to $g \in SL(2,{\bf Z})$ being an isometry on ${\bf H}^2$,
$$
-\Delta \left(\hbox{Im}\, g\cdot z)^s  \right)-s(1-s) \left(\hbox{Im}\, g\cdot z \right)^s=0.
$$
In addition, $\left(\hbox{Im}\, g_0\cdot z \right)^s= \hbox{Im}\, z=y$ for $g_0 \in G_0$.
Therefore, by Lemma 5.5 (2), ${\widetilde E}(z, s)$ satisfies
$$
-\Delta {\widetilde E}(z, s) -s(1-s){\widetilde E}(z, s)=0, \quad \hbox{on}\,\, {\mathcal M}.
$$
By Lemma 5.7, ${\widetilde E}(z, s)-y^s \in L^{\infty}({\mathcal M}) \subset
L^2({\mathcal M})$, in view of
 ${\mathcal M}$ having  finite measure,  $L^{\infty}({\mathcal M}) \subset
L^2({\mathcal M})$. Therefore, for
${\rm Re}\, s > 1$
\begin{equation}
 \widetilde E(z,s) = \chi(y)y^s - R_0(s(1-s))\big([H_0,\chi]y^s\big).
 \nonumber
\end{equation}
Here $R_0(\zeta) = (H_0-\zeta)^{-1}$, and $\chi(y) \in C^{\infty}((0,\infty))$ such that $\chi(y) = 0$ for $y < 2$, $\chi(y) = 1$ for $y > 3$. 
This coincides with the Eisenstein series (\ref{eq:varphixk}) introduced in  \S 4.
By using properties of number theoretic functions and Poisson's summation formula, the S-matrix is computed as follows (see e.g. \cite{Iwa02}, p. 61).


\begin{theorem} For the case of $H_0 = -y^2(\partial_y^2 + \partial_x^2)$, we have
 $$
  \mathcal S(s) = \sqrt{\pi}\,\frac{\Gamma(s-1/2)\, \zeta(2s-1)}{\Gamma(s)\, \zeta(2s)},
 $$
where $\zeta(s)$ is Riemann's zeta function.
\end{theorem}

\begin{remark}
For 3-dimensions, one can define a similar surface by using the Picard group
\begin{equation}
SL(2,{\bf Z} + i{\bf Z}) = 
\left\{
\left(
\begin{array}{cc}
a & b \\
c & d
\end{array}
\right) ; a, b,c,d \in {\bf Z} + i{\bf Z}, \ ad - bc = 1
\right\},
\nonumber
\end{equation}
where the action is defined by quarternios. 
The quotient space $SL(2,{\bf Z} + i{\bf Z})\backslash{\bf H}^3$ is also an orbifold. 
See \cite{EGM98}.
\end{remark}


\chapter{Radon transform and propagation of singularities in ${\bf H}^n$}

The purpose of this chapter is to extend Theorem 1.6.6 to the asymptotically hyperbolic metric on ${\bf R}^n_+$ in the sense of singularity expansion. 


\section{Geodesic coordinates near infinity}


\subsection{Geodesic coordinates}

We shall study  the metric
\begin{equation}
ds^2 = y^{-2}\Big((dx)^2 + (dy)^2 + A(x,y,dx,dy)\Big)
\label{eq:Chap4Sect1metric}
\end{equation}
on ${\bf R}^n_+$ defined in Chapter 2, Subsection 2.1, i.e. the metric satisfying the condition (C) in Chap. 2. Our aim is to transform (\ref{eq:Chap4Sect1metric}) into the following canonical form
\begin{equation}
ds^2 = y^{-2}\Big((dx)^2 + (dy)^2 + B(x,y,dx)\Big)
\label{eq:ds2stnadardform}
\end{equation}
in the region $0 < y < y_0$, $y_0$ being a sufficiently small constant, where $B(x,y,dx)$ is a symmetric covariant tensor of the form
\begin{equation}
B(x,y,dx) = \sum_{i,j=1}^{n-1}b_{ij}(x,y)dx^idx^j.
\nonumber
\end{equation}

Passing to the variable $z = \log y$, we rewrite the Laplace-Beltrami operator $\Delta_g$ associated with (\ref{eq:Chap4Sect1metric}) as
\begin{equation}
\begin{split}
\Delta_g & = \partial_z^2 + e^{2z}\partial_x^2 + \sum_{i,j=1}^{n-1}a^{ij}(x,e^z)e^{2z}\partial_{x_i}\partial_{x_j} \\
& \hskip 23mm + 2\sum_{i=1}^{n-1}a^{in}(x,e^z)e^{z}\partial_{x_i}\partial_z + 
a^{nn}(x,e^z)\partial_z^2
\end{split}
\nonumber
\end{equation}
up to 1st order terms. Then $\left(g^{ij}\right)$ in the variables $x$ and $z$ takes the form
\begin{equation}
g^{ij} = \left\{
\begin{split}
& e^{2z}\big(\delta^{ij} + h^{ij}(x,z)\big), \quad 1 \leq i, j \leq n-1, \\
&e^zh^{in}(x,z), \quad 1 \leq i \leq n-1, \\
& 1 + h^{nn}(x,z), \quad i, j = n,
\end{split}
\right.
\label{C4S1gijform}
\end{equation}
where $h^{ij}(x,z)$ satisfies in the region $z < 0$
\begin{equation}
|\partial_{x}^{\alpha}\partial_z^{\beta} h^{ij}(x,z)| \leq C_{\alpha\beta}
W(x,z)^{-{\rm min}(|\alpha|+\beta,1) - 1 - \epsilon_0},
\label{eq:xzshortrangecond}
\end{equation}
and
\begin{equation}
W(x,z) = 1 + |z| + \log\big(|x| + 1\big).
\nonumber
\end{equation}
We define the Hamiltonian $H(x,z,\xi,\eta)$ by
\begin{equation}
H(x,z,\xi,\eta) = \frac{1}{2}\Big(e^{2z}|\xi|^2 + \eta^2 + h(x,z,\xi,\eta)\Big),\nonumber
\end{equation}
\begin{equation}
h(x,z,\xi,\eta) = \sum_{i,j=1}^{n-1}e^{2z}h^{ij}(x,z)\xi_i\xi_j + 2\sum_{i=1}^{n-1}e^{z}h^{in}(x,z)\xi_i\eta + h^{nn}(x,z)\eta^2.
\nonumber
\end{equation}
The equation of geodesic  is as follows:
\begin{equation}
\left\{
\begin{split}
& \frac{dx}{dt} = \frac{\partial H}{\partial \xi}, \quad  \frac{dz}{dt} = \frac{\partial H}{\partial \eta}, \\
&\frac{d\xi}{dt} = - \frac{\partial H}{\partial x}, \quad \frac{d\eta}{dt} = - \frac{\partial H}{\partial z}.
\end{split}
\right.
\label{eq:HamiltonEquation}
\end{equation}
If $h(x,z,\xi,\eta) = 0$, it has the following solution
$$
x(t) = x_0, \quad \xi(t) = 0, \quad z(t) =  t, \quad \eta(t) = 1.
$$
With this in mind, we seek the solution of the equation (\ref{eq:HamiltonEquation}) which behaves like
\begin{equation}
\left\{
\begin{split}
& x(t) = x_0+ O(W(x_0,t)^{-1-\epsilon}), \quad 
\xi(t) = O(W(x_0,t)^{-1-\epsilon}), \\
& z(t) =  t + O(W(x_0,t)^{-\epsilon}), \quad
\eta(t) = 1 + O(W(x_0,t)^{-1-\epsilon}),
\end{split}
\right.
\nonumber
\end{equation}
as $t \to - \infty$, where $x_0 \in {\bf R}^{n-1}$, $0 < \epsilon < \epsilon_0$. Therefore we put
\begin{equation}
\left\{
\begin{split}
& U_x(x_0,t) = x(t) - x_0, \quad U_z(x_0,t) = z(t) - t, \\
& U_{\xi}(x_0,t) = \xi(t), \quad U_{\eta}(x_0,t,) = \eta(t) - 1,
\end{split}
\right.
\nonumber
\end{equation}
\begin{equation}
U(x_0,t) = \big(U_x(x_0,t), U_z(x_0,t), U_{\xi}(x_0,t), U_{\eta}(x_0,t)\big),
\nonumber
\end{equation}
\begin{equation}
A(U,x_0,t) = \Big(\frac{\partial H}{\partial \xi},\frac{\partial H}{\partial \eta} - 1, - \frac{\partial H}{\partial x}, -\frac{\partial H}{\partial z}\Big)\Big|_{x = U_x+x_0,\xi=U_{\xi},z=U_z+t,\eta=U_{\eta}+1},
\nonumber
\end{equation}
and consider the following non-linear operator
\begin{equation}
\big(B(U(x_0,\cdot);x_0)\big)(t) =  \int_{-\infty}^{t}A(U(x_0,\tau),x_0,\tau)d\tau.
\label{C4S1Bint}
\end{equation}
We shall look for the fixed point of the map : $U \to B(U)$, i.e.
\begin{equation}
U(x_0,t) = \big(B(U(x_0,\cdot);x_0)\big)(t).
\label{eq:UtequalBU}
\end{equation}
We fix $t_0 < 0$, and define the norm 
\begin{equation}
\begin{split}
\|U\|_{t_0} & = \sup_{t<t_0,x_0\in{\bf R}^{n-1}}
\big[|t| + \log(|x_0| + 1)\big]^{\epsilon/2}
|U_z(t)| \\
& \ \ + \sup_{t<t_0,x_0\in{\bf R}^{n-1}}
\big[|t|+ \log(|x_0| + 1)\big]^{1+\epsilon}
\big(|U_{\xi}(t)| + |U_{\eta}(t)| + |U_{x}(t)|\big),
\end{split}
\nonumber
\end{equation}
and the space $\mathcal F_{t_0}$ of functions by
\begin{equation}
\mathcal F_{t_0} \ni U(t) \Longleftrightarrow \|U\|_{t_0} < 1.
\nonumber
\end{equation}
By (\ref{eq:xzshortrangecond}), a simple computation shows 
\begin{equation}
\left|\frac{\partial H}{\partial z}\right| \leq C\|U\|_{t_0}\left(W(x_0,t)^{-2-\epsilon_0} + e^{t}W(x_0,t)^{-1-\epsilon}\right).
\nonumber
\end{equation}
Hence for any $\delta > 0$, there exists $t_0$ such that for $t < t_0$
\begin{equation}
\left|B(U(\cdot),x_0)_{\eta}(t)\right| \leq \int_{-\infty}^t\left|\frac{\partial H}{\partial z}\right|d\tau \leq \delta\|U\|_{t_0}W(x_0,t)^{-1-\epsilon}.
\nonumber
\end{equation}
Using this estimate and  (\ref{eq:xzshortrangecond}), we obtain, taking bigger $|t_0|$ if necessary,
\begin{equation}
\|B(U)(t)\|_{t_0} \leq \delta\|U\|_{t_0}, \quad \forall U \in {\mathcal F}_{t_0}.  
\nonumber
\end{equation}
Similar calculation implies
\begin{equation}
\|B(U)(t) - B(V)(t)\|_{t_0} \leq \delta\|U - V\|_{t_0}, 
\nonumber
\end{equation}
for $U, V \in \mathcal F_{t_0}$.
Then taking $\delta < 1/2$, $B$ maps $\mathcal F_{t_0}$ into $\mathcal F_{t_0}$, and is Lipschitz continuous with Lipschitz constant $ < 1/2$.
Hence, there exists a unique fixed point $U(t) = U(x_0,t) \in {\mathcal F}_{t_0}$ of (\ref{eq:UtequalBU}). 
By differentiating (\ref{C4S1Bint}) with respect to $t$, we see that for some constant $C$
\begin{equation}
\frac{1}{C}W(x_0,t)\partial_tU(x_0,t) \in {\mathcal F}_{t_0}. 
\nonumber
\end{equation}
Differentiating (\ref{eq:UtequalBU}) with respect to $x_0$, we get
$$
(I - B_U(U(x_0,\cdot),x_0))\partial_{x_0}^{\alpha}U = \partial_{x_0}^{\alpha}
B(U,x_0), \quad |\alpha|=1.
$$ 
For $t < |t_0|$, $(I - B_U(U(x_0,\cdot),x_0))$ is invertible, providing
\begin{equation}
\frac{1}{C}W(x_0,t)\partial_{x_0}^{\alpha}U(x_0,t) \in {\mathcal F}_{t_0}, \quad |\alpha|=1. 
\nonumber
\end{equation}
Iterating this procedure, we have the following lemma.


\begin{lemma} Choose $|t_0|$ large enough. Then
 there exists a solution $x(t)$, $z(t)$, $\xi(t)$, $\eta(t)$ of the equation (\ref{eq:HamiltonEquation}) for $(x_0,t) \in {\bf R}^{n-1}\times(- \infty,t_0)$ satisfying
\begin{equation}
\begin{split}
& \big|\partial_{x_0}^{\alpha}\partial_t^{\beta}\big(x(t) - x_0\big)\big| +  
\big|\partial_{x_0}^{\alpha}\partial_t^{\beta}\xi(t)| + 
\big|\partial_{x_0}^{\alpha}\partial_t^{\beta}\big(\eta(t) - 1\big)\big| \\
& \ \ \ \ \ \ \ \ \ \ \ \ \ \ \ \ \ \ \ \ \ \leq C_{\alpha\beta}W(x_0,t)^{-1-\epsilon/2 -{\rm min}(|\alpha|+ \beta,1)},\\
&\big|\partial_{x_0}^{\alpha}\partial_t^{\beta}\big(z(t) - t\big)| 
\leq C_{\alpha\beta}W(x_0,t)^{-\epsilon/2 -{\rm min}(|\alpha|+ \beta,1)}.
\end{split}
\nonumber
\end{equation}
\end{lemma}


\begin{lemma}
As a 2-form on the region ${\bf R}^{n-1}\times(-\infty,t_0)$, we have
\begin{equation}
\sum_{i=1}^{n-1}d\xi_i(x_0,t)\wedge dx^i(x_0,t) + d\eta(x_0,t)\wedge dz(x_0,t) = 0.
\nonumber
\end{equation}
\end{lemma}
Proof. We put $x^n = z$, $\xi_n = \eta$ and $x_0^n = t$. Then we have
\begin{equation}
\sum_{i=1}^nd\xi_i\wedge dx^i = \sum_{j<k}[\xi,x]_{jk}dx_0^j\wedge dx_0^k,
\nonumber
\end{equation}
\begin{equation}
[\xi,x]_{jk} = \frac{\partial \xi}{\partial x_0^j}\cdot
\frac{\partial x}{\partial x_0^k} - \frac{\partial \xi}{\partial x_0^k}\cdot\frac{\partial x}{\partial x_0^j}.
\nonumber
\end{equation}
Noting that
\begin{equation}
\frac{\partial}{\partial t}\left(\frac{\partial\xi}{\partial x_0^j}\cdot\frac{\partial x}{\partial x_0^k}\right) = - \frac{\partial^2H}{\partial x^i\partial x^m}\frac{\partial x^m}{\partial x_0^j}\frac{\partial x^i}{\partial x_0^k} + \frac{\partial^2H}{\partial \xi_i\partial \xi_m}\frac{\partial \xi_i}{\partial x_0^k}\frac{\partial \xi_m}{\partial x_0^j}
\nonumber
\end{equation}
is symmetric with respect to $j$ and $k$, we have
\begin{equation}
\frac{\partial}{\partial t}[\xi,x]_{jk} = 0.
\nonumber
\end{equation}
By Lemma 1.1, $[\xi,x]_{jk} \to 0$ as $t \to - \infty$. Hence $[\xi,x]_{jk} = 0$, which proves the lemma. \qed


\begin{lemma}
For large $|t_0|$, the map 
$$
{\bf R}^{n-1}\times(-\infty,t_0) \ni (x_0,t) \to  (x(x_0,t),z(x_0,t))
$$
 is a diffeomorphism and its image includes
${\bf R}^{n-1}\times(- \infty, 2t_0)$.
\end{lemma}

Proof. We show that this map is locally diffeomorphic and globally injective. Using inverse function theorem, from Lemma 1.1, we have that making $|t_0|$ sufficiently large, there are $r_0, \tilde r_0 > 0$ with the following properties; 

\begin{itemize}
\item
For any $x_0' \in {\bf R}^{n-1}, t_0' < t_0$, the map $(x(x_0,t),z(x_0,t))$ is a diffeomorphism from $B_r(x_0',t_0')$, the ball of radius $r$ with center at $(x_0',t_0')$, onto $U \subset {\bf R}^{n-1}\times (-\infty,t_0)$.
\item 
$B_{\tilde r_0}(x(x_0',t_0'),z(x_0',t_0')) \subset U$.
\end{itemize}
Assume  $x(x_0',t_0')=x(x_0'',t_0'')$, $z(x_0',t_0') = z(x_0'',t_0'')$ for some $(x_0',t_0') \neq (x_0'',t_0'')$. Then by Lemma 1.1, it follows from the 2nd equality that $|t_0' - t_0''| < r/4$ if $|t_0|$ is sufficiently large. Therefore by local injectivity, $|x_0'-x_0''| > 3r/4$. Using again Lemma 1.1, we see that for sufficiently large $|t_0|$, $|x(x_0',t_0')-x_0'| < r/4$, $|x(x_0'',t_0'') - x_0''| < r/4$. This leads to a contradiction. \qed

\medskip

Let $x_0 = x_0(x,z)$, $t = t(x,z)$ be the inverse of the map : $(x_0,t) \to (x,z)$. We put $\xi(x,z) = \xi(x_0(x,z),t(x,z))$, etc.  for the sake of simplicity. Since $\sum_{i=1}^{n-1}\xi_i dx^i + \eta dz$ is a closed 1-form by Lemma 1.2, we have
\begin{equation}
\frac{\partial\xi_j}{\partial x^k} = \frac{\partial\xi_k}{\partial x^j}, \quad
\frac{\partial\xi_j}{\partial z} = \frac{\partial\eta}{\partial x^j}, \quad
1 \leq j, k \leq n-1.
\nonumber
\end{equation}
Recall
\begin{equation}
\begin{split}
& U_{\eta}(x,z) = \eta(x,z) - 1 \\
& = - \int_{-\infty}^{t}\frac{\partial H}{\partial z}
\left(x(x_0,s),z(x_0,s),\xi(x_0,s),\eta(x_0,s)\right)ds\Big|_{x_0=x_0(x,z),t=t(x,z)},
\end{split}
\nonumber
\end{equation}
and define $\Psi(x,z)$ by
\begin{equation}
\Psi(x,z) = z + \int_{-\infty}^{0}U_{\eta}(x,z + \tau)d\tau.
\nonumber
\end{equation}


\begin{lemma} For $z \leq 2t_0$, we have \\
\noindent
(1) $\ \partial_x\Psi(x,z) = \xi(x,z)$, \\
\noindent
(2) $\ \partial_z\Psi(x,z) = \eta(x,z)$, \\
\noindent
(3) $\ H(x,z,\partial_x\Psi(x,z),\partial_z\Psi_z(x,z)) = 1/2$, \\
\noindent
(4) $\ \big|\partial_{x}^{\alpha}\partial_z^{\beta}(\Psi(x,z) - z)\big| \leq 
C_{\alpha\beta}(|z| + \log(|x| + 1))^{-\epsilon/2-{\rm min}(|\alpha|+ \beta,1)}, \quad \forall \alpha, \beta$. \\
\noindent
(5) $\ \Psi(x,z) = t(x,z)$.
\end{lemma}
Proof. We have
\begin{equation}
\begin{split}
\frac{\partial\Psi}{\partial x^j} &=  \int_{-\infty}^{0}\frac{\partial\eta}{\partial x^j}(x,z+\tau)d\tau \\
& =  \int_{-\infty}^{0}\frac{\partial\xi_j}{\partial \tau}(x,z+\tau)d\tau = \xi_j(x,z), 
\end{split}
\nonumber
\end{equation}
\begin{equation}
\begin{split}
\frac{\partial\Psi}{\partial z} &= 1 + \int_{-\infty}^{0}\frac{\partial\eta}{\partial \tau}(x,z+\tau)dt\tau = \eta(x,z), 
\end{split}
\nonumber
\end{equation}
which prove (1) and (2). 

Since $x(t), z(t)$ and $\xi(t), \eta(t)$ are solutions to the equation (\ref{eq:HamiltonEquation}), $H(x(t),p(t),\xi(t),\eta(t))$ is a constant, which turns out to be $1/2$ by letting $t \to - \infty$. This proves (3). (4) follows again from Lemma 1.1 due to the fact that
$$
\left|\partial_{x_0}^{\gamma}\partial_t^{\delta}
\left(\frac{\partial(x,z)}{\partial(x_0,t)}- Id\right)\right| 
\leq C_{\gamma\delta}W(x,z)^{-\epsilon/2-{\rm min}(|\gamma|+\delta,1)}.
$$
 Using (1), (2), we have
\begin{equation}
\begin{split}
\frac{\partial\Psi}{\partial t} & = \frac{\partial\Psi}{\partial x}\cdot\frac{\partial x}{\partial t} + \frac{\partial\Psi}{\partial z}\frac{\partial z}{\partial t} \\
& = \xi(x,z)\cdot\frac{\partial x}{\partial t} + \eta(x,z)\frac{\partial z}{\partial t} \\
& = \xi(x,z)\cdot\frac{\partial H}{\partial \xi} + \eta(x,z)\frac{\partial H}{\partial \eta} \\
&= g^{ij}\partial_i\Psi\partial_j\Psi = 1,
\end{split}
\nonumber
\end{equation}
where the last identity comes from Lemma 1.4 (3). Here
 $\partial_i = \partial/\partial x^i, 1 \leq i \leq n-1$, $\partial_n = \partial/\partial z$.
Therefore $\Psi(x,z) - t$ is independent of $t$. On the other hand, $\Psi - z \to 0$ and $z - t\to 0$ as $t \to -\infty$. Therefore, $\Psi(x,z) =t$. \qed


\begin{lemma}
In the coordinate system $(x_0,t)$, the Riemannian metric (\ref{eq:Chap4Sect1metric}) is written as
\begin{equation}
ds^2 = (dt)^2 + e^{-2t}\Big((dx_0)^2 + \sum_{i,j=1}^{n-1}{\widehat h}_{ij}(x_0,t)dx_0^idx_0^j\Big),
\nonumber
\end{equation}
where $\widehat h_{ij}(x_0,t)$ satisfies
\begin{equation}
\big|\partial_{x_0}^{\alpha}\partial_t^{\beta}{\widehat h}_{ij}(x_0,t)\big| \leq C_{\alpha\beta}W(x_0,t)^{-1-\epsilon/2 - {\rm min}(|\alpha|+ \beta,1)}, \quad \forall \alpha, \beta.
\label{C4S1hij(0)inW}
\end{equation}
\end{lemma}
Proof. We put $y^i = x_0^i, 1 \leq i \leq n-1$, $y^n = t$. Then the associated tensor $\overline{g}^{ij}$ is written as
\begin{equation}
\begin{split}
& \overline{g}^{nn} = g^{ij}\frac{\partial y^n}{\partial x^i}\frac{\partial y^n}{\partial x^j} = g^{ij}(\partial_i\Psi)(\partial_j\Psi) = 1, \\
& \overline{g}^{nk} = g^{ij}\frac{\partial y^n}{\partial x^i}\frac{\partial y^k}{\partial x^j} = g^{ij}(\partial_i\Psi)(\partial_jx_0^k) = 0,
\end{split}
\nonumber
\end{equation}
for $1 \leq k \leq n-1$. Here in the 2nd line, we have used
\begin{equation}
0 = \frac{\partial x_0^k}{\partial t} = \frac{\partial x_0^k}{\partial x^i}
\frac{\partial x^i}{\partial t} = \frac{\partial x_0^k}{\partial x^i}g^{ij}
\partial_j\Psi.
\nonumber
\end{equation}
Therefore the Riemmanian metric has the form
\begin{equation}
ds^2 = (dt)^2 + \sum_{i,j=1}^{n-1}\overline{g}_{ij}dx_0^idx_0^j.\nonumber
\end{equation}
Recall
$$
\overline{g}_{ij}(x_0,t) = g_{kl}\frac{\partial x^k}{\partial x_0^i}\frac{\partial x^l}{\partial x_0^j}
+ 2g_{kn}\frac{\partial x^k}{\partial x_0^i}\frac{\partial z}{\partial x_0^j} 
+ g_{nn}\frac{\partial z}{\partial x_0^i}\frac{\partial z}{\partial x_0^j},
$$
where $1 \leq k, l \leq n-1$, and the right-hand side is evaluated at 
$(x,z) = (x(x_0,t), z(x_0,t))$.
By the formula (\ref{C4S1gijform}), (\ref{eq:xzshortrangecond}) and Lemma 1.1,
the 1st term of the right-hand side is of the form 
$
e^{-2t}\left(\delta_{ij} + {\widehat h}^{(0)}_{ij}\right),
$
where $\widehat{h}^{(0)}_{ij}$ satisfies the estimate (\ref{C4S1hij(0)inW}). By the same reasoning, the 2nd and 3rd terms give rise to $\widehat h^{(1)}_{ij}$ and $\widehat h^{(2)}_{ij}$.
This completes the proof of the lemma. \qed

\bigskip
The coordinates $(x_0,t)$ are actually semi-geodesic coordinates related to 
the boundary at infinity $y = 0$.

\bigskip
Letting $x_0 = \overline{x}$, $t = \log\overline{y}$ in Lemma 1.5 and recalling that $D_{\overline y} = {\overline y}\partial_{\overline y} = \partial_t$, and using Lemma 1.1, we obtain the following theorem.


\begin{theorem}
Choose $y_0 > 0$ sufficiently small. Then there exists a diffeomorphism $(x,y) \to (\overline{x},\overline{y})$ in the region $0 < y < y_0$ such that 
\begin{equation}
\big|\partial_{\overline x}^{\alpha}D_{\overline y}^{\beta}\big(\overline{x} - x\big)| \leq C_{\alpha\beta}(1 + d_h({\overline x},{\overline y}))^{-{\rm min}(|\alpha|+ \beta,1)-1-\epsilon/2}, \quad 
\forall \alpha, \beta,
\nonumber
\end{equation}
\begin{equation}
\big|\partial_{\overline x}^{\alpha}D_{\overline y}^{\beta}\Big(\frac{\overline{y} - y}{\overline y}\Big)| \leq C_{\alpha\beta}(1 + d_h({\overline x},{\overline y}))^{-{\rm min}(|\alpha|+ \beta,1)-1-\epsilon/2}, \quad 
\forall \alpha, \beta,
\nonumber
\end{equation}
and in the $(\overline{x}, \overline{y})$ coordinate system, the Riemannian metric takes the form
\begin{equation}
ds^2 = (\overline y)^{-2}\Big((d\overline{y})^2 + (d\overline{x})^2 + \sum_{i,j=1}^{n-1}{\overline h}_{ij}(\overline{x},\overline{y})d\overline{x}^id\overline{x}^j\Big),
\nonumber
\end{equation}
where 
\begin{equation}
{\overline h}_{ij}(\overline x,\overline y) = {\widehat h}_{ij}(x_0,t), \quad 
x_0 = \overline x, \quad t = \log\overline{y},
\nonumber
\end{equation}
\begin{equation}
\big|\partial_{\overline{x}}^{\alpha}\,D_{\overline{y}}^{\beta}\, {\overline h}_{ij}(\overline{x},\overline{y})\big| \leq C_{\alpha\beta}(1+ d_h(\overline{x},\overline{y}))^{-{\rm min}(|\alpha|+ \beta,1)-1-\epsilon/2}, \quad \forall \alpha, \beta.
\nonumber
\end{equation}
\end{theorem}


\section{Asymptotic solutions to the wave equation}
Theorem 1.6 leads us to consider the metric having the form
\begin{equation}
ds^2 = y^{-2}\Big((dy)^2 + (dx)^2 + \sum_{i,j=1}^{n-1}h_{ij}(x,y)dx^idx^j\Big),
\label{eq:Chap4Sect3Riemmanianmetric}
\end{equation}
in the region ${\bf R}^{n-1}\times(0,y_0)$, 
where $y_0$ is a small constant and $h_{ij}(x,y)$ satisfies 
$$
h_{ij} \in {\mathcal W}^{-1-\epsilon/2}.
$$
 As in Chap. 2,  we consider
\begin{equation}
H = - (y^{2n}g)^{1/4}\Delta_g(y^{2n}g)^{-1/4} - \frac{(n-1)^2}{4} \quad {\rm in}\quad 
L^2\Big({\bf R}^n_+;\frac{dxdy}{y^n}\Big).
\nonumber
\end{equation}
Taking into account that $H$ is self-adjoint, we see that
explicitly, $H$ has the form
\begin{equation}
H = - D_y^2 + (n - 1)D_y  - D_x^2 - \frac{(n-1)^2}{4} - L,
\nonumber
\end{equation}
\begin{equation}
 L = y^2\sum_{|\alpha|\leq2}L_{\alpha}(x,y)\partial_x^{\alpha},
\label{eq:Chap4Sec3FormofL}
\end{equation}
where $D_y = y\partial_y$, $D_x = y\partial_x$. Moreover 
$L_{\alpha} \in {\mathcal W}^{-1-\epsilon/2}$.

It is convenient to rewrite $H$ into the form
\begin{equation}
H = - \Big(D_y - \frac{n-1}{2}\Big)^2  - K,\label{eq:Chap4Sec3rewrintingH}
\end{equation}
\begin{equation}
K = y^2(\partial_x)^2 + 
y^2\sum_{|\alpha|\leq2}L_{\alpha}(x,y)\partial_x^{\alpha}.
\label{eq:Chap4Sect3FormofK}
\end{equation}
Using 
\begin{equation}
\Big(D_y - \frac{n-1}{2}\Big)^m\left(e^{ix\cdot\xi}y^{\frac{n-1}{2}-ik}a\right) = e^{ix\cdot\xi}y^{\frac{n-1}{2}-ik}
(D_y - ik)^ma,
\nonumber
\end{equation}
\begin{equation}
\partial_x^{\alpha}\left(e^{ix\cdot\xi}y^{\frac{n-1}{2}-ik}a\right) = e^{ix\cdot\xi}y^{\frac{n-1}{2}-ik}
(\partial_x + i\xi)^{\alpha}a,
\nonumber
\end{equation}
we have the following identity
\begin{equation}
\begin{split}
& \left(H - k^2\right)\left(e^{ix\cdot\xi}y^{\frac{n-1}{2}-ik}a\right) \\
&\ \ =  e^{ix\cdot\xi}y^{\frac{n-1}{2}-ik}\left\{2ikD_ya 
- \big(D_y^2 + K(\xi)\big)a\right\},
\nonumber
\end{split}
\end{equation}
where $K(\xi)$ is a differential operator of the form
\begin{equation}
K(\xi) = y^2(\partial_x + i\xi)^2 + y^2\sum_{|\alpha|\leq2}L_{\alpha}(x,y)(\partial_x + i\xi)^{\alpha}.
\label{eq:Chap4Sect3FormofKxi}
\end{equation}
We put $a = \sum_{j=0}^Nk^{-j}a_j$. Then the above formula becomes
\begin{equation}
\begin{split}
& e^{-ix\cdot\xi}y^{-\frac{n-1}{2}+ik}\left(H - k^2\right)e^{ix\cdot\xi}y^{\frac{n-1}{2}-ik}a \\
&=  2ikD_y a_0 + \sum_{j=0}^{N-1}k^{-j}\Big\{2iD_ya_{j+1} -
\big(D_y^2 +  K(\xi)\big)a_j\Big\} \\
& -
k^{-N}\big(D_y^2 + K(\xi)\big)a_N.
\label{C4S2(H-k2)a=g}
\end{split}
\end{equation}
We put
\begin{equation}
a_0(x,y) = 1,
\label{eq:Chap4a0}
\end{equation}
and consruct $a_j$ succesively by
\begin{equation}
a_{j+1}(x,y,\xi) = - \frac{i}{2}\int_0^y(D_t^2 + K(\xi))a_j(x,t,\xi)\frac{dt}{t}.
\label{eq:Chap4Sect3jthtermoftransporteq}
\end{equation}
Then we have
\begin{equation}
2iD_ya_{j+1} -
\big(D_y^2 + K(\xi)\big)a_j = 0. 
\label{eq:Chap4Sect3ajtransporteq}
\end{equation}
We put for $p \geq 0$
$$
y^p\,{\mathcal W}^{s} = \{y^p\,w(x,y)\, ; \, w(x,y) \in {\mathcal W}^{s}\}.
$$
Here and what follows, we allow the elements of ${\mathcal W}^s$ to be complex-valued.
Then one can show easily that
\begin{equation}
\int_0^yt^qf(x,t)\frac{dt}{t} \in y^{p+q}\,{\mathcal W}^{s}, \quad {\rm if} \quad f \in y^p\,{\mathcal W}^{s}, \quad  p, q \geq 0, \quad s<0.
\label{C4S2ypqsW}
\end{equation}

In fact, letting $f(x,y) = y^pw(x,y)$, $w \in {\mathcal W}^s$, we are led to estimate
$$
y^{p+q}\int_0^1\tau^{p+q}w(x,y\tau)\frac{d\tau}{\tau}.
$$
Noting that for $0 < y < 1$
$$
\log\langle x\rangle + \langle \log(y\tau)\rangle \geq \log\langle x\rangle + 
\langle\log y\rangle,
$$
we easily get (\ref{C4S2ypqsW}).

\begin{lemma}
For $j \geq 1$, we have
$$
a_j(x,y,\xi) = y^2\xi^2P_{j-1}(y^2\xi^2) + \sum_{p=1}^jy^{2p}\sum_{|\alpha|\leq 2p}A^{(j,p)}_{\alpha}(x,y)\xi^{\alpha},
$$
where $P_{j-1}$ is a polynomial of order $j-1$ with constant coefficients, and 
$A^{(j,p)}_{\alpha}(x,y) \in {\mathcal W}^{-1-\epsilon/2}$.
\end{lemma}

Proof. The proof is by induction using (\ref{C4S2ypqsW}) and the formula
$$
\int_0^y\left(D_t^2t^{\beta}\right)\frac{dt}{t} = \beta y^{\beta}. \qquad \qquad \qquad \qquad \qed
$$

Summing up, we have proven the following theorem.


\begin{theorem}
For any $N > 0$, there exists an asymptotic solution to the equation $(H - k^2)u = 0$ such that in ${\bf R}^{n-1}\times(0,y_0)$
\begin{equation}
(H - k^2)\left(y^{\frac{n-1}{2}-ik}e^{ix\cdot\xi}\sum_{j=0}^Nk^{-j}a_j(x,y,\xi)\right) 
= y^{\frac{n-1}{2}-ik}e^{ix\cdot\xi}k^{-N}g_N(x,y,\xi),
\nonumber
\end{equation}
where $a_j(x,y,\xi)$ has the form in Lemma 2.1. Furthermore $g_N(x,y,\xi)$ has the form
\begin{equation}
g_N(x,y,\xi) = y^2\xi^2Q_{N}(y^2\xi^2) + \sum_{p=1}^{N+1}y^{2p}\sum_{|\alpha|\leq 2p}B^{(N,p)}_{\alpha}(x,y)\xi^{\alpha},
\label{C4S4gNxyxi}
\end{equation}
where $Q_N$ is a polynomial of order $N$ with constant coefficients, and 
$B^{(N,p)}_{\alpha}(x,y) \in {\mathcal W}^{-1-\epsilon/2}$.
\end{theorem}


\section{Mellin transform and pseudo-differential operators}


\subsection{Mellin transform} 
The Mellin transform $U_M$
is defined by
\begin{equation}
\left(U_{M}f\right)(k) = 
\frac{1}{\sqrt{2\pi}}\int_{0}^{\infty}y^{\frac{n-1}{2}+ik}f(y)
\frac{dy}{y^n}, \quad k \in {\bf R}.
\label{eq:MellinTtansform}
\end{equation}
In the following, the Fourier transform and its adjoint are denoted by
\begin{equation}
 F_{k\to z}f(z) = \frac{1}{\sqrt{2\pi}}\int_{-\infty}^{\infty}
 e^{-izk}f(k)dk,
 \label{C4S31dimFourier}
\end{equation}
\begin{equation}
 F_{z\to k}^{\ast}g(k) = \frac{1}{\sqrt{2\pi}}\int_{-\infty}^{\infty}
 e^{izk}g(z)dz.
 \label{C4S31dimFourieradjoint}
\end{equation}
Note that
$$
F_{z\to k}^{\ast} = (F_{k\to z})^{\ast}.
$$
Using the fact that  
$$
T : L^2((0,\infty);dy/y^n) \ni f(y) \to \left(Tf\right)(z) = f(e^z)e^{-(n-1)z/2} \in L^2({\bf R};dz)
$$
 is unitary, 
we have
\begin{equation}
(U_Mf)(k) = \left(F_{z\to k}^{\ast}Tf\right)(k) = \frac{1}{\sqrt{2\pi}}\int_{-\infty}^{\infty}
e^{izk}\left(Tf\right)(z)dz.
\label{eq:MellinandFurier}
\end{equation}
Hence $U_{M} : L^2((0,\infty);dy/y^n) \to 
L^2({\bf R}^1)$ is unitary, and the inversion formula holds:
\begin{equation}
f(y) = \frac{1}{\sqrt{2\pi}}\int_{-\infty}^{\infty}y^{\frac{n-1}{2}-ik}
\left(U_{M}f\right)(k) dk = (U_M)^{\ast}U_Mf.
\nonumber
\end{equation}
We put
\begin{equation}
K_0 = i\left(y\partial_y - \frac{n-1}{2}\right).
\label{eq:DefinitionofK0}
\end{equation}
Then we have for $f \in C_0^{\infty}((0,\infty))$
\begin{equation}
(U_{M}K_0f)(k) = k(U_Mf)(k) = F_{z\to k}^{\ast}\big(i\partial_z\big(Tf\big)\big)(k).
\label{C4S2Correspond}
\end{equation}
Therefore, for a function $\varphi(k)$ on ${\bf R}$, we define the operator $\varphi(K_0)$ by
\begin{equation}
\varphi(K_0) = \big(U_M\big)^{\ast}\varphi(k)U_M.
\label{eq:DefinitionoffunctionofK0}
\end{equation}

By (\ref{C4S2Correspond}), we have the following correspondence between the multiplication operator $k$ and the differential operators $\partial_z$, $y\partial_y$ via the Fourier transform in the $z$-space and the Mellin transform in the $y$-space:
\begin{equation}
i\left(y\partial_y - \frac{n-1}{2}\right) \longleftrightarrow k \longleftrightarrow
i\partial_z.
\label{eq:quantization}
\end{equation}
We also put for $h(x) \in L^2({\bf R}^{n-1})$
\begin{equation}
\left(F_{x\to\xi}h\right)(\xi) = \widehat h(\xi) = (2\pi)^{-(n-1)/2}\int_{{\bf R}^{n-1}}e^{-ix\cdot\xi}h(x)dx.
\nonumber
\end{equation}
Thus we have the following correspondence for the operator $H_0$ on $L^2({\bf H}^n)$ and its symbol:
\begin{equation}
\begin{split}
- D_y^2 + (n-1)D_y - \frac{(n-1)^2}{4}  - y^2\Delta_x  \longleftrightarrow & k^2 + y^2|\xi|^2 \\
=& k^2 + e^{2z}|\xi|^2 
\longleftrightarrow
- \partial_z^2 - e^{2z}\Delta_x.
\end{split}
\label{eq:HamiltonianCorrespondence}
\end{equation}

For $p(x,y,\xi,k) \in C^{\infty}({\bf R}^n_+\times{\bf R}^n)$, we define an operator $p_{FM}$ by
\begin{equation}
\left(p_{FM}f\right)(x,y) = (2\pi)^{-n/2}\int_{{\bf R}^n}e^{ix\cdot\xi}
y^{\frac{n-1}{2}-ik}p(x,y,\xi,k)(U_{M}\widehat f)(\xi,k)d\xi dk.
\label{eq:Definitionofpfm}
\end{equation}
This is rewritten as
\begin{equation}
p_{FM} = T^{\ast}\circ p_T(x,z,-i\partial_x,i\partial_z)\circ T,
\nonumber
\end{equation}
where $P_T := p_T(x,z,-i\partial_x,i\partial_z)$ is a standard pseudo-differential operator ($\Psi$DO) on ${\bf R}^n$:
\begin{equation}
\begin{split}
 \left(P_Th\right)(x,z)  = & (2\pi)^{-n}
\iint_{{\bf R}^n\times{\bf R}^n}e^{i\left((x-x')\cdot\xi - (z-z')k\right)} p_T(x,z,\xi,k)h(x',z')dx'dz'd\xi dk,
\nonumber
\end{split}
\end{equation}
with
\begin{equation}
p_T(x,z,\xi,k) = p(x,e^z,\xi,k).
\label{C4S3pTxzxuk}
\end{equation}
If $p_T(x,z,\xi,k)$ satisfies
\begin{equation}
|\partial_x^{\alpha}\partial_z^m\partial_{\xi}^{\beta}\partial_k^l\,
p_T(x,z,\xi,k)| \leq C_{\alpha\beta m l}, \quad \forall \alpha, \beta, m, l,
\label{C4S4pTestimate}
\end{equation}
$P_T$ is a bounded operator on $L^2({\bf R}^n)$ (see \cite{CalVai}).
Therefore, $p_{FM}$ is a bounded operator on $L^2({\bf H}^n)$.
Note that for the $L^2$-boundedness, it is sufficient to assume (\ref{C4S4pTestimate}) up to some finite order $|\alpha| + |\beta| + m + l \leq \mu(n)$.  

We need the following class of symbols.


\begin{definition}
For $s, t \in {\bf R}$ and $N \geq 0$, let $\widetilde S_{s,t}^N$ be the set of $C^{\infty}$-functions on ${\bf R}^{n}_+\times{\bf R}^{n}$ such that
\begin{equation}
|(\partial_x)^{\alpha}(\partial_{\xi})^{\beta}(y\partial_y)^m(\partial_k)^l \,
p(x,y,\xi,k)| \leq C(1 + |k|)^{s - l}(1 + |\xi|)^{t-\beta}
\nonumber
\end{equation}
holds for $|\alpha| + |\beta| + m + l \leq N$.
\end{definition}

We say that a $\Psi$DO $p_{FM}$ belongs to $\widetilde S_{s,t}^N$ if its symbol belongs to $\widetilde S_{s,t}^N$. We always assume that $N$ is chosen sufficiently large. Standard calculus for $\Psi DO$ applies to $p_{FM}$. For example, 
\begin{equation}
\begin{split}
p \in \widetilde S_{s,t}^{N} & \Longrightarrow
(p_{FM})^{\ast} \in \widetilde S_{s,t}^{N'}, \\
p \in \widetilde S_{s_1,t_1}^{N_1}, \ q \in \widetilde S_{s_2,t_2}^{N_2} & \Longrightarrow
p_{FM}q_{FM} \in \widetilde S_{s_1+s_2,t_1+t_2}^{N'}, \\
p \in \widetilde S_{s_1,t_1}^{N_1}, \ q \in \widetilde S_{s_2,t_2}^{N_2} & \Longrightarrow
[p_{FM},q_{FM}] \in \widetilde S_{s_1+s_2-1,t_1+t_2}^{N'}\cup \widetilde S_{s_1+s_2,t_1+t_2-1}^{N'}
\end{split}
\nonumber
\end{equation}
with suitable $N' > 0$. These can be proven in the same way as in \cite{Hor}, Vol 3, Sect. 18.1.


\subsection{Regularity of the resolvent} 


\begin{lemma} (1) Let $D_x = y\partial_x, \ D_y = y\partial_y$. Then for $N \geq 1$
$$
D_x^{\alpha}D_y^m(H + i)^{-N} \in {\bf B}(L^2({\bf H}^n))
\quad {\rm for} \quad|\alpha| + m \leq 2N.
$$
(2) Let $f \in \mathcal S$. Then we have
$$
D_x^{\alpha}D_y^{m} f(H) \in {\bf B}(L^2({\bf H}^n)), \quad 
\forall \alpha, m.
$$
\end{lemma}

Proof.  For $k \geq 0$, let $\mathcal P_k$ be the elements of $\mathcal P$,  introduced  in Chapter 2, Subsection 2.1, whose order is at most $k$. 

We shall prove (1). The case $N = 1$ is proved in Theorem 2.1.3 (4). 
Assume that the Lemma is true for $N$. Consider $D_x^{\alpha}D_y^m(H+i)^{-N-1}$ where $|\alpha|+m \leq 2(N+1)$. Let first $|\alpha|\geq 2$ so that $\alpha = \alpha' + \alpha''$, where $|\alpha''|=2$. Then 
\begin{equation}
\begin{split}
&D_x^{\alpha}D_y^m(H+i)^{-N-1} \\
&= D_x^{\alpha''}D_x^{\alpha'}D_y^m(H+i)^{-1}(H+i)^{-N} \\
&= D_x^{\alpha''}(H+i)^{-1}D_x^{\alpha'}D_y^m(H+i)^{-1} + 
D_x^{\alpha''}[D_x^{\alpha'}D_y^m,(H+i)^{-1}](H+i)^{-N}.
\end{split}
\nonumber
\end{equation}
The first term is bounded by induction hypothesis. As for the 2nd term, using Lemma 2.1.2 (1) and the definition of $\mathcal W^{-1-\epsilon/2}$, we have
$$
[D_x^{\alpha'}D_y^m,(H+i)^{-1}] = (H+i)^{-1}\Big\{\sum_{i=1}^nD_iA^{(i)} + A^{(0)}\Big\}(H+i)^{-1},
$$
where $A^{(i)} \in \mathcal P_{2N}$, and
$D_i= y\partial_{x_i}$, $1\leq i\leq n-1$, $D_n = D_y$. Thus
\begin{equation}
\begin{split}
& D_x^{\alpha''}[D_x^{\alpha'}D_y^m,(H+i)^{-1}](H+i)^{-N} \\
&=D_x^{\alpha''}(H+i)^{-1}\sum_{i=1}^nD_i(H+i)^{-1}\{A^{(i)}(H+i)^{-N} + [A^{(i)},H](H+i)^{-N}\}\\
& + D_x^{\alpha''}(H+i)^{-1}A^{(0)}(H+i)^{-N-1}. 
\end{split}
\nonumber
\end{equation}
By induction hypothesis, it is sufficient to show that $D_i(H+i)^{-1}[A^{(i)},H](H+i)^{-N}$ is bounded. Note
$$
[A^{(i)},H] = \sum_{j=1}^nD_j\widehat A^{(j)} + \widehat A^{(0)},
$$
where $\widehat A^{(j)} \in \mathcal P_{2N}$. 
However,
\begin{equation}
\begin{split}
D_i(H+i)^{-1}D_j & = D_iD_j(H+i)^{-1} + D_i[(H+i)^{-1},D_j] \\
& = D_iD_j(H+i)^{-1} + D_i(H+i)^{-1}[H,D_j](H+i)^{-1} \in {\bf B}(L^2({\bf H}^n)).
\end{split}
\nonumber
\end{equation}
Thus $D_i(H+i)^{-1}[A^{(i)},H](H+i)^{-N}$ is bounded. The case $|\alpha| < 2$, hence $m\geq 2$, is proved similarly.

Let us prove (2). Take $N$ such that $|\alpha| + m \leq 2N$ and put $g(t) = f(t)(i + t)^{N}$. Let $\widetilde g(z)$ be an almost analytic extension of $g(z)$ defined in Section 3.3.1. Then we have by Lemma 3.3.1
\begin{equation}
D_x^{\alpha}D_y^mg(H) = D_x^{\alpha}D_y^m(i + H)^{-N}
\frac{1}{2\pi i}\int_{{\bf C}}\overline{\partial_z}\widetilde g(z)(i + H)^N
(z - H)^{-1}dzd\overline{z}.
\nonumber
\end{equation}
Since $(i + H)^N(z - H)^{-1} = \sum_{r=-1}^{N-1}c_r(z)(z - H)^{r}$, $c_r(z)$ being a polynomial of $z$ of degree $N-r-1$. Therefore, taking $\sigma = -2N-2$ in Chap. 3 (\ref{C3S3HelfferSjostland}), We see that $D_x^{\alpha}D_y^m g(H)$ is a bounded operator multiplied by a polynomial of $H$ of order $N-1$. By multiplying $(i + H)^{-N}$, we obtain (2).
\qed


\section{Parametrices and regularizers}


\subsection{Wave operators and Mellin transform}
We now introduce wave operators based on the Mellin transform:
\begin{equation}
 W_{M}^{(\pm)} = \mathop{\rm s-lim}_{t\to\pm\infty}
e^{it\sqrt{H_{+}}}e^{\mp itK_0}r_{\pm}(K_0),
\label{eq:Wplusminustilde}
\end{equation}
where $H_{+} = E_H((0,\infty))H = P_{ac}(H)H$, $E_H(\lambda)$ being the spectral resolution for $H$, and $r_{+}(k)$ and $r_-(k)$ are the characteristic function of the interval $(0,\infty)$ and $(-\infty,0)$, respectively (see (\ref{eq:DefinitionoffunctionofK0})). 
Recall $\mathcal F_+$  given in Chap. 2 by formulae (\ref{eq:Chap2Sect7Fplusminusk}),  (\ref{E2.15}) and (\ref{E2.1^}).


\begin{lemma} The strong limits (\ref{eq:Wplusminustilde}) exist and
$$
\mathcal F_+ = \frac{1}{\sqrt2}\left\{r_{+}U_{M}(W_M^{(+)})^{\ast} 
+ r_{-}U_{M}(W_M^{(-)})^{\ast}\right\},
$$
where $r_{\pm}$ is the operator of multiplication by $r_{\pm}(k)$.
\end{lemma}
Proof. Due to formula (\ref{eq:quantization}) and Definition 5.3 of Chap. 1, 
we have
$$
y^{(n-1)/2}\left(\mathcal R_0f\right)(-\log y \mp t,x) = (U_M)^{\ast}\left(e^{\mp ikt}\mathcal F_0(k)f\right)(y,x).
$$
Using again (\ref{eq:quantization}) and Theorem 1.5.5, 
we see that, as $t \to \pm\infty$
\begin{equation}
\left\|e^{-it\sqrt{H_0}}f - \sqrt{2}e^{\mp itK_0}r_{\pm}(K_0)(U_M)^{\ast}
\mathcal F_0f\right\|_{L^2({\bf H}^n)} \to 0. 
\label{C4S4eitsqrtnatcalF}
\end{equation}
By Theorem 2.8.11, the wave operator 
 ${\mathop{\rm s-lim}_{t\to\pm\infty}}\, e^{it\sqrt{H_+}}e^{-it\sqrt{H_0}}$ exists and is equal to $W_{\pm} =  {\mathop{\rm s-lim}_{t\to\pm\infty}}\,e^{itH}e^{-itH_0}$. This and (\ref{C4S4eitsqrtnatcalF}) imply the existence of the  limt $W_{M}^{(\pm)}$  and
\begin{equation}
W_{\pm} = \sqrt{2}\,W_{M}^{(\pm)}(U_M)^{\ast}\mathcal F_0 
= \sqrt{2}\,W_{M}^{(\pm)}r_{\pm}(K_0)(U_M)^{\ast}\mathcal F_0.
\nonumber
\end{equation}
Letting $r_{\pm}$ be the operator of multiplication by $r_{\pm}(k)$ in $L^2({\bf R};L^2({\bf R}^{n-1});dk)$, we then have
\begin{equation}
r_{\pm}\mathcal F_0(W_{\pm})^{\ast} = \sqrt{2}\,r_{\pm}\mathcal F_0(\mathcal F_0)^{\ast}
r_{\pm}U_M(W_{M}^{(\pm)})^{\ast}.
\nonumber
\end{equation}
By Lemma 2.8.3, one can show
\begin{equation}
r_{\pm}\mathcal F_0(\mathcal F_0)^{\ast}r_{\pm} = \frac{1}{2}r_{\pm},
\nonumber
\end{equation}
which together with the formula (\ref{eq:Chap2Sect8FplusWplusminus}) in Lemma 2.8.4 proves the lemma. \qed

\bigskip
Recall that, using the 1-dimensional Fourier transform (\ref{C4S31dimFourieradjoint}), the modified Radon transform is defined by 
$$
\mathcal R_{\pm} = F^{\ast}_{k \to s}\mathcal F_{\pm},
$$
(see Definition 8.5 in Chapter 2).
Then Lemma 4.1 implies


\begin{lemma}
$$
\mathcal R_+ = \frac{1}{\sqrt2}F_{k \to s}^{\ast}
\left(r_+U_M(W_{M}^{(+)})^{\ast} + 
r_-U_M(W_{M}^{(-)})^{\ast}\right).
$$
\end{lemma}


\subsection{Parametrices for the wave equation} 
Let $a_j(x,y,\xi)$  be as in Lemma 2.1. We take $\chi_{\infty}(k) \in C^{\infty}({\bf R})$ such that $\chi_{\infty}(k) = 1 \ (|k| > 2)$, $\chi_{\infty}(k) = 0 \ (|k| < 1)$, and $\widetilde\chi(y) \in C^{\infty}({\bf R})$ such that $\widetilde\chi(y) = 1 \ (y < y_0/2), \ \widetilde\chi(y) = 0 \ (y > y_0)$, $y_0$ being a constant in Theorem 2.2. 
We  define $a^{(\pm)}(x,y,\xi,k)$ by
\begin{equation}
a^{(\pm)}(x,y,\xi,k)  = \chi_{\infty}(k)r_{\pm}(k)\sum_{j=0}^{\infty}
\rho\Big(\frac{\langle\xi\rangle^2}{\epsilon_j\langle k\rangle}\Big)k^{-j}
a_j(x,y,\xi)\widetilde\chi(y).
\label{eq:axyxik}
\end{equation}
Here, $\rho(s) \in C_0^{\infty}({\bf R})$ is such that $\rho(s) = 1$ for $|s|<1/2$, $\rho(s) = 0$ for $|s|>1$, and $\{\epsilon_j\}_{j=0}^{\infty}$ is a sequence such that $\epsilon_0 > \epsilon_1 > \cdots \to 0$. 


\begin{lemma}
For a suitable choice of $\{\epsilon_j\}_{j=0}^{\infty}$, the series (\ref{eq:axyxik}) converges and defines a smooth function having the following properties:
\\
\smallskip
\noindent
(1) $\ {\rm supp}\,a^{(\pm)}(x,y,\xi,k) \subset {\bf R}^{n-1}\times(0,y_0)\times\{(\xi,k)\, ; \, |k| \geq 1, \ \langle \xi\rangle^2\leq \epsilon_0\langle k\rangle\}$. \\
\noindent
(2) If $|\beta| + m + |\gamma| +\ell \leq N$, we have, 
\begin{equation}
\begin{split}
& \Big|\partial_x^{\beta}D_y^m\partial_{\xi}^{\gamma}\partial_k^{\ell}\Big(a^{(\pm)}(x,y,\xi,k) - \chi_{\infty}(k)r_{\pm}(k)\sum_{j=0}^N\rho\Big(\frac{\langle\xi\rangle^2}{\epsilon_j\langle k\rangle}\Big) k^{-j}
a_j(x,y,\xi)\widetilde\chi(y)\Big)\Big|\\
& \ \ \ \ \ \ \ \ \ \  \leq C_{N\beta m \gamma\ell}\,
y^2\left(\frac{\langle\xi\rangle^2}{\langle k\rangle}\right)^{N}\langle\xi\rangle^{-|\gamma|}\langle k\rangle^{-\ell}.
\end{split}
\label{C4S4ajxyxikminusNterms}
\end{equation}
(3) Let $g^{(\pm)}(x,y,\xi,k)$ be defined by 
\begin{equation}
(H - k^2)y^{\frac{n-1}{2}-ik}e^{ix\cdot\xi}a^{(\pm)}(x,y,\xi,k) = 
y^{\frac{n-1}{2}-ik}e^{ix\cdot\xi}g^{(\pm)}(x,y,\xi,k).
\label{eq:Chap4Sect4(L-k2)a=g}
\end{equation}
Then we have for any $N > 0$
\begin{equation}
\Big|\partial_x^{\beta}D_y^m\partial_{\xi}^{\gamma}\partial_k^{\ell}\, g^{(\pm)}(x,y,\xi,k)\Big| \leq 
C_{N\beta m \gamma\ell}\,
y^2\left(\frac{\langle\xi\rangle^2}{\langle k\rangle}\right)^{N}\langle\xi\rangle^{2-|\gamma|}\langle k\rangle^{2-\ell}.
\label{C4S4gpmestimate}
\end{equation}
for $y < y_0/2$ and $\langle\xi\rangle^2 \leq \epsilon_{N+1}\langle k\rangle/2$.
\end{lemma}
 
Proof.  First we derive the following estimate for $j \geq 1$
\begin{equation}
\begin{split}
\left|\partial_x^{\beta}D_y^m\partial_{\xi}^{\gamma}\partial_k^{\ell} 
\left(\rho\Big(\frac{\langle\xi\rangle^2}{\epsilon_j\langle k\rangle}\Big) k^{-j}
a_j(x,y,\xi)\widetilde\chi(y)\right)\right|\\
 \leq C'_{j\beta m \gamma\ell}\,
y^2\left(\frac{\langle\xi\rangle^2}{\langle k\rangle}\right)^{j}\langle\xi\rangle^{-|\gamma|}\langle k\rangle^{-\ell},
\end{split}
\label{C4S4ajxyxikinequality}
\end{equation}
where the constant $C'_{j\beta m \gamma\ell}$ is independent of $\epsilon_j$. 
In fact, by Lemma 2.1,
\begin{equation}
k^{-j}a_j(x,y,\xi)\widetilde\chi(y) = \sum_{|\alpha|\leq 2j}a_{j,\alpha}(x,y)\frac{\xi^{\alpha}}{k^j},
\nonumber
\end{equation}
where $a_{j,\alpha}(x,y) = 0$ for $y > y_0$, and
\begin{equation}
|\partial_x^{\beta}D_y^{m} a_{j,\alpha}(x,y)| \leq C'_{j\beta m}\, y^2, \quad 
\forall \beta, m.
\nonumber
\end{equation}
We define a homegenous polynomial of $(\sigma,\eta) \in {\bf R}^{n}$ by
\begin{equation}
b_j^{(\pm)}(x,y,\sigma,\eta) = (\pm 1)^j\sum_{|\alpha|\leq 2j}
a_{j,\alpha}(x,y)\sigma^{2j-|\alpha|}\eta^{\alpha}.
\nonumber
\end{equation}
We then have
\begin{equation}
k^{-j}a_j(x,y,\xi)\widetilde\chi(y) = b_j^{(\pm)}\big(x,y,\frac{1}{\sqrt{|k|}},\frac{\xi}{\sqrt{|k|}}\big), \quad {\rm for} \quad \pm k >0.
\nonumber
\end{equation}
Put $\Xi = (1/\sqrt{|k|},\xi/\sqrt{|k|})$, and note that
$$
|\partial_{\xi}^{\gamma}\partial_k^{\ell}\Xi| \leq C'_{\beta\ell}\langle \Xi\rangle\langle\xi\rangle^{-|\beta|}| k|^{-|\ell|}
\leq C_{\beta\ell}\langle \xi\rangle^{1-|\gamma|}|k|^{-\ell-1/2}, \ |k| > 1.
$$
Taking into account of the homogeneity of $b_j^{(\pm)}(x,y,\sigma,\eta)$, we then have 
$$
\Big|\partial_x^{\beta}\partial_y^m\partial_{\xi}^{\gamma}\partial_k^{\ell}\,
b_j^{(\pm)}\big(x,y,\frac{1}{\sqrt{|k|}},\frac{\xi}{\sqrt{|k|}}\big)\Big| \leq
C'_{j\beta m \gamma \ell}\, y^2\left(\frac{\langle\xi\rangle^2}{\langle k\rangle}\right)^{j}\langle\xi\rangle^{-\gamma}\langle k\rangle^{-\ell}.
$$
This, together with the inequality,
$$
\left|\partial_x^{\beta}D_y^m\partial_{\xi}^{\gamma}\partial_k^{\ell} \,
\rho\Big(\frac{\langle\xi\rangle^2}{\epsilon_j\langle k\rangle}\Big)\right|
\leq
C'_{\beta m \gamma \ell}\langle\xi\rangle^{-|\gamma|}\langle k\rangle^{-\ell},
$$ 
where the constant $C'_{\beta m \gamma \ell}$ is independent of $\epsilon_j$, gives (\ref{C4S4ajxyxikinequality}). Noting that $\langle\xi\rangle^2/\langle k\rangle \leq \epsilon_j$, we then have
\begin{equation}
\begin{split}
\left|\partial_x^{\beta}D_y^m\partial_{\xi}^{\gamma}\partial_k^{\ell} 
\left(\rho\Big(\frac{\langle\xi\rangle^2}{\epsilon_j\langle k\rangle}\Big) k^{-j}
a_j(x,y,\xi)\widetilde\chi(y)\right)\right|\\
 \leq C'_{j\beta m \gamma\ell}\,
y^2\epsilon_j\left(\frac{\langle\xi\rangle^2}{\langle k\rangle}\right)^{j-1}\langle\xi\rangle^{-|\gamma|}\langle k\rangle^{-\ell},
\end{split}
\label{C4S4ajxyxikinequalitymodify}
\end{equation}
Take $\epsilon_j$ such that
$$
(1 + C'_{j\beta m \gamma\ell})\epsilon_j < 2^{-j}, \quad |\beta| + m + |\gamma| +\ell \leq j.
$$
Then, by (\ref{C4S4ajxyxikinequalitymodify}),  the series (\ref{eq:axyxik}) converges uniformly with all of its derivatives. The inequality (\ref{C4S4ajxyxikminusNterms}) also follows from (\ref{C4S4ajxyxikinequalitymodify}). 
We put
$$
g_{N+1}^{(\pm)} = y^{-\frac{n-1}{2}}e^{-ix\cdot\xi}(H-k^2)y^{\frac{n-1}{2}}e^{ix\cdot\xi}\chi_{\infty}(k)r_{\pm}(k)
\sum_{j=0}^{N}\rho\left(\frac{\langle\xi\rangle^2}{\epsilon_j\langle k\rangle}\right)k^{-j}a_j(x,y,\xi)\widetilde\chi(y),
$$
and $\widetilde g_{N+1}^{(\pm)} = g^{(\pm)} - g^{(\pm)}_{N+1}$. 
Then by (\ref{C4S2(H-k2)a=g}), $g^{(\pm)}_{N+1} = 0$ for $\langle \xi\rangle^2 \leq \epsilon_{N+1}\langle k\rangle/2$ and $y < y_0/2$.
The inequality (\ref{C4S4ajxyxikinequalitymodify}) shows that $\widetilde g_{N+1}^{(\pm)}$ has the estimate in (3). \qed

\bigskip
We define an operator $U_{\pm}(t)$  by
\begin{equation}
U_{\pm}(t) = a^{(\pm)}_{FM}e^{\mp itK_0}\chi(y).
\label{eq:Uplusminust}
\end{equation}
where $\chi(y) \in C^{\infty}({\bf R})$ is such that $\chi(y) = 1 \ (y < y_0/4)$, $\chi(y) = 0 \ (y > y_0/3)$. As in the analysis for the operators $p_{FM}$ (see (\ref{eq:Definitionofpfm}) and thereafter), $a_{FM}^{(\pm)}$ are bounded on $L^2({\bf H}^n)$, and therefore $U_{\pm}(t)$. 
The explicit form of $U_{\pm}(t)$ is as follows:
\begin{equation}
\begin{split}
& \left(U_{\pm}(t)f\right)(x,y) \\
& = (2\pi)^{-\frac{n}{2}}
\int_{{\bf R}^n}e^{ix\cdot\xi}y^{\frac{n-1}{2}-ik}
a^{(\pm)}(x,y,\xi,k)e^{\mp itk}
\big(U_M\chi(y)\widehat f\big)(\xi,k)d\xi dk.
\end{split}
\label{eq:Uplusminusexplicite}
\end{equation}
 We put
\begin{equation}
G_{\pm}(t) = \frac{d}{dt}\Big(e^{it\sqrt{H_+}}U_{\pm}(t)\Big),
\label{eq:Pplusminust}
\end{equation}
and also
\begin{equation}
\Lambda_y = (1 + K_0^2)^{1/2} = (U_{M})^{\ast}(1 + k^2)^{1/2}
U_M,
\label{C4S4Lambday}
\end{equation}
\begin{equation}
\Lambda_x = (1 - \Delta_x)^{1/2} = (F_{x\to\xi}) ^{\ast}(1 + |\xi|^2)^{1/2}F_{x\to\xi}.
\label{C4S4Lambdax}
\end{equation}


\begin{lemma} There exists $N_0 > 0$ such that for any $N > N_0$, there exists a constant $C_{N} > 0$ for which
\begin{equation}
\|G_{\pm}(t)\Lambda_x^{-2N}\Lambda_y^{N/2}\| \leq C_N(1 + |t|)^{-2}, \quad {\rm for} \quad \pm t > 0,
\label{C4S4GtLambda}
\end{equation}
holds, 
where $\|\cdot\|$ denotes the operator norm of $L^2({\bf H}^n)$.
\end{lemma}

Proof.  We consider $G_+(t)$, which is rewritten as
\begin{equation}
G_+(t) = e^{it\sqrt{H_+}}\Big(i\sqrt{H_+}U_+(t) + \frac{d}{dt}U_+(t)\Big).
\nonumber
\end{equation}
Letting $H = \int_{-\infty}^{\infty}\lambda dE_H(\lambda)$, 
we deal with the high energy part and low energy part separately, i.e. 
on the subspace $E_H([1,\infty))L^2({\bf H}^n)$, and $E_H((-\infty,1))L^2({\bf H}^n)$.

\medskip
\noindent
{\it High energy part}. We take $\chi_0(s) \in C^{\infty}_0({\bf R})$ such that
$\chi_0(s) = 1$ for $ -\infty  < s < 1/4$, $\chi_0(s) = 0$ for $s > 1/2$. We  consider $\displaystyle{i\sqrt{H}(1 - \chi_0(H))U_+(t) + \frac{d}{dt}U_+(t)}$. We put $f(s) = s^{-1/2}(1 - \chi_0(s))$. 


\begin{prop}
If $f(s) \in C^{\infty}({\bf R})$ satisfies  for some $\epsilon > 0$, $|f^{(m)}(s)| \leq C_m(1 + |s|)^{-\epsilon-m}, \forall m \geq 0$,  the following formula holds:
\begin{equation}
f(H)\,a^{(\pm)}_{FM} = a^{(\pm)}_{FM}f(K_0^2) + B^{(\pm)},
\nonumber
\end{equation}
\begin{equation}
B^{(\pm)} = \frac{1}{2\pi i}\int_{{\bf C}}\overline{\partial_z}F(\zeta)
\,(\zeta - H)^{-1}g^{(\pm)}_{FM}(\zeta - K_0^2)^{-1}d\zeta d\overline{\zeta},
\label{C4S4Bpm}
\end{equation}
where $F(\zeta)$ is an almost analytic extension of $f$, and $g^{(\pm)}(x,y,\xi,k)$ is defined by (\ref{eq:Chap4Sect4(L-k2)a=g}).
\end{prop}

Proof. Rewriting (\ref{eq:Chap4Sect4(L-k2)a=g}) into the operator form, we have
\begin{equation}
H\,a^{(\pm)}_{FM} = a^{(\pm)}_{FM}K_0^2 + 
g^{(\pm)}_{FM},
\nonumber 
\end{equation}
hence
\begin{equation}
\begin{split}
(\zeta - H)^{-1}\,a^{(\pm)}_{FM} = & a^{(\pm)}_{FM}\,(\zeta - K_0)^{-1}  + 
(\zeta - H)^{-1}g^{(\pm)}_{FM}(\zeta - K_0^2)^{-1}.
\end{split}
\nonumber 
\end{equation}
The proposition then follows from Lemma 3.3.1. \qed

\bigskip
Let us continue the proof for the high energy part. 
We consider the case $t \geq 0$. The case $t\leq 0$ is treated similarly.
Using Proposition 4.4, we have
\begin{equation}
\begin{split}
\sqrt{H}(1 - \chi_0(H))a^{(+)}_{FM} &= f(H)Ha^{(+)}_{FM} \\
& = f(H)a^{(+)}_{FM}K_0^2 + f(H)g^{(+)}_{FM} \\
& = a^{(+)}_{FM}f(K_0^2)K_0^2 + B^{(+)}K_0^2 + f(H)g^{(+)}_{FM}.
\end{split}
\nonumber
\end{equation}
Since $\frac{d}{dt}U_+(t) = - ia_{FM}^{(+)}K_0e^{-itK_0}\chi(y)$, we arrive at
\begin{equation}
\begin{split}
& i\sqrt{H}(1 - \chi_0(H))U_+(t) + \frac{d}{dt}U_+(t) \\
&= iB^{(+)}K_0^2e^{-itK_0}\chi(y) + if(H)g^{(+)}_{FM}e^{-itK_0}\chi(y)\\ 
& \ \ \ - ia^{(+)}_{FM}K_0\chi_0(K_0^2)e^{-itK_0}\chi(y).
\end{split}
\label{eq:HIghenergypartformula}
\end{equation}
Let us note here that
\begin{equation}
a^{(+)}_{FM}K_0\chi_0(K_0^2)=0,
\label{C4S4afMchi0k=0}
\end{equation}
since $|k|\geq 1$ on the support of the symbol of $a^{(+)}_{FM}$, and $\chi_0(k^2) = 0$ if $|k|\geq 1$.

\medskip
Formulae (\ref{C4S4Bpm}) and (\ref{eq:HIghenergypartformula}) contain the operators of the form $g^{(+)}_{FM}e^{-itK_0}\chi(y)$. 
We start with the following result.


\begin{prop}
Assume that $b(x,y,\xi,k) \in C^{\infty}({\bf R}^n_+\times{\bf R}^n)$ have the following properties:
$b(x,y,\xi,k) = 0$ for $ y > y_0$, and there exist $\sigma_0, \tau_0 \in {\bf R}$ such that for any $M, \alpha, m, \beta, \ell$, 
\begin{equation}
 |\partial_x^{\alpha}D_y^m\partial_{\xi}^{\beta}\partial_k^l\, b(x,y,\xi,k)| \leq C_{M\alpha\beta m\ell}\,
 \langle \log y\rangle^{-M}\langle \xi\rangle^{\sigma_0-|\beta|}\langle k\rangle^{\tau_0-\ell},  
\label{C4S4bxyxikestimate}
\end{equation}
for $0 <y < y_0$. 
Let $\chi(y) \in C^{\infty}({\bf R})$ be such that $\chi(y) = 1$ for $0 < y < y_0/4$ and $\chi(y) = 0$ for $y > y_0/3$. Then we have for any $N > 0$, and $\sigma > \sigma_0 + n/2$,
\begin{equation}
\|b_{FM}e^{-itK_0}\chi(y)\Lambda_x^{-\sigma}\Lambda_y^{N}\| \leq C_{\sigma,N}(1 + t)^{-N}, \quad  t > 0.
\label{C4S4bFMe-itK0chiyLambda}
\end{equation} 
\end{prop}

Proof. Take $\psi_0(s) \in C^{\infty}({\bf R})$ such that $\psi_0(s) = 1$ for $|s| < 1$, and $\psi(s) = 0$ for $|s| > 2$, and let for $\epsilon > 0$
$$
b^{(\epsilon)}(x,y,\xi,k) = b(x,y,\xi,k)\psi_0(\epsilon |\xi|)\psi_0(\epsilon k).
$$
Then $b^{(\epsilon)}(x,y,\xi,k)$ satisfies (\ref{C4S4bxyxikestimate}) with constant $C_{M\alpha\beta m\ell}$ independent of $\epsilon > 0$.

 We have, by (\ref{C4S4Lambdax}), (\ref{eq:MellinTtansform}) and (\ref{eq:quantization}),
\begin{equation}
\begin{split}
& b^{(\epsilon)}_{FM}e^{-itK_0}\chi(y)\Lambda_x^{-\sigma}\Lambda_y^{N}f \\
& = (2\pi)^{-\frac{n}{2}}\int_{{\bf R}^n\times{\bf R}_+} e^{ix\cdot\xi}
e^{-ik(t + \log(y/y'))}b^{(\epsilon)}(x,y,\xi,k)\\
& \ \ \ \ \ \ \ \ \ \ \ \ \ \ \ \ \times \chi(y')\langle\xi\rangle^{-\sigma}(yy')^{\frac{n-1}{2}}\Lambda_{y'}^{N}\widehat f(\xi,y')\frac{d\xi dy' dk}{(y')^n} \\
&= \sqrt{2\pi}\left(\,T^{\ast}\circ b^{(\epsilon)}_T(x,z,-i\partial_x,i\partial_z) e^{t\partial_z}\Lambda_x^{-\sigma}\chi(e^z)(1-\partial_z^2)^{N/2}\circ T\right)f.
\end{split}
\label{C4S4IntegralForm}
\end{equation}
Therefore, the estimate of this operator comes down to the calculus of classical, i.e. Euclidean, $\Psi$DO's. For the sake of completeness, we provide a  proof. 

Without loss of generality, we assume that $N/2$ is an integer. Since $(1 - \partial_z^2)^{N/2}$ is a differential operator, commuting $\chi(e^z)$ and $(1-\partial_z^2)^{N/2}$, we see that 
$$
b^{(\epsilon)}_{FM}e^{-itK_0}\chi(y)\Lambda_x^{-\sigma}\Lambda_y^{N/2}
 = T^{\ast}\circ b_T^{\mathcal O,\epsilon}(t,x,z,z',-i\partial_x,i\partial_z)\circ T,
$$
where 
\begin{equation}
\begin{split}
 b_T^{\mathcal O,\epsilon}u & = \left(b_T^{\mathcal O,\epsilon}(t,x,z,z',-i\partial_x,i\partial_z)u\right)(x,z) \\
& = 
\int_{{\bf R}^{n+1}}e^{-ik(t+z-z')}e^{ix\cdot\xi}b_T^{\mathcal O,\epsilon}(x,z,z',\xi,k)\hat u(\xi,z')dz'dkd\xi,
\end{split}
\label{C4S4bTOu}
\end{equation}
Due to (\ref{C4S4bxyxikestimate}), $b^{\mathcal O,\epsilon}_T(x,z,z',\xi,k) \in C^{\infty}({\bf R}^{n+1} \times{\bf R}^n)$ satisfies

\begin{equation}
|\partial_x^{\alpha}\partial_z^m\partial_{z'}^{m'}\partial_{\xi}^{\beta}\partial_k^{\ell}\, b_T^{\mathcal O,\epsilon}(x,z,z',\xi,k)| \leq C_{M\alpha\beta m m'\ell}\,
 \langle z \rangle^{-M}\langle \xi\rangle^{\sigma_0-\sigma -|\beta|}\langle k\rangle^{N + \tau_0-\ell},  
\nonumber
\end{equation}
with constant $C_{M\alpha\beta m m'\ell}$ independent of $\epsilon > 0$, 
and  $ b_T^{\mathcal O,\epsilon}(x,z,z',\xi,k) =0$ when $z' > \log (y_0/3)$.
Since $y_0$ is small enough, $z'< 0$ on the support of the integrand of $b_{T}^{\mathcal O,\epsilon}u$. Hence we have
$$
t  - z' \geq C_0\langle t\rangle, \quad t-z' \geq C_0\langle z'\rangle,  \quad \forall t > 0
$$
for some constant $C_0 > 0$.
Using
$$
e^{-ik(t-z')} = (-i(t-z'))^{-1}\partial_ke^{-ik(t-z')},
\quad
e^{ix\cdot\xi} = (1 +|x|^2)^{-1}(1 - \Delta_
{\xi})e^{ix\cdot\xi},
$$
we integrate $2N + [\tau_0] + 2$ times with respect to $k$ and $n$ times with respect to $\xi$ to have
$$
\Big|\left(b_{T}^{\mathcal O,\epsilon}u\right)(x,z)\Big|
\leq \int_{{\bf R}^{n+1}}A(t,z,z',x,\xi,k)|\widehat u(\xi,z')|dz'd\xi dk,
$$
$$
0 \leq A \leq C\langle t\rangle^{-N}\langle z\rangle^{-1}\langle z'\rangle^{-1}
\langle x\rangle^{-2n}\langle\xi\rangle^{\sigma_0-\sigma}\langle k\rangle^{-1}.
$$
 Then the above estimate together with Cauchy-Schwarz inequality shows that
\begin{equation}
\|b^{\mathcal O,\epsilon}_{T}u\| \leq C(1 + t)^{-N}\|u\|,
\label{C4S4bT2udexay}
\end{equation}
uniformly in $\epsilon > 0$. Letting $\epsilon \to 0$, we have
 (\ref{C4S4bFMe-itK0chiyLambda}). \qed

\medskip
By (\ref{C4S4gNxyxi}), we then see that the 2nd term of the right-hand side of (\ref{eq:HIghenergypartformula}) has the estimate
\begin{equation}
\|f(H)g_{FM}^{(+)}e^{-itK_0}\chi(y)\Lambda_x^{-2N}\Lambda_y^{N/2}\| \leq C_N(1 + t)^{-2}, \quad t \geq 0.
\label{C4S42ndtermdecayestimate}
\end{equation}

To deal with the 1st term, we use the representation (\ref{C4S4Bpm}). To apply Proposition 4.6, we consider
\begin{equation}
\begin{split}
& \ g_{FM}^{(+)}(\zeta - K_0^2)^{-1}K_0^2e^{-itK_0}\chi(y)\Lambda_x^{-2N}
\Lambda_y^{N/2} \\
= & \  g_{FM}^{(+)}K_0^2e^{-itK_0}\Lambda_x^{-2N}
\Lambda_y^{N/2}(\zeta - K_0^2)^{-1}\chi(y) \\
= & \  g_{FM}^{(+)}K_0^2e^{-itK_0}\Lambda_x^{-2N}
\Lambda_y^{N/2}\chi_1(y)(\zeta - K_0^2)^{-1}\chi(y) \\
 & + \  g_{FM}^{(+)}K_0^2e^{-itK_0}\Lambda_x^{-2N}
\Lambda_y^{N/2}\chi_2(y)(\zeta - K_0^2)^{-1}\chi(y),
\end{split}
\label{C4S4gFM(zeta-K02)split}
\end{equation}
where $\chi_1, \chi_2 \in C^{\infty}({\bf R})$, $\chi_1(y) + \chi_2(y) = 1$, $\chi_1(y) = 0$ for $y > y_0$, $\chi_2(y) = 0$ for $y < y_0/2$.
Then, Proposition 4.6 is applicable to the term $g_{FM}^{(+)}K_0^2e^{-itK_0}\Lambda_x^{-2N}
\Lambda_y^{N/2}\chi_1(y)$, and we see that the
1st term of the right-hand side of (\ref{C4S4gFM(zeta-K02)split}) is estimated as
\begin{equation}
 \|g_{FM}^{(+)}K_0^2e^{-itK_0}\Lambda_x^{-2N}
\Lambda_y^{N/2}\chi_1(y)(\zeta - K_0^2)^{-1}\| \leq C|{\rm Im}\, \zeta|^{-1}( 1 + t)^{-2}.
\label{C4S41stterm}
\end{equation}

 The 2nd term of the right-hand side of (\ref{C4S4gFM(zeta-K02)split}) is rewritten as
$$
 g_{FM}^{(+)}K_0^2e^{-itK_0}\Lambda_x^{-2N}
\Lambda_y^{N/2}\langle \log y\rangle^{-2}\cdot \langle \log y \rangle^2\chi_2(y)(\zeta - K_0^2)^{-1}\chi(y).
$$
As in the proof of Proposition 4.6, we represent $g_{FM}^{(+)}K_0^2e^{-itK_0}\Lambda_x^{-2N}
\Lambda_y^{N/2}\langle \log y\rangle^{-2}$ into the integral form like (\ref{C4S4IntegralForm}), and integrate by parts 2 times by using
$e^{-ikt} = (-it)^{-1}\partial_ke^{-ikt}$ and also (\ref{C4S4gpmestimate}).
Then we have
$$
\| g_{FM}^{(+)}K_0^2e^{-itK_0}\Lambda_x^{-2N}
\Lambda_y^{N/2}\langle \log y\rangle^{-2}\| \leq C(1 + t)^{-2}.
$$

Passing to the variable $z = \log y$, the operator $\langle \log y \rangle^2\chi_2(y)(\zeta - K_0^2)^{-1}\chi(y)$ has an integral kernel
$$
K(z,z';\zeta) = - \langle z\rangle^2\chi_2(e^z)\,\frac{\pi i}{2\sqrt{\zeta}}e^{i\sqrt{\zeta}(z-z')}\,\chi(e^{z'}).
$$
Observing the supports of $\chi_2(e^z)$ and $\chi(e^{z'})$, we see that
$z > \log(y_0/2)$, $z' < \log(y_0/3)$. Hence
\begin{equation}
z - z' \geq C\left(\langle z\rangle + \langle z'\rangle\right),
\label{C4S4z-z'>C}
\end{equation}
for a constant $C > 0$. Letting $\sqrt{\zeta} = \sigma + i\tau$, we then have
$$
|K(z,z';\zeta)| \leq \frac{C}{|\sigma| + |\tau|}\langle z\rangle^2\chi_2(e^z)\chi(e^{z'})e^{-\tau(z-z')}.
$$
Using the inequality
$$
e^{-t} \leq C_{\ell}t^{-\ell}, \quad \forall t > 0, \quad \forall \ell \geq 0,
$$
and taking $\ell = 2m + 2$, we have
$$
|K(z,z';\zeta)| \leq \frac{C_m}{\tau^{2m+3}}\langle z\rangle^{-m}\langle z'\rangle^{-m}.
$$
Taking $m > 1$, we then have
$$
\sup_z\int_{{\bf R}}|K(z,z';\zeta)|dz' \leq \frac{C_m}{\tau^{2m+3}}, \quad
\sup_{z'}\int_{{\bf R}}|K(z,z';\zeta)|dz \leq \frac{C_m}{\tau^{2m+3}}.
$$
Noting that
$$
\frac{1}{|\tau|} = \frac{2|\sigma|}{|{\rm Im}\, \zeta|} \leq 
\frac{2|\zeta|^{1/2}}{|{\rm Im}\,\zeta|},
$$
we have obtained the estimate of the operator norm
$$
\|\langle \log y \rangle^2\chi_2(y)(\zeta - K_0^2)^{-1}\chi(y)\| 
\leq C_{p}\Big(\frac{|\zeta|^{1/2}}{|{\rm Im}\,\zeta|}\Big)^p, \quad \forall p > 5.
$$
Therefore, for $p > 5$,
\begin{equation}
\begin{split}
& \|g_{FM}^{(+)}K_0^2e^{-itK_0}\Lambda_x^{-2N-n}
\Lambda_y^{N/2}\chi_2(y)(\zeta - K_0^2)^{-1}\chi(y)\| \\
\leq & C_p|{\rm Im}\,\zeta|^{-p}|\zeta|^{p/2}(1 + t)^{-2}, \quad \forall N > 0.
\end{split}
\label{C4S42ndtermestimate}
\end{equation}
Since
$$
\frac{1}{|{\rm Im}\,\zeta|} \leq \frac{\langle\zeta\rangle^{p-1}}{|{\rm Im}\,\zeta|^p}, \quad \frac{|{\zeta}|^{p/2}}{|{\rm Im}\,\zeta|^p} \leq 
\frac{\langle\zeta\rangle^{p-1}}{|{\rm Im}\,\zeta|^p},
$$
In view of (\ref{C4S41stterm}) and (\ref{C4S42ndtermestimate}), we have, 
for $p>5$,
\begin{equation}
\|g_{FM}^{(+)}(\zeta - K_0^2)^{-1}K_0^2e^{-itK_0}\chi(y)\Lambda_x^{-2N}
\Lambda_y^{N/2}\| \leq C|{\rm Im}\,\zeta|^{-p}\langle \zeta\rangle^{p-1}(1 + t)^{-2}.
\nonumber
\end{equation}
We use Lemma 2.3.1, and take into account that $\sigma$ in Chap. 2 (\ref{C3S3Fzestimate}) is now equal to $-1/2$ to see that
 the 1st term of the righ-hand side of (\ref{eq:HIghenergypartformula}) has the  property 
\begin{equation}
\|B^{(+)}K_0^2e^{-itK_0}\chi(y)\Lambda_x^{-2N}\Lambda_y^{N/2}\| \leq C_N(1 + t)^{-2}, \quad t \geq 0.
\label{C4S41sttermdecay}
\end{equation}

\bigskip
\noindent
{\it Low energy part}. We show
\begin{equation}
\|\chi_0(H)U_+(t)\Lambda_x^{-2N}\Lambda_y^{N/2}\| \leq C(1 + t)^{-2}, \quad \forall t \geq 0.
\label{eq:Chap4Sect5Lowenergypartestimate}
\end{equation}
However, noting that
\begin{equation}
\chi_0(H)a^{(+)}_{FM} = a^{(+)}_{FM}\chi_0(K_0^2) + B^{(+)} = B^{(+)},
\nonumber
\end{equation}
with $B^{(+)}$ given in Proposition 4.4, one can prove (\ref{eq:Chap4Sect5Lowenergypartestimate}) in the same way as above. 

\medskip
By (\ref{C4S42ndtermdecayestimate}), (\ref{C4S41sttermdecay}) and (\ref{eq:Chap4Sect5Lowenergypartestimate}), we have proven Lemma 4.4. \qed


\begin{lemma}
\begin{equation}
\mathop{{\rm s-lim}}_{t\to\pm\infty}e^{it\sqrt{H_{>0}}}U_{\pm}(t) = 
\chi_{\infty}(K_0)
W_{M}^{(\pm)}\chi(y).
\nonumber
\end{equation}
\end{lemma}
Proof. Since $U_{\pm}(t)$ is uniformly bounded in  $t$, we have only to prove the lemma on a dense set of $L^2({\bf H}^n)$. Writing
$$
a^{(\pm)}(x,y,\xi,k) = \chi_{\infty}(k)r_{\pm}(k) + \widetilde a^{(\pm)}(x,y,\xi,k),
$$
the same analysis as in Proposition 4.4 shows that $\|\widetilde a_{FM}^{(\pm)}e^{-itK_0}\chi(y)f\| \to 0$ for $f \in C_0^{\infty}({\bf R}^n)$. Therefore, we have
\begin{equation}
\|U_{\pm}(t)f - (U_M)^{\ast}e^{\mp itk}\chi_{\infty}r_{\pm}U_M\chi(y)f\| \to 0,
\nonumber
\end{equation}
as $t \to \pm\infty$ for any $f \in C_0^{\infty}({\bf H}^n)$. This together with (\ref{eq:Wplusminustilde}) proves the lemma.
\qed

\bigskip
Recall that for any interval  $I \subset (0,\infty)$, $\sigma \in {\bf R}$ and an integer $m \geq 0$,
\begin{equation}
\begin{split}
 & H^{\sigma,m}({\bf R}^{n-1}\times I) \ni f \\
 & \Longleftrightarrow \|f\|_{H^{\sigma,m}({\bf R}^{n-1}\times I)}^2 
 = 
\sum_{0\leq l \leq m}
\int_{{\bf R}^{n-1}\times I}|\langle \xi\rangle^{\sigma}\partial_y^l{\widehat f}(\xi,y)|^2d\xi dy < \infty.
\end{split}
\nonumber
\end{equation}
Using the standard Sobolev space $H^{\sigma,\tau}({\bf R}^n)$, where $\sigma, \tau \in {\bf R}$, we define $H^{\sigma,\tau}({\bf H}^n) = T^{\ast}H^{\sigma,\tau}({\bf R}^n)$. Then 
\begin{equation}
\begin{split}
H^{\sigma,\tau}({\bf H}^n) \ni f \Longleftrightarrow
 \|f\|_{H^{\sigma,\tau}}({\bf H}^n) & = \|Tf\|_{H^{\sigma,\tau}}({\bf R}^n) \\
& = \|\langle\xi\rangle^{\sigma}\langle k\rangle^{\tau}(U_M\hat f)(\xi,k)\|_{L^2({\bf R}^n)} < \infty.
\end{split}
\nonumber
\end{equation}

\medskip
Take $f \in H^{2N,0}$ for  large $N$. By Lemma 4.4,  $\chi_{\infty}(K_0)\int_0^{\pm\infty}G_{\pm}(t)\chi(y)fdt$  converges strongly in $L^2$. Moreover,
by (\ref{eq:Pplusminust}) and Lemma 4.7,
\begin{equation}
\chi_{\infty}(K_0)W_{M}^{(\pm)}\chi(y)f = \chi_{\infty}(K_0)a^{(\pm)}_{FM}\chi(y)f 
+ \chi_{\infty}(K_0)\int_0^{\pm\infty}G_{\pm}(t)f\, dt.
\label{C4S4WpmaFMintG}
\end{equation}
Therefore, the integral of the right-hand side can be extended by continuity as an operator in ${\bf B}(L^2;L^2)$. 

In view of Lemma 4.2 and (\ref{C4S4WpmaFMintG}), we have
\begin{equation}
\begin{split}
\mathcal R_+ & = \frac{1}{\sqrt2}F_{k\to s}^{\ast}
\Big(r_+U_M\chi (a^{(+)}_{FM})^{\ast} + r_-U_M\chi (a^{(-)}_{FM})^{\ast} \\
& \ \ \ \ \ \ \ \ \ \ \ \ \ \ + r_+U_M(1 - \chi)(W_{M}^{(+)})^{\ast} 
+ r_-U_M(1 - \chi)(W_{M}^{(-)})^{\ast}\Big) + R,
\end{split}
\label{eq:Chap4Sect5mathcalRplusexprseeion1}
\end{equation}
where $R$ is written as
\begin{equation}
R = \frac{1}{\sqrt{2}}F^{\ast}_{k\to s}\Big(r_+U_M\int_0^{\infty}G_+(t)^{\ast}\chi_{\infty}(K_0)dt + r_-U_M\int_0^{-\infty}G_-(t)^{\ast}\chi_{\infty}(K_0)dt \Big).
\nonumber
\end{equation}
Observe that since  $\int_0^{\pm\infty}G_{\pm}(t)^{\ast}dt$ enjoys the  property$$
\int_0^{\pm\infty}G_{\pm}(t)^{\ast}dt\,\chi_{\infty}(K_0) \in {\bf B}(L^2;H^{-2N,N/2})\cap {\bf B}(L^2;L^2),
$$
by interpolation,
\begin{equation}
R \in {\bf B}(L^2;H^{-\sigma,\sigma/4}), \quad \forall \sigma \geq 0.
\label{C4S4Restimate}
\end{equation}


\begin{lemma}
Let $s_0 > - \log(y_0/4)$. Then, for any $\tau > 0$, $F_{k\to s}^{\ast}r_{\pm}U_M(1 - \chi)$ is a 
bounded operator from $L^2({\bf H}^n)$ to $H^{0,\tau}({\bf R}^{n-1}\times I)$, where $I = (s_0, \infty)$.
\end{lemma}
Proof. Note $U_M(1-\chi)$ is a bounded operator from $L^2({\bf H}^n)$ to $L^2({\bf R}^n)$. On the support of $1 - \chi(y)$, $\log y > \log y_0/4$. Therefore if $s > s_0 > - \log y_0/4$, 
\begin{equation}
\begin{split}
& F_{k\to s}^{\ast}r_{\pm} U_M(1 - \chi)f  \\
& = F^{\ast}_{k\to s}r_{\pm}(k)F_{z\to k}(1 - \chi(e^z))Tf \\
& =
(2\pi)^{-1}\int_{-\infty}^{\infty}\int_{-\infty}^{\infty}
e^{ik(s + z)}r_{\pm}(k)(1 - \chi(e^z))Tf(x,z)dkdz \\
& = \pm \int_{\bf R}\frac{1}{i(s+z)}(1 - \chi(e^z))Tf(x,z)dz.
\end{split}
\nonumber
\end{equation}
Clearly, the right-hand side is smooth with respect to $s$ with all of its derivatives in $L^2({\bf R}^{n-1}\times I_{\pm})$. \qed

\medskip
Lemma 4.8 and (\ref{eq:Chap4Sect5mathcalRplusexprseeion1}), (\ref{C4S4Restimate}) imply the following lemma.


\begin{lemma}
Let $s _ 0 > - \log y_0/4$, $\sigma \geq 0$. Then we have
\begin{equation}
\begin{split}
\mathcal R_+  - \frac{1}{\sqrt2}F_{k\to s}^{\ast}
\Big(r_+U_M(a^{(+)}_{FM})^{\ast} + r_-U_M(a^{(-)}_{FM})^{\ast}\Big) \\
\in {\bf B}(L^2({\bf H}^n) ; H^{-\sigma,\sigma/4}({\bf R}^{n-1}\times (s_0,\infty))).
\end{split}
\nonumber
\end{equation}
\end{lemma}


\section{Singularity expansion of the Radon transform}
Let us recall the following homogeneous distribution. We define for ${\rm Re}\,\alpha > - 1$
\begin{equation}
h_{\pm}^{\alpha}(s) = \left\{
\begin{split}
&|s|^{\alpha}/\Gamma(\alpha + 1), \quad  \pm s > 0, \\
& 0, \quad \pm s < 0,
\end{split}
\right.
\nonumber
\end{equation}
and, for $n = 1, 2, 3, \cdots$ and ${\rm Re}\,\alpha > -1$,
\begin{equation}
h_{\pm}^{\alpha-n}(s) =  \left(\pm \frac{d}{ds}\right)^n h_{\pm}^{\alpha}(s).
\nonumber
\end{equation}
Thus, $h_{\pm}^{\alpha}(s)$ is analytic with respect to $\alpha$.
Let $\langle \, \, , \, \rangle$ be the coupling of distributions and test functions. Then for any $\alpha, \beta \in {\bf C}$
\begin{equation}
\int_{-\infty}^{\infty}h_{\pm}^{\alpha}(s)h_{\pm}^{\beta}(1-s)ds = 
\langle h_{\pm}^{\alpha}(s) h_{\pm}^{\beta}(1-s),1\rangle = 
\frac{1}{\Gamma(\alpha + \beta + 2)}.
\label{eq;ProductHomogDis}
\end{equation}
In fact, this is true for ${\rm Re}\,\alpha, {\rm Re}\,\beta > - 1$. Let $\chi_0(s), \chi(s) \in C^{\infty}({\bf R})$ be such that
$\chi_0(s) + \chi_1(s) = 1$, $\chi_0(s) = 1 \ (s < 1/3)$, $\chi_0(s) = 0 \ (s > 2/3)$. Then we have
$$
\langle h_{+}^{\alpha}(s)h_{+}^{\beta}(1-s),1 \rangle = 
\langle h^{\alpha}_+(s),\frac{(1 - s)^{\beta}}{\Gamma(\beta+1)}\chi_0(s)\rangle + 
\langle h_+^{\beta}(1 - s),\frac{s^{\alpha}}{\Gamma(\alpha+1)}\chi_1(s)\rangle.
$$
Since $1-s>0$ on ${\rm supp}\,\chi_0$ and $s > 0$ on ${\rm supp}\,\chi_1$, 
 the left-hand side is analytic with respect to $\alpha, \beta$. Hence (\ref{eq;ProductHomogDis}) holds by analytic continuation.

\medskip
The following lemma is well-known (\cite{GeSh64} p.174, \cite{Hor}, Vol 1,  p.167).


\begin{lemma}
For $\alpha \in {\bf R}$
$$
\int_{-\infty}^{\infty}(\pm ik + 0)^{\alpha}e^{iks}dk = 2\pi
h_{\pm}^{-\alpha-1}(s).
$$
\end{lemma}

Let $\chi_{\infty}(k)$ be as in (\ref{eq:axyxik}). Since $1 - \chi_{\infty}(k) \in C_0^{\infty}({\bf R})$, from Lemma 5.1,
\begin{equation}
\frac{1}{2\pi}\int_{-\infty}^{\infty}e^{iks}k^{-j}\chi_{\infty}(k)dk 
- (-i)^jh^{j-1}_-(s) \in C^{\infty}({\bf R}), \quad j = 0, 1, 2, \cdots.
\label{eq:Chap4Sect6Distribution}
\end{equation}

Let $H^{-\sigma,\tau}_{loc}({\bf R}^{n-1}\times (s_0,\infty))$ be the 
set of functions $u$ such that, for any compact interval $I \subset (s_0,\infty)$
$$
u\big|_{{\bf R}^{n-1}\times I} \in H^{-\sigma,\tau}({\bf R}^{n-1}\times I).
$$


\begin{theorem}
Let $s_0 > - \log y_0/4$. Then 
 for any $\sigma > 0$, there is $N = N(\sigma)$ such that 
$$
\mathcal R_+ - \sum_{j=0}^N\mathcal R_j^{(+)} \in 
{\bf B}(L^2({\bf H}^n);H_{loc}^{-\sigma,\sigma/4}({\bf R}^{n-1}\times(s_0,\infty)),
$$
where 
\begin{equation}
\left(\mathcal R_+^{(j)}f\right)(s,x) = \int_0^{\infty}(s + \log y)^{j-1}_-
y^{-\frac{n-1}{2}}P_j(y)f(x,y)\chi(y)\frac{dy}{y},
\nonumber
\end{equation}
\begin{equation}
P_j(y) = \frac{(-i)^j}{\sqrt2}a_j(x,y,-i\partial_x)^{\ast}.
\nonumber
\end{equation}
\end{theorem}
Proof.
Recall from Lemma 4.9, 
$\mathcal R_+f$ is given, up to a smoothening operator, by
\begin{equation}
\frac{1}{\sqrt2}F^{\ast}_{k\to s}\left(F^{\ast}_{z\to k}\left\{(a_T^{(+)})^{\ast} + 
(a_T^{(-)})^{\ast} \right\}\right)Tf.
\label{C4S5formula1}
\end{equation}
Let $M \geq \sigma/4$, and put
$$
a^{(M,\pm)}(x,y,\xi,k) = a^{(\pm)}(x,y,\xi,k) - \chi_{\infty}(k)r_{\pm}(k)\sum_{j=0}^{M}\rho\Big(\frac{\langle\xi\rangle^2}{\epsilon_j\langle k\rangle}\Big)k^{-j}a_{j}^{(\pm)}(x,z,\xi,k).
$$
Denote by $\mathcal R_M$ the operator given by (\ref{C4S5formula1}) with $a^{(\pm)}_T$ replaced by $a^{(M,\pm)}_T$. Letting $a^{(M)} = a^{(M,+)} + a^{(M,-)}$, consider
\begin{equation}
\begin{split}
& \partial_s^p(I - \Delta_x)^{-\ell}\mathcal R_Mf \\
& = \frac{1}{\sqrt2(2\pi)^{n/2}}\int e^{i(x-x')\cdot\xi}e^{-ik(s+z')}\frac{(-ik)^p}{\langle\xi\rangle^{2\ell}}
\overline{a^{(M)}(x',z',\xi,k)}Tf(x',z')dx'dz'd\xi dk.
\end{split}
\nonumber
\end{equation}
By construction of $a_T(x,z,\xi,k)$, $\langle k\rangle \geq \langle\xi\rangle^2/\epsilon_{M+1}$ on ${\rm supp}\,\overline{a^{(M)}}$, and
$$
\Big|\partial_{x'}^{\alpha}\partial_{z'}^m\partial_{\xi}^{\beta}\partial_k^{\gamma}
\left\{(-ik)^p\langle\xi\rangle^{-2\ell}\overline{a^{(M)}(x',z',\xi,k)}\right\}\Big|
\leq C_{\alpha\beta\gamma\delta}\langle\xi\rangle^{2(M-\ell-|\beta|)}\langle k\rangle^{p-M-\gamma}.
$$
The right-hand side is bounded if $p \leq M \leq \ell$, which implies by the $L^2$-boundedness theorem for $\Psi$DO that 
$$
\mathcal R_M \in {\bf B}(L^2({\bf H}^n);H^{-s,\tau}({\bf R}^n)), \quad for \quad s\geq 2\tau, \quad \tau \leq M.
$$
In particular, $\mathcal R_M \in {\bf B}(L^2({\bf H}^n);H^{-\sigma,\sigma/4}({\bf R}^n))$.

By integation by parts using $e^{ix\cdot\xi} = \langle\xi\rangle^2(1 - \Delta_{x'})
e^{ix'\cdot\xi}$, we see that the operator 
$$
\int e^{i(x-x)\cdot\xi}e^{-ik(s+z')}\Big(1 - \rho\big(\frac{\langle\xi\rangle^2}{\epsilon_jk}\big)\Big)\overline{a_{jT}(x',z',\xi,k)}Tf(x',z')dx'dz'd\xi dk
$$
is in ${\bf B}(L^2({\bf H}^n);H^{-\ell,p}({\bf R}^n))$ with $\ell \geq 2p$, 
hence in ${\bf B}(L^2({\bf H}^n);H^{-\sigma,\sigma/4}({\bf R}^n))$.

Therefore, in view of (\ref{eq:axyxik}), we see that $\mathcal R_+f$ is equal to, up to a smoothening operator in ${\bf B}(L^2({\bf H}^n);H^{-\sigma,\sigma/4}({\bf R}^{n-1}\times(s_0,\infty))$,
\begin{equation}
\begin{split}
& \frac{1}{\sqrt2(2\pi)^{n}}\int_{{\bf R}^n\times{\bf R}^n_+}
e^{i(x - x')\cdot\xi}e^{-ik(s + \log y)}y^{\frac{n-1}{2}}
\sum_{j=0}^{M-1}k^{-j}
\overline{a_j(x',y,\xi,k)}f(x',y)\frac{d\xi dkdx'dy}{y^n} \\
&= \frac{1}{\sqrt2}\sum_{j=0}^{M-1}\int_0^{\infty}g_{j}(x,y)y^{-\frac{n-1}{2}}\chi(y)\left(\frac{1}{2\pi}\int_{-\infty}^{\infty}e^{-ik(s + \log y)}k^{-j}\chi_{\infty}(k)dk\right)\frac{dy}{y},
\end{split}
\nonumber
\end{equation}
\begin{equation}
\begin{split}
g_j(x,y) & = \frac{1}{(2\pi)^{(n-1)}}\int_{{\bf R}^{2(n-1)}}
e^{i(x-x')\cdot\xi}\overline{a_j(x',y,\xi)}f(x',y){d\xi dx'} \\
& = a_j(x,y,-i\partial_x)^{\ast}f(x,y).
\end{split}
\nonumber
\end{equation}
This together with (\ref{eq:Chap4Sect6Distribution}) proves the theorem. \qed

\bigskip
Recall that $a_j(x,y,\xi)$ is defined by (\ref{eq:Chap4Sect3jthtermoftransporteq}), and is a polynomial in $\xi$ of order $2j$. Hence $a_j(x,y,-i\partial_x)$ is a differential operator of order $2j$.
 The above theorem in particular yields the following expression
\begin{equation}
\begin{split}
 & \left(\mathcal R_+^{(j)}f\right)(s,x) \\
& =
\left\{
\begin{split}
& \frac{e^{(n-1)s/2}}{\sqrt2}\chi(e^{-s})
f(x,e^{-s}), \quad (j = 0), \\
& \int_0^{e^{-s}}\frac{(s + \log y)^{j-1}}{(j-1)!}y^{-\frac{n-1}{2}}P_j(y)f(x,y)\chi(y)\frac{dy}{y},
\quad (j \geq 1),
\end{split}
\right.
\end{split}
\label{eq:Chap4Sect5Rj}
\end{equation}
where $\chi(y) \in C^{\infty}({\bf R})$ such that $\chi(y) = 1 \ (y < y_0/4)$, $\chi(y) = 0 \ (y > y_0/3)$.
This is a generalization  of Theorem 1.6.6 in the sense of singularity expansion.

\chapter{Introduction to inverse scattering}

Suppose we are given two asymptotically hyperbolic metrics which differ only on a compact set. If the associated scattering operators coincide, one can show that these two metrics coincide up to a diffeomorphism. This result can be extended to manifolds with asymptotically hyperbolic ends when two metrics coincide on one end having a regular infinity. The aim of this chapter is to explain the idea of the proof of these theorems.

\section{Local problem on ${\bf H}^n$}
Recall that in the geodesic polar coordinates centered at $(0,1)$, the metric on ${\bf H}^n$ takes the form
\begin{equation}
ds^2 = (dr)^2 + \sinh^2r\,(d\theta)^2,
\nonumber
\end{equation}
where $(d\theta)^2$ is the standard metric on $S^{n-1}$ (see formula (\ref{NormalCoordinatesmetric}) in Chap. 1).  Letting $y = 2e^{-r}$ and $x = \theta$, one can rewrite the above metric as
$$
ds^2 = \left(\frac{dy}{y}\right)^2 + \left(\frac{1}{y} - \frac{y}{4}\right)^2\,
(dx)^2, \quad y \in (0,2].
$$
Suppose this metric is perturbed so that
$$
ds^2 = \frac{(dy)^2 + (dx)^2 + A(x,y,dx,dy)}{y^2},
$$
with $A(x,y,dx,dy)$ satisfying the assumption (A-4) of Chap. 3, \S 3.
The theorem we are going to prove is as follows.


\begin{theorem}
Suppose we are given two Riemannian metrics $G^{(p)}$, $p = 1,2$, on ${\bf H}^n$ satisfying the above assumption. Suppose their scattering operators coincide. Suppose furthermore $G^{(1)}$ and $G^{(2)}$ coincide except for a compact set. Then $G^{(1)}$ and $G^{(2)}$ are isometric.
\end{theorem}

The proof is done by the following steps. Let $B_a \subset {\bf H}^n$ be a ball of radius $a$ with respect to the unperturbed metric centered at $(0,1)$  such that $G^{(1)} = G^{(2)}$ outside $B_a$. We first take a geodesic sphere $S_a = \partial B_a$, and consider the boundary value problem for the Laplace-Beltrami operators in the interior domain $B_a$. Then the associated Dirichlet-to-Neumann map (or Neumann-to-Dirichlet map) coincide. We use the boundary control method of Belishev-Kurylev to show that $G^{(1)}$ and $G^{(2)}$ are isometric in $B_a$ (see \cite{Be87} and \cite{KKL01}).


\section{Scattering operator and N-D map}


\subsection{Restriction of generalized eigenfunctions to a surface}
For $k > 0$, let $\mathcal F^{(+)}(k)$ be the generalized Fourier transformation defined by Chap. 2 (\ref{eq:Chap2Sect7Fplusminusk}). For a compact hypersurface $S$ in ${\bf H}^n$, we define
\begin{equation}
\langle f,g\rangle_S = \int_Sf(x,y)\overline{g(x,y)}dS_{x,y},
\nonumber
\end{equation}
where $dS_{x,y}$ is the measure induced on $S$.


\begin{lemma}
Let $\Omega$ be a bounded domain in ${\bf H}^n$ with smooth boundary $S = \partial\Omega$. Suppose $k^2 \neq 0$ is not a Neumann eigenvalue for $H$ in $\Omega$. If $f \in L^2(S)$ satisfies
\begin{equation}
\langle f,\partial_{\nu}\mathcal F^{(+)}(k)^{\ast}\phi\rangle_{S} = 0, \quad \forall \phi \in L^2({\bf R}^{n-1}),
\nonumber
\end{equation} 
then $f = 0$, where $\partial_{\nu} = \dfrac{\partial}{\partial\nu}$ is the normal derivative on $S$.
\end{lemma}
Proof. We first study the local regularity of the resolvent. Take $\chi \in C_0^{\infty}({\bf H}^n)$. Then by the well-known elliptic regularity theorem,
$\chi R(k^2 \pm i0)\chi \in {\bf B}(H^s;H^{s+ 2})$, $\forall s \geq0$. By taking the adjoint, we have $\chi R(k^2 \pm i0)\chi \in {\bf B}(H^{-t-2};H^{-t})$, $\forall t \geq0$. By interpolation, we then have 
\begin{equation}
\chi R(k^2 \pm i0)\chi \in {\bf B}(H^m;H^{m+2}), \quad \forall m \in {\bf R}.
\nonumber
\end{equation}

 For $f \in L^2(S)$, we define
$$
( \delta_S'f,g)= \langle f,\partial_{\nu}g\rangle_S, \quad \forall g \in C_0^{\infty}({\bf H}^n).
$$
Then ${\rm supp}\,\delta_S'f \subset S$ and $\delta_S'f \in H^{-3/2}_{comp}$ = the set of $H^{-3/2}$-functions with compact support in ${\bf H}^n$.
For $g \in \mathcal B$, due to Theorem 2.1.3, $\partial_{\nu}R(k^2 + i0)g$ restricted on $S$ is in $H^{1/2}(S)$. Then, for $f \in L^2(S)$, the mapping
\begin{equation}
\mathcal B \ni g \longrightarrow 
\langle \partial_{\nu}R(k^2 - i0)g, f \rangle_{S}
\nonumber
\end{equation}
is a bounded linear functional. 
Using the definition of $\delta_S'f$, we have
\begin{equation}
\langle f, \partial_{\nu}R(k^2 - i0)g\rangle_{S}  = (u,g), \quad \forall g \in \mathcal B,
\nonumber
\end{equation}
where $u = R(k^2 + i0)\delta_S'f \in H^{1/2}_{loc}\cap\mathcal B^{\ast}$.
Using the resolvent equation, we see that
\begin{equation}
u = R_0(k^2 + i0)\delta_S'f - R(k^2 + i0)VR_0(k^2 + i0)\delta_S'f,
\label{ResolventeqdeltaS'}
\end{equation}
where $V = H-H_0$.
Note that $R_0(k^2 + i0)\delta_S'f$ can be written as an integral over $S$
$$
R_0(k^2 + i0)\delta_S'f = \int_S\left(\partial_{\nu'}R_0(k^2 + i0)(x,y,x',y')\right)f(x,y')dS_{x',y'}.
$$
This is an analogue of the classical double layer potential (see e.g. \cite{CoKr83}). 

To understand the properties of this potential, let $S_{\delta}$, where $|\delta|$ is sufficiently small, be an equi-distant surface which lies inside $\Omega$ for positive $\delta$ and inside $\Omega^c$ for negative $\delta$. This defines two types of operators $K_{\delta}$ and  $T_{\delta}$, where
$$
K_{\delta}f = R_0(k^2 + i0)\delta_S'f\big|_{S_{\delta}},
$$
$$
T_{\delta}f = \partial_{\nu}R_0(k^2 + i0)\delta_S'f\big|_{S_{\delta}}.
$$
For $\delta \neq 0$, they are bounded operators on $L^2(S)$, where we use the fact that $S_{\delta}$ is diffeomorphic to $S$. Moreover, $K_{\delta}$ tends to $K_{\pm}$ in the strong operator topology on $L^2(S)$, when $\delta \to \pm 0$, and $K_+ - K_- = Id$. This is proven in ${\bf R}^3$ for the classes of H{\"o}lder continuous functions in Theorem 2.15 and Corollary 2.14 of \cite{CoKr83}. 
However, if we take into account that in the Riemannian normal coordinates, $x = (x_1,\cdots,x_n)$,  $d^2(x,0) = |x|^2 + O(|x|^4)$, the method of  \cite{CoKr83} can be extended to the space $L^2(S)$ and general Riemannian manifold $\mathcal M$. 

Regarding $T_{\delta}$, it is proven in Theorem 2.23,  \cite{CoKr83}, that 
$T_{\delta}$ tends to $T_{\pm}$ in the strong operator topology of bounded operators from $C^{1,\alpha}(S)$ to $C^{\alpha}(S)$, and $T_+ - T_- = 0$. Using duality arguments and the fact that $(T_{\delta})^{\ast}$ has the same structure as $T_{\delta}$, we see that $T_{\delta}$ tends to $T_{\pm}$ in the weak operator topology of boundend operators from $L^2(S)$ to $H^{-s}(S)$, where $s > (n+1)/2$, and $T_+ - T_-=0$.

Extending formula (\ref{eq:Chap2Sect7Fplusminusk}) in Chap. 2, we define 
$\mathcal F^{(+)}(k)$ onto  $H^{-3/2}_{comp}$. Then by Lemma 2.7.3, since $\mathcal G^{(+)}(k)= \mathcal F^{(+)}(k)$, the behavior of $u$ at infinity is given by
\begin{equation}
R(k^2 + i0)\delta_S'f \simeq C(k)\chi(y)y^{\frac{n-1}{2}-ik}
\mathcal F^{(+)}(k)\delta_S'f.
\label{eq:Chap5Sect2Rfdelta}
\end{equation}
However, by the assumption of the lemma
$$
(\delta_S'f,\mathcal F^{(+)}(k)^{\ast}\phi) = (\mathcal F^{(+)}(k)\delta_S'f,\phi)_{L^2({\bf R}^{n-1})} = 0, \quad \forall \phi \in L^2(S^{n-1}).
$$
This, together with (\ref{eq:Chap5Sect2Rfdelta}), implies
\begin{equation}
\lim_{R\to\infty}\frac{1}{\log R}\int_{1/R}^1\|u(y)\|^2_{L^2({\bf R}^{n-1})}\frac{dy}{y^n} = 0.
\nonumber
\end{equation}
Let us note that for any $\varphi \in C_0^{\infty}({\bf H}^n)$
\begin{equation}
\begin{split}
((H - k^2)u,\varphi) &= (u,(H-k^2)\varphi) \\
 & = \langle f, R(k^2 - i0)(H - k^2)\varphi\rangle_S \\
 & = \langle f, \varphi\rangle_S,
\end{split}
\nonumber
\end{equation}
where we have used the fact that $\varphi = R(k^2 - i0)(H - k^2)\varphi$, since $\varphi$ is compactly supported, hence satisfies the radiation condition.
We then have  $(H - k^2)u = 0$  outside and inside $S$. Arguing in the same way as in the proof of Theorem 2.2.10 given in  Subsection 2.3.2, we have $u = 0$ in $\Omega^c := {\bf H}^n\setminus\overline{\Omega}$. 
Thus $T_-f = 0$.

Consider $u_{\Omega} = u\big|_{\Omega}$. Then $(H-k^2)u_{\Omega}=0$ and $\partial_{\nu}u_{\Omega}\big|_{\Gamma}= T_+f=0$. Since $k^2$ is not a Neumann eigenvalue, $u_{\Omega} = 0$ in $\Omega$. Therefore $u = 0$ globally in ${\bf H}^n$, which implies $f = 0$. \qed

\bigskip
By the same arguments, one can prove the following lemma.


\begin{lemma}
Let $\Omega$ be a bounded domain. Suppose $k^2 \neq 0$ is not a Dirichlet eigenvalue for $H$ in $\Omega$. If $f \in L^2(\partial\Omega)$ satisfies
\begin{equation}
\langle f,\mathcal F^{(+)}(k)^{\ast}\phi\rangle_{\partial\Omega} = 0, \quad \forall \phi \in L^2({\bf R}^{n-1}),
\nonumber
\end{equation} 
then $f = 0$.
\end{lemma}

\subsection{Neumann-to-Dirichlet map}
Let $\Omega$ be a bounded domain in ${\bf H}^n$ with smooth boundary $S = \partial\Omega$, and consider the boundary value problem
\begin{equation}
\left\{
\begin{split}
& (H - k^2)u = 0 \quad {\rm in} \quad \Omega, \\
& \partial_{\nu}u = f \in H^{1/2}(S) \quad {\rm on} \quad S. 
\end{split}
\right.
\nonumber
\end{equation}
We denote the corresponding operator in $L^2(\Omega)$ with Neumann boundary condition by $H^N$, keeping the notation $H$ for the operator in ${\bf H}^n$.
If $k^2$ is not an eigenvalue of $H^N$, this problem has a unique solution $u$. The operator
\begin{equation}
\Lambda(k) : f \to u\big|_{S}
\nonumber
\end{equation}
is called the Neumann-to-Dirichlet map, or simply, N-D map. 
 We consider two operators $H_1^N$ and $H_2^N$ associated with two metrics $G^{(1)}$ and $G^{(2)}$. Let $\widehat S_j(k)$ be the S-matrix for $H_j$.


\begin{theorem}
Suppose $k^2 \neq 0$ is not an eigenvalue for both of $H_1^N$ and $H_2^N$. Let $\Lambda_j(k)$ be the N-D map for $H_j^N$, $j = 1,2$. Suppose 
$G^{(1)} = G^{(2)}$ outside $\Omega$. Then 
$\widehat S_1(k) = \widehat S_2(k)$ if and only if $\Lambda_1(k) = \Lambda_2(k)$.
\end{theorem}
Proof. Suppose $\Lambda_1(k) = \Lambda_2(k)$. Let $u_j = \mathcal F^{(+)}_j(k)^{\ast}\phi$ for $\phi \in L^2({\bf R}^{n-1})$. Let $u_{in}$ be the solution to the Neumann problem
\begin{equation}
\left\{
\begin{split}
& (H_2 - k^2)u_{in} = 0 \quad {\rm in} \quad \Omega, \\
& \partial_{\nu}u_{in} = \partial_{\nu}u_1  \quad {\rm on} \quad S. 
\end{split}
\right.
\nonumber
\end{equation}
We define a functon $u_3$ on ${\bf H}^n$ by $u_3 = u_{in}$ on $\Omega$ and $u_3 = u_1$ on $\Omega^c = {\bf H}^n\setminus\overline{\Omega}$. The trace of $u_3$ computed from outside of $S$ is $u_3\big|_S = u_1\big|_S =\Lambda_1(k)\partial_{\nu}u_1$, since $u_1$ satisfies $(H_1-k^2)u_1=0$ in ${\bf H}^n$, hence in $\Omega$. 

On the other hand, the trace computed from inside of $S$ is 
$$
u_{in}\big|_S = \Lambda_2(k)\partial_{\nu}u_{in} =
\Lambda_2(k)\partial_{\nu}u_1 =  \Lambda_1(k)\partial_{\nu}u_1.
$$
 Therefore by our assumption, $u_3$ and $\partial_{\nu}u_3$ are continuous across $S$. Hence $u_3 \in H^2_{loc}$ and satisfies $(H_2 - k^2)u_3 = 0$ on ${\bf H}^n$.

 Let $u_0 = \mathcal F^0(k)^{\ast}\phi$. Then $u_3 - u_0$ satisifies the incoming radiation condition, since so does $u_1 - u_0$. Therefore $v = u_3 - u_2 = (u_3-u_0) - (u_2-u_0)$ is the solution to the equation $(H_2 - k^2)v = 0$ satisfying the radiation condition. By Lemma 2.2.12, $v = 0$. Observing the behavior of $u_1 = u_2$ near infinity and using Theorem 2.7.9, we have $\widehat S_1(k) = \widehat S_2(k)$.

Suppose $\widehat S_1(k) = \widehat S_2(k)$. Let $u_j$ be as above, and put $w = u_1 - u_2$. Then $(H_1 - k^2)w = 0$ in $\Omega^c$. Since $\widehat S_1(k) = \widehat S_2(k)$, $w \simeq 0$ by virtue of Lemma 2.7.2. Consequently, $w = 0$ by Theorem 2.2.10. Then $u_1 = u_2$ and $\partial_{\nu}u_1 = \partial_{\nu}u_2$ on $S$, i.e.
\begin{equation}
\Lambda_1(k)\partial_{\nu}\mathcal F_1^{(+)}(k)^{\ast}\phi = \Lambda_2(k)\partial_{\nu}\mathcal F_2^{(+)}(k)^{\ast}\phi= \Lambda_2(k)\partial_{\nu}\mathcal F_1^{(+)}(k)^{\ast}\phi.
\nonumber
\end{equation}
By Lemma 2.1, $\{\partial_{\nu}\mathcal F^{(+)}_1(k)^{\ast}\phi \, ; \phi \in L^2({\bf R}^{n-1})\}$ is dense in $L^2(S)$, which proves the theorem. \qed


\section{Boundary spectral projection}
 Our inverse problem is now reduced to determining the metric from the N-D map for a bounded domain. Since the following arguments do not rely on individual nature of the metric, we consider in a general situation. Let $\Omega$ be a compact Riemannian manifold with boundary equipped with the metric $ds^2 = g_{ij}(x)dx^idx^j$. Let $\Delta_g$ be the associated Laplace-Beltrami operator, and  $\lambda_1 < \lambda_2 < \cdots$ be the Neumann eigenvalues of $- \Delta_g$. We emphasize that we do not count the multiplicities of eigenvalues here. The N-D map is defined as $\Lambda(\lambda) : f \to u\big|_{\partial\Omega}$, where 
\begin{equation}
\left\{
\begin{split}
& (- \Delta_g - \lambda)u = 0 \quad {\rm in} \quad \Omega, \\
& \partial_{\nu}u = f \in H^{1/2}(\partial\Omega)  \quad {\rm on} \quad \partial\Omega.
\end{split}
\right.
\label{eq:Vhap4Sect10DN}
\end{equation}
Here we are writing $\Lambda(\lambda)$ instead of $\Lambda(\sqrt{\lambda})$. Note that $\Lambda(\lambda)$ is analytic with respect to $\lambda \in {\bf C}\setminus\sigma(- \Delta_g)$.
 Let $\varphi_{i,1}(x), \cdots, \varphi_{i,m(i)}(x)$ be a complete orthonormal system of eigenvectors associated with $\lambda_i$.
We first note that the N-D map $\Lambda(\lambda)$ has the following formal integral kernel
\begin{equation}
\Lambda(\lambda;x,y) =  \sum_{i=1}^{\infty}\sum_{j=1}^{m(i)}
\frac{\varphi_{i,j}(x)\overline{\varphi_{i,j}(y)}}{\lambda_i - \lambda}, \quad
x, y \in \partial\Omega.
\label{eq:Chap4Sect10NDkernel}
\end{equation}
In fact, let $\widetilde f \in H^2(\Omega)$ be such that $\partial_{\nu}\widetilde f = f$ on $\partial\Omega$. Then $v = u - \widetilde f$ solves 
\begin{equation}
\left\{
\begin{split}
& (- \Delta_g - \lambda)v = (\Delta_g + \lambda)\widetilde f =: F \quad {\rm in} \quad \Omega, \\
& \partial_{\nu}v = 0 \in H^{1/2}(\partial\Omega)  \quad {\rm on} \quad \partial\Omega.
\end{split}
\right.
\nonumber
\end{equation}
Therefore, letting $(\; , \, )$ be the inner product of $L^2(\Omega)$
\begin{equation}
v = \sum_{i=1}^{\infty}\frac{1}{\lambda_i - \lambda}\sum_{j=1}^{m(i)}(F,\varphi_{i,j})\varphi_{i,j}(x).
\label{eq:Vhap4Sect1Formofv}
\end{equation}
Letting $\langle \; ,\, \rangle$ be the inner product on $L^2(\partial\Omega)$, we have by integration by parts
$$
(F,\varphi_{i,j}) = \langle f,\varphi_{i,j}\rangle + (\lambda - \lambda_i)(\widetilde f,\varphi_{i,j}),
$$
which proves (\ref{eq:Chap4Sect10NDkernel}).


\begin{definition}
The set $\{\lambda_i,\varphi_{i,j}(x)\big|_{\partial\Omega}\; ;\, j = 1, \cdots, m(i), \, i = 1, 2, \cdots\}$ is called the boundary spectral data ({\bf BSD}) of the 
Neumann problem.
\end{definition}


\begin{lemma} \label{Lemma5.3.2}
Let $\varphi_{i,1}(x), \cdots, \varphi_{i,m(i)}(x)$ be a complete orthnormal system of eigenvectors associated with $\lambda_i$ for the Neumann problem. Then $\varphi_{i,j}(x), 1 \leq j \leq m(i)$, are linearly independent in $L^2(\partial\Omega)$. For another complete orthnormal system $\psi_{i,1}(x),\cdots,\psi_{i,m(i)}(x)$, there is a unitary matrix $U$ such that
$$
\Big(\varphi_{i,1}(x), \cdots,\varphi_{i,m(i)}(x)\Big) = 
\Big(\psi_{i,1}(x),\cdots,\psi_{i,m(i)}(x)\Big)U.
$$ 
\end{lemma}
Proof. Suppose $\sum_{j=1}^{m(i)}c_j\varphi_{i,j}(x) = 0$ on $\partial\Omega$. Then $u = \sum_{j=1}^{m(i)}c_j\varphi_{i,j}(x)$ satisfies $(- \Delta_g - \lambda_i)u = 0$ in $\Omega$, and $u = \partial_{\nu} u = 0$ on $\partial\Omega$. By the uniqueness theorem for the Cauchy problem (see e.g. \cite{Mi73}, p. 373), $u = 0$ in $\Omega$, which implies $c_1 = \cdots = c_{m(i)} = 0$. The 2nd assertion is easy to prove, since $\{\varphi_{i,j}\}$ and $\{\psi_{i,j}\}$ are the orthonomal bases of an $m(i)$-dimensional space. \qed

\bigskip
Let us give an operator theoretical meaning to (\ref{eq:Chap4Sect10NDkernel}). We need the notion of spectral representation. Let $\widehat{\mathcal H} = \oplus_{i=1}^{\infty}{\bf C}^{m(i)}$. We define the (discrete) Fourier transformation $\mathcal F : L^2(\Omega) \to \widehat{\mathcal H}$ by $\mathcal F = (\mathcal F_1,\mathcal F_2,\cdots)$ where
\begin{equation}
\mathcal F_i : L^2(\Omega) \ni u \to \left((u,\varphi_{i,1}),\cdots,(u,\varphi_{i,m(i)})\right) \in {\bf C}^{m(i)}.
\label{eq:discreteFourier}
\end{equation}
$\mathcal F$ is unitary, and diagonalizes the Neumann Laplacian $- \Delta_g$ on $\Omega$ : $\mathcal F_i(- \Delta_gu) = \lambda_i\mathcal F_iu$. 
Let $P_i$ be the eigenprojection associated with the eigenvalue $\lambda_i$. Then, for $z \not\in \sigma(- \Delta_g)$, the resolvent can be written as
\begin{equation}
R_{\Omega}(z) = \sum_{i=1}^{\infty}\frac{P_i}{\lambda_i - z} = 
\sum_{i=1}^{\infty}\frac{\mathcal F_i^{\ast}\mathcal F_i}{\lambda_i - z},
\label{eq:Neumannresolvent}
\end{equation}
which converges in the sense of strong limit in $L^2(\Omega)$.

Let $\Gamma = \partial\Omega$, and $r_{\Gamma} \in {\bf B}(H^{1}(\Omega);H^{1/2}(\Gamma))$ be the trace operator  to $\Gamma$. Define $\delta_{\Gamma} \in {\bf B}(H^{-1/2}(\Gamma);H^{1}(\Omega)^{\ast})$ as its adjoint:
\begin{equation}
(\delta_{\Gamma}f,w)_{L^2(\Omega)} = (f,r_{\Gamma}w)_{L^2(\Gamma)}, \quad  f \in H^{-1/2}(\Gamma), \quad w \in H^{1}(\Omega).
\nonumber
\end{equation}
Accordingly, we write as
\begin{equation}
r_{\Gamma} = \delta_{\Gamma}^{\ast}.
\nonumber
\end{equation}
Then we have
\begin{equation}
\delta_{\Gamma} \in {\bf B}(H^{-1/2}(\Gamma);H^{1}(\Omega)^{\ast}), \quad 
\delta_{\Gamma}^{\ast} \in {\bf B}(H^1(\Omega);H^{1/2}(\Gamma)).
\label{C5S3;deltagamma}
\end{equation}
Then,
\begin{equation}
\Lambda(z) = \delta_{\Gamma}^{\ast}R_{\Omega}(z)\delta_{\Gamma}.
\label{eq:NDand Resolevent}
\end{equation}

Let us prove this formula. We first show that the right-hand side is well-defined. Since $R_{\Omega}(z) \in {\bf B}(L^2(\Omega);H^2(\Omega))$, we have $R_{\Omega}(z) \in {\bf B}(H^{2}(\Omega)^{\ast};L^2(\Omega))$. By an interpolation, we then have $R_{\Omega}(z) \in {\bf B}(H^{1}(\Omega)^{\ast};H^1(\Omega))$. Using (\ref{C5S3;deltagamma}), we see that $\delta_{\Gamma}^{\ast}R_{\Omega}(z)\delta_{\Gamma} \in 
{\bf B}(H^{-1/2}(\Omega);H^{1/2}(\Omega))$.

For $f \in H^{1/2}(\Gamma)$, take $\widetilde f \in H^{3/2}(\Omega)$ such that $\partial_{\nu}\widetilde f = f$ on $\Gamma$. Let $v = R_{\Omega}(z)(\Delta_g + z)\widetilde f$, and put $u = v + \widetilde f$. Then $(- \Delta_g - z)u = 0$ in $\Omega$, and $\partial_{\nu}u = f$ on $\Gamma$. Take $h \in L^2(\Omega)$. Then, by integration by parts,
\begin{equation}
\begin{split}
(P_i(\Delta_g + z)\widetilde f,h)_{L^2(\Omega)} &= (z - \lambda_i)(\widetilde f,P_ih)_{L^2(\Omega)} + (f,r_{\Gamma}P_ih)_{L^2(\Gamma)} \\
&= - (\lambda_i -z)(P_i\widetilde f,h)_{L^2(\Omega)} + 
(P_i\delta_{\Gamma}f,h)_{L^2(\Omega)}.
\end{split}
\nonumber
\end{equation}
This yields
\begin{equation}
P_iu = P_i\widetilde f + R_{\Omega}(z)P_i(\Delta_g + z)\widetilde f = \frac{P_i\delta_{\Gamma}f}{\lambda_i - z}.
\nonumber
\end{equation}
By (\ref{eq:Neumannresolvent}), this implies $u = R_{\Omega}(z)\delta_{\Gamma}f$. By taking the trace to $\Gamma$, we get (\ref{eq:NDand Resolevent}).

By Lemma \ref{Lemma5.3.2}, the operator $\delta_{\Gamma}^{\ast}P_i\delta_{\Gamma}$, whose integral kernel is $\sum_{j=1}^{m(i)}\varphi_{i,j}(x)\overline{\varphi_{i,j}(y)}$ restricted to $\Gamma$, is independent of the choice of the eigenvectors. Let us call the set 
\begin{equation}
\Big\{\big(\lambda_i, \sum_{j=1}^{m(i)}\varphi_{i,j}(x)\overline{\varphi_{i,j}(y)}\Big|_{x\in {\Gamma}, y \in \Gamma}\big)\Big\}_{i=1}^{\infty}, 
\label{eq:BSPdiscrete}
\end{equation} 
{\it boundary spectral projection} ({\bf BSP}).
This is what we actually use in the BC method. BSP is the set of pairs of poles and residues of the N-D map.  We then have the following lemma.


\begin{lemma}
Suppose we are given two metrics on $\Omega$. Then their BSP's coincide if and only if their N-D maps coincide for all $\lambda$ outside the spectrum.
\end{lemma}

In the next chapter, we shall explain how to reconstruct the metric from BSP. 


\section{Inverse problems for hyperbolic ends}


\subsection{Exterior boundary value problem}
 Before entering into the inverse scattering for manifolds with hyperbolic ends, we need to discuss the spectral theory for the exterior boundary value problem.
Let $\Omega$ be a bounded domain in ${\bf H}^n$ with smooth boundary and $\Omega^c := {\bf H}^n\setminus\overline{\Omega}$. Let $H^{N,c}$ be $H$ defined in $\Omega^c$ with Neumann boundary condition. Namely $D(H^{N,c}) = \{u \in H^2(\Omega^c) ; \partial_{\nu}u\big|_{\partial\Omega^c} = 0\}$ and $H^{N,c}u = Hu$ for $u \in D(H^{N,c})$.
Then $H^{N,c}$ is self-adjoint. Let $R^c(z) = (H^{N,c} - z)^{-1}$. The theory developed for $H$ in Chap. 2 can be extended to $H^{N,c}$ without any essential change. In fact, let $u(z) = R^c(z)f$, $f \in L^2(\Omega^c)$, for $z \in {\bf C}\setminus{\bf R}$, and take $\chi \in C^{\infty}({\bf H}^n)$ such that $\chi = 1$ near infinity, and $\chi = 0$ on a bounded open set containing $\overline{\Omega}$.
Then $v(z) = \chi R^c(z)f$ satisfies
\begin{equation}
(H - z)v = [H,\chi]R^c(z)f + \chi f, \quad {\rm in} \quad {\bf H}^n,
\nonumber
\end{equation}
where we use that $\omega := \hbox{supp} [H,\chi] \subset\subset \Omega^c$.  
Let us show that 
\begin{equation}
\|u(z)\|_{\mathcal B^{\ast}} \leq C(\|u\|_{L^2(\Omega_1)} + \|f\|_{\mathcal B}).\label{C5S4Inequalityproof}
\end{equation}
where $\Omega_1$ is a compact set such that $\omega \subset \Omega_1 \subset \Omega^c$. In fact, by elliptic regularity, 
$$
\|u\|_{H^1(\omega)} \leq C(\|u\|_{L^2(\Omega_1)} + \|f\|_{L^2(\Omega_1)}).
$$
The inequality (\ref{C5S4Inequalityproof}) then follows from this and 
(\ref{R(z)UniformEstimate}) in Chap. 2. 

Having inequality (\ref{C5S4Inequalityproof}) in our disposal, we can prove, using the same arguments as for the whole ${\bf H}^n$, Lemma 2.2.13 for $R^c(z)$.


\begin{theorem}
(1) $\ \sigma_{e}(H^{N,c}) = [0,\infty), \ \sigma_p(H^{N,c})\cap(0\
,\infty) = 0$. \\
\noindent
(2) For any $\lambda > 0$, $lim_{\epsilon\to 0}R^c(\lambda \pm i\epsilon) =:R^c(\lambda \pm i0)$ exists in ${\mathcal B}^{\ast}$ in the weak $\ast$-sense.\\
\noindent
(3) For any compact interval $I \subset (0,\infty)$, there exists a constant $C > 0$ such that
\begin{equation}
\|R^c(\lambda \pm i0)f\|_{\mathcal B^{\ast}} \leq C\|f\|_{\mathcal B}, \quad \forall \lambda \in I.
\nonumber
\end{equation}
(4) For any $f, g \in \mathcal B$, $(0,\infty) \ni \lambda \to (R^c(\lambda \pm i0)f,g)$ is continuous. \\
\noindent
(5) For $\lambda > 0$, $R^c(\lambda \pm i0)f$ is a unique solution to the 
equation
\begin{equation}
\left\{
\begin{split}
& (H - \lambda)u = f \in \mathcal B \quad {\rm in} \quad \Omega^c,\\
& \partial_{\nu}u = 0 \quad {\rm on} \quad \partial\Omega
\end{split}
\right.
\nonumber
\end{equation}
satisfying the outgoting (for $+$) or incoming (for $-$) radiation condition.
\end{theorem}

The following lemma can now be proved easily by using Theorem 4.1.


\begin{lemma}
Let $\lambda > 0$ and $f \in H^{1/2}(\partial\Omega)$. Then there exists a unique solution $ u_{\pm} \in \mathcal B^{\ast}$ to the exterior boundary value problem
\begin{equation}
\left\{
\begin{split}
& (H - \lambda)u = 0  \quad {\rm in} \quad \Omega^c,\\
& \partial_{\nu}u = f \quad {\rm on} \quad \partial\Omega
\end{split}
\right.
\nonumber
\end{equation}
satisfying the outgoing or incoming radiation condition.
\end{lemma}
Using the solutions $u_{\pm}$ as above, we define the N-D map by $\Lambda^{(\pm)}(\lambda)f = u_{\pm}\big|_{\partial\Omega}$ in addition to 
 $\Lambda(z)$  for $z \in {\bf C}\setminus\sigma(H^{N,c})$. Note that $\Lambda^{(\pm)}(\lambda)$ is the boundary value of $\Lambda(z)$ as $z \to \lambda \pm i0$. Therefore, $\Lambda^{(\pm)}(\lambda)$ defined for $\lambda > 0$ has a unique analytic continuation to ${\bf C}\setminus{\sigma}(H^{N,c})$.


\subsection{Inverse scattering at regular ends}
Let $\mathcal M$ be a manifold satisfying the assumptions (A.1) $\sim$ (A.4) in Chap. 3, \S 3 with ends of number $N \geq 2$. 
We assume that at least {\it one of the ends has a regular infinity}. Let ${\mathcal M}_1$ be such an end. Namely, in the notation of Chap. 3, \S 2, $\mathcal M_1$ is diffeomorphic to $M_1 \times (0,1)$, in other words, $\mathcal M_1$ is asymptotically equal to a funnel.
Let $\Gamma \subset \mathcal M$ be a compact submanifold of codimension 1 such that $\mathcal M$ splits into 3 parts $\Omega$, $\Omega^c$, $\partial\Omega = \partial\Omega^c = \Gamma$ in the following way :
\begin{equation}
\mathcal M = \Omega\cup\Gamma\cup\Omega^c, \quad
\Omega\cap\Gamma = \Omega^c\cap\Gamma = \emptyset,
\nonumber
\end{equation}
where $\overline{\Omega}$ and $\overline{\Omega^c}$ are assumed to be submanifolds of $\mathcal M$ with boundary $\Gamma$ inheriting the Riemannian metric of $\mathcal M$. Assume also that $\Omega$ is non-compact and has infinity common to $\mathcal M_1$, and has no other infinity, i.e. $\Omega = M_1\times(0,a)$, $0<a<1$. Note that when $N \geq 2$, $\Omega^c$ is also non-compact having a finite number of ends which are either regular or cusps. (The case when $N=1$, which  is equivalent to $\overline{\Omega^c}$ being compact, brings about the inverse boundary spectral problem discussed in \S 3.)

 Let 
$H^N$ be $- \Delta_g - (n-1)^2/4$ in $\Omega$ with Neumann boundary condition, and $H^{N,c}$ be the one on $\Omega^c$.
Then Theorem 4.1 and Lemma 4.2 also hold for $H^N$ and $H^{N,c}$. Note that if all the  ends except for $\mathcal M_1$ have cusps, there may be embedded eigenvalues in the essential spectrum of $H^{N,c}$. However, they are discrete with possible accumulation points only at 0 and infinity with rapidly decreasing eigenvectors.

We generalize Lemma 2.1 to the present case.
Let $\mathcal F^{(\pm)}(k) = (\mathcal F^{(\pm)}_1(k),\cdots,\mathcal F^{(\pm)}_N(k))$ be the generalized Fourier transformation in $\mathcal M$ constructed in Chap. 3, \S 2, and ${\bf h}_{\infty}$ be defined by (\ref{eq:Chap3Sect2hinfty}) in Chap. 3.


\begin{lemma}
Suppose $0 \neq k^2 \not\in \sigma_p(H^{N,c})$. If $f \in L^2(\Gamma)$ satisfies
$$
\langle f,\partial_{\nu}\mathcal F^{(+)}(k)^{\ast}\phi\rangle_{\Gamma} = 0, \quad 
\forall \phi = (\phi_1,0,\cdots,0) \in {\bf h}_{\infty},
$$
then $f = 0$.
\end{lemma}

Proof.  Since (\ref{eq:Chap5Sect2Rfdelta}) holds
in $\mathcal M_1$, arguing in the same way as in Lemma 2.1, we have $u = 0$ in $\Omega$. Consider $u^c = u\big|_{\Omega^c}$. Then we have
$(H - k^2)u^c = 0$ in $\Omega^c$, and similarly to the proof of Lemma 2.1 $\partial_{\nu}u^c = 0$ on $\Gamma$. Since $u^c$ also satisfies the radiation condition, and $k^2 \not\in \sigma_p(H^{N,c})$, we have $u^c = 0$ in $\Omega^c$. This proves the lemma. \qed

\medskip
Recall that $H^{N,c}$ has two parts of spectral representations: the generalized Fourier transform, which we denote by $\mathcal F_{c}^{(+)}$ here, corresponding to the absolutely continuous spectrum for $H^{N,c}$, and the discrete Fourier transform, denoted by $\mathcal F_{p}^c$, corresponding to the point specrum for $H^{N,c}$ defined in the same way as in \S 3. 


\begin{lemma} The N-D map $\Lambda^c(z)$ corresponding to $H^{N,c}$, which is determined  for $z \in {\bf C}\setminus{\bf R}$, is of the form
\begin{equation}
\Lambda^c(z) = \int_0^{\infty}\frac{\delta_{\Gamma}^{\ast}\mathcal F_c^{(+)}(k)^{\ast}\mathcal F_c^{(+)}(k)\delta_{\Gamma}}{k^2 - z}dk + 
\sum_i\frac{\delta_{\Gamma}^{\ast}P_i^c\delta_{\Gamma}}{\lambda_i - z},
\label{eq:Chap5Sec4NDK1}
\end{equation}
where the sum over $i$ may be finite or infinite.
\end{lemma}
Proof. We proceed as in the proof of (\ref{eq:NDand Resolevent}). Take $f \in C^{\infty}(\Gamma)$ and $\widetilde f \in C_0^{\infty}(\overline{\Omega^c})$ such that $\partial_{\nu}\widetilde f = f$ on $\Gamma$.
Let $v$ solve the boundary value problem
\begin{equation}
\left\{
\begin{split}
& (H - z)v = (-H + z)\widetilde f =:F \quad {\rm in } \quad \Omega^c, \\
& \partial_{\nu}v = 0 \quad {\rm on} \quad \Gamma.
\end{split}
\right.
\nonumber
\end{equation}
Then $v$ is represented by eigenvectors $\varphi_{i,j}$ and the generalized Fourier transform $\mathcal F_c^{(+)}$:
$$
v = \int_0^{\infty}\frac{\mathcal F_c^{(+)}(k)^{\ast}\mathcal F_c^{(+)}(k)F}{k^2 - z}dk + \sum_i\frac{\sum_j(F,\varphi_{i,j})\varphi_{i,j}}{\lambda_i - z}.
$$
Take $\phi \in {\bf h}_{\infty}^c$ (see Chap. 3, (\ref{eq:Chap3Sect2hinfty}), where $j$ varies from 2 to $N$). Then we have by integration by parts
\begin{equation}
\begin{split}
(\mathcal F_c^{(+)}(k)F,\phi)_{{\bf h}_{\infty}^c} &= ((-H + z)\widetilde f,\mathcal F_c^{(+)}(k)^{\ast}\phi)_{L^2(\Omega^c)} \\
&= (f,\mathcal F_c^{(+)}(k)^{\ast}\phi)_{L^2(\Gamma)} + 
(z - k^2)(\widetilde f,\mathcal F_c^{(+)}(k^2)^{\ast}\phi)_{L^2(\Omega^c)} \\
&= (\mathcal F_c^{(+)}(k)\delta_{\Gamma}f,\phi)_{{\bf h}_{\infty}^c} + 
(z - k^2)(\mathcal F_c^{(+)}(k^2)\widetilde f,\phi)_{{\bf h}_{\infty}^c}
.
\end{split}
\nonumber
\end{equation}
This implies
$$
\mathcal F_c^{(+)}(k)F = \mathcal F_c^{(+)}(k)\delta_{\Gamma}f + 
(z - k^2)\mathcal F_c^{(+)}(k)\widetilde f.
$$
The term from the point spectrum is dealt with similarly, and the lemma follows from a direct computation. \qed

\medskip
Let us call the set
\begin{equation}
\Big\{\delta_{\Gamma}^{\ast}\mathcal F_c^{(+)}(k)^{\ast}\mathcal F_c^{(+)}(k)\delta_{\Gamma} ; k > 0\Big\} \cup\Big\{\Big(\lambda_i, \delta_{\Gamma}^{\ast}P_i^c\delta_{\Gamma}\Big)\, ; \, i\Big\}
\end{equation}
 the boundary spectral projection ({\bf BSP}) for $H^{N,c}$. By (\ref{eq:Chap5Sec4NDK1}), we have
\begin{equation}
\Lambda^c(z) = \delta_{\Gamma}^{\ast}(H^{N,c} - z)^{-1}\delta_{\Gamma}.
\label{eq:Chap5Sect4NDandResolvent}
\end{equation}


\begin{lemma}
Knowing the N-D map $\Lambda^{(+)}_c(k^2)$ for all $k$ such that $k^2 \not\in \sigma_p(H^{N,c})$ is equivalent to knowing BSP for  $H^{N,c}$.
\end{lemma}
Proof.
$\Lambda^{(+)}_c(k^2)$ has a unique analytic continuation $\Lambda^c(z)$ for $z \in {\bf C}\setminus{\bf R}$, which determines $\Lambda^{(-)}_c(k^2)$ for real $k^2 \not\in \sigma_p(H^{N,c})$. By 
(\ref{eq:Chap5Sect4NDandResolvent}) and Lemma 3.3.11, we have
\begin{equation}
\Lambda^{(+)}_c(k^2) - \Lambda^{(-)}_c(k^2) = \frac{\pi i}{k}\delta_{\Gamma}^{\ast}
\mathcal F_c^{(+)}(k)^{\ast}\mathcal F_c^{(\pm)}(k)\delta_{\Gamma}.
\nonumber
\end{equation}
Therefore we recover 
$\mathcal F_c^{(+)}(k)^{\ast}\mathcal F_c^{(+)}(k)$ for $k^2 \not\in \sigma_p(H^{N,c})$ from $\Lambda^{(+)}_c(k^2)$. By (\ref{eq:Chap5Sec4NDK1}), we also recover $\lambda_i \in \sigma_p(H^{N,c})$ and $\delta_{\Gamma}^{\ast}P_i^c\delta_{\Gamma}$ from the poles and residues of $\Lambda^c(z)$.
 The converse direction is seen by (\ref{eq:Chap5Sec4NDK1}). \qed

\bigskip
Since $\mathcal M$ has $N$-ends, the S-matrix for $\mathcal M$ is an $N\times N$-matrix:
\begin{equation}
\widehat S(k) = \Big(\widehat S_{ij}(k)\Big)_{1\leq i, j \leq N}.
\nonumber
\end{equation}

 Let $\mathcal M^{(j)}, (j=1,2),$ be  manifolds satisfying the assumptions (A.1) $\sim$ (A.4) in Chap. 3, \S 3. Assume that $\mathcal M^{(1)}_1$ and $\mathcal M^{(2)}_1$ are isometric, therefore, $\mathcal M_1^{(1)} = \mathcal M^{(2)}_1 = M_1\times (0,1)$, $M_1$ being a compact manifold of dimension $n-1$.  Letting $\Omega = M_1\times(0,a)$, we construct $\Omega^{c}_j$ and $H^{N,c}_j$ as above.


\begin{theorem}
Suppose $0 \neq k^2 \not\in \sigma_p(H^{N,c}_1)\cup\sigma_p(H^{N,c}_2)$. Let $\Lambda_j^{(+)}(k^2)$ be the N-D map for $H^{N,c}_j$. Then $\widehat S_{11}^{(1)}(k) = \widehat S_{11}^{(2)}(k)$ if and only if $\Lambda_1^{(+)}(k^2) = \Lambda_2^{(+)}(k^2)$.
\end{theorem}
The proof is the same as Theorem 2.3.

\bigskip
We now pass to the boundary control method (BC-method) to show that BSP determines the manifold uniquely. The BC-method works for general Riemannian manifold wih boundary,  if we know the N-D map for all $k$ for the associated Laplace operator. 
The BC-method was first applied to compact manifolds (\cite{BeKu92}), and was extended to non-compact manifolds (see e.g. \cite{KKL04}, \cite{IKL10}). 

\medskip
Let us formulate the inverse problem on non-compact Riemannian manifolds. Let $\mathcal N_1$ and $\mathcal N_2$ be Riemannian manifolds (not necessarily compact) with boundary with metric inherited form the Riemannian metric induced from $\mathcal N_j$.  
We say that $\mathcal N_1$ and $\mathcal N_2$ have  common parts $\Gamma_1 \subset \partial{\mathcal N}_1$ and $\Gamma_2 \subset \partial\mathcal N_2$  if there exists an isometry $\Phi : \Gamma_1 \to \Gamma_2$. Let $\Lambda_j(z)$ be the N-D map for the Laplace operator on $\mathcal N_j$. Then we define
\begin{equation}
\Lambda_1(z)\Big|_{\Gamma_1} = \Lambda_2(z)\Big|_{\Gamma_2} 
\Longleftrightarrow
\Phi\circ\Lambda_1(z)\Big|_{\Gamma_1} = \Lambda_2(z)\Big|_{\Gamma_2}\circ\Phi.
\label{C5S4:NDequiv}
\end{equation}
Here $\Lambda_j(z)\big|_{\Gamma_j}$ is defined by
$$
\Lambda_j(z)\big|_{\Gamma_j}f = \Lambda_j(z)f\big|_{\Gamma_j}, \quad f \in L^2(\Gamma_j).
$$ 
One can then show that (with some additional assumptions) if $\mathcal N_1$ and $\mathcal N_2$ have common parts $\Gamma_1$ and $\Gamma_2$, and (\ref{C5S4:NDequiv}) holds for all $z \not\in {\bf R}$, then $\mathcal N_1$ and $\mathcal N_2$ are isometric. In Chapter 6, we shall give the proof of this theorem (Theorem \ref{main}) for asymptotically hyperbolic manifolds $\Omega_1^c, \Omega_2^c$ under consideration. 
Modulus this {\it theorem}, we have thus proven the following result.


\begin{theorem}
Let $\mathcal M$ be a manifold satisfying the assumptions (A.1) $\sim$ (A.4) in Chap. 3, \S 3. We assume that  one of the ends has a regular infinity, and denote it by $\mathcal M_1$. Suppose we are given two metrics $G^{(j)}$, $j = 1, 2$, on $\mathcal M$ satisfying (A-3) in Chapt. 3, \S 3. 
Assume that $G^{(1)} = G^{(2)}$ on $\mathcal M_1$. If $\widehat S_{11}(k) = \widehat S_{11}(k)$ for all $k > 0$, then $G^{(1)}$ and $G^{(2)}$ are isometric on $\mathcal M$.
\end{theorem}

\medskip
We can actually prove a stronger version of Theorem 4.7, which is valid for two manifolds whose structure, in particular the number of ends, are not known a-priori.


\begin{theorem}
Let $\mathcal M^{(j)}$, $j = 1,2$, be manifolds satisfying the assumptions (A.1) $\sim$ (A.4) in Chap. 3, \S 3 endowed with metric $G^{(j)}$, $j = 1,2$. We assume that for both of $\mathcal M^{(1)}$and $\mathcal M^{(2)}$ one of the ends has a regular infinity, and denote them  by $\mathcal M_1^{(j)}$, $j = 1,2$.  
Assume that $\mathcal M_1^{(1)}$ and $\mathcal M_1^{(2)}$ are isometric, and $\widehat S_{11}(k) = \widehat S_{11}(k)$ for all $k > 0$. Then $\mathcal M^{(1)}$ and $\mathcal M^{(2)}$ are isometric.
\end{theorem}


\subsection{References of inverse scattering on asymptotically hyperbolic manifolds}
Melrose's theory of scattering metric studies the spectral properties of the Laplace-Beltrami operator on manifolds whose ends have the metric of the following type 
$$
ds^2 = \frac{h(x,y,dx,dy)}{y^2}.
$$
Each end is assumed to be isomorphic to $X\times(0,1)$ and $g_0(x,y,dx,dy)$ admits an asymptotic expansion of the form
$$
h(x,y,dx,dy) = (dy)^2 + h_0(x,dx) + y\,h_1(x,dx,dy) + y^2 h_2(x,dx,dy) + \cdots,
$$
$h_0(x,dx)$ being a Riemannian metric on the boundary at infinity, $X$.
Mazzeo and Melrose \cite{MaMe87} developed a pseudo-differential calculus to deal with these manifolds, and proved the existence of analytic continuation of resolvent of the associated Laplace-Beltrami operator into the region ${\bf C}\setminus\{\frac{1}{2}(n - {\bf N}_0)\}$, ${\bf N}_0 = {\bf N}\cup\{0\}$. Borthwick \cite{Bo01} studied the case of variable curvature at the boundary at infinity.
Guillarmou \cite{Gulm05} showed that the resolvent had in general essential singularities at 
$\{\frac{1}{2}(n - {\bf N}_0)\}$. Joshi and S{\'a} Barreto \cite{JoSaBa00} proved that the scattering matrix determined the asymptotic expansion of the metric $h(x,y,dx,dy)$ at infinity. S{\'a} Barreto \cite{SaBa05} proved that the scattering matrix for all energies determined the whole manifold.

Resonance is also an important subject in the inverse scattering theory, and many works are devoted to it. They are summarized in \cite{GuZw97} or in the book of Borthwick \cite{Bo07}.

For the spectral theory of symmetric spaces of higher rank, there ia a work \cite{MaVa07}.

Inverse scattering problem or inverse boundary value problem from a fixed energy is not yet solved completely for the case of the metric. However, in 2-dimensions the inverse boundary value problem is completely solved by Nachman \cite{Na95}, Lassas-Uhlmann \cite{LaUh01}, Astala-Paivarinta \cite{AsPa06} and Astala-Lassas-Paivarinta \cite{AsLaPa05}. For higher dimensions, there is a developed theory for isotropic metrics, 
see the review article of \cite{Uh92}. Morever a method was developed to study anisotropic metrics from a known conformal class. See e.g. \cite{DSKSU09}.

There is a link between the hyperbolic manifolds and the inverse boundary value problems in the Euclidean space. See \cite{Is04a}, \cite{Is04b}, \cite{Is04c}, \cite{Is07a}, \cite{Is07b}. In \cite{IINSU07} an application to the numerical computation is given.


\chapter{Boundary control method}


\section{Brief introduction to the boundary control method}

\subsection{Wave equation and Gel'fand inverse problem}
Let $\mathcal N$ be an $n$-dimensional complete connected Riemannian manifold with boundary $\partial\mathcal N$. We shall consider an IBVP (initial-boundary value problem) for the wave equation
\begin{equation}
\partial_t^2u = \Delta_gu \quad {\rm on} \quad \mathcal N\times(0,\infty),
\nonumber
\end{equation}
where $\Delta_g$ is the Laplace-Beltrami operator. In local coordinates 
\begin{equation}
\Delta_g = g^{-1/2}\partial_i(g^{ij}g^{1/2}\partial_j), \quad
g = \det\,(g_{ij}).
\nonumber
\end{equation}
 We impose the initial condition
\begin{equation}
u\big|_{t=0} = \partial_tu\big|_{t=0} = 0,
\nonumber
\end{equation}
and the boundary condition
\begin{equation}
\partial_{\nu}u\big|_{\partial \mathcal N\times(0,\infty)} = f \in C^{\infty}_0(\partial \mathcal N\times(0,\infty)).
\nonumber
\end{equation}
Here $\nu$ is the outer unit normal to $\partial \mathcal N$. 
Let $u^f(x,t)$ be the solution to the above IBVP. We measure $u^f$ on $\partial \mathcal N\times(0,\infty)$, and call 
\begin{equation}
\Lambda^h : f \to u^f\big|_{\partial \mathcal N\times(0,\infty)}
\label{eq:ChapSect1HypNDmap}
\end{equation}
a {\it hyperbolic Neumann-to-Dirichlet map}.
The basic question we address is the following one.

\medskip
\noindent
{\bf Question} Assume we know $\Lambda^h$. Can we determine $(\mathcal N,g)$, i.e. the manifold $\mathcal N$ and the metric $g$?

\medskip
This is the {\it Gel'fand inverse problem} 
(stated in a slightly different form, \cite{Gel57}). Note that $\Lambda^h$ is an operator defined on $\partial \mathcal N\times(0,\infty)$. Starting from the knowledge on $\partial\mathcal N \times(0,\infty)$, the first issue is the topology of $\mathcal N$, and the second issue is the Riemannian structure.

The answer to the above question is affirmative when $\mathcal N$ is compact, and also for non-compact $\mathcal N$ with some additional geometric assumption. To fix the idea, in this chapter, $\mathcal N$ means either any compact connected Riemannian manifold with boundary, or when dealing with the non-compact case, the manifold $\Omega^c$ discussed in Chap. 5, \S 4. However, the arguments given below also work for non-compact manifolds possesing the spectral representation as in the case of $\Omega^c$. Note that in both cases $\partial \mathcal N$ is compact.


\subsection{Spectral formulation} 
Let us begin with the compact manifold case.
Consider the Neumann Laplacian $H^N$: 
\begin{equation}
H^Nu = - \Delta_gu, \quad u \in H^2(\mathcal N), \quad
\partial_{\nu}u\big|_{\partial \mathcal N} = 0.
\nonumber
\end{equation}
The spectrum of $H^N$ consists of real numbers
\begin{equation}
0 = \lambda_1 < \lambda_2 \leq \lambda_3 \leq \cdots \to \infty.
\nonumber
\end{equation}
 Let $\varphi_k$ be the associated eigenvectors
\begin{equation}
- \Delta_g\varphi_k = \lambda_k\varphi_k, \quad \partial_{\nu}\varphi_k\big|_{\partial \mathcal N} = 0.
\nonumber
\end{equation}
Without loss of generality we can assume $\varphi_k$ to be real-valued.
The set $\{\varphi_k\}_{k=1}^{\infty}$ can be made to
form an orthonormal basis in $L^2(\mathcal N)$ and orthogonal 
basis in $H^1(\mathcal N)$, where 
the inner products of $L^2(\mathcal N)$ and $H^1(\mathcal N)$ are defined by
\begin{equation}
(f,g)_{L^2(\mathcal N)} = \int_{\mathcal N}f(x)\overline{g(x)}dV_g, \quad
dV_g = g^{1/2}dx^1\cdots dx^n,
\nonumber
\end{equation}
\begin{equation}
(f,g)_{H^1(\mathcal N)} = \int_{\mathcal N}g^{ij}\,\partial_if\,\overline{\partial_jg}\,dV_g + 
(f,g)_{L^2}.
\nonumber
\end{equation}
We call $\big\{\big(\lambda_k, \varphi_k\big|_{\partial \mathcal N}\big)\big\}_{k=1}^{\infty}$ 
the {\it boundary spectral data} (BSD). The original 
Gel'fand inverse problem is equivalent to:

\medskip
\noindent
{\bf Question} Given BSD, can we determine $(\mathcal N,g)$?

\medskip
The relation of BSD to the hyperbolic Neumann-to-Dirichlet map is represented by the following (formal) formula:
\begin{equation}
\left(\Lambda^hf\right)(x,t) = \int_{\partial \mathcal N}\int_{{\bf R}_+}
G(x,y,t-s)f(y,s)dS_yds.
\nonumber
\end{equation}
\begin{equation}
G(x,y,t) = \sum_{k=1}^{\infty}\frac{\sin(\sqrt{\lambda_k}t)}{\sqrt{\lambda_k}}
\varphi_k(x)\varphi_k(y)
\big|_{\partial \mathcal N\times\partial \mathcal N}.
\label{eq:NeumannBoundaryOp}
\end{equation}

One can also deal with the Dirichlet Laplacian, i.e. 
\begin{equation}
H^Du = - \Delta_gu,  \quad 
u \in H^2(\mathcal N)\cap H^1_0(\mathcal N).
\nonumber
\end{equation}
Let $0 < \mu_1 < \mu_2 \leq \mu_3 \leq \cdots \to \infty$ be the Dirichlet eigenvalues, and $\psi_k$ the associated eigenvectors. Considering IBVP 
\begin{equation}
\left\{
\begin{split}
& \partial_t^2w = \Delta_gw, \\
& w\big|_{\partial \mathcal N\times{\bf R}+} = 
f \in C_0^{\infty}(\partial \mathcal N\times{\bf R}_+),\\
& w\big|_{t=0} = \partial_tw\big|_{t=0} = 0,
\end{split}
\right.
\nonumber
\end{equation}
we define the hyperbolic Dirichlet-to-Neumann map by
\begin{equation}
R^hf  : f \to \partial_{\nu}w^f\big|_{\partial \mathcal N\times{\bf R}_+}.
\nonumber
\end{equation}
The integral kernel of $R^h$ is formally written as
\begin{equation}
R^h(x,y,t) = \sum_{k=1}^{\infty}\frac{\sin(\sqrt{\mu_k}t)}{\sqrt{\mu_k}}
\partial_{\nu}\psi_k(x)\partial_{\nu}\psi_k(y)
\big|_{\partial \mathcal N\times\partial \mathcal N}.
\nonumber
\end{equation}

The method we are going to talk about is called 
the {\it Boundary Control} (BC) method, whose history goes back to the famous results by M. G. Krein, in the mid-fifties, on the 
$1-$dimensional inverse scattering theory (\cite{Kr51a}, \cite{Kr51b}). Compared with the fundamental methods by Gel'fand-Levitan and Marchenko, the method of Krein is distinguished by the systematic use of the finite propagation speed for the wave equation. However, the ideas based upon the domain of influence, etc. coming from this finite velocity are "disguised" in the work of Krein due to their formulation in the frequency domain (or the stationary equation), where they turn out to be conditions on analyticity of the corresponding Fourier transform of the solution. This principal hyperbolic nature of Krein's method was revealed by Blagovestchenskii who was working in the time-domain (or the time-dependnet equation) using the finite velocity of the wave propagation and ideas of controllability in the filled domain to derive a Volterra-type equation for unknown functions (\cite{Bla71a}). These ideas have become crucial for the extension  of the method to multidimensions pioneered by Belishev \cite{Be87}, see also \cite{KKL01}. One more important ingredient of the BC-method, namely, the possibility to evaluate the inner product of waves sent into $\mathcal N$ 
from $\partial \mathcal N$ also goes back to the 1-dimensional case to the work of Blagovestchenskii
\cite{Bla71b}. See \cite{BeBla92} for the multidimensional case.

The BC method has the following features.

\noindent
{\it (1) BC method is hyperbolic.}

Since the propagation speed of wave motion is finite, and singularities of waves are related with geodesics, this implies the close connection of BC method with geometry.

\noindent
{\it (2) BC method is not  perturbative}.

We do not assume that the given metric is close to some standard one. In this sense, the BC method does not have the character of perturbation theory.


\subsection{Outline of the procedure}  
The crucial tool of the BC-method is the Kuratowski space of boundary distance functions $R(\mathcal N)$ to be defined in \S 5, 
and the reconstruction of the manifold $\mathcal N$ is done by the following 3 steps :
 
\begin{itemize}
\item 
In \S 8, we show that BSP determines $R(\mathcal N)$.

\item
 In \S 5, we show that $R(\mathcal N)$ is topologically isomorphic to
 $\mathcal N$. 

\item
In \S 7, we show that $R(\mathcal N)$ determines the 
Riemannian metric of $\mathcal N$.
\end{itemize}

This is an effective interplay of linear partial differential equations and geometry. The main ingredients of the 1st step are Blagovestchenskii's idenitity, which represents the solution of the initial boundary value problem (IBVP) of the wave equation by BSD, and Tataru's uniqueness theorem, which guarantees the conrollablity of IBVP. The 2nd step is of the character of general topology. The 3rd step is purely from differential geometry, in which the coordinate system of $\mathcal N$ is constructed by $R(\mathcal N)$ and the metric tensor is computed. The analytic and geometric preliminaries are done in \S 2, \S4, and in \S 5, \S 6, respectively.


\section{Blagovestchenskii idenitity}
Given a solution $u^f$ of the wave equation
\begin{equation}
\left\{
\begin{split}
& \partial_t^2u = \Delta_gu, \\
& \partial_{\nu}u\big|_{\partial \mathcal N\times{\bf R}_+} = f, \\
& u\big|_{t=0} = \partial_{t}u\Big|_{t=0} = 0,
\end{split}
\right.
\label{C6S2IBVPEquation}
\end{equation}
we expand it by eigenvectors to get
\begin{equation}
u^f(x,t) = \sum_ku_k^f(t)\varphi_k(x),
\quad
u^f_k(t) = \int_{\mathcal N}u^f(y,t)\varphi_k(y)dV_g.
\nonumber
\end{equation}
Then we have
\begin{equation}
\begin{split}
\frac{d^2}{dt^2}u_k^f(t) &= \int_{\mathcal N}\Delta_gu^f(y,t)\varphi_k(y)dV_g \\
&= \int_{\partial \mathcal N}\left[\partial_{\nu}u^f\varphi_k - u^f\partial_{\nu}\varphi_k\right]dS_g + \int_{\mathcal N}u^f\Delta_g\varphi_kdV_g \\
&= \int_{\partial \mathcal N}f(y,t)\varphi_k(y,t)dS_g - \lambda_k\int_{\mathcal N}u^f(y,t)\varphi_k(y)dV_g.
\end{split}
\nonumber
\end{equation}
We have thus derived
\begin{equation}
\frac{d^2}{dt^2}u_k^f(t) + \lambda_ku_k^f(t) = 
\int_{\partial \mathcal N}f(y,t)\varphi_k(y)dS_g,
\nonumber
\end{equation}
and, due to the initial condition in IBVP,
\begin{equation}
u_k^f(0) = \frac{d}{dt}u^f_k(0) = 0.
\nonumber
\end{equation}
Solving this differential equation, we obtain  {\it Blagovestchenskii idenitity}\begin{equation}
u_k^f(t) = \int_0^tds\int_{\partial \mathcal N}dS_g\frac{\sin(\sqrt{\lambda_k}(t-s))}{\sqrt{\lambda_k}}f(y,s)\varphi_k(y).
\label{eq:ufkrepresentedbyBSD}
\end{equation}
This formula shows that $u_k^f(t)$ is represented by $\lambda_k$ and $\varphi_k\big|_{\partial \mathcal N}$, i.e. BSD.


\begin{lemma} The following  holds:
\begin{equation}
(u^f(t),u^h(s)) = \sum_ku_k^f(t)\overline{u_k^h(s)},
\label{eq:Chap6Sect1Blagoidentity}
\end{equation}
i.e. BSP determines the inner product
$(u^f(t),u^h(s))_{L^2(\mathcal N)}, \ \forall t, s \in {\bf R}, \ \forall f, h \in C_0^{\infty}(\partial \mathcal N \times{\bf R}_+)$.
\end{lemma}
Proof. This follows from (\ref{eq:ufkrepresentedbyBSD}) and the Parseval formula. \qed

\medskip
Lemma 2.1 is the first corner-stone of BC method. We let
\begin{equation}
S(t,\lambda) = \frac{\sin(\sqrt{\lambda}t)}{\sqrt{\lambda}}, \quad
\widetilde S(t,s,\lambda) = S(t,\lambda)S(s,\lambda),
\label{eq:Ch6Sc1Atlambda}
\end{equation}
and use the notation in Chap. 5, \S 3 to rewrite the right-hand side of (\ref{eq:Chap6Sect1Blagoidentity}) as
\begin{equation}
\sum_i \int_0^t\int_0^s dt' ds'\widetilde S(t-t',s-s',\lambda_i)\left(\delta_{\Gamma}^{\ast}P_i\delta_{\Gamma}f(t'),h(s')\right).
\label{eq:Chap6Sect1BSpProjInt}
\end{equation}
This implies the following corollary.


\begin{cor}
The inner product $(u^f(t),u^h(s))$ is written only by BSP. 
\end{cor}

This is also true when $- \Delta_g$ has the continuous spectrum. 
Recall that in \S 4 of Chap. 5, the Laplace-Beltrami operator on $\Omega^c$  admits the spectral representation $\mathcal F_c^{(+)}$. In this case, to modify the formula (\ref{eq:Chap6Sect1Blagoidentity}), we have only to add the integral of $\mathcal F_c^{(+)}(k)^{\ast}\mathcal F_c^{(+)}(k)$ to the right-hand side of (\ref{eq:Chap6Sect1BSpProjInt}): 
\begin{equation}
\begin{split}
&\int_0^{\infty}dk \int_0^t\int_0^sdt'ds'\, \widetilde S(t-t's-s',k^2)\left(\delta_{\Gamma}^{\ast}\mathcal F_c^{(+)}(k)^{\ast}\mathcal F_c^{(+)}(k)\delta_{\Gamma}f(t'),h(s')\right)\\
& + \sum_i\int_0^t\int_0^sdt'ds'\widetilde S(t-t',s-s',\lambda_i)\left(\delta_{\Gamma}^{\ast}P_i^c\delta_{\Gamma}f(t'),h(s')\right).
\end{split}
\label{eq:Chap6Sect1BSpProjIntContSp}
\end{equation}
 Again $(u^f(t),u^h(s))$ is written only by BSP.

 Let us remark that in \cite{KKL01}, p. 214, Lemma 4.9, it is shown that one can construct BSD from BSP up to a multiplication factor if $\mathcal N$ is compact.


\section{Geodesics} 
Let us recall some 
basic notions from Riemannian geometry.
The distance of two points $x, y$ of a Riemannian manifold $\mathcal N$, denoted by $d(x,y)$, is defined by the infimum 
of length of piecewise smooth curves joining $x$ and $y$. This makes $\mathcal N$ a metric space. If $\mathcal N$ is complete 
in this metric, it is said to be {\it metrically complete}. 
When $\partial \mathcal N = \emptyset$, by the theorem of Hopf-Rinow (see e.g. \cite{GaHuLa80}, 
pp. 94, 95),  it is equivalent to that $\mathcal N$ is {\it geodesically complete}, i.e. any solution of the equation 
of geodesics can be extended onto the whole line ${\bf R}$. In this case, again by the theorem of Hopf-Rinow, 
any two points in $\mathcal N$ can be joined by the minimal geodesic (i.e. the shortest curve). 

In local coordinates, the equation of geodesics is written as
\begin{equation}
\frac{d^2x^k}{dt^2} + \Gamma_{ij}^k(x)\frac{dx^i}{dt}\frac{dx^j}{dt} = 0.
\label{C6S2:geodesics}
\end{equation}
Let $x(t,y,v)$ be the solution of (\ref{C6S2:geodesics}) satisfying 
$$
x(0,y,v) = y, \quad \partial_tx(0,y,v) = v \in T_y(\mathcal N), 
$$
where $\partial_t = d/dt$ and $T_y(\mathcal N)$ is the tangent space at $y$. Let $|v|_g$ be the length of $v \in T_y(\mathcal N)$. Then the map defined by
\begin{equation}
\exp_y(v) : T_y(\mathcal N) \ni v \to x(1,y,v) = x(|v|_g,y,\hat v) \in \mathcal N, \quad
\hat v = v/|v|_g
\label{eq:Exponentialmap}
\end{equation}
is called the {\it exponential map}. Using this exponential map, we define the 
{\it Riemannian normal coordinates} centered at $y$ in the following way. Let $B_{y,\rho} = \{v \in T_y(\mathcal N); |v|_g < \rho\}$. Then for $\rho$ sufficiently small, the map
\begin{equation}
\exp_y : B_{y,\rho} \ni v \to \exp_y(v) \in \exp(B_{y,\rho}) \subset \mathcal N
\nonumber
\end{equation}
is a diffeomorphism. Hence $v = (v_1,\cdots,v_n)$ can be used as local coordinates on $\exp_y(B_{y,\rho})$. Note that (\ref{eq:Exponentialmap}) implies that, when dealing with geodesics $x(t,y,v)$, we can always parametrize them so that $|v|_g = 1$. This parametrization is called the arclength parametrization and will be always used in this chapter.

Almost all of the notions from Riemannian geometry can be extended to the manifold with boundary by obvious changes. 
The problem of the existence of the shortest curves, however, is delicate. Think of, for example, non-convex 
domains in ${\bf R}^n$. However,  for any $x, y \in \mathcal N$, there exists
a shortest curve, which is $C^1$-smooth. See e.g. \cite{AlAl81}.
Moreover, the segments of this curve lying inside $\mathcal N$ are (minimal) geodesics in $\mathcal N$,
while the segments of this curve lying on $\partial \mathcal N$ are minimal geodesics on $\partial \mathcal N$.

The following lemma is easy to prove. Let ${ d}(x,y)$ be the distance between $x$ and $y$ with respect to the Riemannian metric $g$, and for a subset $S \subset \mathcal N$, ${d}(x,S) = \inf \{{d}(x,y)\, ; \, y \in S\}$. 


\begin{lemma}\label{Lemma6.3.1}
For any $x \in \mathcal N$, there exists $z \in \partial \mathcal N$ such that $d(x,z) = d(x,\partial \mathcal N)$. Moreover $x = \gamma_z(s)$, where $\gamma_z$ is the geodesic starting from $z$ with initial direction the inner unit normal to $\partial \mathcal N$, and $s = d(x,z)$.
\end{lemma}


\section{Controllabilty and observability}
Two notions in the title of this section are fundamental concepts in control 
theory. They are related to properties of solution operators of dynamical problems.


\subsection{Domains of influence} 
For any set $A \subset \mathcal N$ and  $t_0 > 0$, we define the {\it domain of influence} of $A$ (at time $t_0$) by
\begin{equation}
\mathcal N(A,t_0) = \{x \in \mathcal N \, ; \, d(x,A) \leq t_0\}.
\nonumber
\end{equation}
We introduce the forward, $D_+(A,t_0)$, backward, $D_-(A,t_0)$, and double cones, $D(A,t_0)$,  of dependence by
\begin{equation}
\begin{split}
 D_{\pm}(A,t_0)& = \{(x,t) \, ; \, x \in \mathcal N(A,t_0\mp t), \ 0 \leq \pm t  \leq t_0\}, \\
D(A,t_0)& =  D_+(A,t_0)\cup  D_-(A,t_0).
\end{split}
\nonumber
\end{equation}


\begin{lemma}
Take $t_0 > 0$ and a bounded open set $A \subset \mathcal N$ arbitrarily. 
Let $u$ be a solution to the initial boundary value problem 
\begin{equation}
\left\{
\begin{split}
&\partial_t^2u = \Delta_gu, \quad {\rm in} \quad \mathcal N\times{\bf R}, \\
&u = \partial_tu = 0, \quad {\rm on} \quad \mathcal N(A,t_0) \quad {\rm at} \quad t = 0,\\
&\partial_{\nu}u = 0, \quad {\rm on} \quad  
D(A,t_0)\cap({\partial \mathcal N}\times{\bf R}).
\end{split}
\right.
\label{eq:IBVPsect2}
\end{equation}
Then $u = 0$ in $ D(A,t_0)$. 
\end{lemma}

Proof. We prove this lemma in the case when $\mathcal N$ is a domain in ${\bf R}^n$ and, due to symmetry $t \to -t$, for $t > 0$. The general case can be proved in the same way by taking local coordinates. 

\medskip
First we  recall the well-known energy inequality. Note the identity:
\begin{equation}
\begin{split}
& \frac{1}{2}\partial_t\left((\partial_tv)^2 + g^{ij}\partial_iv\partial_jv\right) - \frac{1}{\sqrt{g}}\partial_i\left(\sqrt{g}g^{ij}\partial_jv\partial_tv\right)\\
& = \Big(\partial_t^2v - \frac{1}{\sqrt{g}}\partial_i(\sqrt{g}g^{ij}\partial_jv)\Big)\partial_tv,
\end{split}
\label{C5S4Identity}
\end{equation}
where $\partial_t = \partial/\partial t$, $\partial_i = \partial/\partial x^i$. Take a time interval I = [0,T], a family of connected open sets $A(t) \subset {\bf R}^n$ $(t \in I)$ and consider a domain $D(T) \subset {\bf R}^n\times{\bf R}^1$ such that
$$
D(T) = \{(x,t)\, ; \, t \in I, \ x \in A(t)\}.
$$
Then $\partial D(t)$ consists of 3 parts:
$$
\partial D(T) = A(T) \cup A(0) \cup S,
$$
where the lateral boundary $S$ consists of 2 parts:
\begin{equation}
S= S_{\partial}\cup S_r,
\quad
S_{\partial} = \overline{D(T)}\cap \big(\partial{\mathcal N}\times[0,T]\big),
\quad
S_r = S\setminus S_{\partial}.
\label{LateralSdelSr}
\end{equation}
Assume that  $S_r$ is piecewise smooth and its unit normal $n = (n_1,\cdots,n_n, n_t)$, with respect to the Euclidean metric, has the property
\begin{equation}
n_t \geq \left(g^{ij}n_in_j\right)^{1/2}, \quad {\rm on} \quad S_r.
\label{ntbiggerthannx}
\end{equation}
Suppose that a real-valued function $v = v(x,t)$ satisfies the wave equation 
\begin{equation}
\left\{
\begin{split}
& \partial_t^2v - \frac{1}{\sqrt{g}}\partial_i(\sqrt{g}g^{ij}\partial_jv) = 0, \quad {\rm in} \quad D(T), \\
& \partial_{\nu}v = 0, \quad {\rm on} \quad S_r.
\end{split}
\right.
\label{C5S4Equatuonforv}
\end{equation}
Mutilplying (\ref{C5S4Identity}) by $\sqrt{g}$ and integrating on $D(T)$, we have
\begin{equation}
\begin{split}
& \frac{1}{2}\Big[\int_{A(t)}\left((\partial_tv)^2+g^{ij}\partial_iv\partial_jv\right)\sqrt{g}\,dx\Big]_{t=0}^{t=T} \\
& = - \frac{1}{2}\int_{S_r}n_t\left((\partial_tv)^2 + g^{ij}\partial_iv\partial_jv\right)\sqrt{g}\, dS + \int_{S_r}n_ig^{ij}\partial_jv\partial_tv\sqrt{g}\,dS,
\end{split}
\end{equation}
where the integral over $S_{\partial}$ disappears due to the boundary condition in (\ref{C5S4Equatuonforv}) and $n_t=0$ on $S_{\partial}$.
The right-hand side is non-positive by (\ref{ntbiggerthannx}), estimate 
$$
|n_ig^{ij}\partial_jv\partial_tv| \leq |\partial_tv|(g^{ij}n_in_j)^{1/2}
(g^{ij}\partial_i\partial_j)^{1/2} \leq n_t|\partial_tv|(g^{ij}\partial_i\partial_j)^{1/2}
$$
and the Cauchy-Schwarz inequality. 
This implies
$$
\int_{A(T)}\left((\partial_tv)^2 + g^{ij}\partial_iv\partial_jv\right)\sqrt{g}\,dx \leq
\int_{A(0)}\left((\partial_tv)^2 + g^{ij}\partial_iv\partial_jv\right)\sqrt{g}\,dx.
$$
This holds with $T$ relplaced by $\tau \in (0,T)$. Therefore, if $v\big|_{t=0}=
\partial_tv\big|_{t=0}=0$ on $A(0)$, we have $\nabla v\big|_{t=\tau}= 0, \partial_tv\big|_{t=\tau} = 0$ on $A(\tau)$, hence $v=0$ on $D(T)$.

\bigskip
We turn to the proof of Lemma 4.1. 
In the following, $C_0$ and $C$ denote constants independent of small $\epsilon > 0$ and $j = (j_1,\cdots,j_n) \in {\bf Z}^n$.  

 For a small $\epsilon > 0$ , we take lattice points 
 $P(j,\epsilon) = (j_1\epsilon/C_0,\cdots,j_n\epsilon/C_0)$,
 where $C_0$ is a large constant. We extend $\big(g^{\alpha\beta}(x)\big)$ smoothly outside $\mathcal N$, and put
\begin{equation}
G^{(j,\epsilon)} = \Big(g^{\alpha\beta}(P(j,\epsilon))\Big) + \epsilon C_0I_n,
\label{MetricGj}
\end{equation}
$I_n$ being the $n\times n$ identity matrix.
Letting $d_{j,\epsilon}(\cdot,\cdot)$ be the distance defined by the  Riemannian metric $G_{j,\epsilon} = (G^{(j,\epsilon)})^{-1}$, we put
$$
B(j,\epsilon) = \{x \in {\bf R}^n \, ; \, d_{j,\epsilon}(x,P(j,\epsilon)) \leq \epsilon\}.
$$
We also let
$$
\mathcal N_{\epsilon}(A,t_0) = \{x \in \mathcal N(A,t_0)\, ; \, d(x,\partial\mathcal N(A,t_0)) > \epsilon\},
$$
where $d(\cdot,\cdot)$ is the distance defined by the Riemannian metric $\big(g_{\alpha\beta}(x)\big)$. Then $\mathcal N_{\epsilon}(A,t_0) \subset \mathcal N(A,t_0)$ and 
$\mathcal N_{\epsilon}(A,t_0) \to \mathcal N(A,t_0)$ as $\epsilon \to 0$.

We now consider a finite set
$$
J(\epsilon) = \{j \, ; \, P(j,\epsilon) \in \mathcal N_{\epsilon}(A,t_0)\},
$$
and for $j \in J(\epsilon)$, we put
$$
D(j,\epsilon) = \Big\{(x,t) \, ; \, x \in\mathcal N, \  d_{j,\epsilon}(x,P(j,\epsilon)) \leq \epsilon + t, \ 0 \leq  t \leq \epsilon/C_0\Big\}.
$$
As above, its lateral boundary consists of 2 parts like (\ref{LateralSdelSr}). We show that the condition (\ref{ntbiggerthannx}) is satisfied on $S_r$. 

For the sake of simplicity, we assume that $P(j,\epsilon) = 0$. The lateral boundary is defined as the zeros of 
$$
\varphi(x,t) = \epsilon + t - (G_{j,\epsilon}x,x)^{1/2}.
$$
Since the Euclidean normal unit of the lateral boundary is given by
$(\nabla_x\varphi,\partial_t\varphi)/(|\nabla_x\varphi|^2 + (\partial_t\varphi)^2)^{1/2}$, we have only to show that for any $x$ on the lateral boundary
\begin{equation}
1 \geq (G(x)^{-1}\nabla_x\varphi,\nabla_x\varphi), \quad G(x) = \big(g_{\alpha\beta}(x)\big).
\label{Lemma 6.4.1Whattoptove}
\end{equation}
Let $G_0 = G(P(j,\epsilon))$. Then $G_{j,\epsilon} = (G_0^{-1} + \epsilon C_0)^{-1}$ and $G(x) = G_0 + O(\epsilon)$.
Since $\nabla_x\varphi = - G_{j,\epsilon}x/(G_{j,\epsilon}x,x)^{1/2}$, we have
\begin{equation}
\begin{split}
(G(x)^{-1}\nabla_x\varphi,\nabla_x\varphi) = 
\frac{(G_0^{-1}G_{j,\epsilon}x,G_{j,\epsilon}x)}{(G_{j,\epsilon}x,x)} + 
\frac{(O(\epsilon)G_{j,\epsilon}x,G_{j,\epsilon}x)}{(G_{j,\epsilon}x,x)}. 
\end{split}
\label{Lemma6.4.1Formula1}
\end{equation}
In the right-hand side, $G_0$ and $G_{j,\epsilon}$ are positive definite, and $O(\epsilon)$ is symmetric. Noting that 
$$
\sqrt{G_{j,\epsilon}}\,O(\epsilon)\sqrt{G_{j,\epsilon}} \leq \epsilon C_1
$$
for some constant $C_1 > 0$, we see that
\begin{equation}
\frac{(O(\epsilon)G_{j,\epsilon}x,G_{j,\epsilon}x)}{(G_{j,\epsilon}x,x)} \leq \epsilon C_1.
\label{Lemma6.4.12ndterm}
\end{equation}
To compute the 1st term of the right-hand side of (\ref{Lemma6.4.1Formula1}), 
we first note $G_0^{-1}G_{j,\epsilon} = (1 + \epsilon C_0G_0)^{-1}$. Letting $\lambda_1$ be the smallest eignvalue of $G_0$, we have
$$
(1 + \epsilon C_0G_0)^{-1} \leq (1 + \epsilon C_0\lambda_1)^{-1}.
$$
 Then, letting $y = \sqrt{G_{j,\epsilon}}x$, and noting that $G_0$ and $G_{j,\epsilon}$ commute, we can estimate the 1st term as
\begin{equation}
\frac{((1 + \epsilon C_0G_0)^{-1}y,y)}{(y,y)}  
\leq \frac{1}{1+ \epsilon C_0\lambda_1}.
\label{Lemma6.4.11stterm}
\end{equation}
 In view of (\ref{Lemma6.4.12ndterm}) and (\ref{Lemma6.4.11stterm}), taking $C_0$ large enough, we see that (\ref{Lemma 6.4.1Whattoptove}) is satisfied.

\medskip
We now put
\begin{equation}
D_1(\epsilon) = {\mathop\cup_{j\in J(\epsilon)}}D(j,\epsilon),
\label{D1epsilon}
\end{equation}
and apply  the energy inequality to have
\begin{equation}
u = 0, \quad {\rm in} \quad D_1(\epsilon).
\label{C6S4u=0inD1}
\end{equation}
Let $D(A,t_0,\tau)$ be the section of $D(A,t_0)$ at time $t = \tau$. We also let $\Sigma_1^{high}(\tau)$ be the boundary of the section of $D_1(\epsilon)$ at time $t = \tau$, and $\Sigma_1^{low}(\tau)$ be the surface such that 
\begin{equation}
\left\{
\begin{split}
& \Sigma_1^{low}(\tau)\  \sqsupset \ \Sigma_1^{high}(\tau), \\ 
& d(\Sigma_1^{low}(\tau),\Sigma_1^{high}(\tau)) = 2\epsilon + C\epsilon \tau,
\end{split}
\right.
\label{C6S4Sigam1define}
\end{equation}
where for 2 compact surfaces $S_1$ and $S_2$, $S_1 \sqsupset S_2$ (or $S_2 \sqsubset S_1$) means that $S_2$ is contained in the bounded domain with boundary $S_1$, and
 where $C$ is chosen large enough.

 The meaning of (\ref{C6S4Sigam1define}) is as follows. At time $t = 0$, we take the surface $\Sigma_1^{high}(0)$ and $\Sigma_1^{low}(0)$ inside and outside of $\partial D(A,t_0)$ with distance $\epsilon$. We then develop them by  speeds higher or lower than that of waves. At time $t$, the distance between $\Sigma_1^{high}(t)$ and $\Sigma_1^{low}(t)$ will increase at most by  $C\epsilon t$.

 Let $\Sigma(\tau)$ be the boundary of $D(A,t_0,\tau)$.
Then we have
\begin{equation}
\Sigma_1^{high}(t) \sqsubset \Sigma(t) \sqsubset \Sigma_1^{low}(t), \quad 0 \leq t \leq \epsilon/C_0.
\label{C6S4SigmamaxandSigmamin}
\end{equation}

The next step starts from the time $t = \epsilon/C_0$ instead of $t = 0$, and $D_1(\epsilon)$ instead of $D(A,t_0)$. One can then construct $D_2(\epsilon)$ and $\Sigma^{high}_2(t)$ as above for $\epsilon/C_0 \leq t \leq 2\epsilon/C_0$.
Then by the energy inequality
\begin{equation}
u = 0, \quad {\rm in} \quad D_2(\epsilon), 
\label{C6S4D=0inD2}
\end{equation}
for the time interval $\epsilon/C_0 \leq t \leq 2\epsilon/C_0$.
The surface $\Sigma_2^{low}(\tau)$ is defined by
\begin{equation}
\left\{
\begin{split}
& \Sigma_2^{low}(\tau) \sqsupset \Sigma_2^{high}(\tau), \\ 
& d(\Sigma_2^{low}(\tau),\Sigma_2^{high}(\tau)) = 2\epsilon + \frac{C}{C_0}\epsilon^2 
+ C\epsilon(\tau - \frac{\epsilon}{C_0}).
\end{split}
\right.
\label{C6S4Sigma2define}
\end{equation}
We continue this procedure. In the $k$-th step, we obtain
\begin{equation}
u = 0, \quad {\rm in} \quad D_k(\epsilon), 
\label{C6S4D=0inDk}
\end{equation}
in the time interval $(k-1)\epsilon/C_0 \leq t \leq k\epsilon/C_0$, and 
\begin{equation}
\left\{
\begin{split}
& \Sigma_k^{low}(\tau) \sqsupset \Sigma_k^{high}(\tau), \\ 
& d(\Sigma_k^{low}(\tau),\Sigma_k^{high}(\tau)) = 2\epsilon + \frac{C}{C_0}(k-1)\epsilon^2 
+ C\epsilon(\tau - \frac{(k-1)}{C_0}\epsilon).
\end{split}
\right.
\label{C6S4Sigmakdefine}
\end{equation}

Now, with a given time $t_0 > 0$ and a large number $N$, we take $\epsilon$ as $N\epsilon/C_0 = t_0$. We put 
$$
D_N = \cup_{j=1}^ND_k(\epsilon).
$$
Then, by the above consideration,
$$
u = 0, \quad {\rm in} \quad D_N.
$$
By our construction, $D_N \subset D(A,t_0)$.
When $N \to \infty$, $D(N)$ tends to $D(A,t_0)$. In fact, by (\ref{C6S4Sigmakdefine}) and $N\epsilon = t_0$,
$$
d(\Sigma_k^{low}(\tau),\Sigma_k^{high}(\tau)) \leq (2 + Ct_0)\epsilon \to 0.
$$
This proves Lemma 4.1. \qed

\medskip
In the proof this lemma, we follow the basic steps of Theorem IV 2.2 of \cite{Lad73}, making them more precise by taking into the account the variable velocity of the wave propagation. 
 
Using Lemma 4.1, we can describe the support of the waves generated by the Neumann boundary sources, namely the solution $u^f$ of the IBVP,
\begin{equation}
\left\{
\begin{split}
& \partial_t^2 u = \Delta_gu, \quad {\rm in} \quad \mathcal N\times (0,\infty) \\
& u\big|_{t=0} = \partial_tu\big|_{t=0} = 0, \quad {\rm on} \quad
\mathcal N, \\
& \partial_{\nu}u\big|_{\partial \mathcal N \times (0,\infty)} = f \in 
C_0^{\infty}(\partial \mathcal N\times (0,\infty)).
\end{split}
\right.
\label{C6S4IBVP}
\end{equation}
To this end, for any subset $A \subset \mathcal N$, we introduce the forward, $C_+(A)$, backward, $C_-(A)$, and the double, $C(A)$, cones of influence 
\begin{equation}
\begin{split}
& C_{\pm}(A) = \{(x,t)\, ; \, d(x,A) \leq \pm t, \ \pm t > 0\},\\
& C(A) = C_+(A)\cup C_-(A).
\end{split}
\label{Coneofinflience}
\end{equation}


\begin{cor} \label{DomainInfCor}
 Let $u^f$ be the solution to IBVP (\ref{C6S4IBVP}). Let, in addition, 
 ${\rm supp}\,f \subset S\times (0,\infty)$, where $S \subset \partial\mathcal N$ is open. Then 
\begin{equation}
{\rm supp}\,u^f \subset C_{+}(S).
\nonumber
\end{equation}
\end{cor}
Proof. Let $t_0 > 0$ and $(y_0,t_0) \not\in C_+(S)$, then for small $r > 0$,
$$
\{(x,t)\, ; \ x \in \mathcal N, \  d(x,S)\leq t,\ 0 \leq t \leq t_0\}\cap D(B_r(y_0),t_0)   = \emptyset,
$$
$B_r(y_0)$ being the ball of radius $r > 0$ centered at $y_0$.
Applying Lemma 6.4.1, we have $u^f(y_0,t_0) = 0$. 
To complete the proof, just note that for $t < 0$, $u^f(x,t) = 0$. 
\qed

\subsection{Unique continuation and controllabilty}
Next we describe the properties of $u^f(\cdot,t)$ in $\mathcal N(S,t)$, when ${\rm supp}\,f \subset 
S\times(0,\infty)$. We start with the following global uniqueness theorem which is essentially due to Tataru (\cite{Ta95}). 


\begin{theorem}\label{Tataruglobal}
Let $u \in H^1_{loc}(\mathcal N\times(-t_0,t_0))$ satisfies 
\begin{equation}
\left\{
\begin{split}
& \partial_t^2u = \Delta_gu \quad {\rm in} \quad \mathcal N \times (-t_0,t_0), \\
& \partial_{\nu}u\big|_{\partial \mathcal N\times(-t_0,t_0)} = 0, \quad
u\big|_{S\times(-t_0,t_0)} = 0.
\end{split}
\right.
\label{TataruEq}
\end{equation}
Then $u = 0$ in $D(S,t_0)$.
\end{theorem}

 For a measurabe subset $D \subset \mathcal N$ and  $v \in L^2(D)$, we define $v = 0$ on $\mathcal N\setminus D$ and regard $L^2(D)$ as a closed subspace of $L^2(\mathcal N)$.


\begin{cor}\label{TataruCor}
 Assume $v$ satisfies 
\begin{equation}
\left\{
\begin{split}
& \partial_t^2v = \Delta_gv \quad {\rm in} \quad \mathcal N \times {\bf R}, \\
&  v\big|_{t = t_0} = 0, 
\quad \partial_tv\big|_{t=t_0}=:\psi \in L^2(\mathcal N(S,t_0)),\\
& \partial_{\nu}v\big|_{\partial \mathcal N\times(0, t_0)} = 0, \quad
v\big|_{S\times(0,t_0)} = 0.
\end{split}
\right.
\label{C6S3:adjointprob}
\end{equation}
Then $\partial_tv\big|_{t = t_0} = 0.$
\end{cor}

Proof. We extend $v(t)$ on the time interval $(t_0,2t_0)$  by $v(t) = - v(2t_0-t)$, and put $w(t) = v(t-t_0)$. Then $w$ satisfies the conditions in Theorem \ref{Tataruglobal}. \qed

\medskip
 Corollary \ref{TataruCor} shows the usefulness of the notion of the 
{\it observability operator},
\begin{equation}
\mathcal O^S_{t_0} : L^2(\mathcal N(S,t_0)) \ni \psi \to v^{\psi}\big|_{S\times(0,t_0)} \in L^2(S\times(0,t_0)),
\nonumber
\end{equation}
where $v^{\psi}$ is the  solution to (\ref{C6S3:adjointprob}).
Note that $v^{\psi}\big|_{\partial \mathcal N\times{\bf R}} \in C({\bf R};H^{1/2}(\partial \mathcal N))$, and
\begin{equation}
\|\mathcal O^S_{t_0}\psi\|_{L^2(S\times(0,t_0))} \leq C\|\psi\|_{L^2(\mathcal N)},
\label{C6S3:Ot0norm}
\end{equation}
where $C = C_{t_0}$ is a constant.

 Corollary \ref{TataruCor} is equivalent to the following fact, called the {\it observabilitry}.


\begin{cor} \label{eq:KernelOt0=0}
For any open set $S \subset \partial\mathcal N$ and $t_0>0$, we have
$$
{\rm Ker}\,\mathcal O^S_{t_0} = \{0\}.
$$
\end{cor}

\medskip
We consider now, the map $\mathcal C^S_{t_0}$ defined by 
\begin{equation}
\mathcal C^S_{t_0} : L^2(S\times(0,t_0)) \ni f \to 
u^f\big|_{t=t_0} \in L^2(\mathcal N(S,t_0)).
\nonumber
\end{equation}
 
 The crucial fact about $\mathcal C^S_{t_0}$ is the following theorem.


\begin{theorem}\label{ApprCont} \  $\overline{{\rm Ran}\,(\mathcal C_{t_0})} = L^2(\mathcal N(S,t_0))$.
\end{theorem}

Proof. Due to Corollary \ref{eq:KernelOt0=0}, it is sufficient to show
\begin{equation}
\mathcal C^S_{t_0} = - \big(\mathcal O^S_{t_0}\big)^{\ast},
\label{C6S3:Wt0norm}
\end{equation}
i.e.
\begin{equation}
(C_{t_0}^Sf,\psi)_{L^2(\mathcal N(S,t_0))} = - (f,\mathcal O^S_{t_0}\psi)_{L^2(S\times(0,t_0))}, 
\end{equation} 
for $f \in L^2(S\times(0,t_0)),\  \psi\in L^2(\mathcal N(S,t_0))$.
Clearly, we can take $f \in C_0^{\infty}(S\times(0,t_0))$ with both $f$ and $\psi$ being real-valued. By integration by parts, we have
\begin{eqnarray*}
0 &=& \int_{\mathcal N}\int_0^{t_0}\left((\partial_t^2u^f - \Delta_gu^f)v^{\psi} - u^f(\partial_t^2v^{\psi} - \Delta_gv^{\psi})\right)dtdV_g \\
&=& \int_{\mathcal N}\left[(\partial_tu^f)v^{\psi} - u^f(\partial_tv^{\psi})\right]_{t=0}^{t=t_0}dV_g \\
& & - \int_{\partial {\mathcal N}}\int_0^{t_0}\left((\partial_{\nu}u^f)v^{\psi} - u^f(\partial_{\nu}v^{\psi})\right)dtdS_g.
\end{eqnarray*}
By the initial conditions, $u^f\big|_{t=0} = \partial_tu^f\big|_{t=0} = 0$, and $v^{\psi}\big|_{t=t_0} = 0,\ \partial_tv^{\psi}\big|_{t=t_0} = \psi$. By the boundary condition, $\partial_{\nu}v^{\psi}\big|_{\partial {\mathcal N}\times{\bf R}} = 0$, and $\partial_{\nu}u^f\big|_{\partial {\mathcal N} \times{\bf R}} = f$. We then have
\begin{equation}
\int_{\mathcal N}\mathcal C^S_{t_0}f\psi dV_g = - \int_{\partial {\mathcal N}}\int_0^{t_0}fv^{\psi}dtdS_g.
\nonumber
\end{equation}
Since $f$ is supported in $S\times(0,t_0)$, the right-hand side is rewritten as 
\begin{equation}
 - \int_{S}\int_0^{t_0}fv^{\psi}\big|_{S\times(0,t_0)}dtdS_g  = 
 - (f,\mathcal O^S_{t_0}\psi)_{L^2(S \times (0, t_0))},
\nonumber
\end{equation}
which proves the lemma. \qed

\medskip
By this theorem, for any $\epsilon > 0$ and $a \in L^2({\mathcal N})$ such that ${\rm supp}\,a \subset {\mathcal N}(S,t_0)$, there exists $f = f_{\epsilon, a} \in C_0^{\infty}(S\times(0,t_0))$ satisfying
$\|u^f(\cdot,t_0) - a\|_{L^2({\mathcal N})} < \epsilon$. 
Therefore the property described in Theorem \ref{ApprCont} should be called {\it approximate controllability}.


\subsection{Further results on uniqueness}

Results of the type of Theorem \ref{Tataruglobal}
 (Holmgren-John type uniqueness theorems) have a long story, starting from the classical result by Holmgren:


\begin{theorem}
Let $u$ be a classical,
i.e. $C^2$,  solution to the partial differential equation $P(x,D_x)u = 0$ with analytic coeffcients. If $u = 0$ in one side of a non-characteristic surface $\Sigma$, then ${\rm supp}\,u\cap\Sigma = \emptyset$, i.e. $u = 0$ near $\Sigma$.
\end{theorem}

For the proof, see e.g. \cite{Hor} Vol 1, p. 309 and \cite{Mi73} p. 250. Recall that for a differential operator $P(x,D_x) = \sum_{|\alpha| \leq m}p_{\alpha}(x)D_x^{\alpha}$ defined on an open set $U$ in ${\bf R}^n$, its {\it principal part} is defined by $P_m(x,\xi) = \sum_{|\alpha|=m}p_{\alpha}(x)\xi^{\alpha}$. A surface $\Sigma$ of co-dimension 1 in $U$ is said to be {\it non-characteristic} to $P(x,D_x)$, if $P_m(x,\nu_x) \neq 0$ for any $x \in \Sigma$ and normal $\nu_x$ to $\Sigma$ at $x$. Theorem 4.6 was first proved by E. Holmgren in 1901 \cite{Hol} and extended by F. John in 1949 \cite{Joh49}. This theorem has been tried to 
be extended to the $C^{\infty}$-coefficient case by Robbiano \cite{Rob91} or H{\"o}rmander \cite{Ho92}, and finally Tataru \cite{Ta95} succeeded 
in obtaining the result in full generality (see also \cite{KKL01}, p. 117). The importance of {\it non-analyticity} should largely 
be emphasized in applications to inverse problems. We formulate Tataru's local uniqueness theorem 
in the form convenient for future applications.


\begin{theorem}
Let $u \in H^1_{loc}(\Omega)$, $\Omega \subset \widetilde{\mathcal N} \times{\bf R}$, be a weak solution to the wave equation $\partial_t^2 u = \Delta_{\tilde g}u$, where $(\widetilde{\mathcal N},\tilde g)$ is a Riemannian manifold. Let $\Sigma \subset \Omega$ be a non-characteristic surface. If $u = 0$ on one side of $\Sigma$, then ${\rm supp}\,u\cap\Sigma = \emptyset$.
\end{theorem}

Actually, this theorem implies Theorem \ref{Tataruglobal} due to the fact that we can continue by 0 untill we hit the chracteristic surface giving rise to the double cone of dependence. Note also that this theorem implies more general version of Theorem \ref{Tataruglobal} where condition $\partial_{\nu}u\big|_{\partial{\mathcal N}\times(-t_0,t_0)}=0$ is changed to  $\partial_{\nu}u\big|_{S\times(-t_0,t_0)}=0$ .


\section{Topological reconstruction of ${\mathcal N}$ by $R({\mathcal N})$}


\subsection{Reconstruction from boundary distance functions}
The key idea of the geometric BC-method is to reconstruct the {\it boundary distance function}, $r_x(z)$,  defined as follows: For any $x \in {\mathcal N}$,  $r_x$ is defined  by
\begin{equation}
r_x(z) = d(x,z), \quad z \in \partial {\mathcal N},
\label{C6S5DefineBoundaryDistFunc}
\end{equation}
$d(x,y)$ being the distance of $x, y \in {\mathcal N}$.
We define the map $R$ by
\begin{equation}
R : {\mathcal N} \ni x \to r_x(\cdot) \in C(\partial {\mathcal N}).
\nonumber
\end{equation}
If $\partial {\mathcal N}$ is compact, $R({\mathcal N})$ becomes a metric space by the distance 
\begin{equation}
d_{\infty}(r_1,r_2) = \|r_1(\cdot) - r_2(\cdot)\|_{L^{\infty}(\partial {\mathcal N})},
\nonumber
\end{equation}
and the following inclusion relation hold
\begin{equation}
R({\mathcal N}) \subset C^{0,1}(\partial {\mathcal N}) \subset L^{\infty}(\partial {\mathcal N}),
\nonumber
\end{equation}
where $C^{0,1}(\partial {\mathcal N})$ is the space of Lipschitz continuous functions on $\partial {\mathcal N}$.
The utility of the boundary distance function is seen in the following lemma.


\begin{lemma} If $\partial {\mathcal N}$ is compact, 
$(R({\mathcal N}),d_{\infty})$ is homeomorphic to $({\mathcal N},d)$.
\end{lemma}
Proof. By the triangle inequality, for any $z \in \partial {\mathcal N}$,
$|d(x,z) -d(y,z)| \leq d(x,y)$. 
Hence $\max_{z \in \partial {\mathcal N}}|d(x,z) -d(y,z)| \leq d(x,y)$. This implies
\begin{equation}
d_{\infty}(r_x,r_y) \leq d(x,y).
\label{eq:dinftyrxrylessdxy}
\end{equation}
Both of $(R({\mathcal N}),d_{\infty})$ and $({\mathcal N},d)$ are complete metric spaces. By
(\ref{eq:dinftyrxrylessdxy}), the map $R : ({\mathcal N},d) \to (R({\mathcal N}),d_{\infty})$ is continuous. Let us show that $R$ is injective. Assume $r_x(z) = r_y(z), \ \forall z \in \partial {\mathcal N}$. Let $z_m$ be a point of minimum of $r_x$ and $r_y$. Then $x$ lies on the geodesic normal to $\partial {\mathcal N}$ from $z_m$ at the arclength $r_x(z_m)$, but also $y$ lies on the geodesic normal at arclength $r_y(z_m) = r_x(z_m)$. Then $x = y$. 

We show that $R^{-1}$ is continuous. Suppose $r_{x_n}(\cdot)$ converges to $r_x(\cdot)$ uniformly on $\partial {\mathcal N}$. Then $\sup_n\|r_{x_n}\|_{L^{\infty}} < \infty$. Since $\min r_{x_n} = d(x_n,\partial {\mathcal N})$, and $\partial{\mathcal N}$ is compact, this means that $\{x_n\}$ is in a compact subset in ${\mathcal N}$. Therefore, for any subsequence of $\{x_n\}$, one can select a sub-subsequnce $\{x_n'\}$ such that $x_n'$ converges to some point $y \in {\mathcal N}$. By (\ref{eq:dinftyrxrylessdxy}), $r_{x_n'}(\cdot)$ converges uniformly to $r_{y}(\cdot)$. However, since $r_{x_n}(\cdot)$ converges to $r_x(\cdot)$, we have $r_x(\cdot) = r_y(\cdot)$. Therefore $x = y$. Since every subsequence of $\{x_n\}$ contains a sub-sub sequence which converges to one and the same limit $x$, $x_n$ converges to $x$. This proves the lemma. \qed


\subsection{Metrics on $R({\mathcal N})$}
$R({\mathcal N})$ is a set of functions indexed by the points $x \in {\mathcal N}$. However in the inverse problem we are now considering, we  know neither ${\mathcal N}$ nor $x$, since they are the objects we are trying to reconstruct. 
So,  changing the notation, we let $r_1 = r_x, r_2 = r_y$, where $x, y \in {\mathcal N}$. 
Now we ask a question: {\it Does $d_{\infty}(r_1,r_2)$ determine $d(x,y)$?} If it is true, it becomes a mile stone for our inverse problem. 

Assume we can find new distance $\widehat d(r_1,r_2)$ from $d_{\infty}(r_1,r_2)$ so that $\widehat d(r_1,r_2) = d(x,y)$ for $x, y$ such that $r_1 = r_x$, $r_2 = r_y$. Then $(R({\mathcal N}),\widehat d)$ becomes isometric,
as a metric space,  to $({\mathcal N},d)$. By the Myers-Steenrod theorem \cite{MySt39} (see e.g. \cite{Cha93}, p. 175), this implies that there is a unique Riemannian manifold structure on $R({\mathcal N})$ such that $R : {\mathcal N} \to R({\mathcal N})$ is isometry. In the following, we give a direct way of reconstructing the 
Riemannian manifold structure on $R({\mathcal N})$
to make $R$ a Riemannian isometry from  ${\mathcal N}$ to $R({\mathcal N})$, without leaning over the abstract nature of the Myers-Steenrod theorem.

To find an isometry from $R({\mathcal N})$ to ${\mathcal N}$, perhaps the simplest case is the {\it simple manifold}. By definition (in the strong sense) simple manifold means that any $x, y \in {\mathcal N}$ are connected by a unique shortest geodesic which continues to both directions to $\partial {\mathcal N}$ as the shortest geodesic, and $\partial{\mathcal N}$ is geodesically convex.


\begin{prop}
If ${\mathcal N}$ is simple, then $d_{\infty}(r_x,r_y) = d(x,y)$.
\end{prop}
Proof. Recall (\ref{eq:dinftyrxrylessdxy}).
Let $z$ be the point on $\partial {\mathcal N}$ lying on the continuation of the geodesic from $x$ to $y$. Then 
$d(x,z) - d(y,z) = d(x,y)$. This proves the proposition. \qed


\begin{remark}
It is known that even in the case of non-simple manifold, there exists a constant $0 < C \leq 1$ such that
\begin{equation}
Cd(x,y) \leq d_{\infty}(r_x,r_y) \leq d(x,y).
\nonumber
\end{equation}
\end{remark}


\begin{remark}  Let $\partial {\mathcal N}_1 = \partial {\mathcal N}_2$, and compare $R({\mathcal N}_1)$ and $R({\mathcal N}_2)$. To this end, we can  take the Hausdorff distance $d_{H}(R({\mathcal N}_1),R({\mathcal N}_2))$. 
 Let us recall that if $\mathcal {\mathcal N}$ be a metric space, $S_1, S_2 \subset \mathcal {\mathcal N}$, then the Hausdorff distance is defined by
\begin{equation}
d_H(S_1,S_2) = \max\{\sup_{x\in S_1}d(x,S_2), \sup_{y\in S_2} d(y,S_1)\}.
\nonumber
\end{equation}

A natural question is, if $d_{H}(R({\mathcal N}_1),R({\mathcal N}_2))$ is small,  does it mean that ${\mathcal N}_1$ and ${\mathcal N}_2$ are close and which sense?

In general, the answer is "No", which is the manifestation of well-known {\it ill-posedness} of the inverse problem.
However, we can add some a-priori conditions, e.g. in terms of Gromov compactness on manifolds $({\mathcal N},g)$, to 
obtain a positive answer. See e.g. \cite{AKKLT04}
\end{remark}


\section{Boundary cut locus}
In this and the next sections, we devote ourselves to geometric preliminaries.  For a Riemannian manifold ${\mathcal N}$, let $T_x({\mathcal N})$ be the tangent space at $x \in {\mathcal N}$. 
Recall that for $\xi, \eta \in T_x({\mathcal N})$, the inner product and the length are defined by
$$
g_x(\xi,\eta) = g_{ij}(x)\xi^i\eta^j = \sum_{i,j=1}^ng_{ij}(x)\xi^i\eta^j, \quad |\xi|_g = \sqrt{g_x(\xi,\xi)}
$$
Put $S_x({\mathcal N}) = \{\xi \in T_x({\mathcal N})\,;\,|\xi|_g = 1\}$. Let $T({\mathcal N})$ and $T^{\ast}({\mathcal N})$ be the tangent bundle and the cotangent bundle of ${\mathcal N}$, respectively.

We are dealing with the manifold ${\mathcal N}$ with boundary. To consider the differential at $z\in \partial\mathcal N$ of a map defined on ${\mathcal N}$, we can extend the manifold ${\mathcal N}$ to a bigger manifold $\widetilde {\mathcal N}$ of the same dimension so that $z$ is in the interior of $\widetilde {\mathcal N}$. 
This defines the tangent space $T_z({\mathcal N})$ at $z$ which is independent of the choice of  ${\mathcal N}$. When we consider the tangent sapce of $\partial {\mathcal N}$ at $z \in \partial {\mathcal N}$, we denote it by $T_z(\partial {\mathcal N})$. Note that $T_z(\partial {\mathcal N})$ is canonically identified with the subspace of codimension 1 in $T_z({\mathcal N})$ whose unit normal is the unit normal to $\partial {\mathcal N}$ at $z$.
 

\subsection{Variation and Jacobi fields} 
Let $c(t)$ be a curve on ${\mathcal N}$. For a vector field $X(t)$ on ${\mathcal N}$, with components $(X^1(t),\cdots,X^n(t))$ in local coordinates, the {\it covariant differential} $\displaystyle{\frac{D}{dt}X(t)}$ along $c(t)$ is defined by
\begin{equation}
\nabla_{\dot c}X(t) = \frac{D}{dt}X^k(t) = \dot X^k(t) + \Gamma^k_{ij}(c(t))\dot c^i(t)X^j(t), 
\nonumber
\end{equation}
where we used the abbreviation $\displaystyle{\dot f(t) = \frac{df(t)}{dt}}$. Note that $\nabla_{\dot c}X(t)$ is independent of local coordinates.
A vector field $Z(t)$ is said to be {\it parallel} along $c(t)$ if it satisfies
$\displaystyle{
\frac{D}{dt}Z(t) = 0}$.
In particular, $c(t)$ is a geodesic if and only if $\dot c(t)$ is parallel along $c(t)$. 
For any $C^{\infty}$-curve $c(t)$ and vector fields $\xi(t)$ and $\eta(t)$ along $c(t)$, we have
$$
\frac{d}{dt}\,g_{c(t)}\big(\xi(t),\eta(t)\big) = g_{c(t)}\left(\frac{D}{dt}\xi(t),\eta(t)\right) + g_{c(t)}\left(\xi(t),\frac{D}{dt}\eta(t)\right).
$$

The {\it energy} of a curve $c(t)$ is defined by
\begin{equation}
E(c) = \frac{1}{2}\int_a^{b}g_{c(t)}(\dot c(t),\dot c(t))dt,
\label{C6S6CurveEnergy}
\end{equation}
and the {\it (arc)length} of $c(t)$ is defined by
\begin{equation}
L(c) = \int_a^{b}\sqrt{g_{c(t)}(\dot c(t),\dot c(t))}dt.
\label{C6S6CurveLength}
\end{equation}
Then by the Cauchy-Schwarz inequality, we have
\begin{equation}
L(c)^2 \leq 2(b-a)E(c),
\label{C6S6LandE}
\end{equation}
where the equality holds only when the speed $\sqrt{g_{c(t)}(\dot c(t),\dot c(t))}$ is constant.

A $C^{\infty}$-map : $[a,b]\times(-\epsilon,\epsilon) \ni (t,s) \to H(t,s) \in {\mathcal N}$ is  said to be a {\it variation} of $c(t)$ if $H(t,0) = c(t) \ (a \leq t \leq b)$.  It is said to be a {\it geodesic variation} if for each $s$, the curve : $t \to H(t,s)$ is a geodesic.

For $p \in \mathcal N$ and $v \in T_p(\mathcal N)$, let $c_p(t,v)$ be the geodesic such that $c_p(0,v) = p$, $\dot c_p(0,v) = v$. The {\it exponential map} is defined by 
$$
\exp_p(v) = c_p(1,v).
$$
For any $v \in T_p(\mathcal N)$, the curve : $t \to \exp_p(tv)$ is a geodesic.

The {\it curvature tensor} $R$ is  defined by
\begin{equation}
\left(R(X,Y)Z\right)^l = R^l_{ijk}X^iY^jZ^k,
\nonumber
\end{equation}
\begin{equation}
R^l_{ijk} = \frac{\partial\Gamma^l_{jk}}{\partial x^i} - \frac{\partial\Gamma^l_{ik}}{\partial x^j} + \Gamma^l_{ir}\Gamma^r_{jk} - \Gamma^l_{jr}\Gamma^r_{ik},
\nonumber
\end{equation}
where $X, Y, Z$ are vector fields on ${\mathcal N}$. 
Note that although we use coordinates to define $R^l_{ijk}$, this is actually a $(1,3)$ tensor. It satisfies
\begin{equation}
R(X,Y)Z = \nabla_X(\nabla_Y Z) - \nabla_Y(\nabla_X Z) - \nabla_{[X,Y]}Z.
\label{C6S6RXYZ}
\end{equation}


\begin{lemma}
Let $H(t,s)$ be a variation of $c(t)$, and put $c_s(t) = H(t,s)$. 
We define the vector field $Y(t)$ along $c(t)$ by 
$$
Y(t) = \frac{\partial}{\partial s}H(t,s)\Big|_{s=0}.
$$
Then the following formulae hold. \\
\noindent
(1) The 1st variation formula:
\begin{equation}
\frac{d}{ds}E(c_s)\Big|_{s=0} = 
g_{c(b)}(Y(b),\dot c(b)) - g_{c(a)}(Y(a),\dot c(a)) - \int_a^bg_{c(t)}\Big(Y(t),\frac{D}{dt}\dot c(t)\Big)dt,
\nonumber
\end{equation}
where $D/dt$ is the covariant differential along $c(t)$. \\
\noindent
(2) The 2nd variation formula:
\begin{equation}
\begin{split}
\frac{d^2}{ds^2}E(c_s)\Big|_{s=0} = & \; g_{c(b)}(S(b),\dot c(b)) - g_{c(a)}(S(a),\dot c(a)) \\
& +\int_a^b\Big\{g_{c(t)}\Big(\frac{D}{dt}Y(t),\frac{D}{dt}Y(t)\Big) - 
g_{c(t)}\big(R(Y(t),\dot c(t))\dot c(t),Y(t)\Big) \\
& \ \ \ \ \ \ \ \ \ \ - g_{c(t)}\Big(S(t),\frac{D}{dt}\dot c(t)\Big)\Big\}dt,
\end{split}
\nonumber
\end{equation}
where, letting $D/ds$ be the covariant differential along the curve $C_t(s) : s \to H(t,s)$,
\begin{equation}
S(t) = \frac{D}{ds}\frac{\partial H(t,s)}{\partial s}\Big|_{s=0}.
\label{C6S6DefineS(t)}
\end{equation}.
\end{lemma}

For the proof of above lemma, see e.g.  \cite{GaHuLa80}, Chap. 3


\begin{lemma} \label{Jacobi filed}
Let $c(t)$ $(a \leq t \leq b)$ be a geodesic on ${\mathcal N}$, and $H(t,s)$ its geodesic variation. Then $Y(t) = \partial H(t,s)/\partial s\big|_{s=0}$ satisfies
\begin{equation}
\left(\frac{D}{dt}\right)^2Y + R(Y,\dot c)\dot c = 0, \quad a \leq t \leq b,
\label{C6S5:EquationJacobi}
\end{equation}
where $D/dt$ is the covariant differential along $c(t)$. Conversely, if a vector field $Y(t)$ along the geodesic $c(t)$ satisfies the equation (\ref{C6S5:EquationJacobi}), there is a geodesic variation $H(t,s)$ such that $H(t,0) = c(t)$ and $Y(t) = \partial H(t,s)/\partial s\big|_{s=0}$.
\end{lemma}

Proof. Direct computation shows that
\begin{equation}
\frac{D}{ds}\frac{\partial}{\partial t}H(t,s) = 
\frac{D}{dt}\frac{\partial}{\partial s}H(t,s).
\nonumber 
\end{equation}
Therefore by (\ref{C6S6RXYZ}), 
$$
\frac{D}{dt}\frac{D}{dt}\frac{\partial H}{\partial s} = \frac{D}{\partial t}\frac{D}{\partial s}\frac{\partial H}{\partial t} = 
\left(\frac{D}{\partial s}\frac{D}{\partial t} + R\big(\frac{\partial H}{\partial t},\frac{\partial H}{\partial s}\big)\right)\frac{\partial H}{\partial t}.
$$
Since $c_s(t)$ are geodesics, $D(\partial H(t,s)/\partial t)/dt = 0$. Thus, letting $s = 0$, we obtain $(D/dt)^2Y = R(\dot c,Y)\dot c$, which proves (\ref{C6S5:EquationJacobi}).

Conversely, suppose $Y(t)$ satisfies (\ref{C6S5:EquationJacobi}). Take a curve $z(s)$ such that $z(0) = c(a)$, $\dot z(0) = Y(a)$. Let $X_0(s)$, $X_1(s)$ are vector fields which are parallel along $z(s)$, and satisfy $X(0) = \dot c(a)$, $X_1(0) = (DY/dt)(a)$. We put
$$
H(t,s) = \exp_{z(s)}\left((t-a)\big(X_0(s) + sX_1(s)\big)\right).
$$
Then the curve : $t \to H(t,s)$ is a geodesic for each $s$, and $H(t,0) = c(t)$. Let $Z(t) = \partial H(t,s)/\partial s\big|_{s=0}$.
Then, as has been shown above, $Z(t)$ satisfies (\ref{C6S5:EquationJacobi}). Moreover, $Z(a) = \dot z(0) = Y(a)$. Then 
\begin{equation}
\begin{split}
\frac{D Z}{dt}(a) &= \frac{D}{dt}\frac{\partial H}{\partial s}\Big|_{t=a,s=0} = \frac{D}{ds}\frac{\partial H}{\partial t}\Big|_{t=a,s=0} \\
&= \frac{D}{ds}\left(X_0(s) + sX_1(s)\right)\big|_{s=0} \\
&= X_1(0) = \frac{DY}{dt}(a),
\end{split}
\nonumber
\end{equation}
where in the last step, we use $X_0(s)$, $X_1(s)$ are parallel along $z(s)$.
Therefore $Y(t) = Z(t)$ by the uniqueness for solutions of differential equations. \qed

\medskip
A solution $Y(t)$ of (\ref{C6S5:EquationJacobi}) is called {\it Jacobi field} along $c(t)$.


\subsection{Focal point}
In the following, we consider the {\it boundary normal geodesic}, denoted by $\gamma_z(t)$ or $\exp_{\partial {\mathcal N}}(z,t)$, starting from $z \in \partial {\mathcal N}$ with initial direction being the inner unit normal at $z$. Explicitly,  take local coordinates $z = (z_1, \cdots,z_{n-1})$ on $\partial \mathcal N$, and $(z_1,\cdots,z_{n-1},x_n)$, where $x_n=0$ is a defining equation of $\partial\mathcal N$, as local coordinates in $\mathcal N$.
Coinsider the equation of geodesics
\begin{equation}
\left\{
\begin{split}
&\frac{d^2 x^k}{dt^2} + \Gamma^{k}_{ij}(x(t))\frac{dx^i}{dt}\frac{dx^j}{dt} = 0, \\
& x(0) = (z,0), \quad \frac{dx}{dt}(0) = \nu(z),
\end{split}
\right.
\nonumber
\end{equation}
where $\nu(z)$ is the unit normal at the boundary. Then, the map $\gamma_{z}(t) : (z,t) \to x(t,z)$ is a diffeomorphism near $\partial {\mathcal N}$, and we use $(z,t)$ as boundary normal coordinates in $\mathcal N$ near $\partial\mathcal N$.


\begin{prop}
In the boundary normal coordinates, the Riemannian metric is written as
$$
ds^2 = (dt)^2 + \sum_{i,j=1}^{n-1}h_{ij}(z,t)dz^idz^j.
$$
\end{prop}

Proof. 
Since $x(t)$ is a geodesic, we have
$$
g_{nn} = g\Big(\frac{\partial x}{\partial t},\frac{\partial x}{\partial t}\Big) = 1.
$$ 
For $1 \leq i \leq n-1$, we have
\begin{equation}
\begin{split}
\frac{d}{dt}g_{ni} & = \frac{d}{dt}g\Big(\frac{\partial x}{\partial t},\frac{\partial x}{\partial z^i}\Big) = g\Big(\frac{\partial x}{\partial t},\frac{D}{dt}\frac{\partial x}{\partial z^i}\Big) \\
&= g\Big(\frac{\partial x}{\partial t},\frac{D}{\partial z^i}\frac{\partial x}{\partial t}\Big) = 
\frac{1}{2}\frac{\partial}{\partial z^i}g\Big(\frac{\partial x}{\partial t},\frac{\partial x}{\partial t}\Big) = 0.
\end{split}
\nonumber
\end{equation}
Since $\dfrac{dx}{dt}(0)=\nu(z)$ is normal to $\partial\mathcal N$, $g_{ni}(z,0)=0$. Therefore, $g_{ni}=0$, and the proof is completed. \qed

\medskip
Fixing $t$, we define the map $\exp_{\partial {\mathcal N}}(\cdot,t)$ by
\begin{equation}
\exp_{\partial {\mathcal N}}(\cdot,t) : \partial \mathcal N \ni z \to \gamma_z(t) \in {\mathcal N}.
\nonumber
\end{equation}
Let $d_{\partial {\mathcal N}}\exp_{\partial {\mathcal N}}(z_0,t) : T_{z_0}(\partial {\mathcal N}) \to T_{\gamma_{z_0}(t)}({\mathcal N})$ be the differential of $\exp_{\partial {\mathcal N}}(\cdot,t)$ evaluated at $z_0$.


\begin{definition}
Let $\gamma_{z_0}(t)$ be the boundary normal geodesic starting from $z_0 \in \partial {\mathcal N}$. The point $\gamma_{z_0}(t_0) = \exp_{\partial {\mathcal N}}(z_0,t_0)$ is called a {\it focal point} along $\gamma_{z_0}(t)$ if
$$
{\rm rank}\left(d_{\partial {\mathcal N}}\exp_{\partial \mathcal N}(z_0,t_0)\right) < n - 1.
$$
\end{definition}


\begin{lemma}\label{FocalpointLemma}
Let $ \gamma_{z_0}(t)$ $(0 \leq t \leq t_0)$ be a boundary normal geodesic starting from 
$z_0 \in \partial {\mathcal N}$. If $\gamma_{z_0}(t_1)$ is a focal point along $\gamma_{z_0}$ for 
some $0 < t_1 < t_0$, then  $\tau=d(\gamma_{z_0}(t_0)<t_0$ and there exist
 $w \in \partial {\mathcal N}$ such that $\gamma_w(\tau)=\gamma_{z_0}(t_0)$. 
\end{lemma}

Note that this lemma is a particular case of Fermi coordinates
associated with $k$-dimensional submanifold in $\mathcal N$, where $k <n$. 
See 
\cite{Cha93},  \S 3.6.  See  \cite {BiCri64}, p. 232, or \cite{Sak96}, Chap. 3, Lemma 2.11 for the complete proof.

We prove this lemma under the following additional assumption.

\medskip
\noindent
{\it Condition (TG)} : In a neighborhood of $z_0$, we can extend $\mathcal N$ to a bigger 
manifold
$\widetilde{\mathcal N}$ so that, in a neighborhood of $z_0$, $\partial{\mathcal N}$ is a 
totally geodesic submanifold of $\widetilde{\mathcal N}$. 

\medskip
Let us recall that, given  a Riemannian manifold  $\widetilde{\mathcal N}$,  its submanifold 
$\mathcal S$  is said to be {\it totally geodesic} if 
any geodesic of 
$\widetilde{\mathcal N}$ starting from a point $z \in \mathcal S$ in a direction tangential
to $\mathcal S$
 lies in 
$\mathcal S$. Note that, if $\hbox{dim}(\mathcal S)=n-1$, which is the case of
$\mathcal S=\partial{\mathcal N}$, this condition is
equivalent to the fact that the second fundamental form (the shape operator) 
of $\mathcal S$ vanishes. In turn, this is equivalent to the fact that $\nu(z)$ is
parallel along $\mathcal S$.

For example, if for some $\epsilon > 0$, 
$ \widetilde{\mathcal N} = \mathcal S\times(-\epsilon,\epsilon)$, 
and the metric of $\widetilde{\mathcal N}$ is of product form:
$$
ds^2 = (dt)^2 + h(\omega,d\omega),
$$
where $h(\omega,d\omega)$ is the positive definite metric on 
$\mathcal S$ induced from that of $\widetilde{\mathcal N}$, then 
$\mathcal S$ is totally geodesic.

\medskip

Proof of Lemma (\ref{FocalpointLemma}). By the assumption, there exists $0 \neq \xi \in T_{z_0}(\partial {\mathcal N})$ such that 
\begin{equation}
\left(d_{\partial {\mathcal N}}\exp_{\partial {\mathcal N}}(z_0,t_1)\right)\xi = 0.
\label{ddNexpt1=xi=0}
\end{equation}
Let $z(s)$ be a geodesic in $\widetilde{\mathcal N}$ such that $z(0) = z_0$, $\dot z(0) = \xi$. 
By the condition (TG), $z(s)$ is also a geodesic in $\partial\mathcal N$.
We put
$$
\widetilde H(t,s) = \left(\exp_{\partial {\mathcal N}}(t)\right)(z(s)) = \gamma_{z(s)}(t),
$$
\begin{equation}
\widetilde Y(t) = \frac{\partial \widetilde H(t,s)}{\partial s}\Big|_{s=0}.
\nonumber
\end{equation}
Then,  by Lemma \ref{Jacobi filed}, $\widetilde Y(t)$ is a Jacobi field along $c(t)$ and satisfies
\begin{equation}
\widetilde Y(0) = \xi, \quad \widetilde Y(t_1) = 0.
\label{C6S5:widetildeYproperty}
\end{equation} 
These facts  follow from $\widetilde H(0,s) = z(s)$, (\ref{ddNexpt1=xi=0}), and
\begin{equation}
\begin{split}
\frac{\partial}{\partial s}\widetilde H(t_1,s)\big|_{s=0}  = 
\frac{\partial}{\partial s}\exp_{\partial N}(t_1)(z(s))\big|_{s=0} = \left(d_{\partial\mathcal N}\exp_{\partial\mathcal N}(z_0,t_1)\right)\xi.
\end{split}
\nonumber
\end{equation}

Take a parallel vector field $Z(t)$ satisfying
\begin{equation}
\left\{
\begin{split}
& \frac{D}{dt}Z(t) = 0, \quad {\rm for} \quad 0<t<t_0, \\
& Z(t_1) = - \frac{D}{dt}\widetilde Y(t_1).
\end{split}
\right.
\label{C6S5:Zt1}
\end{equation}
 Pick $f(t) \in C_0^{\infty}((0,t_0))$ such that $f(t_1) = 1$, and put for 
 $\alpha \in {\bf R}$
\begin{equation}
V_{\alpha}(t) = 
\left\{
\begin{split}
& \widetilde Y(t) + \alpha f(t)Z(t), \quad  0 \leq t \leq t_1, \\
&  \alpha f(t)Z(t), \quad t_1 \leq t \leq t_0.
\end{split}
\right.
\label{C6S5Valphat}
\end{equation}
Note that at $t = t_1$, $V_{\alpha}(t)$ is continuous by (\ref{C6S5:widetildeYproperty}), 
however, $\frac{D}{dt} V_{\alpha}(t)$ is discontinuous. As a variation of 
$c(t)=\gamma_{z_0}(t)$, we 
consider 
\begin{equation}
H_{\alpha}(t,s) = \exp_{c(t)}(sV_{\alpha}(t)).
\label{C6S6DefineHts}
\end{equation}
Let $c_{\alpha,s}(t)$ be the curve : $ t \to H_{\alpha}(t,s)$. 
Then $c_{\alpha,0}(t) = c(t)$ for all $\alpha$.
Define the energy of $c_{\alpha,s}(t)$ by (\ref{C6S6CurveEnergy}).
We can then prove the following formula.


\begin{prop} For small $|\alpha|$, we have
\begin{equation}
 \frac{d^2}{ds^2}E(c_{\alpha,s})\Big|_{s=0} 
=  - 2\alpha \, g_{c(t_1)}\Big(\frac{D\widetilde Y}{dt}(t_1),\frac{D\widetilde Y}{dt}(t_1)\Big) + O(\alpha^2).
\label{C6S6dsdsE(cs)}
\end{equation}
\end{prop}

\medskip
 Granting this proposition for the moment, we complete the proof of Lemma 6.5.
 We have $\displaystyle{\frac{D\widetilde Y}{dt}(t_1) \neq 0}$. In fact, if this vanishes, since $\widetilde Y(t_1) = 0$ and $\widetilde Y(t)$ is a solution of the 2nd order differential equation, $\widetilde Y(t)$ vanishes identically.  Proposition 6.6 then yields 
\begin{equation}
 (d/ds)^2E(c_{\alpha,s})\big|_{s=0} < 0,
\label{C6S6E(Cs)''<0}
\end{equation}
if $\alpha > 0$ is chosen small enough. 
Letting
$$
Y_{\alpha}(t) = \partial H_{\alpha}(t,s)/\partial s\big|_{s=0} = V_{\alpha}(t),
$$
and using
$
Y_{\alpha}(0) = \widetilde Y(0) = \xi, 
Y_{\alpha}(t_0)  = 0, $
 we have by Lemma 6.1 (1), 
$$
(d/ds)E(c_{\alpha,s})\big|_{s=0} = 0.
$$
This, combined with (\ref{C6S6E(Cs)''<0}), implies $E(c_{\alpha,s}) < E(c)$, 
for $0 < s < \epsilon$, if $\epsilon > 0$ is small enough. 
For $0 < s < \epsilon$, we have, by the Cauchy-Schwarz inequality (\ref{C6S6LandE}),
$$
L(c_{\alpha,s})^2 \leq 2t_0 E(c_{\alpha,s}) < 2t_0 E(c) = L(c)^2,
$$
where in the last step we use the fact $c_0(t)$ is a unit speed geodesic.
Therefore, $d(\gamma_{z_0}(t_0), \partial \mathcal N) <t_0$, which implies
an existence of $w \in \partial \mathcal N$ with desired property.
This proves Lemma 6.5. \qed

\medskip
Now we prove Proposition 6.6. We split energy into 2 parts:
\begin{equation}
\begin{split}
E(c_{\alpha,s}) & = 
\frac{1}{2}\int_0^{t_1}g_{c_{\alpha,s}}(t)(\dot c_{\alpha,s}(t),\dot c_{\alpha,s}(t))dt + 
\frac{1}{2}\int_{t_1}^{t_0}g_{c_{\alpha,s}}(t)(\dot c_{\alpha,s}(t),\dot c_{\alpha,s}(t))dt \\
& =: E_1 (c_{\alpha,s}) + E_2(c_{\alpha,s}).
\end{split}
\nonumber
\end{equation}
 Let $S_{\alpha}(t)$ be defined by (\ref{C6S6DefineS(t)}). Then, by Lemma 6.1 (2),
\begin{equation}
\begin{split}
 \frac{d^2}{ds^2}E_1(c_{\alpha,s})\Big|_{s=0} &= g_{c(t_1)}(S_{\alpha}(t_1),\dot c(t_1)) 
- g_{c(0)}(S_{\alpha}(0),\dot c(0))\\
 &+ \int_0^{t_1}
\Big\{g\Big(\frac{D}{dt}V_{\alpha},\frac{D}{dt}V_{\alpha}\Big) - 
g\big(R(V_{\alpha},\dot c)\dot c,V_{\alpha}\big)\Big\}dt.
\end{split}
\nonumber
\end{equation}
Since $DZ/dt = 0$, the integral in the right-hand side is equal to
\begin{equation}
\begin{split}
  &\ \int_{0}^{t_1}
\Big\{g\Big(\frac{D\widetilde Y}{dt} + \alpha\dot fZ,\frac{D\widetilde Y}{dt} + \alpha\dot fZ\Big) - 
g\big(R(\widetilde Y + \alpha fZ,\dot c)\dot c,\widetilde Y + \alpha fZ\big)\Big\}dt \\
= & \int_0^{t_1}\Big\{g\Big(\frac{D\widetilde Y}{dt},\frac{D\widetilde Y}{dt}\Big) - g\big(R(\widetilde Y,\dot c)\dot c,\widetilde Y\big)\Big\}dt \\
& + 2\alpha\int_0^{t_1}\Big\{g\Big(\dot fZ,\frac{D\widetilde Y}{dt}\Big) - 
g\big(R(\widetilde Y,\dot c)\dot c,fZ\big)\Big\}dt + O(\alpha^2).
\end{split}
\nonumber
\end{equation}
Since $\widetilde Y$ is a Jacobi field, it satisfies (\ref{C6S5:EquationJacobi}). This imples
\begin{equation}
\begin{split}
 \frac{d^2}{ds^2}E_1(c_s)\Big|_{s=0} 
&= \  g_{c(t_1)}(S_\alpha(t_1),\dot c(t_1)) 
- g_{c(0)}(S_\alpha(0),\dot c(0))\\
& + 
\int_0^{t_1}\Big\{g\Big(\frac{D\widetilde Y}{dt},\frac{D\widetilde Y}{dt}\Big) + g\big(\frac{D^2\widetilde Y}{dt^2},\widetilde Y\big)\Big\}dt \\
& + 2\alpha\int_0^{t_1}\Big\{g\Big(\dot fZ,\frac{D\widetilde Y}{dt}\Big) +
g\big(\frac{D^2\widetilde Y}{dt^2},fZ\big)\Big\}dt + O(\alpha^2).
\end{split}
\label{C6S6d2ds2E1compute1}
\end{equation}
Then two integrals of the right-hand side are computed as
\begin{equation}
\begin{split}
&  \int_0^{t_1}\frac{d}{dt}g\Big(\frac{D\widetilde Y}{dt},\widetilde Y\Big)dt 
+ 2\alpha\int_0^{t_1}\frac{d}{dt}g\big(\frac{D\widetilde Y}{dt},fZ\big)dt \\
= \ & g_{c(t_1)}\Big(\frac{D\widetilde Y}{dt}(t_1-0),\widetilde Y(t_1)\Big) - g_{c(0)}\Big(\frac{D\widetilde Y}{dt}(0),\widetilde Y(0)\Big) \\
& + 2\alpha\Big\{g_{c(t_1)}\Big(\frac{D\widetilde Y}{dt}(t_1),f(t_1)Z(t_1)\Big) - 
g_{c(0)}\Big(\frac{D\widetilde Y}{dt}(0),f(0)Z(0)\Big)\Big\}.
\end{split}
\label{C6S6d2ds2E1compute2}
\end{equation}
Recall that $\widetilde Y(t_1) = 0$. We also note that the curve : $s \to H(t,s) = \exp_{c(t)}(sV_{\alpha}(t))$ is a geodesic for $t \geq 0$.  Then we have
\begin{equation}
S_\alpha(t) = \frac{D}{ds}\frac{\partial H(t,s)}{\partial s}\Big|_{s=0} = 0, \quad t \geq 0.
\label{C6S6S(t)=0}
\end{equation}
We show that $\frac{D\widetilde Y}{dt}(0) = 0$. In fact, 
since 
\begin{equation}
\frac{D}{dt}\widetilde Y(0) = \frac{D}{dt}\frac{\partial\widetilde H}{\partial s}\Big|_{s=t=0} = \frac{D}{ds}\frac{\partial\widetilde H}{\partial t}\Big|_{s=t=0} = \frac{D}{ds}\nu(z(s))\Big|_{s=0}=0.
\label{C6SDdtWidetildeY=0}
\end{equation}
where the last equation follows from vanishing of the second fundamental form in $z_0$.
Plugging (\ref{C6S6d2ds2E1compute1}) $\sim$ (\ref{C6SDdtWidetildeY=0}), we obtain
\begin{equation}
\frac{d^2}{ds^2}E_1(c_{\alpha,s})\Big|_{s=0} =  2\alpha g_{c(t_1)}\Big(\frac{D\widetilde Y}{dt}(t_1),Z(t_1)\Big) + O(\alpha^2).
\label{C6S6d2ds2E1cs}
\end{equation}

We turn to  $E_2(c_{\alpha,s})$. As above, 
\begin{equation}
\begin{split}
 \frac{d^2}{ds^2}E_2(c_{\alpha,s})\Big|_{s=0} &= g_{c(t_0)}(S_{\alpha}(t_0),\dot c(t_0)) 
- g_{c(t_1)}(S_{\alpha}(t_1),\dot c(t_1))\\
 &+ \int_{t_1}^{t_0}
\Big\{g\Big(\frac{D}{dt}V_{\alpha},\frac{D}{dt}V_{\alpha}\Big) - 
g\big(R(V_{\alpha},\dot c)\dot c,V_{\alpha}\big)\Big\}dt.
\end{split}
\nonumber
\end{equation}
 We compute in the same way as for $E_1(c_{\alpha,s})$. Since 
 $\widetilde Y$ does not appear in this case, we have
\begin{equation}
 \frac{d^2}{ds^2}E_2(c_{\alpha,s})\Big|_{s=0} 
=  O(\alpha^2)
\label{C6S6d2ds2E2cs}
\end{equation}
In view of (\ref{C6S5:Zt1}), (\ref{C6S6d2ds2E1cs}) and (\ref{C6S6d2ds2E2cs}), we have completed the proof. \qed


\begin{remark} The above proof can be immediately extended to the case when the
second fundamental form of $\partial \mathcal N$ vanishes just at the point $z_0$.
Indeed, the above proof shows that, for sufficiently small $\alpha >0$ and $|s|$,
$$
d(z(s), \gamma_{z_0}(t_0)) < t_0- c \alpha s^2.
$$
Since $d(z(s), \partial \mathcal N)=O(|s|^3)$, the result follows.
\end{remark}


\subsection{Boundary cut point}
Let $\gamma_z(\cdot)$ be the boundary normal geodesic starting from $z \in \partial {\mathcal N}$. 
A point $\gamma_z(t)$ is said to be {\it uniquely minimizing} along the geodesic 
$\gamma_z(\cdot)$ if $t = d(\gamma_z(t),\partial {\mathcal N})$ and $t < d(\gamma_z(t),w)$ 
for any $w \in \partial {\mathcal N}$ such that $w \neq z$. Thus, 
$\{\gamma_z(s)\, ; \, 0 \leq s \leq t \}$ is a unique shortest geodesic from 
$\partial {\mathcal N}$ to $\gamma_z(t)$. 


\begin{lemma}\label{Notuniquelymini}
Let $\gamma_z(t)$ $ (0 \leq t \leq t_0)$ be the boundary normal geodesic starting from 
$z \in \partial {\mathcal N}$. If $\gamma_z(t_1)$ is not uniquely minimizing for some 
$0 < t_1 < t_0$, then $d(\gamma_{z}(t_0),\partial {\mathcal N}) < t_0$.
\end{lemma}

Proof.  Since $\gamma_z(t_1)$ is not uniquely minimizing, there exists 
$w \in \partial\mathcal N$ such that $\gamma_w(t) = \gamma_z(t_1)$, $t \leq t_1$. Consider a once broken geodesics $c(s) = \gamma_w([0,t])\cup\gamma_{z}([t_1,t_0])$. Here, for any curve $c(s)$, by $c([a,b])$ we denote the piece of $c(s)$ for $s \in [a,b]$. Then $\gamma_z(t_0) = c(s)$, $s = t_0+ (t-t_1)$. This proves the lemma when $t < t_1$.

For $t=t_1$, consider a curve $c(s)$ which consists of 3 parts: the geodesic $\gamma_w(s), 0 \leq s \leq t-\epsilon$, the minimizing geodesic $c'(\tau)$ connecting $\gamma_w(t-\epsilon)$ and $\gamma_z(t_1+\epsilon)$, and the piece of geodesic $\gamma_z(s)$ for $t_1 + \epsilon \leq s \leq t_0$. Note that, by the short-cut arguments, $L(c') < 2\epsilon$. Therefore,
$$
L(c) = (t-\epsilon) + L(c') + (t_0-(t_1+\epsilon)) 
< t_0-(t_1-t) = t_0,
$$
which proves the lemma. \qed

\medskip
By the above lemma, if $\gamma_z(t)$ is uniquely minimizing along $\gamma_z(\cdot)$, then so is $\gamma_z(s)$ for any $0 < s < t$. We put 
\begin{equation}
\tau(z) = \sup\{t\, ; \, \gamma_z(t) \ {\rm is \ uniquely \ minimizing}\}. 
\label{DefineTau(z)}
\end{equation}
We then have
\begin{equation}
d(\gamma_z(t),\partial {\mathcal N}) < t, \quad {\rm for}\quad \tau(z) < t. 
\nonumber
\end{equation}
In fact, we have only to take $\tau(z) < t_1 < t$ and apply Lemma \ref{Notuniquelymini}.


\begin{definition}\label{Definebdcutpoint}
The function $\tau(z)$ defined by (\ref{DefineTau(z)}) is called the {\it boundary cut function}, and the point $\gamma_{z}(\tau(z))$ for $\tau(z) < \infty$ is called  {\it boundary cut point} of $z$ along $\gamma_z$.  If $\tau(z) = \infty$, we say that there is no boundary cut point along the boundary normal geodesic $\gamma_z$.
\end{definition}


\begin{lemma} \label{2caseatbdcutpt}
For $z_0 \in \partial {\mathcal N}$, let $\tau(z_0)$ be as in Definition \ref{Definebdcutpoint}.
At the boundary cut point, 
$$
d(\gamma_{z_0}(\tau(z_0)),z_0) = \tau(z_0),
$$
 and at least one (possibly both) of the following statements holds:

\noindent
(a) $\gamma_{z_0}(\tau(z_0))$ is an ordinary boundary cut point, i.e. there is $w \in \partial {\mathcal N}$ such that $w \neq z_0$ and $\gamma_{z_0}(\tau(z_0)) = \gamma_{w}(\tau(z_0))$. \\
\noindent
\noindent
(b)  $\gamma_{z_0}(\tau(z_0))$ is the first focal point along $\gamma_{z_0}$, i.e. 
\begin{equation}
\begin{split}
{\rm rank}\left(d_{\partial {\mathcal N}}\exp_{\partial {\mathcal N}}(z_0,t)\right)  = n - 1 \quad  & {\rm if} \quad t < \tau(z_0), \\
{\rm rank}\left(d_{\partial {\mathcal N}}\exp_{\partial {\mathcal N}}(z_0,t)\right)  < n - 1 \quad  & {\rm if} \quad t = \tau(z_0).
\end{split}
\nonumber
\end{equation}
\end{lemma}
Proof. 
By definition, we have $d(\gamma_{z_0}(s),\partial {\mathcal N}) = s$ for $s < \tau(z_0)$. Letting $s \to \tau(z_0)$, we have $d(\gamma_{z_0}(\tau(z_0)),\partial {\mathcal N}) = \tau(z_0)$. This implies, by Lemma 6.5, $\gamma_{z_0}(s)$ is not a focal point for $0 < s < \tau(z_0)$ . 

 There exists $\delta > 0$ such that the geodesic $\gamma_{z_0}(t)$ exists in the interval $[0,\tau(z_0) + \delta]$. Take a sequence $\delta > \epsilon_1 > \epsilon_2 \cdots \to 0$ and put $t_n = \tau(z_0) + \epsilon_n$. Then, by the definition of $\tau(z_0)$, there exists $w_n \in \partial {\mathcal N}$, $w_n\neq z_0$, and $s_n < t_n$ such that $\gamma_{w_n}(s_n) = \gamma_z(t_n)$.
 Since $\partial {\mathcal N}$ is compact, there exists a subsequence $\{w_n, s_n\}$, such that 
 $w_n \to w \in \partial {\mathcal N}$, $s_n\to s$, where $0 \leq s \leq \tau(z_0)$. Then $\gamma_w(s) = \gamma_{z_0}(\tau(z_0))$, which implies $s = \tau(z_0)$. This gives rise to ordinary boundary cut point if $w\neq z_0$.

Suppose $w = z_0$. Let us show that $\gamma_{z_0}(\tau_{z_0})$ is the first focal point along $\gamma_{z_0}$. Assume that ${\rm rank}\left(d_{\partial {\mathcal N}}\exp_{\partial {\mathcal N}}(z_0,\tau(z_0))\right)  = n - 1$. Take a small neighborhood $V$ of $z_0$ in $\partial {\mathcal N}$ and small $\epsilon > 0$.
Then the map : $V\times(\tau(z_0)-\epsilon,\tau(z_0)+\epsilon) \ni (z,t) \to \gamma_{z}(t)$ is a diffeomorphsim.
Therefore, in a small neighborhood $U$ of $\gamma_{z_0}(\tau(z_0))$, $\gamma_{z}(t)^{-1}$ is a diffeomorphism. Since $w_n\to z_0$ and $s_n\to \tau(z_0)$, 
$\gamma_{w_n}(s_n)\in U$.
However, $\gamma_{z_0}(t_n)\in U$, and $\gamma_{z_0}(t_n) = \gamma_{w_n}(s_n)$. We thus arrive at the contradiction. By Lemma \ref{FocalpointLemma}, for $t < \tau_0$, $\gamma_{z_0}(t)$ is not a focal point.
\qed

\medskip
We introduce a topology in ${\bf R}_+\cup\infty$ by taking intervals $(a,b)$ and $(a,\infty] = (a,\infty)\cup\infty$ as basis for the open sets.


\begin{lemma}\label{tau(z)conti}
The function $\tau(z)$ in Definition \ref{Definebdcutpoint} is continuous from $\partial {\mathcal N}$ to ${\bf R}_+\cup\infty$.
\end{lemma}
Proof. Suppose $\tau(z)$ is not continuous at $\overline{z} \in \partial {\mathcal N}$, and let $z_k \in \partial {\mathcal N}$ be such that $z_k \to \overline z$ and $\lim\tau(z_k) \neq \tau(\overline z)$.
 Set $\tau_k = \tau(z_k)$, $\tau_{\infty} = \lim \tau(z_k)$ and $\overline{\tau} = \tau(\overline z)$.

We first consider the case $\overline{\tau} > \tau_{\infty}$.
Since $\overline\tau = \tau(\overline z) > \tau_{\infty}$, then $\tau_{\infty} <\infty$ and by Lemma \ref{FocalpointLemma}, $\exp_{\partial {\mathcal N}}(\tau_{\infty},\overline{z})$ is not a focal point along the boundary normal geodesic $\gamma_{\overline z}(t)$. Therefore, ${\rm rank}(d_{\partial {\mathcal N}}\exp_{\partial {\mathcal N}}(\overline z,\tau_{\infty})) = n - 1$. Then, there is a neighborhood $V$ of $\overline{z}$ in $\partial {\mathcal N}$ and $\epsilon > 0$ such that the map $V\times(\tau_{\infty}-\epsilon,\tau_{\infty}+\epsilon) \ni (z,t) \to \exp_{\partial {\mathcal N}}(t,z)$ is a diffeomorphism.
Since $z_k \to \overline{z}$, $\tau_k \to \tau_{\infty}$, we have $(z_k,\tau_k) \in V\times(\tau_{\infty}-\epsilon,\tau_{\infty}+\epsilon)$ for large $k$.
Therefore,
${\rm rank}(d_{\partial {\mathcal N}}\exp_{\partial {\mathcal N}}(z_k,\tau_k)) = n - 1$ for large $k$.
Then by Lemma \ref{2caseatbdcutpt}, $\exp_{\partial {\mathcal N}}(z_k,\tau_k)$ is not the focal point along the boundary normal geodesic $\exp_{\partial {\mathcal N}}(z_k,t)$, but the ordinary boundary cut point, i.e. there exists $w_k \in \partial {\mathcal N}$ such that $w_k \ne z_k$ and $\exp_{\partial {\mathcal N}}(w_k,\tau_k) = \exp_{\partial {\mathcal N}}(z_k,\tau_k)$. We see that $w_k \not\in V$, since $\exp_{\partial {\mathcal N}}$ is a diffeomorphism on $V\times(\tau_{\infty}-\epsilon,\tau_{\infty}+\epsilon)$. By taking a subsequence if necessary, we can assume that $w_k$ converges to $w \in \partial {\mathcal N}$. By shrinking $V$ if necessary, we have $w \not\in V$. We than have
\begin{equation}
\begin{split}
\exp_{\partial {\mathcal N}}(w,\tau_{\infty}) &= \lim \exp_{\partial {\mathcal N}}(w_k,\tau_k) 
= \lim\exp_{\partial {\mathcal N}}(z_k,\tau_k) \\
&= \exp_{\partial {\mathcal N}}(\overline z,\tau_{\infty}).
\end{split}
\nonumber
\end{equation}
This contradicts Lemma \ref{Notuniquelymini} and the definition of $\overline{\tau} = \tau(\overline{z})$.

Next we assume $\overline{\tau} < \tau_{\infty}$.
Take $\overline{\tau} < \tau < \infty$. Then, there is $w\in\partial\mathcal N$ and $s <\tau$ such that $\gamma_{\overline{z}}(\tau) = \gamma_w(s)$. Since $z_k \to \overline{z}$, $\gamma_{z_k}(\tau) \to \gamma_{\overline{z}}(\tau)$. By the triangle inequlaity, 
$$
d(w,\gamma_{{z_k}}(\tau)) \leq d(w,\gamma_{\overline{z}}(\tau)) +
d(\gamma_{\overline{z}}(\tau),\gamma_{{z_k}}(\tau)) = s + d(\gamma_{\overline{z}}(\tau),\gamma_{{z_k}}(\tau)).
$$
Since $s < \tau$, taking $k$ large enough, we see that $d(w,\gamma_{{z_k}}(\tau)) < \tau$. Since $\tau < \tau_{\infty}$, so that $\tau < \tau(z_k)$ for large $k$, we get the contradiction.
 \qed


\subsection{Boundary cut locus. Boundary normal coordinates}

 
\begin{definition}\label{Sefinrbdcutlocus}
The {\it boundary cut locus} $\omega$ is defined by
\begin{equation}
\omega = \{\gamma_z(\tau(z))\, ;\, z \in \partial {\mathcal N}\},
\nonumber
\end{equation}
where $\gamma_z(\tau(z))$ is the boundary cut point of $z$ along the boundary normal geodesic $\gamma_z(t) = \exp_{\partial {\mathcal N}}(z,t)$ in Definition 6.8. 
\end{definition}
Recall that by Lemma \ref{2caseatbdcutpt}, we have $d(\gamma_z(\tau(z)),z) = \tau(z)$. Let us investigate the structure of $\omega$. We put
\begin{equation}
B(\mathcal N) = \mathop\cup_{z\in\partial {\mathcal N}}\{\gamma_z(t) \, ; \, 0 \leq t < \tau(z)\big\}.
\nonumber
\end{equation}


\begin{lemma} \label{Bnomega} (1) ${\mathcal N} = B({\mathcal N})\cup \omega$, $B({\mathcal N})\cap\omega = \emptyset$. \\
\noindent
(2) $\omega$ is a closed set of measure 0. In particular, it has no interior points.\\
\noindent
(3) $B({\mathcal N})$ is an open set.
\end{lemma}
Proof. For any $x \in {\mathcal N}$, there exists $z_x \in \partial {\mathcal N}$ such that $d(x,z_x) = d(x,\partial {\mathcal N}):=s(x)$. 
Therefore $x = \gamma_{z(x)}(s(x))$ (see Lemma \ref{Lemma6.3.1}).
Let us prove $s(x) \leq \tau(z_x)$, where $\tau(z)$ is boundary cut function, see Definition \ref{Definebdcutpoint}. Indeed, if $s(x) > \tau(z_x)$,  there exists $w \in \partial {\mathcal N}$ such that $d(x,w) < s(x)$, which is a contradiction, since $s(x) = d(x,\partial {\mathcal N})$. 

Therefore, we have shown that, for any $x \in {\mathcal N}$, there exists $z_x \in \partial {\mathcal N}$ such that $x = \exp_{\partial {\mathcal N}}(z_x,d(x,\partial {\mathcal N}))$ and $d(x,\partial {\mathcal N}) \leq \tau(z_x)$. This proves ${\mathcal N} = B({\mathcal N})\cup\omega$. 

The disjointness of $B$ and $\omega$ is obvious. 
Since $\tau(z)$ is continuous, $U := \{(z,\tau(z)) \, ;$ $z \in \partial\mathcal N\} \subset \partial{\mathcal N}\times {\bf R}_+$ has measure 0. Since $\exp_{\partial\mathcal N}(z,t)$ is continuous, $\omega = \exp_{\partial\mathcal N}(U)$ has measure 0. This implies that $\omega$ has no interior points and, since $\partial\mathcal N$ is compact, $\omega$ is compact. \qed


\begin{example}\label{ExampleB1Ellipse}
(1) Let ${\mathcal N} = B^1 = \{|x| < 1\}$ equipped with the Euclidean metric. Then $\omega = \{0\}$, which is both an ordinary boundary cut point and the first focal point.  \\
\noindent
(2) Let ${\mathcal N}$ be the inside of an ellipse : ${\mathcal N} = \{(x,y) \in {\bf R}^2 ; x^2/a^2 + y^2/b^2 < 1\}, \ (a > b > 0)$ equipped with the Euclidean metric. Then $\omega = \{(x,0) ; |x| \leq (a^2 - b^2)/a$\}. The end points $(\pm (a^2-b^2)/a,0)$ are focal points, and all the points in the open interval $\{(x,0); |x| \leq (a^2-b^2)/a\}$ are ordinary boundary cut points.
\end{example}


Based upon Lemma \ref{Bnomega}, we make the following definition.

\begin{definition}\label{Bdnormacoord}
The {\it boundary normal coordinates} is the map, 
\begin{equation}
B(\mathcal N) = \mathcal N\setminus\omega \ni x \to (z(x),s(x)) \in \partial{\mathcal N}\times{\bf R}_+,
\label{C6S5BoundaryNOrmaCoord}
\end{equation}
where $s(x)$ is the distance from $x$ to $\partial\mathcal N$ and $z(x)$ is the unique point on $\partial\mathcal N$ which is the closest to $x$, i.e.
$x = \gamma_{z(x)}( s(x))$. 
\end{definition}


\section{Boundary distance coordinates}

\subsection{Conjugate point}
The boundary cut locus is different from the standard notion of cut locus on the 
manifold without boundary.
 Therefore, we shall assume in this section that the manifold ${\mathcal N}$ is embedded 
 in a complete manifold of the same dimension $\widetilde {\mathcal N}$, 
 where $ \widetilde{\mathcal N}$ has no boundary. 
Note  that we can always construct $\widetilde{\mathcal N}$ taking it to be the Hopf 
double of $\mathcal N$ equipped with metric which is a smooth Seeley-Borel continuation 
across $\partial \mathcal N$.


\begin{definition}\label{DefineConju}
Let $c(t)$ $(a \leq t \leq b)$ be a geodesic on $\widetilde {\mathcal N}$. Two points $c(a)$ and $c(b)$ are said to be {\it conjugate} along $c(t)$ if there exists a non-trivial Jacobi field $Y(t)$ along $c(t)$ such that $Y(a) = 0, Y(b) = 0$. We also say that $c(b)$ is conjugate to $c(a)$ along $c(t)$.
\end{definition}

\medskip
For $y \in \widetilde {\mathcal N}$, let $\gamma_{(y,v)}(t) = \exp_y(tv)$ be the unit speed geodesic starting from $y$ with initial direction $v \in S_y(\widetilde {\mathcal N})$, where $S_y(\widetilde {\mathcal N}) = \{ v \in T_y(\widetilde {\mathcal N}) \, ; \, |v|_g = 1\}$. 


\begin{lemma}\label{c(t0)isconjugate}
Let $c(t) = \gamma_{(y,v)}(t)$ be a unit speed geodesic on $\widetilde {\mathcal N}$. Then $c(t_0)$ is conjugate to $y$ along $c(t)$ if and only if there exists $0 \neq \xi \in T_{y}(\widetilde {\mathcal N}) = T_{t_0v}(T_{y}(\widetilde {\mathcal N}))$ such that
$$
d\exp_{y}\Big|_{t_0v}\xi = 0.
$$
\end{lemma}

For the proof, see e.g. \cite{Au82}, p. 17, or \cite{Cha93}, Theorem 2.16.


\begin{lemma}\label{Comjulemma2}
Let $c(t) \ (a \leq t \leq b)$ be a geodesic on $\widetilde {\mathcal N}$. If there exists $a < \tau < b$ such that $c(\tau)$ is conjugate to $c(a)$ along $c(t)$, there is another 
geodesic with end points $c(a)$ and $c(b)$ which is strictly shorter than the arclength, $b-a$, of the
geodesic $c(t), \, a \leq t \leq b$.
\end{lemma}

For the proof, see e.g. \cite{Cha93}, Theorem 2.11, or \cite{KN69}, p. 87.

Similary to the boundary cut function $\tau(z)$, we introduce (Riemannian) cut function,
$\tau^R$,


\begin{definition}\label{Cutlocusdistance}
The {\it (Riemannian) cut function} $\tau^R,
 : S(\widetilde{\mathcal N}) \to {\bf R}_+$ is given by
\begin{equation}
\tau^R(y,v) = \sup_{t\geq 0}\,\left\{t\, ; \, d(\gamma_{(y,v)}(t),y) = t\right\}.
\label{Rimanniancutfunc}
\end{equation}
\end{definition}

Note that $d(y,\gamma_{(y,v)}(\tau^R(y,v)) = \tau^R(y,v)$. The point $\gamma_{(y,v)}(\tau^R(y,v))$ is called the {\it cut point} for $y$ along the geodesic $\gamma_{(y,v)}(\cdot)$. This should not be confused with the boundary cut point of Definition \ref{Definebdcutpoint}, where we considered the distance to $\partial {\mathcal N}$.

\begin{remark}
Assume that $\mathcal N=\widetilde {\mathcal N} \setminus B(x_0, a)$, where 
$B(x_0, a)$ is a ball of radius $a>0$ centered at $x_0$. Let 
$$
a < \min_{v \in S_{x_0}(\widetilde {\mathcal N})}\tau^R(x_0,v).
$$
Parametrize the points on $\partial \mathcal N=\partial B(x_0, a)$ by $v$
and observe that the normal geodesics to $\partial \mathcal N$,
i.e. $\gamma_v(t)$ are actually the continuations of the geodesics
$\gamma_{x_0,v}(t)$, namely,
$\gamma_v(t)=\gamma_{x_0,v}(t+a).$ Therefore, the focal and boundary cut
points along $\gamma_v$ are actully the conjugate and Riemannian cut
points along $\gamma_{x_0,v}$. This implies, due to Lemma \ref{Comjulemma2},
the validity of Lemma \ref{FocalpointLemma} for 
$\partial \mathcal N=\partial B(x_0, a)$.
\end{remark}


\begin{lemma}\label{TauRconti}
The mapping $\tau^{R}(y,v) : S(\widetilde {\mathcal N}) \to {\bf R}_+\cup{\infty}$ is continuous. 
\end{lemma}

This is proven in the same way as Lemma \ref{tau(z)conti}. See e.g. \cite{Cha93}, Theorem 3.1, or \cite{KN69}, p. 98.


\begin{lemma}\label{Taur>Tauz}
Let $z \in \partial {\mathcal N}$, and $\nu$ be the inner unit normal to $\partial {\mathcal N}$ at $z$. 
Then $\tau^R(z,\nu) > \tau(z)$.
\end{lemma}
Proof. Assume that for some  $z \in \partial\mathcal N$, $\tau^R(z,\nu) \leq \tau(z)$. 
Note that, following our notations for the boundary normal geodesics and geodesics 
starting at $z$, we have $\gamma_z(t) = \gamma_{(z,\nu)}(t)$ for $t > 0$.
Take $x = \gamma_{(z,\nu)}(\tau^R(z,\nu))$ and $\xi = - \dot\gamma_{(z,\nu)}(t)$ 
at $t = \tau^R(z,\nu)$.
By duality, $\tau^R(x,\xi) = \tau^R(z,\nu)$. We extend $\gamma_{(x,\xi)}(t)$ 
 on the interval $[0,\tau^R(x,\xi) + \delta]$ with $\delta > 0$. 
 Since $\dot \gamma_{(x,\xi)}(\tau^R(z,\nu))=-\nu$, by choosing 
 $\delta > 0$ small enough, we  can assume that, if 
 $\tau^R(x,\xi) < s < \tau^R(x,\xi) + \delta$, $\gamma_{(x,\xi)}(s)$ is outside 
 the original ${\mathcal N}$. Let $y(t) = \gamma_{(x,\xi)}(t+\tau^R(x,\xi))$. 
Then, for small $t$, $d(y(t),z) = d(y(t),\partial{\mathcal N}) = t$.

Note that, by the definition of $\tau^R$, for $t > 0$ 
$d(y(t),x) < t + \tau^R(x,\xi)$.
Therefore, there is a shortest geodesic $\mu(s)$ from $y(t)$ to $x$ with 
$\mu(\overline{s}) = x$ and $\overline{s} < \tau^R(x,\xi) + t$. Let $w$ be the 
last point on $\mu$ where $\mu$ crosses $\partial {\mathcal N}$. 

By triangle inequality, 
$$
\overline{s} \geq d(y(t),w) + d(\omega,x) \geq t + d(w,x) \geq t + \tau^R(z,\nu),
$$
where in the last step we use the assumption $\tau^R(z,\nu) \leq \tau(z)$.
This is a contradiction.
 \qed

\medskip
Let $z \in \partial {\mathcal N}$ and $\gamma_z$ be the boundary normal geodesic from $z$. Then, by Lemma \ref{Taur>Tauz}, there exists $\epsilon > 0$ such that for $t < \tau(z) + \epsilon$, $\gamma_z(\cdot)$ is still the shortest geodesic (lying inside ${\mathcal N}$) from $z$ to $\gamma_z(t)$.

\subsection{Hamilton's equation}
Let $(g^{ij}) = (g_{ij})^{-1}$ be the contravariant metric tensor, and define 
a $C^{\infty}$-function on $T^{\ast}(M)$ by $H(x,\xi) = \frac{1}{2}g^{ij}(x)\xi_i\xi_j$. 
As has been mentioned in Subsection 1.4 in Chap. 1, the equation of geodesic can be 
rewritten as Hamiltons's canonical equation
\begin{equation}
\left\{
\begin{split}
\frac{dx^i}{dt} & = \frac{\partial H}{\partial\xi_i} = g^{ij}(x)\xi_j, \\
\frac{d\xi_i}{dt} &= - \frac{\partial H}{\partial x^i} = 
- \frac12 \left(\frac{\partial g^{kl}(x)}{\partial x^i}\right)\xi_k\xi_l.
\end{split}
\right.
\label{C6S6:Hamilton}
\end{equation}
Fix a point $y \in  {\mathcal N}$ and let $x(t), \ \xi(t)$ be the solution 
to (\ref{C6S6:Hamilton}) with initial data $x(0) = y, \ \xi(0) = \xi_0$, 
where $\xi_0$ satisfies $g^{ij}(y)\xi_{0i}\xi_{0j} = 1$. 
Then, by the energy conservation law, 
\begin{equation}
g^{ij}(x(t))\xi_i(t)\xi_j(t) = 1.
\label{C6S6:Energy}
\end{equation}
Let $v^i(t) = dx^i(t)/dt = g^{ij}(x(t))\xi_j(t)$, and put $v(t) = (v^1(t),\cdots,v^n(t))$, 
$v_0 = v(0)$. Then $x(t)$ is a geodesic starting from $y$ with initial direction $v_0$. 
Assume that, for $U \subset S_{y}({\mathcal N}),\, 0<t_1<t_2$, 
the map : $U\times(t_1,t_2) \ni (v_0,t) \to x(t)$ is 
a diffeomorphism. Then $t$ and $v_0$ become smooth
 functions of $x$ depending (smoothly) on the parameter $y$ : $t = t(x,y)$, $v_0 = v_0(x,y)$. 
 Hence, so is $\xi = \xi(x,y)$. Since $t(x,y) = \int_y^x\xi_idx^i$, we have
\begin{equation}
\frac{\partial t(x,y)}{\partial x^i} = \xi_i(x,y).
\label{C6S6:deltatdeltax}
\end{equation}
This equality can be rewritten as
\begin{equation}
\left({\rm grad}_xt(x,y)\right)^i = g^{ij}(x)\frac{\partial t}{\partial x^j}(x,y) = v^i(x,y).
\label{C6S5:gradxt}
\end{equation}
Note also that, if $t_2< \tau^R(y, v_0)$ and $U$ is a small neighborhood of $v_0$, the
above map is, indeed, a diffeomorphism and $t(x, y)=d(x, y)$.


\subsection{Boundary distance coordinates}
Near the cut locus, we cannot use the boundary normal coordinates. However, the {\it boundary distance coordinates} constructed below can be used everywhere on $\mathcal N^{int} = {\mathcal N}\setminus\partial\mathcal N$.


\begin{lemma}\label{BdDistCoord}
For any $x_0 \in {\mathcal N}^{int}$, there exist points 
$z_1,\cdots,z_n \in \partial {\mathcal N}$ such that the functions $(\rho_1(x),\cdots,\rho_n(x))$, where $\rho_i(x) = d(x,z_i)$, give local coordinates 
in a small neighborhhood of $x_0$.
\end{lemma}
Proof. Let $z_0 \in \partial {\mathcal N}$ be a point nearest to $x_0$, i.e. 
$x_0 = \gamma_{z_0}(s_0)$, where $s_0 = d(x_0,z_0) = d(x_0,\partial {\mathcal N})$. 
If there are several such points, one can take any of them. Let 
$v_0 = - \dot\gamma_{z_0}(t)|_{t=t_0} \in S_{x_0}(\mathcal N)$ 
so that $\gamma_{(x_0,v_0)}(s_0) = z_0$. 
 By Lemma \ref{Taur>Tauz}, we have $s_0 < \tau^R(z_0,\nu(z_0)) = \tau^R(x_0,v_0)$. By Lemma \ref{c(t0)isconjugate}, $d\exp_{x_0}\big|_{s_0v_0} : T_{s_0v_0}(T_{x_0}({\mathcal N}))= T_{x_0}(\mathcal N) \to T_{z_0}({\mathcal N})$ is non-singular.

Consider curves $z_i(t), i = 1,\cdots,n-1$, in $\partial {\mathcal N}$ such that 
$z_i(0) = z_0$ and the vectors $\dot z_i(0)$, $i = 1, \cdots, n-1$, form an 
orthonormal basis of $T_{z_0}(\partial {\mathcal N})$. Let $v_i = (d\exp_{x_0}\big|_{s_0v_0})^{-1}\dot z_i(0)$ for $i = 1,\cdots,n-1,$ and 
$v_n = v_0$. Then $v_i$, $i = 1,\cdots,n$, form a basis of 
$T_{x_0}({\mathcal N})$. Furthermore, 
$c_i(s) := (\exp_{x_0})^{-1}(z_i(s)) \in T_{x_0}({\mathcal N})$, $i = 1,\cdots,n-1$, 
satisfy 
$c_i(0) = s_0v_0$ and $\dot c_i(0) = v_i$. For $i = 1,\cdots,n-1$, let 
$z_i = z_i(\epsilon)$ for a sufficiently small $\epsilon$ and $z_n = z_0$. We 
define $\rho_i(x) = d(x,z_i)$, $i = 1,\cdots,n$. Then, by (\ref{C6S5:gradxt}), 
${\rm grad}_x\rho_i(x_0) = - \dot c_i(\epsilon)/|\dot c_i(\epsilon)|_g$, $i = 1,\dots,n$, 
are linearly independent. The inverse function theorem completes the proof. \qed


\begin{example}\label{BallBddist}
Let ${\mathcal N}$ be a Euclidean sphere : ${\mathcal N} = \{|x| \leq 1\}$. Then the 
boundary normal coordinates are essentially polar coordinates with center at 
the with $r \rightarrow 1-r,\, r \leq 1$ . 
The center is the cut locus. To define the local coordinate around the origin, we have only to take $n$ points $w_1, \cdots, w_n$ on $\partial {\mathcal N}$ which are linearly independent, and $\rho_i(x) = |x-w_i|$.
\end{example}


\subsection{Reconstruction of the metric}
The following lemma is a key trick to reconstruct the Riemannian metric.


\begin{lemma}\label{tesorrecover}
Let 
 $x_0 \in {\mathcal N}$. Then we can recover the metric tensor $g_{ij}(x)$ 
 from the boundary distance functions $\partial {\mathcal N} \ni w \to d(x,w)$,
 where $x \in U$, $U$ being a neighborhood of $x_0$.
\end{lemma}
Proof.
For $x_0 \in {\mathcal N}$, let $z_0 \in \partial {\mathcal N}$ be such that $d(x_0,z_0) = d(x_0,\partial {\mathcal N})$. Then there is a small open cone of directions $C \subset S_{x_0}({\mathcal N})$ such that the geodesic starting from $x_0$ with initial direction in $C$ hits  $\partial {\mathcal N}$ transversally in a neighborhood $W_0$ of $z_0$.
Using the proof of Lemma \ref{BdDistCoord}, this means that
the directions of the shortest geodesics from $z \in W_0$ to $x_0$ form the cone $-C$ in $S_{x_0}({\mathcal N})$. 

Let $U$ be a small neighborhood of $x_0$.  For $x \in U$ and $z \in W_0$, we consier $d(x,z)$. Passing to Hamilton's equation, we have $d(x,z) = t(x,z)$, where $t(x,z)$ is defined in Subsection 7.2. By (\ref{C6S6:Energy}), we have
$$
g^{ij}(x_0)\xi_i(x_0,z)\xi_j(x_0,z) = 1.
$$
We can compute $\xi_i(x_0,z)$ from (\ref{C6S6:deltatdeltax}): $\displaystyle{\xi_i(x_0,z) = \frac{\partial d}{\partial x^i}(x_0,z)}$. Let $z$ vary on $W_0$. Then, since $\xi(x_0,z)$ varies over an open set in $S^*_{x_0}({\mathcal N})$, the unit sphere in the cotangent space $T^{\ast}_{x_0}(\mathcal N)$, we can recover the contravariant metric tensor $g^{ij}(x_0)$. \qed


\section{Reconstruction of $R({\mathcal N})$ from BSP}
In this section, we shall prove that if two manifolds ${\mathcal N}^{(1)}$ and ${\mathcal N}^{(2)}$ have the same BSP, the space of boundary distance functions $R({\mathcal N}^{(1)})$ and $R({\mathcal N}^{(2)})$ coincide. We use the expression "BSP determines the quantity $A$" to mean the following: Let $A^{(1)}$ and $A^{(2)}$ be the quantities associated to the manifolds ${\mathcal N}^{(1)}$ and ${\mathcal N}^{(2)}$, respectively. Then if ${\mathcal N}^{(1)}$ and ${\mathcal N}^{(2)}$ have the same BSP,  $A^{(1)} = A^{(2)}$ holds.


\subsection{Projection to the domain of influence}
Recall that, for a subset $\Gamma \subset \partial {\mathcal N}\subset \mathcal N$ and $\tau > 0$, we put
\begin{equation}
{\mathcal N}(\Gamma,\tau) = \{x \in \mathcal N\, ; \, d(x,\Gamma) \leq \tau\}.
\nonumber
\end{equation}
We also define for $z \in \partial {\mathcal N}$
\begin{equation}
{\mathcal N}(z,\tau) = \{x \in {\mathcal N}\,;\, d(x,z) \leq \tau\}.
\nonumber
\end{equation}
Let $\chi_{{\mathcal N}(\Gamma,\tau)}(x)$ be the characteristic function of ${\mathcal N}(\Gamma,\tau)$. We define a projection on $L^2({\mathcal N})$ by
\begin{equation}
P_{\Gamma,\tau}f(x) = \chi_{{\mathcal N}(\Gamma,\tau)}(x)f(x) \in L^2({\mathcal N}(\Gamma,\tau)), \quad 
f \in L^2({\mathcal N}).
\label{C6S8Projection}
\end{equation}
 Let $u^f(t)$ be the solution to IBVP (\ref{C6S2IBVPEquation}).


\begin{lemma}
Let $f \in C_0^{\infty}(\partial \mathcal N\times (0,\infty))$ and $\tau, t> 0$. Let $\Gamma \subset \partial {\mathcal N}$ be an open set. Then $BSP$ determines a sequence $f_j \in C_0^{\infty}(\Gamma \times (0,\tau))$ such that 
$u^{f_j}(t) \to P_{\Gamma,\tau}u^{f}(t)$.
\end{lemma}
Proof.  Let us recall an elementary fact on the projection in a Hilbert space $\mathcal H$. Let $P$ be a projection onto a closed subspace $S$ of $\mathcal H$. For $u \in \mathcal H$, take $v_n \in S$ such that $\lim_{n\to\infty}\|u - v_n\| = \inf_{v \in S}\|u - v\| = \|(1 - P)u\|$. Then $v_n \to Pu$.

Using Theorem 4.6, we have
\begin{equation}
\begin{split}
\|u^f(t)\|^2 - \|P_{\Gamma,\tau}u^f(t)\|^2 &= 
\|(1 - P_{\Gamma,\tau})u^f(t)\|^2 \\
&= \mathop{\rm inf}_{\eta \in C_0^{\infty}(\Gamma\times(0,\tau))}
\|u^f(t) - u^{\eta}(\tau)\|^2.
\end{split}
\label{C6S8:ProjectionInf}
\end{equation}
Noting that
\begin{equation}
\|u^f(t) - u^{\eta}(\tau)\|^2 = \|u^f(t)\|^2 - 2{\rm Re}(u^f(t),u^{\eta}(\tau)) + \|u^{\eta}(\tau)\|^2,
\nonumber
\end{equation}
one can compute the right-hand side of (\ref{C6S8:ProjectionInf}) by Corollary 2.2. We then choose a sequence $f_j \in C_0^{\infty}(\Gamma\times(0,\tau))$ which attains the infimum of (\ref{C6S8:ProjectionInf}). Then $u^{f_j}(\tau) \to P_{\Gamma,\tau}u^f(t)$. This procedure depends only on BSP. \qed


\begin{lemma} \label{BSPdetermines123}
Let $f, h \in C_0^{\infty}(\partial {\mathcal N}\times (0,\infty))$ and $\tau_1, \tau_2, t, s > 0$. \\
\noindent
(1) Let $\Gamma_1, \Gamma_2 \subset \partial {\mathcal N}$ be open sets. Then $BSP$ determines the inner product
\begin{equation}
\left(P_{\Gamma_1,\tau_1}u^f(t),P_{\Gamma_2,\tau_2}u^h(s)\right)_{L^2({\mathcal N})}.
\nonumber
\end{equation}
(2) Let $z_1, z_2 \in \partial {\mathcal N}$. Then $BSP$ determines the inner product
\begin{equation}
\left(P_{z_1,\tau_1}u^f(t),P_{z_2,\tau_2}u^h(s)\right)_{L^2({\mathcal N})}.
\nonumber
\end{equation}
\end{lemma}
\medskip
Proof. (1) is an obvious consequence of Lemma 8.1.
Taking open sets $\Gamma_1, \Gamma_2 \subset \partial {\mathcal N}$ shrinking to $z_1, z_2 \in \partial {\mathcal N}$, and applying Lebesgue's convergence theorem, we obtain (2). \qed


\subsection{Domain of influence and $R({\mathcal N})$}
Following \cite{KKL04}, we can identify 
the boundary normal geodesic from BSP.


\begin{lemma}\label{gammaz(s)Equiv}
Let $\gamma_z(\cdot)$ be the boundary normal geodesic starting from $z \in \partial {\mathcal N}$, and $s > 0$. Then the following 3 assertions are equivalent. \\
\noindent
(1) $d(\gamma_z(s),z) = d(\gamma_z(s),\partial {\mathcal N})$. \\
\noindent
(2) For any $\epsilon > 0$  and any neighborhood $\Gamma \subset \partial {\mathcal N}$ of $z$, the interior of $\big({\mathcal N}(\Gamma,s)\setminus {\mathcal N}(\partial {\mathcal N},s-\epsilon)\big) \neq \emptyset$. \\
\noindent
(3) For any neighborhood $\Gamma \subset \partial {\mathcal N}$ of $z$, there exists $h \in C_0^{\infty}(\Gamma\times(0,s))$ such that $\|u^h(s)\| > \|P_{\partial {\mathcal N},s-\epsilon}u^h(s)\|$.
\end{lemma}
Proof. Suppose (1) holds, and consider the open ball $B_{\epsilon/2}(x_{\epsilon})$, where $x_{\epsilon} = \gamma_z(s - \epsilon/2)$. Clearly $B_{\epsilon/2}(x_{\epsilon}) \subset {\mathcal N}(\Gamma,s)$. Let us show $B_{\epsilon/2}(x_{\epsilon}) \cap {\mathcal N}(\partial {\mathcal N},s-\epsilon) = \emptyset$. Indeed, if there exists $x \in B_{\epsilon/2}(x_{\epsilon}) \cap {\mathcal N}(\partial {\mathcal N},s-\epsilon)$, Then 
$$
d(x_{\epsilon},\partial\mathcal N) \leq d(x_{\epsilon},x) + d(x,\partial\mathcal N) < \epsilon/2 + (s-\epsilon) = s-\epsilon/2,
$$
which contradicts (1).
Hence (2) holds. 

Suppose (2) holds. Take a sequence $\epsilon_n \to 0$ and a neighborhood $\Gamma_n \subset \partial {\mathcal N}$ of $z$ of $\hbox{diam}\,(\Gamma_n) < \epsilon_n$. There exists a sequence $x_n, \delta_n \in (0,\epsilon_n/2)$ such that $B_{\delta_n}(x_n) \subset {\mathcal N}(\Gamma_n,s)\setminus {\mathcal N}(\partial {\mathcal N},s-\epsilon_n)$. Up to taking a subsequence, $x_n \to \overline{x} \in {\mathcal N}$. Since $s - \epsilon_n < d(x_n,\partial {\mathcal N}) \leq d(x_n,\Gamma_n) \leq s$, we have $d(\overline{x},\partial {\mathcal N}) = d(\overline{x},z) = s$. This implies that $\overline{x} = \gamma_{z}(s)$, hence (1) holds.

Suppose (2) holds. Let $\chi$ be the characteristic function of ${\mathcal N}(\Gamma,s)\setminus {\mathcal N}(\partial {\mathcal N},s-\epsilon)$.  Then $\|\chi\|_{L^2(\mathcal N)} > 0$. Approximating $\chi$ by $u^h(s)$, where $h \in C_0^{\infty}(\Gamma\times(0,s))$, we get (3). 

Evidently, (3) implies (2). \qed


\begin{lemma}\label{C6S6gammaw(s)Lemma}
Let $\gamma_w(\cdot)$ be the boundary normal geodesic starting from $w \in \partial {\mathcal N}$, and $s > 0$ be such that $d(\gamma_w(s),w) = d(\gamma_w(s),\partial {\mathcal N})$. Let $z \in \partial {\mathcal N}$ and $t > 0$. Then the following 3 assertions are equivalent. \\
\noindent
(1) $t > d(\gamma_w(s),z)$. \\
\noindent
(2) There exist a neighborhood $\Gamma \subset \partial\mathcal N$ of $w$ and  $\epsilon > 0$ such that
$$
{\mathcal N}(\Gamma,s) \subset {\mathcal N}(\partial {\mathcal N},s-\epsilon)\cup {\mathcal N}(z,t-\epsilon).
$$
(3)  There exist a neighborhood $\Gamma \subset \partial\mathcal N$ of $w$ and $\epsilon > 0$ such that for any $h \in C_0^{\infty}(\Gamma\times(0,s))$
$$
\|u^h(s)\|^2 = \|P_{\partial {\mathcal N},s-\epsilon}u^h(s)\|^2 + \|P_{z,t-\epsilon}u^h(s)\|^2 - (P_{\partial {\mathcal N},s-\epsilon}u^h(s),P_{z,t-\epsilon}u^h(s)).
$$
\end{lemma}
Proof.
Assume (1) holds.  
If (2) does not hold, there exist a sequence $\Gamma_n \subset \partial {\mathcal N}$ shrinking to $\{w\}$ and $\epsilon_n \to 0$, such that ${\mathcal N}(\Gamma_n,s) \not\subset {\mathcal N}(\partial {\mathcal N},s-\epsilon_n)\cup {\mathcal N}(z,t-\epsilon_n)$. Then there exists $x_n \in {\mathcal N}$ such that $d(x_n,\partial {\mathcal N}) > s - \epsilon_n$, $d(x_n,z) > t-\epsilon_n$, and  $d(x_n,\Gamma_n) \leq s$. Then, up to subsequence, $x_n \to \overline{x}$, with $d(\overline{x},\partial {\mathcal N}) = d(\overline{x},w)=s$, and $d(\overline{x},z) \geq t$. Therefore $\overline{x} = \gamma_w(s)$, which by (1) implies $d(\gamma_w(s),z) = d(\overline{x},z) < t$. This contradiction shows that (1) implies (2).

Suppose (2) holds. Since the condition $d(\gamma_w(s),w) = d(\gamma_w(s),\partial {\mathcal N})$ implies that $\gamma_{w}(s) \not\in {\mathcal N}(\partial {\mathcal N},s-\epsilon)$, then $\gamma_{w}(s) \in  {\mathcal N}(z,t-\epsilon)$. Thus, $d(\gamma_{w}(s),z) \leq t-\epsilon$, proving (1).

Let $P = P_{\partial {\mathcal N},s-\epsilon}$, $Q = P_{z,t-\epsilon}$. Using (\ref{C6S8Projection}), we see that $R = P + Q - PQ$ is a projection onto $L^2(\mathcal N(\partial {\mathcal N},s-\epsilon)\cup {\mathcal N}(z,t-\epsilon))$. Then (2) is equivalent to
$$
u^h(s) = Ru^h(s), \quad \forall h \in C_0^{\infty}(\Gamma\times(0,s)).
$$
Since $R$ is a projection, this is equivalent to 
$$
\|u^h(s) \|^2 = \|Ru^h(s)\|^2, \quad \forall h \in C_0^{\infty}(\Gamma\times(0,s)).
$$
which is equivalent to (3). \qed


\subsection{Main theorem} 
We are now in a position to prove the following 
theorem.


\begin{theorem} \label{main} 
Let $({\mathcal N},g)$ be a connected Riemannian manifold with compact boundary.
Suppose we are given the boundary spectral projections of the Neumann Laplacian on ${\mathcal N}$. Then these data determine $({\mathcal N}, g)$ uniquely.
\end{theorem}
Proof. We take $w \in \partial {\mathcal N}$. By Lemma \ref{BSPdetermines123} and Lemma \ref{gammaz(s)Equiv} (3), we can determine, by using BSP,  whether or not $\gamma_w([0,s])$ is a shortest geodesic to $\partial {\mathcal N}$. In particular, this detemines the boundary cut function $\tau(w)$. 

By Lemma \ref{C6S6gammaw(s)Lemma}, for $s \leq \tau(w)$, we can compute, by using BSP, $d(\gamma_w(s),z)$ for any $z \in \partial {\mathcal N}$. 
Thus, for any $w \in \partial\mathcal N$ and $s \leq \tau(w)$, we associate, using BSD, a function $r^{(w,s)}(\cdot) \in C(\partial\mathcal N)$:
$$
 r^{(w,s)}(z) = d(\gamma_w(s),z), \quad z \in\partial\mathcal N.
$$
Note, see (\ref{C6S5DefineBoundaryDistFunc}), that $r^{(w,s)}(\cdot)$ is the boundary distance function corresponding to $x=\gamma_{w}(s)$.

Lemma \ref{Bnomega} 
shows that, when $w$ runs over $\partial {\mathcal N}$ and $s$ runs over $[0,\tau(w)]$, then $ r^{(w,s)}(z)$ runs over the whole $R(\mathcal N) \subset C(\partial\mathcal N)$. Thus, BSP determines $R(\mathcal N)$.

We then recover the topology of ${\mathcal N}$  by Lemma 5.1. 
By Lemma 7.10, we recover the metric by BSP. \qed

\medskip
We note that the uniqueness in the above Theorem means "up to an isometry". 
We have used the generalized Fourier transform to represent BSP. However, in the above proof, we have actually used the hyperbolic Neumann-to-Dirichlet map and this can be controlled under milder assumptions. In fact, the BC-method also works for 
the manifold of bounded geometry, i.e. with the assumption of uniform injective radius of Riemannian normal coordinates, and the boundedness of curvature tensor. See \cite{KKL04}.


\section{Wave fronts and $R({\mathcal N})$}
As has been seen above, the construction of boundary distance functions from BSP is the step where the geodesic is traced using Blagovestchenski  identity for the solutions to IBVP, providing an interplay between geometry and  partial differential equations. 
Therefore, it is of interest to try other ideas. In this section, we explain the method which deals with the wave front of solution $u^f(t)$ to IBVP  (\ref{eq:IBVPsect2}).

\medskip
\noindent
(i) {\it Controlled subspaces}. By the finite propagation property, we have 
\begin{equation}
{\rm supp}\,u^f(\cdot,t) \subset {\mathcal N}(\Gamma,t) := \{x \in {\mathcal N}\, ; \, d\,(x,\Gamma) \leq t\}.
\nonumber
\end{equation}
Recall that the closure in $L^2({\mathcal N})$ of $\{u^f(\cdot,t)\,;\,f \in C_0^{\infty}(\Gamma\times(0,t)\}$ is $L^2({\mathcal N}(\Gamma,t))$.

We define a unitary operator
\begin{equation}
\mathcal F = (\mathcal F_c^{(+)},\mathcal F_p) : L^2({\mathcal N}) \to L^2((0,\infty);{\bf h};dk)\oplus {\bf C}^d,
\nonumber
\end{equation}
where $\mathcal F_c^{(+)}$ is the generalized Fourier transform, and $\mathcal F_p$ is the spectral representation associated with the point spectrum for $H$:
\begin{equation}
\mathcal F_p : L^2({\mathcal N}) \ni u = \sum_ia_i\varphi_i(x) \to (a_1,a_2,\cdots) \in {\bf C}^d,
\nonumber
\end{equation} 
where $d$ is the dimension of the point spectral subspace of $H$. If $d = \infty$, ${\bf C}^d = l^2$.
If $\mathcal N$ is compact, $\mathcal F_c^{(+)}$ is absent.

\medskip
\noindent
(ii) {\it Projections}.
Let $P_{\Gamma,t}$ be the orthogonal projection
\begin{equation}
P_{\Gamma,t} : L^2({\mathcal N}) \ni u \to \chi_{{\mathcal N}(\Gamma,t)}(x)u(x) \in L^2({\mathcal N}(\Gamma,t)),
\nonumber
\end{equation}
$\chi_{{\mathcal N}(\Gamma,t)}(x)$ being the characteristsic function of the set ${\mathcal N}(\Gamma,t)$. Passing to the Fourier transform, we have
\begin{equation}
\mathcal FP_{\Gamma,t} = \mathcal P_{\Gamma,t}\mathcal F,
\nonumber
\end{equation}
where $\mathcal P_{\Gamma,t}$ is the orthogonal projection : 
\begin{equation}
\mathcal P_{\Gamma,t} : L^2((0,\infty);{\bf h};dk)\oplus {\bf C}^d \to \mathcal L^2(\Gamma,t).
\nonumber
\end{equation}
(iii) {\it Layers}. 
It is obvious that 
\begin{equation}
L^2({\mathcal N}(\Gamma,t_-)) \subset L^2({\mathcal N}(\Gamma,t_+)), \quad 
0 \leq t_- < t_+,
\nonumber
\end{equation}
\begin{equation}
\mathcal L^2(\Gamma,t_-) \subset \mathcal L^2(\Gamma,t_+), \quad 
0 \leq t_- < t_+.
\nonumber
\end{equation}
Take $\mathcal L^2(\Gamma,t_+,t_-) = \mathcal L^2(\Gamma,t_+)\ominus \mathcal L^2(\Gamma,t_-)$, which are the Fourier transforms of functions with 
support in the {\it shell type layer} or {\it approximate wave front}
\begin{equation}
{\mathcal N}(\Gamma,t^+)\setminus {\mathcal N}(\Gamma,t^-): =
{\mathcal Sh}(\Gamma,t^+,t^-).
\nonumber
\end{equation}
Take $(\Gamma_1,t_1^+,t_1^-)$ and $(\Gamma_2,t_2^+,t_2^-)$. Then 
\begin{equation}
\begin{split}
& \mathcal L^2(\Gamma_1,t_1^+,t_1^-)\cap \mathcal L^2(\Gamma_2,t_2^+,t_2^-) \\
& = \mathcal F\{a \, ; \, {\rm supp}\,a \subset {\mathcal Sh}(\Gamma_1,t_1^+,t_1^-)\cap {\mathcal Sh}(\Gamma_2,t_2^+,t_2^-)\}.
\end{split}
\label{eq:l2gamma1ticapl2Gamma2t2}
\end{equation}

\medskip
\noindent
(iv) {\it Approximate distance functions}. We take $\Gamma_i, \ t_i^{\pm},\ i = 1, \cdots, N$, and consider $\cap_{i=1}^N\mathcal L^2(\Gamma_i,t_1^+,t_1^-)$, which is the Fourier image of functions with support in the intersection of layers. If the intersection of layers has measure $0$, then  $\cap_{i=1}^N\mathcal L^2(\Gamma_i,t_1^+,t_1^-) = \{0\}$. If this intersection has positive measure, then $\hbox{dim}\,\big( \cap_{i=1}^N\mathcal L^2(\Gamma_i,t_1^+,t_1^-)\big) =\infty$. In particular, there is $x \in {\mathcal N}$ such that $t_i^- \leq d(x,\Gamma_i) \leq t_i^+$. 

 Divide $\partial {\mathcal N}$ into a large number, which is denoted by $N(\epsilon)$, of $\Gamma_i$ with ${\rm diam}\,\Gamma_i < \epsilon$. For any vector ${\bf n} = (n_1,\cdots,n_{N(\epsilon)}) \in {\bf Z}_+^{N(\epsilon)}$, put $t_i^- = (n_i - 1)\epsilon$, $t_i^+ = n_i\epsilon$. 
Construct $\cap_{i}\mathcal L^2(\Gamma_i,t_i^+,t_i^-)$. We call $\bf n$ admissible, if $\cap_{i}\mathcal L^2(\Gamma_i,t_i^+,t_i^-)\neq \{0\}$. 
For any admissible $\bf n$, we associate a function 
$$
\kappa_{\bf n}  \in L^{\infty}(\partial {\mathcal N}), \quad
\kappa_{\bf n}(z) =n_i\epsilon, \quad {\rm for} \quad z \in \Gamma_i.
$$ 
 Take all these $\kappa_{\bf n}(z)$ for all admissible $\bf n$, and get a finite number of $L^{\infty}(\partial {\mathcal N})$ functions. They are roughly distances from various points in ${\mathcal N}$ to $\partial {\mathcal N}$. Let us denote the set of these functions as $R^{\epsilon}({\mathcal N})$.
 
\medskip
\noindent
(v) {\it Boundary distance representation of ${\mathcal N}$}. Recall that, see \S 5.1, for any $x \in {\mathcal N}$, there is the boundary distance function $r_x(z), \ z \in \partial {\mathcal N}$, 
\begin{equation}
r_x(z) = d(x,z).
\nonumber
\end{equation}
This defines the map
\begin{equation}
R : {\mathcal N} \to C^{0,1}(\partial {\mathcal N}) \subset L^{\infty}(\partial {\mathcal N}), \quad R(x)=r_x(\cdot).
\nonumber
\end{equation}
Let $R({\mathcal N})$ be the image of ${\mathcal N}$ by this map. Then
the Hausdorff distance
in $L^\infty(\partial {\mathcal N})$ between $R({\mathcal N})$ and $R^{\epsilon}({\mathcal N})$ is estimated as
\begin{equation}
d_H(R({\mathcal N}),R^{\epsilon}({\mathcal N})) < 3\epsilon.
\label{eq:Hausdorffdistance3epsilon}
\end{equation}
In fact, since $(n_i - 1)\epsilon \leq d(x,\Gamma_i) \leq n_i\epsilon$ and ${\rm diam}\,\Gamma_i \leq \epsilon$, we have
$$
|d(x,z) - n_i\epsilon| \leq 2\epsilon, \quad z \in \Gamma_i,
$$
for all $x \in \cap {\mathcal Sh}(\Gamma_i,n_i \epsilon, (n_i-1)\epsilon)$.
As, for any $x \in {\mathcal N}$, there is $\tilde x \in \cap {\mathcal Sh}(\Gamma_i,n_i \epsilon, (n_i-1)\epsilon)$
with $d(x, \tilde x) < \epsilon$, this
 proves (\ref{eq:Hausdorffdistance3epsilon}).

\medskip
In summary, we have shown the following lemma.


\begin{lemma} For any $\epsilon > 0$, we can construct, from BSP, a finite set $R^{\epsilon}({\mathcal N}) \subset L^{\infty}(\partial {\mathcal N})$, such that $d_H(R({\mathcal N}),R^{\infty}({\mathcal N})) < 3\epsilon$.
Taking $\epsilon \to 0$, we obtain the {\it boundary distance representation}
$R({\mathcal N})$ of ${\mathcal N}$.
\end{lemma}


\section{Propagation of singularities and $R({\mathcal N})$}
The singularities of solutions to the wave equation on Riemannian manifolds propagate along the geodesics. Using this property, we can determine the boundary distance function from BSP. The tool we use is the Gaussian beams which are complex valued asymptotic solutions to the wave equation in $\mathcal N \times{\bf R}$ having the following property: A Gaussian beam is concentrated near a light ray $(\gamma(t),t)$, where $\gamma(t)$ is a unit speed geodesic. For any $t$, the profile of the Gaussian beam is close to Gaussian, 
with its peak at $x = \gamma(t)$. Therefore, it is a wave packet moving along the geodesic. Since whole procedure requires  long computations, we only give the sketch here. The details can be found in \cite{KKL01}. The exposition of \cite{Ral82} is a good introduction to the theory of Gaussian beams.

 The {\it Gaussian beam} is an asymptotic solution to the wave equation of the form
\begin{equation}
U_{\epsilon}(x,t) = (\pi\epsilon)^{-n/4}\exp\left(-\frac{\theta(x,t)}{i\epsilon}\right)
\sum_{j=0}^{\infty}(i\epsilon)^ju_j(x,t),
\label{C6S10GausssianBeam}
\end{equation}
where the phase function has the following property: 
\begin{equation}
{\rm Im}\,\theta(\gamma(t),t) = 0, \quad
 {\rm Im}\,\theta(x,t) \geq C_0d(x,\gamma(t))^2,
\label{C6S10thetaandgeo}
\end{equation}
where $\gamma(t)$ is a geodesic associated with $U_{\epsilon}$.
 The fact that $U_{\epsilon}$ is an asymptotic solution means that, if we take a finite sum,
$$
U_{\epsilon}^{(N)}(x,t) = (\pi\epsilon)^{-n/4}\exp\left(-\frac{\theta(x,t)}{i\epsilon}\right)
\sum_{j=0}^{N}(i\epsilon)^ju_j(x,t),
$$
then, for any given time interval $[0,T]$, there exists a constant $C_T > 0$ such that $U_{\epsilon}^{(N)}(x,t)$ satisfies
\begin{equation}
\left|(\partial_t^2 - \Delta_g)U_{\epsilon}^{(N)}(x,t)\right| \leq C_T\epsilon^{\alpha(N)}, \quad
{\rm on} \quad \mathcal N\times [0,T],
\label{C6S10AsumptSol}
\end{equation}
$$
\alpha(N) \to \infty, \quad {\rm for}\quad N \to \infty.
$$

Fixing  boundary normal coordinates, we consider in the half-space ${\bf R}^n_+ = \{x = (z,x_n)\, ; \, z \in {\bf R}^{n-1}, x_n > 0\}$. For $z_0 \in {\bf R}^{n-1}$ and $t_0 > 0$, and we put the following highly oscillatory data on the boundary:
\begin{equation}
f_{\epsilon}(z,t) = (\pi\epsilon)^{-n/4}\chi_0(z,t)\exp\left(-\frac{\Theta(z,t)}{i\epsilon}\right),
\label{C6S10BoundaryData}
\end{equation}
where $\epsilon > 0$ is a small parameter, $\chi_0(z,t)$ is a smooth cut-off function near $(z_0,t_0)$ and
\begin{equation}
\Theta(z,t) = - ( t - t_0) + \frac{1}{2}\big(H_0(z-z_0),z-z_0)\big) + \frac{i}{2}(t-t_0)^2,
\label{C6S10InitialPhase}
\end{equation}
$(\;,\;)$ being the Euclidean inner product, $H_0$ a complex symmetric matrix with a positive definite imaginary part.

Since we are taking boundary normal coordinates, the Riemannian metric becomes $ds^2 = g_{ij}(x)dz^idz^j + (dx^n)^2$, and the boundary normal geodesic emanating from $z_0$ at time $t = t_0$ is $\gamma_{z_0}(t) = (z_0,t-t_0)$. Then for any given $z_0, t_0, H_0$ and $V$, one can construct the Gaussian beam (\ref{C6S10GausssianBeam}) as follows: 

\medskip
\noindent
(i) Let $l(z_0)$ be the time when the normal geodesic starting from $z_0$ at time $0$ hits the boundary. Then the Gaussian beam is constructed on the time interval $I(z_0) = [0, t_0 + l(z_0))$.

\noindent
(ii) It concentrates along the geodesic $\gamma_{z_0}(t) = (z_0,t - t_0)$, i.e. (\ref{C6S10thetaandgeo}) is satisfield for $\gamma(t) = \gamma_{z_0}(t)$ on $I(z_0)$.

\noindent
(iii) Its phase function and the amplitude functions satisfy
\begin{equation}
\theta(z,0,t) \approx \Theta(z,0), \quad
u_j(z,0,t) \approx \delta_{j0},
\nonumber
\end{equation}
where $f(z) \approx g(z)$ means $\partial_z^{\alpha}(f(z) - g(z)) = 0, \ \forall \alpha$, at $z = z_0$, and 
\begin{equation}
(\partial_t\theta)^2 - g_{ij}(x)(\partial_i\theta)(\partial_j\theta) \asymp 0,
\nonumber
\end{equation} 
\begin{equation}
L_{\theta}u_n \asymp (\partial_t^2 - \Delta_g)u_{n-1}, \quad u_{-1} = 0,
\nonumber
\end{equation}
where $L_{\theta} = 2(\partial_t\theta)\partial_t - 2g^{ij}(\partial_i\theta)\partial_j + (\partial_t^2 - \Delta_g)\theta$, $\partial_j = \partial/\partial x^j$, and $f(x) \asymp g(x)$ means
$\partial_x^{\alpha}(f(x) - g(x)) = 0$, $ \forall \alpha$, at $x =  \gamma_{z_0}(t)$ on $I(z_0)$.

\medskip
Let $u_{\epsilon}(t)$ be the solution to IBVP (\ref{eq:IBVPsect2}) with $f$ repalced by $f_{\epsilon}(z,t)$ of (\ref{C6S10BoundaryData}). Then as can be checked easily
\begin{equation}
\|u_{\epsilon}(t) - U_{\epsilon}^{(N)}(t)\| \leq C_N\epsilon^{\alpha(N)}.
\nonumber
\end{equation}
Using this Gaussian beam one can prove the following lemma (see Corollary 3.25 of \cite{KKL01}).

\begin{lemma}
For any $z_0 \in \partial {\mathcal N}$, $t_0 < t < t_0 + l(z_0)$ and $\tau > 0$, we have
\begin{equation}
\lim_{\epsilon\to0}\left(P_{y,\tau}u_{\epsilon}(t),u_{\epsilon}(t)\right) = 
\left\{
\begin{split}
 \alpha(t), \quad &{\rm if} \quad d(\gamma_{z_0}(t),y) < \tau, \\
0, \quad &{\rm if} \quad d(\gamma_{z_0}(t),y) > \tau,
\end{split}
\right.
\nonumber
\end{equation}
where $\alpha(t) > 0$.
\end{lemma}

Therefore we can compute $d(\gamma_{z_0}(t),y)$ from BSP.


\section{Eigenfunction coordinates}

\subsection{Regularity of the metric}
Let us discuss regularity problems for the metric. For the details, see \cite{AKKLT04}. 
If $g_{ij} \in C^{k,\alpha}$, the distance is locally $C^{k-1,\alpha}$. Then $g_{ij}$ is only
in distance coordinates is $C^{k-2,\alpha}$, since the Jacobian is involved. As regard 
to this regularity loss problem, a nice choice is the 
{\it harmonic coordinates} 
$X^i(x)$, $i = 1, \cdots,n$, such that $\Delta_g X^i = 0$. The feature of these 
harmonic coordinates is that they are the best possible for smoothness. In fact, assume 
that, in some coordinates $(x^1,\cdots,x^n)$, $g_{ij}$ is $C^{k,\alpha}$. Then 
$X^j(x),\, j=1, \dots, n,$ are $C^{k+1,\alpha}$, which implies that $g_{ij}$ 
is $C^{k,\alpha}$ in the coordinates $(X^1,\cdots,X^n)$. Another important feature 
is that, in the harmonic coordinates, the following equation holds:
\begin{equation}
\Delta_gg_{ij} = -2{\rm Ric}_{ij} + \mathcal F_{ij}(g,\nabla g),
\nonumber
\end{equation}
where ${\rm Ric}_{ij}$ is the Ricci curvature.
For the proof, see \cite{DeKa81}, Lemma 4.1. See also \cite{Heb96} for 
harmonic 
coordinates. 

We should also remark that eigenfunctions of $\Delta_g$ are good candidates of coordinates. In this section, we only consider the case of compact manifold. 


\begin{lemma}
Let $\varphi_j(x), \ j = 1, 2, \cdots,$ be a complete orthonormal system of 
eigenfunctions 
of $\Delta_g$ with Neumann boundary condition. Then, for any 
$x_0 \in \mathcal N^{int}$, there exists a neighborhood 
of $x_0$ and $j_1, \cdots, j_n$ such that $\varphi_{j_1}(x), \cdots, \varphi_{j_n}(x)$ 
form local coordinates on $U$.
\end{lemma}
Proof. By the Fourier expansion for any $a \in C_0^{\infty}({\mathcal N})$, 
$a(x) = \sum a_k\varphi_k(x)$, where the series converges in $C^{\infty}(\mathcal N)$. 
From this one 
can show that, for any $x_0 \in \mathcal N^{int}$,
${\rm Sp}\{\nabla\varphi_k(x_0)\}_{k=1}^{\infty} =T_{x_0}( \mathcal N):={\bf R}^n$, 
where 
${\rm Sp}(A)$ means the linear span of the set $A$. 
In fact, 
take some local coordinates near
$x_0$  and let  $a(x)$  be a smooth function which is 
linear around $x_0$. Then $\nabla a(x) = \sum a_k\nabla \varphi_k(x)$ near $x_0$. 
This 
means that the direction $\nabla a(x_0)$ is approximated by a linear 
combination 
of $\nabla\varphi_k(x_0)$. Therefore, one can choose $n$ 
functions 
$\varphi_{j_i}(x), i = 1, \cdots, n$, such that 
${\rm Sp}\{\nabla\varphi_k(x_0) ; k = j_1,\cdots, j_n\} = {\bf R}^n$. \qed

\medskip
Note that, since $\Delta_g\varphi_k = \lambda_k\varphi_k$, we have,
by elliptic regularity, that $\varphi_k \in C^{k+1,\alpha}$ if $g_{ij} \in C^{k,\alpha}$.

Suppose we can find $\varphi_k(x), k = 1, 2, \cdots$, in $R({\mathcal N})$. Then, we can reconstruct the distance on ${\mathcal N}$ by looking at the heat kernel
\begin{equation}
\begin{split}
h(x,y,t) & = \sum e^{-\lambda_kt}\varphi_k(x)\varphi_k(y).
\end{split}
\nonumber
\end{equation}
In fact, we have as $t \to 0$
\begin{equation}
h(x,y,t) \sim \frac{C_n}{t^{n/2}}e^{-\frac{d^2(x,y)}{4t} }.
\nonumber
\end{equation}
Therefore,
\begin{equation}
\left(- \lim_{t\to 0}4t\log h(x,y,t)\right)^{1/2} = d(x,y).
\nonumber
\end{equation}
This is another way of reconstructing the distance on $R({\mathcal N})$.


\subsection{Spectral map} 
From $R({\mathcal N})$, we have reconstructed the differential structure of ${\mathcal N}$ by finding boundary normal coordinates and boundary distance coordinates. However, the distance coordinates have the disadvantage that we lose 2 orders of regularity, say, of $g_{ij}$. As for the regularity problem, the best choice is the coordinate system made of eigenfunctions. Let
\begin{equation}
 \mu_1, \mu_2, \mu_3 \cdots \quad 
{\rm and} \quad \psi_1(x), \psi_2(x), \psi_3(x) \cdots
\nonumber
\end{equation}
be the eigenvalues and  eigenfunctions of Dirichlet problem, and 
\begin{equation}
\lambda_0, \lambda_1, \lambda_2, \cdots \quad {\rm and} \quad \varphi_0(x), \varphi_1(x), \varphi_2(x), \cdots
\nonumber
\end{equation}
those of Neumann problem.


\begin{lemma}
Having BSD for, say, Neumann problem, we can find BSD for Dirichlet proplem.
\end{lemma}
Proof. Let $\Delta^N$ and $\Delta^D$ be Neumann and Dirichlet Laplacians on ${\mathcal N}$, and $\{\lambda_i, \varphi_i\big|_{\partial {\mathcal N}}; i = 0,1,2,\cdots\}$ and 
$\{\mu_i, \partial\psi_i/\partial\nu\big|_{\partial {\mathcal N}}; i = 1,2,\cdots\}$ 
be the boundary spectral data for Neumann and Dirichlet problem, respectively.
Take $z \not\in \sigma(- \Delta^N)\cup\sigma(- \Delta^D)$.
The Neumann-to-Dirichlet map is defined to be $R^N(z) : f \to u\big|_{\partial M}$, where
\begin{equation}
\left\{
\begin{split}
& (- \Delta_g - z)u = 0 \quad {\rm in} \quad {\mathcal N}, \\
& \frac{\partial u}{\partial\nu} = f \quad {\rm on} \quad \partial {\mathcal N}.
\end{split}
\right.
\nonumber
\end{equation}
and the Dirichlet-to-Neumann map is defined to be $R^D(z) : f \to \partial v/\partial\nu\big|_{\partial {\mathcal N}}$, where
\begin{equation}
\left\{
\begin{split}
& (- \Delta_g - v)u = 0 \quad {\rm in} \quad {\mathcal N}, \\
& v = f \quad {\rm on} \quad \partial {\mathcal N}.
\end{split}
\right.
\nonumber
\end{equation}
As is seen before, $R^N(z)$ has an integral kernel
\begin{equation}
R^N(z;x,y) = \sum_{i=0}^{\infty}\frac{\varphi_i(x)\varphi_i(y)}{z - \lambda_i}, \quad x, y \in \partial \mathcal N.
\nonumber
\end{equation}
By definition, one can easily see that $(R^N(z))^{-1} = R^D(z)$, and $R^N(z)$ is determined by the Neumann spectral data. Therefore, $R^D(z)$ is determined by the Neumann spectral data. Now $R^D(z)$ has the following formal integral kernel
\begin{equation}
R^D(z;x,y) = \sum_{i=1}^{\infty}\frac{\partial_{\nu}\psi_i(x)\partial_{\nu}\psi_i(y)}{z - \mu_i}, \quad x, y \in \partial \mathcal N.
\nonumber
\end{equation}
Actually this sum does not converge. However, $R^D(z)$ is known to be an operator-valued meromorphic function of $z$ with simple poles at $z = \mu_i$ and its
residue is given by $\sum_{\mu_k = \mu_i}\partial_{\nu}\psi_{\mu_k}(x)\partial_{\nu}\psi_{\mu_k}(y)$, which proves the lemma. \qed

\bigskip
By the same argument as in the proof of Lemma 11.1, one can show the following lemma.


\begin{lemma}
Let $x \in \partial{\mathcal N}$. Then there are $n-1$ eigenfunctions of Neumann problem, and one eigenfunction of the Dirichlet problem such that $\{\varphi_{i_1},\cdots,\varphi_{i_{n-1}}, \psi_{i_n}\}$ form a coordinate system near $x$.
\end{lemma}

Now we define the spectral map $S : {\mathcal N} \to {\bf R}^{\infty}$ by
\begin{equation}
S(x) = \{\varphi_0(x), \psi_1(x), \varphi_1(x), \psi_2(x), \varphi_2(x), \cdots\}.
\nonumber
\end{equation}
Since these eigenfunctions satisfy $- \Delta_g\varphi_i = \lambda_i\varphi_i$, 
$- \Delta_g\psi_i = \mu_i\psi_i$, they can be used to find coefficients of 
$\Delta_g$ in "eigenfunction coordinates", i.e. the metric tensor. 
This is now an well-known idea in geometry, see e.g
\cite{BBG94}, \cite{KaKu9496}.

The problem is how to find these eigenfunction coordinates.


\begin{lemma}
BSD determines $S({\mathcal N}) \subset {\bf R}^{\infty}$.
\end{lemma}

Proof. Let us recall the slicing procedure in \S 9. There, by solving the initial boundary value problem for the 
wave equation, we have constructed a layer $Sh(\Gamma,t^+,t^-)$. By taking the intersection of these layers 
in a generic position, we can find a region of positive measure in ${\mathcal N}$. Let us call it "a pixel", and denote by 
$P_X$. Passing to the Fourier transforms $\mathcal F^N$ (Neumann case) or $\mathcal F^D$ (Dirichlet case), we then find
\begin{equation}
l^{2,N}(P_X) := \mathcal F^N(L^2(P_X)), \quad 
l^{2,D}(P_X) := \mathcal F^D(L^2(P_X)).
\nonumber 
\end{equation}
Observe that 
\begin{equation}
\mathcal F^D\psi_i = e_i = (0,\cdots,0,1,0,\cdots,0,\cdots),
\nonumber
\end{equation}
\begin{equation}
\mathcal F^N\varphi_i = f_i = (0,\cdots,0,1,0,\cdots,0,\cdots),
\nonumber
\end{equation}
Let 
\begin{equation}
\begin{split}
& Q^D(P_X) : l^2 \to l^{2,D}(P_X), \\
& Q^N(P_X) : l^2 \to l^{2,N}(P_X)
\end{split}
\nonumber
\end{equation}
be the associated orthogonal projections. We then have
\begin{equation}
(Q^D(P_X)e_i,e_j) = \int_{P_X}\psi_i(x)\psi_j(x)dV,
\nonumber
\end{equation}
\begin{equation}
(Q^N(P_X)f_0,f_0) = \frac{1}{{\rm Vol}(({\mathcal N})}\int_{P_X}dV.
\nonumber
\end{equation}
We now let $P_X$ shrink to a point : $P_X \to \{x\}$. Then we have
\begin{equation}
\frac{(Q^De_i,e_j)}{(Q^Nf_0,f_0)} \to {\rm Vol}({\mathcal N})\psi_i(x)\psi_j(x),
\nonumber
\end{equation}
\begin{equation}
\frac{(Q^Nf_i,f_0)}{(Q^Nf_0,f_0)} \to {\rm Vol}^{1/2}({\mathcal N})\varphi_i(x), \quad
{\rm Vol}^{-1/2}({\mathcal N})=\varphi_0|_{\partial {\mathcal N}} .
\nonumber
\end{equation}
We thus find a map
\begin{equation}
\widetilde S : {\mathcal N} \ni x \to \{\varphi_0(x), \psi_1(x)^2, \varphi_1(x), 
\psi_2(x)\psi_1(x), \cdots\}.
\nonumber
\end{equation}
Since $\psi_1(x) > 0$, one can find $\psi_1(x)$ from $\psi_1(x)^2$ on ${\mathcal N}$. Therefore by dividing by $\psi_1(x)$, we get $\{\varphi_0(x), \psi_1(x), \varphi_1(x), \psi_2(x), \cdots\} = S(x)$. \qed


\appendix


\chapter{Radon transform and propagation of singularities in ${\bf R}^n$}
In this appendix, we study the relation between the Radon transform and the propagation of singularties of solutions to the wave equation. This is basically well-known, however, we feel it important to study this problem in a general Riemannian metric, and it is not an obvious task even in the asymptotically Euclidean metric.

The Radon transform associated with the Euclidean metric is defined by
\begin{equation}
\left(\mathcal R_0f\right)(s,\theta) = \int_{s = x\cdot\theta}f(x)d\Pi_x, 
\quad s \in {\bf R}, \quad \theta \in S^{n-1},
\nonumber
\end{equation}
$d\Pi_x$ being the measure induced on the hyperplane $\{x \in {\bf R}^n ; \, s = x\cdot\theta\}$ from the Lebesgue measure $dx$ on ${\bf R}^n$. This is rewritten as
\begin{equation}
\left(\mathcal R_0f\right)(s,\theta) = (2\pi)^{(n-1)/2}\int_{-\infty}^{\infty}e^{isk}\widehat f(k\theta)dk,
\nonumber
\end{equation}
where $\widehat f$ is the Fourier transform:
\begin{equation}
\widehat f(\xi) = (2\pi)^{-n/2}\int_{{\bf R}^n}e^{-ix\cdot\xi}f(x)dx.
\nonumber
\end{equation}
Let us  consider the Riemannian metric on ${\bf R}^n$ satisfying the following condition:
\begin{equation}
|\partial_x^{\alpha}(g_{ij}(x) - \delta_{ij})| \leq C_{\alpha}
(1 + |x|)^{-1-\epsilon_0-|\alpha|}, \quad \forall \alpha,
\label{eq:DecayAssumption}
\end{equation}
where $\epsilon_0 > 0$ is a constant. 
In Chap. 2, \S 7, we have already constructed a generalized Fourier transformation $\mathcal F^{(\pm)}$ for $\Delta_g$. As in Chap. 2, \S 7, we construct $\mathcal F_{\pm}$ from $\mathcal F^{(\pm)}$, and 
define the modified Radon transform $\mathcal R_{\pm}$ by
$$
\mathcal R_{\pm}f(s,\theta) = \frac{1}{\sqrt{2\pi}}\int_{-\infty}^{\infty}
e^{isk}(\mathcal F_{\pm}f)(k,\theta)dk.
$$
 For the Euclidean Laplacian in ${\bf R}^n$ this turns out to be
$$
\mathcal R_{\pm} = \big(\mp \partial_s + 0)^{\frac{n-1}{2}}\mathcal R_0.
$$

The main issue of this chapter is  the {\it singular support theorem} for $\mathcal R_{\pm}$. We construct $\varphi(x,\theta) \in C^{\infty}({\bf R}^n\times S^{n-1})$ such that
$$
|\partial_{\theta}^{\alpha}\partial_{x}^{\beta}(\varphi(x,\theta) - x\cdot\theta)| \leq C_{\alpha\beta}(1 + |x|)^{-|\beta| - \epsilon_0},
$$
and it solves the eikonal equation
$$
g^{ij}(x)(\partial_i\varphi(x,\theta))(\partial_j\varphi(x,\theta)) = 1, \quad
\partial_i = \partial/\partial x_i,
$$
 in an appropriate region in ${\bf R}^n$. We put $\Sigma(s,\theta) = \{x \in {\bf R}^n ; s = \varphi(x,\theta)\}$, which describes a wave front of a plane wave solution to the wave equation $\partial_t^2u = \Delta_gu$. Then by observing the propagation of singularities, we obtain the following theorem: 
Let $\mathcal R_+(s,\theta,x)$ be the distribution kernel of $\mathcal R_+$. Then if we fix $s > 0$ large enough, we have the following singularity expansion:
\begin{equation}
\mathcal R_+(s,\theta,x) \sim \sum_{j=0}^{\infty}
(s - \varphi(x,\theta))^{-\frac{n+1}{2}+j}_-r_j(x,\theta).
\nonumber
\end{equation}
Let $\Sigma(s)$ be the envelope of the family of hypersurfaces $\{\Sigma(s,\theta)\; ; \; \theta \in S^{n-1}\}$, which describes a spherical wave front. We then show that $f$ 
 (satisfying a suitable condition on the wave front set) is piecewise smooth near $\Sigma(\sigma)$ with interface $\Sigma(\sigma)$ if and only if 
 $\big(\mathcal R_+f\big)(s)$ is piecewise smooth near $\{s = \sigma\}$ with interface $s = \sigma$. Moreover we also obtain the singularity expansion of $\mathcal R_+f$ in terms of spherical wave solution to the eikonal equation.


\section{Fourier and Radon transforms for perturbed metric}


\subsection{Spectral properties}
The  Laplace-Beltrami operator $\Delta_g$ is  symmetric in $L^2({\bf R}^n ; \sqrt{g(x)}dx)$. To avoid the denstity $\sqrt{g(x)}$, we apply a unitary transformation : $u \to ug(x)^{1/4}$, and consider the differential operator
\begin{equation}
H = - g(x)^{1/4}\Delta_gg(x)^{-1/4} = - \sum_{i,j=1}^na_{ij}(x)\partial_{i}\partial_{j} + \sum_{i=1}^nb_i(x)\partial_i + c(x)
\nonumber
\end{equation}
on $L^2({\bf R}^n;dx)$. Note that $a_{ij}(x) = g^{ij}(x)$ and $a_{ij}(x) - \delta_{ij}, b_i(x), c(x)$ satisfy
\begin{equation}
|\partial^{\alpha}_xa(x)| \leq C_{\alpha}(1 + |x|)^{-|\alpha|  - 1 - \epsilon_0}, \quad \forall \alpha.
\nonumber
\end{equation}
We put
\begin{equation}
H_0 = - \Delta = - \sum_{i=1}^n(\partial/\partial x_i)^2, \quad V = H - H_0,
\nonumber
\end{equation}
\begin{equation}
R_0(z) = (H_0 - z)^{-1}, \quad R(z) = (H - z)^{-1}.
\nonumber
\end{equation}


\begin{theorem}
(1) $\sigma(H) = \sigma_{ac}(H) = [0,\infty)$. \\
\noindent
(2) $\sigma_p(H) = \sigma_{sc}(H) = \emptyset$. \\
\noindent
(3) For any $\lambda > 0$ and $f, g \in {\mathcal B}$, there exists a limit
$$
\lim_{\epsilon\to 0}\big(R(\lambda \pm i\epsilon)f,g\big) 
=: \big(R(\lambda \pm i0)f,g\big).
$$
(4) For any $0 < a < b < \infty$, there exists a constant $C > 0$ such that
$$
\|R(\lambda \pm i0)f\|_{{\mathcal B}^{\ast}} \leq C\|f\|_{\mathcal B}, \quad 
a < \lambda < b.
$$
(5) For any $f, \, g \in {\mathcal B}$, $(R(\lambda \pm i0)f,g)$ is a continuous function of $\lambda > 0$.
\end{theorem}

The proof is omitted. The limiting absorption principle in weighted $L^2$ spaces was proved in, e.g., \cite{IkSa72}, and in $\mathcal B-\mathcal B^{\ast}$ spaces by Agmon and Agmon-H{\"o}rmander \cite{Hor}, and \cite{JePe85}.


\subsection{Generalized Fourier transform} Let us recall the notation in Chap. 2, \S 7.
For $k \in {\bf R}\setminus\{0\}$ and $f \in {\mathcal B}$, we define
\begin{equation}
\left(\mathcal F^0(k)f\right)(\omega) = (2\pi)^{-n/2}\int_{{\bf R}^n}
e^{- ik\omega\cdot x}f(x)dx. 
\nonumber
\end{equation}
It has the following properties
\begin{equation}
\mathcal F^0(k) \in {\bf B}({\mathcal B} ; L^2(S^{n-1})),
\nonumber
\end{equation}
\begin{equation}
\mathcal F^0(-k) = J\mathcal F^0(k),
\end{equation}
$J$ being the anti-podal operator defined by 
\begin{equation}
\big(J\psi\big)(\omega) = \psi(-\omega).
\label{eq:Chap4Sec1Antipodal}
\end{equation}
We put
\begin{equation}
\widehat{\mathcal H}_{>0} = L^2((0,\infty) ; L^2(S^{n-1});k^{n-1}dk),
\nonumber
\end{equation} 
\begin{equation}
\widehat{\mathcal H}_{<0} = L^2((-\infty,0) ; L^2(S^{n-1});|k|^{n-1}dk).
\nonumber
\end{equation} 
Then the operator $(\mathcal F^0f)(k) := \mathcal F^0(k)f$ is uniquely extended to a unitary operator from $L^2({\bf R}^n)$ to $\widehat{\mathcal H}_{>0}$. It is also extended to a unitary operator from $L^2({\bf R}^n)$ to $\widehat{\mathcal H}_{<0}$. With these properties in mind, we define the generalized Fourier transform $\mathcal F^{(\pm)}(k)$ by the following formula:
\begin{equation}
\mathcal F^{(\pm)}(k) = \mathcal F^0(k)\big(1 - VR((k \pm i0)^2)\big).
\nonumber
\end{equation}
Note that $(k + i0)^2 = k^2 + i0$ for $k > 0$ and $(k + i0)^2 = k^2 - i0$ for $k < 0$. By (\ref{eq:Chap4Sec1Antipodal}) we have
\begin{equation}
\mathcal F^{(+)}(-k) = J\mathcal F^{(-)}(k).
\label{eq:FplusandJFminus}
\end{equation}
By Theorem 2.7.11, $\mathcal F^{(\pm)}$ is uniquely extended to a unitary operator from $L^2({\bf R}^n)$ to $\widehat{\mathcal H}_{>0}$ and diagonalizes $H$, and $\mathcal F^{(\pm)}$ is also  unitary from $L^2({\bf R}^n)$ to $\widehat{\mathcal H}_{<0}$.

\bigskip
\noindent
{\bf Remark}. One can also prove that $\big(\mathcal F^{(\pm)}f\big)(k,\theta)$ is smooth with respect to $k$ and $\theta$. In fact, let $\varphi(\lambda) \in C_0^{\infty}((0,\infty))$, $f(x) \in C^{\infty}_0({\bf R}^n)$ and put $g(\xi) = \big(\mathcal F^{(\pm)}(k)\varphi(L)f\big)(\omega)$ with $k = |\xi|$, $\omega = \xi/|\xi|$. Then $g(\xi) \in C^{\infty}({\bf R}^n)$. For the case of the Schr{\"o}dinger operator $ - \Delta + V$ where $V$ is a real-valued potential, we have proven this property in \cite{Is85} by using a parametrics at infinity of the time evlolution equation. One can repeat the same argument by using the geometrical optics solutions to be constructed in \S3 of this chapter. 

\medskip
The following theorem is proved in the same way as in \cite{Yaf91}.


\begin{theorem}
For $k \in {\bf R}\setminus\{0\}$ and $f \in {\mathcal B}$
$$
R((k + i0)^2)f(x) \simeq C_0(k)r^{-(n-1)/2}e^{ikr}
\left({\mathcal F}^{(+)}(k)f\right)(\omega),
$$
 where $r = |x|, \ \omega = x/r$, and
$$
C_0(k) = \sqrt{\frac{\pi}{2}}(-ik + 0)^{(n-3)/2}.
$$
\end{theorem}

Here $f \simeq g$ means that
\begin{equation}
\lim_{R\to\infty}\frac{1}{R}\int_{|x|<R}|f(x) - g(x)|^2dx = 0.
\nonumber
\end{equation}


\subsection{Wave operators and scattering matrix} 
The wave operator $W_{\pm}$ for the Schr{\"o}dinger equation is defined by the following strong limit in $L^2({\bf R}^n)$:
\begin{equation}
W_{\pm} = \mathop{\rm s-lim}_{t\to\pm \infty}e^{itH}e^{-itH_0}.
\nonumber
\end{equation}
It is well-known that this limit exists and regarding $\mathcal F^0$ and $\mathcal F^{(\pm)}$ as unitary from $L^2({\bf R}^n)$ to 
$\widehat{\mathcal H}_{>0}$, we have the following relation
\begin{equation}
W_{\pm} = \big(\mathcal F^{(\pm)}\big)^{\ast}\mathcal F^0.
\label{eq:WaveOpandFourier}
\end{equation}
The wave operator for the wave equation is usually defined by the energy norm. 
We can also employ the following equivalent operator
\begin{equation}
\mathop{\rm s-lim}_{t\to\pm \infty}e^{it\sqrt{H}}e^{-it\sqrt{H_0}} = W_{\pm} = 
\big(\mathcal F^{(\pm)}\big)^{\ast}\mathcal F^0.
\label{eq:Waveopforwaveeq}
\end{equation}
The point is that the limit in the left-hand side exists, and coincides with the wave operator for the Schr{\"o}dinger equation.
This fact, called the invariance principle, is known to hold in a broad situations (see e.g. \cite{Ka76}, p. 579). The equality (\ref{eq:Waveopforwaveeq}) can of course be proved directly by using $\mathcal F^{(\pm)}$ (see e.g. \cite{Moc83}). 

As a by-product, one can show that the solution $u(t)$ of the wave equation
\begin{equation}
\left\{
\begin{split}
& \partial_t^2u = - Hu, \\
& u(0) = f, \quad \partial_tu(0) = - i\sqrt Hf
\end{split}
\right.
\nonumber
\end{equation}
behaves as follows
\begin{equation}
\|u(t) - e^{-it\sqrt{H_0}}f_{\pm}\|_{L^2} \to 0 \quad 
{\rm as} \quad t \to \pm \infty,
\nonumber
\end{equation}
where $f_{\pm} = \big(\mathcal F^0\big)^{\ast}\mathcal F^{(\pm)}f$.
Therefore $\mathcal F^{(\pm)}$ represents the far field behavior of waves. The same fact can be proven for the Schr{\"o}dinger equation.


\begin{definition}
Regarding $\mathcal F^0$ and $\mathcal F^{(\pm)}$ as unitary from $L^2({\bf R}^n)$ to 
$\widehat{\mathcal H}_{>0}$, we define the scattering operator $S$, its Fourier transform $\widehat S$, and the physical S-matrix $\widehat S_{phy}(k)$  by
$$
S = \big(W_+\big)^{\ast}W_-, \quad
\widehat S = \big(\mathcal F^0\big)^{\ast}S\mathcal F^0 = \mathcal F^{(+)}\big(\mathcal F^{(-)}\big)^{\ast}.
$$
$$
\widehat S_{phy}(k) = I - \pi ik^{n-2}\mathcal F^{(+)}(k)V\mathcal F^0(k)^{\ast}, \quad k > 0.
$$
\end{definition}


\begin{lemma} 
$\widehat S_{phy}(k)$ is unitary on  $L^2(S^{n-1})$ for any $k > 0$, and 
$$
\big(\widehat Sf\big)(k) = \widehat S_{phy}(k)f(k), \quad 
\forall f \in \widehat{\mathcal H}_{>0}, \quad {\rm a.e.} \ \  k > 0,
$$
$$
{\mathcal F}^{(+)}(k) = \widehat S_{phy}(k){\mathcal F}^{(-)}(k), \quad \forall k > 0.
$$
\end{lemma}


\begin{definition}
For $k > 0$, we define the geometric scattering matrix $\widehat S_{geo}(k)$ by
$$
\widehat S_{geo}(k) = \widehat S_{phy}(k)J.
$$
\end{definition}

The following theorem is proved in the same way as in \cite{Yaf91}, (see also \cite{Is01}, \cite{Is04a}).

\begin{theorem} Let $k > 0$, and put  
\begin{equation}
{\mathcal N}(k) = \{u \in {\mathcal B}^{\ast} ; (H - k^2)u = 0\}.
\nonumber
\end{equation}
(1) We have
\begin{equation}
\mathcal N(k) = 
{\mathcal F}^{(\pm)}(k)^{\ast}\big(L^2(S^{n-1})\big).
\nonumber
\end{equation}
(2) For any $u \in \mathcal N(k)$ there exist $\varphi_{\pm} \in L^2(S^{n-1})$ such that
\begin{equation}
u(x) \simeq \frac{e^{i(kr - (n-1)\pi/4)}}{r^{(n-1)/2}}\varphi_+(\omega) + 
\frac{e^{-i(kr - (n-1)\pi/4})}{r^{(n-1)/2}}\varphi_{-}(\omega),
\label{eq:AsympHelmholtz}
\end{equation}
where $r = |x|, \ \omega = x/r$. \\
\noindent
(3) For any $\varphi_- \in L^2(S^{n-1})$, there exist a unique $u 
\in {\mathcal N}(k)$ and $\varphi_+ \in L^2(S^{n-1})$ such that the expansion (\ref{eq:AsympHelmholtz}) holds. Moreover they are related as follows :
\begin{equation}
\varphi_+ = \widehat S_{geo}(k)\varphi_-.
\nonumber
\end{equation}
\end{theorem}


\subsection{Modified Radon transform}
It is convenient to change the definition of the generalized Fourier transform slightly. 
For $k \in {\bf R}\setminus\{0\}$, we define
\begin{equation}
\mathcal F_{\pm}(k)  = \frac{1}{\sqrt2}
(\mp ik + 0)^{(n-1)/2}\mathcal F^{(\pm)}(k),
\nonumber
\end{equation}
\begin{equation}
\mathcal F_{0}(k)  = \frac{1}{\sqrt2}
(- ik + 0)^{(n-1)/2}\mathcal F^{0}(k),
\nonumber
\end{equation}
 and put $\big(\mathcal F_{\pm}f\big)(k) = 
\mathcal F_{\pm}(k)f$, $\big(\mathcal F_{0}f\big)(k) = 
\mathcal F_{0}(k)f$. Note that by (\ref{eq:FplusandJFminus})
\begin{equation}
\mathcal F_+(-k) = J{\mathcal F}_-(k).
\label{eq:LowerplusminusF}
\end{equation}


\begin{theorem}
(1) $\mathcal F_{\pm} : L^2({\bf R}^n) \to L^2({\bf R};L^2(S^{n-1});dk)$ is an isometry. Moreover we have
$$
\left(\mathcal F_{\pm} Hf\right)(k) = k^2\left(\mathcal F_{\pm}f\right)(k).
$$
(2) For $k > 0$, we have
\begin{equation}
\mathcal F_+(k) = (-i)^{n-1}\widehat S_{phy}(k)J\mathcal F_{+}(-k).
\nonumber
\end{equation}
Consequently, the range of $\mathcal F_{\pm}$ has the following characterization: 
\begin{equation}
g \in {\rm  Ran}\,\mathcal F_{+} \  \Longleftrightarrow \ 
(-i)^{n-1}\widehat S_{phy}(k)Jg(-k) = g(k), \ k > 0,
\nonumber
\end{equation}
\begin{equation}
g \in {\rm  Ran}\,\mathcal F_{-} \  \Longleftrightarrow \ 
(-i)^{n-1}\widehat S_{phy}(k)g(k) = Jg(-k), \ k > 0.
\nonumber
\end{equation}
(3) Let $r_{+}$ $(r_-)$ be the projection onto $\widehat{\mathcal H}_{>0}$ $(\widehat{\mathcal H}_{<0})$. Then we have
\begin{equation}
W_+ = 2\big(\mathcal F_+\big)^{\ast}r_{+}\mathcal F_{0}, \quad 
W_- =  2\big(\mathcal F_+\big)^{\ast}r_{-}\mathcal F_{0},
\label{eq:WpluminusFplus}
\end{equation}
\begin{equation}
W_+ = 2(-i)^{n-1}\big(\mathcal F_-\big)^{\ast}r_{-}\mathcal F_{0}, \quad 
W_- =  2i^{n-1}\big(\mathcal F_-\big)^{\ast}r_{+}\mathcal F_{0}.
\label{eq:WplusminusFminus}
\end{equation}
\end{theorem}
Proof. Theorem 2.7.11 proves (1). Lemma 1.4 and (\ref{eq:FplusandJFminus}) imply $\widehat S_{phy}(k)J{\mathcal F}^{(+)}(-k)$ $ = {\mathcal F}^{(+)}(k)$ for $k > 0$, which proves (2).
The formula (\ref{eq:WaveOpandFourier}) proves (\ref{eq:WpluminusFplus}) for $W_+$.  For $f, g \in \mathcal B$, we have by using (\ref{eq:FplusandJFminus}) and (\ref{eq:WaveOpandFourier}) for $W_-$
\begin{eqnarray*}
\big(W_-f,g) &=& (\mathcal F^0f,\mathcal F^{(-)}g) \\
&=& \int_0^{\infty}(\mathcal F^0(k)f,\mathcal F^{(-)}(k)g)k^{n-1}dk \\
&=& \int_{-\infty}^0(J\mathcal F^0(k)f,J\mathcal F^{(+)}(k)g)|k|^{n-1}dk \\
&=& \int_{-\infty}^0((-ik + 0)^{(n-1)/2}\mathcal F^0(k)f,(-ik + 0)^{(n-1)/2}\mathcal F^{(+)}(k)g)dk \\
&=& 2\big((\mathcal F_+)^{\ast}r_-\mathcal F_0f,g).
\end{eqnarray*}
This proves (\ref{eq:WpluminusFplus}) for $W_-$. By a similar computation using $$
(\mp ik + 0)^{\alpha} = e^{\mp {\rm sgn}(k)\alpha\pi i/2}|k|^{\alpha}, \quad
{\rm sgn}(k) = k/|k|,
$$
we have
\begin{eqnarray*}
\big(W_+f,g)
&=& \int_0^{\infty}(\mathcal F^0(k)f,\mathcal F^{(+)}(k)g)k^{n-1}dk \\
&=& \int_{-\infty}^0(J\mathcal F^0(k)f,J\mathcal F^{(-)}(k)g)|k|^{n-1}dk \\
&=& \int_{-\infty}^0((ik + 0)^{(n-1)/2}\mathcal F^0(k)f,(ik + 0)^{(n-1)/2}\mathcal F^{(-)}(k)g)dk \\
&=& (-i)^{n-1}\int_{-\infty}^0((-ik + 0)^{(n-1)/2}\mathcal F^0(k)f,(ik + 0)^{(n-1)/2}\mathcal F^{(-)}(k)g)dk\\
&=& 2(-i)^{n-1}\big((\mathcal F_-)^{\ast}r_-\mathcal F_0f,g),
\end{eqnarray*}
which proves (\ref{eq:WplusminusFminus}) for $W_+$. Finally by (\ref{eq:WaveOpandFourier})
\begin{eqnarray*}
\big(W_-f,g)
&=& \int_0^{\infty}(\mathcal F^0(k)f,\mathcal F^{(-)}(k)g)k^{n-1}dk \\
&=& \int_{0}^{\infty}((ik +0)^{(n-1)/2}\mathcal F^0(k)f,(ik + 0)^{(n-1)/2}\mathcal F^{(-)}(k)g)dk \\
&=&i^{n-1} \int_{0}^{\infty}((-ik + 0)^{(n-1)/2}\mathcal F^0(k)f,(ik + 0)^{(n-1)/2}\mathcal F^{(-)}(k)g)dk \\
&=& 2i^{n-1}\big((\mathcal F_-)^{\ast}r_+\mathcal F_0f,g),
\end{eqnarray*}
which proves (\ref{eq:WplusminusFminus}) for $W_-$.
 \qed

\bigskip
As a consequence of Theorem 1.7 (2), we have
\begin{equation}
g \in {\rm Ran}\,\mathcal F_0 \Longleftrightarrow 
g(-k,-\omega) = i^{n-1}g(k,\omega), \quad k > 0.
\nonumber
\end{equation}
The projection onto the range of $\mathcal F_0$ is written as follows.


\begin{lemma} We define the operator $\widetilde J$  by $(\widetilde Jf)(k,\omega) = f(-k,-\omega)$. Then 
\begin{equation}
 \mathcal F_0(\mathcal F_0)^{\ast} 
=  \frac{1}{2} + \frac{1}{2}\big((-i)^{n-1} r_+ + i^{n-1} r_-\big)\widetilde J.
\nonumber
\end{equation}
\end{lemma}
Proof. We put $\big(U_0f\big)(k,\omega) = \frac{1}{\sqrt2}|k|^{(n-1)/2}\widehat f(k\omega)$.
Then $U_0$ is an isometry from $L^2({\bf R}^n)$ to $L^2({\bf R};L^2(S^{n-1});dk)$ and
$$
g \in {\rm Ran}\,U_0 \Longleftrightarrow g = \widetilde Jg.
$$
Since $U_0(U_0)^{\ast}$ is the projection onto the range of $U_0$, we have
$$
U_0(U_0)^{\ast} = \frac{1}{2}(1 + \widetilde J).
$$
Let $h = \zeta^{1/2} r_+ + \overline{\zeta}^{1/2}r_-$, $\zeta = e^{-(n-1)\pi i/2}$. Then we have $\mathcal F_0 = hU_0$, hence
$$
\mathcal F_0(\mathcal F_0)^{\ast} = hU_0(U_0)^{\ast}h^{\ast}.
$$
As can be checked easily
\begin{equation}
\widetilde Jr_{\pm} = r_{\mp}\widetilde J.
\label{eq:Chap2Sect7widetildeJandrplusminus}
\end{equation}
Using these formulas we obtain the lemma by a direct computation. \qed

\begin{cor} 
\begin{equation} 
\mathcal F_+ = r_+\mathcal F_0(W_+)^{\ast} + r_-\mathcal F_0(W_-)^{\ast},
\label{eq:FplusandWplusminus}
\end{equation}
\begin{equation}
\mathcal F_- = i^{n-1}r_+\mathcal F_0(W_-)^{\ast} + (-i)^{n-1} r_-\mathcal F_0(W_+)^{\ast}.
\label{eq:FminusandWplusminus}
\end{equation}
\end{cor}
Proof. By (\ref{eq:WpluminusFplus}), $\mathcal F_0\big(W_{\pm}\big)^{\ast} = 2\mathcal F_0({\mathcal F_0})^{\ast}r_{\pm}\mathcal F_+$. By Lemma 1.8 and (\ref{eq:Chap2Sect7widetildeJandrplusminus}) we have
$$
r_{\pm}\mathcal F_0(\mathcal F_0)^{\ast}r_{\pm} = \frac{1}{2}r_{\pm}.
$$
 This proves (\ref{eq:FplusandWplusminus}).
By (\ref{eq:WplusminusFminus}), we have $\mathcal F_0(W_+)^{\ast} = 2i^{n-1}\mathcal F_0(\mathcal F_0)^{\ast}r_-\mathcal F_-$, and $\mathcal F_0(W_-)^{\ast} = 2(-i)^{n-1}\mathcal F_0(\mathcal F_0)^{\ast}r_+\mathcal F_-$. Therefore 
$$
r_-\mathcal F_0(W_+)^{\ast} = i^{n-1}r_-\mathcal F_-, \quad 
r_+\mathcal F_0(W_-)^{\ast} = (-i)^{n-1}r_+\mathcal F_-.
$$
Hence (\ref{eq:FminusandWplusminus}) follows. \qed


\begin{definition}
The modified Radon transform $\mathcal R_{\pm}$ is defined  by 
$$
\left({\mathcal R}_{\pm} f\right)(s) = \frac{1}{\sqrt{2\pi}}\int_{-\infty}^{\infty}e^{iks}\left(\mathcal F_{\pm} f\right)(k)dk.
$$
\end{definition}

By (\ref{eq:LowerplusminusF}) and Theorem 1.7, we have


\begin{theorem}
$\mathcal R_{\pm} : L^2({\bf R}^n) \to L^2({\bf R};L^2(S^{n-1});dk)$ is an isometry and
$$
\left(\mathcal R_{\pm} H f\right)(s) = - \partial_s^2\left(\mathcal R_{\pm}f\right)(s).
$$
 Moreover
\begin{equation}
\left(\mathcal R_+f\right)(-s) = J\left(\mathcal R_-f\right)(s).
\nonumber
\end{equation}
\end{theorem}


\begin{definition} For an open interval $I \subset {\bf R}$,
let $\widehat H^m(I)$ be the set of functions $\phi(s,\omega)$ satisfying
$$
\sum_{0\leq j\leq m}\int_{I}\big\|\partial_s^j\phi(s,\cdot)\big\|^2_{L^2(S^{n-1})}ds < \infty. 
$$
If $I = {\bf R}$, we simply write $\widehat H^m$, in which case $m$ can be any real number by passing to the Fourier transformation. 
\end{definition}


\begin{lemma}
For any  $m \geq 0$ we have
$$
f \in H^m \Longleftrightarrow \mathcal R_{\pm}f  \in \widehat H^m.
$$
\end{lemma}
Proof. A direct consequence of Theorem 1.11. \qed


\subsection{Asymptotic profiles of solutions to the wave equation} 
The following theroem is proved in the same way as Theorem 2.8.9.


\begin{theorem}
For $x \in {\bf R}^n$, we write $r = |x|, \omega = x/r$. Then for $f  \in L^2({\bf R}^n)$, we have as $t \to  \infty$
\begin{equation}
\Big\|\left(\cos(t\sqrt{H})f\right)(x) - \frac{r^{-(n-1)/2}}{\sqrt 2}\left({\mathcal R}_{+} f\right)(r - t,\omega)\Big\| \to 0,
\nonumber
\end{equation}
\begin{equation}
\Big\|\left(\sin(t\sqrt{H})f\right)(x) - \frac{i\,r^{-(n-1)/2}}{\sqrt 2}\left(h\big(-i\frac{\partial}{\partial s}\big){\mathcal R}_{+}f\right)(r - t,\omega)\Big\|
 \to 0,
 \nonumber
\end{equation}
where $\|\cdot\|$ is the $L^2({\bf R}^n)$-norm, and $h(k) = 1 \ (k > 0), \ h(k) = - 1 \ (k < 0)$.
\end{theorem}


\subsection{Relation between scattering operators}
The scattering operator is also defined by the Radon transform, namely


\begin{definition}
$ \hskip 1cm \mathcal S_R = \mathcal R_+\left(\mathcal R_-\right)^{\ast}.$
\end{definition}

The following lemma follows easily from Theorem 1.11 and Lemma 1.13.


\begin{lemma}
(1) $\mathcal S_R$ is a partial isometry with initial set Ran$\,\big(\mathcal R_-\big)$ and final set Ran$\,\big(\mathcal R_+\big)$. \\
\noindent
(2)
$\displaystyle{\ 
\partial_s^2\mathcal S_R = \mathcal S_R \partial_s^2
\label{eq:Smatrixandds}
}$.\\
\noindent
(3) $\ {\mathcal S}_R \widehat H^m \subset \widehat H^m, \ \forall m \geq 0.$
\end{lemma}

The relation to the scattering operator $S$ in Definition 1.3 is as follows.


\begin{lemma}
Let $\mathcal F_1$ be the 1-dimensional Fourier transform, $r_{\pm}$ the projection in Theorem 1.7 (3) and $\widetilde J$ as in Lemma 1.8. Then we have
\begin{equation}
\mathcal F_1{\mathcal S}_R(\mathcal F_1)^{\ast}  = 
(-i)^{n-1}r_+\mathcal F_0S(\mathcal F_0)^{\ast}r_+ + 
i^{n-1}r_-\mathcal F_0S^{\ast}(\mathcal F_0)^{\ast}r_- + \frac{1}{2}\widetilde J.
\nonumber
\end{equation}
\end{lemma}
Proof. Since $\mathcal F_1\mathcal S_R(\mathcal F_1)^{\ast} = \mathcal F_+(\mathcal F_-)^{\ast}$, the lemma follows from Corollary 1.9. \qed


\section{Asymptotic solutions}

\subsection{Geometrical optics} In this section we construct an asymptotic solution to the equation
$$
- \Delta_g\left(e^{ik\varphi}a\right) = k^2 e^{ik\varphi}a,
$$
$k \in {\bf R}$ being a large parameter. We put
$
a = \sum_{j=0}^N k^{-j}a_j.
$ 
Then we have
\begin{equation}
\begin{split}
 e^{-ik\varphi}(- \Delta_g - k^2)e^{ik\varphi}a  & =  \; k^2\left[g^{\alpha\beta}(\partial_\alpha\varphi)(\partial_{\beta}\varphi) - 1\right]a 
- ik Ta - \Delta_ga \\
& = \; k^2\left[g^{\alpha\beta}(\partial_\alpha\varphi)(\partial_{\beta}\varphi) - 1\right]a -ikTa_0 \\
& \hskip 0.5cm - i\sum_{j=0}^{N-1}k^{-j}(Ta_{j+1} - i\Delta_ga_{j}) - 
ik^{-N}\Delta_ga_N,
\end{split}
\label{eq:EikonalIdentity}
\end{equation}
where $T$ is the following differential operator
\begin{equation}
T = 2g^{\alpha\beta}(\partial_{\alpha}\varphi)\partial_{\beta} + \Delta_g\varphi.
\nonumber
\end{equation}
We define the Hamiltonian $h(x,p)$ by
\begin{equation}
h(x,p) = \frac{1}{2} g^{ij}(x)p_ip_j.
\nonumber
\end{equation}

Our aim is to constrcut a real function $\varphi(x,\theta) \in C^{\infty}({\bf R}^{n}\times S^{n-1})$ which behaves like $x\cdot\theta + O(|x|^{-\epsilon_0})$ as $|x| \to \infty$, and solves the eikonal equation $h(x,\nabla_x\varphi) = 1/2$ in the region $\{x\cdot\theta + |x_{\perp}|/\epsilon > R\}$, where $x_{\perp} = x - (x\cdot\theta)\theta$, and $R$, $1/\epsilon$ are sufficiently large constants. We shall parametrize the bicharacteristics by the asymptotic data at infinity. 

We fix $\theta \in S^{n-1}$ arbitrarily. We seek a solution 
$x(t), p(t)$ of the Hamilton-Jacobi equation
\begin{equation}
\frac{dx}{dt} = \frac{\partial h}{\partial p}, \quad  
\frac{dp}{dt} = - \frac{\partial h}{\partial x},
\label{eq:HamiltonJacobi}
\end{equation}
having the following asymptotics:
\begin{equation}
x(t) = t\theta + y + O(t^{-\epsilon_0}), \quad p(t) = \theta + 
O(t^{-1-\epsilon_0}), \quad (t \to \infty)
\nonumber
\end{equation}
for some $y \in {\bf R}^n$. A simple calculation shows that $x(t)$ satisfies the following integral equation
\begin{equation}
x(t) = t\theta + y + \int_t^{\infty}(s - t)\frac{d^2x(s)}{ds^2}ds.
\nonumber
\end{equation}
Since Hamilton's equation (\ref{eq:HamiltonJacobi}) coincides with the equation of geodesic, 
we have
$$
\frac{d^2x^k}{dt^2} = - \Gamma^k_{ij}\frac{dx^i}{dt}\frac{dx^j}{dt} = 
- \Gamma^k_{ij}g^{i\alpha}g^{j\beta}p_{\alpha}p_{\beta},
$$
$\Gamma^k_{ij}$ being Christoffel's symbol. In view of these formulas, we  put
\begin{eqnarray*}
z(t) &=& x(t) - t\theta - y,\\
A^k(t,s,y,\theta;z,p) 
& =&  (t - s)\Gamma^k_{ij}
(s\theta + y + z)g^{i\alpha}(s\theta + y + z)g^{j\beta}(s\theta + y + z)p_{\alpha}p_{\beta},\\
B^k(s,y,\theta;z,p) &=& \frac{1}{2}\frac{\partial g^{ij}}{\partial x^k}
(s\theta + y + z)p_ip_j,\\
A &=& (A^1,\cdots,A^n), \quad
B = (B^1,\cdots,B^n), 
\end{eqnarray*}
and consider the integral equation
\begin{equation}
\left\{
\begin{split}
& \ z(t) = \int_t^{\infty}A(t,s,y,\theta;z(s),p(s))ds, \\
& \ p(t) = \theta + \int_t^{\infty}B(s,y,\theta;z(s),p(s))ds.
\end{split}
\right.
\label{eq:IntegralEquation}
\end{equation}
We fix a sufficiently small $\epsilon > 0$. For a sufficiently large $R > 0$, let $\Omega_{R,\epsilon}(\theta)$ be the region defined by
\begin{equation}
\Omega_{R,\epsilon}(\theta) = \{(t,y,z) \, ; \, t + |y|/\epsilon > R, \, y\cdot\theta = 0, 
\, |z| < 3\}.
\nonumber
\end{equation}
Then taking $R$ large enough we have by a simple computation  
\begin{equation}
|t\theta + y + z| \geq C(|t| + |y| + R), \quad
 \forall (t,y,z) \in \Omega_{R,\epsilon}(\theta),
 \label{eq:tthetaplusyfrombelow}
\end{equation}
where the constant $C$ is independent of $(t,y,z) \in \Omega_{R,\epsilon}(\theta)$ and $R > 0$.
We  put
\begin{equation}
X(t) = (z(t),p(t)),
\nonumber
\end{equation}
and define the non-linear map $\mathcal L(X)$ by
\begin{equation}
\mathcal L(X)(t,y,\theta) = 
\left(\int_t^{\infty}A(t,s,y,\theta;z(s),p(s))ds,
\int_t^{\infty}B(s,y,\theta;z(s),p(s))ds\right).
\nonumber
\end{equation}
We parametrize $y$ in the following way. Take vectors $e_1(\theta),\cdots,e_{n-1}(\theta)$ so that $e_1(\theta),\cdots,e_{n-1}(\theta)$ and $\theta$ form an orthonormal basis of ${\bf R}^n$. Then if $y\cdot\theta = 0$, $y$ is written as $y = \sum_{i=1}^{n-1}y_ie_i(\theta)$. This $(y_1,\cdots,y_{n-1})$ gives the desired parametrization. Note that $e_1(\theta),\cdots,e_{n-1}(\theta)$ can be chosen to be smooth with respect to $\theta \in S^{n-1}$ (at least locally). We put
\begin{equation}
|X|_{\infty} = \sup_{(t,y,z)\in\Omega_{R,\epsilon}(\theta)}|X(t)|.
\nonumber
\end{equation}


\begin{lemma} Suppose $|X|_{\infty} < 2, \ |\widetilde X|_{\infty} < 2$. Then the following inequalities hold:
$$
\left|\partial_t^m\partial_y^{\alpha}\mathcal L(X)(t,y,\theta)\right| 
\leq C_{m\alpha}
(|t| + |y| + R)^{-\epsilon_0 - m - |\alpha|}, \quad
\forall m, \alpha,
$$
$$
\left|\mathcal L(X)(t,y,\theta) - \mathcal L(\widetilde X)(t,y,\theta)\right| 
\leq C
(|t| + |y| + R)^{-\epsilon_0 }|X - \widetilde X|_{\infty}.
$$
\end{lemma}
Proof. This is a direct consequence of (\ref{eq:tthetaplusyfrombelow}) and the estimate 
$\partial_x^{\alpha}\Gamma^k_{ij}(x) = O(|x|^{-2-\epsilon_0 - |\alpha|})$,
which follows from (\ref{eq:DecayAssumption}).
\qed

\medskip
We now put $X_0 = (0,\theta)$ and take $R > 0$ large enough. Then by Lemma 2.1 and the standard method of iteration, there exists a unique solution $X(t,y,\theta)$ of the integral equation
\begin{equation}
X = X_0 + \mathcal L(X)
\nonumber
\end{equation}
in the region $\{t + |y|/\epsilon > R, \ y\cdot\theta = 0\}$ satisfying
\begin{equation}
\left|\partial_t^m\partial_y^{\alpha} (X(t,y,\theta) - X_0)\right| 
\leq C_{m\alpha}(|t| + |y| + R)^{-\epsilon_0 - m - |\alpha|}, \quad 
\forall m, \ \alpha.
\nonumber
\end{equation}

Returning back to the equation (\ref{eq:HamiltonJacobi}), we have proven the following lemma.

\begin{lemma}
Take $\theta \in S^{n-1}$ arbitrarily and $R > 0$ large enough. Then there exists a unique solution $x(t,y,\theta), p(t,y,\theta)$ of the equation (\ref{eq:HamiltonJacobi}) such that in the region $\{t + |y|/\epsilon > R, \ y\cdot\theta = 0\}$ it satisfies
$$
\left|\partial_t^m\partial_y^{\alpha}(x(t,y,\theta) - t\theta - y)\right| \leq 
C_{m\alpha}(|t| + |y| + R)^{-\epsilon_0 - m - |\alpha|}, \quad 
\forall m,\ \alpha,
$$
$$
\left|\partial_t^m\partial_y^{\alpha}(p(t,y,\theta) - \theta)\right| \leq 
C_{m\alpha}(|t| + |y| + R)^{-1-\epsilon_0 - m - |\alpha|}, \quad 
\forall m,\ \alpha.
$$
\end{lemma}
Proof. By differentiating the integral equation (\ref{eq:IntegralEquation}), we have
\begin{equation}
\frac{dx^k}{dt} = \theta^k + \int_t^{\infty}\Gamma_{ij}^kg^{i\alpha}g^{j\beta}
p_{\alpha}p_{\beta}ds,
\nonumber
\end{equation}
\begin{equation}
\frac{dp_k}{dt} = - \frac{1}{2}\frac{\partial g^{\alpha\beta}}{\partial x^k}
p_{\alpha}p_{\beta} = - \frac{\partial h}{\partial x^k}.
\label{eq:dpkdtequal-partialpartialx}
\end{equation}
Therefore we have to show that
\begin{equation}
g^{k\alpha}p_{\alpha} = \theta^k + \int_t^{\infty}\Gamma_{ij}^kg^{i\alpha}g^{j\beta}
p_{\alpha}p_{\beta}ds.
\nonumber 
\end{equation}
Since both sides tend to $\theta^k$ as $t \to \infty$, we have only to show that their time derivatives coincide. By (\ref{eq:dpkdtequal-partialpartialx}), the formula to show is
\begin{equation}
\frac{\partial g^{k\alpha}}{\partial x^i}g^{i\beta} - \frac{1}{2}g^{ki}\frac{\partial g^{\alpha\beta}}{\partial x^i} = - \Gamma^k_{ij}g^{i\alpha}g^{j\beta},
\nonumber
\end{equation}
which follows from a direct computation and the formula
\begin{equation}
\frac{\partial g^{ij}}{\partial x^m} = - g^{ik}\left(\frac{\partial g_{kr}}{\partial x^m}\right)g^{rj}.
\nonumber
\end{equation}
The estimates of $x(t), p(t)$ are easy to derive. \qed


\begin{lemma} As a 2-form on the region $\{(t,y) \, ; \, t + |y|/\epsilon > R, \ y\cdot\theta = 0\}$, we have
$$
\sum_{i=1}^ndp_i(t,y,\theta)\wedge dx^i(t,y,\theta) = 0.
$$
\end{lemma}
Proof. Without loss of generality we assume $\theta =(0,\cdots,0,1)$ and put
$y = (u^1,\cdots,u^{n-1},0)$, $t = u^n$.  Then we have
$$
\sum_idp_i\wedge dx^i = \sum_{j<k}[p,x]_{jk}du^j\wedge du^k,
$$
$$
[p,x]_{jk} = \frac{\partial p}{\partial u^j}\cdot\frac{\partial x}{\partial u^k} - \frac{\partial p}{\partial u^k}\cdot\frac{\partial x}{\partial u^j}.
$$
Noting that
$$
\frac{\partial}{\partial t}\left(\frac{\partial p}{\partial u^j}\cdot\frac{\partial x}{\partial u^k}\right) = 
- \frac{\partial^2h}{\partial x^i\partial x^m}
\frac{\partial x^m}{\partial u^j}\frac{\partial x^i}{\partial u^k} + 
 \frac{\partial^2h}{\partial p_i\partial p_m}
\frac{\partial p_i}{\partial u^k}\frac{\partial p_m}{\partial u^j}
$$
is symmetric with respect to $j$ and $k$, we have
$$
\frac{\partial}{\partial t}[p,x]_{jk} = 0.
$$
By Lemma 2.2, $[p,x]_{jk} \to 0$ as $t \to \infty$. Hence $[p,x]_{jk} = 0$, which proves the lemma. \qed

\medskip
For $x \in {\bf R}^n$, we put $x_{\perp} = x - (x\cdot\theta)\theta$ and define the region $\triangle_{R,\epsilon}(\theta)$ by
\begin{equation}
\triangle_{R,\epsilon}(\theta) = \{x \in {\bf R}^n \, ; \, x\cdot\theta + |x_{\perp}|/\epsilon > 
R\}.
\nonumber
\end{equation}
In the coordinates with basis $\theta, e_1(\theta), \cdots, e_{n-1}(\theta)$, the differential of the map $(t,y) \to x(t,y,\theta)$ is $I + O(R^{-\epsilon_0})$. Therefore the following lemma holds.


\begin{lemma}
For large $R > 0$, the map $(t,y) \to x(t,y,\theta)$ is a diffeomorphism and its image includes $\triangle_{2R,\epsilon}(\theta)$.
\end{lemma}

Let $t = t(x,\theta), y = y(x,\theta)$ be the inverse of the map : $(t,y) \to x(t,y,\theta)$. We put $p(x,\theta) = p(t(x,\theta),y(x,\theta),\theta)$ for the sake of simplicity. Lemma 2.3 implies $d(\sum_jp_j(x,\theta)dx^j) = 0$, which shows
\begin{equation}
\frac{\partial p_j(x,\theta)}{\partial x^i} = 
\frac{\partial p_i(x,\theta)}{\partial x^j}.
\label{eq:IntegrabilityCond}
\end{equation}
We put
\begin{equation}
f(x,\theta) = p(x,\theta) - \theta = \int_t^{\infty}
\frac{\partial h}{\partial x}\big(x(s,y,\theta),p(s,y,\theta)\big)ds\Big|_{t = t(x,\theta), y = y(x.\theta)}, 
\nonumber
\end{equation}
and define $\Psi(x,\theta)$ by
\begin{equation}
\Psi(x,\theta) = x\cdot\theta - \int_0^{\infty}f(x + t\theta,\theta)\cdot\theta dt.
\nonumber
\end{equation}


\begin{lemma}
On $\triangle_{2R,\epsilon}(\theta)$, we have
\begin{equation}
\nabla_x\Psi(x,\theta) = p(x,\theta),
\label{eq:nablaPsi}
\end{equation}
\begin{equation}
h(x,\nabla_x\Psi(x,\theta)) = 1/2, 
\label{eq:EnergyCons}
\end{equation}
\begin{equation}
|\partial_x^{\alpha}(\Psi(x,\theta) - x\cdot\theta)| \leq 
C_{\alpha}(1 + |x|)^{-\epsilon_0 - |\alpha|}, \quad \forall \alpha.
\label{eq:EstimatePsi}
\end{equation}
\begin{equation}
\Psi(x,\theta) = t(x,\theta).
\label{eq:traveltime}
\end{equation}
\end{lemma}
Proof. Letting $f = (f_1,\cdots,f_n)$, we have
${\displaystyle
\frac{\partial
 f_j}{\partial x^i}(x,\theta) = 
\frac{\partial f_i}{\partial x^j}(x,\theta)}$
by (\ref{eq:IntegrabilityCond}).
We then have
\begin{eqnarray*}
\frac{\partial\Psi}{\partial x^i} &=& \theta_i - 
\int_0^{\infty}\sum_j\frac{\partial f_j}{\partial x^i}(x + t\theta,\theta)\theta_jdt \\&=& \theta_i - 
\int_0^{\infty}\sum_j\frac{\partial f_i}{\partial x^j}(x + t\theta,\theta)\theta_jdt \\&=& \theta_i - 
\int_0^{\infty}\frac{d}{dt}f_i(x + t\theta,\theta)dt \\
&=& \theta_i + f_i(x,\theta) \\
&=& p_i(x,\theta),
\end{eqnarray*} 
which proves (\ref{eq:nablaPsi}). Since $x(t), p(t)$ solve the equation (3.2), $h(x(t),p(t))$ is a constant. Letting $t \to \infty$, this constant is seen to be equal to 1/2, which together with (\ref{eq:nablaPsi}) proves (\ref{eq:EnergyCons}). The estimate (\ref{eq:EstimatePsi}) follows from Lemma 2.1. By (\ref{eq:nablaPsi}), we have
$$
\frac{\partial\Psi}{\partial t} = \left(\partial_i\Psi\right)\frac{\partial x^i}{\partial t} = g^{ij}\left(\partial_i\Psi\right)\left(\partial_j\Psi\right) = 1.
$$
Therefore $\Psi = t + t_0(y,\theta)$ for some $t_0(y,\theta)$. However by Lemma 3.2, $x(t,y,\theta)\cdot\theta = t + O(t^{-\epsilon_0})$, which implies $t_0(y,\theta) = \Psi - x\cdot\theta + O(t^{-\epsilon_0}) = O(t^{-\epsilon_0})$. Therefore $t_0(y,\theta) = 0$, which proves (\ref{eq:traveltime}).
\qed 

\medskip
The equality (\ref{eq:IntegrabilityCond}) yields the following corollary.


\begin{cor}
For any smooth function $f(x)$ on ${\bf R}^n$, we have
$$
\frac{\partial}{\partial t} f(x(t,y,\theta))\Big|_{t=t(x,\theta),y=y(x,\theta)}  = g^{ij}(x)\frac{\partial \Psi(x,\theta)}{\partial x^j}\frac{\partial f(x)}{\partial x^i}.
$$
\end{cor}

By the above construction, $\Psi(x,\theta)$ is actually a function on the fibered space $\{(\theta,x) \; ;\; \theta \in S^{n-1}, x \in \Delta_{2R,\epsilon}(\theta)\}$ and satisfies
\begin{equation}
|\partial_{\theta}^{\alpha}\partial_{x}^{\beta}(\Psi(x,\theta) - x\cdot\theta)| \leq 
C_{\alpha\beta}(1 + |x|)^{-\epsilon_0 - |\beta|}, \quad \forall \alpha, \ \beta.\nonumber
\end{equation}

\begin{definition}
We take $\chi_{\infty}(t) \in C^{\infty}({\bf R})$ and $\chi(t) \in C^{\infty}({\bf R})$ such that $\chi_{\infty}(t) = 1, \ ( t > 3R)$, \
$\chi_{\infty}(t) = 0, \ (t < 2R)$, \
$\chi(t) = 1, \ (t > -1 +  2\epsilon)$,\ $\chi(t) = 0, (t  < - 1 + \epsilon)$, 
where $1/R$ and $\epsilon > 0$ are sufficiently small constants.
We define 
$$
 \varphi(x,\theta)   =  x\cdot\theta 
 + \chi_{\infty}(|x|)\chi(\widehat x\cdot\theta)\Big(\Psi(x,\theta) - x\cdot\theta\Big),
$$
$$
\varphi_{\pm}(x,\xi) = \pm |\xi|\,\varphi(x,\pm\widehat\xi\,), \quad
\widehat\xi = \xi/|\xi|.
$$
\end{definition}

The following lemma is a direct consequence of the above definition.


\begin{theorem} (1)
$\varphi_{\pm}(x,\xi) \in C^{\infty}({\bf R}^n\times\left({\bf R}^n\setminus\{0\}\right))$ and
$$
\Big|\partial_{\xi}^{\alpha}\partial_{x}^{\beta}\Big(\varphi_{\pm}(x,\xi) - x\cdot\xi\Big)\Big| 
\leq C_{\alpha\beta}|\xi|^{1-|\alpha|}(1 + |x|)^{- \epsilon_0 - |\beta|}.
$$
(2) If $|x| > 3R$ and $\pm \widehat x\cdot\widehat\xi > - 1 + 2\epsilon$, it satisfies the eikonal equation
$$
h(x,\nabla_x\varphi_{\pm}(x,\xi)) = |\xi|^2/2. 
$$
(3) $\ \varphi_-(x,\xi) = - \varphi_+(x,-\xi)$.
\end{theorem}


\subsection{Asymptotic solutions} 
We employ the above $\varphi(x,\theta)$ as $\varphi$ in (\ref{eq:EikonalIdentity}). Letting
\begin{equation}
a_0(x,\theta) = \exp\left(\int_t^{\infty}\frac{1}{2}
(\Delta_g\varphi)(x(s,y,\theta),\theta)ds\right)
\Big|_{t=t(x,\theta), y = y(x,\theta)},
\label{eq:A0}
\end{equation}
and using Corollary 2.6, we have
$$
Ta_0(x,\theta) = 0 \quad {\rm for} \quad
 |x| > 3R, \quad \widehat x\cdot\theta > - 1 + 2\epsilon.
$$
By Theorem 2.8 (1),
$a_0(x,\theta)$ satisfies
$$
|\partial_{\theta}^{\alpha}\partial_x^{\beta}\left(a_0(x,\theta) - 1\right)| \leq C_{\alpha\beta}(1 + |x|)^{-|\beta|-\epsilon_0}.
$$
We integrate the higer order transport equation
$$
Ta_j - i\Delta_ga_{j-1} = 0, \quad j \geq 1
$$
in a similar manner, and obtain
$$
|\partial_{\theta}^{\alpha}\partial_x^{\beta}\, a_j(x,\theta)| \leq C_{\alpha\beta}(1 + |x|)^{-j -|\beta|-\epsilon_0}.
$$
Let $\chi(t), \chi_{\epsilon}(t) \in C^{\infty}({\bf R})$ be such that 
$\chi(t) = 1 \ (t > 4), \ \chi(t) = 0 \ (t < 3), \ \chi_{\epsilon}(t) = 1 \ (t > - 1 + 3\epsilon), \ \chi_{\epsilon}(t) = 0 \ (t < - 1 + 2\epsilon)$.
We put
\begin{equation}
a(x,k,\theta) = g(x)^{1/4}
\chi_{\epsilon}(\widehat x\cdot\theta)
\sum_{j=0}^{\infty}k^{-j}a_j(x,\theta)\chi(\epsilon_j|x|)
\chi(\epsilon_j|k|).
\label{eq:AmplitudeWithK}
\end{equation}
By a suitable choice of the sequence $\epsilon_0 > \epsilon_1 >\cdots \to 0$, 
this series converges and defines a smooth function.
We finally define
\begin{equation}
a_{\pm}(x,\xi) = a(x,\pm|\xi|,\pm\widehat\xi\,).
\nonumber
\end{equation}
 The following lemma holds.


\begin{lemma}  
(1) On ${\bf R}^n\times{\bf R}^n$, $a_{\pm}(x,\xi)$ satisfies
$$
|\partial_{\xi}^{\alpha}\partial_x^{\beta}\, a_{\pm}(x,\xi)| \leq
C_{\alpha\beta}(1 + |\xi|)^{-|\alpha|}(1 + |x|)^{-|\beta|}.
$$
(2) Let $g_{\pm}(x,\xi) = e^{-i\varphi_{\pm}(x,\xi)}(L - |\xi|^2)e^{i\varphi_{\pm}(x,\xi)}a_{\pm}(x,\xi)$. Then  it satisfies
$$
\big|\partial_{\xi}^{\alpha}\partial_x^{\beta}\, g_{\pm}(x,\xi)\big| 
\leq C_{\alpha\beta N}(1 + |\xi|)^{-N}(1 + |x|)^{- N}
$$
 for any $N > 0$ in the region $ |x| > 4R, \ \pm 
\widehat x\cdot\widehat\xi > - 1 + 3\epsilon$.
\end{lemma}

\section{Fourier integral operators and functional calculus} 

\subsection{Product formula for FIO} Lets us recall the theory of FIO's. Since we need precise product formulas, we employ the computation  by \cite{Kum76}, \cite{Kum81}.
For $m \in {\bf R}$, let $S^{m}$ be the class of symbols defined by
$$
S^{m} \ni p(x,\xi) \Longleftrightarrow 
\big|\partial_{\xi}^{\alpha}\partial_x^{\beta}\,p(x,\xi)\big| \leq C_{\alpha\beta}
(1 + |\xi|)^{m - |\alpha|}, \quad 
\forall \alpha,\ \beta.
$$
The phase function $\varphi(x,\xi) \in 
C^{\infty}({\bf R}^n\times{\bf R}^n)$ is assumed to be real-valued and satisfy the following conditions $(3.1) \sim (3.4)$ for a sufficiently small constant $\delta_0 > 0$: 
\begin{equation}
\varphi(x,\xi) - x\cdot\xi \in S^{1}, 
\end{equation}
\begin{equation}
 \big|\nabla_{\xi}\left(\varphi(x,\xi) - x\cdot\xi\right)\big| < \delta_0, 
\end{equation}
\begin{equation}
\big|\nabla_{x}\left(\varphi(x,\xi) - x\cdot\xi\right)\big| < \delta_0(1 + |\xi|),
\end{equation}
\begin{equation}
\left|\frac{\partial^2}{\partial x\partial\xi}\varphi(x,\xi) - I\right| 
< \delta_0.
\end{equation}

\medskip
We define FIO's $I_{\varphi,a}, \ I_{\varphi^{\ast},a}$ by
$$
I_{\varphi,a}u(x) = (2\pi)^{-n}\iint\limits_{{\bf R}^n\times{\bf R}^n}
e^{i(\varphi(x,\xi) - y\cdot\xi)}a(x,\xi)u(y)dyd\xi,
$$
$$
I_{\varphi^{\ast},a}u(x) = (2\pi)^{-n}\iint\limits_{{\bf R}^n\times{\bf R}^n}
e^{i(x\cdot\xi - \varphi(y,\xi))}a(y,\xi)u(y)dyd\xi.
$$
We put $D_x = - i\partial_x$ and define the $\psi$DO $p(x,D_x)$ with symbol $p(x,\xi)$ by
$$
p(x,D_x)u(x) = (2\pi)^{-n}\iint\limits_{{\bf R}^n\times{\bf R}^n}
e^{i(x-y)\cdot\xi}p(x,\xi)u(y)dyd\xi.
$$
Using the conditions (3.1) $\sim$ (3.4) we can prove the following lemma.

\begin{lemma}
(1) The map ${\bf R}^n \ni \xi \to \eta = \nabla_x\varphi(x,\xi) \in {\bf R}^n$ is a global diffeomorphism on ${\bf R}^n$. Letting its inverse by $\xi(x,\eta)$, we have
$$
\xi(x,\eta) - \eta \in S^{1},
$$
$$
C^{-1}(1 + |\eta|) \leq 1 + |\xi| \leq C(1 + |\eta|).
$$
(2) The map ${\bf R}^n \ni x \to y = \nabla_{\xi}\varphi(x,\xi)$ is a global diffeomorphism on  ${\bf R}^n$. Letting $x(y,\xi)$ be its inverse, we have
$$
x(y,\xi) - y \in S^{0},
$$
$$
C^{-1}(1 + |y|) \leq 1 + |x| \leq C(1 + |y|).
$$
\end{lemma}

In the following Theorem 3.2, all symbols 
$c(x,\xi)$ belong to $S^{s_1+s_2}$ and have the following asymptoic expansion: 
\begin{equation}
c(x,\xi) \sim \sum_{j=1}^{\infty}c_j(x,\xi),
\quad c_j(x,\xi) \in S^{s_1 + s_2 - j}.
\label{eq:Sect4SymbolAsymptotic Expansion}
\end{equation}


\begin{theorem}
Let $a \in S^{s_1}, \ b \in S^{s_2}$. Then we have the 
following formulas.  
\begin{equation}
\left\{
\begin{split}
& I_{\varphi,a}I_{\varphi^{\ast},b} = c(x,D_x), \\
& c(x,\eta) \sim  a(x,\xi)b(x,\xi)\det\left(\frac{\partial^2}{\partial x\partial\xi}\varphi(x,\xi)\right)^{-1}\Bigg|_{\xi=\xi(x,\eta)} + \cdots,
\end{split}
\right.
\label{eq:PhiAstPhi}
\end{equation}
where $\xi(x,\eta)$ is the inverse map of $\eta = \nabla_x\varphi(x,\xi)$,
\begin{equation}
\left\{
\begin{split}
& I_{\varphi^{\ast},a}I_{\varphi,b} = c(x,D_x), \\
&
c(y,\xi) \sim a(x,\xi)b(x,\xi)\det\left(\frac{\partial^2}{\partial x\partial\xi}\varphi(x,\xi)\right)^{-1}\Bigg|_{x = x(y,\xi)} + \cdots,
\end{split}
\right.
\label{eq:AstPhiPhi}
\end{equation}
where $x(y,\xi)$ is the inverse map of $y = \nabla_{\xi}\varphi(x,\xi)$,
\begin{equation}
\left\{
\begin{split}
& I_{\varphi,a}b(x,D_x) = I_{\varphi,c},\\
& c(x,\xi) \sim a(x,\xi)b(\nabla_{\xi}\varphi(x,\xi),\xi) + \cdots,
\end{split}
\right.
\end{equation}
\begin{equation}
\left\{
\begin{split}
& a(x,D_x)I_{\varphi,b} = I_{\varphi,c},\\
&
c(x,\xi) \sim a(x,\nabla_x\varphi(x,\xi))b(x,\xi) + \cdots.
\end{split}
\right.
\end{equation}
\end{theorem}

For the proof, see \cite{Kum76}, Theorems 2.1 $\sim$ 2.4. We need the following explicit form of the asymptotic expansion (3.5) later. We put
$$
\widetilde\nabla_{\xi}\varphi(x,\xi,\eta) = \int_0^1\left(\nabla_{\xi}\varphi\right)(x,t\xi + (1 - t)\eta)dt,
$$
$$
\widetilde\nabla_x\varphi(x,y,\xi) = \int_0^1\left(\nabla_x\varphi\right)(tx + (1 - t)y,\xi)dt.
$$
Then $c(x,\xi)$ in (3.8) has the following asymptotic expansion:
\begin{equation}
c(x,\eta) \sim \sum_{\alpha}\frac{1}{\alpha!}\partial_{\xi}^{\alpha}
\left\{a(x,\xi)(D_x^{\alpha}b)(\widetilde\nabla_{\xi}\varphi(x,\xi,\eta),\eta)\right\}\Big|_{\xi = \eta},
\label{eq:AsympSymbPhiB}
\end{equation}
and $c(x,\xi)$ in (3.9) has the following asymptotic expansion:
\begin{equation}
c(x,\xi) \sim \sum_{\alpha}\frac{1}{\alpha!}D_{y}^{\alpha}
\left\{\big(\partial_{\xi}^{\alpha}a\big)(x,\widetilde\nabla_x\varphi(x,y,\xi))b(y,\xi)\right\}\Big|_{y=x},
\end{equation}
(see \cite{Kum76}, (2.41), (2.57)).


\subsection{Functional calculus}
In Chap. 3, \S 2, we have introduced the almost analytic extension $F(z)$ of $f(t)$. By the construction procedure, we see that $\partial_tF(t + is)$ is an almost analytic extension of $f'(t)$.
 Let 
\begin{equation}
X = (1 + |x|^2)^{1/2}, \quad \Lambda = (1 + |D_x|^2)^{1/2}.
\label{eq:EksLambda}
\end{equation}


\begin{lemma}
Let $f(t) \in C_0^{\infty}({\bf R})$. Then we have for any $N > 0$
\begin{equation}
f(H) = f(H_0) + \sum_{n=1}^Np_n(x,D_x)f^{(n)}(H_0) + R_N,
\label{eq:Lemma4.4Expansionoff(L)}
\end{equation}
where $p_n(x,D_x) = \sum_{|\alpha| \leq \mu(n)}a_{\alpha}^{(n)}(x)D_x^{\alpha}$ such that $|\partial_{x}^{\beta}\,a_{\alpha}^{(n)}(x)| \leq C_{\alpha\beta}(1 + |x|)^{-|\beta|- 1 - \epsilon_0}$, 
and $R_N$ satisfies
\begin{equation}
X^{N}\Lambda^{N}R_N\Lambda^{N}X^{N} \in {\bf B}(L^2({\bf R}^n)).
\label{eq:Lemma4.4EstimatesofRN}
\end{equation}
\end{lemma}
Proof. We first prove the lemma with the property (\ref{eq:Lemma4.4EstimatesofRN}) replaced by
\begin{equation}
X^NR_NX^N \in {\bf B}(L^2({\bf R}^n)).
\label{eq:Lemma4.4WeakEstimate}
\end{equation}
 We prove the case $N = 1$. By the resolvent equation, we have
\begin{eqnarray*}
(z - H)^{-1} - (z - H_0)^{-1} &=& (z - H)^{-1}V(z - H_0)^{-1} \\
&=& V(z - H)^{-1}(z - H_0)^{-1} + [(z - H)^{-1},V](z - H_0)^{-1}\\
&=& V(z - H_0)^{-2} + K(z),
\end{eqnarray*}
\begin{eqnarray*}
K(z) &=& V(z - H)^{-1}V(z - H_0)^{-2} \\
& & + \ (z - H)^{-1}[H,V](z - H)^{-1}(z - H_0)^{-1}.
\end{eqnarray*}
Therefore by virtue of Lemma 3.2.1
\begin{equation}
\begin{split}
f(H) - f(H_0) &= V\frac{1}{2\pi i}\int_{{\bf C}}\overline{\partial_z}F(z)
(z - H_0)^{-2}dzd\overline{z} \\
&  + \frac{1}{2\pi i}\int_{{\bf C}}\overline{\partial_z}F(z)
K(z)dzd\overline{z}.
\end{split}
\label{eq:f(L)-f(L_0)}
\end{equation}
Since $\partial_tF(t + is)$ is an almost analytic extension of $f'(t)$, we have  by integration by parts 
\begin{equation}
\begin{split}
f'(H_0) & = \frac{1}{2\pi i}\int_{{\bf C}}\overline{\partial_z}\partial_tF(z)
(z - H_0)^{-1}dzd\overline{z} \\
 & = \frac{1}{2\pi i}\int_{{\bf C}}\overline{\partial_z}F(z)
(z - H_0)^{-2}dzd\overline{z}.
\end{split}
\nonumber
\end{equation}
Therefore the 1st term of the right-hand side of (\ref{eq:f(L)-f(L_0)}) is equal to $Vf'(H_0)$. If $P_j$ is a differential operator of order $j = 1, 2$ with bounded coefficients, we have by passing to the spectral decomposition 
$$
\|P_j(z - H)^{-1}\| \leq 
C\,|{\rm Im}\,z|^{-1}(1 + |z|)^{j/2}.
$$
We then have
$$
\|XK(z)X\| \leq C|{\rm Im}\,z|^{-p}(1 + |z|)^{p},
$$
for some $p \geq 2$. Since $F(z)$ satisfies $|\overline{\partial_z}F(z)| \leq C
|{\rm Im}\,z|^{p}(1 + |z|)^{s-p-1}$ for any $s < 0$, the remainder term has the desired estimate (\ref{eq:Lemma4.4WeakEstimate}). The proof for $N \geq 2$ is similar.

 Now for $f \in C_0^{\infty}({\bf R})$ we take $\chi \in C_0^{\infty}({\bf R})$ such that $\chi(t) = 1$ on ${\rm supp}\,\chi$. We multiply (\ref{eq:Lemma4.4Expansionoff(L)}) by the expansion
\begin{equation}
\chi(H) = \chi(H_0) + \sum_{j=1}^N\chi^{(j)}(H_0)q_j(x,D_x) + (\widetilde R_N)^{\ast},
\nonumber
\end{equation}
with $q_{j}(x,D_x)$ and $\widetilde R_N$ having the above mentioned properties.
We then have
$$
f(H_0)\chi(H) = f(H_0) + f(H_0)(\widetilde R_N)^{\ast}.
$$
Since $\widetilde R_N$ satisfies (\ref{eq:Lemma4.4WeakEstimate}), one can prove that $f(H_0)(\widetilde R_N)^{\ast}$ satisfies (\ref{eq:Lemma4.4EstimatesofRN}). One can deal with $p_n(x,D_x)f^{(n)}(H_0)\chi(H)$ and $R_N\chi(H)$ in a similar manner. 
\qed


\section{Parametrices and regularizers}
We construct parametrices for the wave equation in the form of a FIO using $\varphi_{\pm}$ and $a_{\pm}$ in \S 2. Recall that $\varphi_{\pm}, a_{\pm}$ contain cut-off functions. Here we need another cut-off function which restricts $x$ and $\xi$ in a smaller region.
Let $R$ and $\epsilon$ be as in Definition 2.7. Take $\chi_{\infty}(t), \ \chi(t) \in C^{\infty}({\bf R})$ such that
$\chi_{\infty}(t) = 1 \ (t > 10R)$, $\chi_{\infty}(t) = 0 \ (t < 9R)$, $\chi(t) = 1\ (t > - 1 + 5\epsilon), \ \chi(t) = 0 \ (t < - 1 + 4\epsilon)$, and put
\begin{equation}
\chi_{\pm}(x,\xi) = \chi_{\infty}(|x|)\chi_{\infty}(|\xi|)\chi(\pm \widehat x\cdot\widehat\xi\,).
\label{eq:cutoffSection5}
\end{equation}


\begin{definition}
Let $\varphi_{\pm}$, $a_{\pm}$ be as in Theorem 2.8 and Lemma 2.9, and $\chi_{\pm}$ as in (\ref{eq:cutoffSection5}). We define a 
FIO $U_{\pm}(t)$ by
\begin{equation}
U_{\pm}(t) = I_{\varphi_{\pm},a_{\pm}}e^{- it\sqrt{H_0}}I_{\varphi_{\pm}^{\ast},\chi_{\pm}}.
\nonumber
\end{equation}
\end{definition}

In the following, $\|\cdot\|$ denotes either the operator norm $\|T\|_{{\bf B}(L^2({\bf R}^n))}$ of a bounded operator $T$ on $L^2({\bf R}^n)$ or the $L^2$-norm $\|u\|_{L^2({\bf R}^n)}$ of a vector $u \in L^2({\bf R}^n)$. There will be no fear of confusion. 
We put
\begin{equation}
G_+(t) = \frac{d}{dt}\left(e^{it\sqrt{H}}U_+(t)\right).
\nonumber
\end{equation}
Let $X$ and $\Lambda$ be as in (\ref{eq:EksLambda}). 


\begin{lemma} For any $N > 0$, there exists a constant $C_N > 0$ such that
$$
\|\Lambda^NG_{+}(t)\Lambda^NX^N\| \leq C_N(1 + t)^{-N}, \quad  t > 0.
$$
\end{lemma}
Proof. We have
$$
G_+(t) = 
e^{it{\sqrt H}}\left(i{\sqrt H}U_+(t) + 
\frac{d}{dt}U_+(t)\right).
$$
We decompose this operator into two parts and make use of the tools in \S 3.

\medskip
\noindent
{\it Low energy part}. 
First we deal with the low energy part. We take $\chi_0(t) \in C^{\infty}({\bf R})$ such that $\chi_0(t) = 1 \ (t < 1)$, $\chi_0(t) = 0 \ (t > 2)$ and consider $\Lambda^Ne^{it\sqrt{H}}{\sqrt H}\chi_0(H)U_+(t)$. Noting that
$$
\Lambda^Ne^{it\sqrt{H}}\sqrt{H}\chi_0(H)U_+(t) = 
\Lambda^N(1 + H)^{-N/2}e^{it\sqrt{H}}(1 + H)^{N/2}\sqrt{H}\chi_0(H)U_+(t),
$$
we have only to show 
\begin{equation}
\|\chi_0(H)U_+(t)\Lambda^NX^N\| \leq C_N(1 + t)^{-N}, \quad \forall t, N > 0.
\label{eq:LowEnergy}
\end{equation}
We decompose $\chi_0(H)U_+(t)$ into two parts:
\begin{equation}
\chi_0(H)U_+(t) = \chi_0(H)I_{\varphi_+,a_+}\cdot e^{-it\sqrt{H_0}}I_{\varphi_+^{\ast},\chi_+}.
\label{eq:OpDecompo}
\end{equation}


\begin{prop} 
$ \ \chi_0(H)I_{\varphi_+,a_+}\Lambda^{N}X^N \in {\bf B}(L^2({\bf R}^n)), \quad
\forall N > 0$. 
\end{prop}
Proof.
Lemma 3.3 entails the asymptotic expansion
\begin{equation}
\chi_0(H) = \chi_0(H_0) + \sum_{n=1}^{N}p_n(x,D_x) + R_N,
\label{eq;ki0expansion}
\end{equation}
\begin{equation}
p_n(x,\xi) = 0 \quad {\rm for} \quad |\xi| > 2, \quad
X^{N}\Lambda^{N}R_N\Lambda^{N}X^{N} \in {\bf B}(L^2({\bf R}^n)).
\label{eq:chi0symbol}
\end{equation}  
By the construction of $a_+(x,\xi)$ in \S 2 (see (\ref{eq:AmplitudeWithK})),
$|\xi| \geq 1/\epsilon_0$ and $|x| \geq 1/\epsilon_0$ on supp$\,a_+(x,\xi)$. Therefore in the expression
\begin{equation}
\iint e^{-ix\cdot\eta}\chi_0(|\eta|^2) e^{i\varphi_+(x,\xi)}a_+(x,\xi)
(1 + |\xi|^2)^{N/2}(1 - \Delta_{\xi})^{N/2}\widehat f(\xi)d\xi dx,
\label{eq:Intchi0}
\end{equation}
which is the Fourier transform of $\chi_0(L_0)I_{\varphi_+,a_+}\Lambda^NX^Nf$,
the phase has the following estimate
$$
\big|\nabla_x\big(x\cdot\eta - \varphi_+(x,\xi)\big)\big| \geq C(1 + |\xi|), \quad C > 0.
$$
Using the differential operator
$$
P = i\big|\eta - \nabla_x\varphi_+(x,\xi)\big|^{-2}\big(\eta - \nabla_x\varphi_+(x,\xi)\big)\cdot\nabla_x,
$$
and integration by parts, we can then rewrite (\ref{eq:Intchi0}) as
$$
\iint e^{-i(x\cdot\eta - \varphi_+(x,\xi))}\chi_0(|\eta|^2)\big(P^{\ast}\big)^{2N}a_+(x,\xi)(1 + |\xi|^2)^{N/2}(1 - \Delta_{\xi})^{N/2}\widehat f(\xi)d\xi dx.
$$
Since $|\big(P^{\ast}\big)^{2N}a_+(x,\xi)| \leq C_N(1 + |x|)^{-2N}(1 + |\xi|)^{-2N}$, by integrating by parts with respect to $\xi$, the proposition is proved if $\chi_0(H)$ is replaced by $\chi_0(H_0)$ .
By (\ref{eq:chi0symbol}) one can prove the same result if $\chi_0(H_0)$ is replaced by $p_n(x,D_x)$ or $R_N$. This proves the above proposition. \qed

\medskip
By (\ref{eq:OpDecompo}) and Proposition 4.3, the proof of (\ref{eq:LowEnergy}) is reduced to the following Proposition.


\begin{prop}
$$
\|X^{-N}\Lambda^{-N}e^{-it\sqrt{H_0}}I_{\varphi_+^{\ast},\chi_+}\Lambda^NX^N\| \leq C_N(1 + t)^{-N}, \quad \forall t, N > 0.
$$
\end{prop}
Proof. We estimate the phase function of
 \begin{equation}
e^{-it\sqrt{H_0}}I_{\varphi_+^{\ast},\chi_+}f = 
(2\pi)^{-n}\iint e^{i(x\cdot\xi - t|\xi| - \varphi_+(y,\xi))}
\chi_+(y,\xi)f(y)dyd\xi.
\nonumber
\end{equation}
First we have
$$
|\nabla_{\xi}(t|\xi| + \varphi_+(y,\xi))| \geq |t\widehat\xi + y| - C|y|^{-\epsilon_0}.
$$
Here the localization with respect to the directions of $y$ and $\xi$ plays an important  role. Since $\widehat\xi\cdot\widehat y > - 1 + 4\epsilon$ on supp$\,\chi_+(y,\xi)$, we have
\begin{equation}
\begin{split}
|t\widehat\xi + y|^2 & = t^2 + 2t|y|\widehat\xi\cdot\widehat y + |y|^2 \\
&\geq t^2 - 2t|y|(1 - 4\epsilon) + |y|^2 \\
&\geq 4\epsilon(t^2 + |y|^2).
\end{split}
\nonumber
\end{equation}
By choosing $R$ large enough, we have
\begin{equation}
|\nabla_{\xi}(t|\xi| + \varphi_+(y,\xi))| \geq 
C(t + |y|)
\label{eq:Propo5.4Estimateofphase}
\end{equation}
with a constant $C > 0$ independent of $y$ and $t > 0$. Integration by parts then proves the proposition. \qed

\medskip
\noindent
{\it High energy part}.
Next we consider
$\displaystyle{i\sqrt{H}(1 - \chi_0(H))U_+(t) + \frac{d}{dt}U_+(t)}$. By the definition of $g_+$ in Lemma 2.9, we have
\begin{equation}
HI_{\varphi_+,a_+} - I_{\varphi_+,a_+}H_0 = I_{\varphi_+,g_+},
\label{eq:LandL0}
\end{equation} 
which implies
$$
I_{\varphi_{+},a_{+}}(H_0 - z)^{-1} - (H - z)^{-1}I_{\varphi_{+},a_{+}} = 
(H - z)^{-1}I_{\varphi_{+},g_{+}}(H_0 - z)^{-1}.
$$
We put $f(t) = t^{-1/2}(1 - \chi_{0}(t))$ and let $F(z)$ be its almost analytic extension. 
Then we  have by virtue of Lemma 4.3
\begin{equation}
f(H)I_{\varphi_{+},a_{+}} - I_{\varphi_{+},a_{+}}f(H_0) = B,
\label{eq:ef(L)-ef(L0)}
\end{equation}
\begin{equation}
B = \frac{1}{2\pi i}\int\limits_{{\bf C}}\overline{\partial_z}F(z)(H - z)^{-1}
\nonumber
I_{\varphi_{+},g_{+}}(H_0 - z)^{-1}dzd\overline z.
\end{equation}
Using this formula, we then have 
\begin{eqnarray*}
\sqrt{H}(1 - \chi_0(H))I_{\varphi_+,a_+} &=& 
f(H)HI_{\varphi_+,a_+} \\
&=& f(H)I_{\varphi_+,a_+}H_0 + f(H)I_{\varphi_+.g_+} \\
&=& I_{\varphi_+,a_+}f(H_0)H_0 + BH_0 + f(H)I_{\varphi_+,g_+},
\end{eqnarray*}
where we have used (\ref{eq:LandL0}), (\ref{eq:ef(L)-ef(L0)}) in the first and second lines.
Therefore we have
\begin{equation}
\begin{split}
& i\sqrt{H}(1 - \chi_0(H))U_+(t) + \frac{d}{dt}U_+(t) \\
&= 
iBH_0e^{-it\sqrt{H_0}}I_{\varphi_+^{\ast},\chi_+}  + if(H)I_{\varphi_+,g_+}e^{-it\sqrt{H_0}}I_{\varphi_+^{\ast},\chi_+} \\
& \ \ - iI_{\varphi_+,a_+}\sqrt{H_0}\chi_0(H_0)e^{-it\sqrt{H_0}}I_{\varphi^{\ast},\chi_+}.
\end{split}
\label{eq:ThreeTerms}
\end{equation}
The third term of the right-hand side vanishes, since $\chi_0(|\xi|^2)\chi_+(y,\xi) = 0$. 
 Let us consider the second term. Taking notice of the relation
 $$
 \Lambda^Ne^{it\sqrt{H}}f(H) = \Lambda^N(1 + H)^{-N/2}\cdot e^{it\sqrt{H}}\cdot f(H)(1 + H)^{N/2}\Lambda^{-N}\cdot \Lambda^N,
 $$
we have only to show the following

 
\begin{prop} 
$$
 \|\Lambda^NI_{\varphi_+,g_+}e^{-it\sqrt{H_0}}I_{\varphi_+^{\ast},\chi_+}\Lambda^NX^N\| \leq C_N(1 + t)^{-N}, \ \forall t, N > 0. 
 $$
\end{prop}
Proof. We choose $\psi_1(t), \psi_2(t) \in C^{\infty}({\bf R})$ such that 
$\psi_1(t) + \psi_2(t) = 1 \ (t \in {\bf R})$, $\psi_1(t) = 1 \ (t < - 1 + 3\epsilon)$, $\psi_1(t) = 0 \ (t > - 1 + 7\epsilon/2)$, and put
$$
J_k(t)f = (2\pi)^{-n}\iint e^{i(\varphi_+(x,\xi) - t|\xi| - \varphi_+(y,\xi))}
\psi_k(\widehat x\cdot\widehat\xi\,)g_+(x,\xi)\chi_+(y,\xi)f(y)dyd\xi.
$$
Then $I_{\varphi_+,g_+}e^{-it\sqrt{L_0}}I_{\varphi_+^{\ast},\chi_+} = J_1(t) + J_2(t)$. Note that $\widehat x\cdot\widehat \xi > - 1+ 3\epsilon$ on the support of $\psi_2(\widehat x\cdot\widehat\xi)$, on which region $g_+(x,\xi)$ decays rapidly in $x$ and $\xi$ by Lemma 3.9. 
Using (\ref{eq:Propo5.4Estimateofphase}) and integrating by parts, we then have
$$
\|\Lambda^NJ_2(t)\Lambda^NX^N\| \leq C_N(1 + t)^{-N}, \quad 
\forall t, N > 0.
$$
  
  We next show that on the support of the integrand of $J_1(t)$
\begin{equation}
|\nabla_{\xi}(\varphi_+(x,\xi) - t|\xi| - \varphi_+(y,\xi))| \geq 
C(t + |x| + |y|)
\label{eq:PhiXtY}
\end{equation}
for a constant $C > 0$. Once this is proved, one can prove
$$
\|\Lambda^NJ_1(t)\Lambda^NX^N\| \leq C_N(1 + t)^{-N}, \quad 
\forall t, N > 0
$$
by integration by parts. To prove (\ref{eq:PhiXtY}), we put
$$
D_+ = \{y \in {\bf R}^n\,;\,\widehat y\cdot\widehat\xi > - 1 + 4\epsilon\}, \quad
D_- = \{x \in {\bf R}^n\,;\,\widehat x\cdot\widehat\xi < - 1 + 7\epsilon/2\}.
$$
Then there exists $0 < c_0 < 1$ such that
$$
y\cdot x \leq c_0|y||x| \quad {\rm if} \quad y \in D_+,\quad x \in D_-.
$$
We also see that $y + t\widehat\xi \in D_+$ if $y \in D_+$, $t \geq 0$.
Therefore
$$
|y + t\widehat\xi - x|^2 \geq (1 - c_0)(|y + t\widehat\xi\,|^2 + |x|^2).
$$
In the proof of Proposition 5.4, we have already seen that $|y + t\widehat\xi| \geq C(t + |y|)$ for some $C > 0$. This proves (\ref{eq:PhiXtY}). \qed

\medskip
It remains to consider the first term of the right-hand side of (\ref{eq:ThreeTerms}). 


\begin{prop}
$$
 \|\Lambda^NBH_0e^{-it\sqrt{H_0}}I_{\varphi^{\ast}_+,\chi_+}\Lambda^NX^N\| \leq C_N(1 + t)^{-N}, \quad \forall t, N > 0.
 $$
\end{prop}
Proof. We rewrite $BH_0e^{-it\sqrt{H_0}}I_{\varphi_+^{\ast},\chi_+}$ as
\begin{equation}
\begin{split}
& \frac{1}{2\pi i}\int_{\bf C}\left(\overline{\partial_z}F(z)\right)
|{\rm Im}\,z|^{-m}(1 + |z|)^{m - 1} \cdot 
|{\rm Im}\,z|(H - z)^{-1} \\
&\ \ \ \ \ \ \cdot I_{\varphi_+,g_+}\cdot\left(\frac{|{\rm Im}\,z|}{1 + |z|}\right)^{m-1}(H_0 - z)^{-1}L_0e^{-it\sqrt{H_0}}
I_{\varphi_+^{\ast},\chi_+}dzd\overline{z},
\end{split}
\nonumber
\end{equation}
$m$ being an arbitrily chosen integer. By the property of almost analytic extension, $\left(\overline{\partial_z}F(z)\right)
|{\rm Im}\,z|^{-m}(1 + |z|)^{m - 1}$ is integrable, and $\||{\rm Im}\,z|(H - z)^{-1}\|$ is uniformly bounded on ${\bf C}$. We show that by taking $m$ large enough, one can deal with $|{\rm Im}\,z|^{m-1}(1 + |z|)^{-m+1}(H_0 - z)^{-1}L_0$ like a $\psi$DO with smooth symbol whose operator norm is uniformly bounded in $z$. To show this, we have only to prove
\begin{equation}
\left(\frac{|{\rm Im}\,z|}{1 + |z|}\right)^{|\alpha|+1}
|\partial_{\xi}^{\alpha}(|\xi|^2 - z)^{-1}| \leq C(1 + |z|)^{-1},
\label{eq:EstimatesofQxiz}
\end{equation}
where $C$ is a constant independent of $\xi \in {\bf R}$ and $z \in {\bf C}\setminus{\bf R}$. In fact, one can show by induction that
$$
\partial_{\xi}^{\alpha}(|\xi|^2 -z)^{-1} = 
\sum_{n=1}^{|\alpha|}\frac{P_n(\xi)}{(|\xi|^2 - z)^{n+1}},
$$
where $P_n(\xi)$ is a polynomial of order $n$. Using the inequality $|\xi| \leq C(1 + |z| + ||\xi|^2 - z|)$, we have $|P_n(\xi)| \leq C((1 + |z|)^n + ||\xi|^2 -z|^n)$, which implies
$$
|\partial_{\xi}^{\alpha}(|\xi|^2 - z)^{-1}| \leq 
C\sum_{n=1}^{|\alpha|+1}\frac{(1+|z|)^{n-1}}{||\xi|^2-z|^{n}}.
$$
This proves (\ref{eq:EstimatesofQxiz}). Then by the same computation as in the proof of Proposition 4.5, we can prove the desired estimate. \qed

\medskip
The proof of Lemma 4.2 is now completed. \qed


\begin{lemma}
For any $f \in L^2({\bf R}^n)$ we have in the sense  of $L^2({\bf R}^n)$
$$
U_{\pm}(t)f = e^{- it\sqrt{H_0}}I_{\varphi_{\pm}^{\ast},\chi_{\pm}}f 
 + o(1), \quad t \to \pm \infty.
$$
\end{lemma}
Proof. We have only to prove that
$$
I_{\varphi_{\pm},a_{\pm}}e^{- it\sqrt{H_0}}g = e^{- it\sqrt{H_0}}g  + o(1), \quad {\rm as} \quad t \pm \infty
$$
for $g$ satisfying $\widehat g(\xi) = \chi_{\infty}(\xi)\widehat g(\xi) \in 
C_0^{\infty}({\bf R}^n)$. We prove the case as $t \to \infty$. Take $\chi_0(t), \chi_1(t) \in C^{\infty}({\bf R})$ such that $\chi_0(t) + \chi_1(t) = 1 \ (t \in {\bf R})$, $\chi_0(t) = 1 \ (t < 1/3)$, $\chi_0(t) = 0 \ (t > 2/3)$. Then we have
$$
\chi_0\big(\frac{|x|}{t}\big)I_{\varphi_{+},a_{+}}e^{- it\sqrt{H_0}}
g = (2\pi)^{-n/2}\int_{{\bf R}^n}e^{i(\varphi_{+}(x,\xi) - t|\xi|)}\chi_0\big(\frac{|x|}{t}\big)a_{+}(x,\xi)\widehat g(\xi)d\xi. 
$$
Since $\nabla_{\xi}(\varphi_+(x,\xi) - t|\xi|) = x - t\widehat\xi + O(|x|^{-\epsilon_0})$, we have
$$
|\nabla_{\xi}(\varphi_+(x,\xi) - t|\xi|)| \geq Ct
$$
for some constant $C > 0$ on the support of the integrand. By integration by parts, we then have
$$
\|\chi_0\big(\frac{|x|}{t}\big)I_{\varphi_{+},a_{+}}e^{- it\sqrt{H_0}}
g\| \leq C_Nt^{-N}, \quad \forall N, t >0.
$$
We rewrite $\chi_1(\frac{|x|}{t})I_{\varphi_{+},a_{+}}e^{- it\sqrt{H_0}}g$ as above. Since $a_+(x,\xi) = $ $ \chi(\epsilon_0|\xi|)\chi_{\epsilon}(\widehat x\cdot\widehat\xi\,) + O(|x|^{-\epsilon_0})$ (see (\ref{eq:AmplitudeWithK})), and the integral over the region $\{\widehat x \cdot\xi < 0\}$ disappears (which is proven by the same method of integration by parts), we have
$$
\chi_1\big(\frac{|x|}{t}\big)I_{\varphi_{+},a_{+}}e^{- it\sqrt{H_0}}
g = (2\pi)^{-n/2}\int_{{\bf R}^n}e^{i(\varphi_{+}(x,\xi) - t|\xi|)}\chi_1\big(\frac{|x|}{t}\big)\chi(\epsilon_0|\xi|)\widehat g(\xi)d\xi + o(1). 
$$
In (\ref{eq:cutoffSection5}), we take $R$ large enough so that $\chi_{\infty}(|\xi|) = \chi_{\infty}(|\xi|)\chi(\epsilon_0|\xi|)$. Then we have $\chi(\epsilon_0|\xi|)\widehat g(\xi) = \chi_{\infty}(|\xi|)\widehat g(\xi) = \widehat g(\xi)$. Therefore
\begin{equation}
\begin{split}
\chi_1\big(\frac{|x|}{t}\big)I_{\varphi_{+},a_{+}}e^{- it\sqrt{H_0}}
g & = \chi_1\big(\frac{|x|}{t}\big)e^{- it\sqrt{H_0}}g + o(1) \\
&= e^{-it\sqrt{H_0}}g + o(1),
\end{split} 
\nonumber
\end{equation}
which proves the lemma \qed

\medskip
Let $\widehat H^m$ be the Sobolev space in Definition 1.12.


\begin{definition}
(1) An operator $R$  is  called a {\it regularizer of order $N$} if it satisfies \begin{equation}
R \in \mathop\cap_{m=-\infty}^{\infty}{\bf B}(H^{m}\,;\,H^{m + N}) \quad 
{\rm or } \quad R \in \mathop\cap_{m=-\infty}^{\infty}{\bf B}(H^{m}\,;\,\widehat H^{m+N}).
\nonumber
\end{equation}
If $N$ can be taken arbitarily large, $R$ is simply called a regularizer. \\
\noindent
(2) A $\psi$DO $P_+$ $(P_-)$ is called an {\it approximate outgoing (incoming) projection} if its symbol $p_+(x,\xi)$ $(p_-(x,\xi))$ has the form
$$
p_{\pm}(y,\xi) = \chi_{\pm}(x,\xi)\Big|_{x = x_{\pm}(y,\xi)},
$$
where $\chi_{\pm}(x,\xi)$ is specified in (\ref{eq:cutoffSection5}), and $x_{\pm}(y,\xi)$ is the inverse function of $y = \nabla_{\xi}\varphi_{\pm}(x,\xi)$. 
\end{definition}

Let $W_{\pm}$ be the wave operator defined in Subsection 1.3.

\begin{theorem}
For any $N > 0$, there exist an approximate outgoing (incoming) projection $P_+$ ($P_-$) and a regularizer of order $N$, which is denotede by $R_{\pm}^N$, such that
$$
W_{\pm}P_{\pm}= I_{\varphi_{\pm},a_{\pm}}
P_{\pm} + R^N_{\pm}.
$$
\end{theorem}
Proof. We consider $W_+$. 
Lemmas 4.2 and 4.7 imply
\begin{equation}
W_{+}I_{\varphi_{+}^{\ast},\chi_{+}} = 
I_{\varphi_{+},a_{+}}I_{\varphi_{+}^{\ast},\chi_{+}} + 
\int_0^{\infty}G_{+}(t)dt,
\end{equation}
the 2nd term of the right-hand side being a regularizer. 
In the following we use the abbreviation
$$
 b\Big|_{x_+(y,\xi)} = b(x,\xi)\Big|_{x = x_+(y,\xi)}.
$$
We now put $b_0(x,\xi) = \det\Big({\partial^2\varphi_+}/{\partial x\partial \xi}\Big)\Big|_{x_+(y,\xi)}$, and let
$$
I_{\varphi_+^{\ast},\chi_+}I_{\varphi_+,b_0} = c_+(x,D_x).
$$
Then we have modulo a regularizer
$$
W_+c_+(x,D_x) \equiv I_{\varphi_+,a_+}c_+(x,D_x).
$$
By virtue of (4.7), $c_+(x,\xi)$ has an asymptotic expansion
$$
c_+(y,\xi) \sim \chi_+\Big|_{x_+(y,\xi)} + c_1(y,\xi) + \cdots, \quad c_1 \in S^{-1}.
$$
Let $\widetilde\chi_+(x,\xi)$ be a function similar to $\chi_+(x,\xi)$ such that $\chi_+(x,\xi) = 1$ on supp$\,\widetilde\chi_+(x,\xi)$. Namely, we slightly 
shrink the support of $\chi_+$. Let $q_1 \in S^{-1}$ and $Q_1$ be a $\psi$DO with symbol $\widetilde\chi\Big|_{x_+(y,\xi)} + q_1(y,\xi)$. Then the symbol of $c_+(x,D_x)Q_1$ has an asymptotic expansion
\begin{gather}
\chi_+\Big|_{x_+(y,\xi)}\widetilde\chi_+\Big|_{x_+(y,\xi)} + \chi_+\Big|_{x_+(y,\xi)}q_1 + c_1\widetilde\chi_+\Big|_{x_+(y,\xi)}
\nonumber \\
 + \sum_{|\alpha|=1}\partial_{\xi}^{\alpha}\chi_+\Big|_{x_+(y,\xi)}\cdot D_y^{\alpha}\widetilde\chi_+\Big|_{x_+(y,\xi)} \quad 
{\rm mod} \quad S^{-2}.
\nonumber
\end{gather}
We choose $q_1$ as follows:
$$
q_1 = - \frac{1}{\chi_+\Big|_{x_+(y,\xi)}}\left(c_1\widetilde\chi_+\Big|_{x_+(y,\xi)} + 
\sum_{|\alpha|=1}\partial_{\xi}^{\alpha}\chi_+\Big|_{x_+(y,\xi)}\cdot D_y^{\alpha}\widetilde\chi_+
\Big|_{x_+(y,\xi)} \quad \right).
$$
Since $\chi_+ = 1$ on supp$\,\widetilde\chi_+$, $q_1(y,\xi)$ is smooth and
$$
I_{\varphi^{\ast}_+,\chi_+}I_{\varphi_+,c_+}Q_1 = \widetilde c_+(x,D_x),
$$
$$
\widetilde c_+(y,\xi) \sim \chi_+\Big|_{x_+(y,\xi)} + c_2(y,\xi) + \cdots, \quad
c_2 \in S^{-2}.
$$
Repeating this procedure, we complete the proof the theorem. 
\qed


\section{Propagation of singularities}

\subsection{Singularity expansions I}
We show how $\mathcal R_+$ describes the singularities of solutions to the wave equation.  We start with the following lemma, which can be proved easily by integration by parts. 


\begin{lemma}
The integral operator defined by
\begin{equation}
\left(Af\right)(s,\omega) = \int_{-\infty}^{\infty}\int_{{\bf R}^n}
e^{ik(s-\omega\cdot y)}a(s,\omega,k,y)f(y)dkdy
\nonumber
\end{equation}
($s \in {\bf R}^1, \ \omega \in S^{n-1}$) is a regularizer if there exist  constants $\nu \in {\bf R}$ and $C_0 > 0$ such that
\begin{equation}
\left|\partial_s^{\alpha}\partial_k^{\beta}\partial_y^{\gamma}\,a(s,\omega,k,y)\right| \leq 
C_{\alpha\beta\gamma}(1 + |k|)^{\nu -\beta}, \quad \forall \alpha, \beta, \gamma, 
\label{eq:Section6Estimatesasomegaky}
\end{equation}
\begin{equation}
|s - \omega\cdot y| \geq C_0(1 + |s| + |y|)
\label{eq:sminusomegayfrombelow}
\end{equation}
on the support of $a(s,\omega,k,y)$. 
\end{lemma}

 By Corollary 1.9, we have the following expression:
\begin{equation}
\begin{split}
\big({\mathcal R}_+f\big)(s) &= 
\frac{1}{\sqrt{2\pi}}\int_0^{\infty}e^{iks}
\left(\mathcal F_0(W_+)^{\ast}f\right)(k)dk \\
&+
\frac{1}{\sqrt{2\pi}}\int_{-\infty}^0e^{iks}
\left(\mathcal F_0(W_-)^{\ast}f\right)(k)dk. 
\end{split}
\label{eq:RadonPlus}
\end{equation}
We take $\chi_R(s) \in C^{\infty}({\bf R})$  such that $\chi_R(s) = 0 \ (s < 15R)$, $\chi_R(s) = 1 \ (s > 20R)$, and study the singularity of $\chi_R(s)\mathcal R_+f(s)$ with respect to $s$.


\begin{lemma}
We take $N > 0$ large enough. Then there exist approximate outgoing, incoming projections  $P_+$, $P_-$ such that
\begin{equation}
\chi_R(s)\int_0^{\infty}e^{iks}\mathcal F_0(k)(W_+)^{\ast}dk \equiv
\chi_R(s)\int_0^{\infty}e^{iks}\mathcal F_0(k)P_+^{\ast}I_{\varphi_+^{\ast},\overline{a_+}}dk,
\label{eq:Lemma6.2chirsF0kWplusasteqiv}
\end{equation}
\begin{equation}
\chi_R(s)\int_{-\infty}^{0}e^{iks}\mathcal F_0(k)(W_-)^{\ast}dk \equiv 
\chi_R(s)\int_{-\infty}^{0}e^{iks}\mathcal F_0(k)P_-^{\ast}I_{\varphi_-^{\ast},\overline{a_-}}dk
\label{eq:Lemma6.2chiRsminusinfty0F0kWminusast}
\end{equation}
modulo regularizers of order $N$.
\end{lemma}
Proof.
We compute the first term of the right-hand side of (\ref{eq:RadonPlus}). Let $\chi_{\infty}(t)$ and $\chi(t)$ be as in (\ref{eq:cutoffSection5}). Modulo a regularizer, we can insert $\chi_{\infty}(|D_x|)$ between 
$\mathcal F_0(k)$ and $(W_+)^{\ast}$. 
Let $Q_0$ and $Q_{\infty}$ be defined by
$$
Q_0f(x) = (2\pi)^{-n}\int\!\!\!\int e^{i(x-y)\cdot\xi}\big(1 - \chi_{\infty}(|x_+(y,\xi)|)\big)\chi_{\infty}(|\xi|)f(y)dyd\xi,
$$
$$
Q_{\infty}f(x) = (2\pi)^{-n}\int\!\!\!\int e^{i(x-y)\cdot\xi}\chi_{\infty}(|x_+(y,\xi)|)\chi_{\infty}(|\xi|)f(y)dyd\xi,
$$
where $x_{\pm}(y,\xi)$ is the inverse function of $y = \nabla_{\xi}\varphi_{\pm}(x,\xi)$. Then we have
\begin{equation}
\chi_R(s)\int_0^{\infty}e^{iks}\mathcal F_0(k)Q_0fdk = 
\int_0^{\infty}\!\!\int_{{\bf R}^n}e^{ik(s - \omega\cdot y)}a(s,\omega,k,y)f(y)dydk,
\label{Sect6chiRsF0kQofdk}
\end{equation}
\begin{equation}
a(s,\omega,k,y) = \frac{\chi_R(s)}{\sqrt{2}(2\pi)^{n/2}}(-ik+0)^{(n-1)/2}(1 - \chi_{\infty}(|x_+(y,k\omega)|))\chi_{\infty}(k).
\nonumber
\end{equation}
 Since $|y| \leq 11R$  on the support of $a(s,\omega,k,y)$, the condition (\ref{eq:sminusomegayfrombelow}) is satisfied. Moreover by differentiating $y = \nabla_{\xi}\varphi_+(x,\xi)$, we have
$$
\left|\partial_k^m\partial_y^{\gamma} \, x_+(y,k\omega)\right| \leq C_{m\gamma}(1 + |k|)^{-m}, \quad \forall m \geq 1, \quad \forall \gamma,
$$
from which one can show that the condition (\ref{eq:Section6Estimatesasomegaky}) is also satisfied. Hence by Lemma 5.1, (\ref{Sect6chiRsF0kQofdk}) is a regularizer.

 Therefore we have only to consider
\begin{equation}
\chi_R(s)\int_0^{\infty}e^{iks}\big(\mathcal F_0Q_{\infty}(W_+)^{\ast}f\big)(k)dk.
\label{Sect6chiRsF0QinftyWplusast}
\end{equation}
We put $\chi_-(t) = 1 - \chi(t)$ and let $Q_-$ be defined by
\begin{equation}
Q_-f(x) = (2\pi)^{-n}\int\!\!\!\int e^{i(x-y)\cdot\xi}\chi_{\infty}(|x_+(y,\xi)|)\chi_{\infty}(|\xi|)\chi_-\big(\frac{x_+(y,\xi)}{|x_+(y,k\omega)|}\cdot\frac{\xi}{|\xi|}\big)f(y)dyd\xi.
\nonumber
\end{equation}
Then the operator (\ref{Sect6chiRsF0QinftyWplusast}) is split into two parts:
\begin{equation}
\chi_R(s)\int_0^{\infty}e^{iks}\big(\mathcal F_0P_+^{\ast}(W_+)^{\ast}f\big)(k)dk + 
\chi_R(s)\int_0^{\infty}e^{iks}\big(\mathcal F_0Q_-(W_+)^{\ast}f\big)(k)dk.
\nonumber
\end{equation}
The second term is rewritten as, up to a constant,
\begin{equation}
\chi_R(s)\int_0^{\infty}\!\!\!\int_{{\bf R}^n} e^{ik(s - \omega\cdot y)}
\chi_-\big(\frac{x_+(y,k\omega)}{|x_+(y,k\omega)|}\cdot\frac{k\omega}{|k\omega|}\big)\cdots dkdy,
\nonumber
\end{equation}
which is a regularizer by virtue of Lemma 5.1, since $s > 15R$ and $\omega\cdot y \leq - |y|/2$  on the support of the integrand. 
By Theorem 4.9, 
\begin{equation}
P_+^{\ast}(W_+)^{\ast} \equiv P_+^{\ast}I_{\varphi_+^{\ast},\overline{a_+}}
\nonumber
\end{equation}
modulo a regularizer of order $N$. We have thus proved (\ref{eq:Lemma6.2chirsF0kWplusasteqiv}).

Next we consider the second term of the right-hand side of (\ref{eq:RadonPlus}). We repeat the same arguments as above with $x_+(y,\xi)$ replaced by $x_-(y,\xi)$ and $\int_0^{\infty}\cdots dk$ by $\int_{-\infty}^0\cdots dk$.
Let $\chi_+(t) = 1 - \chi(- t)$ and $Q_+$ be defeined by
\begin{equation}
Q_+f(x) = (2\pi)^{-n}\int\!\!\!\int e^{i(x-y)\cdot\xi}\chi_{\infty}(|x_-(y,\xi)|)\chi_{\infty}(|\xi|)\chi_+\big(\frac{x_-(y,\xi)}{|x_-(y,k\omega)|}\cdot\frac{\xi}{|\xi|}\big)f(y)dyd\xi.
\nonumber
\end{equation}
Then as above, we are led to consider
\begin{equation}
\chi_{R}(s)\int_{-\infty}^0e^{iks}\big(\mathcal F_0P_-^{\ast}(W_-)^{\ast}f\big)(k)dk + 
\chi_{R}(s)\int_{-\infty}^0e^{iks}\big(\mathcal F_0Q_+(W_-)^{\ast}f\big)(k)dk
\nonumber
\end{equation}
modulo a regularizer. Since $k < 0$ this time, we have
\begin{equation}
\chi_+\big(\frac{x_-(y,k\omega)}{|x_-(y,k\omega)|}\cdot\frac{k\omega}{|k\omega|}\big) = \chi_+\big(- \frac{x_-(y,k\omega)}{|x_-(y,k\omega)|}\cdot\omega\big), 
\nonumber
\end{equation}
on which support, we have $\omega\cdot y \leq - |y|/2$. Therefore the second term is a regularizer. Again using Theorem 4.9, we have
\begin{equation}
P_-^{\ast}(W_-)^{\ast} \equiv P_-^{\ast}I_{\varphi_-^{\ast},\overline{a_-}}
\nonumber
\end{equation}
modulo a regularizer of order $N$. We have thus derived (\ref{eq:Lemma6.2chiRsminusinfty0F0kWminusast}) \qed

\bigskip
Let $(s)^{\alpha}_-$ be the homogeneous distribution defined in Chap.4, \S 5.


\begin{lemma}
Let $\chi_{\infty}(k)$ be as in (\ref{eq:cutoffSection5}), and put
\begin{equation}
D_j(s) = \frac{1}{2\pi}\int_{-\infty}^{\infty}e^{iks}\left(-ik+0\right)^{\frac{n-1}{2}-j}\chi_{\infty}(|k|)dk.
\label{eq:Dj(s)}
\end{equation}
Then we have
\begin{equation}
D_j(s) = \big(s\big)^{-\frac{n+1}{2}+j}_- + \Psi_0(s),
\nonumber
\end{equation}
where $\Psi_0(s)$ is a polynomially bounded smooth function on ${\bf R}$.
\end{lemma}
Proof. Letting $\psi_0(t)$ be the Fourier transform of $1 - \chi_{\infty}(|k|)$, we have
\begin{equation}
D_j(s) = \big(s\big)^{-\frac{n+1}{2}+j}_- - \frac{1}{\sqrt{2\pi}}\int_{-\infty}^{\infty}\big(s + t\big)_-^{-\frac{n+1}{2}+j}\psi_0(t)dt,
\nonumber
\end{equation}
from which the lemma follows immediately. \qed

\medskip

In the following we use the notation $\sim$ in the same meaning as in (\ref{eq:Sect4SymbolAsymptotic Expansion}). Namely
\begin{equation}
c(x,\xi) \sim \sum_{j=0}^{\infty}|\xi|^{-j}c_j(x,\widehat \xi)
\nonumber
\end{equation}
if and only if
\begin{equation}
\big|\partial_{\xi}^{\alpha}\partial_x^{\beta}\big(c(x,\xi) - \sum_{j=0}^{N-1}|\xi|^{-j}c_j(x,\widehat\xi)\big)\big| \leq C_{\alpha\beta N}|\xi|^{-N-|\alpha|}, \quad |\xi| > 1
\nonumber
\end{equation}
holds for any $\alpha, \beta$ and $N$. Note that this asymptotic expansion can be differentiated term by term with respect to $x$ and $\xi$.

By Theorem 3.2, we have for some $b_{\pm}(x,\xi) \in S^0$,
\begin{equation}
I_{\varphi_{\pm},a_{\pm}}P_{\pm} = I_{\varphi_{\pm},b_{\pm}}.
\label{eq:Bplusminus}
\end{equation}


\begin{lemma}
There exist $b_j(x,\theta) \ (j = 0, 1, 2, \cdots)$ such that $b_{\pm}(x,\xi)$ have the following asymptotic expansions as $|\xi| \to \infty$: 
\begin{equation}
b_{\pm}(x,\xi) \sim \sum_{j=0}^{\infty}(\pm|\xi|)^{-j}b_j(x,\pm\widehat\xi\,),
\label{eq:Basymp}
\end{equation}
\begin{equation}
b_0(x,\theta) = g(x)^{1/4}a_0(x,\theta)\chi_{\infty}(|x|)\chi(\widehat x\cdot\theta),
\label{eq:FormofB0pm}
\end{equation}
where $a_0(x,\theta)$ is given in (\ref{eq:A0}) and $\chi_{\infty}$, $\chi$ are given in (\ref{eq:cutoffSection5}).
\end{lemma}

Granting this lemma for the moment, we state the main theorem of this section.


\begin{theorem}
Let $\mathcal R_+(s,\theta,x)$ be the distribution kernel of $\mathcal R_+$. Then there exist $s_0 > 0$ such that  for any  $N > (n+1)/2$, the follolwing expansion holds for $s > s_0$:
$$
\mathcal R_+(s,\theta,x) = 
\sum_{j=0}^{N-1}
(s - \varphi(x,\theta))^{-\frac{n+1}{2}+j}_-r_j(x,\theta) + r^{(N)}(s,\theta,x),
$$
where $(s_0,\infty) \ni s \to r^{(N)}(s,\theta,x) \in \mathcal D'(S^{n-1}\times{\bf R}^n)$ is in $C^{\mu(N)}$, $\mu(N)$ is the greatest integer $\leq N - (n+1)/2$, $\varphi(x,\theta)$ is given by Definition 2.7, and
\begin{equation}
r_j(x,\theta)  = 2^{-1/2}(2\pi)^{(1-n)/2}i^{-j}\overline{b_j(x,\theta)},
\label{eq;rjxtheta}
\end{equation}
$b_j(x,\theta)$ being given in Lemma 5.5.
\end{theorem}
Proof. First let us note that
\begin{equation}
\varphi_-(x,k\theta) = k\varphi_+(x,\theta) \quad
{\rm for} \quad k < 0,
\label{eq:Phipm}
\end{equation}
\begin{equation}
b_-(x,k\theta) \sim \sum_{j=0}^{\infty}k^{-j}b_j(x,\theta) \quad
{\rm as} \quad k \to - \infty.
\label{eq:AsympBminus}
\end{equation}
In fact by Theorem 2.8 (3) we have for $k < 0$
\begin{equation}
\varphi_-(x,k\theta) = - \varphi_+(x,-k\theta) 
= - \varphi_+(x,|k|\theta)  = - |k|\varphi_+(x,\theta) = k\varphi_+(x,\theta),
\nonumber
\end{equation}
which proves (\ref{eq:Phipm}). By (\ref{eq:Basymp}) we have as $k \to - \infty$ $$
b_-(x,k\theta) \sim \sum_{j=0}^{\infty}(-|k|)^{-j}b_j\big(x,-\frac{k\theta}{|k\theta|}\big) = \sum_{j=0}^{\infty}k^{-j}b_j(x,\theta)
$$
which proves (\ref{eq:AsympBminus}).

Take $f \in C_0^{\infty}({\bf R}^n)$. Since $\varphi_+(x,\theta) = \varphi(x,\theta)$ by Definition 2.7, using (\ref{eq:Basymp}) we have as $k \to \infty$
\begin{eqnarray*}
& & \mathcal F_0(k)\big(I_{\varphi_+,b_+}\big)^{\ast}f \\
&=& \frac{1}{\sqrt2(2\pi)^{n/2}}\big(-ik + 0\big)^{(n-1)/2}
\int_{{\bf R}^n}e^{-i\varphi_+(x,k\theta)}
\overline{b_+(x,k\theta)}f(x)dx \\
&\sim& \frac{1}{\sqrt2(2\pi)^{n/2}}\sum_{j=0}^{\infty}
\int_{{\bf R}^n}e^{-ik\varphi(x,\theta)}\big(-ik + 0\big)^{\frac{n-1}{2} - j}
\chi_{\infty}(k)
i^{-j}\overline{b_j(x,\theta)}f(x)dx,
\end{eqnarray*}
where $\chi_{\infty}(k)$ is as in (\ref{eq:cutoffSection5}). Here we have used the fact that
\begin{equation}
(-ik+0)^{\alpha}(-ik)^m = (-ik + 0)^{\alpha+m} \quad {\rm if}\quad
0 \neq k \in {\bf R}, \quad \alpha \in {\bf R}, \quad m \in {\bf Z}.
\nonumber
\end{equation}
By (\ref{eq:Phipm}) and (\ref{eq:AsympBminus}), we have as  $k \to - \infty$
\begin{eqnarray*}
& & \mathcal F_0(k)\big(I_{\varphi_-,b_-}\big)^{\ast}f \\
&=& \frac{1}{\sqrt2(2\pi)^{n/2}}\big(-ik + 0\big)^{(n-1)/2}
\int_{{\bf R}^n}e^{-i\varphi_-(x,k\theta)}
\overline{b_-(x,k\theta)}f(x)dx \\
&\sim& \frac{e^{-(n-1)\pi i/4}}{\sqrt2(2\pi)^{n/2}}\sum_{j=0}^{\infty}
\int_{{\bf R}^n}e^{-ik\varphi(x,\theta)}\big(-ik + 0\big)^{\frac{n-1}{2} - j}
\chi_{\infty}(k)
i^{-j}\overline{b_j(x,\theta)}f(x)dx.
\end{eqnarray*}
Using (\ref{eq:RadonPlus}), Lemma 5.2 and (\ref{eq:Bplusminus}), we have
\begin{equation}
\begin{split}
 \chi_R(s)\mathcal R_+f(s) 
& \equiv \frac{\chi_R(s)}{\sqrt{2\pi}}
\int_{0}^{\infty}e^{iks}
\mathcal F_0(k)\left(I_{\varphi_+,b_+}\right)^{\ast}fdk \\
& + \frac{\chi_R(s)}{\sqrt{2\pi}}
\int_{-\infty}^0e^{iks}
\mathcal F_0(k)\left(I_{\varphi_-,b_-}\right)^{\ast}fdk
\end{split}
\nonumber
\end{equation}
modulo a regularizer of order $N$. We replace $\mathcal F_0(k)
\left(I_{\varphi_{\pm},b_{\pm}}\right)^{\ast}$ by the above asymptotic expansion to obtain 
\begin{equation}
\begin{split}
 \chi_R(s)\mathcal R_+f(s) 
& \equiv \frac{\chi_R(s)}{\sqrt{2}(2\pi)^{(n+1)/2}}
\sum_{j=0}^{N}
\int_{-\infty}^{\infty}\int_{{\bf R}^n}e^{ik(s - \varphi(x,\theta))}\\
& \hskip15mm \cdot
(-ik + 0)^{\frac{n-1}{2} - j}\chi_{\infty}(k)i^{-j}\overline{b_j(x,\theta)}f(x)dxdk
\end{split}
\nonumber
\end{equation}
modulo a term sufficiently regular in $s$. Performing the integral in $k$ and using Lemma 5.3, we have
$$
 \chi_R(s)\mathcal R_+f(s) 
 \equiv \frac{\chi_R(s)}{\sqrt{2}(2\pi)^{(n+1)/2}}
\sum_{j=0}^{N}
\int_{{\bf R}^n}\big(s - \varphi(x,\theta)\big)_-^{-\frac{n+1}{2}+j}i^{-j}\overline{b_j(x,\theta)}f(x)dx,
$$
modulo a term sufficiently regular in $s$, which proves the asymptotic expansion of $\mathcal R_+(s,\theta,x)$.
\qed

\bigskip
It remains to prove Lemma 5.4. Let $(\nabla_{\xi}\varphi_{\pm}\big)^{-1}(x,\xi)$  the inverse of the map $: x \to \nabla_{\xi}\varphi_{\pm}(x,\xi)$. Then by (\ref{eq:cutoffSection5}), the symbol $p_{\pm}(x,\xi)$ of $P_{\pm}$ is written as
\begin{equation}
p_{\pm}(x,\xi) = 
\chi_{\pm}\circ\big(\nabla_{\xi}\varphi_{\pm}\big)^{-1}(x,\xi).
\label{eq:PandChi}
\end{equation}
Now in view of (\ref{eq:AsympSymbPhiB}), we have
\begin{equation}
b_{\pm}(x,\eta) \sim \sum_{\alpha}\frac{1}{\alpha!}\;\partial_{\xi}^{\alpha}\left\{a_{\pm}(x,\xi)\big(D_x^{\alpha}p_{\pm}\big)\big(\widetilde\nabla_{\xi}\varphi_{\pm}(x,\xi,\eta),\eta\big) \right\}\Big|_{\xi=\eta}.
\label{eq:ExpandBpm1}
\end{equation}
Each term of the right-hand side consists of a sum of functions homogeneous in $\eta$. 
We rearrange them as
\begin{equation}
b_{\pm}(x,\eta) \sim \sum_{j=0}^{\infty}|\eta|^{-j}b_{\pm}^{(j)}(x,\widehat\eta),
\label{eq:ExpandBpm2}
\end{equation}
and compare (\ref{eq:ExpandBpm1}) and  (\ref{eq:ExpandBpm2}) to obtain
$$
b_{\pm}^{(0)}(x,\theta) = g(x)^{1/4}\chi_{\epsilon}(\pm x\cdot\theta)a_0(x,\pm \theta)p_{\pm}\big(\widetilde\nabla_{\xi}\varphi_{\pm}(x,\xi,\eta),\eta\big)\Big|_{\xi=\eta=\theta},
$$
where we have used (\ref{eq:AmplitudeWithK}).
Since
$$
\widetilde\nabla_{\xi}\varphi_{\pm}(x,\xi,\eta)\Big|_{\xi=\eta} = \big(\nabla_{\xi}\varphi_{\pm}\big)(x,\eta),
$$
we have by (\ref{eq:PandChi})
$$
p_{\pm}\big(\widetilde\nabla_{\xi}\varphi_{\pm}(x,\xi,\eta),\eta\big)\Big|_{\xi=\eta=\theta} = \chi_{\pm}(x,\theta),
$$
which proves (\ref{eq:FormofB0pm}).

To prove (\ref{eq:Basymp}), we make the following definition. Two functions $f_+(x,\xi)$ and $f_-(x,\xi)$ are said to be {\it compatible} if there exist $f_j(x,\theta) \ (j = 0,1,2,\cdots)$ such that $f_{\pm}(x,\xi)$ have the following asymptotic expansion as $|\xi| \to \infty$:
$$
f_{\pm}(x,\xi) \sim \sum_{j=0}^{\infty}\big(\pm|\xi|\big)^{-j}f_j(x,\pm\widehat\xi\,).
$$


\begin{lemma}
(1) If $f_+(x,\xi)$ and $f_-(x,\xi)$ are compatible, so are $\partial_{\xi}^{\alpha}f_+(x,\xi)$ and $\partial_{\xi}^{\alpha}f_-(x,\xi)$. \\
\noindent
(2) If $f_+(x,\xi)$ and $f_-(x,\xi)$ as well as $g_+(x,\xi)$ and $g_-(x,\xi)$ are compatible, so are $f_+(x,\xi)g_+(x,\xi)$ and $f_-(x,\xi)g_-(x,\xi)$. \\
\noindent
(3) $\partial_{\xi}^{\beta}\big(D_x^{\alpha}p_+\big)\big(\widetilde\nabla_{\xi}\varphi_+(x,\xi,\eta),\eta\big)\Big|_{\xi=\eta}$ and $\partial_{\xi}^{\beta}\big(D_x^{\alpha}p_-\big)\big(\widetilde\nabla_{\xi}\varphi_-(x,\xi,\eta),\eta\big)\Big|_{\xi=\eta}$ are compatible.
\end{lemma}
Proof. The assertions follow from a direct computation. In order to prove (1), we let $\partial_i = \partial/\partial\xi_i$ and take notice of
$$
\partial_if_+(x,\xi) \sim \sum_{m=0}^{\infty}|\xi|^{-m-1}\big\{-m\widehat\xi_i f_m(x,\widehat\xi\,) + \sum_{j=1}^n\Big(\partial_jf_m\Big)(x,\widehat\xi\,)(\delta_{ij} - \widehat\xi_i\widehat\xi_j)\Big\},
$$
$$
\partial_if_-(x,\xi) \sim \sum_{m=0}^{\infty}(-|\xi|)^{-m-1}\big\{m\widehat\xi_i f_m(x,-\widehat\xi\,) + \sum_{j=1}^n\Big(\partial_jf_m\Big)(x,-\widehat\xi\,)(\delta_{ij} - \widehat\xi_i\widehat\xi_j)\Big\}.
$$
The assertion (2) is obvious. To show (3), note that by Definition 3.7
$$
\partial_i\varphi_-(x,\xi) = - \widehat\xi_i\varphi(x,-\widehat\xi) + 
\sum_{j=1}^n\left(\frac{\partial\varphi}{\partial\xi_j}\right)(x,-\widehat\xi)(\delta_{ij} - \widehat\xi_i\widehat\xi_j) = 
\left(\partial_i\varphi_+\right)(x,-\xi).
$$
Since $\nabla_{\xi}\varphi_{\pm}$ are homogeneous of degree 0, this means that $\nabla_{\xi}\varphi_+(x,\xi)$ and $\nabla_{\xi}\varphi_-(x,\xi)$ are compatible. Since $D_x^{\alpha}p_+(x,\xi)$ and $D_x^{\alpha}p_-(x,\xi)$ are compatible, one can prove (3) inductively.
\qed

\medskip
By Lemma 5.6 and (\ref{eq:ExpandBpm1}), $b_{\pm}(x,\xi)$ are compatible. This proves Lemma 5.4. 


\subsection{Recovering partial regularities near infinity} 
Let us rewrite Theorem 5.6 in the operator form.
Let $D_j(s)$ and $r_j(x,\theta)$ be as in (\ref{eq:Dj(s)}) and (\ref{eq;rjxtheta}), rspectively. We put
\begin{equation}
\left(\mathcal R_+^{(j)}f\right)(s,\theta) = 
\int_{{\bf R}^n}
D_j(s - \varphi(x,\theta))
r_j(x,\theta)f(x)dx.
\nonumber
\end{equation}


\begin{lemma}
(1) For any $j, m \geq 0$, we have $
\mathcal R_+^{(j)} \in {\bf B}(H^m;\widehat H^{j+m})$.\\
\noindent
(2) Let $\chi_R(s)$ be as in Lemma 5.2. Then
for any $N$ 
$$
\chi_R(s)\mathcal R_+ \equiv \chi_R(s)\sum_{j=0}^{N-1}\mathcal R_{+}^{(j)}
$$
modulo a regularizer of order $N$. 
\end{lemma}

Proof. To prove the assertion (1), we have only to note that the operator
$$
\int_{{\bf R}^n}e^{-i\varphi(x,\xi)}\overline{r_j(x,\xi/|\xi|)}\chi_{\infty}(|\xi|)f(x)dx
$$
is $L^2$-bounded.
The assertion (2) has been proven in Theorem 5.5.  \qed

\medskip
The purpose of this section is to prove Lemma 1.13 in a localized form. Let us recall that the stationary phase method shows the scattered waves propagate to infinity along the directions close to $\widehat\xi = \pm \widehat x$. With this in mind, we prepare the following notion.


\begin{definition}
 For a constant $0 < \delta < 1$, let $S(\delta)$ be the set of symbols $p(x,\xi) \in S^{0}$ such that
${\rm supp}\,p \subset \{(x,\xi) \, ; \, |\widehat x\cdot\widehat\xi\,| < \delta\}$. 
We say that $f \in L^2({\bf R}^n)$ is {\it regular in  non-scattering region} if there exists $0 < \delta < 1$ such that $p(x,D_x)f \in H^{\infty}({\bf R}^n), \ \forall p(x,\xi) \in S(\delta)$.
\end{definition}

If $f$ is regular in non-scattering region, its wave front set, denoted by ${\rm WF}\,(f)$, satisfies ${\rm WF}\,(f)\cap\{|\widehat x\cdot\widehat\xi| < \delta\} = \emptyset$. As an example, let $B_R = \{x \in {\bf R}^n ; |x| < R\}$. If $f \in H^{\infty}(B_R)$ and $f(x) = 0$ for $|x| > R$, by the stationary phase method, $f(x)$ is shown to be regular in non-scattering region (see Lemma 6.8). The necessity of this notion will be made clear in the proof of Lemma 5.9.
We put $\widehat H^m(s > \sigma) = \widehat H^m(I_{\sigma})$ and $H^m(|x| > \rho) = H^m(B_{\rho}^c)$, where $I_{\sigma} = (\sigma,\infty)$ and $B_{\rho}^c = \{x \in {\bf R}^n ; |x| > \rho\}$.


\begin{lemma}
There exist constants $\rho > \sigma> 0$ such that the following assertion holds: If $f \in L^2({\bf R}^n)$ is regular in  non-scattering region and $\mathcal R_+^{(0)}f \in \widehat H^{m}(s>\sigma)$ for some $m \geq 0$, then $f \in H^{m}(|x|>\rho)$. Moreover $\rho $ can be chosen arbitrarily close to $\sigma$.
\end{lemma}
Proof. The proof is complicated and is split into several parts. Let $\chi(s) \in C^{\infty}({\bf R})$ be such that $\chi(s) = 1 \ (s > \sigma + 2)$, $\chi(s) = 0 \ (s < \sigma + 1)$, where $\sigma > 0$ will be determined later. We put
\begin{equation}
u(s,\theta) = \chi(s)\int_{-\infty}^{\infty}\!\int_{{\bf R}^n}
e^{ik(s - \varphi(x,\theta))} (-ik + 0)^{\frac{n-1}{2}}\chi_{\infty}(|k|)
r_0(x,\theta)f(x)dxdk,
\nonumber
\end{equation}
and assume that $u \in \widehat H^m$. We take $\psi_0(t), \psi_{\infty}(t) \in C^{\infty}({\bf R})$ such that $\psi_0(t) + \psi_{\infty}(t) = 1 \ (t \in {\bf R})$, $\psi_{\infty}(t) = 1 \ (t > 2)$, $\psi_{\infty}(t) = 0 \ (t < 1)$, and $c_0(t), c_1(t) \in C^{\infty}({\bf R})$ such that $c_0(t) + c_1(t) = 1 \ (t \in {\bf R})$, $c_1(t) = 1 \ (|t| > \delta/2)$, $c_1(t) = 0 \ (|t| < \delta/4)$, where $\delta$ is the constant appearing in the assumption of regularity in non-scattering region for $f$.
We split $f(x)$ into 3 parts :
$$
f(x) = \psi_{\infty}(|x|)c_1(\widehat x\cdot\theta)f(x) + \psi_0(|x|)f(x) +
\psi_{\infty}(|x|)c_0(\widehat x\cdot\theta)f(x).
$$
{\it 1st Step}. We put
\begin{equation}
\begin{split}
u_1(s,\theta) = \chi(s)\int_{-\infty}^{\infty}\int_{{\bf R}^n} e^{ik(s - \varphi(x,\theta))}(-ik + 0)^{\frac{n-1}{2}}
\chi_{\infty}(|k|)\\
\cdot r_0(x,\theta)
\psi_{\infty}(|x|)c_0(\widehat x\cdot\theta)f(x)dxdk,
\end{split}
\nonumber
\end{equation}
and show that $u_1 \in \widehat H^{\infty}$. This is proved if we show
\begin{equation}
v_1(x) := (2\pi)^{-n}\iint_{{\bf R}^n\times{\bf R}^n}e^{i(x\cdot\xi -\varphi(y,\xi))}\chi_{\infty}(|\xi|)r_0(y,\pm \widehat\xi\,)\psi_{\infty}(|y|)c_0(\pm \widehat y\cdot\widehat \xi\,)f(y)dyd\xi
\nonumber
\end{equation}
is in $H^{\infty}$. In view of (\ref{eq:PhiAstPhi}), we have
$$
w_1 := I_{\varphi,1}v_1 = Pf,
$$
where, modulo a regularizer, $P$ is a $\psi$DO whose symbol is supported in the region $\{|\widehat x\cdot\widehat\xi| < \delta\}$. Therefore $w_1 \in H^{\infty}$, since $f$ is regular in non-scattering region. Computing $I_{\varphi^{\ast},1}w_1$ and using (\ref{eq:AstPhiPhi}), we then have
$$
(1 + P_1 + P_2 + \cdots)v_1 = g,
$$
where $P_i \in S^{-i}$ and $g \in H^{\infty}$. By multiplying suitable $\psi$DO's, we have $v_1 \in H^{\infty}$.

\medskip
\noindent
{\it 2nd Step}.  Next we consider
\begin{equation}
\begin{split}
 \chi(s)\int_{-\infty}^{\infty}\int_{{\bf R}^n} e^{ik(s - \varphi(x,\theta))}(-ik + 0)^{\frac{n-1}{2}}
\chi_{\infty}(|k|)r_0(x,\theta)\\
\cdot \big[\psi_{\infty}(|x|)c_1(\widehat x\cdot\theta) + \psi_0(|x|)\big]f(x)dxdk.
\end{split}
\label{eq:ScattRegInt}
\end{equation}
Let $\widetilde\chi(s) \in C^{\infty}({\bf R})$ be such that $\widetilde\chi(s) = 1 \ (s > \sigma)$, $\widetilde\chi(s) = 0 \ (s < \sigma -1)$. By integration by parts, the operator
$$
\chi(s)\iint e^{ik(s - \varphi(x,\theta))}\left(1 - \widetilde\chi(\varphi(x,\theta))\right)\cdots dxdk
$$
is a regularizer. In fact, since $\varphi(x,\theta) < \sigma$, we have $|s - \varphi(x,\theta)| \geq C(s + |x|)$ for a constant $C > 0$ thanks to the factor $\psi_{\infty}(|x|)c_1(\widehat x\cdot\theta) + \psi_0(|x|)$. 

We are thus led to consider
\begin{equation}
\begin{split}
u_{2}(s,\theta) = & \chi(s)\int_{-\infty}^{\infty}\int_{{\bf R}^n} e^{ik(s - \varphi(x,\theta))}(-ik + 0)^{\frac{n-1}{2}}
\chi_{\infty}(|k|)\\
 & \cdot  \widetilde\chi(\varphi(x,\theta))r_0(x,\theta)
\big[\psi_{\infty}(|x|)c_1(\widehat x\cdot\theta) + \psi_0(|x|)\big]f(x)dxdk,
\end{split}
\nonumber
\end{equation}
which belongs to $\widehat H^{m}$. Here we choose $\sigma$ large enough so as to be able to apply Lemma 2.4, and make the change of variables $x \to (t,y) = (t(x,\theta), y(x,\theta))$. Since $t(x,\theta) = \varphi(x,\theta)$ by virtue of Lemma 2.5, the above integral is rewritten as 
\begin{equation}
\frac{1}{2\pi}\chi(s)\iint e^{ik(s-t)}q_0(t,k,y,\theta)\widetilde f(t,y,\theta)dkdtdy =: v_{2}(s,\theta),
\label{eq:vinfty}
\end{equation}
\begin{equation}
\begin{split}
q_0(t,k,y,\theta) = &2\pi (-ik + 0)^{\frac{n-1}{2}}\chi_{\infty}(|k|)\\
& \cdot \widetilde\chi(t)r_0(x,\theta)
\left[\psi_{\infty}(|x|)c_1(\widehat x\cdot\theta) + \psi_0(|x|)\right]J(t,y,\theta),
\end{split}
\label{eq:q0}
\end{equation}
$J(t,y,\theta)$ being the Jacobian of the map : $x \to (t,y)$, and in the expression of $q_0$, $x$ should be read as $x(t,y,\theta)$, $\widetilde f(t,y,\theta) = f(x)$. 
This reduces the problem to the 1-dimensional $\psi$DO calculus. 

Let $Q_0$ be the 1-dimensional $\psi$DO with symbol $\overline{q_0(t,k,y,\theta)}$, where $y, \theta$ are regarded as parameters. Then (\ref{eq:vinfty}) reads
$$
\int\chi(s)\left(Q_0^{\ast}\widetilde f(\cdot,y,\theta)\right)(s)dy = v_2(s,\theta),
$$
where $v_2 \in \widehat H^m$. By $\psi$DO calculus, we  have modulo $\widehat H^{m+1}$
\begin{equation}
\int\chi(s)\left(Q_0^{\ast}\widetilde f(\cdot,y,\theta)\right)(s)dy \equiv
\int\left(P_0^{\ast}\widetilde f(\cdot,y,\theta)\right)(s)dy 
 \in \widehat H^m,
 \label{eq:PAst}
\end{equation}
where the symbol of $P_0$ is the product of $\chi(t)$ and $q_0(t,k,y,\theta)$, namely, it is obtained with $\widetilde\chi(t)$ replaced by $\chi(t)$ in (\ref{eq:q0}). Passing to the Fourier transformation with respect to $s$ in (\ref{eq:PAst}), we get
\begin{equation}
\begin{split}
& \iint e^{-ikt}(-ik + 0)^{\frac{n-1}{2}}\chi_{\infty}(|k|)
\chi(t)r_0(x,\theta) \\
& \cdot \left[\psi_{\infty}(|x|)c_1(\widehat x\cdot\theta) + \psi_0(|x|)\right]J(t,y,\theta)\widetilde f(t,y,\theta)dtdy =: w(k,\theta),
\end{split}
\nonumber
\end{equation}
where $w(k,\theta)$ satisfies
$$
\int (1 + |k|)^{2m}\|w(k,\cdot)\|^2_{L^2(S^{n-1})}dk < \infty.
$$
Transforming back to the original variable $x$, we get
\begin{equation}
\begin{split}
& (-ik + 0)^{\frac{n-1}{2}}\chi_{\infty}(|k|)\int e^{-ik\varphi(x,\theta)}\chi(\varphi(x,\theta))r_0(x,\theta) \\
& \cdot \left[\psi_{\infty}(|x|)c_1(\widehat x\cdot\theta) + \psi_0(|x|)\right]f(x)dx = w(k,\theta).
\end{split}
\label{eq:IPhiAstf}
\end{equation}
We try to regard (\ref{eq:IPhiAstf}) as a FIO putting $\xi = k\theta$. 
Here we must note that the term $\chi(\varphi(x,\theta))$ behaves like
$$
|\partial_{\theta}^{\alpha}\chi(\varphi(x,\theta))| \leq C_{\alpha}(1 + |x|)^{|\alpha|},
$$
which seems to cause a trouble in defining a suitable class of symbols. However thanks to the locaization factor $\psi_{\infty}(|x|)c_1(\widehat x\cdot\theta) + \psi_0(|x|)$,
the amplitude $b(x,\theta)$ of (\ref{eq:IPhiAstf}) has the estimate
\begin{equation}
|\partial_{\theta}^{\alpha}\partial_x^{\beta}b(x,\theta)| \leq C_{\alpha\beta}(1 + |x|)^{-|\beta|}.
\nonumber
\end{equation}
In fact, by the estimate (\ref{eq:EstimatePsi}), on the support of $\chi'(\varphi(x,\theta))$, $|x\cdot\theta|$ is bounded. Due to the locaization factor $\psi_{\infty}(|x|)c_1(\widehat x\cdot\theta) + \psi_0(|x|)$, if $|x\cdot\theta|$ is bounded so is $x$. Therefore, the derivatives of $\chi(\varphi(x,\theta))$ does no harm to our analysis. This is the reason why we have introduced the notion of regularity in non-scattering region.

\medskip
\noindent
{\it 3rd Step}. 
We consider (\ref{eq:IPhiAstf}) separately in the region $k > 0$ and $k < 0$. For $\pm k > 0$, we put $k = \pm |\xi|$ and $\theta = \pm \widehat\xi$. Then we can rewrite (\ref{eq:IPhiAstf}) as
$$
\int e^{-i\varphi_{\pm}(x,\xi)}p_{\pm}(x,\xi)f(x)dx = g_{\pm}(\xi),
$$
where $p_{\pm}(x,\xi) \in S^0$ has its support in the region $\pm \widehat x\cdot\widehat\xi > \delta/3$ and $g_{\pm}(\xi)$ satisfies $(1 + |\xi|)^mg_{\pm}(\xi) \in L^2({\bf R}^n)$. We now mulitiply $e^{i\varphi_{\pm}(x,\xi)}$ and integrate in $\xi$. Then we have by FIO calculus
$$
q_{\pm}(x,D_x)\chi(|x|)f \in H^{m}, 
$$
where $q_{\pm}(x,\xi) \in S^0$, \ $q_{\pm}(x,\xi) = 1$ for $\pm \widehat x\cdot\widehat\xi > \delta$ and $|x| > 1$, $q_{\pm}(x,\xi) = 0$ for $\pm \widehat x\cdot\widehat\xi < \delta/5$ and $|x| > 1$, and $\chi(t) \in C^{\infty}({\bf R})$ such that $\chi(t) = 1$ \ $(t > \sigma + 2)$, $\chi(t) = 0$ \ $(t < \sigma + 1)$.
Taking into account that $f$ is regular in non-scattering region, we finally prove that $f \in H^{m}(|x|> \rho)$ for $\rho = s + 2$. By examining the proof, we see that $\rho$ can be chosen arbitrarily close to $\sigma$. \qed


\begin{theorem}
There exist $\rho > \sigma > 0$ such that if $f$ is regular in non-scattering region and $\mathcal R_+f \in 
\widehat H^m(s>\sigma)$ for some $m \geq 1$, then $f \in H^m(|x|>\rho)$. Moreover $\rho$ can be chosen arbitrarily close to $\sigma$.
\end{theorem}
Proof. If $\mathcal R_+f \in \widehat H^1(s>\sigma)$, we have $\mathcal R_+^{(0)}f \in \widehat H^1(s>\sigma)$ by Lemma 5.6 (1). Therefore the case $m = 1$ follows from Lemma 5.9. Let us assume the theorem when $m = k-1$. Then if $\mathcal R_+f \in \widehat H^k(s>\sigma)$, we have $f \in H^{k-1}(|x|>\rho)$. Therefore if $j \geq 1$, we have $\mathcal R_+^{(j)}f \in \widehat H^k(s>\sigma)$, which implies that $\mathcal R_+^{(0)}f \in \widehat H^k(s>\sigma)$. By Lemma 5.9, we have $f \in H^k(|x|>\rho)$, which completes the proof. \qed


\section{Singular support theorem}


\subsection{Envelope} Let us first recall the classical notion of envelope. Let $U$ and $\Omega$ be open sets in ${\bf R}^n$ and ${\bf R}^{n-1}$, respectively. Suppose a real-valued function $\phi(x,\omega) \in C^{\infty}(U\times\Omega)$ satisfies
\begin{equation}
\det\left(\nabla_x\phi,\frac{\partial}{\partial\omega_1}\nabla_x\phi,\cdots,
\frac{\partial}{\partial\omega_{n-1}}\nabla_x\phi\right) \neq 0, \quad
x \in U, \quad \omega \in \Omega,
\label{eq:Det1}
\end{equation}
\begin{equation}
\det\left(\frac{\partial^2\phi}{\partial\omega_i\partial\omega_j}\right)_{1\leq i, j \leq n-1} \neq 0, \quad
x \in U, \quad \omega \in \Omega.
\label{eq:Det2}
\end{equation}
Given an interval $I \subset {\bf R}$, we consider a family of surfaces
$$
\Sigma(s,\omega) = \left\{x \in U \; ; \; \phi(x,\omega) = s \right\}, \quad s \in I, \quad \omega \in \Omega.
$$
Assume that for $x \in U$ there exists a unique solution $\omega(x)$ to the system of equations
\begin{equation}
\frac{\partial\phi}{\partial\omega_1}(x,\omega) = \cdots = \frac{\partial\phi}{\partial\omega_{n-1}}(x,\omega)
= 0.
\label{eq:Phidiffomega}
\end{equation}
Then the envelope $\Sigma(s)$ of $\big\{\Sigma(s,\omega)\big\}_{\omega\in\Omega}$ is defined by
$$
\Sigma(s) = \{x \in U \; ; \; \phi(x,\omega(x)) = s\}.
$$
We put $y = (s,\omega)$ and $f(x,y) = \big(f_1(x,y),\cdots,f_n(x,y)\big)$, where$$
f_i(x,y) = \partial \phi(x,\omega)/\partial \omega_i, \quad (1 \leq i \leq n-1), \quad 
f_n(x,y) = \phi(x,\omega) - s.
$$
 Then the equation for the envelope  and the conditions (\ref{eq:Det1}), (\ref{eq:Det2}) are rewritten as
$$
f(x,y) = 0, \quad \det\left(\frac{\partial f}{\partial x}\right) \neq 0, \quad
\det\left(\frac{\partial f}{\partial y}\right) \neq 0.
$$
Hence by the implicit function theorem the map : $U \ni x \to y(x) = (s(x),\omega(x)) \in I\times\Omega$ is a diffeomorphism. Let $X(s,\omega)$ be its inverse.


\begin{lemma}
Let $g_{ij}(x)dx^idx^j$ be a Riemannian metric on $U$ and put $h(x,\xi) = \frac{1}{2}g^{ij}(x)\xi_i\xi_j$. Assume that $\phi(x,\omega)$ satisfies the eikonal equation
\begin{equation}
h(x,\nabla_x\phi(x,\omega)) = 1/2, \quad x \in U, \quad \omega \in \Omega.
\label{eq:EikonalSection7}
\end{equation}
(1) We put $\Phi(x) = \phi(x,\omega(x))$. Then $\Phi(x)$ also satisfies the eikonal equation
$$
h(x,\nabla_x\Phi(x))= 1/2, \quad x \in X.
$$
(2) Let $P(s,\omega) = \left(\nabla_x\Phi\right)(X(s,\omega))$. Then we have for $s \in I$ and $\omega \in \Omega$,
\begin{equation}
\left\{
\begin{split}
&\frac{\partial }{\partial s}X(s,\omega) = \left(\frac{\partial h}{\partial\xi}\right)(X(s,\omega),P(s,\omega)), \\
&\frac{\partial }{\partial s}P(s,\omega) = - \left(\frac{\partial h}{\partial x}\right)(X(s,\omega),P(s,\omega)). 
\end{split}
\right.
\label{eq:HamiJaco}
\end{equation}
\end{lemma}

Proof. By virtue of (\ref{eq:Phidiffomega}), we have
\begin{equation}
\nabla_x\Phi(x) = \left(\nabla_x\phi\right)(x,\omega(x)),
\label{eq:Phieqphi}
\end{equation}
which implies (1). We let $k(x,\omega) = \left(\nabla_x\phi\right)(x,\omega)$ and differentiate (\ref{eq:EikonalSection7}) by $\omega_j$ to have
\begin{equation}
\left(\frac{\partial h}{\partial\xi}\right)(x,k(x,\omega))\cdot
\frac{\partial k}{\partial\omega_j}(x,\omega) = 0, \quad
1 \leq j \leq n-1.
\nonumber
\end{equation}
Using (\ref{eq:Phieqphi}), we have $P(s,\omega) = k(X(s,\omega),\omega)$, hence
\begin{equation}
\left(\frac{\partial k}{\partial\omega_j}\right)(P(s,\omega),\omega)\cdot
\left(\frac{\partial h}{\partial\xi}\right)(X(s,\omega),P(s,\omega)) = 0, 
\quad 1 \leq j \leq n-1.
\label{eq:hkEq1}
\end{equation}
On the other hand, we have by differentiating $\big(\partial\phi/\partial\omega_j\big)(X(s,\omega),\omega) = 0$ by $s$
\begin{equation}
\left(\frac{\partial k}{\partial\omega_j}\right)(X(s,\omega),\omega)\cdot
\frac{\partial X}{\partial s}(s,\omega) = 0, \quad 
1 \leq j \leq n-1.
\label{eq:hkEq2}
\end{equation}
By (\ref{eq:Det1}), $\partial k/\partial\omega_1, \cdots, \partial k/\partial\omega_{n-1}$ are linearly independent. Therefore by (\ref{eq:hkEq1}) and 
(\ref{eq:hkEq2}) we have
$$
\frac{\partial X}{\partial s}(s,\omega) = \lambda(s,\omega)
\left(\frac{\partial h}{\partial \xi}\right)(X(s,\omega),P(s,\omega))
$$
for some scalar function $\lambda(s,\omega)$. Differentiating $s = \phi(X(s,\omega),\omega)$ with respect to $s$, we then have
$$
1 = k\cdot\frac{\partial X}{\partial s} = 
\lambda k\cdot\left(\frac{\partial h}{\partial\xi}\right)(X,k) = 2\lambda h(X,P)= \lambda.
$$
Finally by differentiating $P_i(s,\omega) = \left(\partial\phi/\partial x_i\right)(X(s,\omega),\omega)$ we have
\begin{eqnarray*}
\frac{\partial}{\partial s}P_i(s,\omega) &=& 
\sum_j\left(\frac{\partial^2\phi}{\partial x_i\partial x_j}\right)
(X(s,\omega),\omega)\frac{\partial X_j}{\partial s}(s,\omega) \\
&=& \left(\frac{\partial k}{\partial x_i}\right)
(X(s,\omega),\omega)\cdot \left(\frac{\partial h}{\partial \xi}\right)(X(s,\omega),P(s,\omega)) \\
&=& - \left(\frac{\partial h}{\partial x_i}\right)(X(s,\omega),P(s,\omega)),
\end{eqnarray*}
since by differentiating $h(x,k(x,\omega)) = 1/2$, we get
$$
\left(\frac{\partial h}{\partial x_i}\right)(x,k(x,\omega)) + 
\left(\frac{\partial h}{\partial\xi}\right)(x,k(x,\omega))\cdot
\frac{\partial k}{\partial x_i}(x,\omega) = 0.
\qed
$$
Let us note that by (\ref{eq:Phieqphi}), $\Sigma(s,\omega)$ is tangent to $\Sigma(s)$ at $X(s,\omega)$.

\medskip
We now put
\begin{equation}
\Sigma^{(\pm)}(s,\theta) = \left\{x \in {\bf R}^n ; \varphi_{\pm}(x,\theta) = s \right\},
\nonumber
\end{equation}
and construct the envelope of $\big\{\Sigma^{(\pm)}(s,\theta)\big\}_{\theta \in S^{n-1}}$. Since $\varphi_+(x,\theta) = - \varphi_-(x,-\theta)$ by Theorem 2.8 (3), we have
\begin{equation}
\Sigma^{(+)}(s,\theta) = \Sigma^{(-)}(-s,-\theta).
\nonumber
\end{equation}
Therefore we have only to consider $\varphi_+(x,\theta) = \varphi(x,\theta)$. 
For $\varphi(x,\theta)$, the assumptions (\ref{eq:Det1}), (\ref{eq:Det2}) are satisfied on the region $\{|x| > r_0\}\times S^{n-1}$, where $r_0 > 0$ is chosen largre enogh. We consider the equation
\begin{equation}
\nabla_{\theta}\varphi(x,\theta) = 0, \quad x\cdot\theta > 0,
\label{eq:NablaPhi}
\end{equation} 
 $\nabla_{\theta}$ being the gradient on $S^{n-1}$, which corresponds to (\ref{eq:Phidiffomega}). If $\varphi(x,\theta) = x\cdot\theta$, the solution is unique and given by $\theta = \widehat x$. Since
$\partial_x^{\alpha}(\varphi(x,\theta) - x\cdot\theta) = O(|x|^{-|\alpha|-\epsilon_0})$, we see that (\ref{eq:NablaPhi}) has a unique solution $\theta(x) = \widehat x + O(|x|^{-\epsilon_0})$.  Let $s(x) = \varphi(x,\theta(x))$ and $X(s,\theta)$ be the inverse of the map : $x \to (s(x),\theta(x))$. We summarize the properties of these diffeomorphisms in the following theorem. We put $\Sigma(s,\theta) = \Sigma^{(+)}(s,\theta)$.


\begin{theorem}
There exist $r_0 > 0$ and $s_0 > 0$ for which the following assertions hold. \\
\noindent
(1) For any $x \in {\bf R}^n$ such that $|x| > r_0$, there exists a unique 
$\theta(x) \in S^{n-1}$ satsifying $\big(\nabla_{\theta}\varphi\big)(x,\theta(x)) = 0$ and $\theta(x)\cdot x > 0$.
We define
\begin{equation}
\Phi(x) = \varphi(x,\theta(x)) \quad 
{\rm for} \quad |x| > r_0,
\nonumber
\end{equation}
and extend it smoothly for $|x| \leq r_0$ so that $\Phi(x)$ is monotone increasing with respect to $|x|$. Then $\Phi(x) \sim  |x|$ as $|x| \to \infty$ and satisfies the eikonal equation
$$
g^{ij}(x)(\partial_i\Phi(x))(\partial_j\Phi(x)) = 1
\quad {\rm for} \quad |x| > r_0.
$$
(2) For any $s > s_0$, the set
\begin{equation}
\Sigma(s) = \{x \in {\bf R}^n ; \Phi(x) = s\}
\nonumber
\end{equation}
is a strictly convex  compact hypersurface. \\
\noindent
(3) For any $ s > s_0$ and $x \in \Sigma(s)$, $\Sigma(s)$ is tangent to  $\Sigma(s,\theta(x))$ at $x$. Moreover $\theta(x)$ is a unique point $\theta$ in $S^{n-1}$ such that $\Sigma(s)$ is tangent to $\Sigma(s,\theta)$ at $x$. We also have for $|x| > r_0$
\begin{equation}
\max_{\theta\in S^{n-1}}\varphi(x,\theta) = \Phi(x), 
\label{eq:maxphiPhi}
\end{equation}
and the maximum is attained if and only if $\theta = \theta(x)$.
\\
\noindent
(4) For any $ s > s_0$ and $\theta \in S^{n-1}$, there exists a unique $X(s,\theta) \in \Sigma(s)$ such that $\Sigma(s,\theta)$ is tangent to $\Sigma(s)$ at $X(s,\theta)$. We also have for any $\theta \in S^{n-1}$
\begin{equation}
\max_{x\in\Sigma(s)}\varphi(x,\theta) = s = \Phi(X(s,\theta)), 
\label{eq:MaxphisPhi}
\end{equation}
and the maximum  is attained if and only if $x = X(s,\theta)$.
\\
\noindent
(5) For any $ s > s_0$, the map
\begin{equation}
S^{n-1} \ni \theta \to X(s,\theta) \in \Sigma(s)
\nonumber
\end{equation}
is a diffeomorphism and its inverse is given by
\begin{equation}
\Sigma(s) \ni x \to \theta(x) \in S^{n-1}.
\nonumber
\end{equation}
(6) The map
\begin{equation}
 X : (s_0,\infty)\times S^{n-1} \ni (s,\theta) \to X(s,\theta) \in {\bf R}^n
 \nonumber
\end{equation}
is a diffeomorphism whose image contains the region $\{x \, ;\, |x| > r_0\}$.
The inverse of this map is 
\begin{equation}
X^{-1} : x \to \big(\Phi(x),\theta(x)\big).
\nonumber
\end{equation}
It has the following estimates $(\widehat x = x/|x|)$
\begin{equation}
|\partial_x^{\alpha}(\Phi(x) - |x|)| \leq C_{\alpha}(1 + |x|)^{-\epsilon_0-|\alpha|}, \quad \forall \alpha,
\label{eq:EstimatesPHI}
\end{equation}
\begin{equation}
|\partial_x^{\alpha}(\theta(x) - \widehat x)| \leq C_{\alpha}(1 + |x|)^{-1-\epsilon_0-|\alpha|}, \quad \forall \alpha.
\label{eq:EstimatesTHETA}
\end{equation}
(7) The diffeomorphism $X^{-1}$ gives the geodesic polar coordinates in a neighborhood of infinity, and in this coordinate system the Riemannian metric 
$G = g_{ij}(x)dx^idx^j$ takes the following form
\begin{equation}
X^{\ast}G = (ds)^2 + \sum_{i,j=1}^{n-1}h_{ij}(s,\theta)d\theta^id\theta^j.
\nonumber
\end{equation}
\end{theorem}
Proof.
As is noted above $\varphi(x,\theta) = x\cdot\theta$ for the Euclidean metric, hence $\theta(x) = \widehat x$, $\Phi(x) =  |x|$, and the theorem is obvious. The assertion (1) follows from Lemma 6.1. Since $\Sigma(s)$ is a slight perturbation of sphere, (2) follows. The first part of the assertion (3) is obvious. We shall prove (\ref{eq:maxphiPhi}). If $\varphi(x,\theta)$ attains its maximum at $\theta$, $(\nabla_{\theta}\varphi)(x,\theta) = 0$ holds. This equation has two solutions $\widetilde\theta_{\pm}$ such that $\pm x\cdot\widetilde\theta_{\pm} > 0$. The Hessian matrix of $\varphi(x,\theta)$ at $\widetilde\theta_+ \ (\widetilde\theta_-)$ is negative (positive) definite. Hence the maximum is attained at $\widetilde\theta_+$, furthermore, $\widetilde\theta_+ = \theta(x)$. The first part of (4) is obvious. At the point $x$ where $\varphi(x,\theta)$ attains its maximum on $\Sigma(s)$, $\nabla_x\Phi(x)$ and $\nabla_x\varphi(x,\theta)$ are propotional. This is just the point on which two surfaces $\Sigma(s)$ and $\Sigma(s,\theta)$ are tangent each other, hence (\ref{eq:MaxphisPhi}) holds.
The mapping properties in (5) and (6) are clear.
From the equation $\nabla_{\theta}\varphi(x,\theta) = 0$, we get $\nabla_{\theta}\widehat x\cdot\theta =  O(|x|^{-1-\epsilon_0})$, 
from which (\ref{eq:EstimatesTHETA}) follows. The estimate (\ref{eq:EstimatesPHI}) then follows from Theorem 2.8 (1).
Let us prove (7). By the equation (\ref{eq:HamiJaco}), $X(s,\theta)$ is a geodesic. Hence $(s(x),\theta(x))$ are geodesic polar coordinates. We put $\overline{x}^i = \theta_i(x) \ (1 \leq i \leq n-1)$, 
$\overline{x}^n = \Phi(x)$. Then the associated Riemannian metric $\overline{g}_{ij}$ is computed as follows :
$$
\overline{g}^{nn} = g^{ij}\frac{\partial\overline{x}^n}{\partial x^i} 
\frac{\partial\overline{x}^n}{\partial x^j} = 
g^{ij}\left(\partial_i\Phi\right)\left(\partial_j\Phi\right)
= 1, 
$$
$$
\overline{g}^{nk} = g^{ij}\frac{\partial\overline{x}^n}{\partial x^i} 
\frac{\partial\overline{x}^k}{\partial x^j} = 
g^{ij}\left(\partial_i\Phi\right)\left(\partial_j\theta_{ k}\right)
= 0, 
$$
for $1 \leq k \leq n-1$. Here we have used the equation (7.5) and 
$$
0 = \frac{\partial\theta_k}{\partial s} = 
\frac{\partial\theta_k}{\partial x^m}\frac{\partial X^m}{\partial s} 
= \left(\partial_m\theta_k\right)g^{im}P_i = 0.
$$
This proves (7). \qed


\begin{cor}
For large $|x|$, we have $\varphi(x,\theta) \leq \Phi(x)$,
and the equality holds if and only if $\theta = \theta(x)$, equivalently, $x = X(s,\theta)$ for some $ s > s_0$.
\end{cor}


\subsection{Singularity expansions II}
Our next aim is to compute an asymptotic expansion around $s = \sigma$ of the integral (coupling of distribution and test function, actually)
\begin{equation}
\int_{{\bf R}^n} (s - \varphi(x,\theta))^{\alpha}_-(\sigma - \Phi(x))^{\beta}_+
f(x)dx, \quad f \in C_0^{\infty}({\bf R}^n).
\label{eq:intalphabeta}
\end{equation} 
For any $\theta \in S^{n-1}$, we have constructed a bicharacteristic $x(t,y,\theta), \, p(t,y,\theta)$ having the properties in Lemma 2.2. We use the variables $t,y $ to calculate (\ref{eq:intalphabeta}), which is possible by virtue of Lemma 2.4. In perfoming the computation below it will be helpful to recall that for the Euclidean metric $\sum_{i=1}^n(dx^i)^2$, $x(t,y,\theta) = t\theta + y$, $\theta\cdot y = 0$, $\varphi(x,\theta) = x\cdot\theta$ and $\Phi(x) = |x|$.

Let
$\widetilde\Phi(t,y,\theta) = \Phi(x(t,y,\theta))$. Then since $t = \varphi(x,\theta)$ by Lemma 2.5 
we have by Corollary 6.3
$$
\widetilde\Phi(t,y,\theta) - t = \Phi(x) - \varphi(x,\theta) \geq 0,
$$
and for a fixed $t$ the last equality holds only at one point, which we denote by $y(t,\theta)$. At $y(t,\theta)$ the surface $t = \Phi(x)$ is tangent to the surafce $t = \varphi(x,\theta)$. Therefore  $(t,y(t,\theta))$ is the coordiante of $X(t,\theta)$ given in Theorem 6.2 (4). By the Taylor expansion with respect to $y$ we have
$$
\widetilde\Phi(t,y,\theta) - t = \frac{1}{2}\langle A(y - y(t,\theta)),y - y(t,\theta)\rangle + O(|y - y(t,\theta)|^3),
$$
as $y \to y(t,\theta)$,
where 
$$
A = A(t,\theta) = \left(\frac{\partial^2\widetilde\Phi}{\partial y_i\partial y_j}
(t,y(t,\theta),\theta)\right)
$$
 is a positive definite matrix and $\langle \,\,,\,\rangle$ is the Euclidean inner product of ${\bf R}^{n-1}$. By the Morse lemma, one can find a function $z = z(t,y,\theta)$ defined in a neighborhood of $y(t,\theta)$ such that
$$
\Phi(x) = \widetilde\Phi(t,y,\theta) = t + \frac{1}{2}\langle A(t,\theta)z,z\rangle,
$$
and $z = y - y(t,\theta) + O(|y - y(t,\theta)|^2)$.
We now make a new change of variables: $x \to (t,z)$ and put $\widetilde f(t,z,\theta) = f(x)$. We denote by 
$$
J_P(t,z,\theta) = |\det\left(\partial x/\partial(t,z)\right)|
$$
the associated Jacobian. (Here the subscript 
$P$ menas that we are using the plane wave like characteristic surface $t = \varphi(x,\theta)$). 
 Then we have
\begin{equation}
\begin{split}
&  \displaystyle\int\big(s - \varphi(x,\theta)\big)_-^{\alpha}\big(\sigma - \Phi(x)\big)^{\beta}_+f(x)dx \\
& = \displaystyle \iint \big(s - t\big)_-^{\alpha}\Big(\sigma - t - \frac{1}{2}
\langle A(t,\theta)z,z\rangle\Big)^{\beta}_+\widetilde f(t,z,\theta)J_P(t,z,\theta)dtdz.
\end{split}
\label{eq:Changevari}
\end{equation}

We say that $g(s,\theta)$ admits the asymptotic expansion
\begin{equation}
g(s,\theta) \sim \sum_{k=0}^{\infty}(\sigma - s)^{\lambda+k}_+g_k(\theta), \quad g_k \in C^{\infty}(S^{n-1})
\nonumber
\end{equation}
around $s= \sigma$, if there exists $\epsilon_0 > 0$ with the following property. For any $N > 0$, there exist $G_N(s,\theta), H_N(s,\theta) \in C^{\infty}({\bf R};L^2(S^{n-1}))$ such that
\begin{equation}
g(s,\theta) = \sum_{k=0}^{N-1}(\sigma - s)^{\lambda+k}_+g_k(\theta) + 
(\sigma - s)^{\lambda+N}_+G_N(s,\theta) + H_N(s,\theta)
\nonumber
\end{equation}
holds for $|s - \sigma| < \epsilon_0$. Similarly, we say that $f(x)$ admits the asymptotic expansion
\begin{equation}
f(x) \sim \sum_{k=0}^{\infty}(\sigma - \Phi(x))^{\lambda+k}_+f_k(\theta), \quad f_k(\theta) \in C^{\infty}(\Sigma(\sigma)) 
\nonumber
\end{equation}
around $\Phi(x) = \sigma$, where $\Sigma(\sigma) = \{\sigma = \Phi(x)\}$ and $\theta$ denotes the local coordinate on $\Sigma(\sigma)$, if there exists $\epsilon_0 > 0$ with the following property. For any $N > 0$, there exist $G_N(x), H_N(x) \in C^{\infty}({\bf R}^n)$ such that
\begin{equation}
f(x) = \sum_{k=0}^{N-1}(\sigma - \Phi(x))^{\lambda+k}_+f_k(\theta) + 
(\sigma - \Phi(x))^{\lambda+N}_+G_N(x) + H_N(x)
\nonumber
\end{equation}
holds when $|\Phi(x) - \sigma| < \epsilon_0$.


\begin{lemma}
Let $g(t,z) \in C^{\infty}_0({\bf R}\times{\bf R}^{n-1})$, and $\sigma > 0$ be a sufficiently large constant. Then if $\beta > -1$, we have the following asymptotic expansion around $s = \sigma$
\begin{equation}
\begin{split}
& \iint\big(s - t\big)^{\alpha}_-\left(\sigma - t - \frac{1}{2}\langle A(t,\theta)z,z\rangle\right)^{\beta}_+g(t,z)dtdz \\
& \sim \sum_{k=0}^{\infty}(\sigma - s)_+^{\alpha + \beta + \frac{n+1}{2} + k}
\left(P^{(\alpha,\beta)}_kg\right)(\sigma,0),
\end{split}
\label{eq:SingExpa1}
\end{equation}
where $P^{(\alpha,\beta)}_k$ is a differential operator having the following form
\begin{equation} 
P^{(\alpha,\beta)}_k = 
\sum_{\begin{subarray}{c}
m + |\gamma|/2 \, \leq \, k, \\|\gamma|= {\rm even}
\end{subarray}}
C_{km\gamma}(\alpha,\beta)p_{km\gamma}(\sigma,\theta)\partial_t^m\partial_z^{\gamma}.
\label{eq:CoeffDiff1}
\end{equation} 
 If $|\gamma| = m = k = 0$, we have
\begin{equation}
 C_{000}{(\alpha,\beta)}p_{000}(\sigma,\theta) = (2\pi)^{\frac{n-1}{2}}\det A(\sigma,\theta)^{-1/2}.
\label{eq:SingExpa2}
\end{equation}
\end{lemma}
Proof. First let us note that the left-hand side of (\ref{eq:SingExpa1}) vanishes if $s > \sigma$. For $s < \sigma$, we put $\epsilon = \sigma - s$, $s - t = \epsilon \rho$, $z = \sqrt{2\epsilon(1+ \rho)}A(t,\theta)^{-1/2}w$ and
\begin{equation}
\begin{split}
g_{\epsilon}(\rho,w) = g\big(\sigma - \epsilon(1 + \rho),\sqrt{2\epsilon(1 + \rho)}A(\sigma - \epsilon(1 + \rho),\theta)^{-1/2}w\big) \\ \cdot
\det A(\sigma - \epsilon(1 + \rho),\theta)^{-1/2}.
\end{split}
\nonumber
\end{equation}
Note that since $\sigma \geq t + \frac{1}{2}\langle Az,z\rangle \geq t$, we have $\sigma - t = \epsilon(1 + \rho) \geq 0$.
Then the left-hand side of (\ref{eq:SingExpa1}) is rewritten as
\begin{equation}
\begin{split}
& 2^{\frac{n-1}{2}}\frac{\Gamma(\beta + \frac{n+1}{2})}{\Gamma(\beta + 1)}\epsilon^{\alpha + \beta + \frac{n+1}{2}}\\
& \times \int_{-1}^0\int_{|w|<1} (\rho)^{\alpha}_-(1 + \rho)^{\beta + \frac{n-1}{2}}_+(1 - |w|^2)^{\beta}\, g_{\epsilon}(\rho,w)\,d\rho \,dw.
\label{eq:intrhow}
\end{split}
\end{equation}
Since $A(t,\theta)$ is a positive definite matrix and smooth in $t$, so is $A(t,\theta)^{-1/2}$. This follows from the well-known Dunford-Taylor integral of bounded operators (see e.g. p. 44 of \cite{Ka76}). 
We put $\delta = \sqrt{\epsilon(1 + \rho)}$ and expand $g_{\epsilon}(\rho,w)$ into  a Taylor series with respect to $\delta$ to see that each term of the expansion consists of the product of a function of $\sigma, \theta$ and 
\begin{equation}
\delta^{2p + |\gamma|}w^{\gamma}\left(\partial_s^{m}\partial_z^{\gamma}g\right)(\sigma,0), \quad m \leq p.
\label{eq:TermsOfgepsilon}
\end{equation}
In fact, we first expand $g\big(\sigma - \delta^2,\delta y\big)$ to obtain terms like $\delta^{2m+|\gamma|}y^{\gamma}(\partial_s^m\partial_z^{\gamma}g)(\sigma,0)$, and next expand $y = \sqrt2A(\sigma-\delta^2,\theta)^{-1/2}w$ and $\det A(\sigma - \delta^2,\theta)^{-1/2}$ to have (\ref{eq:TermsOfgepsilon}).  We 
replace $g_{\epsilon}(\rho,w)$ in (\ref{eq:intrhow}) by this asymptotic expansion. If $|\gamma|$ is odd, $\int(1 - |w|^2)^{\beta}w^{\gamma}dw = 0$. Therefore, letting $k = p + |\gamma|/2$ and rearranging the terms, we obtain (\ref{eq:SingExpa1}).  To compute (\ref{eq:SingExpa2}), we have only to use (\ref{eq;ProductHomogDis}) and the formula
\begin{equation}
\int_{|w|<1}(1 - |w|^2)^{\beta}dw = \pi^{\frac{n-1}{2}}
\frac{\Gamma(\beta+1)}{\Gamma(\beta + \frac{n+1}{2})}. 
\nonumber
\end{equation}
Here we have assumed $\beta > -1$ to guarantee the convergence of the integral \qed


\begin{lemma}
Let $\sigma > 0$ be sufficiently large, and assume that $\beta > -1$.  Then for any $f(x) \in C_0^{\infty}({\bf R}^n)$, we have the following asymptotic expansion around $s =  \sigma$:
\begin{equation}
 \int\big(s - \varphi(x,\theta)\big)^{\alpha}_-\big(\sigma - \Phi(x)\big)^{\beta}_+f(x)dx 
\sim 
\sum_{k=0}^{\infty}\big(\sigma - s\big)_+^{\alpha + \beta + \frac{n+1}{2} + k}
g_{k}^{(\alpha,\beta)}(\sigma,\theta).
\label{eq:SingExpa3}
\end{equation}
Each term of the expansion (\ref{eq:SingExpa3}) is represented by a differential operator $M^{(\alpha,\beta)}_k$ on ${\bf R}\times S^{n-1}$ in the following way: 
$$
g_k^{(\alpha,\beta)}(\sigma,\theta) = \left(M_k^{(\alpha,\beta)}f\circ X\right)(\sigma,\theta),
$$
where $X(s,\theta)$ is the diffeomorphism in Theorem 6.2 (6). In the local coordinates $M_k^{(\alpha,\beta)}$ has the following expression
\begin{equation}
M_k^{(\alpha,\beta)} = \sum_{j+|\gamma|/2\leq k}C_{kj\gamma}(\alpha,\beta)m_{kj\gamma}(s,\theta)\partial_s^j\partial_{\theta}^{\gamma}.
\label{eq:MkAlphaBeta}
\end{equation}
In particular, 
\begin{equation}
M_0^{(\alpha,\beta)} = 
(2\pi)^{\frac{n-1}{2}}\det\left(A(\sigma,\theta)\right)^{-1/2}J_P(\sigma,0,\theta).
\label{eq:TopTerm0}
\end{equation}
\end{lemma}
Proof. We plug (\ref{eq:Changevari}) with (\ref{eq:SingExpa1}). Let $X : (s,\theta) \to X(s,\theta)$ be the diffeomorphism in Theorem 6.2 (6). In the $(t,y)$ coordinate system employed to derive (\ref{eq:Changevari}), the condition $z = 0$ and $t = \sigma$ means that $y = y(\sigma,\theta)$ and $\varphi(x(\sigma,y,\theta),\theta) = \sigma$, which represents the point $X(\sigma,\theta)$. Therefore each term of the asymptotic expansion (\ref{eq:SingExpa3}) is a derivative of $f(x)$ evaluated at $x = X(\sigma,\theta)$. Moreover
\begin{eqnarray*}
\partial_t\Big|_{t=s,y=y(s,\theta)} &=& 
\sum_{i,j=1}^ng^{ij}(X(s,\theta))\left(\frac{\partial\varphi}{\partial x_j}\right)(X(s,\theta),\theta)\frac{\partial}{\partial x_i} \\
&=& 
\sum_{i,j=1}^ng^{ij}(X(s,\theta))\left(\frac{\partial\Phi}{\partial x_j}\right)(X(s,\theta))\frac{\partial}{\partial x_i},
\end{eqnarray*}
which is equal to $\partial_s$ in the coordinate system $(s,\theta) = X^{-1}(x)$.
Thus we have the asymptotic expansion (\ref{eq:SingExpa3}). The formulas (\ref{eq:MkAlphaBeta}) and (\ref{eq:TopTerm0}) follow from ({\ref{eq:CoeffDiff1}) and (\ref{eq:SingExpa2}). \qed

\medskip
The first term $M_0^{(\alpha,\beta)}$ is written by geometric quantities. By a simple computation one can show that
\begin{equation}
\left(\det A(\sigma,\theta)\right)^{-1/2} = |\nabla_x\Phi(x)|^{-(n-1)/2}\left(\det {\mathcal H}_{PS}\big(\frac{\partial x}{\partial y_i},\frac{\partial x}{\partial y_j}\big)\right)^{-1/2}\Big|_{x = X(\sigma,\theta)},
\nonumber
\end{equation}
\begin{equation}
{\mathcal H}_{PS} = {\mathcal H}_P - {\mathcal H}_S,
\nonumber
\end{equation}
where $x = x(t,y,\theta)$,  ${\mathcal H}_P$ and ${\mathcal H}_S$ are second fundamental forms on $\{\sigma = \varphi(x,\theta)\}$ and $\{\sigma = \Phi(x)\}$ induced from the Euclidean metric, and
\begin{equation}
J_P(\sigma,0,\theta) = \big|G(x)^{-1}\nabla_x\Phi(x)\big|\left(\det G_S(x)\right)^{1/2}\Big|_{x = X(\sigma,\theta)},
\nonumber
\end{equation}
where $G(x) = \big(g_{ij}(x)\big)$, and $G_S(x)$ is the matrix of first fundamental form on $\{\sigma = \Phi(x)\}$ induced from the Euclidean metric. 


\begin{theorem}
Let $\sigma > 0$ be sufficiently large and $\lambda > - 1/2$. Then for any $f \in C_0^{\infty}({\bf R}^n)$, we have the following asymptotic expansion around $s = \sigma$
\begin{equation}
\left(\mathcal R_+(\sigma - \Phi(x))^{\lambda}_+f\right)(s,\theta) \sim
\sum_{k=0}^{\infty}(\sigma - s)^{\lambda + k}_+g_{k}^{(\lambda)}(\sigma,\theta).\nonumber
\end{equation}
\end{theorem}
Proof. This follows from Theorem 5.5 and Lemma 6.5. Note that $(\sigma - \Phi(x))^{\lambda}_+f \in L^2({\bf R}^n)$ if $\lambda > -1/2$. \qed

\medskip
In order to prove the converse of Theorem 6.6, we expand $(\sigma - \Phi(x))^{\lambda}_+f(x)$ into an asymptotic series $\sum_{k=0}^{\infty}(\sigma - \Phi(x))^{\lambda + k}_+f_k(x)$ and study the relations between $f_k$ and $g_k$. We compute in the following way. For $f(x) \in C_0^{\infty}({\bf R}^n)$, take $\chi(x) \in C_0^{\infty}({\bf R}^n)$ such that $\chi(x) = 1$ on ${\rm supp}\,f$.
 Then by Taylor expansion
\begin{equation}
(\sigma - \Phi(x))^{\lambda}_+f(x) = \sum_{j=0}^N(\sigma - \Phi(x))^{\lambda+j}_+f_j^{(\sigma)}\chi(x) + F_N(x),
\nonumber
\end{equation}
where $f_j^{(\sigma)}$ is a smooth function on $\{\sigma = \Phi(x)\}$ and $F_N(x)$ is a compactly supported $C^{\mu(N)}$-function, where $\mu(N) \to \infty$ as $N \to \infty$. This implies modulo $C^{\mu(N)}$-function
\begin{equation}
\left(\mathcal R_+\big((\sigma - \Phi(x))^{\lambda}_+f(x)\big)\right)(s,\theta) \equiv 
\sum_{j=0}^N\left(\mathcal R_+\big((\sigma - \Phi(x))^{\lambda+j}_+f_j^{(\sigma)}\chi(x)\big)\right)(s,\theta),
\nonumber
\end{equation}
and up to a smooth function the right-hand side is equal to
\begin{equation}
\sum_{i,j}\int (s - \varphi(x))_-^{-\frac{n+1}{2}+i}(\sigma - \Phi(x))^{\lambda+j}_+r_if_j^{(\sigma)}dx
\nonumber
\end{equation}
near $s = \sigma$, since $\chi(x) \equiv 1$ near $\{\sigma = \Phi(x)\}$. Omitting the cut-off function $\chi(x)$, we express this computation as 
\begin{equation}
\left(\mathcal R_+\big((\sigma - \Phi(x))^{\lambda}_+f(x)\big)\right)(s,\theta) \sim
\sum_{j=0}^{\infty}\left(\mathcal R_+\big((\sigma - \Phi(x))^{\lambda+j}_+f_j^{(\sigma)}\big)\right)(s,\theta),
\nonumber
\end{equation}
which will not give a confusion.

 In order to write down the expansion it is convenient to use the diffeomorphism $X(s,\theta)$ in Theorem 6.2 (6). We insert the asymptotic expansion
\begin{equation}
\big((\sigma - \Phi(x))^{\lambda}_+f\circ X\big)(s,\theta) \sim \sum_{k=0}^{\infty}(\sigma - s)_+^{\lambda + k}f_k(\sigma,\theta)
\nonumber
\end{equation}
into the formula in Theorem 6.6 and obtain
\begin{equation}
\Big(\mathcal R_+\Big(\sum_{k=0}^{\infty}(\sigma - \Phi(x))^{\lambda+k}_+f_k^{\ast}\Big)\Big)(\tau,\theta) \sim
\sum_{k=0}^{\infty}(\sigma - \tau)^{\lambda + k}_+g_k(\lambda,\sigma,\theta),
\nonumber
\end{equation}
where $f_k^{\ast} = f_k\circ X^{-1}$. Note that we fix $\sigma$ and regard $f_k^{\ast}$ as a function on $\{\sigma = \Phi(x)\}$. Let us look at $g_k(\lambda,\sigma,\theta)$ more precisely. Using Theorem 5.5 and Lemma 6.5, we have
\begin{equation}
\begin{split}
& \Big({\mathcal R}_+\Big(\sum_{\alpha=0}^{\infty}(\sigma - \Phi)^{\lambda+\alpha}_+f_{\alpha}^{\ast}\Big)\Big)(\tau,\theta) \\
& \sim
\sum_{k=0}^{\infty}(\sigma-\tau)^{\lambda+k}_+\sum_{\alpha+\beta+\gamma=k}
g_{\gamma}^{(-\frac{n+1}{2}+\beta,\lambda+\alpha)}(\sigma,\theta),
\end{split}
\nonumber
\end{equation}
\begin{equation}
g_{\gamma}^{(-\frac{n+1}{2}+\beta,\lambda+\alpha)}(\sigma,\theta) = M_{\gamma}^{(-\frac{n+1}{2}+\beta,\lambda+\alpha)}r_{\beta}f_{\alpha}^{\ast}\circ X.
\nonumber
\end{equation}
Therefore we have
\begin{equation}
g_k(\lambda,\sigma,\theta) = 
\sum_{\alpha=0}^k\left(\sum_{\beta+\gamma=k-\alpha}M_{\gamma}^{(-\frac{n+1}{2}+\beta,\lambda+\alpha)}r_{\beta}\right)f_{\alpha}^{\ast}\circ X.
\nonumber
\end{equation}
 Hence we have the following formula
\begin{equation}
\begin{split}
& g_k(\lambda,\sigma,\theta)  = P_{0}^{(k)}(\lambda)f_k(\sigma,\theta) + P_{2}^{(k-1)}(\lambda)f_{k-1}(\sigma,\theta) \\
 & \hskip 20mm + \cdots 
 +  P_{2k}^{(0)}(\lambda)f_0(\sigma,\theta),
 \end{split}
\label{eq:fntogn}
\end{equation}
where $P_{2(k-j)}^{(j)}(\lambda)$ is a differential operator with respect to $\theta$, and $P_0^{(k)}$ is the operator of multiplication by
\begin{equation}
P_0^{(k)}(\sigma,\theta) = (2\pi)^{\frac{n-1}{2}}
\det A(\sigma,\theta)^{-1/2}J_P(\sigma,0,\theta)r_0(X(\sigma,\theta),\theta).
\label{eq:NonzeroCond1}
\end{equation}
Using (\ref{eq:NonzeroCond1}), one can solve (\ref{eq:fntogn}) with respect to $f_j$ to have
\begin{equation}
\begin{split}
& f_k(\lambda,\sigma,\theta)  = Q_{0}^{(k)}(\lambda)g_k(\sigma,\theta) + Q_{2}^{(k-1)}(\lambda)g_{k-1}(\sigma,\theta) \\
 & \hskip 20mm + \cdots 
 +  Q_{2k}^{(0)}(\lambda)g_0(\sigma,\theta),
 \end{split}
\label{eq:gntofn}
\end{equation}
where $Q_{2(k-j)}^{(j)}(\lambda)$ is a differential operator with respect to $\theta$, and
\begin{equation}
Q_0^{(k)}(\sigma,\theta) = 1/P_0^{(k)}(\sigma,\theta).
\nonumber
\end{equation}


\begin{theorem}
Let $\sigma > 0$ be sufficiently large and $\lambda > - 1/2$. Given any $g(s,\theta)$ having the following asymptotic expansion around $s = \sigma$
\begin{equation}
g(s,\theta) \sim \sum_{k=0}^{\infty}(\sigma - s)^{\lambda + k}_+g_k(\theta)
\nonumber
\end{equation}
with $g_k(\theta) \in C^{\infty}(S^{n-1})$, there exists $f(x)$ such that around $s = \sigma$
\begin{equation}
\left(\mathcal R_+f\right)(s,\theta) \sim  \sum_{k=0}^{\infty}(\sigma - s)^{\lambda + k}_+g_k(\theta),
\nonumber
\end{equation}
and $f(x)$ admits the asymptotic expansion
\begin{equation}
f(x) \sim \sum_{k=0}^{\infty}(\sigma - \Phi(x))^{\lambda + k}_+f_k(\theta)
\label{eq:AsympExpaf}
\end{equation}
around $\Sigma(\sigma)$, $\theta$ being the local coordinates on $\Sigma(\sigma)$. Furthermore
\begin{equation}
g_0(\theta) = 
N(\sigma,\theta)f_0(X(\sigma,\theta)),
\nonumber
\end{equation}
$N(\sigma,\theta)$ being given by (\ref{eq:NonzeroCond1}). This $f(x)$ is unique in the sense that if there exist two such $f^{(1)}(x)$ and $f^{(2)}(x)$, $f^{(1)}(x) - f^{(2)}(x)$ is smooth. In particular, $f^{(1)}(x)$ and $f^{(2)}(x)$ have the asymptotic expansion as in (\ref{eq:AsympExpaf}) with the same $f_k(\theta)$.
\end{theorem}
Proof. By (\ref{eq:gntofn}), one can construct $f_k(\theta)$. Using Borel's procedure we then construct $f(x)$ having the asymptotic expansion $f(x) \sim 
\sum_{k=0}^{\infty}(\sigma - \Phi(x))_+^{\lambda+k}f_k(\theta)$. Suppose there exist two such $f^{(1)}$ and $f^{(2)}$. As is seen by the lemma below, $f^{(1)} - f^{(2)}$ is regular in non-scattering region, hence it is in $H^{\infty}$ by Theorem 5.11. \qed


\begin{lemma}
For $\sigma > 0$ large enough, let $u(x) = (\sigma - \Phi(x))^{\mu}_+f(x)$, where $f(x) \in C^{\infty}({\bf R}^n)$ wwhose support is sufficiently close to $\{\sigma = \Phi(x)\}$, and $\mu > -1/2$. Then $u(x)$ is regular in non scattering region.
\end{lemma}
Proof. Let $P$ be the $\psi$DO with symbol $p(x,\xi) \in S^0$ such that for some $0 < \delta < 1$, ${\rm supp}\,p(x,\xi) \subset \{|\widehat x\cdot\widehat\xi| < \delta\}$. Then by using the polar coordinates $(s,\theta)$ in Theorem 6.2 (6), 
\begin{equation}
\begin{split}
\widehat{Pu}(\xi) & = (2\pi)^{-n/2}\int_{\Phi(x) < \sigma}e^{-ix\cdot\xi}
\overline{p(x,\xi)}u(x)dx \\
& = \int_0^{\sigma}\int_{S^{n-1}}e^{-iX(s,\theta)\cdot\xi}(\sigma - s)^{\mu}p(X(x,\theta),\xi)g(s,\theta)dsd\theta,
\end{split}
\nonumber
\end{equation}
with suitable $g(s,\theta) \in C^{\infty}$. We apply the stationary phase method (as $|\xi| \to \infty)$ to the integral on $S^{n-1}$. Since $X(s,\theta)$ is close to $s\theta$, the critical points are close to $\pm \widehat \xi$, on which $p(X(s,\theta),\xi)$ vanishes. Therefore above integral is rapidly decreasing in $\xi$. \qed


\subsection{Singular support theorem}
The following Theorem 6.10  will elucidate how the modified Radon transform describes the propagation of singularities for the wave equation. 


\begin{definition} Assume $\Sigma(t) \subset \{|x| > r_0\}$.
A function $f(x) \in L^2({\bf R}^n)$ is said to be piecewise $H^{\infty}(|x|>r_0)$ with interface $\Sigma(t)$ if there exist $f_1, f_2 \in H^{\infty}(|x|>r_0)$ such that $f = \big(t - \Phi(x)\big)^0_+f_1 + \big(t - \Phi(x)\big)^0_-f_2$ on $|x|>r_0$. 
Similarly a function $f(s) \in L^2({\bf R};L^2(S^{n-1}))$ is said to be piecewise $\widehat H^{\infty}(s>s_0)$ with interface $s = t \ (> s_0)$ if there exist $f_1, f_2 \in \widehat H^{\infty}(s>s_0)$ such that $f = (t - s)^0_+f_1 + (t - s)^0_-f_2$ for $s > s_0$. 
\end{definition}


\begin{theorem}
Pick $r_0, s_0 > 0$ large enough, and let $t > $ max $\{r_0 + 1, s_0 +1\}$. Assume that $f \in L^2({\bf R}^n)$ is regular in non-scattering region. Then $f$ is piecewise $H^{\infty}(|x|>r_0)$ with interface $\Sigma(t)$ if and only if $\mathcal R_+ f$ is piecewise $\widehat H^{\infty}(s>s_0)$ with interface $s = t$.
\end{theorem}
Proof. Suppose $f$ is piecewise $H^{\infty}(|x|>r_0)$ with interface $\Sigma(t)$. Up to an $H^{\infty}$-function, $f$ is equal to $(t - \Phi(x))^0_+\widetilde f(x)$ with $\widetilde f \in H^{\infty}({\bf R}^n)$. By Theorem 5.5, $(\mathcal R_+f)(s,\theta)$ is smooth with respect to $s$ if $s \neq t$. By Theorem 6.6, $({\mathcal R}_+f)(s,\theta) \sim \sum_{k\geq0}(t - s)^k_+g_k(\theta)$ around $s = t$. Therefore $\mathcal R_+ f$ is piecewise $\widehat H^{\infty}(s>s_0)$ with interface $s = t$.

Conversely, suppose $\mathcal R_+ f$ is piecewise $\widehat H^{\infty}(s>s_0)$ with interface $s = t$. Up to an $\widehat H^{\infty}$-function, $(\mathcal R_+f)(s,\theta) \equiv (t - s)^0_+g(s,\theta)$ with $g \in \widehat H^{\infty}(s > s_0)$. By Theorem 6.7, there exists $\widetilde f$ such that $({\mathcal R}_+\widetilde f)(s,\theta) \sim (t - s)^0_+g(s,\theta)$ around $s = t$. Then $\mathcal R_+(f - \widetilde f) \in \widehat H^{\infty}(s > s_0)$. By Theorem 5.10, $f - \widetilde f \in H^{\infty}(|x|>r_0)$. This shows that $f$ is piecewise $H^{\infty}(|x|>r_0)$ with interface $\Sigma(t)$. \qed

\bigskip
The meaning of Theorem 6.10 in propagation of singularities is as follows.  We put $v(t,s) = \big(\mathcal R_+ \partial_t u(t)\big)(s)$  for the solution $u(t)$ to the wave equation $\partial_t^2u = Hu$ with initial data $u(0) = 0, \ \partial_tu(0) = f$. Then $v(t,s)$ solves the 1-dimensional wave equation
\begin{equation}
\left\{
\begin{split}
&(\partial_t^2 - \partial_s^2)v(t,s) = 0, \\
& v(0,s) = \big(\mathcal R_+f\big)(s), \quad \partial_tv(0,s) = 0,
\end{split}
\right.
\nonumber
\end{equation}
hence is written as
$$
v(t,s) = \frac{1}{2}\Big((\mathcal R_+f)(s + t) + (\mathcal R_+f)(s - t)\Big).
$$
If $\sigma$ is sufficiently large, $t \geq 0$ and $f$ is regular in non-scattering region, we then see that $f$ is piecewise $H^m(|x|>r_0)$ with interface $\Sigma(\sigma)$ if and only if  $(\mathcal R_+\partial_tu(t))(s)$ is piecewise $\widehat H^m(s>s_0)$ with interface $s = t + \sigma$, which is equivalent to that $\partial_tu(t)$ is piecewise $H^m(|x|>t+r_0)$ with interface $\Sigma(t+\sigma)$.

\end{document}